UCLouvain, Faculté des sciences
Institut de Recherche en Mathématique et Physique

# Odd Khovanov homology, higher representation theory and higher rewriting theory

Léo Schelstraete

Thesis submitted in partial fulfillment
of the requirements for the degree of
*Docteur en Sciences*

Thesis committee:

| | |
|---|---|
| Prof. Ben Elias | University of Oregon |
| Prof. Marino Gran *(President)* | UCLouvain |
| Prof. Pascal Lambrechts | UCLouvain |
| Prof. Louis-Hadrien Robert | Université Clermont-Auvergne |
| Prof. Tim Van der Linden *(Secretary)* | UCLouvain |
| Prof. Pedro Vaz *(Supervisor)* | UCLouvain |

Odd Khovanov homology, higher representation theory
and higher rewriting theory

by Léo Schelstraete

This work was supported by the Fonds de la Recherche Scientifique - FNRS
under the Aspirant Fellowship FC 38559.

# Remerciements

Une thèse est un objet étrange: un long projet solitaire qu'on ne peut mener à bien seul. Pour garder la cadence, on a autant besoin de meneurs d'allure que d'un public en folie. À l'occasion de ce dernier kilomètre, je voudrais remercier celles et ceux qui ont fait de ces quatre années un savoureux mélange d'émerveillement et de doux instants.

Comme tout, cette aventure commence par hasard: un groupe de travail sur les groupes quantiques, un étudiant de première master qui s'y invite, et des particules qui, en quelques coups de craie, deviennent des ficelles qui s'emmêlent. Pedro, cela fait déjà cinq ans et demi que tu m'as insuflé la passion pour la topologie quantique. Depuis, tu as été un mentor investi, disponible, passionné, humain. Je ne compte plus la myriade de discussions que nous avons eu, et ce même lors des débuts confinés. Ta porte toujours ouverte a été un luxe inestimable. Grâce à ta confiance sans cesse réitérée, je me suis construit en tant que mathématicien indépendant. Merci pour tout, cette thèse n'existerait pas sans toi.

Loulou, tu m'as suivi ces quatre années en chassé-croisé, au gré des conférences. Merci pour tes nombreux conseils et tes commentaires assidus sur ma thèse, même sur la partie réécriture.

Ben, thank you for your warm welcome in Eugene and in your home. Your interest and enthusiasm for rewriting theory gave me the courage and confidence for the laborious months that followed, carefully writing down the theory.

Marino, Pascal, Tim, merci d'avoir accepté de faire partie de mon jury de thèse, mais aussi de suivre mon parcours depuis neuf ans maintenant.

Je ne ferai pas des maths sans sa communauté et celles et ceux qui la font. Merci à tous les grands frères de l'IRMP: Fathi, François, Geoffrey et Jacques. Vous m'avez appris autant de mathématiques que ce qui fait la vie d'un mathématicien. Merci aux invétérés des pauses midi, berceaux des ragots et des micro-révolutions; à Carine, Cathy, Élodie et Martine, pour leur aide précieuse; à toutes les équipes de recherche, qui m'ont accepté quand je m'incrustais dans leurs groupes de travail; à Benoît, pour toutes ces discussions mathématico-politiques; à Pierre, pour ta passion contagieuse.

Merci aussi à toute la communauté de topologie de basses dimensions française. Au fil des années, les conférences ressemblaient de plus à plus à des retrouvailles entre vieux amis; et les workshops "Between The Waves", à des réunions familiales. Merci Delphine, Jules, Marco, Léo et Loulou pour l'organisation: c'était toujours avec impatience que je descendais de ma petite Belgique, retrouver tout le monde. Une pensée particulière pour Yohan, et nos discussions qui se teintaient si vite de métaphysique. Merci aussi à mes deux grands frères quantiques, Jules et Marco, pour les rires et les conseils.

Merci enfin à celles et ceux qui ont rendu ces quatre années inoubliables. Avant tout, merci à la coloc des Pommiers, Cam, Ed, Seb et Vic. J'ai une pensée déjà nostalique de ces premiers mois confinés et leur gym matinale, des soirées films-fondants, des jeux sociétés, des repas du soir et leurs histoires, des randonnées, des apéros dans le potager. Merci pour ces moments paisibles, de rire et de tendresse, de camaraderie et d'entraide.

Ma thèse a été parsemée de divers projets musicaux. Étrangement, je me rappelle mieux des répétitions que des représentations… Merci à tous les musicos, et en particulier à Yann, la fondamentale de tous ces morceaux. Merci aussi aux amis du LSO: que de séances, que de guindailles, sont passées à s'amuser autour d'un verre!

Sandrine, merci d'être dans ma vie: les mots sont trop fades pour encrer ce qui a déjà été dit.

Enfin, merci à Papa, Maman et Emma. Votre amour et votre soutien continuel m'ont armé d'une certitude lucide, de celle qui permet de tout réussir, tout en définissant soi-même ce que cela veut dire.

Over the course of this four-year thesis, my professional travel has resulted in the emission of 5.27 tons of $CO_2$–0.23 tons from train travel and 5.04 tons from air travel. This averages out to 1.32 tons per year. For context, to align with the Paris Agreement goals, average individual total $CO_2$ emissions should be reduced to 2 tons per year by 2030 (see the IPCC report for a more precise statement). Sources: *ecopassenger* for train consumption and *myclimate.org* for plane consumption.

# Contents

# I Odd Khovanov homology and higher representation theory

# II Rewriting theory and higher algebra

# Introduction to part I

## *odd Khovanov homology and higher representation theory*

*Quantum topology* studies the interplay between low-dimensional topology and representation theory; *rewriting theory* studies algebraic structures with an algorithmic perspective. This thesis is devoted to both fields, and their surprising interconnections.

The common thread is the realisation that certain commutative theories admit *odd* (or *super*) analogues, where the commutative rule only holds up to signs. In topology, the last decades have seen the emergence of powerful homology theories of knots, starting with the pioneer work of Khovanov and his celebrated *Khovanov homology.* Their study, together with their connections to the representation theory of quantum groups, form the field of modern quantum topology. Surprisingly, Khovanov homology admits an odd analogue, known as *odd Khovanov homology.* While the former builds upon polynomial algebras, the latter has anti-commutative features. This led to the discovery of similar odd analogues in related fields. Despite these successes, a grasp on odd Khovanov homology has remained elusive. In part I of this thesis, we relate odd Khovanov homology to representation theory, leading to a natural framework for its understanding.

The main characters in our construction are certain singular surfaces called *super* $\mathfrak{gl}_2$*-foams.* Their structural properties are governed by what is known as a *super-2-category*, the categorical analogue of a super algebra. Unfortunately, the increased complexity of working with such structures exhausts the standard toolkit of the higher algebraist. Recently however, techniques from rewriting theory have started to emerge. In part II of this thesis, we developed the theory of *linear Gray rewriting modulo*, a general set of tools to study presented weakly commutative 2-categories. We then apply these methods to give a basis theorem for super $\mathfrak{gl}_2$-foams.

This thesis is organised such that both parts can be read independently, including distinct introductions. Then, the dependency is as follows:

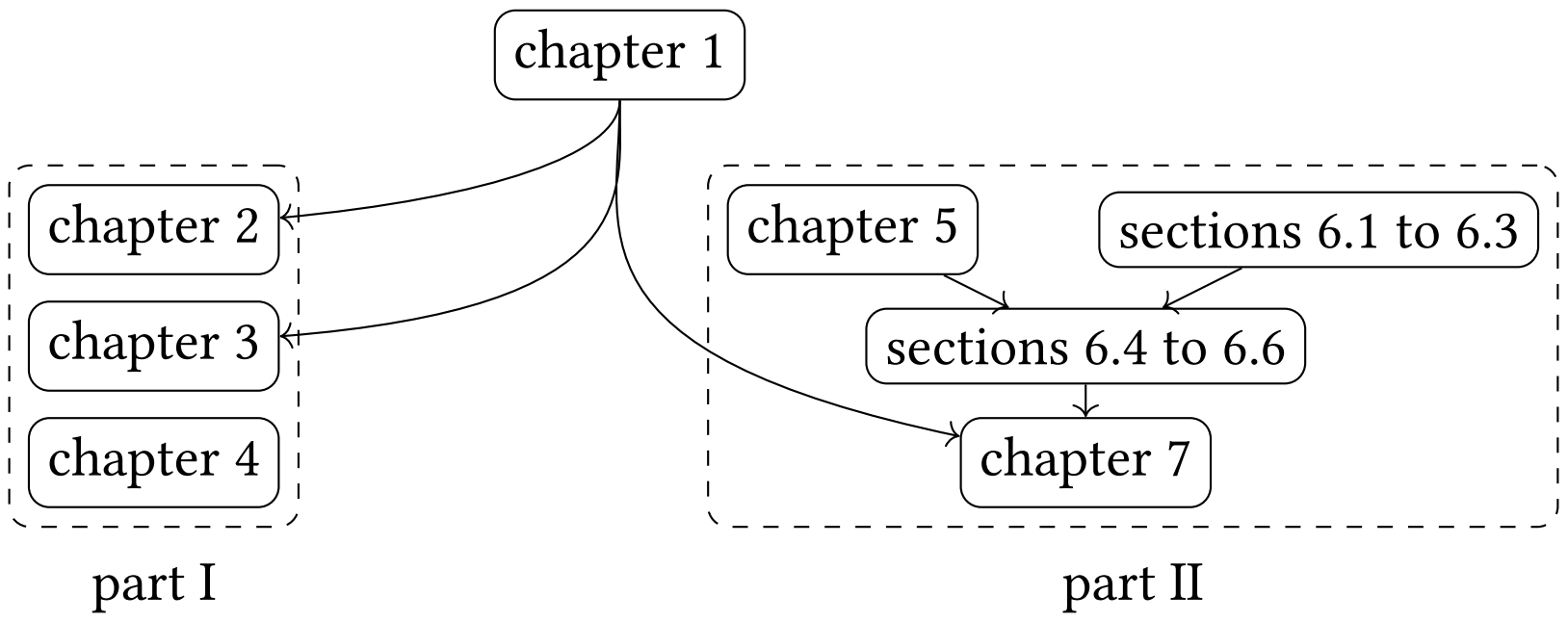

This is only a reading dependency; mathematically, the main results of chapter 2 depend on chapter 4, and the main results of part I depend on part II.

## i.1 Overview of the field

### i.1.1 Quantum topology

The field of quantum topology was initiated in 1984, when Vaughan Jones discovered a new (Laurent) polynomial invariant of links[1], nowadays known as the *Jones polynomial* [94]. This was followed by rapid developments connecting low-dimensional topology with representation theory, culminating in a general method to construct invariants of links from algebraic data, called the *Reshetikhin–Turaev construction* [160, 161]:

$$\left\{ \begin{array}{c} \text{representation} \\ \text{theory of } U_q(\mathfrak{g}) \end{array} \right\} \xrightarrow{\substack{\text{Reshetikhin–Turaev} \\ \text{construction}}} \left\{ \begin{array}{c} \text{invariants} \\ \text{of links} \end{array} \right\}$$

Here $U_q(\mathfrak{g})$ is a certain algebra associated to a simple finite-dimensional complex Lie algebra $\mathfrak{g}$, called a *quantum group*; more details below. Invariants of links arising from the Reshetikhin–Turaev construction are called *quantum*

[1] A knot is an embedding of the circle in $\mathbb{R}^3$ regarded up to ambient isotopies, while a link generalizes to an embedding of an arbitrary union of circles. In this thesis, everything that holds for knots also holds for links.

*polynomial invariants.* The Jones polynomial is realised as the simplest quantum polynomial invariant, associated to the quantum group $U_q(\mathfrak{sl}_2)$, where $\mathfrak{sl}_2$ is the Lie algebra of traceless 2-by-2 matrices.

From this breakthrough arose a new categorical perspective on knot theory. Intuitively, an invariant of links is *local* if it can be computed by chopping off the link in small pieces, assigning algebraic data to these pieces, and glueing them back to obtain the invariant of the whole link. Formally, a piece of link, or *tangle*, is a union of embedded circles and intervals in $\mathbb{R}^2 \times [0,1]$, regarded up to boundary-preserving isotopies. Furthermore, each endpoint of an interval must either lie on the bottom line $\mathbb{R} \times \{0\} \times \{0\}$ or the top line $\mathbb{R} \times \{0\} \times \{1\}$; see Fig. i.1 on the left-hand side, where a link is decomposed into tangles. To formalize the glueing, tangles are gathered into a monoidal category[2] $(Tang, \sqcup, \emptyset)$, with composition and monoidal product $\otimes = \sqcup$ pictured below:

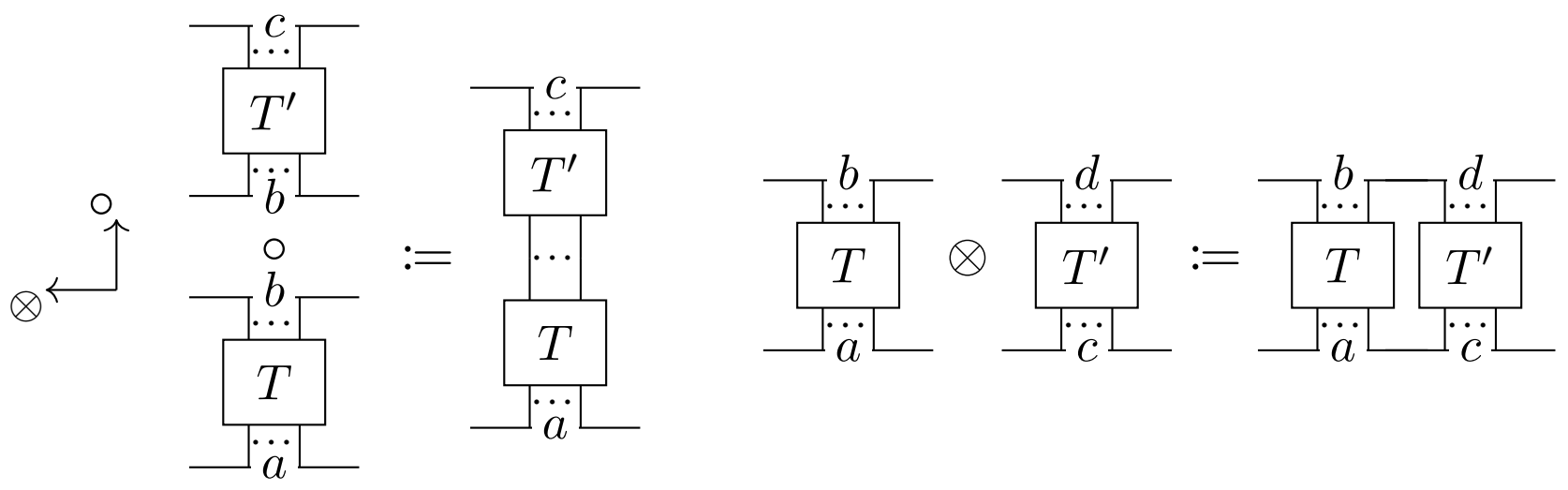

In other words, we view a tangle $T$ as a morphism $T\colon a \to b$ from its bottom endpoints $a$ to its top endpoints $b$. Two tangles $T\colon a \to b$ and $T'\colon b \to c$ with compatible source and target can be composed $T' \circ T\colon a \to c$ by stacking them atop each other, and vertically rescaling; and two tangles $T\colon a \to b$ and $T'\colon c \to d$ can be tensored $T \sqcup T'\colon a \sqcup c \to b \sqcup d$ using the disjoint union, putting them next to each other horizontally. An *invariant of tangles* is defined to be a monoidal functor

$$F\colon Tang \to C,$$

where $(C, \otimes, \mathbb{1})$ is another monoidal category. A similar construction also holds for *oriented* tangles; we restrict to the unoriented case for simplicity.

In this categorical paradigm for knot theory, the Reshetikhin–Turaev construction corresponds to the choice $C = A\text{-Mod}$, where $A$ is a bialgebra (say over the field $\mathbb{C}(q)$) and $A\text{-Mod}$ denotes its category of finite-dimensional representations. The monoidal structure on $A$-Mod is induced by the co-algebra structure on $A$. The ground field $\mathbb{C}(q)$ is the monoidal unit.

[2]More precisely, $Tang$ should be taken as a choice of strictification of what our intuitive definition suggests.

So what does it take to define a quantum invariant? Objects in $Tang$ are monoidally generated by the point; similarly, morphisms are monoidally generated by three elementary tangles, the *cup*, the *cap* and the *crossing* (pictured below). This means that a monoidal functor $F\colon Tang \to A\text{-Mod}$ is determined by the choice of a representation $V$ such that $F(\bullet) = V$, and intertwiners of representations $R$, $\epsilon$ and $\eta$, such that:

$F_{A,V}$ ( 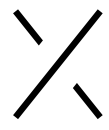 ) $= \begin{matrix} V \otimes V \\ \uparrow R \\ V \otimes V \end{matrix}$ ,

$F_{A,V}$ ( 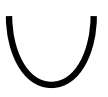 ) $= \begin{matrix} V \otimes V \\ \uparrow \epsilon \\ \mathbb{C}(q) \end{matrix}$ and $F_{A,V}$ ( 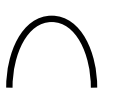 ) $= \begin{matrix} \mathbb{C}(q) \\ \uparrow \eta \\ V \otimes V \end{matrix}$ .

See the middle part of Fig. i.1 for an example. One can describe isotopies in $Tang$ with a finite set of relations, in analogy with the classical Reidemeister theorem for links. Transferring these relations through $F$ provides a finite set of conditions on the triple $(R, \epsilon, \eta)$; we shall say that this data makes $V$ a *tangle object*[3] in $A$-Mod. In that way, we reduce the problem of defining invariants of tangles with value in $A$-Mod to the problem of finding tangle objects in $A$-Mod.

Of course, the theory would all be in vain if tangle objects did not exist. This is where quantum groups come in, which we now properly introduce. Given a simple finite-dimensional complex Lie algebra $\mathfrak{g}$, one can construct its universal enveloping algebra $U(\mathfrak{g})$. The quantum group $U_q(\mathfrak{g})$ is a $q$-deformation of $U(\mathfrak{g})$; that is, an algebra defined over the field $\mathbb{C}(q)$, such that taking the limit $q \to 1$ (in some suitable sense) recovers the $\mathbb{C}$-algebra $U(\mathfrak{g})$. To motivate the terminology, note that the respective categories of representations of the associative algebra $U(\mathfrak{g})$, the Lie algebra $\mathfrak{g}$ and the associated simply-connected Lie group $G$ are all equivalent. Hence, from the point of view of representation theory, $U_q(\mathfrak{g})$ may also be thought as a deformation of the group $G$, although $U_q(\mathfrak{g})$ itself is not a group. Quantum groups were independently introduced by Drinfel'd [53] and Jumbo [93] in the 1980s, in their study of the Yang–Baxter equation in quantum integral systems. Quantum groups also play a prominent role in noncommutative geometry, as quantization deformation of algebras of functions on Lie groups.

[3] A more standard terminology would be *braided object with duals*.

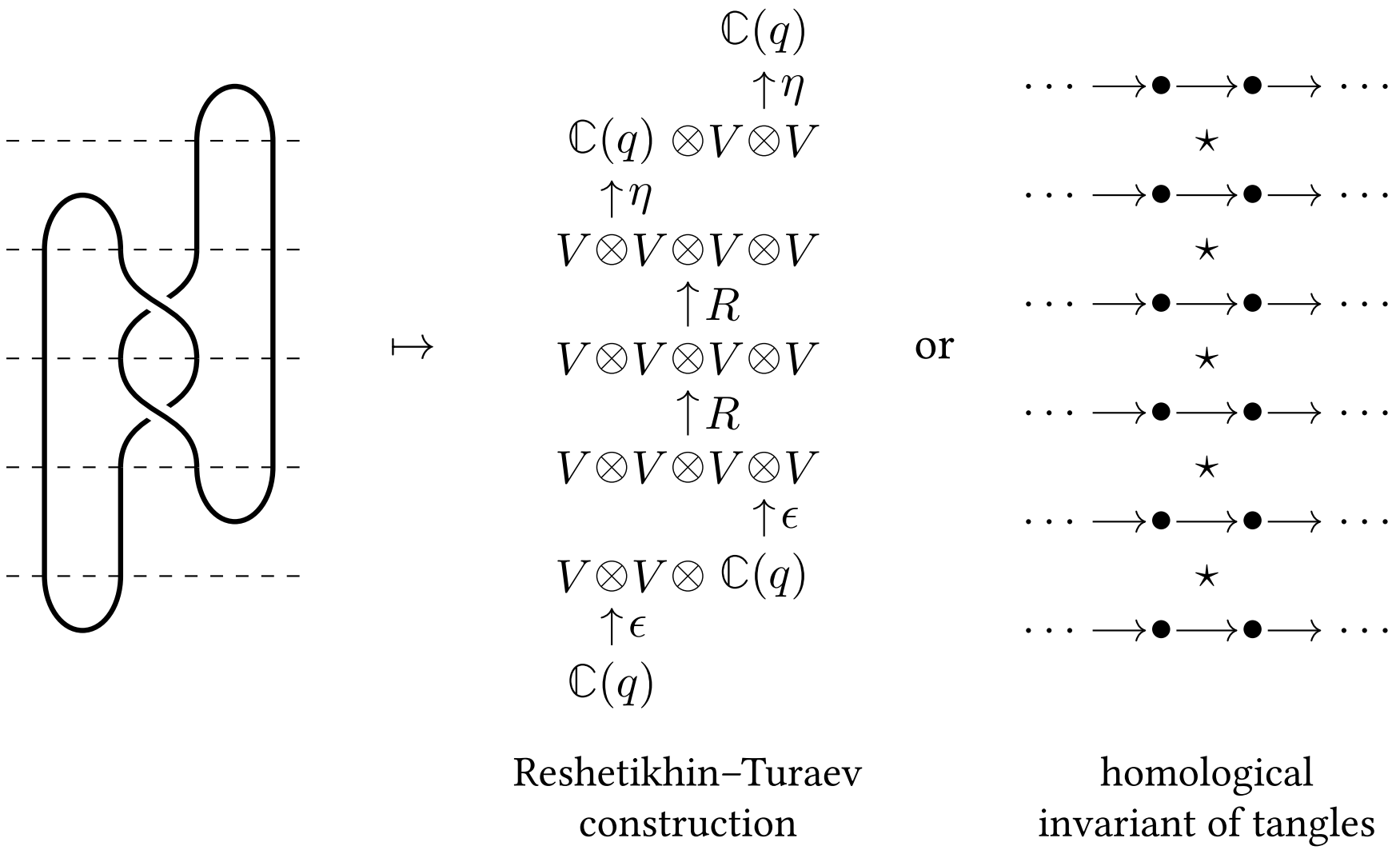

**Fig. i.1** The left-hand side pictures the Hopf link as a composition of elementary tangles; the middle part pictures the classical Reshetikhin–Turaev construction; and the right-hand its higher analogue, replacing composition of maps with composition (tensor product) of chain complexes.

So why quantum groups? Because they have the extremely rich structure of a *quasi-triangular Hopf algebra.* This turns out to be precisely the required structure to turn every finite-dimensional representation $V$ of $U_q(\mathfrak{g})$ into a tangle object. In other words, any choice of pair $(\mathfrak{g}, V)$ provides an invariant of tangles $F_{\mathfrak{g},V}\colon Tang \to U_q(\mathfrak{g})\text{-Mod}$. Note that if $L$ is a link, $F_{\mathfrak{g},V}(L)$ is an endomorphism of the ground ring $\mathbb{C}(q)$, and hence an element of $\mathbb{C}(q)$. In fact, one can show that $F_{\mathfrak{g},V}(L) \in \mathbb{Z}[q, q^{-1}]$, giving the aforementioned family of quantum polynomial invariants.

### i.1.2 Crane and Frenkel's conjecture

By the time quantum topology had emerged as a new field of study, tools from gauge theory (notably Donaldon's theory [50]) and sympletic geometry (notably Floer homology [78, 79, 80]) had been successfully applied to low-dimensional topology, particularly to the study of exotic structures in dimension four. Two smooth manifolds form an *exotic pair* if they are homeomorphic but not diffeomorphic. Dimension four has shown to be especially challenging, being the only dimension for which exotic $\mathbb{R}^4$s exist, and the only dimension for which the smooth Poincaré conjecture, stating that the four-sphere $S^4$ has a unique smooth structure, remains unsolved. While powerful, gauge theoretic and sympletic invariants are knowingly hard to compute, in sharp contrast with quantum polynomial invariants.

In 1989, Witten proposed a construction of the Jones polynomial using gauge theory [190]. This inspired Reshetikhin and Turaev to extend their method to 3-manifolds [160], leading to a 3-dimensional *Topological Quantum Field Theory* (TQFT); the reader should think of the latter as another terminology for a *local* invariant of 3-manifolds, similar to the role played by tangles for links. This is known as the *Witten–Reshetikhin–Turaev TQFT*.

In 1994, Crane and Frenkel conjectured that the whole construction could be generalized to four dimension [48], based on the concept of *categorification* [47]. More precisely, they conjectured that quantum groups were only the shadow of certain "Hopf categories" from which one could define a 4-dimensional TQFT, expected to be both computable and capable of detecting exotic structures.

To best explain these ideas, we now turn to the more algebraic side of our exposition. We will return to topology when discussing Khovanov homology.

### i.1.3 Higher algebra and string calculus

From a modern perspective, classical algebra is the study of objects with a categorical structure; for instance an algebra is a one-object linear category. *Higher algebra* extends this study to higher categorical structures, starting with the 2-dimensional case. Recall that a 2-category[4] has objects, morphisms, called *1-morphisms*, and morphisms between morphisms, called *2-morphisms*. The latter carry two distinct associative structures, the 0-composition $\star_0$ (composing along objects) and the 1-composition $\star_1$ (composing along 1-morphisms), compatible via the *interchange law* (see Fig. i.2). The 2-category of small categories, functors and natural transformations provides a stereotypical example, with the 0-composition and the 1-composition being respectively the horizontal and vertical compositions of natural transformations. Note that a one-object 2-category is precisely a strict monoidal category, with $\circ = \star_1$ and $\otimes = \star_0$.

[4] In this thesis, an $n$-category always refers to a strict $n$-category.

The 2-dimensional aspect of higher algebra, along with the additional complexity it brings, calls for a new form of calculus that would subsume the ubiquitous symbolic calculus in classical algebra. This calculus is known as *string calculus* (or *diagrammatic calculus*) to category theorists (following the work of Joyal and Street [98]) and as *tensor networks* to physicists (following the work of Penrose [151]). One constructs the *string diagram representation* of a 2-morphism by taking the graph Poincaré dual of its globular representation:

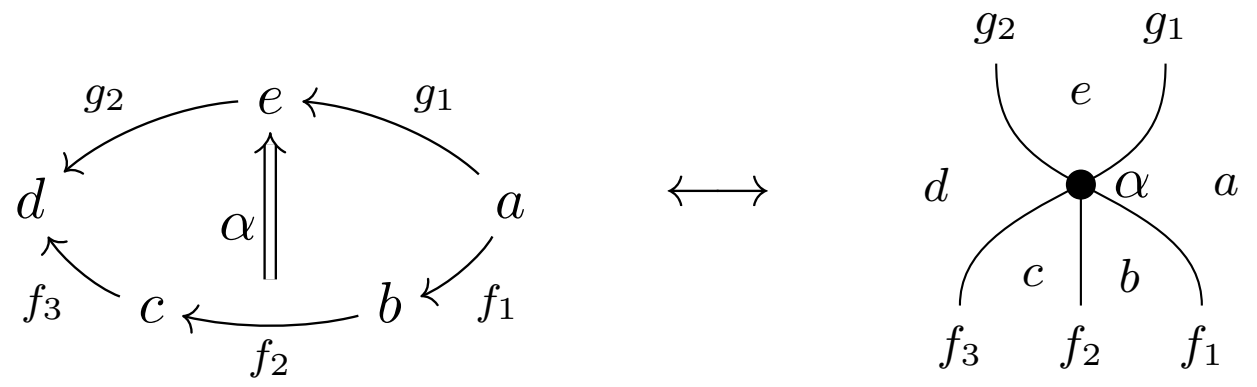

Note that we read string diagrams from the bottom-right to the top-left. With these diagrammatics, the interchange law amounts to "sliding two dots past one another" (see Fig. i.2). In fact, one can show that if there exists a planar isotopy preserving the labels between two string diagrams, then these string diagrams represent the same 2-morphism [98]. In practice, many categorical properties can be translated into meaningful topological interpretations, provided one uses appropriate conventions. For instance, an identity morphism is usually not pictured, be it a 1-morphism (a line) or a 2-morphism (a dot). As an application of this convention, consider an adjunction of categories

$$C \underset{\mathcal{G}}{\overset{\mathcal{F}}{\underset{\leftarrow}{\overset{\rightarrow}{\bot}}}} D,$$

with unit $\epsilon\colon \mathcal{F}\mathcal{G} \to 1_D$ and counit $\eta\colon 1_C \to \mathcal{G}\mathcal{F}$. The string diagrammatics of each triangular relation amounts to a simple "zigzag isotopy", as shown in Fig. i.3.

String diagrams are now perversive in many fields; we refer the reader to [9, 173], the nLab or Wikipedia for further comments on the history and applications.

## i.1.4 Categorification

*Categorification* is the process of enhancing an object with an additional, higher, structure. The inverse process, *de*categorification, forgets the higher structure and returns only a shadow of the original object. The former is

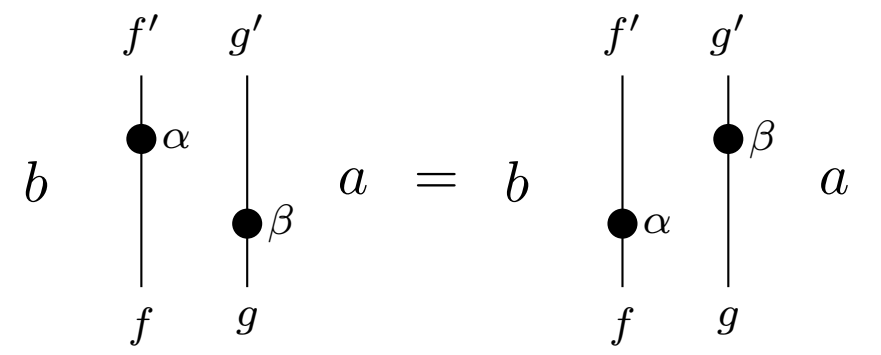

$$(\alpha \star_0 \mathrm{id}_{g'}) \star_1 (\mathrm{id}_f \star_0 \beta) = (\mathrm{id}_{f'} \star_0 \beta) \star_1 (\alpha \star_0 \mathrm{id}_g)$$

**Fig. i.2** interchange law in 2-categories

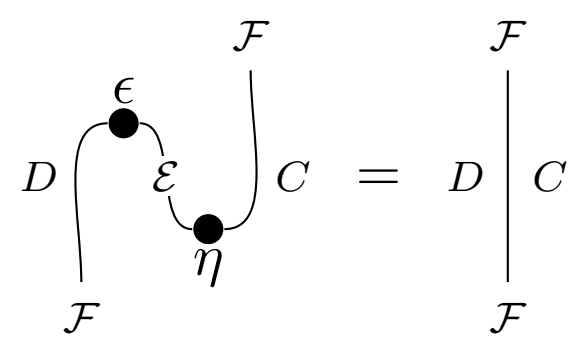

**Fig. i.3** string diagrammatics for adjunction of categories

a highly non-canonical process; the latter is context-dependent, but always well-defined. Intuitively, categorification "add morphisms between objects", lifting equalities to isomorphisms. The added data of morphisms provides finer ways to compare objects, away from the dichotomy of the equality. Let us motivate the idea through examples; we refer the reader to [118] for a great introduction.

- *Example* 1. If $X$ is a set, its decategorification is the cardinal $\mathrm{card}(X)$; conversely a cardinal $\aleph$ is categorified by a *choice* of set $X$ such that $\mathrm{card}(X) = \aleph$. An equality between cardinals is categorified by a *choice* of bijection between two sets. Note that in combinatorics, one often explains a numerical identity as stemming from a bijection between two finite sets.

- *Example* 2. If $V$ is a finite-dimensional vector space, its decategorification is the natural number $\dim(V)$. This process is compatible with the additive and monoidal structure on finite-dimensional vector spaces, in the sense that

$$\begin{aligned} &\dim(V \oplus W) = \dim V + \dim W \\ \text{and} \quad &\dim(V \otimes W) = \dim(V) \cdot \dim(W). \end{aligned}$$

  In that sense, linear algebra categorifies arithmetic.

- *Example* 2'. If now $V = \bigoplus_{i\in\mathbb{Z}} V_i$ is a *graded* finite-dimensional vector space, its decategorification is its *graded* dimension

$$\mathrm{gdim}(V) = \sum_{i\in\mathbb{Z}} \dim(V_i) q^i,$$

  and similar compatibilities hold.

- *Example* 3. If $V_\bullet$ is a bounded chain complex of finite-dimensional vector spaces, its decategorification is its Euler characteristic

$$\chi(V_\bullet) = \sum_{i\in\mathbb{Z}} \dim(V_i)(-1)^i.$$

In that sense, (co)homology theories categorify their Euler characteristics. Conversely, given an alternating sum, one can wonder whether it categorifies into an (interesting) homology.

- *Example* 4. Recall the monoidal category $\mathcal{T}ang$ from subsection i.1.1 and the notion of 2-category from subsection i.1.3. One can analogously define a monoidal 2-category 2-$\mathcal{T}ang$ (the reader should not focus on the precise definition of a monoidal 2-category) of embedded surfaces in $\mathbb{R}^4$ up to isotopies, viewed as 2-morphisms, or *cobordisms*, between tangles (see the schematic on the right)[a]. In this construction, tangles are no longer considered up to isotopies; instead, two isotopic tangles are related by an embedded surface, invertible in 2-$\mathcal{T}ang$. The monoidal 2-category 2-$\mathcal{T}ang$ is a categorification of the monoidal category $\mathcal{T}ang$.

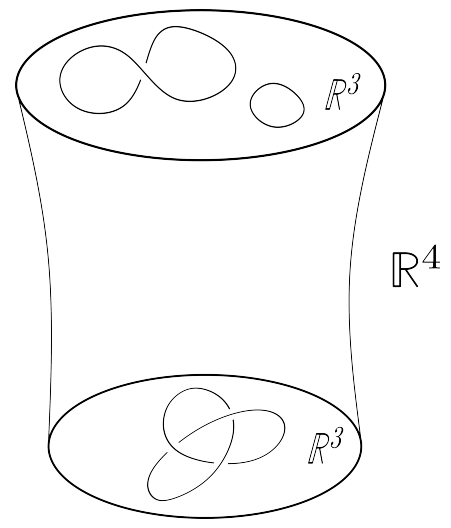

**Fig. i.4** Schematic for link cobordism

[a] As for $\mathcal{T}ang$, 2-$\mathcal{T}ang$ should be taken as a choice of strictification of what our intuitive definition suggests.

Each of these examples fits into a well-defined notion of decategorification, for a given choice of categorical structure. For instance, if $\mathcal{C}$ is an additive monoidal category, we define its (split) *Grothendieck ring* $K_0(\mathcal{C})$ as generated as an abelian group by elements $[A]$ for each object $A \in \mathrm{ob}(\mathcal{C})$, quotiented by the relation $[C] = [A] \oplus [B]$ whenever $C \cong A \oplus B$, and with the induced multiplicative structure $[A]\cdot[B] := [A \otimes B]$. In particular, if $A \cong B$ then $[A] = [B]$. By definition, the additive (resp. the multiplicative) structure on $K_0(\mathcal{C})$ is induced by the additive (resp. the monoidal) structure on $\mathcal{C}$. Schematically:

$$\mathrm{K}_0\colon (\mathcal{C},\oplus,\otimes) \rightsquigarrow (K_0(\mathcal{C}),+,\cdot).$$

When $\mathcal{C}$ is the category of finite-dimensional vectors spaces, we recover Example 2, with $(K_0(\mathcal{C}),+,\cdot) = (\mathbb{Z},+,\cdot)$ and $[V] = \dim(V)$.

On the other hand, categorification is always a *choice*; it needs not be unique, and may not even exist. Categorification adds morphisms between objects, and it is in these new morphisms that may lie exiting new mathematics; either as a mean to better understand its shadow, and for its own sake.

### i.1.5 Higher representation theory

What does it take then to categorify a $\mathbb{C}$-algebra $A$? Recall the Grothendieck ring $K_0(\mathcal{A})$ of an additive monoidal category $\mathcal{A}$ is a ring; that is, a $\mathbb{Z}$-algebra. We say that $\mathcal{A}$ *categorifies* $A$ if

$$K_0(\mathcal{A}) \otimes_{\mathbb{Z}} \mathbb{C} \cong A.$$

If $\mathcal{A}$ is sufficiently nice[5], every object can be uniquely written as a finite sum of indecomposable objects. Let $\{B_i\}_{i\in I}$ be a set of indecomposable objects, unique and complete up to isomorphism. Then, the tensor product $B_j \otimes B_k$ can be written as a direct sum of $B_i$'s, and:

$$B_j \otimes B_k \cong \bigoplus_{i\in I} B_i^{\oplus n_{jk}^i} \rightsquigarrow [B_j]\cdot[B_k] = \sum_{i\in I} n_{jk}^i [B_i]$$

where the $n_{jk}^i$'s are *positive* integers. In other words, the categorification $\mathcal{A}$ implies the existence of an integral form $K_0(\mathcal{A})$ for $A$, together with a basis $[B_i]$ whose multiplication rule has positive coefficients!

A similar remark applies for algebras over the field $\mathbb{C}(q)$. Recall that an algebra is the same as a one-object linear category; it is categorified by a monoidal category, that is, a one-object additive 2-category. Similarly, an *idempotented algebra* is the same as a linear category, categorified by an additive 2-category. It turns out that each quantum group $U_q(\mathfrak{g})$ admits an idempotented integral form $\dot{U}_q(\mathfrak{g})_{\mathbb{Z}}$ admitting a basis with positive multiplicative coefficients, known as the Luzstig–Kashiwara canonical basis. This suggests that a categorification could exist. Indeed, *categorified quantum groups* were ultimately constructed by Khovanov–Lauda [106, 116] and Rouquier [166], building on important previous work [37, 82]; see [22] for the equivalence between the two constructions. Categorified quantum groups *à la* Khovanov–Lauda are defined by generators and relations, heavily relying on string calculus (subsection i.1.3).

A representation of an algebra $A$ is an algebra morphism $F : A \to \mathrm{End}(V)$, where $V$ is a vector space. Viewing both $A$ and $\mathrm{End}(V)$ as linear categories makes $F$ into a linear functor. More generally, a representation of $A$ is the data of a linear category $C$ and a linear functor $F\colon A \to C$. Similarly, a

[5] $\mathcal{A}$ is a Krull–Schmidt category.

*2-representation of an additive 2-category* $\mathcal{A}$ is the data of an additive 2-category $\mathcal{C}$ and an additive 2-functor $\mathcal{F}\colon \mathcal{A} \to \mathcal{C}$. If $\mathcal{A}$ categorifies $A$, taking the Grothendieck ring on both sides induces a representation of (an integral form of) $A$:

$$\begin{array}{ccc} \mathcal{A} & \xrightarrow{\mathcal{F}} & \mathcal{C} \\ {\scriptstyle K_0}\wr\downarrow & \circlearrowleft & \wr\downarrow{\scriptstyle K_0} \\ A & \xrightarrow[K_0(F)]{} & C \end{array}$$

In the opposite direction, given a representation of $A$, one can wonder wether it lifts to a 2-representation of a certain categorification of $A$.

The study of 2-representations is the field of *higher representation theory*. It was initiated in seminal papers such as [134], [37] and [77, 147], and further developed in [133, 134, 135, 136, 137, 138]. As its classical counterpart, it serves as a unifying principle to understand higher symmetries, leading for instance to a proof of Broué's conjecture for the symmetric group by Chuang and Rouquier [37].

Let us finally mention the parallel story on the Coxeter side, historically preceding the categorification of quantum groups. Recall the notion of Coxeter groups, an abstraction of the notion of reflection groups. They are closely related to Lie algebras, as to every simple finite-dimensional complex Lie algebra $\mathfrak{g}$ corresponds a canonical reflection group known as its Weyl group. Similarly to quantum groups, the group algebra of a Coxeter group admits a $q$-deformation, known as its *Hecke algebra*, equipped with an integral basis having positive multiplicative coefficients. A candidate categorification arose from Soergel's approach to the Kazhdan–Lusztig conjectures [103, 176], whose statements have no reference to categorification. The proof was originally given by Soergel in the case of Weyl groups, and completed by Elias and Williamson [67, 68, 70, 71] for all Coxeter groups, using diagrammatic methods.

### i.1.6 Khovanov homology

In 1999, Mikhail Khovanov announced the discovery of a new homological invariant of links $\mathrm{Kh}(L)$, now called *Khovanov homology*, whose Euler characteristic is the Jones polynomial $J(L)$. (In addition to its homological grading, Khovanov homology has an extra grading, which explains why its Euler characteristic is a polynomial rather than a number; see Example 2' and Example 3 above.) This initiated a renewal of the field of quantum topology, focusing on link homologies rather than link polynomials.

As most (co)homology theories, Khovanov homology's significance lies in its *functorial* nature. Recall the proper notions of morphisms of links, known as link cobordisms, that we sketched in Example 4. Khovanov homology is functorial in the sense that to any link cobordism $S\colon L_1 \to L_2$ corresponds a linear map

$$\mathrm{Kh}(S)\colon \mathrm{Kh}(L_1) \to \mathrm{Kh}(L_2)$$

between the respective Khovanov homologies. This essential property has been used to either recover results only known via gauge theory, or achieve new results in knot theory [111, 153, 158]. Recent results have extended Khovanov homology to an invariant of smooth four-manifolds [143], with promising properties regarding detection of exotic pairs [159].

Other Khovanov-like link homologies [105, 126, 187], constituting the family of *quantum link homologies*, have been constructed. Many admit spectral sequences to gauge or sympletic counterparts; see [108, section 3] for a survey. Typically, quantum invariants are easier to compute, but with obscure topological interpretations; and vice-versa for gauge or symplectic analogues. Understanding their relationship remains an active field of research.

Khovanov homology serves as a nexus between a wide range of mathematical ideas, as testified by the numerous approaches that were developed: singular surfaces [14, 19, 29, 39, 105, 126, 167], categorical skew Howe duality [115, 156], matrix factorizations [110], categorification of tensor products of representations [187], category $\mathcal{O}$ [132], algebraic geometry [32], symplectic geometry [172], gauge theory [189], and mirror symmetry [2]. Although each has its own flavor, most have been shown to be equivalent [31, 128]. See [108] for a survey.

### i.1.7 Foams and functoriality

We now describe three approaches to Khovanov homology: Khovanov's original approach, Bar-Natan's approach leading to a *tangle* version of Khovanov homology, and Blanchet's approach leading to a *functorial* version of Khovanov homology.

*Khovanov's original approach*

The representation theory involved in the definition of the Jones polynomial admits a diagrammatic description, given by the Temperley–Lieb category **TL**. Kauffman used this fact to give a combinatorial description of the Jones polynomial. The main tool is so-called the Kauffman bracket, which "resolves" a crossing into a linear combination of planar curves; see Fig. i.5. If $L$ is a link and $D_L$ is a choice of link diagram for $L$ (a planar depiction of $L$ as in Fig. i.1), one can "resolve" the diagram into a linear combination of diagrams consisting only in circles. Evaluating a union of $n$ circles to $(q+q^{-1})^n$ defines a polynomial $\langle D_L \rangle \in \mathbb{Z}[q,q^{-1}]$, known as the *Kauffman bracket of* $L$. It coincides, up to normalization[6], with the Jones polynomial, and thus provides a simple combinatorial definition leaving the representation theory in the background. This can be phrased as a monoidal functor

$$\langle - \rangle \colon \mathit{Tang} \to \mathbf{TL}$$

as in subsection i.1.1.

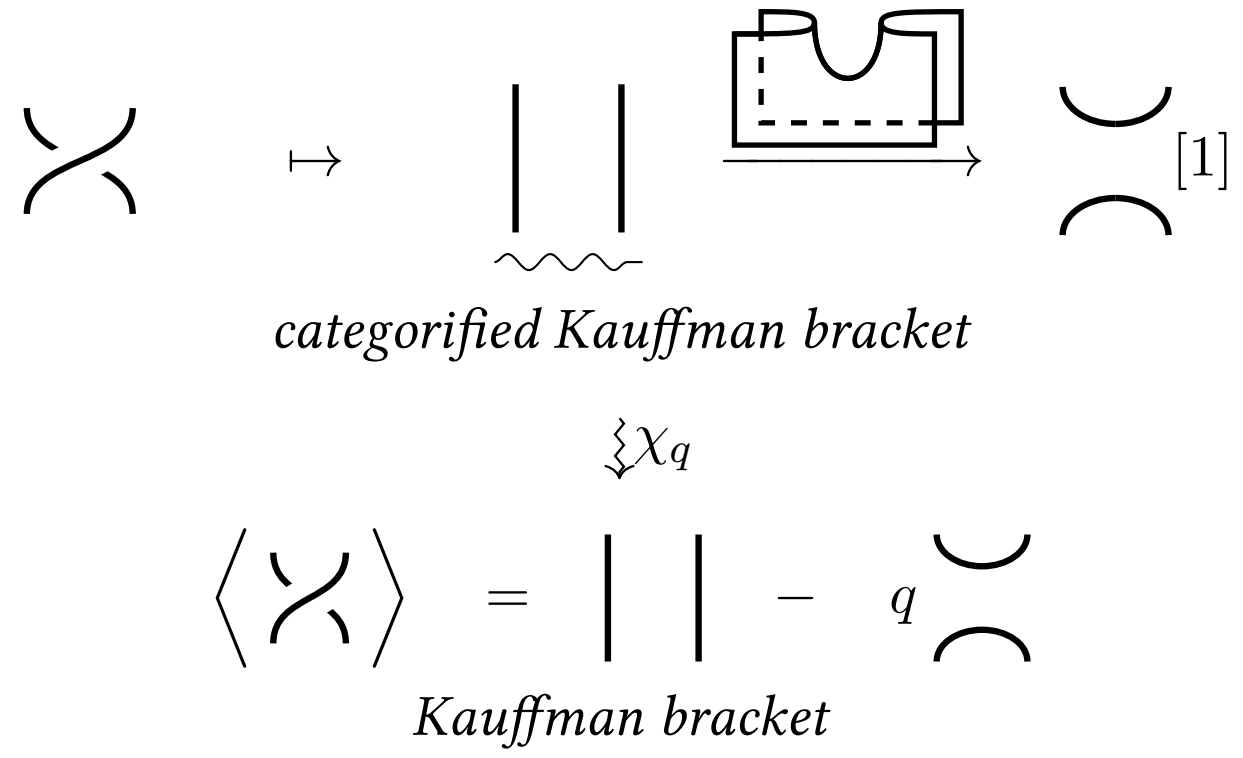

**Fig. i.5** The definition of Khovanov homology in a nutshell. The homological degree zero is underlined, and the symbol [1] denotes a shift in the quantum grading.

Khovanov's brilliant insight was to lift the Kauffman bracket to a length-two complex, replacing the minus sign with a saddle cobordism; see Fig. i.5. Doing so for every crossing in the diagram $D_L$ defines a hypercube of dimension $n$, whose vertices are union of circles and edges are saddles, either merging two circles into one ("merge saddle") or splitting a circle into two ("split saddle"). See Fig. i.6 for an example.

To turn these pictures into something computable, one replaces every circle by a (graded) vector space $V$, and each merge saddle (resp. split saddle) with a linear map $m\colon V\otimes V\to V$ (resp. $\Delta\colon V\to V\otimes V$); formally, this amounts to applying a 2-dimensional TQFT (see subsection i.1.2) to the hypercube.

[6]As before, we ignore orientation for simplicity.

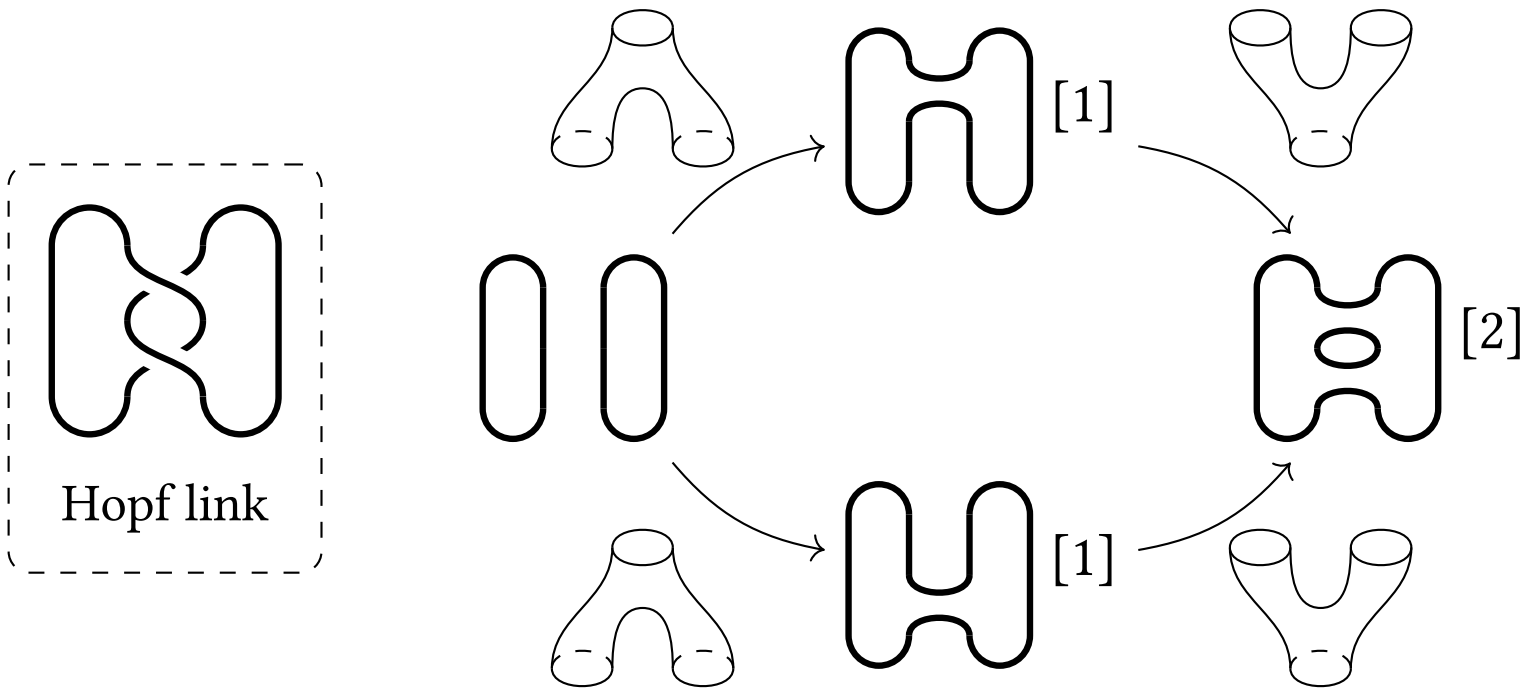

**Fig. i.6** Khovanov's hypercube (here a square) for the Hopf link.

Adding suitable signs following the Koszul rule turns the hypercube into a chain complex. Khovanov showed that, for a wisely chosen 2-dimensional TQFT, the homotopy type of the resulting complex is an invariant of links which categorifies the Jones polynomial.

*A local combinatorial approach to Khovanov homology*

Khovanov's original definition is not local: whether a saddle is a merge saddle or a split saddle is a global property. One solution, pursued by Bar-Natan [14], is to lift the relations imposed by the 2-dimensional TQFT directly on the cobordisms. More precisely, he defined a linear 2-category $\mathbf{BN}$ whose 2-morphisms are linear combinations of cobordisms up to some relations. One then view the hypercube directly as a hypercube in the linear 2-category $\mathbf{BN}$, which can now be understood as arising from a tensor product of length-two complexes, each associated to a crossing; see the right-hand side of Fig. i.1. Bar-Natan showed that the homotopy type of this complex is an invariant of tangles, in the sense given in subsection i.1.1. This leads to a monoidal functor $\mathrm{Kh}\colon \mathit{Tang} \to \mathcal{K}(\mathbf{BN})$ lifting the Kauffman monoidal functor $\langle - \rangle \colon \mathit{Tang} \to \mathbf{TL}$, where $\mathcal{K}(\mathbf{BN})$ denotes the homotopy category of $\mathbf{BN}$:

$$\begin{array}{ccc} & & \mathcal{K}(\mathbf{BN}) \\ & \overset{\mathrm{Kh}}{\nearrow} & \big\downarrow K_0 \\ \mathit{Tang} & \xrightarrow{\langle - \rangle} & \mathbf{TL} \end{array}$$

*local* Khovanov homology

$$\begin{array}{ccc} 2\text{-}\mathcal{T}\mathit{ang} & \xrightarrow{?} & \mathrm{Ch}_\bullet(\mathbf{BN}) \\ \big\downarrow K_0 & & \big\downarrow K_0 \\ \mathit{Tang} & \xrightarrow{\langle - \rangle} & \mathbf{TL} \end{array}$$

*functorial* local Khovanov homology?

In the above construction, if two tangles are isotopic, then the associated chain complexes are homotopic. One can wonder whether for each isotopy, there exists a *choice* of homotopy between the associated chain complexes, natural in some ways. Recall that in the monoidal 2-category 2-$\mathcal{T}ang$ introduced in Example 4, isotopies of tangles lift to 2-isomorphisms. Better still, one can wonder whether data can be assigned to every 2-morphism of 2-$\mathcal{T}ang$, that is, to every tangle cobordism. This would define a monoidal 2-functor 2-$\mathcal{T}ang \to \mathrm{Ch}_\bullet(\mathbf{BN})$ whose decategorification is the Kauffman bracket, where $\mathrm{Ch}_\bullet(\mathbf{BN})$ is monoidal 2-category of chain complexes in $\mathbf{BN}$. In general, a *functorial invariant of tangles* is the data a monoidal 2-functor $\mathcal{F}\colon$ 2-$\mathcal{T}ang \to \mathcal{C}$ for some monoidal 2-category $\mathcal{C}$. Bar-Natan showed his construction to be functorial *up to signs*.

*A functorial local approach to Khovanov homology*

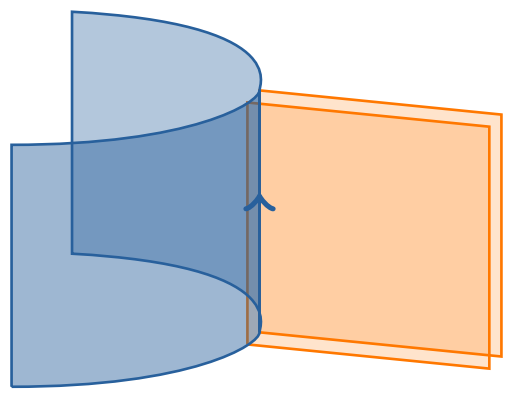

Solutions to this sign problem were provided by numerous authors [29, 39, 167, 185]. A solution introduced by Blanchet [19] used *foams*, certain decorated singular surfaces first introduced by Khovanov in his definition of $\mathfrak{sl}_3$-Khovanov homology [105] and generalized to $\mathfrak{gl}_n$ for $n \geq 3$ in [126]. Adapting this construction to the $n = 2$ case, Blanchet defined a functorial version of Khovanov homology using $\mathfrak{gl}_2$-*foams*[7] (see on the right). The proof of functoriality was later generalized to all $\mathfrak{gl}_n$-link homologies in [65]. In practice, working with $\mathfrak{gl}_2$ instead of $\mathfrak{sl}_2$ leads to better-behaved constructions. This comes down to the fact that in the former case, the trivial representation and the determinant representation are not isomorphic, and keeping track of this distinction leads to better control on signs (see [115, 1F] for a discussion).

At the decategorified level, moving from $\mathfrak{sl}_2$ to $\mathfrak{gl}_2$ does not affect the topology: the associated invariant remains the Jones polynomial. The Temperley–Lieb category is replaced by the category $\mathbf{Web}_d$[8] of $\mathfrak{gl}_2$-*webs*. The latter are trivalent graphs providing a combinatorial description of the representation theory of $U_q(\mathfrak{gl}_2)$:

[7]Blanchet called them *enhanced* $\mathfrak{sl}_2$-*foams*, but as they were later understood to be related to $\mathfrak{gl}_2$ rather than $\mathfrak{sl}_2$, we call them $\mathfrak{gl}_2$-foams.

[8]The parameter $d \in \mathbb{N}$ is a certain width parameter. Fixing $d$ has the effect of forgetting the monoidal structure, which we implicitly do below.

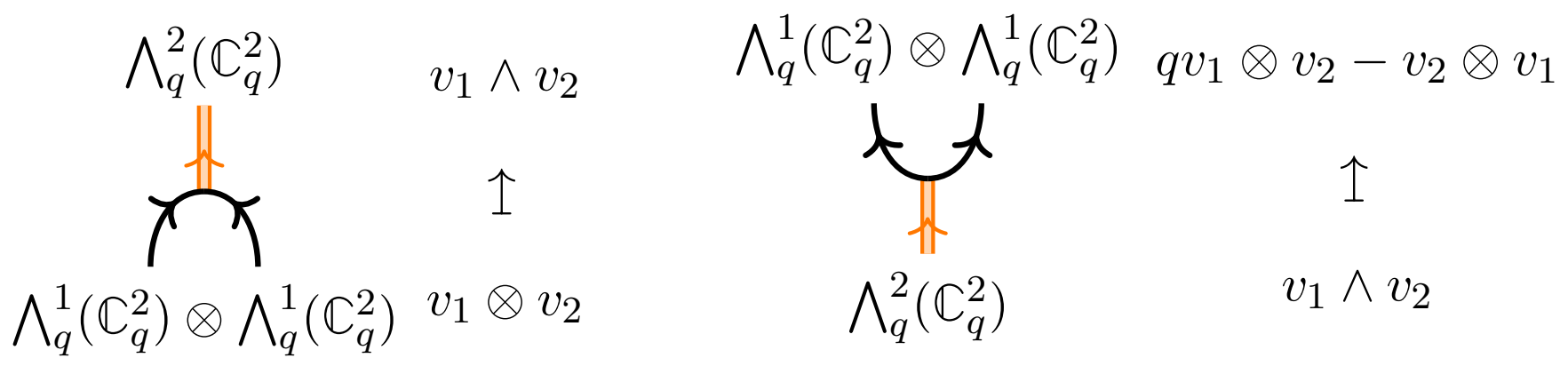

**Fig. i.7** Generating $\mathfrak{gl}_2$-webs.

Here $\mathbb{C}_q := \mathbb{C}(q)$ and $\bigwedge_q^0(\mathbb{C}_q^2)$ is the trivial representation, $\bigwedge_q^1(\mathbb{C}_q^2)$ is the standard representation, and $\bigwedge_q^2(\mathbb{C}_q^2)$ is the determinant representation. At the categorified level, the Bar-Natan 2-category $\mathbf{BN}$ is replaced by the 2-category $\mathbf{Foam}_d$ of $\mathfrak{gl}_2$-foams, which categorifies $\mathbf{Web}_d$.

### i.1.8 Categorical skew Howe duality

Let us return to Crane and Frenkel's conjecture (subsection i.1.2). Recall the notion of a functorial invariant of tangles $\mathcal{F}\colon 2\text{-}\mathcal{T}ang \to \mathcal{C}$ from the previous section. The monoidal 2-category $2\text{-}\mathcal{T}ang$ admits a presentation by generators and relations [8, 30], so that $\mathcal{F}$ is determined by a *2-tangle object* in $\mathcal{C}$, similarly to the classical case (see subsection i.1.1). Then one may hope that categorified quantum groups (subsection i.1.5) carry a certain rich categorical structure ("Hopf category") that naturality induces a structure of 2-tangle object on any of its 2-representations. Making this precise is a possible interpretation of Crane and Frenkel's conjecture, in the case of links (see also [179]).

A major difficulty in this statement comes from defining the tensor product of 2-representations. Indeed, recall that in the classical case, the tensor product of representations is induced by the coproduct on $U_q(\mathfrak{g})$. However, *by design* categorification only sees the product, and it is not clear what the higher analogue of the coproduct would be; see however [129] for recent progress. To circumvent this problem, Webster [187] has given an approach to link homologies by directly categorifying tensor product of representations. His construction provides a link homology for any pair $(\mathfrak{g}, V)$ categorifying the associated quantum polynomial, vastly generalizing Khovanov homology. Figure i.8 gives an overview of what we discussed so far.

While conceptually appealing, Webster's approach is technically demanding. For instance, functoriality in this setting is yet to be achieved. It turns out that the foam construction can also be understood via higher representation theory, albeit with a different taste. This goes via *categorical skew Howe duality*. We start by recalling classical skew Howe duality.

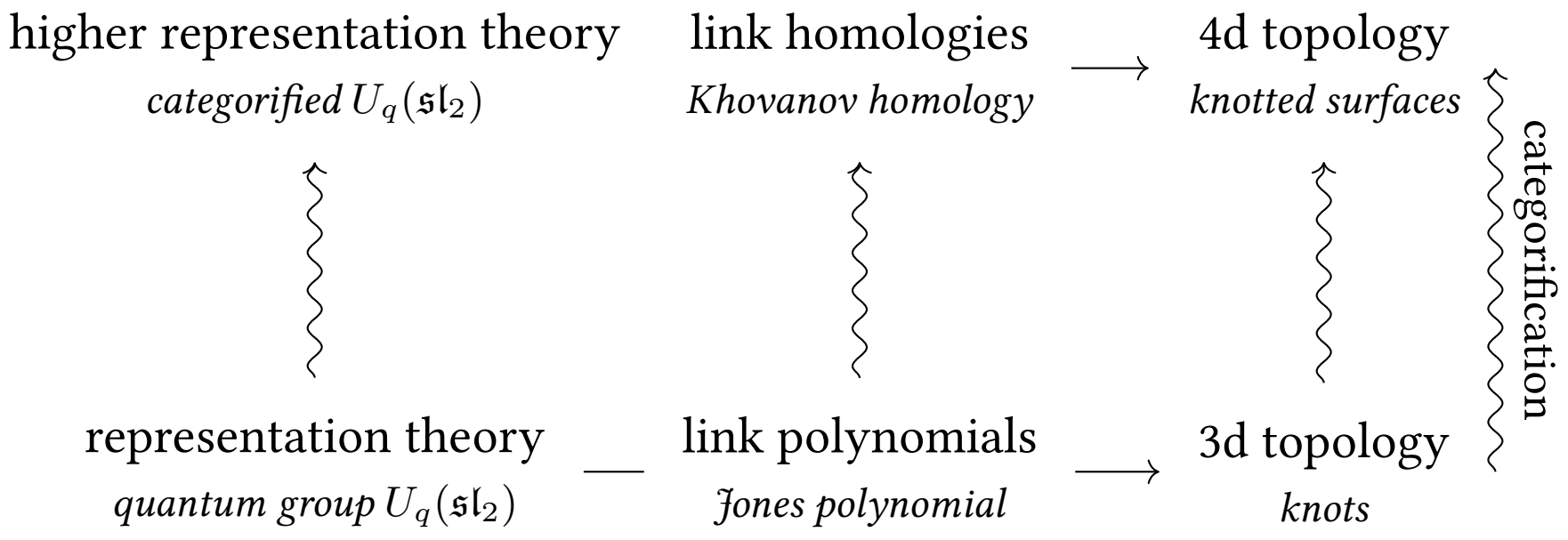

**Fig. i.8** A biased view on the field of quantum topology.

Recall the notation $\mathbb{C}_q := \mathbb{C}(q)$, and consider the (quantum) exterior power $\bigwedge_q^d(\mathbb{C}_q^n \otimes \mathbb{C}_q^2)$ equipped with the actions of $U_q(\mathfrak{gl}_n)$ and $U_q(\mathfrak{gl}_2)$. These actions commute, and in fact describe each other intertwiners. This fact is known as *skew Howe duality*, and was used (generalizing from $\mathfrak{gl}_2$ to $\mathfrak{gl}_k$) by Cautis, Kamnitzer and Morisson to describe $U_q(\mathfrak{gl}_k)$-intertwiners on fundamental representations by generators and relations [33]. In particular, the map

$$\phi\colon U_q(\mathfrak{gl}_n) \to \mathrm{End}_{U_q(\mathfrak{gl}_2)}\left(\textstyle\bigwedge_q^d(\mathbb{C}_q^n \otimes \mathbb{C}_q^2)\right)$$

is surjective. Taking the quotient by the kernel of this action defines the *$q$-Schur algebra of level two*:

$$S_{n,d} := U_q(\mathfrak{gl}_n)/\ker\phi \cong \mathrm{End}_{U_q(\mathfrak{gl}_2)}\left(\textstyle\bigwedge_q^d(\mathbb{C}_q^n \otimes \mathbb{C}_q^2)\right).$$

While the above isomorphism holds over $\mathbb{C}(q)$, the two sides have distinct interpretations, and hence distinct natural choices of integral and idempotented forms. On the one hand, $U_q(\mathfrak{gl}_2)$-intertwiners of fundamental representations are generated by the intertwiners presented in Fig. i.7, leading to the integral and idempotented form $\mathbf{Web}_d$ (recall that a small linear category is the same data as an idempotented algebra). On the other hand, Recall from subsection i.1.5 that $U_q(\mathfrak{gl}_n)$ admits an integral and an idempotented form denoted $\dot{U}_q(\mathfrak{gl}_n)$, which also provides an integral and idempotented form $\dot{S}_{n,d}$ for the $q$-Schur algebra of level two via the above quotient. Switching back to the categorical terminology, these $\mathbb{Z}[q,q^{-1}]$-linear categories are related by

$\mathbb{Z}[q, q^{-1}]$-linear functors:

$$\dot{U}_q(\mathfrak{gl}_n) \to \dot{S}_{n,d} \to \mathbf{Web}_d.$$

The relationship between $\dot{S}_{n,d}$ and $\mathbf{Web}_d$ is diagrammatically described by so-called *ladder webs*; see section 3.1 for details.

As shown Lauda, Queffelec and Rose [115, 156], skew Howe duality can be categorified:

$$\begin{array}{ccccc} \mathcal{U}(\mathfrak{gl}_n) & \longrightarrow & \mathcal{S}_{n,d} & \longrightarrow & \mathbf{Foam}_d \\ \downarrow K_0 & & \downarrow K_0 & & \downarrow K_0 \\ \dot{U}_q(\mathfrak{gl}_n) & \longrightarrow & \dot{S}_{n,d} & \longrightarrow & \mathbf{Web}_d \end{array}$$

The top arrows are 2-functors, $\mathcal{U}(\mathfrak{gl}_n)$ is the categorified quantum group discussed in subsection i.1.5, and $\mathcal{S}_{n,d}$ is the *2-Schur algebra* as introduced by Mackaay–Stošić–Vaz [127]. The latter is obtained by taking an appropriate quotient of $\mathcal{U}(\mathfrak{gl}_n)$. It categorifies the Schur algebra $\dot{S}_{n,d}$. The composition of the two top arrows is called the *foamation 2-functor* in [115], exhibiting $\mathbf{Foam}_d$ as a 2-representation of $\mathcal{U}(\mathfrak{gl}_n)$. In that sense, it provides a higher representation theoretic understanding of Khovanov homology.

## i.1.9 An odd analogue to Khovanov homology

Given the central role played by Khovanov homology, it is surprising to discover that the Jones polynomial admits another[9] categorification, called *odd Khovanov homology* [149]. It agrees with Khovanov homology modulo 2, but the two homologies are distinct in the sense that one can find pairs of knots distinguished by one but not the other, and vice-versa [175]. Heuristically, while (even) Khovanov homology is related to commutative algebra, odd Khovanov homology should be thought as an anti-commutative, hence *odd*, analogue:

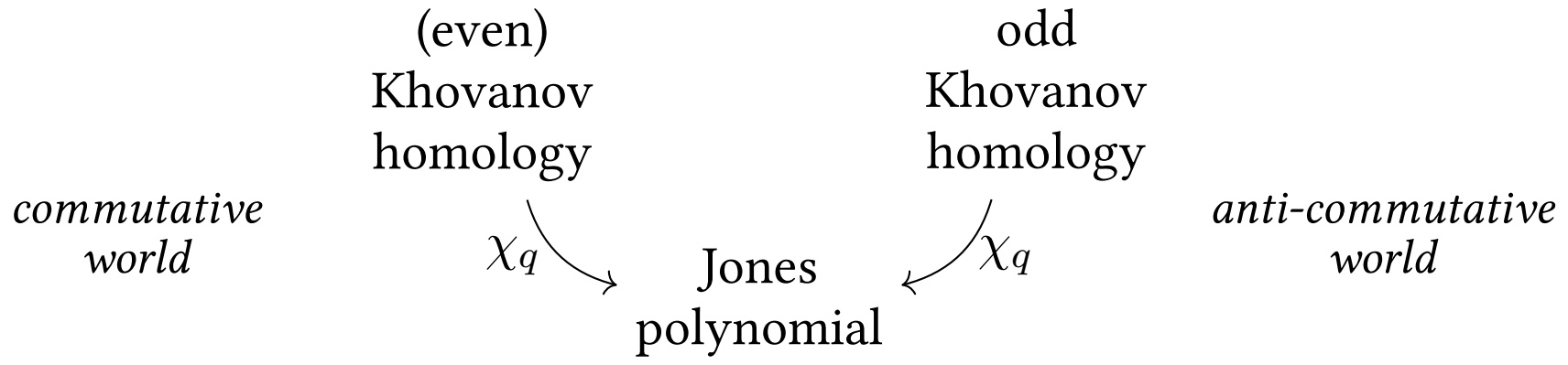

[9] In another direction, symmetric Khovanov homology [157, 165] is yet another categorification of the Jones polynomial.

Odd Khovanov homology does not immediately fit into Webster's general theory; it is the only known case of an even and odd pair of link homologies. Let us emphasize that this is a strictly higher phenomena, as both even and odd Khovanov homology categorify the Jones polynomial.

The definition of odd Khovanov homology is based on Khovanov's original definition sketched in subsection i.1.6. Recall that in a nutshell, it consists in introducing a *hypercube of resolutions* associated with every link diagram, using a suitable 2-dimensional TQFT to algebrize the hypercube, and turning the hypercube into a chain complex by assigning signs to its edges following the Koszul rule. In that case, the 2-dimensional TQFT associate to a union of $n$ circles the polynomial algebra $\mathbb{Z}[x_1,\ldots,x_n]/(x_1^2,\ldots,x_n^2)$.

Odd Khovanov homology is constructed similarly, substituting the polynomial algebra $\mathbb{Z}[x_1,\ldots,x_n]/(x_1^2,\ldots,x_n^2)$ for the exterior algebra $\bigwedge(x_1,\ldots,x_n)$. Since the exterior algebra is not commutative, the associated TQFT is only a *projective* TQFT in the sense that it is only functorial up to signs. Somehow miraculously, it was shown in [149] that this defect in functoriality can be balanced out when assigning signs to the hypercube. This however requires a much more intricate sign assignment than the Koszul rule, based on a case-by-case analysis of possible squares in the hypercube.

## i.1.10 The oddification program

The appearance of odd Khovanov homology sparked interest in finding odd analogues to known categorified and geometric structures [24, 25, 26, 59, 61, 62, 74, 75, 76, 109, 117, 145], motivated by their relation with (even) Khovanov homology. These investigations, known as the *oddification program* [74], were paralleled by independent developments on categorifications of super Kac–Moody algebras [100, 101, 102]. Similar behaviours have appeared in related fields, including condensed matter physics [1, 181] and gauge theory [140].

Underlying many of these constructions is the notion of a *super-2-category* [24] (or 2-supercategory), a certain categorical structure akin to a linear 2-category, where the interchange law is twisted by the extra data of a $\mathbb{Z}/2\mathbb{Z}$-grading, or *parity*, on the 2Hom-spaces. Diagrammatically, this is pictured as follows:

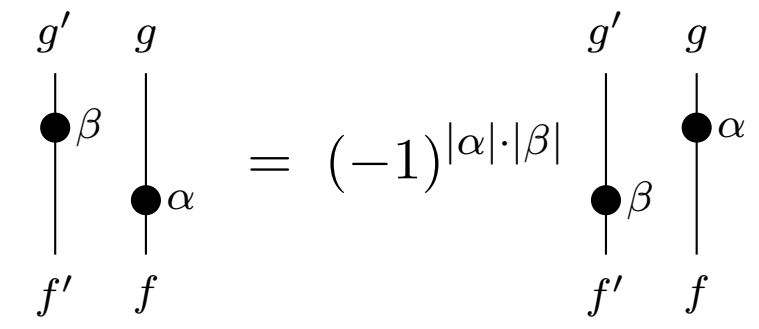

Here $|\alpha|$ and $|\beta|$ denote the respective parities of the 2-morphisms $\alpha$ and $\beta$. Note that in particular, a super-2-category is *not* a 2-category endowed with extra structure. We call *supercategorification* the process of categorifying a category with a super-2-category. Heuristically, an *odd analogue* (or *super analogue*) is an algebraic structure demonstrating anti-commutative behaviours, with the same graded rank as its even (commutative) counterpart and such that the even and odd constructions agree when reduced modulo two.

We conclude this overview of the field with a (biased) portrait of odd analogues, established or conjectural. See section i.2 for further discussions.

(i) *arc algebras.* Putyra gave a topological interpretation of odd Khovanov homology using so-called *chronological cobordisms* [155]. An odd analogue of the arc algebra was defined in [145], later used in [144] to give an extension of odd Khovanov homology to tangles. Interestingly, the non-commutativity translates to non-associativity for the odd arc algebra.

(ii) *combinatorics.* Khovanov homology is closely related to the combinatorics of symmetric functions. Ellis and Lauda introduced *odd symmetric functions* [72, 73], on which an odd analogue of the nilHecke algebra acts [74].

(iii) *Lie super theory.* The quantum group $U_q(\mathfrak{g})$ admits an odd analogue, associated to a certain Lie superalgebra counterpart of $\mathfrak{g}$. Both fit into a single covering quantum group $U_{q,\pi}(\mathfrak{g})$, where $\pi$ is an extra parameter such that setting $\pi = 1$ recovers the even case, while setting $\pi = -1$ recovers the super case [38, 41, 42, 43, 44, 45].

(iv) *categorification.* Recall from subsection i.1.5 that quantum groups can be categorified. Analogously, Brundan and Ellis [25] gave a construction of supercategorified quantum groups (or super 2-Kac–Moody algebras), building on previous work [75, 100, 101, 102]. Conjecturally, their construction provide a categorification of the quantum covering group $U_{q,\pi}(\mathfrak{g})$. Outside $\mathfrak{sl}_2$ [26, 57], this conjecture remains open.

(v) *Floer and gauge theories.* Odd Khovanov homology was discovered by Ozsváth, Rasmussen and Szabó in an attempt to lift to the integers the Ozsváth–Szabó's spectral sequence [148] from reduced Khovanov homology to the Heegaard–Floer homology of the branched double cover. This remains conjectural. Analogous spectral sequences from Khovanov homology to symplectic and gauge theoretic homologies, originally defined over $\mathbb{F}_2$, have conjectural lifts to odd Khovanov homology; see [108] for a survey and [16] for an established lift.

(vi) *algebraic geometry.* Khovanov homology is closely related to the cohomology of complex Springer varieties [104, 180]. Lauda and Russel defined an oddification of this cohomology ring [117], which was related to *real* Springer varieties in [59].

## i.2 Contributions

In part I of this thesis, we show, joint with Pedro Vaz [170], that odd Khovanov homology can be constructed from an odd analogue in categorified representation theory. This gives a natural framework for its understanding, together with an extension to tangles.

### i.2.1 Graded $\mathfrak{gl}_2$-foams

Recall that the construction of odd Khovanov homology is based on the exterior algebra. As the anti-commutativity in the exterior algebra is controlled by a $\mathbb{Z}/2\mathbb{Z}$-grading, it is natural to wonder whether one could give a construction of odd Khovanov homology using a super-2-category. Ideally, the superstructure would control all signs appearing in odd Khovanov homology, that is, all interchanges of saddles. One solution, pursued by Putyra [155], is to pull back the TQFT to linear relations on cobordisms. This provides an analogue for odd Khovanov homology of Bar-Natan's "picture-world" construction of (even) Khovanov homology, restricted to links. In this context, a merge is even and a split is odd:

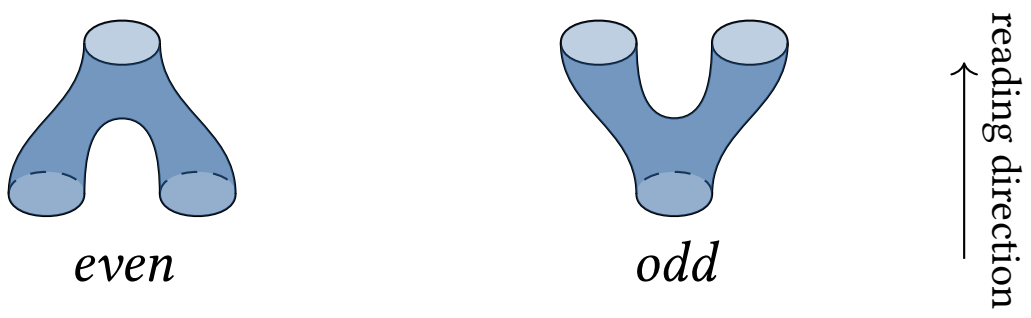

However, the superstructure only partially controls interchanges of saddles, and one still needs to use the artificial sign assignment from the original construction. Also, it does not generalize in an obvious way to tangles, as whether a saddle is a split or a merge is a global property. See however [144] for an answer to this question.

This suggests that one should look further away from the original construction. Heuristically, it is plausible that different choices of TQFTs, identical up to signs, lead to the same invariant. After all, this is the take-home message of odd Khovanov homology: sign issues on the level of the TQFT can be balanced out when choosing a sign assignment on the hypercube.

As we noted already in subsection i.1.7, an important early-day problem on Khovanov homology was also about signs. Namely, Khovanov's original construction is not properly functorial under link cobordisms, but only so up to signs. One solution to this problem was given by Blanchet [19], using $\mathfrak{gl}_2$-foams. In this thesis, we show that the same heuristics give control on signs in odd Khovanov homology. By fixing a (sort of) Morse decomposition, every $\mathfrak{gl}_2$-foams can be decomposed into a composition of the following local pictures:

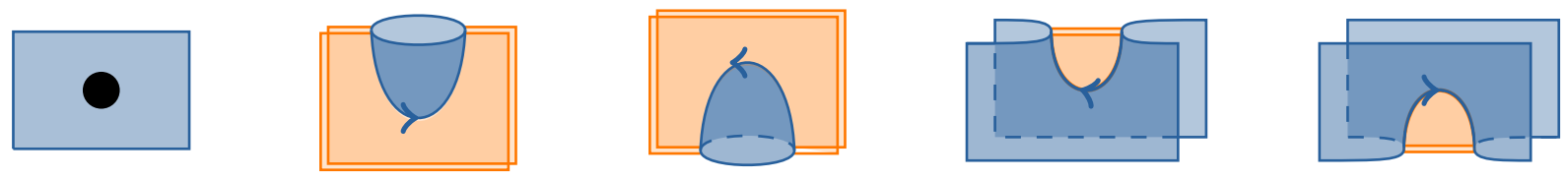

dot $(odd)$ cup $(odd)$ cap $(even)$ zip $(even)$ unzip $(odd)$

Assigning parities to these local pictures turns a $\mathfrak{gl}_2$-foam into a *super* $\mathfrak{gl}_2$-foam. The super-2-category $\mathbf{SFoam}_d$ is generated by linear combination of super $\mathfrak{gl}_2$-foams modulo relations, some of which encoding *super* isotopies. Note that (super) $\mathfrak{gl}_2$-foams can be decorated with dots: in the super case, this dot is odd, as the variable $x_i$ in the exterior algebra $\bigwedge(x_1, \ldots, x_n)$.

The super-2-category $\mathbf{SFoam}_d$ provides an odd analogue of $\mathbf{Foam}_d$ in the sense described in subsection i.1.10. In fact, both constructions fit together in a single graded-2-category $\mathbf{GFoam}_d$ of *graded* $\mathfrak{gl}_2$*-foams*, which we introduce in chapter 1. Here a *graded-2-category* is defined similarly as the super-2-category, replacing signs by generic scalars. The graded-2-category $\mathbf{GFoam}_d$ is defined over the generic ring $\mathbb{Z}[X, Y, Z^{\pm 1}]/(X^2 = Y^2 = 1)$. Setting $X = Y = Z = 1$ recovers $\mathbf{Foam}_d$, the category of $\mathfrak{gl}_2$-foams described in subsection i.1.6, while setting $X = Z = 1$ and $Y = -1$ gives a super-2-category $\mathbf{SFoam}_d$ of super $\mathfrak{gl}_2$-foams. In part II of this thesis, we use rewriting theory to describe a basis of $\mathbf{GFoam}_d$. This is summarized as follows:

**Main theorem A** (see chapter 1). *There exists a graded-2-category* $\mathbf{GFoam}_d$ *which graded-categorifies the integral form* $\mathbf{Web}_d$ *of the representation theory of* $U_q(\mathfrak{gl}_2)$.

### i.2.2 Covering $\mathfrak{gl}_2$-Khovanov homology

In chapter 2, we define an invariant of oriented tangles as a certain tensor product of complexes in $\mathbf{GFoam}_d$, following the usual scheme of a categorified Kauffman bracket (see subsection i.1.6 and Fig. i.1; recall from Fig. i.5 that $[\cdot]$ denotes the shift in quantum grading and wiggly lines the zero homological degree):

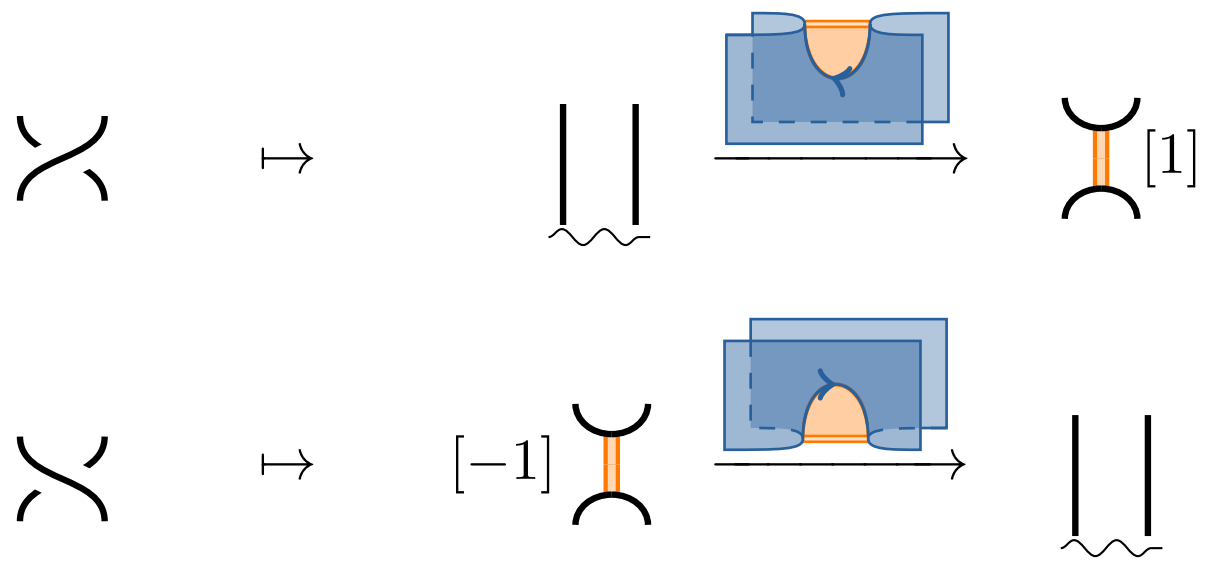

In contrast with the standard construction however, here the tensor product of complexes is understood in a *graded* sense. We define this notion properly in chapter 4 for a subclass of chain complexes called *homogeneous polycomplexes*. A *homogeneous complex* is a chain complex in a graded-2-category such that each differential is homogeneous (although the grading can differ at distinct homological degrees), and a homogeneous polycomplex is a tensor product of homogeneous complexes. This tensor product is coherent with homotopies in the following sense:

**Main theorem B** (Theorem 4.1.9)**.** *In any super-2-category, there exists a well-defined tensor product on homogeneous polycomplexes such that if $A_1^\bullet$ and $A_2^\bullet$ (resp. $B_1^\bullet$ and $B_2^\bullet$) are homotopic homogeneous polycomplexes, then so are $A_1^\bullet \otimes B_1^\bullet$ and $A_2^\bullet \otimes B_2^\bullet$.*

This result first appeared in the author's Master thesis [169]. If all differentials have a trivial grading, this recovers the usual Koszul rule for chain complexes in linear 2-categories. Although not conceptually difficult, this chapter is technical. For that reason we also give a minimal version in subsection 2.1.1, sufficient for the purpose of defining the invariant.

We then show that this invariant recovers odd Khovanov homology:

**Main theorem C** (Theorem 2.0.3)**.** *Setting $X = Z = 1$ and $Y = -1$ at the level of chain complexes, our construction coincides with odd Khovanov homology when restricted to links.*

More generally, our invariant recovers *covering Khovanov homology* when restricted to links, a homological invariant of links defined by Putyra [155] over the ring $\mathbb{Z}[X, Y, Z^{\pm 1}]/(X^2 = Y^2 = 1)$, such that setting $X = Y = Z = 1$ and $X = Z = 1$ and $Y = -1$ *at the level of chain complexes* respectively recovers (even) Khovanov homology and odd Khovanov homology:

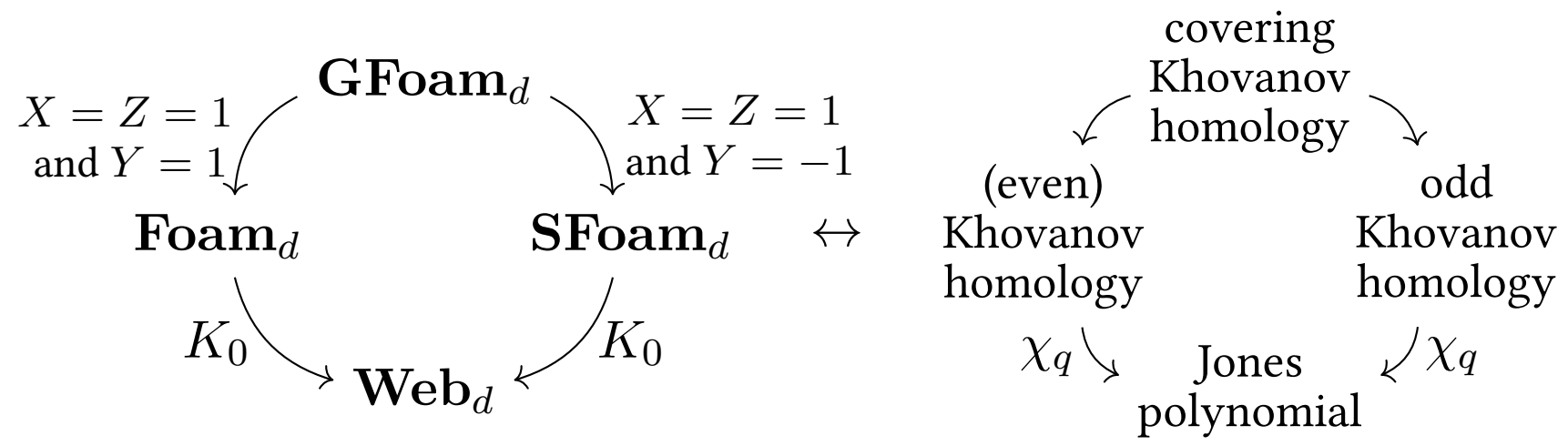

Our result may be summarized as follows:

> Slogan: *Odd Khovanov homology arises from the super-interplay of two categorified Kauffman brackets, one even and the other odd.*

In particular, this gives an extension of odd Khovanov homology to oriented tangles (another extension was given in [144]; see section i.3). This also gives a *sign-coherent* definition, where signs naturally arise from the super-2-categorical structure. This contrasts with the original definition, where signs are fixed in a somehow ad-hoc way on the hypercube of resolutions. We hope that this new definition will open the way to further connections and applications.

### i.2.3 A graded-categorification of the $q$-Schur algebra of level two

Finally, in chapter 3 we relate our construction to higher representation theory. In [183], the second author defined a superalgebra which can be seen as an odd analogue of the KLR algebra [107, 166] of level two for the $A_n$ quiver. Taking a cyclotomic quotient of this construction leads to a supercategorification of the negative half of the $q$-Schur algebra of level two. In chapter 3, we extend it to a graded-categorification of the $q$-Schur algebra of the level two (subsection i.1.8), giving a graded analogue of [127] in the level two case. We then define a *graded-foamation 2-functor*, fitting into the following commutative square:

$$\begin{array}{ccc} \mathcal{GS}_{n,d} & \longrightarrow & \mathbf{GFoam}_d \\ \big\downarrow{\scriptstyle K_0} & & \big\downarrow{\scriptstyle K_0} \\ \dot{S}_{n,d} & \longrightarrow & \mathbf{Web}_d \end{array}$$

In that way, $\mathbf{GFoam}_d$ is realised as graded-2-representation of $\mathcal{GS}_{n,d}$. In particular, the invariant defined in [183] coincides with ours, and hence with odd Khovanov homology when restricted to links. This gives a partial graded analogue to the work of Lauda, Queffelec and Rose [115]; see subsection i.1.8.

In a nutshell:

**Main theorem D.** *Setting $X = Z = 1$ and $Y = -1$ to go from "graded" to "super", The supercategorification of the $q$-Schur algebra of level two together with the super foamation 2-functor provides a representation theoretic construction of odd Khovanov homology.*

See the beginning of chapter 3 for more references.

## i.3 Perspectives

We discuss possible future directions of research. See also subsection 4.1.1 for chain complexes in graded-monoidal categories.

*Relationship to odd arc algebras*

In [144], Naisse and Putyra gave another extension of odd Khovanov homology to tangles based on arc algebras, building on previous work of Putyra [155] and Naisse–Vaz [145]. They conjectured that their construction coincides with the construction of the second author in [183]. By construction, our tangle invariant coincides with the tangle invariant defined by the second author in [183]. Following their conjecture, Naisse and Putyra's construction should coincide with ours. This remains an open question. See also the introduction of chapter 3 for further connections with their work.

*Relationship to super Lie theory*

Supercategorification is known to be related to super Lie theory. For instance, consider the case of $\mathfrak{sl}_2$, associated with the Cartan datum consisting of a single vertex. The Lie superalgebra $\mathfrak{osp}_{1|2}$ similarly arises from the Cartan super datum consisting of a single odd vertex. Their quantizations are unified by the *covering quantum group* $U_{q,\pi}(\mathfrak{sl}_2)$ defined over $\mathbb{Q}(q)[\pi](\pi^2-1)$: setting $\pi=1$ recovers $U_q(\mathfrak{sl}_2)$ while setting $\pi=-1$ recovers $U_q(\mathfrak{osp}_{1|2})$ ([90]; see also the series of papers [38, 41, 42, 43, 44, 45]). On the other hand, Ellis and Lauda [75] constructed a supercategorification of $U_{q,\pi}(\mathfrak{sl}_2)$, later reformulated (and extended to other super Cartan data) by Brundan and Ellis [25] (see also the beginning of chapter 3 for further references). Given those interactions, it was conjectured that an odd homology should correspond to a covering quantum group (resp. a Lie superalgebra), with odd Khovanov homology corresponding to $U_{q,\pi}(\mathfrak{sl}_2)$ (resp. $\mathfrak{osp}_{1|2}$). Work of Mikhaylov and Witten [140] on odd $\mathfrak{so}_{2n+1}$-link homologies arising from gauge theory further supports this claim, considering the exceptional isomorphism $\mathfrak{sl}_2 \cong \mathfrak{so}_3$. However, an explicit connection between odd Khovanov homology and $\mathfrak{osp}_{1|2}$ remains an open problem (see however [40, 60]). We conjecture that a further careful study of our construction can lead to such a connection.

*Relationship to supercategorified quantum groups*

Contrary to the classical case, our graded-categorification of the $q$-Schur algebra of level two is *not* defined as quotient of some supercategorified quantum groups. If such a connection were to exist, it would most likely arise from a certain duality on the decategorified level.

*Relationship to Heegaard–Floer homology*

Odd Khovanov homology was discovered by Ozsváth, Rasmussen and Szabó in an attempt to lift to the integers the Ozsváth–Szabó's spectral sequence from reduced Khovanov homology $\widetilde{\mathrm{Kh}}(K,\mathbb{Z}/2\mathbb{Z})$ to the Heegaard–Floer homology of the branched double cover of the knot (with orientation reversed) $\widehat{\mathrm{HF}}(-\Sigma(K),\mathbb{Z}/2\mathbb{Z})$. They conjectured that in order to lift the spectral sequence to the integers, one should replace reduced Khovanov homology $\widetilde{\mathrm{Kh}}(K,\mathbb{Z})$ with reduced *odd* Khovanov homology $\widetilde{\mathrm{OKh}}(K,\mathbb{Z})$:

$$\begin{array}{ccc} \widetilde{\mathrm{OKh}}(K,\mathbb{Z}) & \overset{?}{\rightsquigarrow} & \widehat{\mathrm{HF}}(-\Sigma(K),\mathbb{Z}) \\ \downarrow{\scriptstyle \otimes_{\mathbb{Z}}\mathbb{Z}/2\mathbb{Z}} & & \downarrow{\scriptstyle \otimes_{\mathbb{Z}}\mathbb{Z}/2\mathbb{Z}} \\ \widetilde{\mathrm{Kh}}(K,\mathbb{Z}/2\mathbb{Z}) & \rightsquigarrow & \widehat{\mathrm{HF}}(-\Sigma(K),\mathbb{Z}/2\mathbb{Z}) \end{array}$$

We expect that our construction can shed new lights on this question.

*Relationship to geometry*

The geometry underlying odd analogues is barely understood, and in particular the one underlying super $\mathfrak{gl}_2$-foams, if there exists any. In particular, there is no known super-2-representation of super $\mathfrak{gl}_2$-foams, from which one would deduce a simple proof of non-degeneracy. (This is somehow reminiscent of the fact that, as we learned from [69], there is no known geometry associated to Soergel bimodules outside of the crystallographic case.)

On a related note, Robert and Wagner [164] gave an evaluation formula for $\mathfrak{gl}_n$-foams, from which one can deduce non-degeneracy. There is no known odd analogue of an evaluation formula for super $\mathfrak{gl}_2$-foams.

*Functoriality*

We expect that our construction can lead to some sort of functoriality for odd Khovanov homology, albeit probably not in the usual sense. Indeed, isotopies already only hold up to signs in the super-2-category $\mathbf{SFoam}_d$ of super $\mathfrak{gl}_2$-foams.

# Introduction to part II

## *higher rewriting theory and categorification*

We give an introduction to part II; see the beginning of part I for a quick overview of the whole thesis. We recommend casual readers to start with the overview of the field (section ii.1), while experts on rewriting theory may prefer to read the contributions and perspectives directly (sections ii.2 and ii.3).

## ii.1 Overview of the field

The casual reader should read subsection i.1.3 on higher algebra first. References will also be made to subsections i.1.4, i.1.5, i.1.10 and i.2.1.

### ii.1.1 Rewriting theory

The word problem for monoids asks the following question: given a presented monoid $G$ with generators in the set $\mathsf{X}$, is there an algorithm that decides whether two words with letters in $\mathsf{X}$ are equal as elements of $G$? While known to be undecidable in general [131, 154], one can hope to solve the word problem in practical cases. *Rewriting theory* suggests the following method, which we illustrate with the symmetric group on three strands (understood as a monoid), defined using its Coxeter presentation:

$$\mathfrak{S}_3 = \langle \sigma, \tau \mid \sigma\sigma = 1, \tau\tau = 1, \sigma\tau\sigma = \tau\sigma\tau \rangle \,.$$

A *rewriting system* consists of a choice of orientations on the defining relations. For instance:

$$\mathsf{P}_{\mathfrak{S}_3} := (\mathsf{X}, \mathsf{R}) : \quad \mathsf{X} = \{\sigma, \tau\} \quad \mathsf{R} = \{\sigma\sigma \to 1, \tau\tau \to 1, \sigma\tau\sigma \to \tau\sigma\tau\} \quad \text{(ii.1)}$$

Let $\mathsf{X}^*$ denote the set of words in letters in $\mathsf{X}$. Any rewriting system defines a non-deterministic algorithm, where for each oriented relation $A \to B$ and words $x, y \in \mathsf{X}^*$, the algorithm may perform the reduction $xAy \to xBy$, but not the reduction $xBy \to xAy$. In the terminology of rewriting theory, $xAy \to xBy$ is called a *rewriting step*, and a successive composition of rewriting steps, denoted $a \overset{*}{\to} b$, is called a *rewriting sequence*.

A pair of co-initial rewriting sequences $(f : a \overset{*}{\to} b, g : a \overset{*}{\to} b')$ is called a *branching*, and a pair of a co-terminal rewriting sequences $(f' : b \overset{*}{\to} c, g' : b' \overset{*}{\to} c)$ is called a *confluence*. A branching that admits a confluence is said to be *confluent*. This is illustrated in the diagrams below, where plain arrows (resp. dotted arrows) denote branchings (resp. confluences):

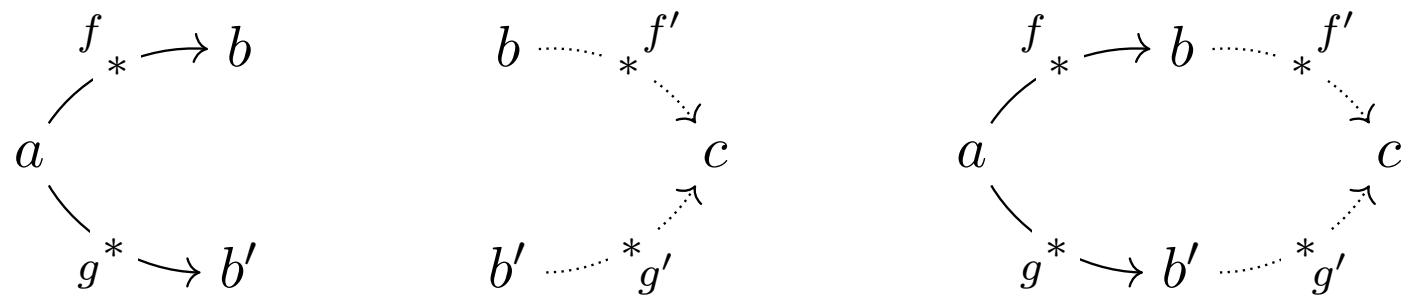

a branching a confluence a confluent branching

Figure ii.1 gives examples of confluent branchings in $\mathsf{P}_{\mathfrak{S}_3}$.

To solve the word problem, we require two key properties: *termination*, which postulates that any rewriting sequence terminates, and *confluence*, which postulates that every branching is confluent. A rewriting system which is both terminating and confluent is called *convergent*. In that case, the non-deterministic algorithm always produces an output, and this output is always the same. Given a word as input, we call the output its *reduced expression*. It is not so hard to prove that under convergence, two words are equal in the associated monoid if and only if they have identical reduced expressions. This provides a solution to the word problem.

We are left with the problem of showing that indeed, the rewriting system $\mathsf{P}_{\mathfrak{S}_3}$ is convergent. Termination is not hard: it can be shown using a suitable partial order. On the other hand, confluence must in principle be checked for *every* branching. A branching $(f, g)$ for which $f$ and $g$ are rewriting *steps* (in contrast to rewriting *sequences*) is said to be *local*. For instance, the three branchings in Fig. ii.1, depicted in plain arrows, are local. *Newmann's*

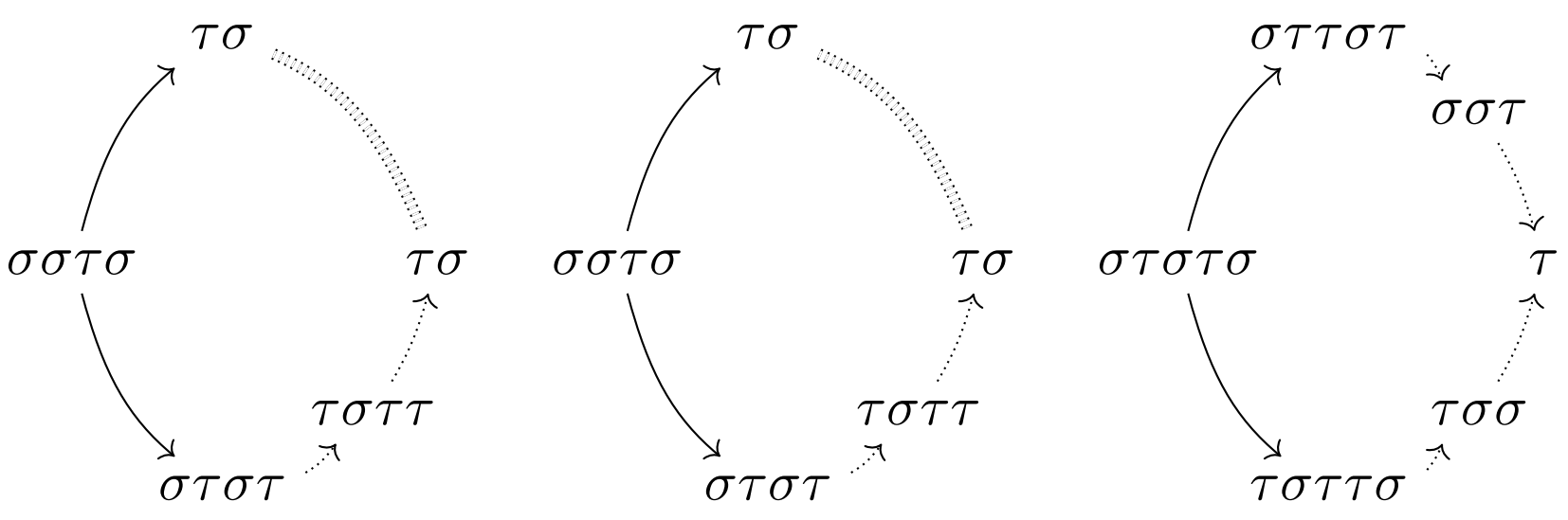

**Fig. ii.1** Critical branchings in $\mathsf{P}_{\mathfrak{S}_3}$.

*lemma* (Lemma 6.2.5) states that assuming termination, confluence follows from confluence of local branchings.

Still, many local branchings remain, as the following ones (for any two words $x, y \in \mathsf{X}^*$):

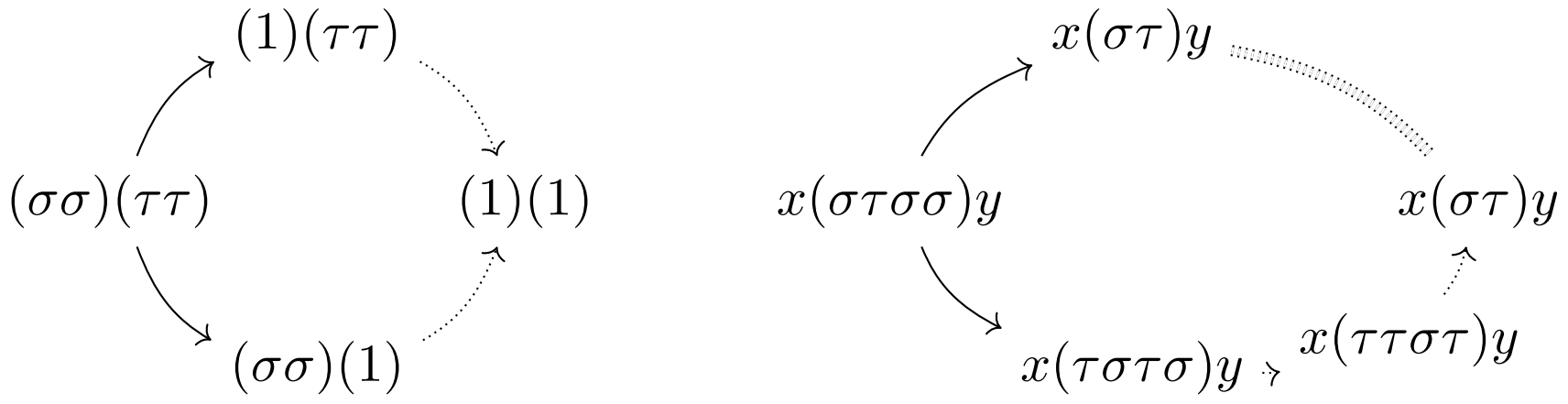

independent branching

contextualized branching

Each of them admits a “somehow canonical” confluence. The first branching is an *independent branching*: intuitively, it consists of two rewriting steps that do not interact with each other. The second branching is a *contextualized branching*: it is of the form $x(f, g)y$, for $(f, g)$ the first local branching in Fig. ii.1, and the confluence of $x(f, g)y$ is canonically induced from the confluence of $(f, g)$. A *critical branching* is a local branching which is neither an independent branching nor a contextualized branching. As suggested by our discussion, confluence of local branchings follows from confluence of critical branchings. The latter constitute the minimal amount of computations one needs to perform in a given situation; the rest follows from general considerations. The reader may convince themselves that Fig. ii.1 gives a complete list of critical branchings in $\mathsf{P}_{\mathfrak{S}_3}$. It follows that the rewriting system $\mathsf{P}_{\mathfrak{S}_3}$ is both terminating and locally confluent, and hence convergent by Newmann’s lemma.

The rewriting approach can be considered in any algebraic structure, with both general and context-dependent methods to find convergent rewriting systems. While Newmann's lemma applies to all situations, one must leverage the specifics of the given algebraic structure to further simplify the study of confluence. This often leads to a structure-dependent notion of critical branchings, together with a "critical branching lemma" that reduces local confluence to confluence of critical branchings.

For commutative algebras, convergent rewriting systems are known as Gröbner bases [27]. Shirshov [174] used rewriting techniques to give an explicit construction of Poincaré–Birkhoff–Witt bases in Lie algebras. Bokut' and Bergman then generalized these works to associative algebras [18, 20], where the critical branching lemma is known as the composition lemma and the diamond lemma, respectively. See also [86] for a modern approach on rewriting in associative algebras.

### ii.1.2 Polygraphs

What is a presentation of a (strict and small) $n$-category? As a first example, consider again the rewriting system $\mathsf{P}_{\mathfrak{S}_3} = (\mathsf{X}, \mathsf{R})$ defined in (ii.1). Viewing a monoid as a category with a single object $\{*\}$, the set $\mathsf{X}$ consists of generating 1-cells with source and target the object $\{*\}$. We encapsulate the latter fact with (trivial) source and target maps $s_0, t_0 \colon \mathsf{X} \to \{*\}$. In this perspective, the set of words $\mathsf{X}^*$ is the free category generated by $\mathsf{X}$. Similarly, the set $\mathsf{R}$ consists of generating 2-cells, and we define maps $s_1, t_1 \colon \mathsf{R} \to \mathsf{X}^*$ setting $s_1(r) = A$ and $t_1(r) = B$ for each oriented relation $r \colon A \to B$. Reformulated in this way, $\mathsf{P}_{\mathfrak{S}_3}$ defines the data of a *2-polygraph.* More generally, an $n$-category can be presented by an $(n+1)$*-polygraph*, with generating $(k+1)$-cells defined on the free $k$-category generated by the lower cells:

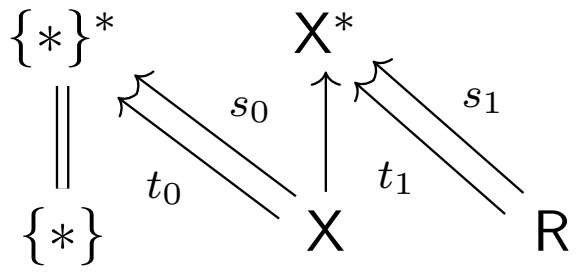

the 2-polygraph $\mathsf{P}_{\mathfrak{S}_3}$ presenting $\mathfrak{S}_3$

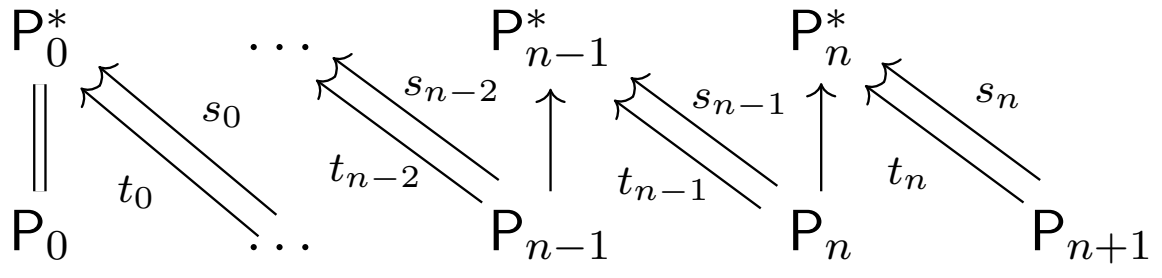

an $(n+1)$-polygraph $\mathsf{P}$ presenting an $n$-category

See Table ii.1 for a summary of low-dimensional $n$-polygraphs. The $n$-category presented by a $(n+1)$-polygraph $\mathsf{P}$ is obtained by quotienting $\mathsf{P}_n^*$ by the relation $s_n(r) = t_n(r)$ for each $r \in \mathsf{P}_{n+1}$.

| *structure* | *presentation* |
|---|---|
| set | 1-polygraph |
| category (monoid) | 2-polygraph (with one object) |
| 2-category | 3-polygraph |

**Table ii.1** Low-dimensional $n$-polygraphs

| *structure* | *presentation* |
|---|---|
| module | linear 1-polygraph |
| linear category (associative algebra) | linear 2-polygraph (with one object) |
| linear 2-category | linear 3-polygraph |

**Table ii.2** Low-dimensional linear $n$-polygraphs

Polygraphs were first introduced by Street [178], under the name of *computads*; the term *signatures* also appears in the literature. Polygraphs were independently introduced by Burroni [28] to study generalizations of the word problem. He introduced the polygraphic terminology, now standard in the rewriting community.

### ii.1.3 Higher rewriting theory

Polygraphs can be thought as higher-dimensional rewriting systems. Consider once again the 2-polygraph $\mathsf{P}_{\mathfrak{S}_3} = (\mathsf{X}, \mathsf{R})$ presenting $\mathfrak{S}_3$ from (ii.1). A word in $\mathsf{X}$ is nothing else than a path of generating 1-cells in $\mathsf{X}$:

$$\sigma\tau\sigma \in \mathsf{X}^* \quad \leftrightarrow \quad * \xrightarrow{\sigma}_{\mathsf{X}} * \xrightarrow{\tau}_{\mathsf{X}} * \xrightarrow{\sigma}_{\mathsf{X}} *$$

In this way, each element of $\mathsf{X}$ is thought as a 1-dimensional rewriting step, and each word in $\mathsf{X}^*$ as a rewriting sequence. With this point of view, a generating 2-cell in $\mathsf{R}$ is nothing else than a generating 2-dimensional rewriting step between 1-dimensional rewriting sequences:

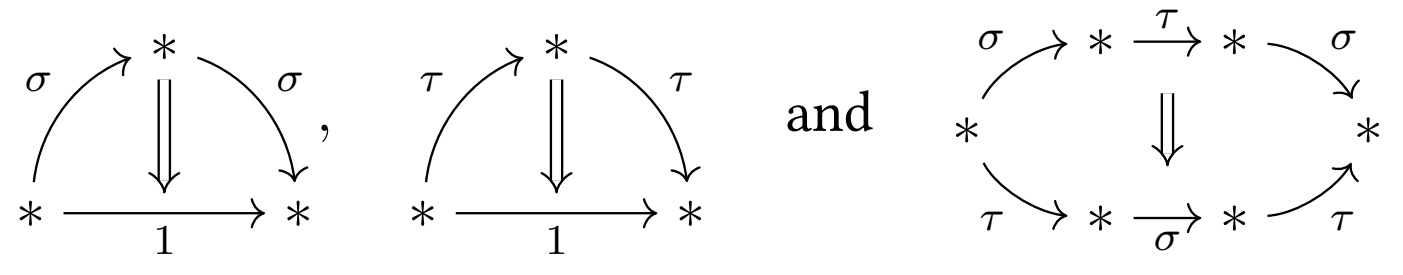

Similarly, in an $n$-polygraph one thinks of a $(k+1)$-cell as a generating $(k+1)$-dimensional rewriting step between $k$-dimensional rewriting sequences.

Higher rewriting theory has witnessed increasing interests in the last two decades, prominently by the French school [85, 87, 88, 89, 112, 142]. As a major application, let us note the construction of a homology theory for $\omega$-categories, based on polygraphic resolutions [84, 113, 114, 139, 141]. For further details, we refer the reader to the recent monograph on the subject [5].

In general, the rewriting theory associated to an $(n+1)$-polygraph is called $(n+1)$*-dimensional rewriting theory*; it presents an $n$-category. The 1-dimensional case corresponds to rewriting in sets; it is also known as *abstract rewriting*. As we have seen, 2-dimensional rewriting theory corresponds to rewriting in categories, and monoids in particular. In this case, a relation $r\colon A \to B$, i.e. a generating 2-cell, can always be composed on the left and on right with 1-cells $x$ and $y$ respectively, leading to a new relation

$$x \star_0 r \star_0 y\colon x \star_0 A \star_0 y \to x \star_0 B \star_0 y.$$

This process is called *contextualization*. We call the data of $x$ and $y$ a *context*, denoted with the letter $\Gamma$, and write $\Gamma[r]$ for the relation $x \star_0 r \star_0 y$. In monoids, contextualization amounts to multiplying a relation on the left and on the right with words, as we have seen already in subsection ii.1.1.

Three-dimensional rewriting is rewriting in 2-categories. In this case, a relation $r\colon A \to B$ is a generating 3-cell, and contextualization amounts to first composing horizontally with 1-cells $x$ and $y$, and then vertically with 2-cells $\alpha$ and $\beta$:

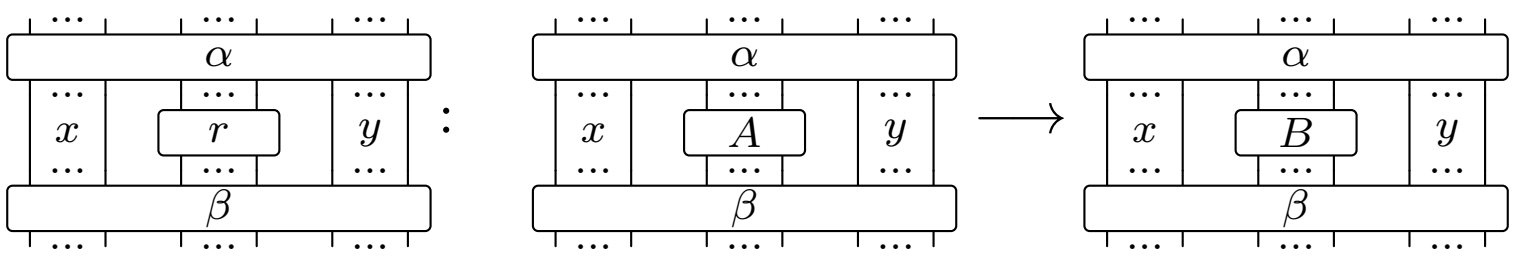

We similarly write $\Gamma[r]$ for a contextualization of the relation $r$. In what follows, we will rarely discuss the cases $n > 3$, and often say *higher rewriting theory* to refer to the case $n = 3$. The toolbox of higher rewriting theory resembles the one of rewriting in monoids: Newmann's lemma is (as always) applicable, and one has similar notions of independent and contextualized branchings.

### ii.1.4 Linear rewriting

Linear rewriting is the algorithmic study of presented modules. Given a commutative ring $\Bbbk$, a $\Bbbk$-module $\mathtt{M}$ is presented by a set $\mathcal{B}$, together with a set $\mathtt{R}$ of oriented relations on $\langle \mathcal{B} \rangle_{\Bbbk}$, the free $\Bbbk$-module generated by $\mathcal{B}$. This is encapsulated by the data of a *linear 1-polygraph*:

$$\langle \mathcal{B} \rangle_{\Bbbk} \overset{s}{\underset{t}{\leftleftarrows}} \mathtt{R}.$$

Similarly, *linear 2-polygraphs* present linear categories, and in particular associative algebras [86]; and linear $(n+1)$-polygraphs present linear $n$-categories [3]. See Table ii.1 for a summary of low-dimensional linear $n$-polygraphs. The associated linear $n$-dimensional rewriting theory was studied by Guiraud, Hoffbeck and Malbos in the case $n = 2$ [86], and by Alleaume in the case $n = 3$ [3], with application to the oriented affine Brauer algebra.

The most classical setting for linear rewriting is commutative algebras. For that reason, we call *monomials* the elements of $\mathcal{B}$. An oriented relation $r\colon s(r) \to t(r) \in \mathtt{R}$ is assumed to be of the form "rewrite a monomial into a linear combination of monomials", that is:

$$r\colon b \to_{\mathtt{R}} \lambda_1 b_1 + \ldots + \lambda_n b_n \qquad \text{for } \lambda_1, \ldots, \lambda_n \in \Bbbk, b, b_1, \ldots, b_n \in \mathcal{B} \quad \text{(ii.2)}$$

We say that relations in $r$ are *left-monomial*.

A generic rewriting step is of the form

$$\lambda r + v\colon \lambda s(r) + v \dashrightarrow_{\mathtt{R}} \lambda t(r) + v \quad \lambda \in \Bbbk \setminus \{0\}, v \in \langle \mathcal{B} \rangle_{\Bbbk},$$

and it is said to be *positive* if the monomial $s(r) \in \mathcal{B}$ does not appear in the linear decomposition of $v$. For example, if $\mathcal{B} = \{a, b, c\}$ and $r\colon a \to b + c \in \mathtt{R}$ is an oriented relation, then both

$$2a + b = a + (a + b) \dashrightarrow_{\mathtt{R}} (b + c) + (a + b) \quad \text{and} \quad 2a + b \to_{\mathtt{R}} 2(b + c) + b$$

are rewriting steps, but only the latter is positive (note the use of dashed and plain arrows to differ between the two). To avoid rewriting loops, one must restrict to positive rewriting steps. For instance, in the previous example:

$$0 = a - a \dashrightarrow_{\mathtt{R}} (b + c) - a \dashrightarrow_{\mathtt{R}} (b + c) - (b + c) = 0.$$

Positivity may look like a minor modification of the theory; in fact, it constitutes the main difficulty of the linear setting. We shall say more about that in subsection ii.2.2.

## ii.1.5 Rewriting modulo

In many practical situations, enforcing all relations to be oriented is too restrictive. Instead, one may wish to rewrite with a set of oriented relations $\mathbb{R}$,[1] modulo another set of *un*oriented relations $\mathbb{E}$. More precisely, the working data of *abstract rewriting modulo* is given by two 1-polygraphs

$$X \underset{t}{\overset{s}{\leftleftarrows}} \mathbb{R} \quad \text{and} \quad X \underset{t}{\overset{s}{\leftleftarrows}} \mathbb{E},$$

defined on the same underlying set $X$. Denote by $\mathbb{E}^\top = (\mathbb{E} \cup \mathbb{E}^{-1})^*$ the free groupoid generated by $\mathbb{E}$. Intuitively, relations in $\mathbb{E}^\top$ are *un*oriented rewriting sequences in $\mathbb{E}$. In this context, a *rewriting step modulo* is a composition $e' \circ r \circ e$ with $r \in \mathbb{R}$ and $e, e' \in \mathbb{E}^\top$:

$$\bullet \overset{e}{\underset{\mathbb{E}}{\rightsquigarrow}} \bullet \overset{r}{\underset{\mathbb{R}}{\longrightarrow}} \bullet \overset{e'}{\underset{\mathbb{E}}{\rightsquigarrow}} \bullet$$

In other words, in between a rewriting step in $\mathbb{R}$ one can apply an arbitrary number of relations in $\mathbb{E}$, in any direction. The data $\mathbb{S} = (\mathbb{R}, \mathbb{E})$ defines an *abstract rewriting system modulo*, and a rewriting step modulo as above is called an $\mathbb{S}$-*rewriting step.*

Typically, the modulo data $\mathbb{E}$ will consists of relations thought as being "structural", in the sense of being part of some underlying algebraic structure. For instance, rewriting in commutative algebras is implicitly rewriting in associative algebras modulo commutativity ($\mathbb{E}$ consists of relations $xy \to yx$ for all monomials $x$ and $y$).

Rewriting modulo allows an inductive approach to the word problem. Indeed, instead of trying to fit the relations $\mathbb{R} \sqcup \mathbb{E}$ into a single convergent rewriting system, one can instead show that $\mathbb{R}$ is convergent "modulo $\mathbb{E}$" on one hand, and that $\mathbb{E}$ is convergent on the other hand. This is especially important when dealing with intricate algebraic structures such as the ones appearing in higher algebra.

Different variants of rewriting modulo have been developed, often with more restrictive modulo rules than the ones described above [92, 95, 96, 130, 152, 184]. Rewriting modulo is used in [35] to study confluence in Lawvere theories. A higher analogue to rewriting modulo was introduced in [58], using the formalism of double categories. In the linear setting, Dupont extended

[1] We use blackboard font to refer to abstract rewriting; this is unrelated to the set of real numbers.

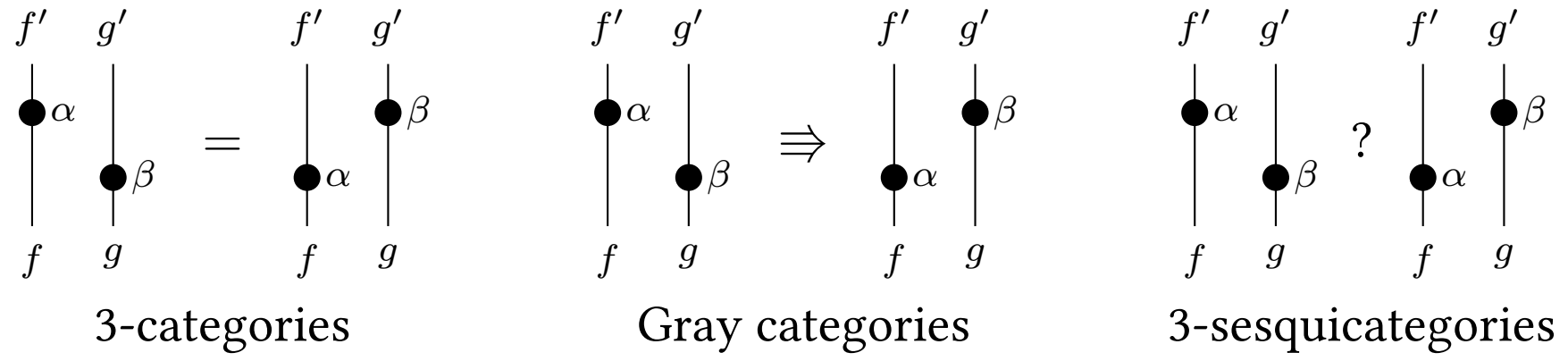

**Fig. ii.2** In a 3-category, a Gray category or a 3-sesquicategory, the interchange law for 2-morphisms respectively holds strictly, holds weakly via interchangers, or does not hold (a priori).

Alleaume's approach [3] modulo in order to rewrite modulo pivotality in 2-categories [56]. He then applied his theory to categorified quantum groups [55]. See also [66] for related investigations. An approach to rewriting modulo in super-2-categories was also proposed in [57], motivated by the study of supercategorified quantum $\mathfrak{sl}_2$ (see also subsection ii.1.7).

## ii.1.6 Gray rewriting theory

Starting with $n = 3$, the strict and weak notions of an $n$-category start to diverge: while a bicategory is alway equivalent to a 2-category, not every tricategory is equivalent to a 3-category. However, every tricategory is equivalent to a *Gray category* [83] (see Remark A.1.1 for an understanding of this fact). In that sense, Gray categories provide a simpler notion than tricategories, while retaining their expressivity:

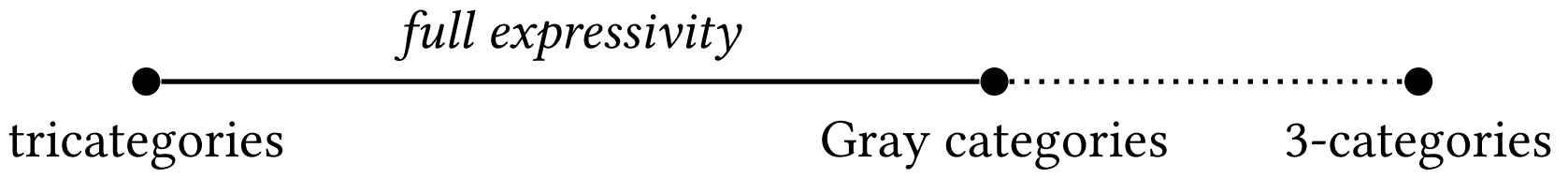

In Gray categories, associativity and unitality hold strictly, as for 3-categories. However, the interchange law for 2-morphisms only holds up to certain coherent 3-morphisms, called *interchangers*; see Fig. ii.2.

In [81], Forest and Mimram initiated the study of rewriting theory for *weak* $n$-categories, starting with Gray categories. It turns out to be easier to consider the more general framework of *$n$-sesquicategories* [6, 81, 177] (called *$n$-precategories* in [81]). As for Gray categories, $n$-sesquicategories are strictly associative and unital structures. However, they do not contain any coherence data for interchange laws; see Fig. ii.2. One can understand an $n$-sesquicategory as an "unfinished definition" for what a semistrict $n$-category, analogous to the notion of Gray category for $n = 3$, would be. Sesquicategories play a prominent role in the graphical proof-assistant *Globular* and its successor *homotopy.io* [11, 12, 51, 162].

Recall that an $n$-category is presented by an $(n+1)$-polygraph (see subsection ii.1.2); similarly, an $n$-sesquicategory is presented by an $(n+1)$-sesquipolygraph [81]. A 0-, 1- and 2-sesquipolygraph is the same as 0-, 1- and 2-polygraph; the distinction only starts to appear with 3-sesquipolygraphs. In particular, a 2-sesquicategory is presented by a 3-sesquipolygraph. A 3-sesquipolygraph which contains generating interchangers is called a *Gray polygraph*; in this case, it presents a 2-category. Note that both Gray polygraphs and 3-polygraphs present 2-categories. However, the former explicitly contains interchangers as generating 3-cells, which alters the associated rewriting theory.

## ii.1.7 Non-degeneracy of graded-2-categories

Recall the notion of super-2-categories from subsection i.1.10. Graded-2-categories are defined analogously, replacing signs with scalars. More precisely, let $\Bbbk$ be a commutative ring with $\Bbbk^\times$ its group of invertible scalars, $G$ an abelian group and $\mu\colon G \times G \to \Bbbk^\times$ a $\mathbb{Z}$-bilinear form. A $(G,\mu)$-*graded-2-category* is a certain categorical structure akin to a $G$-graded linear 2-category, but where the interchange law is twisted by $\mu$ (see section 1.1):

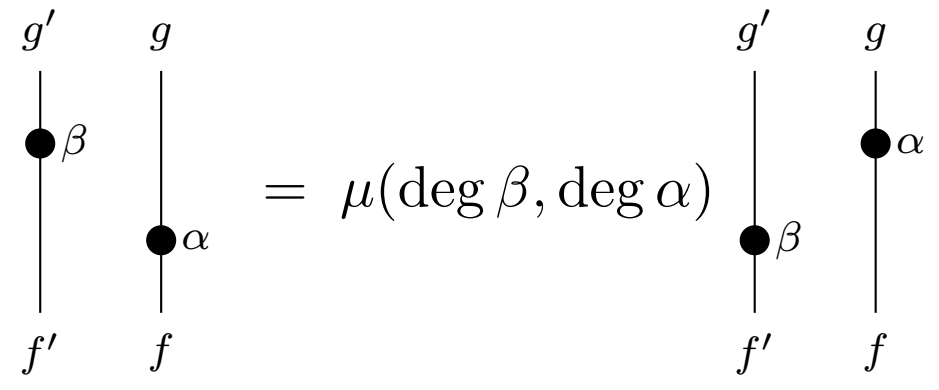

(ii.3)

Here $\beta$ and $\alpha$ are assumed to be homogeneous, and $\deg(-)$ denotes the $G$-grading.

Recall from subsections i.1.3 to i.1.5 the idea of categorification in representation theory. The complexity of the presented 2-categories appearing in these contexts makes the computation of bases challenging. More precisely, while one often has a candidate basis, proving its linear independence is difficult, and sometimes left as a "non-degeneracy conjecture". Computing the Grothendieck ring, and hence showing that the given 2-category categorifies the expected algebra, usually relies on this non-degeneracy conjecture. For instance, while Khovanov and Lauda gave a definition of $\mathcal{U}(\mathfrak{g})$ that holds in any type $\mathfrak{g}$ [106], they could only show its non-degeneracy in type $\mathfrak{g} = \mathfrak{sl}_n$. The general case was shown a decade later by Webster [188] (building on [99, 187]), thereby showing that, as expected, $\mathcal{U}(\mathfrak{g})$ categorifies $U_q(\mathfrak{g})$. For $\mathfrak{g} = \mathfrak{sl}_n$, Khovanov and Lauda's proof of non-degeneracy relies on a certain 2-representation on the cohomology ring of flag varieties. In most cases, non-degeneracy results similarly rely on a suitable 2-representation.

As we noted already in subsection i.1.10, super-2-categories, and to a lesser extent graded-2-categories, have started to appear in categorification. In this new setting, non-degeneracy is an even harder challenge. On the one hand, one rarely has a 2-representation to rely on (see however [26]). On the other hand, while graded analogues have similar string diagrammatics with interpretations as appropriate isotopies, these isotopies only hold *up to scalars*. This greatly complicates a formal analysis of their defining relations.

## ii.2 Contributions

In part II of this thesis, we develop the theory of *linear Gray rewriting modulo*, extending the work of Forest and Mimram [81] to incorporate both linear rewriting and rewriting modulo. This is done in chapters 5 and 6. Our approach is also inspired by Alleaume's approach to rewriting in linear 2-categories [3, 4], and Dupont's approach to rewriting modulo in linear 2-categories [54, 55, 56], although with significant differences in the implementation. Indeed, a major part of our work consists in defining a suitable theory of (1-dimensional) linear rewriting modulo.

Our theory applies in particular to graded-2-categories. In chapter 1 (see also subsection i.2.1), a certain graded-2-category of $\mathfrak{gl}_2$-foams $\mathbf{GFoam}_d$ is introduced, used throughout part I of this thesis. In chapter 7, we use Gray rewriting theory modulo to show that $\mathbf{GFoam}_d$ is non-degenerate (see Corollary 1.6.4).

### ii.2.1 Linear Gray polygraphs

In chapter 5, we define *linear $n$-sesquicategories* and *linear $n$-sesquipolygraphs* as direct linear analogues of $n$-sesquicategories and $n$-sesquipolygraphs [81]. We call *graded interchangers* the 3-cells capturing the relation (ii.3). A *linear*

*Gray polygraph* is defined as a linear 3-sesquipolygraph which contains its own graded interchangers. We get a notion of presentation for a graded-2-category, suitable for rewriting theory:

**Main definition A** (Definition 5.4.4)**.** *A* presentation *of a graded-2-category is a linear Gray polygraph.*

The chapter starts with a thorough example (section 5.1): we expect it to be sufficient for the impatient reader.

## ii.2.2 Higher linear rewriting modulo

Chapter 6 is the technical heart of part II, where the rewriting theory is developed. Although our goal contains the respective goals of [3, 55, 56, 57], our approach has important differences at various stage of the theory, even when restricting to their respective settings. For that reason, chapter 6 starts from first principles, and is essentially self-contained:

- In section 6.1, we define convergence modulo and explain how it relates to the word problem. We also define the notion of *coherence modulo* and relate it to convergence modulo. This motivates the study of convergence modulo.
- In section 6.2, we give the foundations of *abstract* (i.e. 1-dimensional) *rewriting modulo*. The working data $\mathbb{S} = (\mathbb{R}, \mathbb{E})$ is given by two 1-polygraphs $\mathbb{R}$ and $\mathbb{E}$ defined on the same underlying set $X$, as in subsection ii.1.5.
- In section 6.3, we give the foundations of *linear* (1-dimensional) *rewriting modulo*. The working data $\mathtt{S} = (\mathtt{R}, \mathtt{E})$ is given by two linear 1-polygraphs $\mathtt{R}$ and $\mathtt{E}$ defined on the same underlying set of monomials $\mathcal{B}$, with the notations of subsection ii.1.4. There are additional conditions, including that $\mathtt{R}$ is left-monomial (see (ii.2)), and that relations in $\mathtt{E}$ are of the form $b \to \lambda b'$ for some $\lambda$ an invertible scalar; we say that $\mathtt{E}$ is *monomial-invertible*.
- In section 6.4, we give the foundations of weak 3-dimensional rewriting modulo, simply called *higher rewriting modulo*. The working data $\mathsf{S} = (\mathsf{R}, \mathsf{E})$ is given by two 3-sesquipolygraphs $\mathsf{R}$ and $\mathsf{E}$ defined on the same underlying 2-polygraph.

- In section 6.5 we combine all of the above to define weak 3-dimensional linear rewriting modulo, simply called *higher linear rewriting modulo.* The working data $\mathsf{S} = (\mathsf{R}, \mathsf{E})$ is given by two linear 3-sesquipolygraphs $\mathsf{R}$ and $\mathsf{E}$ defined on the same underlying 2-polygraph, with further conditions similar to the ones for 1-dimensional linear rewriting modulo. If it contains graded interchangers, we call it *linear Gray rewriting modulo.* The latter allows rewriting in graded-2-categories.

Finally, section 6.6 is a "user's guide to applications", which we illustrate with a simple example (subsection 6.6.3).

*Why combining higher, linear and modulo is hard*

As discussed in subsection ii.1.3, the difference between 1-dimensional rewriting and higher rewriting lies in the existence of contexts. For instance, in monoids one can multiply an element $A$ on the left and on the right, giving a contextualization $\Gamma[A]$ of $A$. One can think of the contextualization $\Gamma[-]$ as some sort of action. Note that for monoids, it acts faithfully: $\Gamma[A] = \Gamma[B]$ implies $A = B$, for any two words $A$ and $B$. The same is true for contextualization in $n$-sesquipolygraphs.

Now assume we work modulo $\mathsf{E}$; for instance, rewriting in monoids modulo commutativity. We write $A \sim_{\mathsf{E}} B$ if the cells $A$ and $B$ are related by unoriented relations in $\mathsf{E}$. A priori, there is no guaranty that contextualization remains faithful modulo, i.e. that $\Gamma[A] \sim_{\mathsf{E}} \Gamma[B]$ implies $A \sim_{\mathsf{E}} B$. The reader can convince themselves that in the case of rewriting in monoids modulo commutativity, contextualization is in fact faithful.

However, *contextualization modulo interchangers is not faithful.* Let us explain this statement in more details. The weak approach to higher rewriting emphasizes that rewriting in 2-categories is best understood as rewriting in 2-sesquicategories modulo interchangers. Consider a 2-polygraph $\mathsf{R} = (\mathsf{R}_0, \mathsf{R}_1, \mathsf{R}_2)$ with $\mathsf{R}_0 = \{*\}$, $\mathsf{R}_1 = \{ | \}$ and $\mathsf{R}_2 = \{ \bullet , \frown \}$ and write $\mathsf{RGray}$ the 3-sesquipolygraph of interchangers on $\mathsf{R}$. Then:

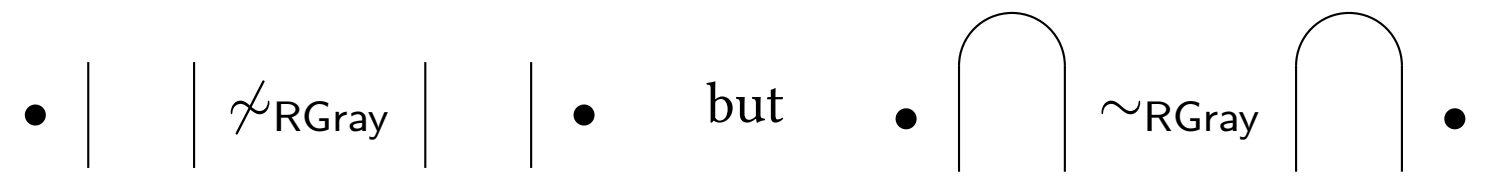

(ii.4)

This elementary example shows that when rewriting modulo interchangers, and hence in many practical cases, one cannot expect contextualization to act faithfully.

This fact is especially important in the linear setting. Indeed, recall from subsection ii.1.4 the importance of considering positive rewriting steps. If contextualization is not faithful, then *the contextualization* $\Gamma[r]$ *of a positive rewriting step* $r$ *needs not be positive itself.* To solve this problem, we introduce the notion of tamed congruence.

### *Tamed congruence*

The instability of the positivity condition is a common issue in linear rewriting. For instance, the canonical confluence of an additive branching (a canonical type of branchings in linear rewriting) needs not be positive; and similarly for the canonical confluences of an independent branching. In standard approaches (e.g. [86]), these issues are dealt with “inside proofs”, keeping positive confluence as the central concept of the theory.

However, as we generalize to higher settings, it becomes increasingly difficult to hide these issues in “black-boxed” statements. This is especially the case when adding non-faithful contextualization into the mix, as now the canonical confluence of a contextualized positively confluent branching needs not be positive. For that reason, we introduce the notion of *tamed congruence* (Definition 6.2.6) as a replacement for confluence, and build the theory of linear rewriting modulo around it. As we demonstrate, tamed congruence behaves better than confluence with respect to positivity: in each of the three cases given above, the canonical confluence is a tamed congruence (see respectively Theorem 6.3.24, Lemma 6.5.10 and Corollary 6.5.8).

In fact, tamed congruence happens to also behave better abstractly. Roughly speaking, the study of tamed congruence allows more flexibility than the study of confluence. Precisely, given four rewriting sequences defined on the same source, if any three pairs are tamed congruent, then any pair is tamed congruent (see Lemma 6.2.11). This flexibility is particularly useful for higher (linear) rewriting modulo.

### *Context-dependent termination*

In the higher setting, one often needs a *context-dependent* restriction on rewriting steps to ensure termination. For instance, Consider again the example (ii.4), together with the rewriting step $r$ below (this situation appears in practice [3]):

$$r : \bullet \; \Big| \; \to_{\mathsf{R}} \; \Big| \; \bullet$$

Allowing $\Gamma[r]$ as a rewriting rule for all contexts $\Gamma$ will result in rewriting loops. To avoid this obstruction, one needs a context-dependent termination, that is, a choice of rewriting rules that is sensitive to contexts. As one can expect, this complicates the study of confluence. Indeed, if a branching $(f, g)$ admits a confluence $(f', g')$, this does not imply that $\Gamma[f', g']$ is confluent, since $\Gamma[f']$ and $\Gamma[g']$ may not be valid rewriting sequences.

In section 6.4 and section 6.5, we develop the relevant techniques to deal with context-dependent termination.

*Rewriting modulo invertible scalar*

To account for rewriting in graded-2-categories, one requires modulo rules that can add invertible scalars. In particular, the modulo data may not be coherent; that is, two parallel paths of modulo rules may lead to distinct scalar. This happens in practice; for instance, it is the case of graded $\mathfrak{gl}_2$-foams (see Remark 6.3.5).

One must be extremely careful when working with non-coherent modulo rules. Indeed, the way the modulo rule is applied matters. In other words, rewriting with a non-coherent modulo cannot be reduced to rewriting on modulo equivalence classes.

For that reason, a careful understanding of (non-)coherence of the modulo data is imperative. In subsection 6.6.2, we explain how to proceed in practice, and detail the approach is several examples (subsection 6.6.3, section 7.1), including graded $\mathfrak{gl}_2$-foams (section 7.3).

### ii.2.3 Special cases

Our theory has three main features: it can be *(weakly) higher*, it can be *linear* and it can be *modulo*. In other words, it allows to rewrite modulo in 2-sesquicategories and in linear 2-sesquicategories. We use the terminology *strict* to emphasize that a rewriting theory does not allow modulo. The following are special cases of higher (linear) rewriting modulo:

- Strict rewriting in 2-sesquicategories, and in particular strict rewriting containing interchangers, recovers Forest and Mimram's Gray rewriting theory [81].
- Strict rewriting in 2-categories (see subsection ii.1.3) is rewriting in 2-sesquicategories modulo interchangers.

- Strict rewriting in linear 2-categories is rewriting in linear 2-sesquicategories modulo interchangers. Alleaume [3] also developed a rewriting theory for linear 2-categories.
- Rewriting modulo in linear 2-categories is rewriting modulo in linear 2-sesquicategories with a modulo containing interchangers. Dupont [56] also developed a rewriting modulo theory for linear 2-categories.
- Rewriting modulo in graded-2-categories is rewriting modulo in linear 2-sesquicategories with a modulo containing graded interchangers. Dupont, Ebert and Lauda [57] also developed a rewriting modulo theory for super-2-categories.
- Although we do not describe it explicitly, our work can be used to define rewriting modulo in monoids and associative algebras. In the latter case and without the modulo, this recovers the work of Guiraud, Hoffbeck and Malbos [86].

We note that linear rewriting modulo holds over any commutative ring; one needs not restrict to a field.

*Limitation*

We stress that in our approach, the modulo data in linear rewriting modulo must be monomial-invertible. In other words, while it incorporates multiplication by an invertible scalar, a relation such as $b \sim_{\mathrm{E}} b_1 + b_2$ for distinct monomials $b$, $b_1$ and $b_2$ is not a valid modulo rule. Extending our work to this more general setting is an important open problem.

### ii.2.4 Non-degeneracy of graded $\mathfrak{gl}_2$-foams

In chapter 7 we prove the non-degeneracy of the graded-2-category of graded $\mathfrak{gl}_2$-foams $\mathbf{GFoam}_d$ (Corollary 1.6.4). This uses the full-strength of linear Gray rewriting modulo, including context-dependent termination and rewriting with a non-coherent modulo. The proof illustrates general techniques to study tamed congruence, such as characterization of branchings up to tamed congruence (see e.g. Proposition 7.3.8) or independent rewriting (subsections 6.2.4, 6.4.5 and 6.5.4).

## ii.3 Perspectives

A rewriting approach can be heavy: finding a practical convergent rewriting system may require a lot of trial and error, and classifying critical branchings can be laborious. Moreover:

> Obstruction to rewriting theory: there is no guaranty that a convergent rewriting system exists, and if so, that it is sufficiently reasonable to be used in practice.

However, once established the rewriting perspective provides a rich understanding of the combinatorial structure of the presentation. We list some of these (sometimes conjectural) advantages below, with special emphasizes on higher algebra:

(i) *independence on background:* rewriting theory only relies on the data of the presentation. One needs not know whether an action exists, or any far-reaching connections with other fields, specific to the presented higher algebra at hand. In particular, rewriting theory provides a self-contained approach to non-degeneracy conjectures in categorification (see subsection ii.1.7).

(ii) *computation heavy, but trouble-free:* once the rewriting foundations are well-established, applying it to a given example is a matter of choosing the right rewriting modulo data. Then, the proof of confluence is reduced to (lengthy) computations.

(iii) *computer implementation:* by design, rewriting theory provides an algorithm, which in principle can be implemented on a computer. In turn, such an implementation can be used to perform large computations. Computer implementation, and in particular machine learning, may also be used to find (efficient) convergent presentations (see e.g. [150]).

   In the context of higher algebra, this remains a largely unexplored direction of research, although the graphical proof-assistants *Globular* and *homotopy.io* [11, 12, 51, 162] constitute important first steps.

(iv) *coherence:* by considering how rewriting sequences forms unoriented cycles, one can extract an understanding of relations between relations, or *coherence data*; see section 6.1. In the linear setting, this coherence data is known as *syzygies*. Typically, critical branchings of a convergent rewriting system lead to an explicit description of this coherence data;

this has for instance been used to compute homological invariants and study Koszulness [86].

In the recent years, stable homotopy theory and $\infty$-categories have become more and more prevalent in link homologies and higher representation theory [52, 91, 119, 121, 123, 124, 129, 168], inspired by the construction of four-dimensional TQFTs and related constructions in Floer-like homologies [46, 120]. At present however, the complexity of the $\infty$-setting remains a major obstacle to exploration beyond the simplest cases. Explicit presentations of the coherence data, obtained via rewriting methods, could be an important step forward.[2]

This is a wide open direction of research.

(v) *deformation theory:* as a byproduct of coherence, rewriting theory can be used to study deformations of linear structures: we discuss it in subsection 6.6.3 and section 7.4. To our knowledge, this perspective has not explicitly appeared in the literature; see however the very recent work of Barbier on deformations of Brauer categories [15], where rewriting theory is implicit.

As we discussed in subsections i.1.3 and i.1.5, diagrammatics is the main technical tool that allowed the development of 2-dimensional algebra, especially in the context of categorification and higher representation theory. Arguably, the diagrammatics of 3-dimensional algebra is given by paths (or movies) of diagrammatics; that is, rewriting sequences. This leads us to formulate the following bold claim:

> Bold claim. *Rewriting is the new higher diagrammatics; that is, diagrammatical rewriting theory will be as significant to 3-dimensional algebra as diagrammatics is to 2-dimensional algebra.*

At present, there is not enough evidence to validate that claim. Indeed, only a couple of higher linear algebras have been considered using rewriting theory. As more and more examples are studied, we hope that the theory of linear Gray rewriting modulo will evolve into a standard set of tools, both easy-to-master and powerful, fostering the exploration of still finer and richer higher symmetries.

[2] See Remark 7.4.1 for an interesting hint in that direction.

# 1

# Graded $\mathfrak{gl}_2$-foams

This chapter defines the graded-2-category $\mathbf{GFoam}_d$ of graded $\mathfrak{gl}_2$-foams that categorifies the category $\mathbf{Web}_d$ of $\mathfrak{gl}_2$-webs. This can be viewed as a graded analogue of $\mathfrak{gl}_2$-foams *à la* Blanchet [19]. See also [17, 63, 64] for further studies of $\mathfrak{gl}_2$-foams. More precisely, the graded-2-category $\mathbf{GFoam}_d$ gives a graded analogue of the category of *directed* $\mathfrak{gl}_2$-foams as defined in [156] (see Remark 1.3.7).

Section 1.1 defines the notion of graded-2-categories. The category $\mathbf{Web}_d$ and the graded-2-category $\mathbf{GFoam}_d$ are then respectively defined in sections 1.2 and 1.3. In sections 1.4 and 1.5, we define two diagrammatics for graded $\mathfrak{gl}_2$-foams, using strings and shadings respectively. Section 1.6 describes a candidate basis for graded $\mathfrak{gl}_2$-foams; the fact that this is indeed a basis is the subject of chapter 7. Finally, based on this result we show in section 1.7 that $\mathbf{GFoam}_d$ indeed categorifies $\mathbf{Web}_d$.

## 1.1 Graded-2-categories

We define the notion of graded-2-categories to allow gradings with generic abelian groups. Taking this abelian group to be $\mathbb{Z}/2\mathbb{Z}$ recovers the notion of super-2-categories [25]. Graded-2-categories appeared in [23][1]. A similar definition also appeared in [155] in the context of algebras, under the name of

[1]We thank James MacPherson for pointing out this reference to us.

*chronological algebras*. See also [24] for an in-depth study of superstructures, including non-strict versions. Our exposition is a direct generalization to the graded case of the exposition given in [25] for super-2-categories (which they call *2-supercategories*).

Throughout the section we fix an abelian group $G$ and a unital commutative ring $\Bbbk$. We denote $\Bbbk^\times$ the abelian group of invertible elements in $\Bbbk$ equipped with the multiplicative structure. We also fix a $\Bbbk^\times$-valued bilinear form $\mu$ on $G$, that is a $\mathbb{Z}$-bilinear map

$$\mu\colon G \times G \to \Bbbk^\times.$$

Recall that a $\Bbbk$-module $V$ is said to be *$G$-graded* if it is equipped with the data of a direct sum decomposition $V = \bigoplus_{g\in G} V_g$. If $v \in V_g$ for some $g \in G$, the vector $v$ is said to be *homogeneous*; if furthermore $v \neq 0$, *it has degree* $g$, which we write as $\deg(v) = g$. Note that the zero vector is homogeneous, but does not have a well-defined degree: however, it is sometimes convenient to set $\deg(0) = 0$.

The hom-space $\mathrm{Hom}_\Bbbk(V, W)$ between two $G$-graded $\Bbbk$-modules inherits a structure of $G$-graded $\Bbbk$-module, stating that a non-zero $\Bbbk$-linear map

$$f\colon \bigoplus_{g\in G} V_g \to \bigoplus_{g\in G} W_g$$

is of degree $\deg f = h$ if $f(V_g) \subset W_{g+h}$ for all $g \in G$. If $f = 0$ or if $\deg f = 0$, we say that $f$ is *degree-preserving*.

Denote by $\Bbbk\text{-}\mathrm{Mod}_G$ the category of $G$-graded $\Bbbk$-modules, and by $\Bbbk\text{-}\mathrm{Mod}_G^0$ the wide subcategory of $\Bbbk\text{-}\mathrm{Mod}_G$ consisting only of degree-preserving $\Bbbk$-linear maps. The category $\Bbbk\text{-}\mathrm{Mod}_G$ admits a standard monoidal structure, defined on objects $V = \bigoplus_{g\in G} V_g$ and $W = \bigoplus_{g\in G} W_g$ by the formula

$$(V \otimes W)_g = \bigoplus_{\substack{(g_1,g_2)\in G\times G\\ g_1+g_2=g}} V_{g_1} \otimes_\Bbbk W_{g_2},$$

and on morphisms $f$ and $g$ by the formula $(f \otimes g)(v \otimes w) = f(v) \otimes f(w)$. However, one can give an alternative definition using the data of the grading and the bilinear form $\mu$:

**Definition 1.1.1.** *The* $(G,\mu)$-graded tensor product *is defined on objects as* $V \otimes_{G,\mu} W := V \otimes W$ *and on morphisms as*

$$(f \otimes_{G,\mu} g)(v \otimes_{G,\mu} w) = \mu(\deg v, \deg g) f(v) \otimes_{G,\mu} g(w).$$

Equipped with this graded tensor product, $\Bbbk\text{-Mod}_G$ is *not* in general a monoidal category. Indeed, morphisms respect the following *graded interchange law*:

$$(f \otimes_{G,\mu} g) \circ (h \otimes_{G,\mu} k) = \mu(\deg g, \deg h)(f \circ h) \otimes_{G,\mu} (g \circ k). \tag{1.1}$$

We now define the proper categorical structures that encompass this behaviour.

**Definition 1.1.2.** *A* $G$-graded $\Bbbk$-linear category *is a* $(\Bbbk\text{-Mod}^0_G, \otimes)$*-enriched category. A* $G$-graded $\Bbbk$-linear functor *is a* $(\Bbbk\text{-Mod}^0_G, \otimes)$*-enriched functor.*

In other words, A $G$-graded $\Bbbk$-linear category is a category such that each Hom is a $G$-graded $\Bbbk$-module, and such that composition is $\Bbbk$-bilinear and preserves the grading in the sense that $\deg(f \circ g) = \deg f + \deg g$. A $G$-graded $\Bbbk$-linear functor is a functor between two $G$-graded $\Bbbk$-linear categories that restricts to a degree-preserving $G$-graded $\Bbbk$-linear map on Homspaces.

Denote by $\Bbbk\text{-}\mathcal{C}at_G$ the category of small $G$-graded $\Bbbk$-linear categories and $G$-graded $\Bbbk$-linear functors. For $\mathcal{A}$ and $\mathcal{B}$ two $G$-graded $\Bbbk$-linear categories, their $(G,\mu)$*-graded-cartesian tensor product* is the $G$-graded $\Bbbk$-linear category $\mathcal{A} \boxtimes_{G,\mu} \mathcal{B}$ such that $\mathrm{ob}(\mathcal{A} \boxtimes_{G,\mu} \mathcal{B}) := \mathrm{ob}(\mathcal{A}) \times \mathrm{ob}(\mathcal{B})$ and

$$\mathrm{Hom}_{\mathcal{A}\boxtimes_{G,\mu}\mathcal{B}}((a,b),(c,d)) := \mathrm{Hom}_{\mathcal{A}}(a,c) \otimes_{G,\mu} \mathrm{Hom}_{\mathcal{B}}(b,d).$$

Most importantly, the composition in $\mathcal{A} \boxtimes_{G,\mu} \mathcal{B}$ is given by the graded interchange relation (1.1). This makes $\mathcal{A} \boxtimes_{G,\mu} \mathcal{B}$ a $G$-graded $\Bbbk$-linear category. Moreover, one can check that $\boxtimes_{G,\mu}$ is associative and unital, giving $\Bbbk\text{-}\mathcal{C}at_G$ the structure of a monoidal category.

**Definition 1.1.3.** *A* $(G,\mu)$-graded-2-category *is a* $\Bbbk\text{-}\mathcal{C}at_G$*-enriched category. A* $(G,\mu)$-graded-2-functor *is a* $\Bbbk\text{-}\mathcal{C}at_G$*-enriched functor.*

In particular:

**Definition 1.1.4.** *A* $(G,\mu)$-graded-monoidal category *is a* $G$*-graded* $\Bbbk$*-linear category* $\mathcal{C}$ *together with the data of a unit object and a* $G$*-graded* $\Bbbk$*-linear functor* $\otimes_{G,\mu} \colon \mathcal{C} \boxtimes_{G,\mu} \mathcal{C} \to \mathcal{C}$ *satisfying the same unital and associativity axioms as for a strict monoidal category.*

For instance, the $(G,\mu)$-graded tensor product from Definition 1.1.1 assembles into a $G$-graded $\Bbbk$-linear functor

$$\otimes_{G,\mu} \colon \Bbbk\text{-Mod}_G \boxtimes_{G,\mu} \Bbbk\text{-Mod}_G \to \Bbbk\text{-Mod}_G,$$

making $\Bbbk\text{-Mod}_G$ a $(G,\mu)$-graded-monoidal category.

Any 2-category that is both $\Bbbk$-linear and $G$-graded as a $\Bbbk$-linear 2-category is a $(G,\mu)$-graded-2-category, with $\mu$ the trivial map. In general though, the converse is not true, as the interchange law (1.1) of a generic $(G,\mu)$-graded-2-category is graded.

*Remark* 1.1.5. Graded-2-categories are closely related to Gray categories, certain 3-dimensional categorical structures where the interchange law for 2-morphisms holds weakly. We use this point of view in part II to develop a rewriting theory suitable for graded-2-categories.

*Remark* 1.1.6. If $\mu$ is symmetric, the monoidal category $(\Bbbk\text{-Mod}_G^0, \otimes)$ has as symmetric structure given by $v \otimes w \mapsto \mu(\deg(v), \deg(w)) w \otimes v$. In general, if $\mathcal{V}$ is a symmetric monoidal category, then $\mathcal{V}$-$\mathcal{C}at$, the category of small $\mathcal{V}$-enriched categories and $\mathcal{V}$-enriched functors, is itself symmetric monoidal (see e.g. [21, Proposition 6.2.9]). This allows an inductive definition of $\mathcal{V}$-enriched $n$-categories. In that general framework, $(G,\mu)$-graded-2-categories are precisely $(\Bbbk\text{-Mod}_G^0)$-enriched 2-categories.

### 1.1.1 Diagrammatics

Let $\mathcal{C}$ be a $(G,\mu)$-graded-2-category. Unpacking Definition 1.1.3, $\mathcal{C}$ consists of objects $\mathcal{C}_0$ together with a $G$-graded $\Bbbk$-linear category $\mathcal{C}(a,b)$ for each pair of objects $(a,b)$. We denote $\mathcal{C}_1(a,b)$ its objects and for $f,g \in \mathcal{C}_1(a,b)$, we denote $\mathcal{C}_2(f,g)$ the Hom-space of morphisms from $f$ to $g$. Elements of $\mathcal{C}_1(a,b)$ and $\mathcal{C}_2(f,g)$ are respectively called *1-morphisms* and *2-morphisms*. A 2-morphism $\alpha \in \mathcal{C}_2(f,g)$ can be pictured using *string diagrams*, akin to the string diagrammatics of 2-categories (see subsection i.1.3):

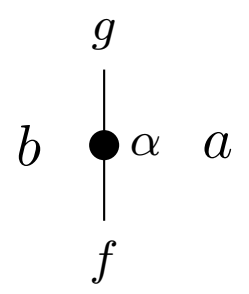

Note that we read from right to left and from bottom to top. 2-morphisms in $\mathcal{C}$ also comes equipped with two compositions, the *horizontal composition* $\star_0$ and the *vertical composition* $\star_1$. The vertical composition denotes all the compositions in the $G$-graded $\Bbbk$-linear categories $\mathcal{C}(a,b)$. It is pictured by stacking the 2-morphisms atop each other:

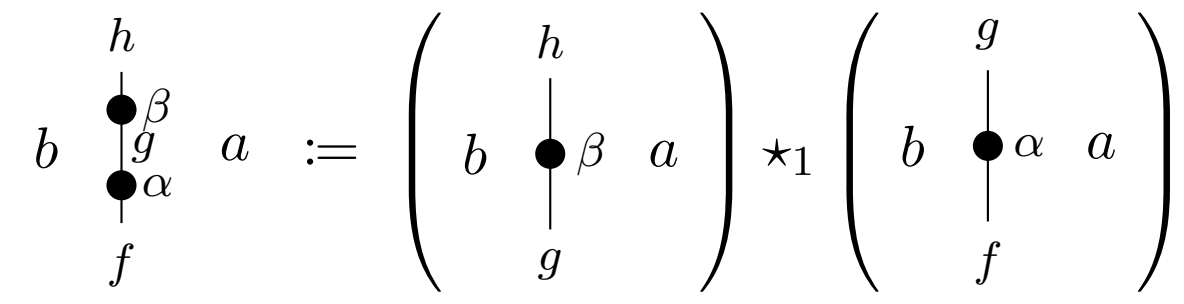

On the other hand, the horizontal composition is pictured by putting the two 2-morphisms side-by-side:

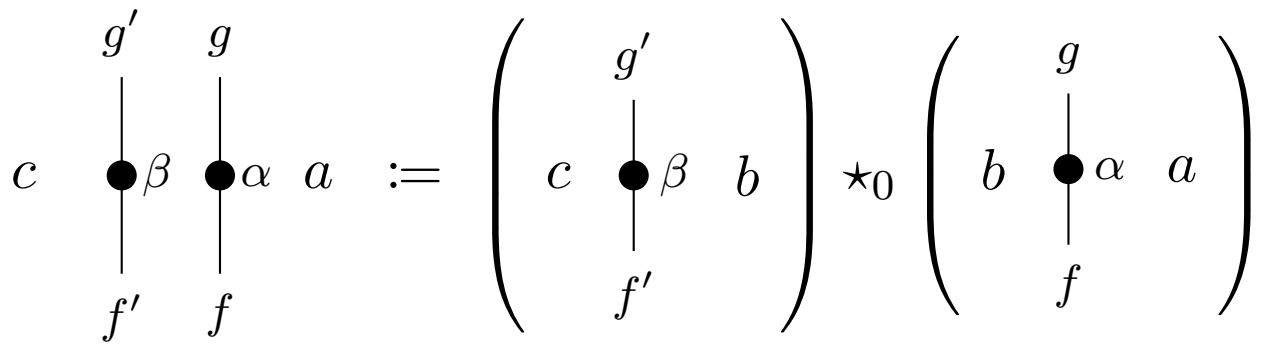

To avoid clutter, we leave regions and lines unlabelled for now. The rule "compose first horizontally, then vertically" resolves the ambiguity for the order of compositions:

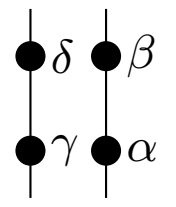

$$:= (\delta \star_0 \beta) \star_1 (\gamma \star_0 \alpha)$$

Contrary to 2-categories, it is in general *not* the same as "compose first vertically, then horizontally". Indeed, in a graded-2-category sliding a 2-morphism past another vertically comes at the price of an invertible scalar. This is the graded interchange law:

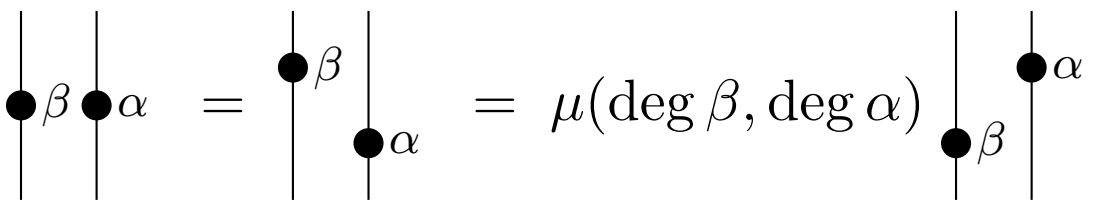

(1.2)

In particular, one must be careful with the relative vertical positions of 2-morphisms. To avoid confusion, in this thesis we always draw string diagrams such that no two 2-morphisms lie on the same vertical level. One can make an exception for 2-morphisms with trivial grading, as those can slide vertically without adding scalars. Finally, as customary already in the string diagrammatics of 2-categories, we usually do not picture identities.

**Definition 1.1.7.** *We say that $\mu$ is* symmetric *if $\mu(g,h)\mu(h,g) = 1$ for all $g, h \in G$.*

Note that $\mu$ is automatically symmetric in the super case $G = \mathbb{Z}/2\mathbb{Z}$. If $\mu$ is symmetric, the scalar appearing in the graded interchange law depends *only* on the relative vertical positions of the 2-morphisms:

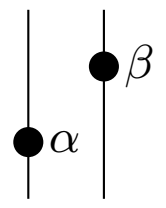

$= \mu(\deg\beta, \deg\alpha)$
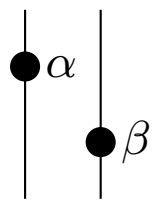

if $\mu$ is symmetric.

Compare with Eq. (1.2). In other words, passing $\beta$ *down* $\alpha$ always adds the scalar $\mu(\deg\beta, \deg\alpha)$, regardless of their respective horizontal positions.

### 1.1.2 Grothendieck ring

Let $H$ be an abelian group. Recall that the Grothendieck ring of an $H$-graded 2-category is a certain $\mathbb{Z}[H]$-linear category $K_0(\mathcal{C})|_H$, capturing the isomorphism classes of 1-morphisms in $\mathcal{C}$. More precisely, one defines

$$K_0(\mathcal{C})|_H := K_0(\underline{\mathcal{C}}_H^{\oplus}),$$

where $\mathcal{C}_H$ is the $H$-*envelope* of $\mathcal{C}$ (formally adjoining $H$-shifts), $\underline{\mathcal{C}}_H$ is the *underlying category* of $\mathcal{C}_H$ (restricting to grading-preserving 2-morphisms), and $\underline{\mathcal{C}}_H^{\oplus}$ is the *additive closure* of $\underline{\mathcal{C}}_H$ (formally adjoining direct sums):

$$\mathcal{C} \quad\rightsquigarrow\quad \mathcal{C}_H \quad\rightsquigarrow\quad \underline{\mathcal{C}}_H \quad\rightsquigarrow\quad \underline{\mathcal{C}}_H^{\oplus}.$$

One could give a similar definition for $(G,\mu)$-graded-2-categories; in that case, the Grothendieck ring would be a $\mathbb{Z}[G]$-linear category. This generalizes the notion of Grothendieck ring for super-2-categories [24]. The main difficulty lies in properly defining the horizontal composition in the $G$-envelope, ensuring that the graded interchange law holds (see Remark 1.1.8). More generally, if a $(G,\mu)$-graded-2-category is equipped with an extra $H$-grading, one can define the Grothendieck ring to be $\mathbb{Z}[G \times H]$-linear.

In this thesis, we will deal with $H$-graded $(G,\mu)$-graded-2-categories where $G = \mathbb{Z} \times \mathbb{Z}$ and $H = \mathbb{Z}$; in this context, the $H$-grading is called the *quantum grading*. In fact, the quantum grading will be derived from the $G$-grading, *but we shall consider them as distinct gradings.* Moreover, *we shall ignore the $G$-grading when taking Grothendieck rings.* To emphasize this point, we write the Grothendieck ring $K_0(\mathcal{C})|_q$, where "$q$" refers to the quantum grading. We call it the *quantum Grothendieck ring*. In particular, the *$q$-envelope* $\mathcal{C}_q$ is a certain $(G, \mu)$-graded-2-category defined by adjoining $q$-shifts; the latter do *not* change the $G$-grading, and hence do not interact with the $(G, \mu)$-graded-2-categorical structure.

We give below the formal definitions, starting with the 1-dimensional case; see also [69, chap. 11] for a review of the basic notions.

*Remark* 1.1.8. For completeness and future reference, we give below the definition of the horizontal composition in the $G$-envelope of a $(G, \mu)$-graded-2-category, following the conventions of [24]:

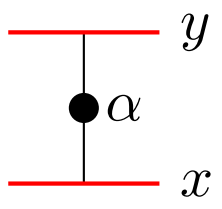

$\star_0$

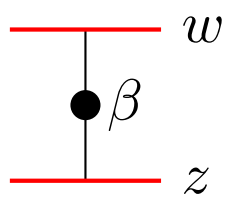

$$= \mu(-x, \deg_G \beta - z + w)\mu(\deg_G \alpha, w) \tag{1.3}$$

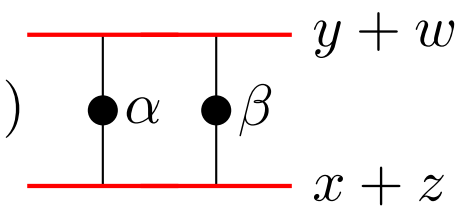

*1-dimensional case: $H$-graded categories*

Let $H$ be an abelian group. An *$H$-graded category* is a category enriched over $H$-graded abelian groups. If $H = \{*\}$ is the trivial group, the category is said to be *pre-additive*. A pre-additive category $C$ is a *category with $H$-shifts* if it is equipped with a group morphism $H \to \mathcal{A}ut(C)$ where $\mathcal{A}ut(C)$ is the group of invertible endofunctors; that is, $C$ is equipped with an action of $H$ by autofunctors. We denote $f\{x\}$ the action of $x \in H$ on an object $f \in \mathrm{ob}(C)$. When $H = \mathbb{Z}$, we also write $q^n f := f\{n\}$.

Let $C$ be a pre-additive category. The category $C$ naturally embeds in an additive category $C^{\oplus}$, its *additive closure*. Objects in $C^{\oplus}$ are formal direct sums of objects in $C$, and morphisms are matrices whose entries are morphisms in $C$. If $C$ is a category with $H$-shifts, then so is $C^{\oplus}$ by extending the $H$-action additively.

Let $C$ be an additive category. The *(split) Grothendieck ring of $C$* is the abelian group $K_0(C)$ generated by elements $[f]$ for each object $f \in \mathrm{ob}(C)$, and subject to the relations $[f \oplus g] = [f] + [g]$; in particular, if $f \cong g$ then $[f] = [g]$. If $C$ is with $H$-shifts, then $K_0(C)$ has the structure of a $\mathbb{Z}[H]$-module, where the $H$-action is given by $x \cdot_H [f] := [f\{x\}]$.

Let $C$ be a pre-additive category. The *Grothendieck ring of $C$* is the Grothendieck ring of its additive closure: $K_0(C) := K_0(C^{\oplus})$.

Let $C$ be an $H$-graded category. The *$H$-envelope of $C$* is the category $C_H$ whose objects are formal $H$-shifts $f\{x\}$, where $x \in H$ and $f \in \mathrm{ob}(C)$, and there is a morphism

$$\alpha\colon f\{x\} \to g\{y\}$$

in $C_H$ for each morphism $\alpha\colon f \to g$ in $C$ and each pair of elements $x, y \in H$. If $\deg_H \alpha$ denotes the $H$-degree of $\alpha$, we set $\deg_H(\alpha\colon f\{x\} \to g\{y\}) = \deg_H \alpha - x + y$. The $H$-envelope is both $H$-graded and with $H$-shifts.

Let $C$ be an $H$-graded category. The *underlying category of $C$* is the category $\underline{C}$ consisting of degree-preserving morphisms. The *$H$-shifted closure of $C$*, denoted by $\underline{C}_H$, is the underlying category of its $H$-envelope; note that it has $H$-shifts. Finally, the *$H$-Grothendieck ring of $C$* is the Grothendieck ring of its $H$-shifted closure $\underline{C}_H$:

$$K_0(C)|_H \coloneqq K_0(\underline{C}_H) = K_0((\underline{C}_H)^{\oplus}).$$

The $H$-Grothendieck ring is a $\mathbb{Z}[H]$-module.

All of the above extend to the case where $C$ is an $(G \times H)$-graded category, ignoring the $G$-grading throughout. Note that the $H$-envelope is canonically $G$-graded, setting

$$\deg_G(\alpha\colon f\{x\} \to g\{y\}) = \deg_G \alpha.$$

*2-dimensional case: $(G, \mu)$-graded-2-categories*

Let $\mathcal{C}$ be a $(G, \mu)$-graded-2-category and $H$ an abelian group. We say that $\mathcal{C}$ is *additive* if each hom-category $\mathcal{C}(a, b)$ is additive and horizontal composition is bilinear. We say that $\mathcal{C}$ is *with $H$-shifts* if it is equipped with an action of $H$ by automorphisms which is the identity on objects. In particular, for each pair of objects $a, b \in \mathrm{ob}(\mathcal{C})$ the $G$-graded category $\mathcal{C}(a, b)$ is with $H$-shifts. Its Grothendieck ring is the $\mathbb{Z}[H]$-linear category obtained by taking the Grothendieck ring of each hom-categories $\mathcal{C}(a, b)$, and inducing a bilinear composition by setting $[f] \circ [g] \coloneqq [f \star_0 g]$.

Let $\mathcal{C}$ be an $H$-graded $(G, \mu)$-graded-2-category. The *$H$-envelope of $\mathcal{C}$* is the $(G, \mu)$-graded-2-category with $H$-shifts defined by taking the $H$-envelope of each hom-category $\mathcal{C}(a, b)$. The *underlying linear 2-category of $\mathcal{C}$* is the sub-2-category $\underline{\mathcal{C}}$ obtained by taking the underlying category of each hom-category $\mathcal{C}(a,b)$. Notions of $H$-shifted closure and additive closure are defined analogously to the 1-dimensional case. Finally:

**Definition 1.1.9.** *Let $\mathcal{C}$ be an $H$-graded $(G,\mu)$-graded-2-category. Its $H$-*Grothendieck ring *is the Grothendieck ring of the additive closure of its $H$-shifted closure:*

$$K_0(\mathcal{C})|_H := K_0((\underline{\mathcal{C}}_H)^{\oplus}).$$

*The $H$-Grothendieck ring is a $\mathbb{Z}[H]$-linear category.*

## 1.2 $\mathfrak{gl}_2$-webs

In our context, a $\mathfrak{gl}_2$-web is a certain trivalent graph, smoothly embedded in $\mathbb{R}\times[0,1]$, that we view as a morphism in a certain category defined below. Objects, called *weights*, are elements of the following set:

$$\underline{\Lambda}_d := \bigsqcup_{k\in\mathbb{N}} \{\lambda\in\{1,2\}^k \mid \lambda_1+\ldots+\lambda_k = d\}.$$

For each $\lambda\in\underline{\Lambda}_d$ with $k$ coordinates, we can define a label on its coordinates

$$l_\lambda\colon \{1,\ldots,k\}\to\{1,\ldots,d\}$$

by setting $l_\lambda(i) = \sum_{j<i}\lambda_i + 1$. For instance, $l_{(1,1,2,1)} = (1,2,3,5)$. Foreseeing the string diagrammatics, we call this label the *colour* of the coordinate. The identity web of a weight is pictured as a juxtaposition of straight vertical lines in $\mathbb{R}\times[0,1]$, decorated as *single* (black) or *double* (orange) lines:

$\mathrm{id}_{(1,1,2,1)} =$

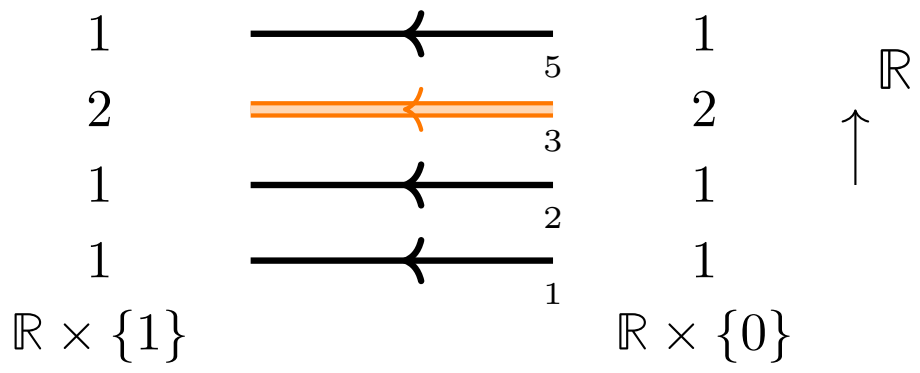

Note that we read webs from right to left. Note also that each line has a colour, given by the colour of the corresponding coordinate (this colour is unrelated to whether the line is pictured "black" or "orange"). Effectively, the colour of a line counts the number of preceding lines (plus one), counting twice double lines. A generic $\mathfrak{gl}_2$-*web* is either an identity web or a composition of the following generating webs:

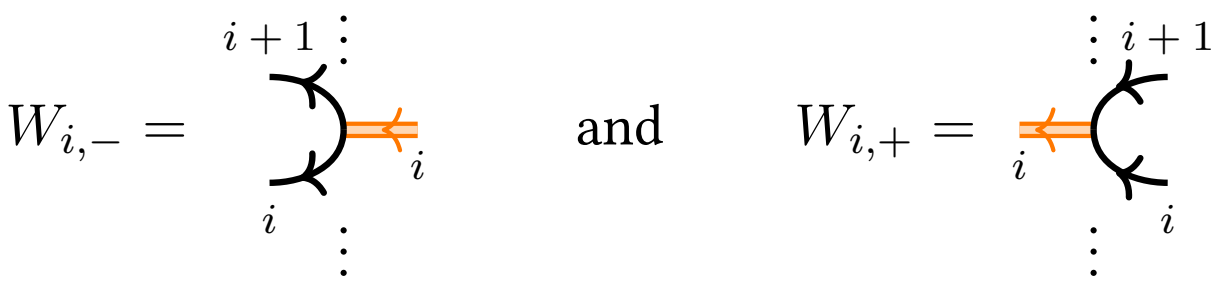

Hereabove the dots denote possibly additional vertical single or double lines, and composition is given by stacking webs atop each other (reading from right to left). Note that a web $W$ has an underlying unoriented flat tangle diagram, denoted $\mathfrak{sl}(W)$, given by forgetting the double lines and the orientations.

*Remark* 1.2.1 (Orientation). Our webs are given an orientation "flowing" from $\mathbb{R} \times \{0\}$ to $\mathbb{R} \times \{1\}$, which is more restrictive than some definitions in the literature. Such webs are sometimes called *acyclic* or *left-directed* (e.g. in [156]). As this orientation is canonical, we often omit it.

**Definition 1.2.2.** *The category* $\mathbf{Web}_d$ *has objects* $\underline{\Lambda}_d$*, and morphisms are* $\mathbb{Z}[q, q^{-1}]$*-linear combinations of* $\mathfrak{gl}_2$*-webs, up to the following* web relations*:*

$$W_{i,s_1} W_{j,s_2} = W_{j,s_2} W_{i,s_1}$$
$$(\textit{for all } s_1, s_2 \in \{-,+\} \textit{ and } |i-j| > 1)$$

*interchange relations*

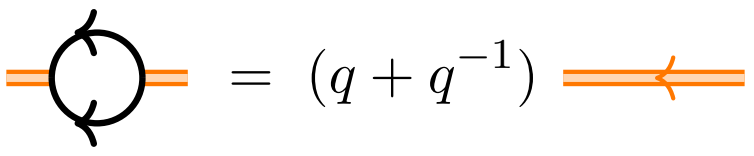

*circle evaluation*

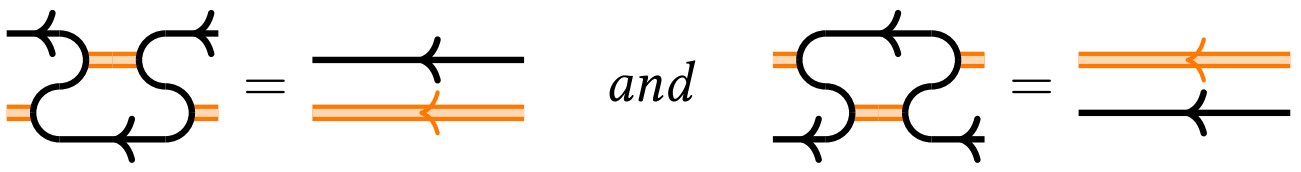

*isotopies of flat tangle diagrams*

We shall need the following notion:

**Definition 1.2.3.** *A* spatial isotopy[2] *of flat tangle diagrams is a usual isotopy with the addition of the following* spatial move*:*

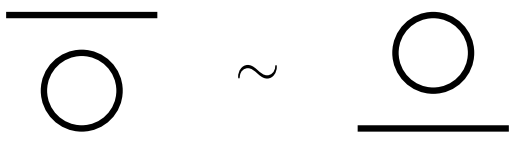

The relations in $\mathbf{Web}_d$ fully capture the spatial isotopy classes of the underlying flat tangle diagrams, in the sense of the following lemma.

[2]The terminology is taken from [173].

**Lemma 1.2.4.** *Let $W$ and $W'$ be two $\mathfrak{gl}_2$-webs with the same domain and codomain. Then $W$ and $W'$ are equal in* $\mathbf{Web}_d$ *if and only if there exists a spatial isotopy between their underlying flat tangle diagrams $\mathfrak{sl}(W)$ and $\mathfrak{sl}(W')$.*

Using relations in Definition 1.2.2, any closed strand in $\mathfrak{sl}(W)$ evaluates to $q+q^{-1}$ in $W$. Moreover, if $\mathfrak{sl}(W)$ does not have any closed strand, then the last two relations in Definition 1.2.2 are enough to capture all isotopies of $\mathfrak{sl}(W)$. These two facts essentially constitute the proof of Lemma 1.2.4; a formal proof can be found in subsection 6.6.3 using rewriting theory.

## 1.3 Graded $\mathfrak{gl}_2$-foams

Foams provide a suitable notion of cobordisms between webs. They are certain singular surfaces locally modelled on the product of the interval $[0,1]$ with the letter "Y" (see Fig. 1.1), embedded in $(\mathbb{R}\times[0,1])\times[0,1]$. In Fig. 1.1, $\mathbb{R}$ is pictured from front to back, while the last interval is pictured from bottom to top. In this context, we refer to the singular curves as *seams*, and call *facets* the components of the complement of the set of seams. Facets have a thickness, either single or double, and we refer to them respectively as *1-facets* (or *single facets*, shaded blue) and *2-facets* (or *double facets*, shaded orange). We refer to *cross-sections* as cross-sections of the projection $\pi\colon(\mathbb{R}\times[0,1])\times[0,1]\to[0,1]$ onto the last coordinate (pictured vertically).

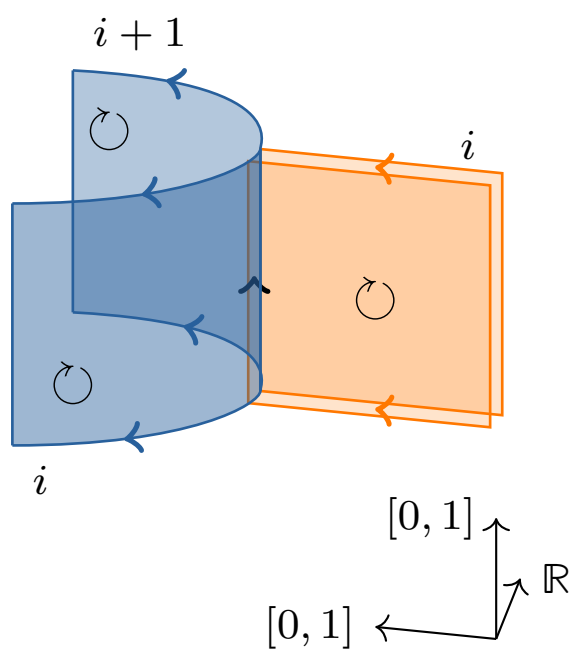

**Fig. 1.1** Local model for foams.

**Definition 1.3.1.** *A* graded $\mathfrak{gl}_2$-foam, *or simply* foam, *is a topological singular surface embedded in $(\mathbb{R}\times[0,1])\times[0,1]$, such that generic cross-sections are webs and all non-generic cross-sections respect one the following local behaviours:*

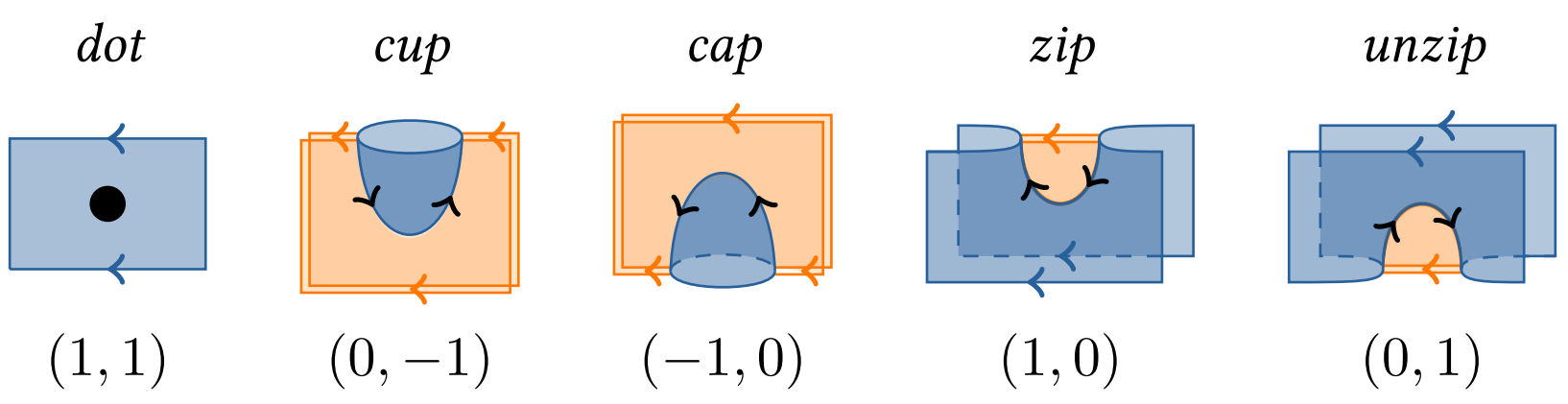

*where the axes are oriented as in Fig. 1.1. Each such local behaviour $F$ is endowed with a $\mathbb{Z}^2$-degree $\deg_{\mathbb{Z}^2}(F)$, denoted below the associated picture. Extending additively, this induces a $\mathbb{Z}^2$-grading on graded $\mathfrak{gl}_2$-foams.*

We also call the above local behaviours *generating foams*, or simply *generators*. A graded $\mathfrak{gl}_2$-foam $F \subset (\mathbb{R} \times [0,1]) \times [0,1] \to [0,1]$ is viewed as a morphism from the web $F \cap \pi^{-1}(\{0\})$ to the web $F \cap \pi^{-1}(\{1\})$, reading from bottom to top. The *dot* in the first picture of Definition 1.3.1 is a formal decoration that any 1-facet may carry. By removing the 2-facets, a foam $F$ has an underlying surface, denoted $\mathfrak{sl}(F)$. We assume that $\mathfrak{sl}(F)$ is smooth and that the vertical projection $\pi$ defines a separative Morse function for $\mathfrak{sl}(F)$, considering dots as critical points. Note that the local behaviours in Definition 1.3.1 dictate the local behaviours around critical points of $\pi$. Because $\pi$ is separative, each such local behaviour lies on a distinct vertical position.

The facets of a foam $F$ admit a canonical orientation, induced by the canonical orientation on the web $F \cap \pi^{-1}(\{0\})$. That is, we endow the facets incident to $F \cap \pi^{-1}(\{0\})$ with an orientation compatible with the canonical orientation of $F \cap \pi^{-1}(\{0\})$, and extend globally with the condition that at a given seam, the orientation induced by the 2-facet is opposite to the orientation induced by the two 1-facets (see Fig. 1.1). In particular, the facets incident to $F \cap \pi^{-1}(\{1\})$ induce an orientation on $F \cap \pi^{-1}(\{1\})$ opposite to its canonical orientation. Similarly, facets are endowed with a colour induced by the colour on webs; see also the string diagrammatics in section 1.4. This also defines an orientation and colour on each seam, induced from the orientation and colour of its incident 2-facet.

As in the non-graded case, we will consider graded $\mathfrak{gl}_2$-foams up to isotopies. In our context, an *isotopy* is a boundary-preserving isotopy which generically preserves the cross-section condition in Definition 1.3.3. Sliding a dot along its 1-facet is also considered as an isotopy. Moreover, we assume that the restriction of an isotopy to the underlying surface is a diffeotopy.

We distinguish different kinds of isotopies, in analogy with the property of the corresponding diffeotopy for the underlying surface:

- If the isotopy preserves the relative vertical positions of the generating foams, we say that it is a *Morse-preserving isotopy*.
- If the isotopy only interchanges the vertical positions of generating foams, we say that it is a *Morse-singular isotopy*.

- If the isotopy is such the underlying diffeotopy is a birth-death diffeotopy,[3] we say that the isotopy is a *birth-death isotopy*.

We refer to [155] for relevant details on Morse theory.

In the graded case, isotopies only hold up to invertible scalars. Some are controlled by a graded-2-categorical structure:

**Definition 1.3.2.** *Set the commutative ring*

$$\Bbbk := \mathbb{Z}[X, Y, Z^{\pm 1}]/(X^2 = Y^2 = 1),$$

*the abelian group* $G := \mathbb{Z}^2$ *and the bilinear form* $\mu$ *as*

$$\begin{aligned} \mu\colon \mathbb{Z}^2 \times \mathbb{Z}^2 &\to \Bbbk^\times, \\ ((a,b),(c,d)) &\mapsto X^{ac} Y^{bd} Z^{ad-bc}. \end{aligned}$$

Note that $\mu$ is symmetric in the sense of Definition 1.1.7.

**Definition 1.3.3.** *The* $(\mathbb{Z}^2, \mu)$*-graded-2-category* $\mathbf{GFoam}_d$ *has objects* $\underline{\Lambda}_d$*, its 1-morphisms are* $\mathfrak{gl}_2$*-webs and its 2-morphisms are* $\Bbbk$*-linear combinations of graded* $\mathfrak{gl}_2$*-foams, regarded up to the following relations:*

*(i) If* $\varphi\colon F_1 \to F_2$ *is a Morse-preserving isotopy, then* $F_1 = F_2$ *in* $\mathbf{GFoam}_d$*.*

*(ii) If* $\varphi\colon F_1 \to F_2$ *is a Morse-singular isotopy interchanging the vertical positions of exactly two critical points* $p$ *and* $q$*, with* $p$ *vertically above* $q$ *in* $F_1$*, then* $F_1 = \mu(\deg p, \deg q) F_2$ *in* $\mathbf{GFoam}_d$*.*

*(iii) All the local relations in Fig. 1.2 below.*

Setting $X = Y = Z = 1$ recovers (even) $\mathfrak{gl}_2$-foams (see Remark 1.3.7), while setting $X = Z = 1$ and $Y = -1$ gives so-called *super* $\mathfrak{gl}_2$*-foams* as discussed in the introduction. One could also set $Y = Z = 1$ and $X = -1$, leading to an essentially identical theory.

*Remark* 1.3.4. The graded-2-categories $\mathbf{GFoam}_d$ for $0 \leq d < \infty$ could be gathered together into some sort of "monoidal graded-2-category" $\mathbf{GFoam}_\otimes$. We do not pursue this direction here; see however section 5.1 where we use the analogue $\mathbf{Web}_\otimes$ for $\mathfrak{gl}_2$-webs.

[3] That is, collapsing the two singularities associated to a saddle and a cup or cap by "smoothing out" the surface, or the converse; see the zigzag relations in Fig. 1.2 for the underlying surfaces.

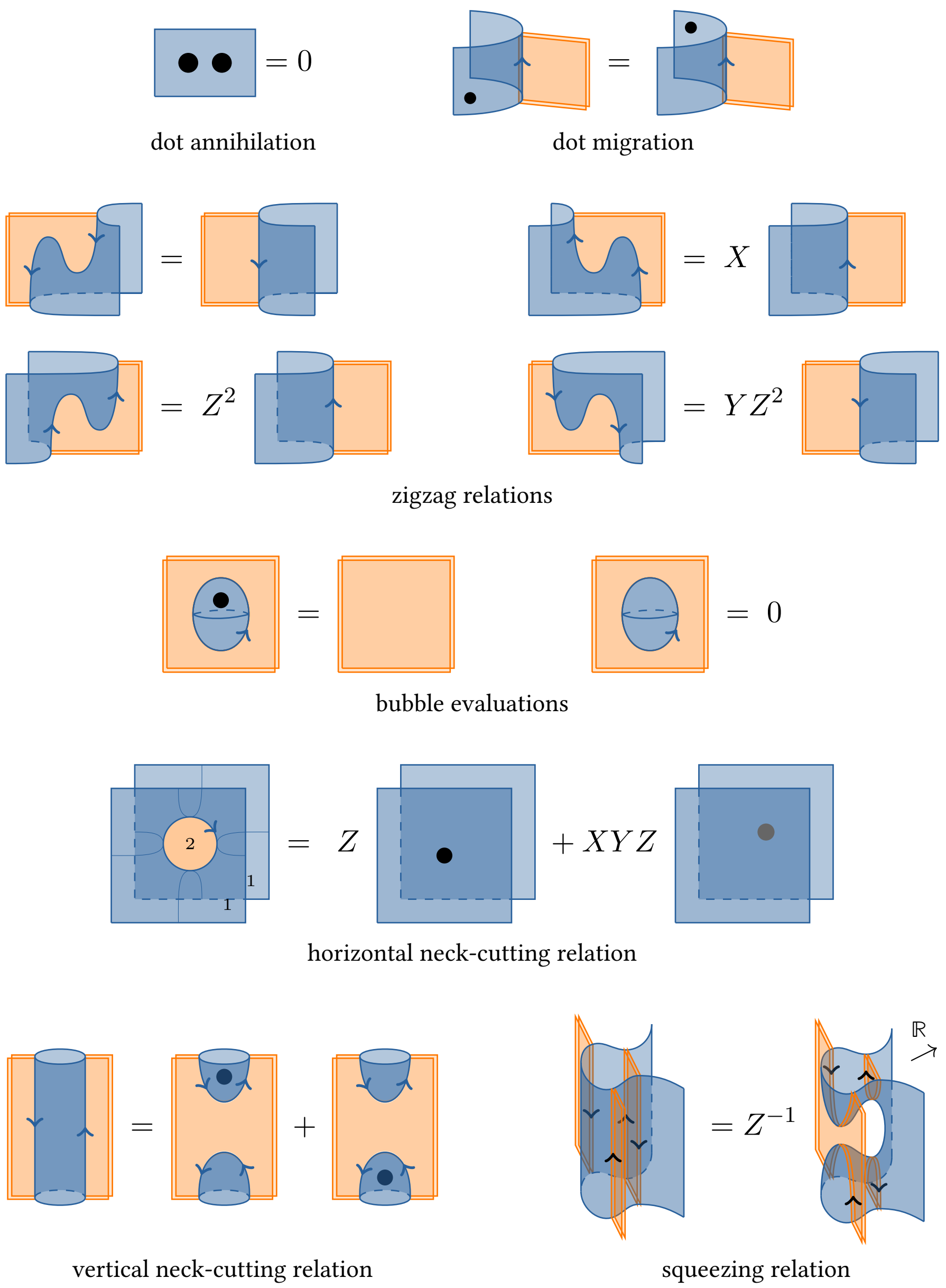

**Fig. 1.2** Relations in $\mathbf{GFoam}_d$. A grey dot denotes a dot behind a facet. The $\mathbb{R}$ axis is pictured from front to back, except for the squeezing relation for which it is pictured from left to right for better readability.

*Remark* 1.3.5. The $\mathbb{Z}^2$-grading induces a $\mathbb{Z}$-grading on graded $\mathfrak{gl}_2$-foams, the *quantum grading* (or *$q$-grading*). If $q\colon \mathbb{Z}^2 \to \mathbb{Z}$ denotes the map $q(a,b) = a+b$, we set for each foam $F$:

$$\operatorname{qdeg}(F) := q(\deg_{\mathbb{Z}^2}(F)).$$

Following subsection 1.1.2, we write $\underline{(\mathbf{GFoam}_d)}_q^{\oplus}$ the additive $q$-shifted closure of $\mathbf{GFoam}$: it allows formal shifts of webs in the quantum grading, restrict graded foams to those preserving the quantum grading, and allows formal direct sums. In this construction, we view the quantum grading as independent from the $\mathbb{Z}^2$-grading; hence $\underline{(\mathbf{GFoam})}|_q^{\oplus}$ is still a $(\mathbb{Z}^2, \mu)$-graded-2-category, and shifts in quantum degree do not affect the $\mathbb{Z}^2$-degree.

Similar to webs (Lemma 1.2.4), the following lemma shows that relations on foams capture the diffeotopy classes of the underlying surfaces; or rather, the underlying dotted surfaces, where sliding a dot along a connected component is considered to be a diffeotopy. Below we write $F \dot{\sim} F'$ whenever there exists an invertible scalar $r \in \mathbb{k}^{\times}$ such that $F = rF'$. Recall that $\mathfrak{sl}(F)$ denotes the underlying surface of a foam $F$.

**Lemma 1.3.6.** *Let $F$ and $F'$ be two foams in $\mathbf{GFoam}_d$ with the same domain and codomain. If $\mathfrak{sl}(F)$ and $\mathfrak{sl}(F')$ are isotopic, then $F \dot{\sim} F'$ in $\mathbf{GFoam}_d$.*

*Proof.* By Cerf theory ([34]; see [155, appendix A] for a review), diffeotopic surfaces are related by diffeotopies preserving the relative vertical positions of critical points, Morse-singular diffeotopies interchanging two critical points, and birth-death diffeotopies. These correspond to the relations (i), (ii) and the zigzag relations in $\mathbf{GFoam}_d$, respectively. Finally, if $\mathfrak{sl}(F)$ and $\mathfrak{sl}(F')$ are diffeotopic by sliding a dot along a connected component, then $F$ and $F'$ are related by sliding the dot along 1-facets, possibly crossings seams using dot migration. □

*Remark* 1.3.7. Recall the quantum grading from Remark 1.3.5, and $\underline{(\mathbf{GFoam}_d)}_q^{\oplus}$ the additive $q$-shifted closure of $\mathbf{GFoam}_d$. One can check that the $\mathbb{k}$-linear 2-category

$$\underline{(\mathbf{GFoam}_d)}_q^{\oplus}|_{X=Y=Z=1}$$

is precisely the category $n\mathbf{Foam}(N)^{\bullet}$ from [156, p. 1322] with $n = 2$ and $N = d$, with the same quantum grading. This follows from Lemma 1.3.6 and renormalizing the $i$-labelled dot, the $i$-labelled cap and the $i$-labelled zip by $(-1)^i$. In particular, we warn the reader familiar with the even setting that, while the dot migration relation in Fig. 1.2 has no sign, this is *not* an essential feature of our construction, but rather a choice of normalization, allowed by the restriction to directed $\mathfrak{gl}_2$-foams (see Remark 1.2.1).

## 1.4 String diagrammatics

We introduce string diagrammatics for graded $\mathfrak{gl}_2$-foams. This is based on the observation that *a foam is fully described by its seams*; more precisely, by its domain object, its seams, and their orientations and colours. For identity foams, this holds since webs are generated by $W_{i,\pm}$ and $\mathrm{id}_{W_{i,\pm}}$ is determined by its domain and the seam, together with its orientation and colour. Thus, we can represent a web with an *identity foam diagram* (or simply *identity diagram*), a horizontal juxtaposition of oriented vertical strands coloured with elements in $\{1, 2, \dots, d-1\}$. For instance:

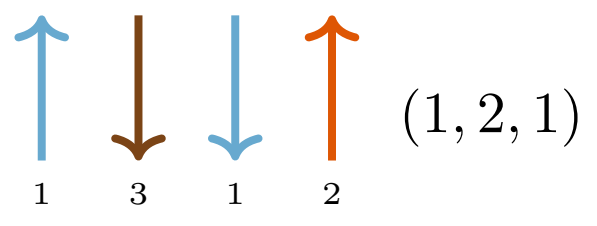

a labelled identity diagram

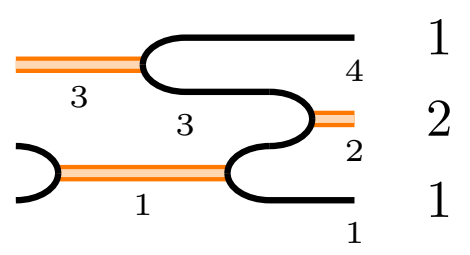

the web it represents

where we remind the reader that webs are read and oriented from left to right. In the example above, we have labelled its rightmost region to specify the domain of the web. This induces a label on all of its regions thanks to the following rule (here the value $l_\lambda$ (see section 1.2) is given below the corresponding coordinate):

$$(\dots,\underset{i}{1},\underset{i+1}{1},\dots)\ \underset{i}{\big\uparrow}\ (\dots,\underset{i}{2},\dots) \quad\text{and}\quad (\dots,\underset{i}{2},\dots)\ \underset{i}{\big\downarrow}\ (\dots,\underset{i}{1},\underset{i+1}{1},\dots) \tag{1.4}$$

A label of the regions of an identity diagram with elements in $\underline{\Lambda}_d$ is said to be *legal* if it locally satisfies (1.4). In general, if the regions of an identity diagram are labelled with a legal label, then it is only necessary to label one region to specify the label. Note that the converse does not hold: the data of a label on a region cannot necessarily be extended to a legal label on the entire diagram. An identity diagram equipped with the data of a legal label is called a *labelled* *identity (foam) diagrams.* By the above discussion, labelled identity diagrams are in one-to-one correspondence with webs.

Given these diagrammatics for webs, it is not difficult to extend it to foams. Indeed, local behaviours in Definition 1.3.1 are determined by the seam, together with its orientation and colour.

**Definition 1.4.1.** *A (generic)* foam diagram, *or simply* diagram, *is a diagram obtained by vertical and horizontal juxtapositions of identity diagrams and generators given below:*

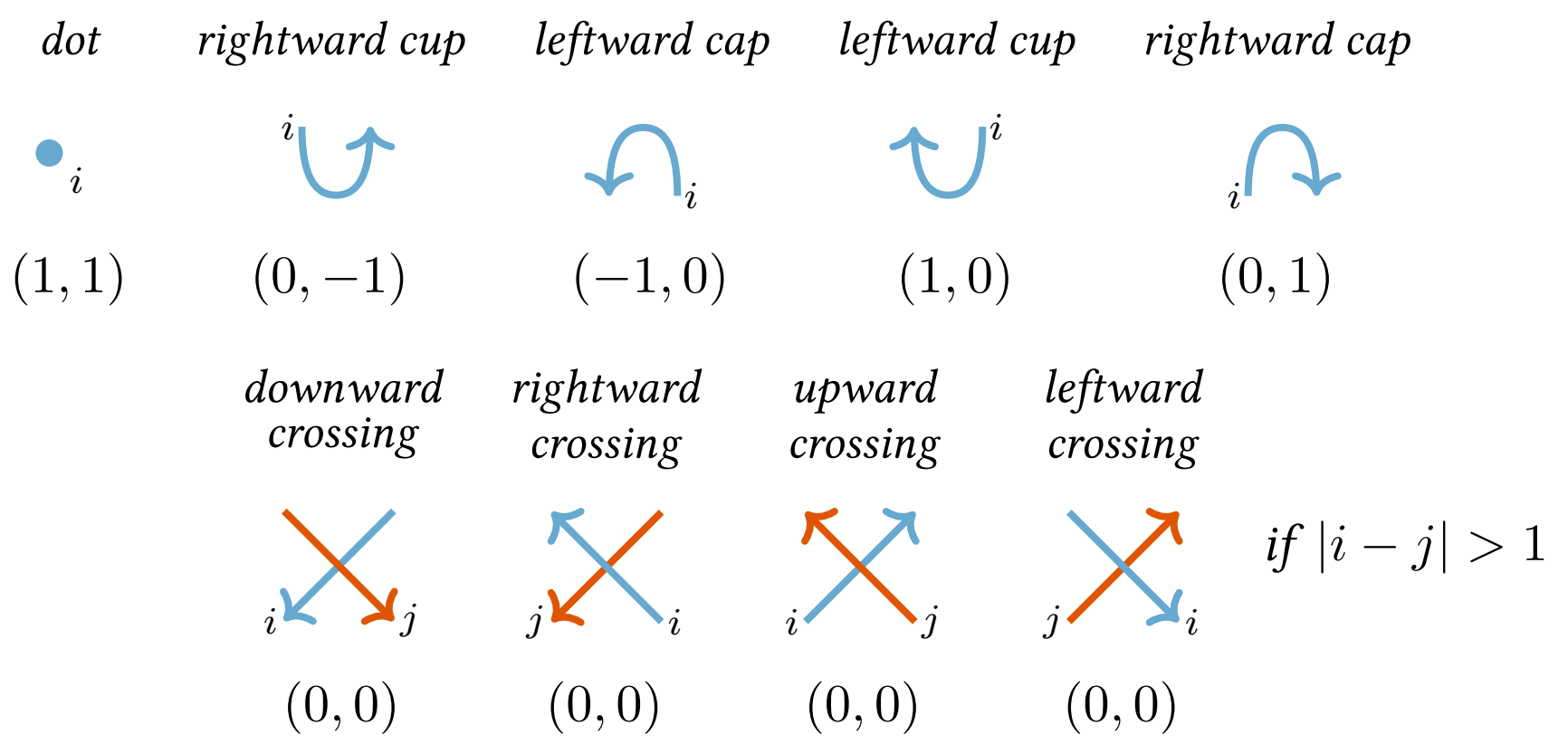

*Each generator is equipped with a $\mathbb{Z}^2$-degree, which extends additively to generic diagrams. We assume that in a generic diagram, generators are in general position with respect to the vertical direction.*

In a foam diagram, we refer to strands coloured $i$ as $i$-*strands*, and dots coloured $i$ as $i$-*dots*. As for an identity diagram, we can label the regions of a generic diagram with elements of $\underline{\Lambda}_d$. The label is said to be *legal* if it satisfies condition (1.4) as before, and if for each $i$-dot contained in a region labelled with $\lambda$, we have $\lambda_i = 1$. This latter condition corresponds to the fact that dots can only sit on 1-facets. A *labelled (foam) diagram* is a foam diagram equipped with a legal label. By the above discussion, labelled diagrams are in one-to-one correspondence with foams.

This provides a string diagrammatics for the graded-2-category $\mathbf{GFoam}_d$:

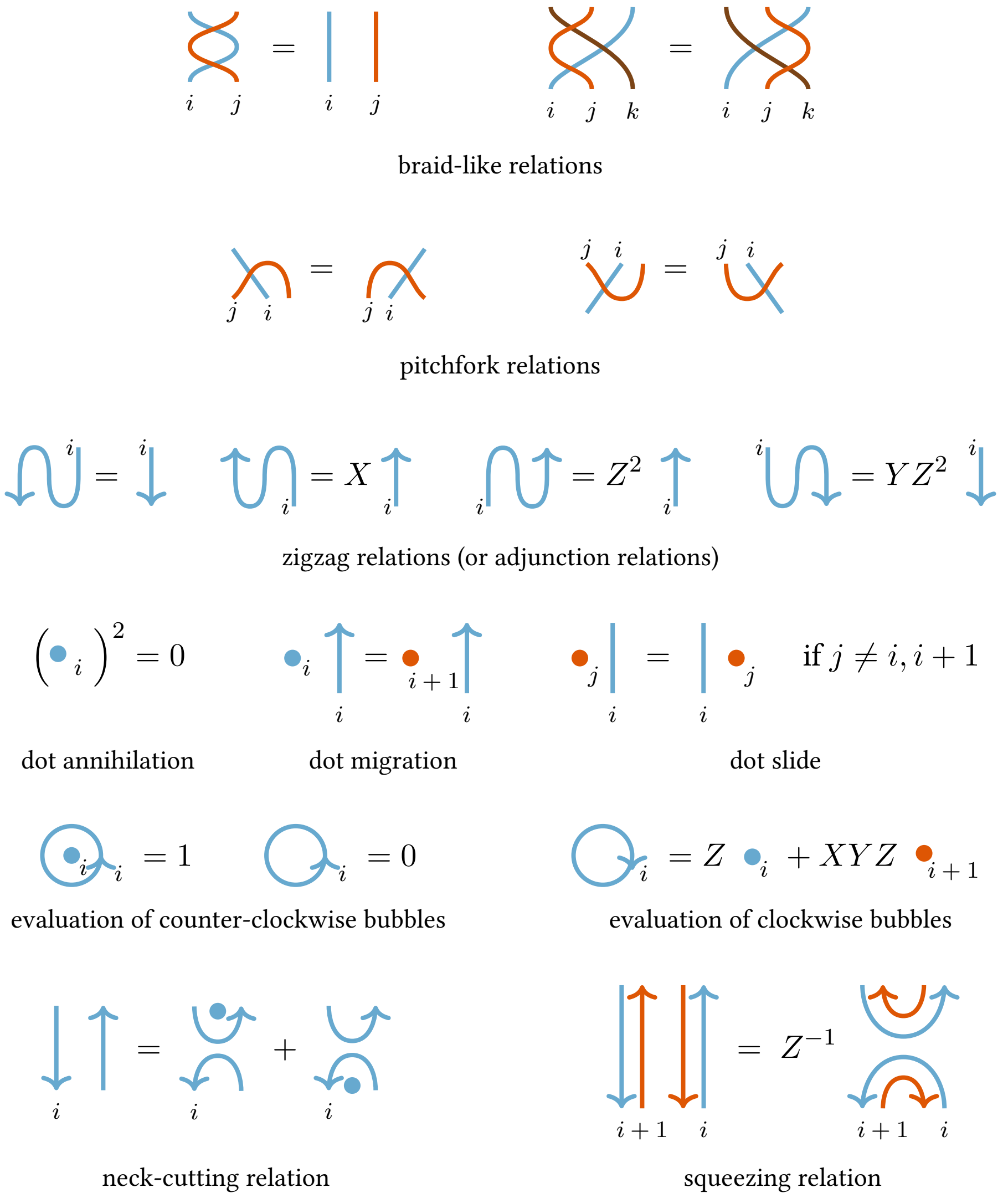

**Fig. 1.3** Relations in $\mathbf{Diag}_d^{\Lambda}$. We omitted the objects labelling the regions of each diagram: this avoids clutter and emphasizes that relations are independent on the ambient object. If no orientation is given, the relation holds for all orientations. In the case of the braid-like and pitchfork relations, colours should be so that the crossings exist.

**Definition 1.4.2.** *The $(\mathbb{Z}^2, \mu)$-graded-monoidal category* $\mathbf{Diag}_d^\Lambda$ *has objects* $\underline{\Lambda}_{n,d}$*, 1-morphisms labelled identity foam diagrams, and 2-morphisms labelled foam diagrams, regarded up to the following relations:*

*(i) If* $\varphi\colon D_1 \to D_2$ *is a planar isotopy such that* $\varphi$ *preserves the relative vertical positions of the generators, then* $D_1 = D_2$ *in* $\mathbf{Diag}_d^\Lambda$*.*

*(ii) If* $\varphi\colon D_1 \to D_2$ *is a planar isotopy as above* except *that it interchanges the vertical positions of two generators* $p$ *and* $q$*, with* $p$ *vertically above* $q$*, then* $D_1 = \mu(\deg p, \deg q) D_2$ *in* $\mathbf{Diag}_d^\Lambda$*.*

*(iii) All the local relations in Fig. 1.3 above, viewed with a legal label.*

**Proposition 1.4.3.** $\mathbf{GFoam}_d$ *and* $\mathbf{Diag}_d^\Lambda$ *are isomorphic as* $(\mathbb{Z}^2, \mu)$*-graded-2-categories.*

*Proof.* The $\mathbb{Z}^2$-grading is preserved by correspondence between foams and diagrams. It remains to check that foam relations in Definition 1.3.3 correspond to diagrammatic relations in Definition 1.4.2.

Relations (i) in Definition 1.3.3 correspond to relations (i) in Definition 1.4.2, together with graded interchange laws that involve at least one crossing, and braid-like relations, pitchfork relations and dot slide.

Relations (ii) in Definition 1.3.3 correspond to graded interchange laws that only involve cups, caps and dots.

Relations in Fig. 1.2 correspond to relations in Fig. 1.3, except braid-like relations, pitchfork relations and dot slide. For instance, the horizontal neck-cutting corresponds to the evaluation of clockwise bubbles. □

Note that none of the relations in Definition 1.4.2 depends on the objects, that is, the label of the foam diagrams. In the sequel, we mostly leave foam diagrams unlabelled, as the label is either irrelevant to the discussion or understood from the context. In particular, this applies when computing secondary relations from the defining relations in $\mathbf{Diag}_d^\Lambda$, as in the following lemma:

**Lemma 1.4.4.** *The following local relations hold in* $\mathbf{Diag}_d^\Lambda$ *for any choice of legal label (omitted below):*

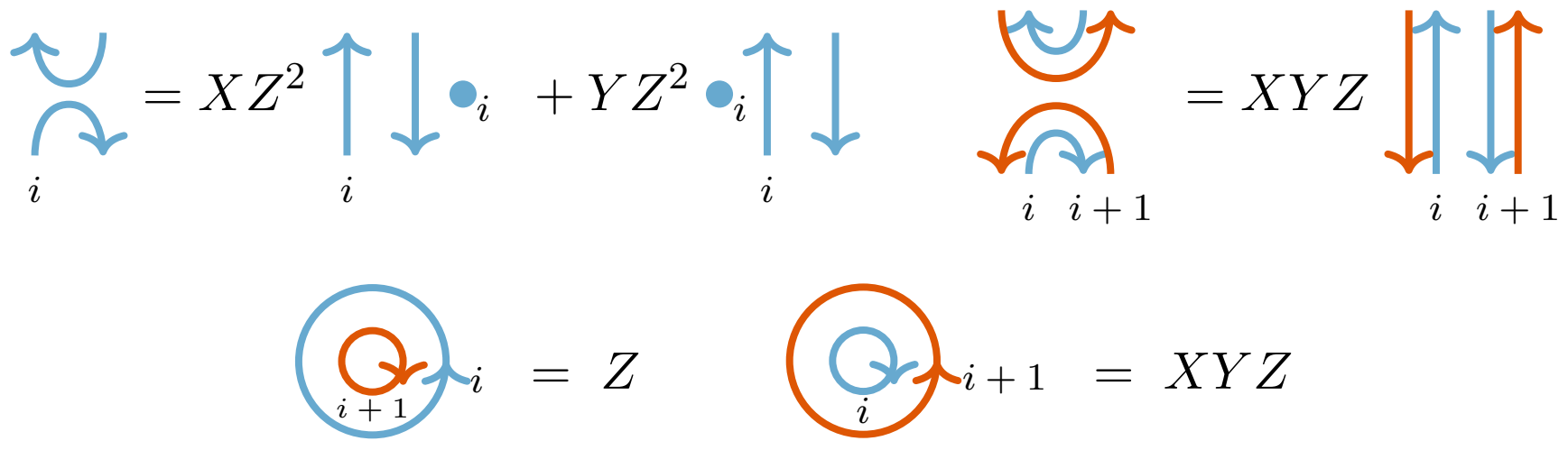

## 1.5 Shading diagrammatics

As we often work without labels, it is useful to know when a foam diagram *could* be labelled with a legal label. For that, we introduce an alternative diagrammatics, the *shading diagrammatics*. This corresponds to recording 2-facets instead of seams. For example:

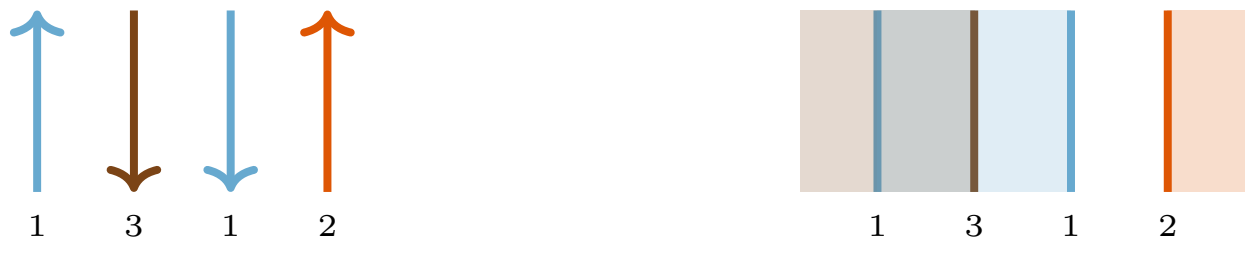

an identity diagram that admits a legal label

the same identity diagram using shading diagrammatics

We get the following characterisation:

**Lemma 1.5.1.** *A foam diagram admits a legal label if and only if the following conditions hold:*

1. *EXISTENCE OF SHADINGS: for all $i$, we can partially shade the regions delimited by the $i$-strand such that if we follow any $i$-strand in the direction of its orientation, the right-hand side is shaded while the left-hand side is not. We call this data a* shading.
2. *SUPERPOSITION OF SHADINGS: the $i$-shaded regions and the $j$-shaded regions do not overlap, for all $i$ and $j$ such that $|i-j|=1$.*
3. *DOTS AND SHADINGS: $i$-dots never sits in $j$-shaded regions, for all $i$ and $j=i-1,i$.*

*We refer to the above conditions as* shading conditions.

Here are some of the defining relations using shading diagrammatics:

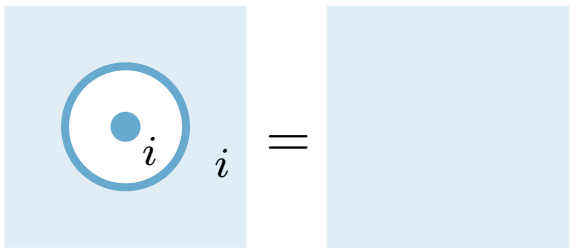

evaluation of dotted counter-clockwise bubble

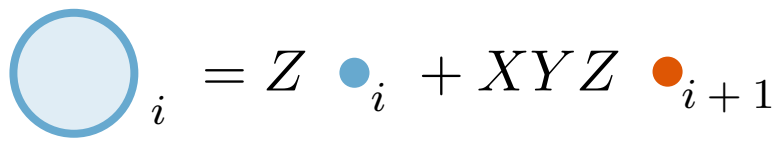

evaluation of clockwise bubble

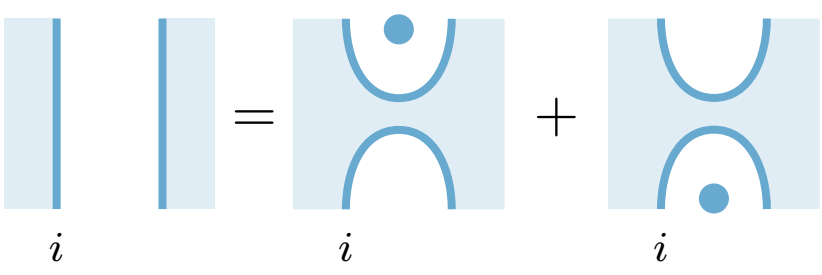

neck-cutting relation

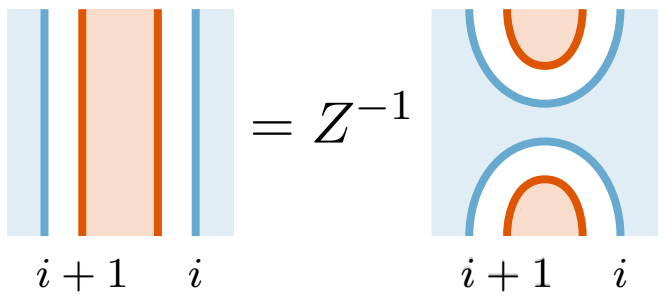

squeezing relation

Shading diagrammatics stand closer to the original foam. However, it is heavier than string diagrammatics, and we shall prefer the latter in most part of the text. However, we encourage the reader to think with both diagrammatics.

## 1.6 A basis for graded $\mathfrak{gl}_2$-foams

Each 2Hom-space of $\mathbf{GFoam}_d$ admits a basis with the following one-line description:

> SLOGAN. *A hom-basis in* $\mathbf{GFoam}_d$ *is given by picking a foam whose underlying surface is a union of disks, and considering all the ways to put dots on the disks.*

This kind of family of foams is formalized below as a *reduced family*, and foams appearing in such families as *reduced foams*. In this section, we describe reduced families and sketch a candidate reduction algorithm to rewrite any foam as a linear combination of reduced foams. In part II, this algorithm is formalized as a *Gray linear rewriting system modulo*. It is shown to terminate (each foam can be written as a linear combination of reduced foams) and to confluate (this decomposition is unique). In other words:

**Theorem 1.6.1** (Basis theorem for graded $\mathfrak{gl}_2$-foams)**.** *Let* $W, W' \colon \mu \to \lambda$ *be two parallel webs.* $W$ *and* $W'$ *admit a reduced family and every reduced family constitutes a basis of* $\mathrm{Hom}_{\mathbf{GFoam}_d}(W, W')$.

### 1.6.1 Reduced families

Let $W, W' \colon \mu \to \lambda$ be two parallel webs. We denote $\mathfrak{sl}(W) \sqcup_\partial \mathfrak{sl}(W')$ the closed 1-manifold obtained by glueing $\mathfrak{sl}(W)$ and $\mathfrak{sl}(W')$ along their common boundary points. Note that if $F \colon W \to W'$ is a foam, then $\partial(\mathfrak{sl}(F))$ is homeomorphic to $\mathfrak{sl}(W) \sqcup_\partial \mathfrak{sl}(W')$. A *reduced foam* is a foam $F$ such that $\mathfrak{sl}(F)$ is a (dotted) union of disks with at most one dot on each disk. In that case, there is a bijection between $\pi_0(\mathfrak{sl}(F))$, the connected components of $\mathfrak{sl}(F)$, and $\pi_0(\mathfrak{sl}(W) \sqcup_\partial \mathfrak{sl}(W'))$, the connected components of $\mathfrak{sl}(W) \sqcup_\partial \mathfrak{sl}(W')$, mapping a disk to its boundary. For $\delta \subset \pi_0(\mathfrak{sl}(W) \sqcup_\partial \mathfrak{sl}(W'))$, we say that $\mathfrak{sl}(F)$ (or abusing terminology, $F$) is *$\delta$-dotted* if the following holds:

> *$\delta$-dotted:* there is a dot (resp. no dot) sitting on a disk $d$ in $\mathfrak{sl}(F)$ if and only if its boundary $\partial d \in \pi_0(\mathfrak{sl}(W) \sqcup_\partial \mathfrak{sl}(W'))$ is in $\delta$ (resp. is not in $\delta$).

We say that $F$ is *undotted* whenever it is $\emptyset$-dotted.

**Definition 1.6.2.** *Let $W, W' \colon \mu \to \lambda$ be two parallel webs. A* reduced family *is a family of foams $F_\delta \colon W \to W'$, indexed by subsets $\delta \subset \pi_0(\mathfrak{sl}(W) \sqcup_\partial \mathfrak{sl}(W'))$, such that each $F_\delta$ is a $\delta$-dotted reduced foam.*

Thanks to Lemma 1.3.6, reduced families are essentially unique, in the sense that if $\{F_\delta\}_\delta$ and $\{F'_\delta\}_\delta$ are two reduced families for the same hom-space, then $F_\delta \mathrel{\dot\sim} F'_\delta$ for all $\delta$. In particular:

**Lemma 1.6.3.** *Let $W, W' \colon \mu \to \lambda$ be two parallel webs. If a reduced family is a basis of the $\hom_{\mathbf{GFoam}_d}(W, W')$, then every reduced family is a basis of $\hom_{\mathbf{GFoam}_d}(W, W')$.* □

A choice of reduced family can be given by first choosing an undotted reduced foam $F$ and defining the $F_\delta$'s by adding dots sitting on the disks, as in the slogan above. Furthermore, $\mathfrak{sl}(F)$ can be chosen with the cellular form given in Fig. 1.4. The cross-section of the middle part consists of a union of closed intervals, without closed components. To each closed interval corresponds a disk whose boundary intersects both $\mathfrak{sl}(W)$ and $\mathfrak{sl}(W')$; call this kind of disks *through disks*. The rest of the disks in $\mathfrak{sl}(F)$ are cups or caps, whose boundary only intersect either $\mathfrak{sl}(W')$ or $\mathfrak{sl}(W)$. In the example of Fig. 1.4, there are exactly one through disk, one cup disk and one cap disk.

We can count the number of through disks. If we denote $\#(1:\lambda)$ the number of ones in an object $\lambda \in \underline{\Lambda}_d$, then:

$$\#\{\text{through disks}\} = [\#(1:\lambda) + \#(1:\mu)]/2.$$

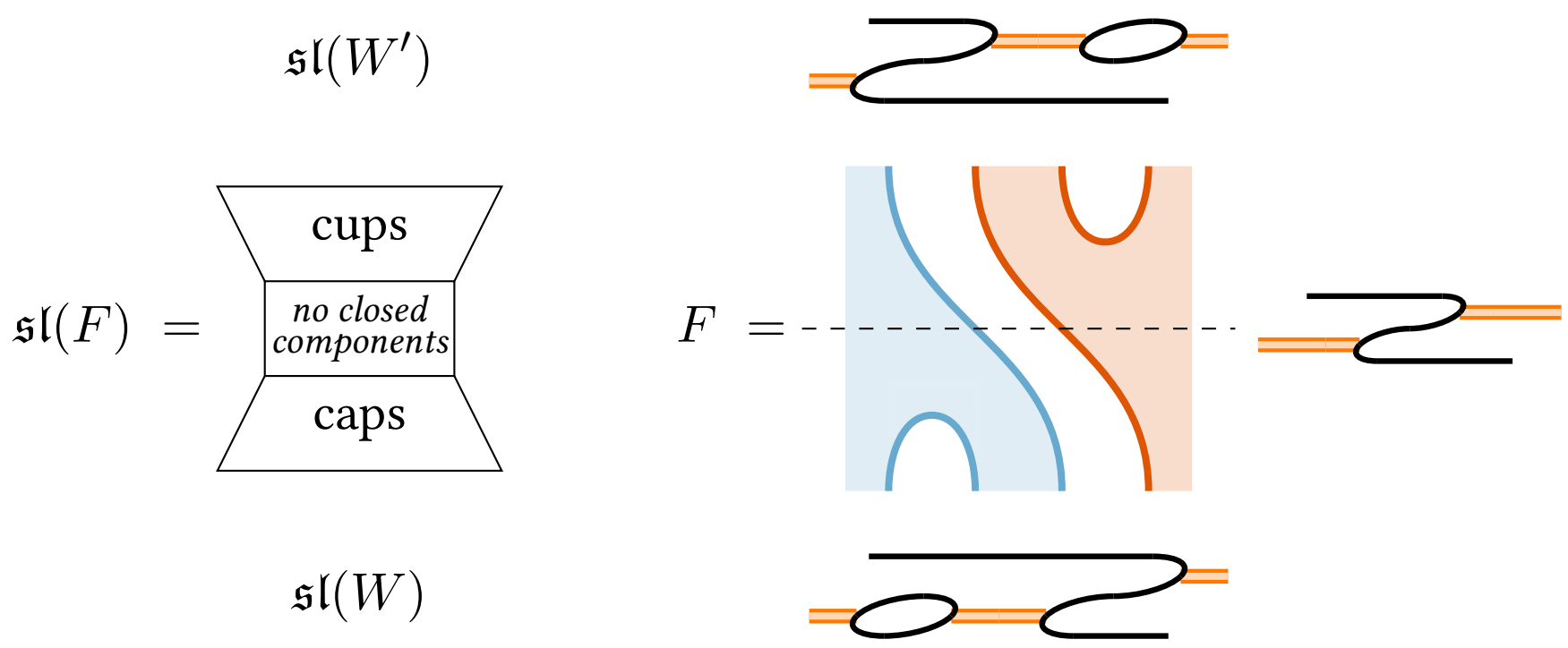

**Fig. 1.4** Undotted reduced foams, schematic (on the left) and example (on the right).

As shown in [17, Lemma 2.5], counting the number of closed components in $\mathfrak{sl}(W) \sqcup_\partial \mathfrak{sl}(W')$ defines a non-degenerate pairing on $\mathfrak{gl}_2$-webs:

$$\langle W, W' \rangle := (q + q^{-1})^{\#[\mathfrak{sl}(W) \sqcup_\partial \mathfrak{sl}(W')]}.$$

This pairing coincides with the web evaluation formula given in [163]. See also [17, p. 1315].

We now have all the ingredients to compute the graded dimension[4] of

$$\mathrm{Hom}_{\mathbf{GFoam}_d}(W, W').$$

It factors into the contributions of each disk in $\mathfrak{sl}(F)$: either dotted or undotted. If the disk is a cup or a cap, then these two possibilities contribute with $q + q^{-1}$. If the disk is a through disk, they contribute with $(q^2 + 1)$ to the graded rank.

This shows the following corollary:

**Corollary 1.6.4** (Non-degeneracy of graded $\mathfrak{gl}_2$-foams)**.** *Let $W$ and $W'$ be two webs with boundary points $\lambda, \mu$ in $\underline{\Lambda}_d$. Then the space* $\mathrm{Hom}_{\mathbf{GFoam}_d}(W, W')$ *is a free* $\Bbbk$*-module and:*

$$\mathrm{gdim}_q \mathrm{Hom}_{\mathbf{GFoam}_d}(W, W') = q^{[\#(1:\lambda) + \#(1:\mu)]/2} \langle W, W' \rangle,$$

*where* $\mathrm{gdim}_q$ *denotes the graded rank with respect to the $q$-grading.* □

[4] See subsection i.1.4, Example 2', for the definition.

### 1.6.2 A reduction algorithm

Given two parallel webs $W, W' \colon \lambda \to \mu$, we define a non-deterministic algorithm, denoted $\mathtt{T}^+$, on $V := \mathrm{Hom}_{\mathbf{GFoam}_d}(W, W')$. Recall the notion of foam isotopies, which in the string diagrammatics correspond to compositions or graded interchanges, braid-like relations, pitchfork relations, zigzag relations and dot slide. Each step of the algorithm has the following form:

INPUT: $u \in V$. Then:

(1) apply an arbitrary foam isotopy to the foams in the linear decomposition of $u$, returning a vector $u' \in V$;
(2) apply one of the following relations, *from left to right*, on one of the foams appearing in the linear decomposition of $u'$:
   - dot annihilation,
   - dot migration,
   - evaluation of a counter-clockwise bubble,
   - evaluation of a clockwise bubble,
   - neck-cutting relation coloured $i$, *provided the two pieces of $i$-strands belong to distinct $i$-strands*,
   - squeezing relation coloured $(i, i+1)$, *provided the two pieces of $i$-strands belong to distinct $i$-strands*,

   returning a vector $v' \in V$;
(3) apply, as in step (1), an arbitrary foam isotopy to the foams in the linear decomposition of $v'$, returning a vector $v \in V$;

OUTPUT: $v$. ■

If no relation in (2) can be applied, no matter the foam isotopy one applies in (1), then the algorithm terminates. In that case, we say that the foam cannot be further reduced by $\mathtt{T}^+$. We shall see below that:

**Proposition 1.6.5.** *A foam is reduced in the sense of subsection 1.6.1 if and only if it cannot be further reduced by* $\mathtt{T}^+$.

In part II, we formalize the algorithm $\mathtt{T}^+$ as a linear Gray rewriting system modulo, whose monomial normal forms are precisely reduced foams. In particular, we show that $\mathtt{T}^+$ terminates (Proposition 7.2.4), and that if $W, W' \colon \mu \to \lambda$ are two parallel webs, and $B$ is a set of unique isotopy representative for reduced foams, then $B$ is a basis of $\mathrm{Hom}_{\mathbf{GFoam}_d}(W, W')$ (Theorem 7.2.5). Together with Lemma 1.6.3, this shows Theorem 1.6.1.

The remaining of this section is the proof of Proposition 1.6.5.

**Definition 1.6.6.** *Let $F$ be foam. A* generalized bubble $\psi$ *is the data of dots and closed strands in $F$ (viewed in its string diagrammatics), such that if a dot or a closed strand $x$ is in the interior of a closed strand in $\psi$, then $x$ is in $\psi$.*

**Lemma 1.6.7.** *Let $F$ be a foam. The algorithm $\mathtt{T}^+$ evaluates every generalized bubble $\psi$ to a linear combination of dots:*

$$F \quad \to_{\mathtt{T}^+} \quad \sum_{dots} F \setminus \psi \cup \{dots\}.$$

*More precisely, $F$ reduces either to zero or to a linear combination of diagrams, obtained from $F$ by removing the data of $\psi$ and (possibly) adding dots.*

The above lemma only says that an evaluation exists: it does not imply that it is unique.

*Proof.* Pick a closed strand $S$ in $\psi$. If $S$ is a bubble with an arbitrary number of dots in its interior, it can be evaluated to dots using a bubble evaluation, together with dot migration, dot annihilation and dot slide. If not, the only situation that could prevent $S$ to be isotoped to a bubble is if some braid-like relation has a generalized bubble $\psi'$ in its center. By assumption, $\psi' \subset \psi$, so that we can inductively assume that $\psi'$ only consists of dots. Using dot migration and dot slide, these dots can be slid to the boundary of the braid-like relation. Hence, we can inductively reduce any closed strand in $\psi$ to a bubble, which can then be evaluated to dots. □

**Corollary 1.6.8.** *A reduced foam does not have any closed seam.* □

**Lemma 1.6.9.** *In a reduced foam, two $i$-shadings are always separated by at least one $(i-1)$-shading. In particular, a reduced foam has at most one $1$-shading.*

*Proof.* Denote by $S_1$ and $S_2$ two $i$-shadings in $F$. If $S_1$ and $S_2$ can be isotoped next to each other, they can be rewritten using a neck-cutting relation. Hence, they must be separated by a $(i-1)$-shading or by a $(i+1)$-shading. If $i = d-1$, the statement holds as $d$-shadings do not exist. Assume by induction that the statement holds for $j > i$, and that $S_1$ and $S_2$ are only separated by $(i+1)$-shadings. By the induction hypothesis, these $(i+1)$-shadings are themselves separated by $i$-shadings. But that implies that a squeezing relation can be applied. We conclude that $S_1$ and $S_2$ are separated by at least one $(i-1)$-shading. □

*Proof of Proposition 1.6.5.* If a foam is reduced then it cannot be further reduced by $\mathtt{T}^+$. Consider then a foam $F$ that cannot be further reduced by $\mathtt{T}^+$. Denote $\mathfrak{sl}(F)_{\leq i}$ the surface obtained from $\mathfrak{sl}(F)$ by discarding all 1-facets labelled $j$ for $j > i$. We proceed by induction to show the following:

> $P(i)$: $i$-labelled seams all belong to distinct connected components when view ed as sitting on $\mathfrak{sl}(F)_{\leq i}$.

$P(1)$ follows from the fact that $F$ has at most one 1-shading, and hence at most one 1-labelled seam. Assume $P(j)$ holds for $j < i$. Let $S_1$ and $S_2$ be two $i$-labelled seams, and $F_1$ (resp. $F_2$) the $i$-labelled 1-facet incident to $S_1$ (resp. $S_2$). Since $S_1$ and $S_2$ are separated by some $(i-1)$-shading (Lemma 1.6.9), $F_1$ and $F_2$ are incident to (necessarily distinct) $(i-1)$-labelled seams. $P(i-1)$ implies that $F_1$ and $F_2$ belong to distinct connected components.

Finally, assume $\mathfrak{sl}(F)$ has a connected component which is not a disk. This component cannot be a sphere, because $F$ does not have any closed seam (Corollary 1.6.8). Consider then a non-contractible loop on this component. There exists a label $i$ such that twice, the loop intersects some $i$-labelled seam. We can inductively assume that the two intersection points belong to distinct $i$-seams; otherwise, we can cut the loop into two loops, such that one of them is non-contractible and intersect strictly fewer seams. We conclude that when viewed as sitting on $\mathfrak{sl}(F)_{\leq i}$, these two $i$-seams belong to the same component. This contradicts $P(i)$. □

## 1.7 The categorification theorem

Recall the definition of the Grothendieck ring of a graded-2-category and related notions from subsection 1.1.2. We decategorify $\mathbf{GFoam}_d$ with respect to the quantum grading (see Remark 1.3.5), that is:

$$K_0(\mathbf{GFoam}_d)|_q = K_0(\underline{(\mathbf{GFoam}_d)_q^{\oplus}}).$$

As explained in subsection 1.1.2, we call it the *quantum Grothendieck ring* of $\mathbf{GFoam}_d$. It has the structure of a $\mathbb{Z}[q, q^{-1}]$-linear category.

**Theorem 1.7.1** (Categorification theorem)**.** *The graded-2-category of foams categories the category of webs:*

$$K_0(\mathbf{GFoam}_d)|_q \cong \mathbf{Web}_d,$$

*where the isomorphism is an isomorphism of* $\mathbb{Z}[q, q^{-1}]$*-linear categories.*

This generalizes the analogous statement in the non-graded case, see [17, Theorem 2.11]. The proof follows the same line of thought, borrowing the general strategy from [107].

*Proof of Theorem 1.7.1.* Let $\gamma\colon \mathbf{Web}_d \to K_0(\mathbf{GFoam}_d)|_q$ mapping a web $W$ to its image $[W]$ in $K_0(\mathbf{GFoam}_d)|_q$. The neck-cutting relation, the squeezing relation and the Reidemeister II relation in Fig. 1.3 show that the web relations lift to isomorphisms in $\underline{(\mathbf{GFoam}_d)}_q^{\oplus}$, so that $\gamma$ is well-defined. By Corollary 1.6.4, $\gamma$ preserves a non-degenerate pairing, so it must be injective. □

# PART I
# Odd Khovanov homology
# and
# higher representation theory

# 2

# Covering Khovanov homology

In this chapter, we define an invariant of oriented tangles and show that it coincides with odd Khovanov homology when restricted to links. More generally, we define an invariant that coincides with *covering Khovanov homology*, an invariant of links defined by Putyra [155]. Both constructions are defined over the ring $\Bbbk = \mathbb{Z}[X, Y, Z^{\pm}](X^2 = Y^2 = 1)$ (see Definition 1.3.2). Setting $X = Y = Z = 1$ at the level of chain complexes recovers (even) Khovanov homology, while setting $X = Z = 1$ and $Y = -1$ (still at the level of chain complexes) recovers odd Khovanov homology.

To distinguish the two constructions, we call Putyra's construction *covering* $\mathfrak{sl}_2$*-Khovanov homology* and denote it $\mathrm{CKh}_{\mathfrak{sl}_2}(L)$ for an oriented link $L$, and we call our construction *covering* $\mathfrak{gl}_2$*-Khovanov homology* and denote it $\mathrm{CKh}_{\mathfrak{gl}_2}(T)$ for an oriented tangle $T$. The latter coincides with *not even Khovanov homology* [183] once we set $X = Z = 1$ and $Y = -1$ at the level of chain complexes. See also section i.3 for connections with the work of Naisse and Putyra [144].

To state our claim precisely, we introduce the following completions:

**Definition 2.0.1.** *The set* $\underline{\Lambda}$*, the* $\mathbb{Z}[q, q^{-1}]$*-linear category* $\mathbf{Web}$ *and the* $(\mathbb{Z}^2, \mu)$*-graded-2-category* $\mathbf{GFoam}$ *are respectively defined as:*

$$\underline{\Lambda} \coloneqq \mathrm{colim}(\ldots \hookrightarrow \underline{\Lambda}_d \hookrightarrow \underline{\Lambda}_{d+2} \hookrightarrow \ldots),$$
$$\mathbf{Web} \coloneqq \mathrm{colim}(\ldots \hookrightarrow \mathbf{Web}_d \hookrightarrow \mathbf{Web}_{d+2} \hookrightarrow \ldots),$$
$$\mathbf{GFoam} \coloneqq \mathrm{colim}(\ldots \hookrightarrow \mathbf{GFoam}_d \hookrightarrow \mathbf{GFoam}_{d+2} \hookrightarrow \ldots),$$

*where the embeddings denote the addition of a double point, a double line and a double facet on top.*

The fact that the above indeed are embeddings follows from Lemma 1.2.4 and Theorem 1.6.1.

As an example, in the category $\mathbf{Web}$ the following identity webs are identified:

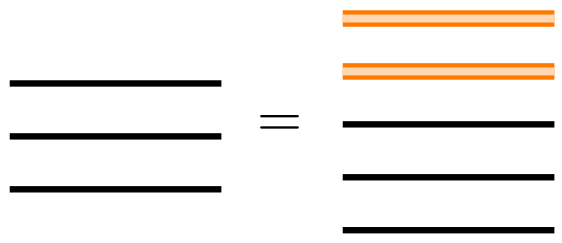

in $\mathbf{Web}$.

Informally, working in $\underline{\Lambda}$, $\mathbf{Web}$ and $\mathbf{GFoam}$ means that one can always "add a double point, a double line or a double facet on top".

For $\mathcal{C}$ a $(G,\mu)$-graded-2-category, we denote by $\mathrm{Ch}_\bullet(\mathcal{C})$ the $\Bbbk$-linear category of chain complexes in $\mathcal{C}$ and chain morphisms. We say that a chain morphism is *degree-preserving* if each of its components is degree-preserving. Note that a priori, $\mathrm{Ch}_\bullet(\mathcal{C})$ does not inherit a $G$-grading from $\mathcal{C}$.

Recall from subsection 1.1.2 that $\underline{(\mathbf{GFoam})|_q^\oplus}$ denotes the additive $q$-shifted closure of $\mathbf{GFoam}$: it allows formal shifts of webs in the quantum grading, restrict graded foams to those preserving the quantum grading, and allows formal direct sums. As before, we view the quantum grading as independent from the $\mathbb{Z}^2$-grading; hence $\underline{(\mathbf{GFoam})|_q^\oplus}$ is still a $(\mathbb{Z}^2,\mu)$-graded-2-category (for $\mu$ as in Definition 1.3.2), and shifts in quantum degree do not affect the $\mathbb{Z}^2$-degree.

For each sliced oriented tangle diagram $D_T$ representing an oriented tangle $T$, we define in subsection 2.1.2 a chain complex

$$\mathrm{Kom}_{\mathfrak{gl}_2}(D_T) \in \mathrm{Ch}_\bullet(\underline{(\mathbf{GFoam}_d)|_q^\oplus}).$$

Then:

**Theorem 2.0.2.** *Let $D_T$ be a sliced oriented tangle diagram presenting an oriented tangle $T$, and denote $N_+$ and $N_-$ the number of positive and negative crossings, respectively. The homotopy type of $q^{-2N_+ + N_-}t^{N_+}\mathrm{Kom}_{\mathfrak{gl}_2}(D_T)$, denoted $\mathrm{CKh}_{\mathfrak{gl}_2}(T)$, is an invariant of the oriented tangle $T$.*

The proof of Theorem 2.0.2 is given in subsection 2.1.3. This construction is a graded analogue of the construction given in [115]. Crucially, it uses a graded horizontal composition (or "object-adapted" tensor product) of chain complexes, for which we give a minimal introduction in subsection 2.1.1.

In section 2.2, we review the definition of covering $\mathfrak{sl}_2$-Khovanov homology. For each oriented link diagram $D_L$, it associates a complex $\mathrm{Kom}_{\mathfrak{sl}_2}(D_L)$ in $\mathrm{Ch}_\bullet(\Bbbk\text{-Mod}_{\mathbb{Z}})$, the category of chain complexes in $\mathbb{Z}$-graded $\Bbbk$-modules. The homotopy type of $q^{-2N_+ + N_-}t^{N_+}\mathrm{Kom}_{\mathfrak{sl}_2}(D_L)$ is an invariant of the oriented link $L$.

Finally, section 2.3 shows the equivalence between the two constructions, when restricted to links. To state it, denote by $\emptyset \in \underline{\Lambda}$ the empty weight and by $\emptyset := \mathrm{id}_\emptyset$ its identity, the empty web. Recall that in $\underline{\Lambda}$ (resp. in $\mathbf{Web}$), the empty weight (resp. the empty web) is the same as an arbitrary vertical juxtaposition of double points (resp. double lines). Denote by $\mathbf{GFoam}(\emptyset,\emptyset)$ the $\mathbb{Z}$-graded $\Bbbk$-linear category obtained by restricting $\mathbf{GFoam}$ to the object $\emptyset$ (the $\mathbb{Z}$-grading being the quantum grading), and let

$$\mathcal{A}_{\mathfrak{gl}_2} : \mathbf{GFoam}(\emptyset,\emptyset) \to \Bbbk\text{-Mod}_{\mathbb{Z}}$$

be the representable functor $\mathcal{A}_{\mathfrak{gl}_2} := \mathrm{Hom}_{\mathbf{GFoam}(\emptyset,\emptyset)}(\emptyset, -)$. It canonically extends to a functor

$$\mathcal{A}_{\mathfrak{gl}_2} : \mathrm{Ch}_\bullet(\underline{(\mathbf{GFoam})|_q^{\oplus}})(\emptyset,\emptyset) \to \mathrm{Ch}_\bullet(\Bbbk\text{-Mod}_{\mathbb{Z}}).$$

We can now state the main result of this chapter:

**Theorem 2.0.3.** *Let $D_L$ be a sliced link diagram presenting an oriented link $L$. We have the following degree-preserving isomorphism of chain complexes of $\Bbbk$-modules:*

$$\mathcal{A}_{\mathfrak{gl}_2}(\mathrm{Kom}_{\mathfrak{gl}_2}(D_L)) \cong \mathrm{Kom}_{\mathfrak{sl}_2}(D_L).$$

This theorem is the content of Main theorem C when setting $X = Z = 1$ and $Y = -1$ at the level of chain complexes.

## 2.1 A covering $\mathfrak{gl}_2$-Khovanov homology for oriented tangles

### 2.1.1 Composition of hypercubic chain complexes

We describe the horizontal composition of two *hypercubic complexes* in a given graded-2-category. Hypercubic complexes are special cases of homogeneous polycomplexes, introduced in full generality in Definition 4.1.1. The horizontal composition that we describe is the specialization of the definitions Definition 4.1.3 and Definition 4.1.6 (adapted to the context of graded-2-categories). This is the minimal description necessary for the construction of covering $\mathfrak{gl}_2$-Khovanov homology.

*Notation* 2.1.1. Fix $N \in \mathbb{N}$. We use the shorthand $\mathcal{I} := \{1, \dots, N\}$ for the set of indices $1 \leq i \leq N$. We view $\{0,1\}^N$ as a hypercubic lattice and denote $(\boldsymbol{e}_i)_{i \in \mathcal{I}}$ the canonical basis of $\mathbb{Z}^N$. For each $\boldsymbol{r} \in \{0,1\}^N$ and each $i \in \mathcal{I}$, we write $\boldsymbol{r} \to \boldsymbol{r} + \boldsymbol{e}_i$ the corresponding edge in the hypercube $\{0,1\}^N$.

Fix $\mathcal{C}$ a $(G, \mu)$-graded-2-category. Whenever we write a composition in $\mathcal{C}$, it is tacitly assumed that the 1-morphisms or 2-morphisms involved are composable.

**Definition 2.1.2.** *A* hypercubic complex *in $\mathcal{C}$ and of dimension $N$ is a pair $\mathbb{A} = (A, \alpha)$ consisting of the following data:*

*(i) for each vertex $\boldsymbol{r} \in \{0,1\}^N$, a 1-morphism $A^{\boldsymbol{r}}$ in $\mathcal{C}$,*

*(ii) for each edge $\boldsymbol{r} \to \boldsymbol{r} + \boldsymbol{e}_i$, a homogeneous 2-morphism $\alpha_i^{\boldsymbol{r}} \colon A^{\boldsymbol{r}} \to A^{\boldsymbol{r}+\boldsymbol{e}_i}$ in $\mathcal{C}$, such that each square anti-commutes:*

$$\alpha_{i_2}^{\boldsymbol{r}+\boldsymbol{e}_{i_1}} \star_1 \alpha_{i_1}^{\boldsymbol{r}} = -\alpha_{i_1}^{\boldsymbol{r}+\boldsymbol{e}_{i_2}} \star_1 \alpha_{i_2}^{\boldsymbol{r}}$$

*for all suitable $\boldsymbol{r} \in \{0,1\}^N$ and $i_1, i_2 \in \mathcal{I}$.*

*(iii) The grading is constant in a given direction, in the sense that for any $i \in \mathcal{I}$, either $\alpha_i^{\boldsymbol{r}} = \alpha_i^{\boldsymbol{s}} = 0$ for all $\boldsymbol{r}, \boldsymbol{s} \in \{0,1\}^N$, or $\deg \alpha_i^{\boldsymbol{r}} = \deg \alpha_i^{\boldsymbol{s}}$ for all $\boldsymbol{r}, \boldsymbol{s} \in \{0,1\}^N$.*

Given such a hypercubic complex $\mathbb{A} = (A, \alpha)$, we define the following element of $G$:

$$|\alpha|(\boldsymbol{r}) := \sum_{i \colon \boldsymbol{r}_i = 1} \deg_G(\alpha_i^{\boldsymbol{0}}),$$

where $\boldsymbol{0} := (0, \dots, 0) \in \{0,1\}^N$ and we used the convention $\deg_G(0) = 0$. Alternatively, the element $|\alpha|(\boldsymbol{r})$ is the sum of the $G$-degrees along a path from $\boldsymbol{0}$ to $\boldsymbol{r}$.

**Definition 2.1.3.** *Let $\mathbb{A} = (A, \alpha)$ and $\mathbb{B} = (B, \beta)$ be two hypercubic complexes of dimensions $N$ and $M$ respectively. The* horizontal composition $\mathbb{A} \star_0 \mathbb{B}$ *of $\mathbb{A}$ and $\mathbb{B}$ is the hypercubic chain complex of dimension $N + M$ defined by the following data:*

*(i) on each vertex $(\boldsymbol{r}, \boldsymbol{s}) \in \{0,1\}^{n+m}$, the 1-morphism $A^{\boldsymbol{r}} \star_0 B^{\boldsymbol{s}}$;*

*(ii) on each edge* $(\boldsymbol{r},\boldsymbol{s}) \to (\boldsymbol{r},\boldsymbol{s}) + \boldsymbol{e}_k$*, the homogeneous 2-morphism*

$$(\alpha \star_0 \beta)_k^{(\boldsymbol{r},\boldsymbol{s})} := \begin{cases} \alpha_i^{\boldsymbol{r}} \star_0 \mathrm{id}_{B^{\boldsymbol{s}}} & \text{if } k = i \in \{1,\ldots,N\}, \\ (-1)^{|\boldsymbol{r}|}\mu\left(|\alpha|\,(\boldsymbol{r}), \beta_j^{\boldsymbol{s}}\right) \mathrm{id}_{A^{\boldsymbol{r}}} \star_0 \beta_j^{\boldsymbol{s}} & \text{if } k = j \in \{N+1,\ldots,N+M\}, \end{cases}$$

The sign appearing above is the *graded Koszul rule.* By Theorem 4.1.9, this horizontal composition is coherent with homotopies (see also Main theorem B).

Note that a length-two chain complex whose differential is homogeneous is exactly a hypercubic complex of dimension one. In particular, if $\mathbb{A}_1,\ldots,\mathbb{A}_N$ is a family of horizontally composable length-two chain complexes with homogeneous differentials, Definition 2.1.3 defines their $N$-fold horizontal composition.

### 2.1.2 Definition of covering $\mathfrak{gl}_2$-Khovanov homology

We now define a chain complex $\mathrm{Kom}_{\mathfrak{gl}_2}(D) \in \mathrm{Ch}_\bullet(\underline{(\mathbf{GFoam})}|_q^{\oplus})$ for every sliced tangle diagram $D$. In fact, $\mathrm{Kom}_{\mathfrak{gl}_2}(D)$ is independent of the orientation of $D$, and we shall not mention it again in this section. The reader can follow the procedure on the example given in Fig. 2.1, with $D$ pictured at step ①.

We shall need the following definitions of "mixed crossings", between a single line and a double line:

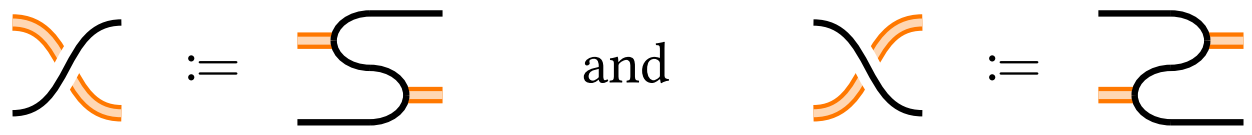

(2.1)

These mixed crossings satisfy the following relations in $\mathbf{Web}$:

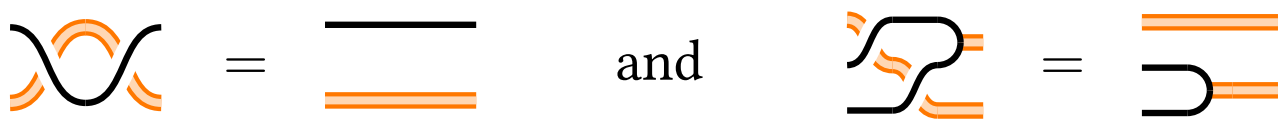

(2.2)

as well as all the relations obtained from the above by reflecting vertically and horizontally.

The procedure starts by telling how to assign a web to an elementary flat tangle diagram, that is, to a cap and a cup. There are more than one web $W$ such that $\mathfrak{sl}(W)$ is a cup (or a cap), but we fix a choice by fixing the endpoints. To enforce the given endpoints, we use the mixed crossings (2.1).

Say that $\lambda \in \underline{\Lambda}$ is *antidominant* if it is antidominant as a $\mathfrak{gl}_2$-weight, that is, if it is (non-strictly) increasing. To any set of $n$ points on a line corresponds

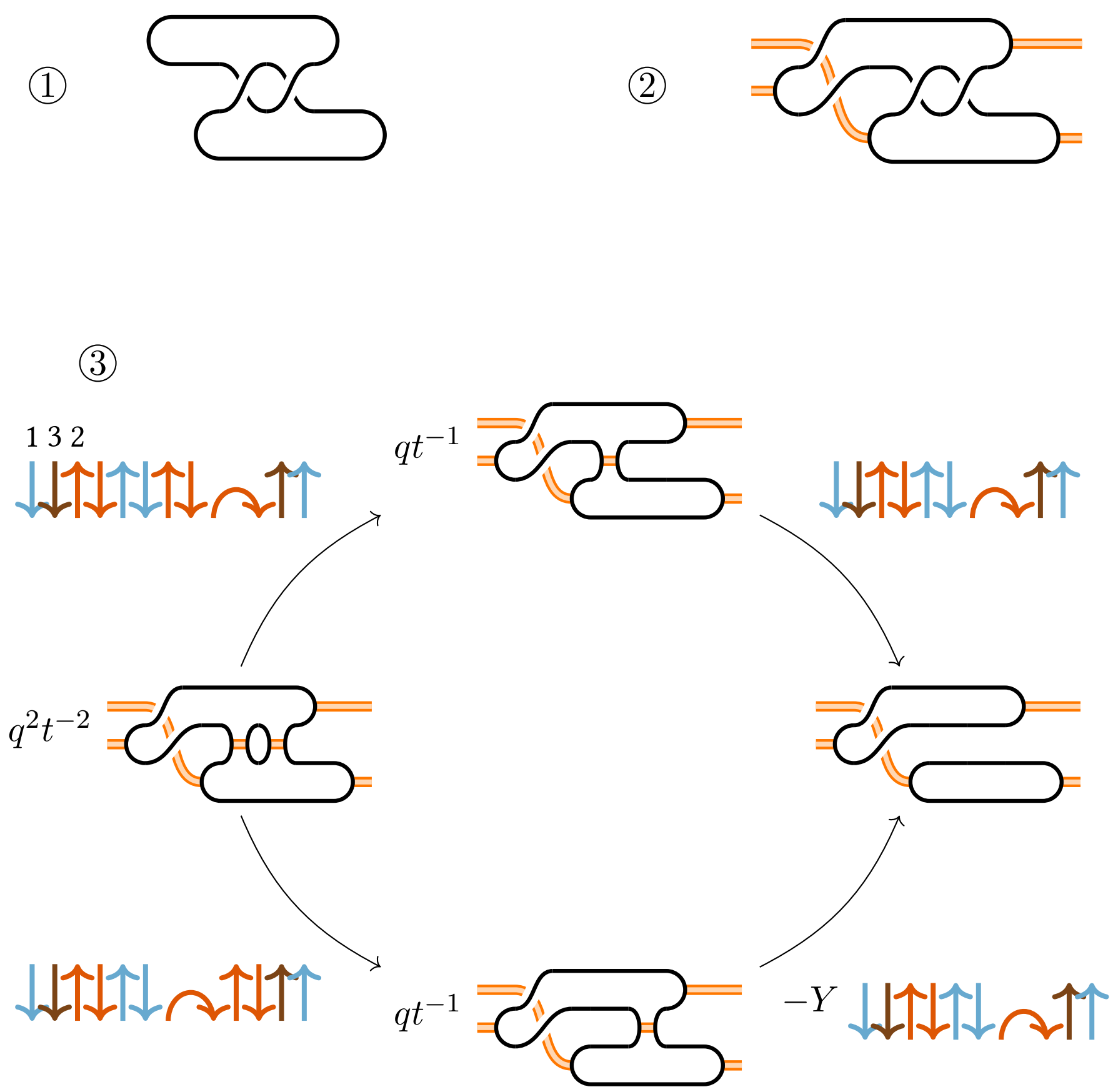

**Fig. 2.1** Defining procedure for $\mathrm{Kom}_{\mathfrak{gl}_2}$, in the case of a sliced tangle diagram presenting the Hopf link. The variables $q$ and $t$ respectively denote shift in quantum and homological degrees. Both differentials have $\mathbb{Z}^2$-degree $(0,1)$: one checks that the graded Koszul rule in Definition 2.1.3 adds the scalar $-Y$, as pictured in the figure.

a unique antidominant weight $\lambda \in \underline{\Lambda}_d$ for $n \leq d$ and $n = d \mod 2$. In turn, those antidominant weights define a unique element in $\underline{\Lambda}$. Given any elementary flat tangle diagram, we pick a web representative whose endpoints are antidominant by "adding a double line to the cup or cap and sliding it to the top". For instance:

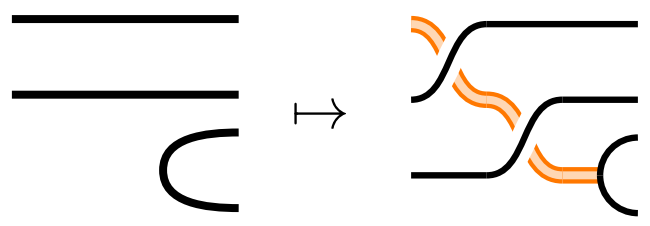

Note that fixing a choice for the endpoints ensures that if two elementary flat tangle diagrams are composable, then so are the corresponding webs in **Web**. We may extend this assignment to non-flat tangle diagrams by formally adjoining crossings to our web diagrammatics (see step ② in Fig. 2.1).

The procedure extends to chain complexes in $\mathrm{Ch}_\bullet(\underline{(\mathbf{GFoam})|_q^{\oplus}})$. For cups and caps, it is the chain complex concentrated in homological degree 0 corresponding to the previously assigned web. For crossings, this is given by the Khovanov–Blanchet bracket, generalized to the graded case:

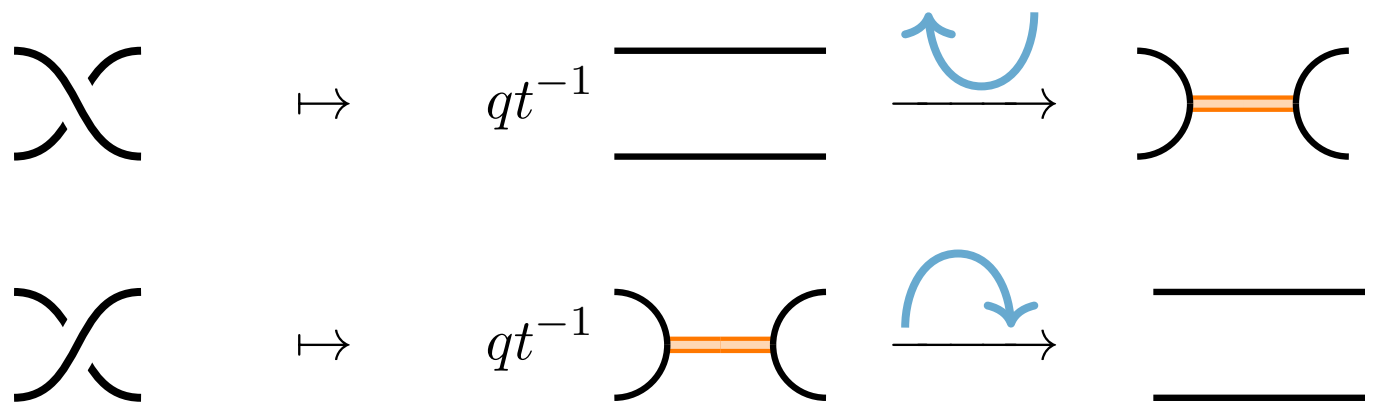

The variables $q$ and $t$ respectively denote shift in quantum and homological degrees. Note that the differentials are degree-preserving with respect to the quantum degree. However, they are not degree-preserving with respect to the $\mathbb{Z}^2$-degree: the former has $\mathbb{Z}^2$-degree $(1, 0)$, while the later has $\mathbb{Z}^2$-degree $(0, 1)$. If one restricts to the odd case ($X = Z = 1$ and $Y = -1$), the former has even parity and the latter has odd parity.

Finally, let $D$ be a sliced tangle diagram. Then $\mathrm{Kom}_{\mathfrak{gl}_2}(D)$ is defined as the horizontal composition (see Definition 2.1.3) of the chain complexes assigned to each slice of $D$. This ends the definition of $\mathrm{Kom}_{\mathfrak{gl}_2}(D)$ (see step ③ in Fig. 2.1). ◇

### 2.1.3 Proof of invariance

In this subsection, we prove Theorem 2.0.2.

Since the horizontal composition of chain complexes is coherent with homotopies (Theorem 4.1.9), the proof can be done locally. It suffices to check invariance under the Reidemeister–Turaev moves for sliced oriented tangle diagrams (see e.g. [146]). We split the proof in two lemmas, the first one (Lemma 2.1.4) dealing with planar isotopies and the second one (Lemma 2.1.6) dealing with Reidemeister moves.

**Lemma 2.1.4.** *Let $D_1$ and $D_2$ be two sliced tangle diagrams. If $D_1$ and $D_2$ are planar isotopic, then* $\mathrm{Kom}_{\mathfrak{gl}_2}(D_1)$ *and* $\mathrm{Kom}_{\mathfrak{gl}_2}(D_2)$ *are isomorphic.*

*Proof.* It suffices to check invariance under elementary planar isotopies for sliced tangle diagrams, as given by Fig. 2.2.

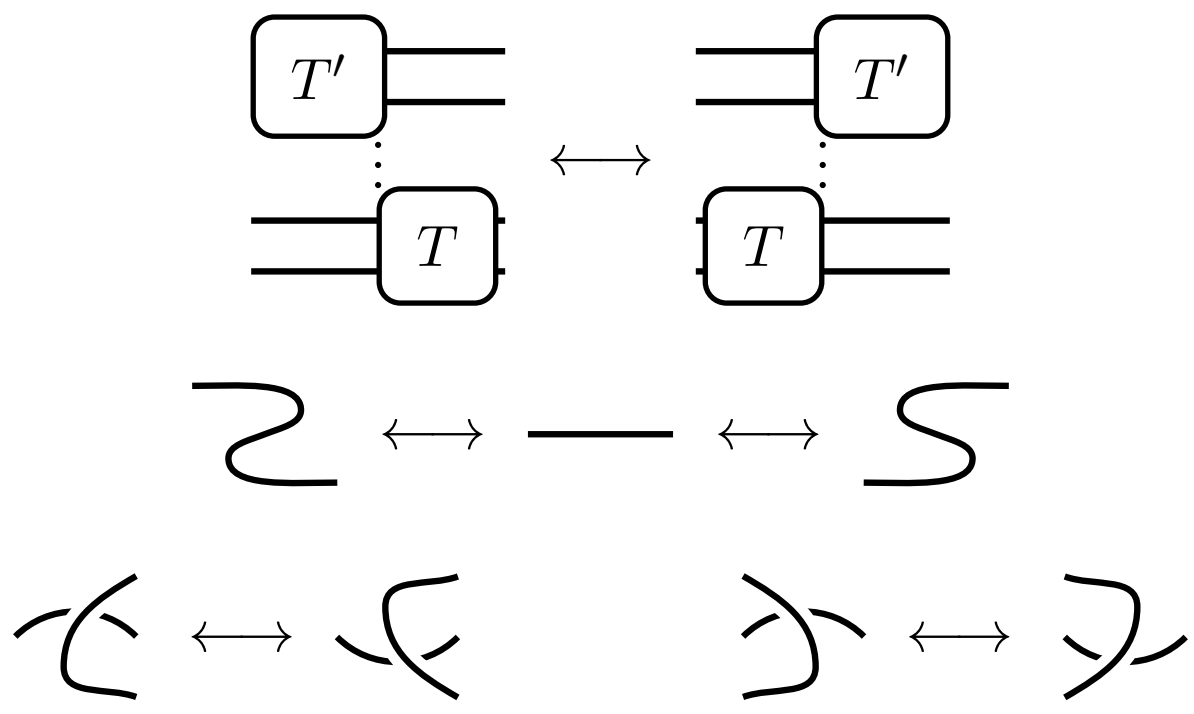

**Fig. 2.2** Elementary planar moves for sliced tangle diagrams. Here $T$ (resp. $T'$) denotes a crossing, an cup or an cap.

If the elementary planar isotopy does not involve any crossing, the complex is concentrated in a single homological degree. Finding an isomorphism of complexes reduces to finding an isomorphism of webs. Thanks to the categorification theorem (Theorem 1.7.1) and Lemma 1.2.4, two webs $W_1$ and $W_2$ are isomorphic in $\underline{(\mathbf{GFoam})}|_q^{\oplus}$ precisely if $\mathfrak{sl}(W_1)$ and $\mathfrak{sl}(W_2)$ are isotopic.

If the elementary planar isotopy contains exactly one crossing, the assigned complexes are length-two complexes with isotopic webs as vertices. Fix a pair of isomorphisms between these webs. Thanks to Lemma 1.3.6 (if two foams have the same underlying surface, there are equal up to invertible scalar), this pair defines a chain morphism up to invertible scalar. Renormalizing provides a genuine chain isomorphism.

Finally, the only elementary planar isotopy with at least two crossings is the planar isotopy interchanging two crossings. The associated chain complexes have the form $F \star_1 F'$ and $F' \star_1 F$ respectively, for some pair of length-two complexes $F$ and $F'$. A chain isomorphism $F \star_1 F' \cong F' \star_1 F$ is given by suitable pointwise compositions of foam crossings, exchanging $F$ with $F'$. □

To show invariance under Reidemeister moves, we follow Bar-Natan's strategy for the even case [13, 14], using delooping and gaussian elimination. *Delooping* denotes using circle evaluation in **Web** to remove a circle, while *gaussian elimination* is the following general homological fact:

**Lemma 2.1.5.** *In any additive category, if $\alpha$ is an isomorphism in the complex $C_\bullet$ given by*

$$\begin{array}{ccccc} & & X & & \\ & \overset{\alpha}{\nearrow} & & \overset{\gamma}{\searrow} & \\ W & & & & Z\,, \\ & \underset{\beta}{\searrow} & & \underset{\delta}{\nearrow} & \\ & & Y & & \end{array}$$

*then there exists a homotopy equivalence of complexes $C_\bullet \to C'_\bullet$ where $C'_\bullet$ is the complex $Y \xrightarrow{\delta} Z$. Moreover, this homotopy equivalence is a strong deformation retract.* □

One needs not know the definition of a strong deformation retract, except for its appearance in Lemma 2.1.7 below.

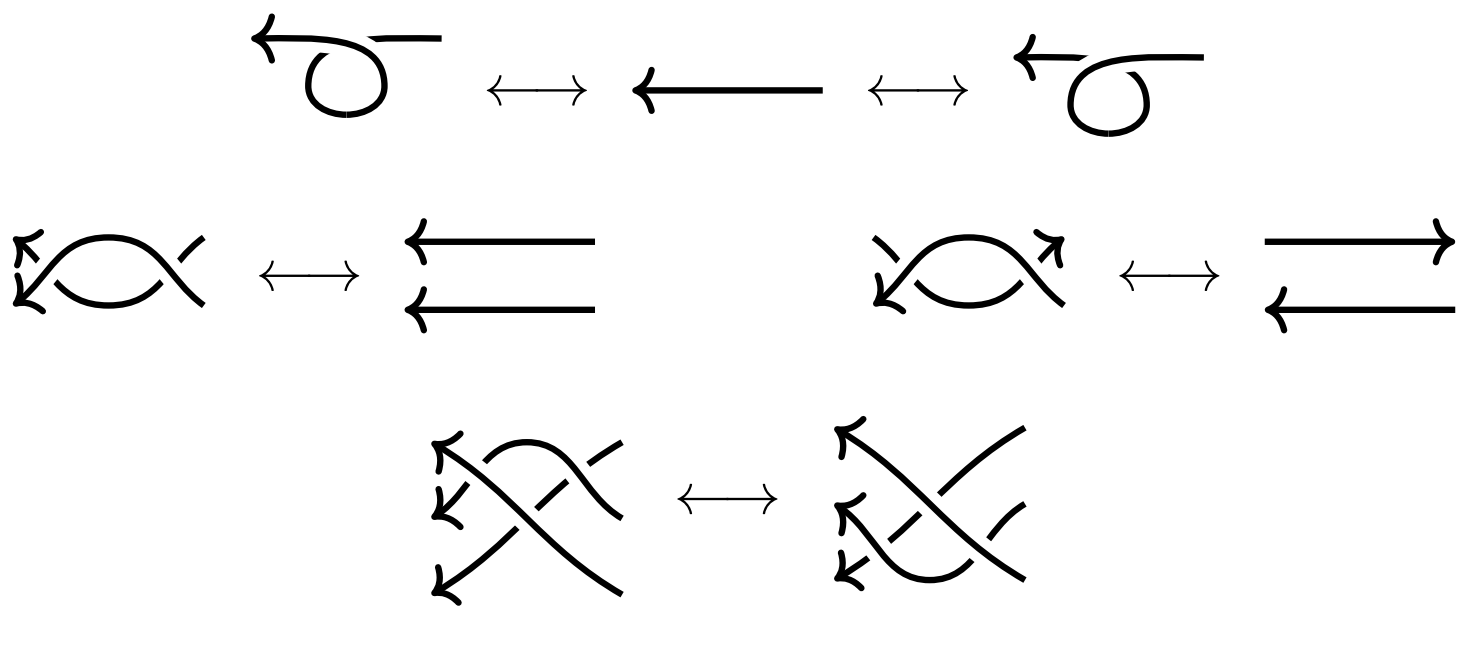

**Fig. 2.3** Reidemeister moves for sliced oriented tangle diagrams.

**Lemma 2.1.6.** *Let $D_1$ and $D_2$ be two oriented sliced tangle diagrams. If $D_1$ and $D_2$ are related by a Reidemeister move (see Fig. 2.3), then*

$$q^{-2N_+ + N_-} t^{N_+} \mathrm{Kom}_{\mathfrak{gl}_2}(D_1) \text{ and } q^{-2N_+ + N_-} t^{N_+} \mathrm{Kom}_{\mathfrak{gl}_2}(D_2)$$

*are homotopic.*

*Proof.* Consider first the Reidemeister I move. The proof of invariance is essentially contained in Fig. 2.4. On the left, the complex associated to the left-hand-side of the move. On the right, the same complex after delooping in homological degree zero, and simplifying in homological degree one. By Lemma 1.4.4, the bottom arrow is an isomorphism, so gaussian elimination gives a homotopy equivalence with the left-top web. It is shifted in quantum and homological degrees, but this is fixed by the renormalization. Invariance under the other Reidemeister I move is proved similarly.

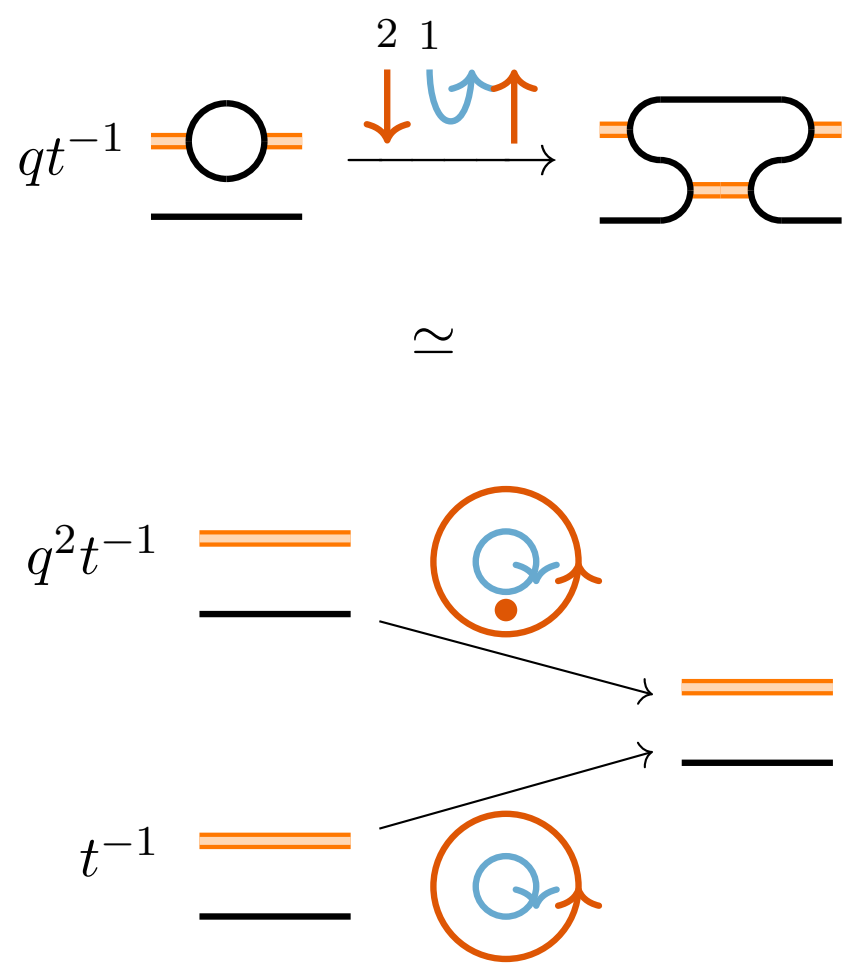

**Fig. 2.4** proof of invariance under Reidemeister I

Reidemeister II is proved similarly. The complex associated to the non-trivial side of this move is pictured in Fig. 2.5. Delooping can be applied to the bottom web, and gaussian elimination (together with zigzag relations) shows that this complex is homotopy equivalent to the top web. Renormalization concludes.

Following Bar-Natan [14], we use cones to simplify the proof of invariance under Reidemeister III. Denote $\Gamma(\psi)$ the cone associated to a morphism of complex $\psi$. Any hypercube $T$ of dimension $n$ can be seen as cone. Indeed, choose a direction $k$ in the hypercube. Ignoring the $k$-edges, $T$ breaks in two hypercubes $T_1$ and $T_2$ of dimension $n-1$. Then, switch the signs of all differentials in $T_1$. the hypercube $T_1$ is still a complex, but the $k$-faces of $T$ now commute instead of anti-commuting: the $k$-edges form a morphism

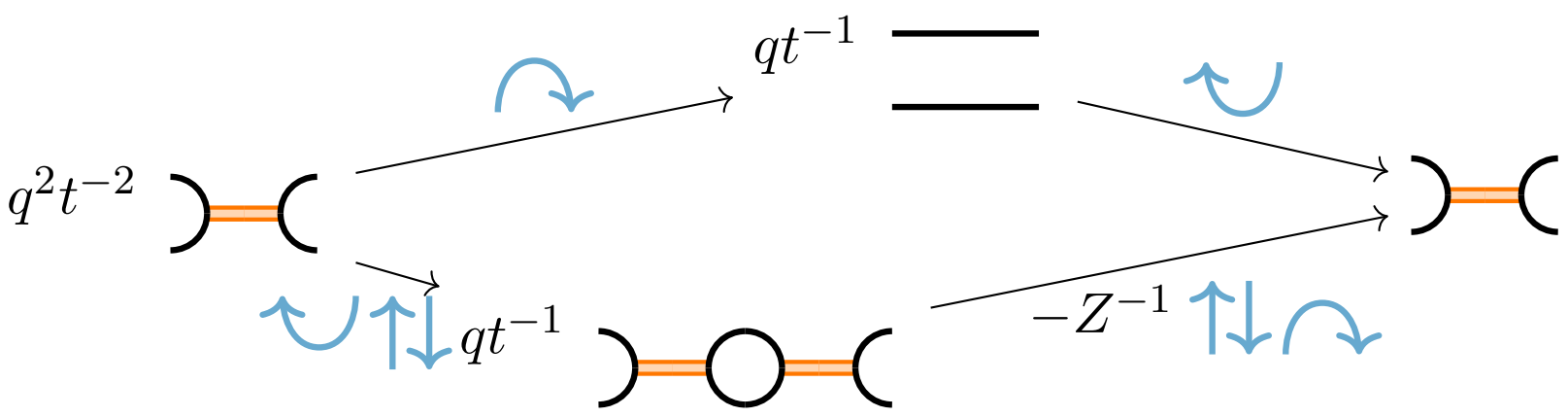

**Fig. 2.5** proof of invariance under Reidemeister II

$\psi\colon T_1 \to T_2$. It is then easy to see that $\Gamma(\psi) = T$. We shall use the following lemma:

**Lemma 2.1.7** ([14, Lemma 4.5]). *The cone construction is invariant, up to homotopy, under compositions with strong deformation retracts. That is, if $\psi$ and $F$ are composable morphisms of complexes and $F$ is a deformation retract, then $\Gamma(F\psi) \simeq \Gamma(\psi)$.* □

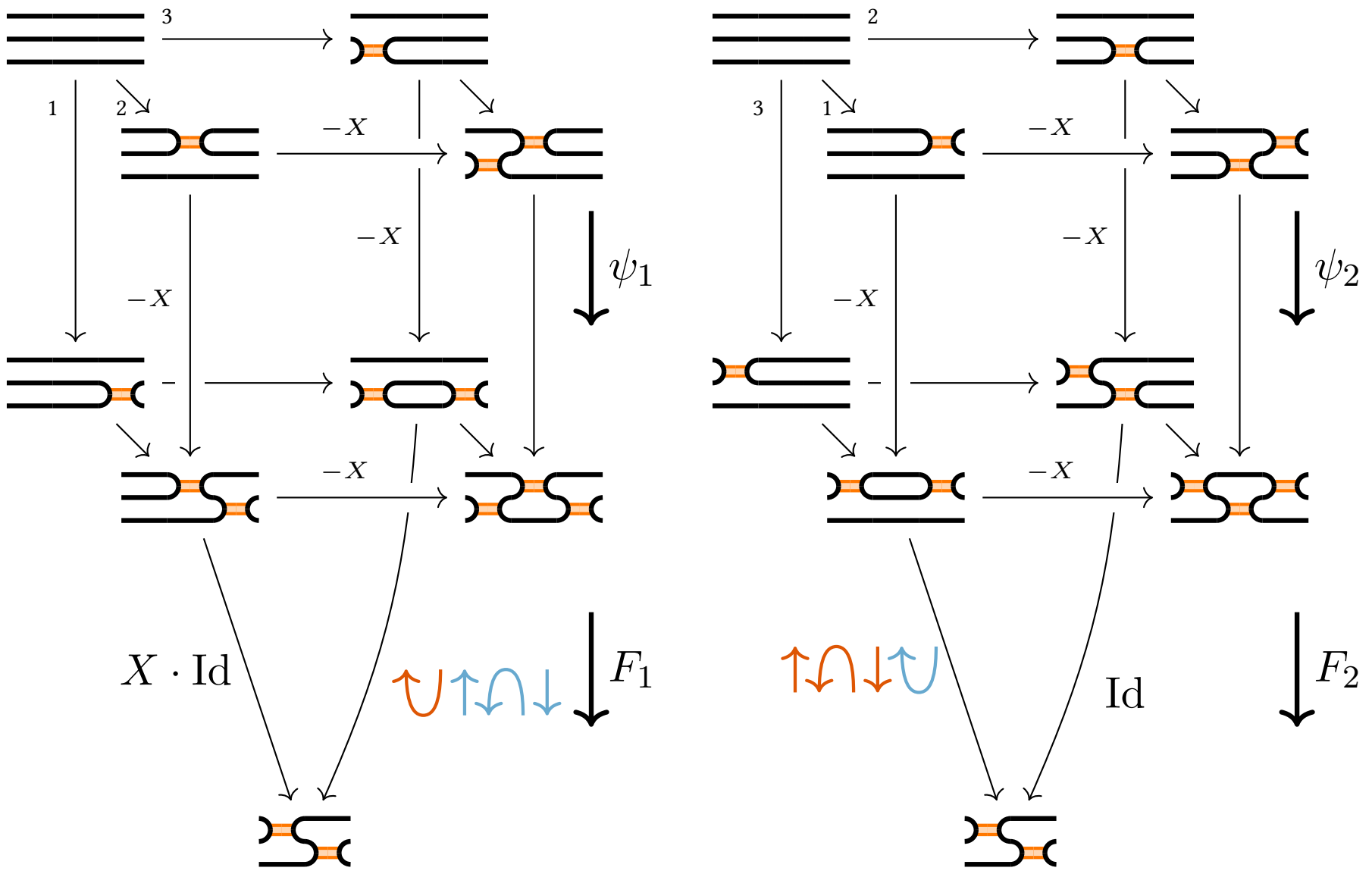

**Fig. 2.6** proof of invariance under Reidemeister III

The proof of invariance is then essentially contained in Fig. 2.6, where we discarded quantum and homological shifts for clarity. On the top, one sees the hypercubes associated to each side of the Reidemeister III move, viewed as cones over morphisms respectively denoted $\psi_1$ and $\psi_2$. The ordering 1, 2 and 3 of the three directions corresponds to the ordering of the crossings, reading from right to left. Note that the top faces are identical. The bottom of the picture shows how the bottom faces are simplified using delooping and gaussian elimination, giving homotopy equivalences $F_1$ and $F_2$.

Thanks to Lemma 2.1.5, $F_1$ and $F_2$ are strong deformation retracts, so by the above lemma, $\Gamma(F_i\psi_i) \simeq \Gamma(\psi_i)$ for $i = 1, 2$. To conclude, it only remains to compare $F_1\psi_1$ and $F_2\psi_2$. Delving into the proof of Lemma 2.1.5 (or computing by hand), we get explicit $F_1$ and $F_2$; the result is shown in Fig. 2.6. We get that $F_1\psi_1 = F_2\psi_2$, and hence $\Gamma(F_1\psi_1) \simeq \Gamma(F_2\psi_2)$. □

*Proof of Theorem 2.0.2.* The renormalization with $q^{-2N_+ + N_-}$ does not affect the result of Lemma 2.1.4: if $D_1$ and $D_2$ are two planar isotopic oriented sliced tangle diagrams, then $q^{-2N_+ + N_-}\mathrm{Kom}_{\mathfrak{gl}_2}(D_1)$ and $q^{-2N_+ + N_-}\mathrm{Kom}_{\mathfrak{gl}_2}(D_2)$ are isomorphic, and in particular homotopic. We conclude with invariance under Reidemeister moves, given by Lemma 2.1.6. □

## 2.2 Review of covering $\mathfrak{sl}_2$-Khovanov homology for links

We review the construction of covering $\mathfrak{sl}_2$-Khovanov homology as defined by Putyra [155]. His construction uses a 2-category of “chronological cobordisms”, but for our purpose, we give here a “low-tech” definition of covering $\mathfrak{sl}_2$-Khovanov homology, directly generalizing the original definition of odd Khovanov homology of Ozsváth, Rasmussen and Szabó [149].

We first give some preliminary definitions. Recall the ring

$$\Bbbk = \mathbb{Z}[X, Y, Z^{\pm}]/(X^2 = Y^2 = 1)$$

from Definition 1.3.2. For $n \in \mathbb{N}$, let $\wedge_{\Bbbk}(a_1, \dots, a_n)$ be the $\Bbbk$-algebra generated by variables $a_1, \dots, a_n$ and subject to the following relations:

$$\begin{aligned} a_i a_j &= XY a_j a_i && \text{for } 1 \leq i, j \leq n, \\ a_i^2 &= 0 && \text{for } 1 \leq i \leq n. \end{aligned}$$

Denote by $\wedge^r_{\Bbbk}(a_1, \dots, a_n)$ the $\Bbbk$-submodule of $\wedge_{\Bbbk}(a_1, \dots, a_n)$ generated by words of length $r$ in the letters $a_1, \dots, a_n$. We endow $\wedge_{\Bbbk}(a_1, \dots, a_n)$ with a $\mathbb{Z}$-grading, the *q-grading*, setting $\operatorname{qdeg} p = 2r - n$ whenever $p \in \wedge^r_{\Bbbk}(a_1, \dots, a_n)$. Define also the following linear maps:

$$\begin{aligned} m_{a_1,a_2;a} \colon \wedge_{\Bbbk}(a_1, a_2, x_1, \dots, x_n) &\to \wedge_{\Bbbk}(a, x_1, \dots, x_n) \\ p &\mapsto p|_{a_1,a_2 \mapsto a} \\ \Delta_{a;a_1,a_2} \colon \wedge_{\Bbbk}(a, x_1, \dots, x_n) &\to \wedge_{\Bbbk}(a_1, a_2, x_1, \dots, x_n) \\ p &\mapsto (a_1 + XYa_2)p|_{a \mapsto a_1} \\ &\quad = (a_1 + XYa_2)p|_{a \mapsto a_2} \end{aligned}$$

Here $a_1, a_2 \mapsto a$ means that one should replace every instance of $a_1$ and $a_2$ by $a$ in $p$, and similarly for $a \mapsto a_1$ and $a \mapsto a_2$. With respect to the $q$-grading, these maps are graded maps with $q$-degree

$$\operatorname{qdeg}(m_{a_1,a_2;a}) = \operatorname{qdeg}(\Delta_{a;a_1,a_2}) = 1.$$

Note that one recovers the algebra $\mathbb{Z}[a_1, \dots, a_n]/(a_1^2 = \dots = a_n^2 = 0)$ with its product and coproduct by setting $X = Y = Z = 1$, and the exterior algebra in variables $a_1, \dots, a_n$ by setting $X = Z = 1$ and $Y = -1$.

Recall Notation 2.1.1. For $\boldsymbol{r} \in \{0,1\}^N$ and $k, l \in \{1, \dots, N\}$ where $k < l$, the square

$$\begin{array}{ccc} \boldsymbol{r} & \longrightarrow & \boldsymbol{r} + \boldsymbol{e}_k \\ \downarrow & \circlearrowleft & \downarrow \\ \boldsymbol{r} + \boldsymbol{e}_l & \longrightarrow & \boldsymbol{r} + \boldsymbol{e}_k + \boldsymbol{e}_l \end{array} \tag{2.3}$$

is given an orientation as depicted, and we denote it by $\square^r_{k,l}$.

Let then $D$ be a link diagram with $N$ crossings. The hypercubic complex $\operatorname{Kom}_{\mathfrak{sl}_2}(D)$ is constructed through the following steps:

(i) *Hypercube of resolutions:* fix an arbitrary order on the crossings of $D$. Each crossing can be *resolved* into two possible planar diagrams, respectively the *0-resolution* (on the left) or the *1-resolution* (on the right):

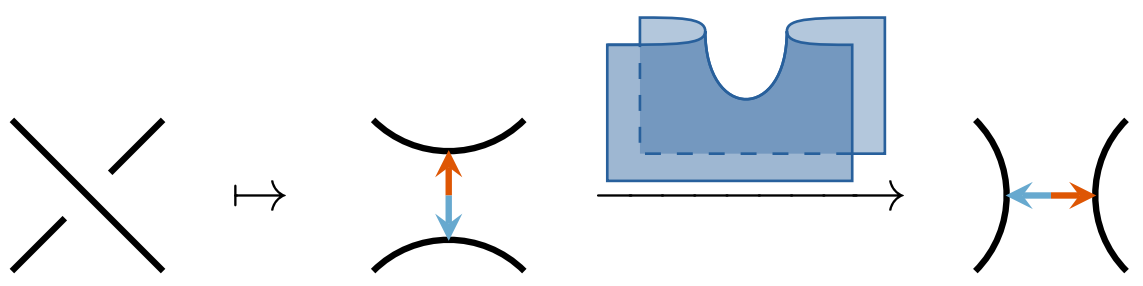

A *resolution* of $D$ is a choice of resolutions for each crossing. The resolutions of $D$ can be pictured as sitting on the vertices of a hypercube $\{0,1\}^N$, where for $\boldsymbol{r} \in \{0,1\}^N$ the binary $\boldsymbol{r}_i$ encodes the chosen resolution for the $i$-th crossing. Each edge $\boldsymbol{r} \to \boldsymbol{r} + \boldsymbol{e}_i$ of the hypercube connects two resolutions that only differ at the $i$-crossing. This edge is decorated with a saddle cobordism, which can either be a merge or a split depending on the global context. Finally, for each crossing one must choose an orientation on the arcs of the two resolutions: the red or the blue orientation. We call it the *arc orientation.* Equivalently, an arc orientation is a choice of arc orientation for the 0-resolution, which induces an arc orientation for the 1-resolution by rotating a quarter of a turn clockwise.

(ii) *Algebrization:* we turn the hypercube of resolutions into a hypercube in the category of $\mathbb{Z}$-graded $\Bbbk$-modules. To each vertex $\boldsymbol{r} \in \{0,1\}^N$ we associate the $\Bbbk$-module

$$V_{\boldsymbol{r}} := q^{N-|\boldsymbol{r}|}t^{-(N-|\boldsymbol{r}|)} \wedge_{\Bbbk} (a_1, \dots, a_n),$$

where $q$ and $t$ denote shifts in quantum and homological degree respectively, and $n$ is the number of connected components in the corresponding resolution. One should think of each variable as attached to one connected component. In addition, each edge is replaced by an $\Bbbk$-linear map between relevant $\Bbbk$-modules:

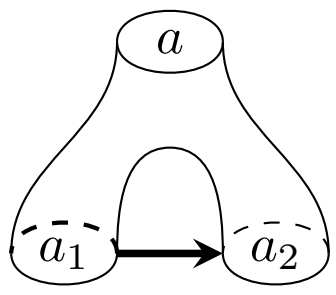

$\mapsto m_{a_1,a_2;a}$ and

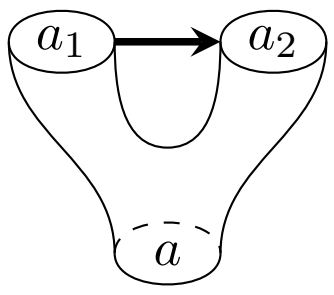

$\mapsto \Delta_{a;a_1,a_2}$

Note the importance of the extra arrows, which give a preferred choice of ordering between the two circles corresponding to the variables $a_1$ and $a_2$. We denote $H_{\mathfrak{sl}_2}(D)$ the resulting hypercube.

(iii) *Commutativity:* As defined, squares in the algebrized hypercube do not necessarily commute. In fact, if we consider a generic square

$$\begin{array}{ccc} \boldsymbol{r} & \xrightarrow{F_{*0}} & \boldsymbol{r}+\boldsymbol{e}_k \\ {\scriptstyle F_{0*}}\downarrow & \circlearrowright & \downarrow{\scriptstyle F_{1*}} \\ \boldsymbol{r}+\boldsymbol{e}_l & \xrightarrow[F_{*1}]{} & \boldsymbol{r}+\boldsymbol{e}_k+\boldsymbol{e}_l \end{array},$$

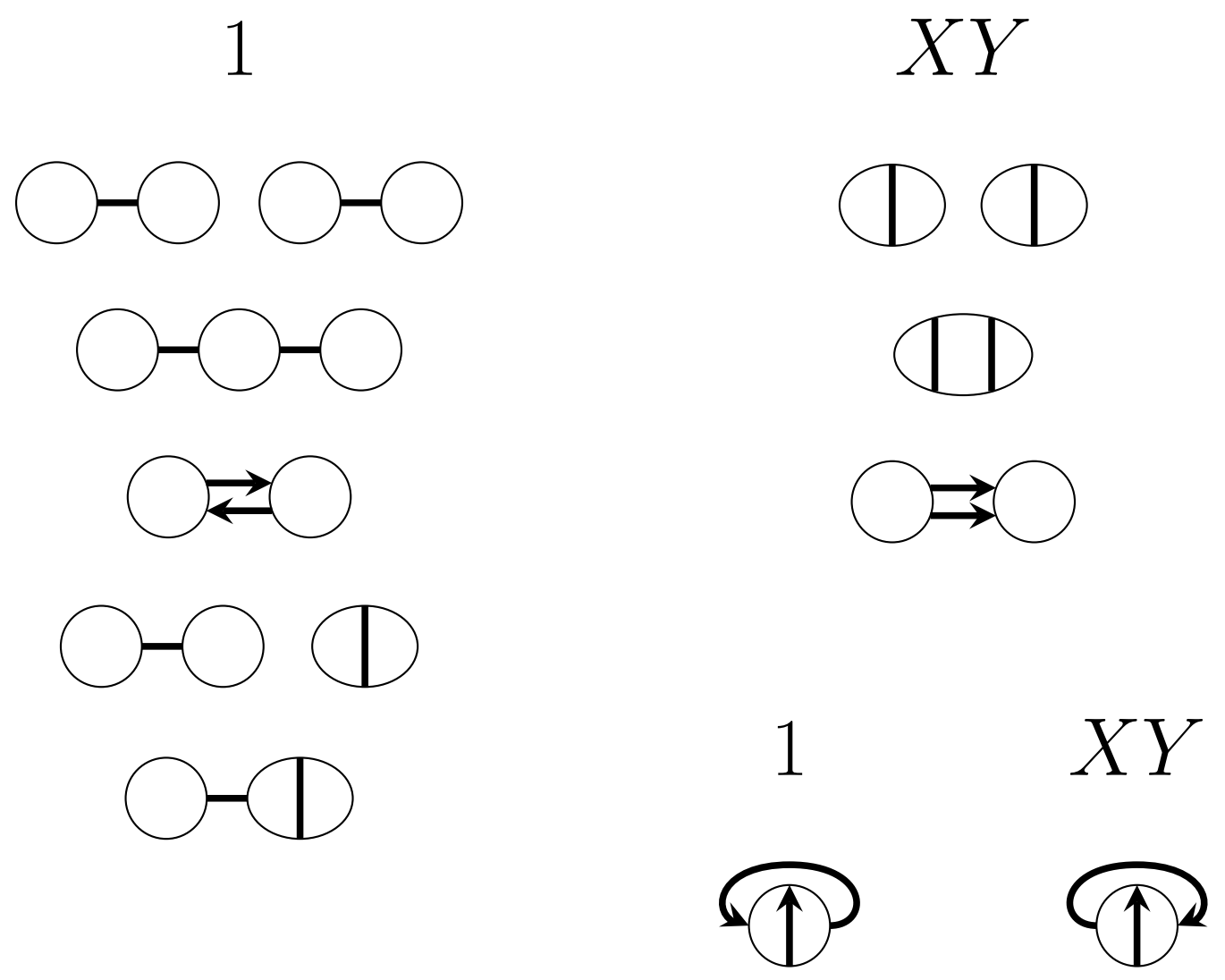

**Table 2.1** Definition of $\psi_{\mathfrak{sl}_2}$ for covering $\mathfrak{sl}_2$-Khovanov homology. Each square is uniquely represented by the (relevant local piece of) resolution at the initial point. If no orientation on the arrows is given, then the value of $\psi_{\mathfrak{sl}_2}$ is independent of the choice of orientations. The last two cases are called the *ladybugs*.

we have:

$$F_{1*} \circ F_{*0} = \psi_{\mathfrak{sl}_2}(\square^r_{k,l}) F_{*1} \circ F_{0*}.$$

Here $\psi_{\mathfrak{sl}_2}$ is the $\Bbbk^\times$-valued 2-cochain on the hypercube defined by Table 2.1. As shown in [149, 155], $\psi_{\mathfrak{sl}_2}$ is a cocycle. An $\mathfrak{sl}_2$-*scalar assignment* is a choice of a 1-cochain $\epsilon_{\mathfrak{sl}_2}$ such that $\partial\epsilon_{\mathfrak{sl}_2} = \psi_{\mathfrak{sl}_2}$. Such a choice always exists: by contractibility of the hypercube, a 2-cocycle is always a 2-coboundary. Given a choice of $\mathfrak{sl}_2$-scalar assignment $\epsilon_{\mathfrak{sl}_2}$, we multiply each edge $e$ of the hypercube by $\epsilon_{\mathfrak{sl}_2}(e)$: this makes each square commute. We denote $H_{\mathfrak{sl}_2}(D, \epsilon_{\mathfrak{sl}_2})$ the resulting hypercube.

(iv) *Koszul rule:* we apply the Koszul rule (multiplying each edge $\boldsymbol{r} \to \boldsymbol{r} + \boldsymbol{e}_i$ by $(-1)^{\#\{r_j=1 \mid j<i\}}$) to turn every commutative square into an anti-commutative square. We denote $\mathrm{Kom}_{\mathfrak{sl}_2}(D)$ the resulting hypercubic complex.

**Theorem 2.2.1** ([149, 155]). *Let $D$ be an oriented link diagram with respectively $N_+$ and $N_-$ positive and negative crossings, presenting an oriented link $L$. The isomorphism class of* $\mathrm{Kom}_{\mathfrak{sl}_2}(D)$ *is independent of the choice of ordering on crossings, the choice of arc orientations, and the choice of* $\mathfrak{sl}_2$*-scalar assignment.*

*Moreover, the homotopy type of* $q^{-2N_+ + N_-} t^{N_+} \mathrm{Kom}_{\mathfrak{sl}_2}(D)$*, denoted* $\mathrm{CKh}_{\mathfrak{sl}_2}(L)$*, is an invariant of* $L$*.*

*Remark* 2.2.2. The 2-cocycle $\psi_{\mathfrak{sl}_2}$ is not the only choice that makes the construction above work. Indeed, for the last two cases of Table 2.1, called the *ladybugs*,[1] We both have

$$F_{1*} \circ F_{*0} = \psi_{\mathfrak{sl}_2}(\square^r_{k,l}) F_{*1} \circ F_{0*} \text{ and } F_{1*} \circ F_{*0} = XY \psi_{\mathfrak{sl}_2}(\square^r_{k,l}) F_{*1} \circ F_{0*}$$

Define $\overline{\psi}_{\mathfrak{sl}_2}$ to be the 2-cochain defined as $\psi_{\mathfrak{sl}_2}$ except for the ladybugs, where instead we set

$$\overline{\psi}_{\mathfrak{sl}_2}\left( \right) = XY \qquad \text{and} \qquad \overline{\psi}_{\mathfrak{sl}_2}\left( \right) = 1.$$

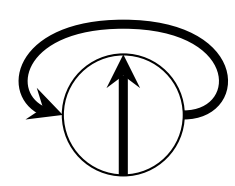

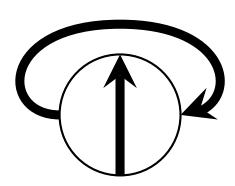

The results of Theorem 2.2.1 still hold in this case. Let us use the notations $\mathrm{CKh}_{\mathfrak{sl}_2}(L, \psi_{\mathfrak{sl}_2})$ and $\mathrm{CKh}_{\mathfrak{sl}_2}(L, \overline{\psi}_{\mathfrak{sl}_2})$ to distinguish the two constructions. In [149], they are respectively called "type X" and "type Y" (setting $X = Z = 1$ and $Y = -1$; no analogy between the scalar $X$ and "type X" intended). It is shown in Putyra [155] that $\mathrm{CKh}_{\mathfrak{sl}_2}(L, \psi_{\mathfrak{sl}_2})$ and $\mathrm{CKh}_{\mathfrak{sl}_2}(L, \overline{\psi}_{\mathfrak{sl}_2})$ are in fact isomorphic.

## 2.3 Covering $\mathfrak{sl}_2$- and $\mathfrak{gl}_2$-Khovanov homologies are isomorphic

This section is devoted to the proof of Theorem 2.0.3. Here is a quick summary:

(i) In subsection 2.3.2, we restate the definition of covering $\mathfrak{gl}_2$-Khovanov homology using a $\mathfrak{gl}_2$*-hypercube of resolutions*. The rest of the proof consists in comparing this $\mathfrak{gl}_2$-hypercube with the $\mathfrak{sl}_2$-hypercube defined above.

(ii) To compare the hypercubes, we need to compare the $\Bbbk$-modules at each vertex. This requires a choice of basis for each $\mathrm{Hom}_{\mathbf{GFoam}}(\emptyset, W)$, called *cup foams*, that we describe in subsection 2.3.1.

(iii) In subsection 2.3.3, we use the above basis to define a family of isomorphisms on the level of vertices. This defines a proper morphism of hypercubes only *up to invertible scalar*; we call it a *projective morphism*. We state a certain 2-cocycle condition such that, if satisfied, the aforementioned family of isomorphisms can be rescaled into a genuine isomorphism of hypercubes.

[1]This terminology is borrowed from [121].

(iv) The proof of Theorem 2.0.3 then reduces to the analysis of this 2-cocycle condition. Subsection 2.3.4 shows that this can be done locally, looking only at the cases pictured in Table 2.1. In most cases, general considerations show that the 2-cocycle condition is necessarily verified.

(v) However, these general considerations do not work for the ladybugs (see Remark 2.2.2). To deal with these two cases, we require finer results on the independence on all choices involved in the above family of isomorphisms. This is done in subsection 2.3.5, which concludes the proof.

### 2.3.1 Cup foams

Call a web $W$ *closed* if $\mathfrak{sl}(W)$ is a closed 1-manifold. Recall the basis described in section 1.6. In the special case where the domain is the empty web and the codomain is a closed web $W$ (recall that in $\underline{(\mathbf{GFoam})}|_q^{\oplus}$, the empty web is equal to a juxtaposition of double lines), an undotted reduced foam as the follow form:

$\mathfrak{sl}(F) =$ 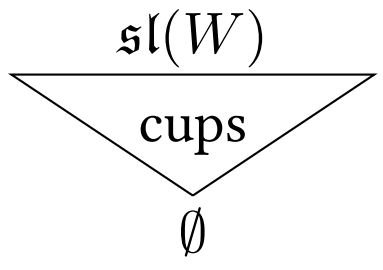

We call such an $F$ an *undotted cup foam on* $W$. If we wish to allow $F$ to (possibly) carry dots, we simply say that $F$ is a *cup foam on* $W$. We write $\pi_0(\mathfrak{sl}(W))$ for the set of closed components of $\mathfrak{sl}(W)$. In this context, Theorem 1.6.1 is restated as follows:

**Proposition 2.3.1.** *Let $W$ be a closed web. Let $B$ be a set containing precisely one cup foam*

$$\beta_\delta \colon \emptyset \to W$$

*for each subset $\delta \subset \pi_0(\mathfrak{sl}(W))$, so that for each $c \in \pi_0(\mathfrak{sl}(W))$, the corresponding disk in $\mathfrak{sl}(\beta_\delta)$ is dotted if and only if $c \in \delta$. Then $B$ is basis for the $\Bbbk$-module* $\mathrm{Hom}_{\mathbf{GFoam}}(\emptyset, W)$. □

Fix an undotted cup foam $\beta^W$ for $W$ and pick a total order on $\pi_0(\mathfrak{sl}(W))$. For each subset $\delta \subset \mathfrak{sl}(W)$, denote by $\mathrm{id}_W^\delta \colon W \to W$ the foam identical to $\mathrm{id}_W$ except for an additional dot on the connected component $c$ for each $c \in \delta$, ordering the dots increasingly with respect to the total order on $\pi_0(\mathfrak{sl}(W))$, reading from bottom to top. This defines $\mathrm{id}_W^\delta$ uniquely. We denote $\beta_\delta^W := \mathrm{id}_W^\delta \circ \beta^W$. Schematically:

$\mathfrak{sl}(\beta_\delta^W) =$ 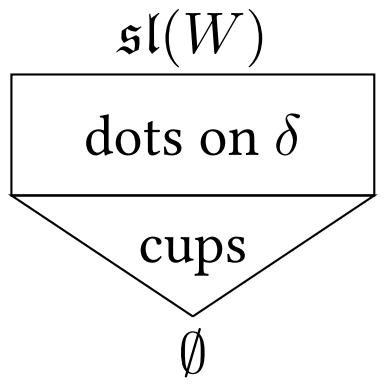

It follows from Proposition 2.3.1 that:

**Corollary 2.3.2.** *For every choices of undotted cup foam $\beta^W$ and total order on $\pi_0(\mathfrak{sl}(W))$, the family $\{\beta_\delta^W\}_{\delta \subset \mathfrak{sl}(W)}$ defines a basis for the $\Bbbk$-module* $\mathrm{Hom}_{\mathbf{GFoam}}(\emptyset, W)$.

Let $V$ and $W$ be two closed webs. If $F\colon V \to W$ is a zip or an unzip, then $\mathfrak{sl}(F)$ is either a merge or a split. Below we sometimes speak of the closed components in the domain and codomain of $\mathfrak{sl}(F)$ to refer only to the closed components making up the boundary of the closed component in $\mathfrak{sl}(F)$ containing the saddle. The distinction should be clear by the context.

Recall the symbol $\dot{\sim}$ to mean "equal up to multiplying by an invertible scalar", defined in the paragraph before Lemma 1.3.6.

**Proposition 2.3.3.** *Let $V$ and $W$ be closed webs and $\beta^V$ (resp. $\beta^W$) be a choice of undotted cup foam for $V$ (resp. $W$). Let $F\colon V \to W$ be a zip or an unzip. Then:*

*(i) If $\mathfrak{sl}(F)$ is a merge, then (as depicted in* (2.4)*)*

$$F \circ \beta^V \mathrel{\dot{\sim}} \beta^W.$$

*(ii) If $\mathfrak{sl}(F)$ is split, then (as depicted in* (2.5)*)*

$$F \circ \beta^V \mathrel{\dot{\sim}} \beta_{i_1}^W + XY\beta_{i_2}^W,$$

*where $i_1$ and $i_2$ are the connected components of $W$ corresponding to the codomain of $\mathfrak{sl}(F)$ (note that the ordering of $i_1$ and $i_2$ is irrelevant in that statement).*

Here is the schematic for Proposition 2.3.3:

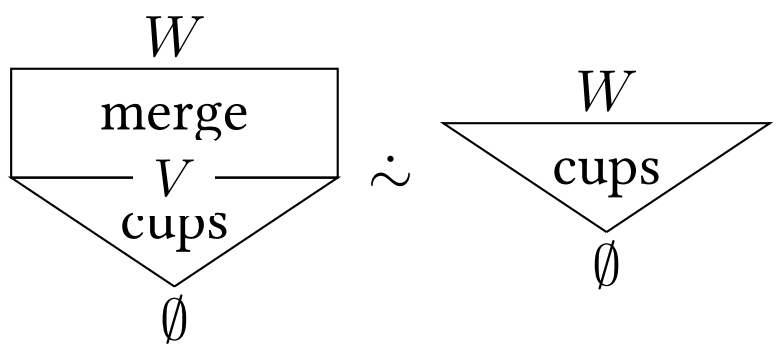

(2.4)

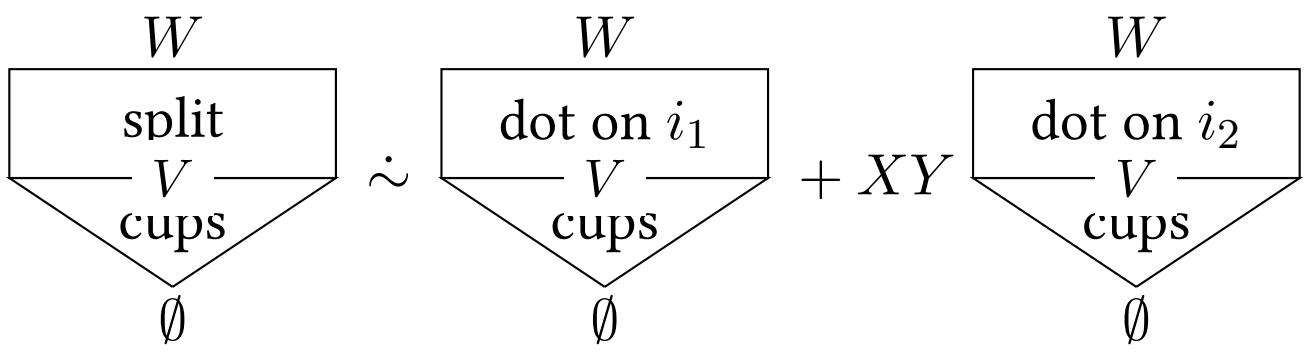

(2.5)

To prove the proposition, we will need the following lemma:

**Lemma 2.3.4.** *Let $F\colon W \to W'$ and $F'\colon W' \to W''$ be two foams and let $D_1$ and $D_2$ be two foams identical to $\mathrm{id}_{W'}$ except for a dot sitting on a 1-facet. If the two dots belong to the same closed component of $\mathfrak{sl}(F' \circ F)$, then*

$$F' \circ D_1 \circ F = F' \circ D_2 \circ F$$

*in* **GFoam**.

*Proof.* This is a consequence of the fact that a dot can slide across 2-facets at no cost of scalar, and along 1-facets past generators depending on $\mu$. Because $\mu$ is symmetric (see Definition 1.1.7), the latter scalar only depends on the relative vertical position of the dot. □

*Proof of Proposition 2.3.3.* Fix a total order on $\pi_0(\mathfrak{sl}(V))$ and $\pi_0(\mathfrak{sl}(W))$. Using the neck-cutting successively on the identity of $W$, one can decompose it as

$$\mathrm{id}_W = \sum_{\delta \subset \pi_0(\mathfrak{sl}(W))} \beta_\delta^W \circ \beta_{\delta^c}^c,$$

where $\delta^c := \pi_0(\mathfrak{sl}(W)) \setminus \delta$, and each $\beta_\delta^c\colon W \to \emptyset$ is a *cap foam* in the sense that $\mathfrak{sl}(\beta_\delta^c)$ is a union of disks, each disk being dotted depending on $\delta$ as in Proposition 2.3.1. This allows us to write:

$$F \circ \beta^V = \sum_{\delta \subset \pi_0(\mathfrak{sl}(W))} \beta_\delta^W \circ \left(\beta_{\delta^c}^c \circ F \circ \beta^V\right),$$

where each $\beta^{c}_{\delta^c} \circ F \circ \beta^V$ is a *closed foam*, that is a foam with domain and codomain the empty web. Then, we apply the following result on the evaluation of closed foams, which is a consequence of Theorem 1.6.1:

**Lemma 2.3.5.** *Let* $U\colon \emptyset \to \emptyset$ *be a closed foam in* **GFoam**. *Then:*

$$U \mathrel{\dot\sim} \begin{cases} \mathrm{id}_\emptyset & \text{if each closed component of } \mathfrak{sl}(U) \text{ is a sphere with a single dot,} \\ 0 & \text{otherwise.} \end{cases}$$

□

By the above lemma, there exist invertible scalars $\tau$, $\tau_1$ and $\tau_2$ such that:

$$F \circ \beta^V = \tau\beta^W \quad \text{or} \quad F \circ \beta^V = \tau_1\beta^W_{i_1} + \tau_2\beta^W_{i_2},$$

depending on whether $\mathfrak{sl}(F)$ is a merge or a split. It remains to show that $\tau_1/\tau_2 = XY$ in the latter case. For that, we use Lemma 2.3.4. Assume $i_1 < i_2$ for the purpose of the computation, so that $\mathrm{id}_{i_1,i_2} = \mathrm{id}_{i_2} \circ \mathrm{id}_{i_1}$:

$$\begin{aligned} \tau_1\beta^W_{i_1,i_2} &= \mathrm{id}^W_{i_2} \circ (\tau_1\beta^W_{i_1} + \tau_2\beta^W_{i_2}) = \mathrm{id}^W_{i_2} \circ F \circ \beta^V \\ &\overset{2.3.4}{=} \mathrm{id}^W_{i_1} \circ F \circ \beta^V = \mathrm{id}^W_{i_1} \circ (\tau_1\beta^W_{i_1} + \tau_2\beta^W_{i_2}) = XY\tau_2\beta^W_{i_1,i_2}, \end{aligned}$$

where in the last equality we used that two dots interchange at the cost of the scalar $XY$. Since $\beta^W_{i_1,i_2} \neq 0$ belongs to a free family by Proposition 2.3.1, we must have $(\tau_1 - XY\tau_2) = 0$, which concludes. □

### 2.3.2 The $\mathfrak{gl}_2$-hypercube of resolutions

We reformulate the definition of covering $\mathfrak{gl}_2$-Khovanov homology to emphasize the similarities with covering $\mathfrak{sl}_2$-Khovanov homology. Recall

$$\mathcal{A}_{\mathfrak{gl}_2} := \mathrm{Hom}_{\mathbf{GFoam}(\emptyset,\emptyset)}(\emptyset, -)$$

defined in the introduction of this chapter.

Let $D$ be a sliced oriented link diagram with $N$ crossings. As described in subsection 2.1.2, we can associate with $D$ a knotted web $W_D$. Then, starting with $W_D$, the complex $\mathcal{A}_{\mathfrak{gl}_2}(\mathrm{Kom}_{\mathfrak{gl}_2}(D))$ can be defined as follows:

(i) *Hypercube of resolutions:* each crossing can be resolved into a *web 0-resolution* or a *web 1-resolution*:

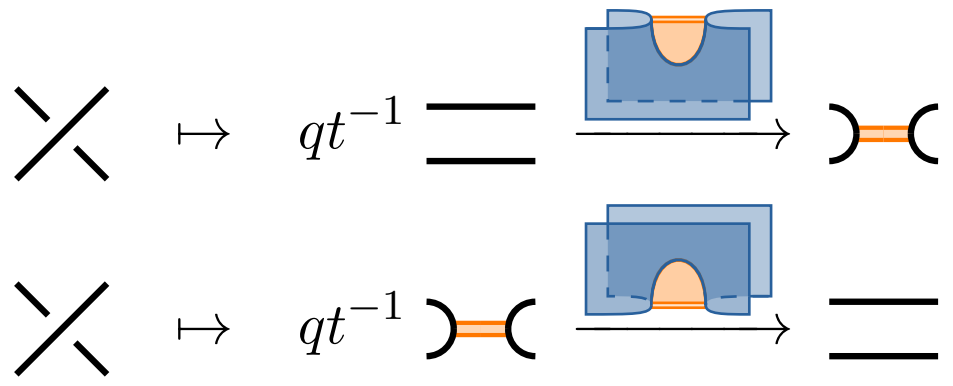

A *web resolution* is a choice of web resolutions for each crossing. Fixing an ordering on the crossings, they can be pictured as sitting on the vertices of the hypercube $\{0,1\}^N$, whose edges are decorated with a zip or an unzip, depending on whether the associated crossing is positive or negative.

(ii) *Algebrization:* we apply the functor $\mathcal{A}_{\mathfrak{gl}_2}$ to the hypercube. Denote $H_{\mathfrak{gl}_2}(D)$ the decorated hypercube so obtained.

(iii) *Commutativity:* a $\mathfrak{gl}_2$*-scalar assignment*[2] is an $\Bbbk^\times$-valued 1-cochain $\epsilon_{\mathfrak{gl}_2}$ on the hypercube $\{0,1\}^N$, such that $\partial\epsilon_{\mathfrak{gl}_2} = \psi_{\mathfrak{gl}_2}$ where $\psi_{\mathfrak{gl}_2}$ is a 2-cocycle defined as

$$\psi_{\mathfrak{gl}_2}(\square^r_{k,l}) := \mu(\deg F_{0*}, \deg F_{*0})^{-1} = \mu(\deg F_{*1}, \deg F_{1*}).$$

We multiply each edge $e$ by $\epsilon_{\mathfrak{gl}_2}(e)$. This makes each square commutes. This defines a hypercube $H_{\mathfrak{gl}_2}(D;\epsilon_{\mathfrak{gl}_2})$.

(iv) *Koszul rule:* we apply the Koszul rule to turn every commutative square into an anti-commutative square.

As a consequence of Lemma 4.1.7, the isomorphism class of the complex obtained does not depend on the choice of $\mathfrak{gl}_2$-scalar assignment. It is easily checked that the construction above coincides with the definition of covering $\mathcal{A}_{\mathfrak{gl}_2}(\mathrm{Kom}_{\mathfrak{gl}_2}(D))$ given in subsection 2.1.1. In particular, it is shown in Definition 4.1.6 that the graded Koszul rule is a $\mathfrak{gl}_2$-scalar assignment.[3]

### 2.3.3 A projective isomorphism of hypercubes

Let then $D$ be a sliced link diagram with $N$ crossings, with a fixed choice of ordering on crossings and arc orientations. To show Theorem 2.0.3, it suffices to exhibit an isomorphism between the two hypercubes $H_{\mathfrak{sl}_2}(D;\epsilon_{\mathfrak{sl}_2})$ and $H_{\mathfrak{gl}_2}(D;\epsilon_{\mathfrak{gl}_2})$, where an *isomorphism of hypercubes* is a family of isomorphisms at each vertex, such that all squares involved commute. We abuse notation and denote $H_{\mathfrak{sl}_2}(D)$ (resp. $H_{\mathfrak{gl}_2}(D)$) both the hypercube of resolutions (step (i)) and its algebrization (step (ii)): the relevant hypercube should be clear by the context. We also use $\mathfrak{sl}_2$ and $\mathfrak{gl}_2$ as subscripts to distinguish features of the two constructions. For instance, we write $\boldsymbol{r}_{\mathfrak{gl}_2}$ to denote the decoration on the vertex $\boldsymbol{r}\in\{0,1\}^N$ of the hypercube $H_{\mathfrak{gl}_2}(D)$ (that is, a web $W$, or the graded $\Bbbk$-module $\mathrm{Hom}_{\mathbf{GFoam}}(\emptyset,W)$ depending on the context).

[2]The notion of scalar assignment here is slightly different from chapter 4, as the latter already includes the Kozsul rule.

[3]With the caveat of the previous footnote.

Looking at step (i) in the construction of the hypercubes, it is clear that $\mathfrak{sl}(H_{\mathfrak{gl}_2}(D))=H_{\mathfrak{sl}_2}(D)$: that is, $\mathfrak{sl}(\boldsymbol{r}_{\mathfrak{gl}_2})=\boldsymbol{r}_{\mathfrak{sl}_2}$ for each $\boldsymbol{r}\in\{0,1\}^N$, and similarly for edges. For each vertex $W=\boldsymbol{r}_{\mathfrak{gl}_2}$ of $H_{\mathfrak{gl}_2}(D)$, fix a choice of undotted cup foam $\beta^W$ and a choice of total ordering on the set of connected components $\pi_0(\mathfrak{sl}(W))$. By Corollary 2.3.2, $\boldsymbol{r}_{\mathfrak{gl}_2}$ has basis given by the set $\{\beta^W_\delta\}_{\delta\subset\pi_0(\mathfrak{sl}(W))}$. On the other hand, $\boldsymbol{r}_{\mathfrak{sl}_2}$ has basis given by the set $\{a_\delta\}_{\delta\subset\pi_0(\mathfrak{sl}(W))}$, where $a_\delta=a_{i_k}\ldots a_{i_1}$ with $\delta=\{i_k,\ldots,i_1\}$ and $i_k>\ldots>i_1$. Hence, the map

$$\iota_W\colon \beta^W_\delta\mapsto \mu(|\delta|\,(1,1),\beta^W)a_\delta$$

defines an isomorphism of graded $\Bbbk$-modules.[4] Note that $\iota_W$ depends on the choice of $\beta^W$, but not on the choice of ordering on $\pi_0(\mathfrak{sl}(W))$. Moreover, by Proposition 2.3.3 the family of those isomorphisms at each vertex defines a *projective* isomorphism of hypercubes, in the sense that for each edge $e\colon V\to W$ in $H_{\mathfrak{gl}_2}(D)$, the square

$$\square_e := \begin{array}{ccc} V & \xrightarrow{\;e\;} & W \\ {\scriptstyle\iota_V}\downarrow & \circlearrowleft & \downarrow{\scriptstyle\iota_W} \\ \mathfrak{sl}(V) & \xrightarrow[\mathfrak{sl}(e)]{} & \mathfrak{sl}(W) \end{array} \qquad \begin{array}{c} H_{\mathfrak{gl}_2}(D) \\ \\ H_{\mathfrak{sl}_2}(D) \end{array} \tag{2.6}$$

commutes up to invertible scalar (the symbol $\circlearrowleft$ denotes a choice of orientation). This scalar is computed in the following lemma:

**Lemma 2.3.6.** *Denote by $\tau\in\Bbbk^\times$ the invertible scalar given by Proposition 2.3.3. Then either*

$$e\circ\beta^V=\tau\beta^W \quad \textit{or} \quad e\circ\beta^V=\tau\left(\beta^W_{i_1}+XY\beta^W_{i_2}\right),$$

[4]One can think of $a_\delta$ as corresponding to $\beta^W\circ\mathrm{id}^W_\delta$: this makes no sense when it comes to composition, but this is coherent with the interchange relation. Thinking this way makes some of the identities below clearer.

*depending on whether* $\mathfrak{sl}(e)$ *is a merge or a split. We distinguish* $i_1$ *from* $i_2$ *using the arc orientation: it goes from* $i_1$ *to* $i_2$*. Define*

$$\psi_\beta(e) := \begin{cases} \tau & \text{if } \mathfrak{sl}(e) \text{ is a merge,} \\ \tau\mu((1,1),\beta^W) & \text{if } \mathfrak{sl}(e) \text{ is a split.} \end{cases}$$

*Then* $\iota_W \circ e = \psi_\beta(e)(\mathfrak{sl}(e) \circ \iota_V)$*.*

Note that $\psi_\beta(e)$ depends in general on the choices of arc orientation on $e$ and undotted cup foams $\beta^V$ and $\beta^W$, but not on the choice of ordering on closed components. We postpone the proof of Lemma 2.3.6 until after this discussion.

Ultimately, we are interested in comparing the hypercubes $H_{\mathfrak{sl}_2}(D;\epsilon_{\mathfrak{sl}_2})$ and $H_{\mathfrak{gl}_2}(D;\epsilon_{\mathfrak{gl}_2})$. The above suggests the strategy of finding a 0-cochain $\varphi$ on $\{0,1\}^N$ such that the square

$$\begin{array}{ccccc} \mathfrak{sl}(V) & \xrightarrow{\epsilon_{\mathfrak{gl}_2}(e)e} & W & \qquad & H_{\mathfrak{gl}_2}(D;\epsilon_{\mathfrak{gl}_2}) \\ {\scriptstyle \varphi(V)\iota_V}\big\downarrow & \circlearrowleft & \big\downarrow{\scriptstyle \varphi(W)\iota_W} & & \\ \mathfrak{sl}(V) & \xrightarrow[\epsilon_{\mathfrak{sl}_2}(e)\mathfrak{sl}(e)]{} & \mathfrak{sl}(W) & & H_{\mathfrak{sl}_2}(D;\epsilon_{\mathfrak{sl}_2}) \end{array}$$

commutes. That would define an isomorphism of hypercubes, and prove Theorem 2.0.3.

We can rephrase the problem as follows. Denote by $H_\iota(D)$ the $(N+1)$-dimensional hypercube decorated as $H_{\mathfrak{gl}_2}(D)$ on $\{0,1\}^N \times \{0\}$, as $H_{\mathfrak{gl}_2}(D)$ on $\{0,1\}^N \times \{1\}$, and decorated with $\iota$ on edges $\boldsymbol{r} \times \{0\} \to \boldsymbol{r} \times \{1\}$. Define also the following 2-cochain $\psi$ on $H_\iota(D)$:

$$\psi := \begin{cases} \psi_{\mathfrak{gl}_2} & \text{on } \{0,1\}^N \times \{0\}, \\ \psi_{\mathfrak{sl}_2} & \text{on } \{0,1\}^N \times \{1\}, \\ \psi_\beta(e)^{-1} & \text{on } \square_e. \end{cases}$$

where we recall $\square_e$ from (2.6).

Assume that $\psi$ is a 2-cocycle. Then by contractibility of the hypercube, it is a coboundary, and there exists some 1-cochain $\epsilon$ such that $\partial\epsilon = \psi$. Denote by $H_\iota(D;\epsilon)$ the hypercube $H_\iota(D)$ obtained by multiplying each edge by its value on $\epsilon$. By definition, $\epsilon_{\mathfrak{gl}_2} := \epsilon|_{\{0,1\}^N\times\{0\}}$ (resp. $\epsilon_{\mathfrak{sl}_2} := \epsilon|_{\{0,1\}^N\times\{1\}}$) is a $\mathfrak{gl}_2$-scalar assignment (resp. a $\mathfrak{sl}_2$-scalar assignment). In other words, the hypercube $\{0,1\}^N \times \{0\}$ (resp. $\{0,1\}^N \times \{1\}$) in $H_\iota(D;\epsilon)$ coincides with $H_{\mathfrak{gl}_2}(D;\epsilon_{\mathfrak{gl}_2})$ (resp $H_{\mathfrak{sl}_2}(D;\epsilon_{\mathfrak{sl}_2})$). Moreover, by Lemma 2.3.6 all squares of the kind $\square_e$ commute, so that the $(N+1)$-direction in the hypercube $H_\iota(D;\epsilon)$ defines an isomorphism of hypercubes between $H_{\mathfrak{gl}_2}(D;\epsilon_{\mathfrak{gl}_2})$ and $H_{\mathfrak{sl}_2}(D;\epsilon_{\mathfrak{sl}_2})$.

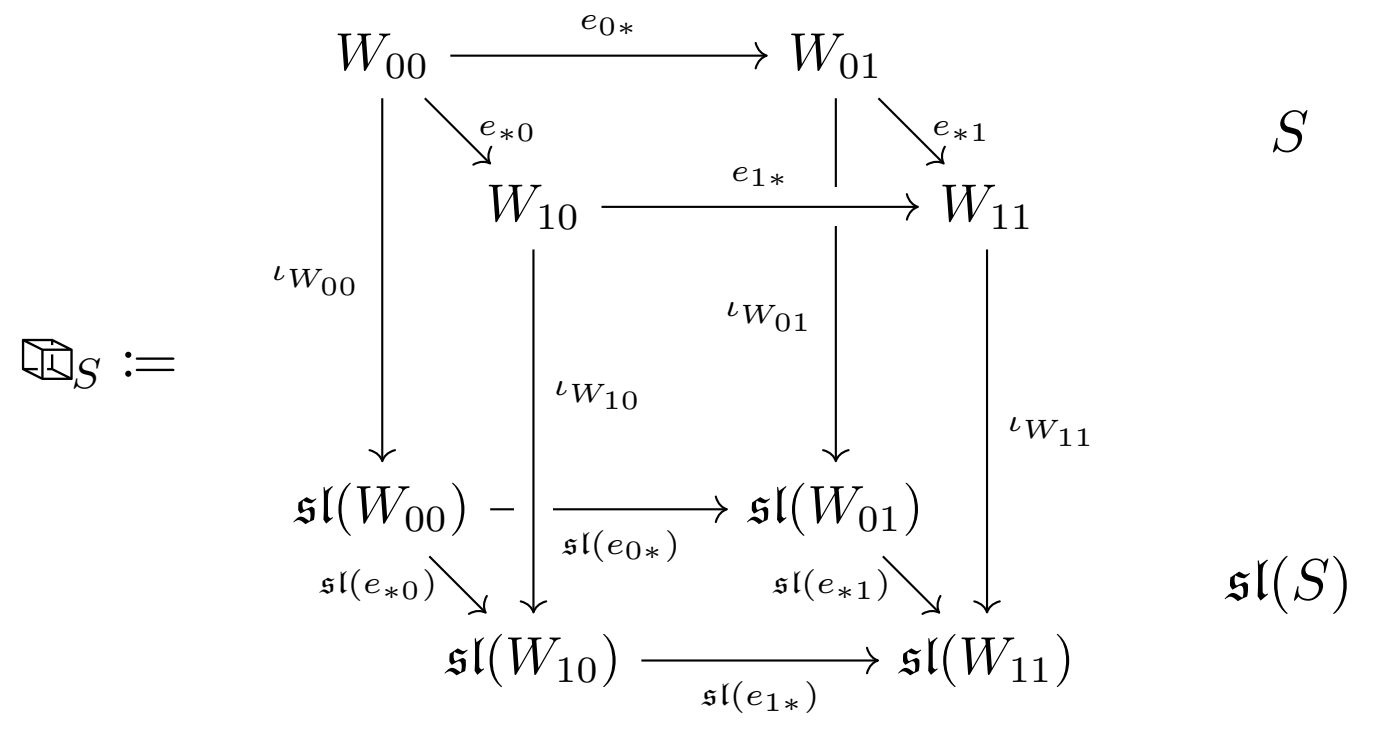

**Fig. 2.7** The 3-dimensional cube $⧉_S$, where $S$ is a square of $H_{\mathfrak{gl}_2}(D)$, pictured on the top.

When is $\psi$ a 2-cocycle? Note that it suffices that $\psi$ is a 2-cocycle on every 3-dimensional cube of $H_\iota(D)$. As $\psi_{\mathfrak{gl}_2}$ and $\psi_{\mathfrak{sl}_2}$ are already 2-cocycles, it is only necessary that $\psi$ is a 2-cocycle on every 3-dimensional cube $⧉_S := S\times\{0,1\}$ for $S$ a square of $H_{\mathfrak{gl}_2}(D)$ (see Fig. 2.7). We give an orientation on $⧉_S$ such that it agrees with the orientation of $e_{0*}$.

To sum up, we have shown that:

**Proposition 2.3.7.** *Let $D$ be a sliced link diagram with $N$ crossings. Assume given an ordering on crossings, a choice of arc orientations, and a choice of undotted cup foams on the webs decorating $H_{\mathfrak{gl}_2}(D)$. If for every square $S$ of $H_{\mathfrak{gl}_2}(D)$, the identity*

$$\partial\psi(⧉_S) = 1$$

*holds, then there exist scalar assignments $\epsilon_{\mathfrak{gl}_2}$ and $\epsilon_{\mathfrak{sl}_2}$ such that $H_{\mathfrak{gl}_2}(D;\epsilon_{\mathfrak{gl}_2})$ and $H_{\mathfrak{sl}_2}(D;\epsilon_{\mathfrak{sl}_2})$ are isomorphic.* □

We end this subsection with the proof of Lemma 2.3.6.

*Proof of Lemma 2.3.6.* We have two cases: either $\mathfrak{sl}(e)$ is a merge, or it is a split. In both cases, we compare where the basis element $\beta_\delta^W$ is mapped through the two paths defining the square. The first case is depicted below, where the two end-results are separated by a dashed line:

$$
\begin{array}{ccc}
\beta^V_\delta & \longrightarrow & \tau\mu(\deg e, |\delta|\,(1,1))\ \beta^W_\delta \\
\Big\downarrow & & \Big\downarrow \\
 & & \tau\mu(\deg e, |\delta|\,(1,1))\mu(|\delta|\,(1,1), \deg\beta^W)\ a^W_\delta \\
 & & \text{- - - - - - - - - - - - - - - - - -} \\
\mu(|\delta|\,(1,1), \deg\beta^V)\ a^V_\delta & \longrightarrow & \mu(|\delta|\,(1,1), \deg\beta^V)\ a^W_\delta
\end{array}
$$

We use symmetry of $\mu$ and the identity $\deg\beta^W = \deg\beta^V + \deg e$ to conclude. Similarly, the computation for the second case gives:

$$
\begin{array}{ccc}
\beta^V_\delta & \longrightarrow & \tau\mu(\deg e, |\delta|\,(1,1))\ \mathrm{id}^\delta_W \circ \left(\beta^W_{i_1} + XY\beta^W_{i_2}\right) \\
\Big\downarrow & & \Big\downarrow \\
 & & \tau\mu(\deg e, |\delta|\,(1,1))\mu((|\delta|+1)(1,1), \deg\beta^W)\ a^W_\delta \left(a^W_{i_1} + XYa^W_{i_2}\right) \\
 & & \text{- - - - - - - - - - - - - - - - - -} \\
\mu(|\delta|\,(1,1), \deg\beta^V)\ a^V_\delta & \longrightarrow & \mu(|\delta|\,(1,1), \deg\beta^V)\ \left(a^W_{i_1} + XYa^W_{i_2}\right) a^W_\delta
\end{array}
$$

and the identity $\deg\beta^W + (1,1) = \deg\beta^V + \deg e$ concludes. □

### 2.3.4 Local analysis

For a sliced link diagram $D$ together with choices as in Proposition 2.3.7, we need to verify

$$\partial\psi(\text{⧉}_S) = 1$$

for every square $S$ in $H_{\mathfrak{gl}_2}(D)$. Note that this certainly does not depend on the choice of ordering on crossings. We implicitly assume such a choice in the sequel.

Note also that given a square $S$ in $H_{\mathfrak{gl}_2}(D)$, the value of $\partial\psi(\text{⧉}_S)$ only depends on *local* choices: a choice of arc orientations for the two crossings involved and choices of undotted cup foams for the four webs involved. (This is essentially because $\iota$ only depends on local choices.) In other words, we can restrict our analysis to the set of squares that can appear in a $\mathfrak{gl}_2$-hypercube of resolutions, that is, to diagrams with exactly two crossings:

**Lemma 2.3.8.** *Assume that for every sliced link diagram $S$ with exactly two crossings, and for all choices of arc orientations and undotted cup foams, the identity $\partial\psi(\square_S) = 1$ holds. Then Theorem 2.0.3 holds.* □

Recall the pictures of Table 2.1: they describe each possible isotopy class, together with the data of arc orientations, associated with such a sliced link diagram $S$. They are obtained by recording only the 0-resolution for both crossings together with their arc orientation. More precisely, pictures in Table 2.1 only picture the non-trivial local part, obtained by removing the simple closed loops that do not contain the boundary of an arc.

Recall the ladybug local arc presentation from Remark 2.2.2: this was the only case where the value of the $\mathfrak{sl}_2$-2-cocycle $\psi_{\mathfrak{sl}_2}$ could be set differently. For the other cases, a generic argument is sufficient:

**Proposition 2.3.9.** *Let $S$ be a sliced link diagram with exactly two crossings, together with choices of arc orientations and undotted cup foams. Then:*

*(i) if the local arc presentation of $S$ is not a ladybug, then $\partial\psi(\square_S) = 1$,*

*(ii) if the local arc presentation of $S$ is a ladybug, then $\partial\psi(\square_S) = 1$ or $\partial\psi(\square_S) = XY$.*

*Proof.* Recall the notations of Fig. 2.7. As $\psi$ is the 2-cochain controlling the commutativity in $\square_S$, by definition we have that $p = \partial\psi(\square_S)p$ for $p = \mathfrak{sl}(e_{*1}) \circ \mathfrak{sl}(e_{0*})$ (see Lemma 2.3.6). In particular:

$$(1 - \partial\psi(\square_S))p(1) = 0. \tag{2.7}$$

the element $p(1)$ admits a unique decomposition into basis elements: $p(1) = \sum_{i=1}^{n} \lambda_i a_{\delta_i}$ for some scalars $\lambda_i \in \Bbbk$ and some subsets $\delta_i \subset \mathfrak{sl}(W_{11})$. The above relation and the unicity of the decomposition implies that $(1 - \partial\psi(\square_S))\lambda_i = 0$ for all $i = 1, \ldots, n$. Hence, if any of the $\lambda_i$'s is invertible, we automatically get that $\partial\psi(\square_S) = 1$.

The only case where we get non-invertible coefficients is when $\mathfrak{sl}(e_{0*})$ is a split and $\mathfrak{sl}(e_{*1})$ is a merge, that is, in the ladybug cases. Instead, we have $p(1) = (1 + XY)a_j$ for some $j \in \mathfrak{sl}(W_{11})$, which forces either $\partial\psi(\square_S) = 1$ or $\partial\psi(\square_S) = XY$. □

If $\partial\psi(\square_S) = 1$ holds for all choices with a ladybug local arc presentation, then Theorem 2.0.3 holds. If on the contrary $\partial\psi(\square_S) = XY$ holds for all choices with a ladybug local arc presentation, then Theorem 2.0.3 still holds, as we can apply the same reasoning using instead the $\mathfrak{sl}_2$-2-cocycle $\overline{\psi}_{\mathfrak{sl}_2}$ defined in Remark 2.2.2. This amounts to compare our $\mathfrak{gl}_2$-construction with the $\mathfrak{sl}_2$-construction of type Y. In either case, we construct an isomorphism between the $\mathfrak{gl}_2$-hypercube and the $\mathfrak{sl}_2$-hypercubes of type X or type Y, the latter two being isomorphic.

In other words, what matters is that, whatever the value of $\partial\psi(\square_S)$, it remains the same for all the choices involved. This is given by the following proposition:

**Proposition 2.3.10.** *Let $S$ be a sliced link diagram with exactly two crossings, together with choices of arc orientations and undotted cup foams. Then the value of $\partial\psi(\square_S)$ only depends on the local arc presentation of $S$.*

The proof of Proposition 2.3.10 is given in subsection 2.3.5.

*Remark* 2.3.11. A direct computation shows that in fact, we do have $\partial\psi(\square_S) = 1$ even in the case of a ladybug local arc presentation. It is interesting to note that, if we carry out the proof for the $(\mathbb{Z}^2,\mu)$-graded-2-category $\mathbf{GFoam}'_d$ (see section 7.4 in part II) instead of $\mathbf{GFoam}_d$, then we get $\partial\psi(\square_S) = XY$ in the ladybug cases; see Remark 7.4.1 for further comments.

### 2.3.5 Independence on choices

We conclude the proof of Theorem 2.0.2 by proving Proposition 2.3.10. With the notation of Fig. 2.7, recall that $S$ is oriented as follows:

$$S = \begin{array}{ccc} W_{00} & \xrightarrow{e_{*0}} & W_{10} \\ {\scriptstyle e_{0*}}\downarrow & \circlearrowright & \downarrow{\scriptstyle e_{1*}} \\ W_{01} & \xrightarrow[e_{*1}]{} & W_{11} \end{array}$$

and similarly for $\mathfrak{sl}(S)$. We compute that:

$$\partial\psi(\square_S) = \psi_{\mathfrak{gl}_2}(S)\psi_{\mathfrak{sl}_2}(\mathfrak{sl}(S))^{-1}\partial\psi_\beta(S).$$

**Lemma 2.3.12.** *Let $S$ be a sliced link diagram with exactly two crossings, together with choices of arc orientations and undotted cup foams. Then the value of $\partial\psi(⧉_S)$ does not depend on the choices of arc orientations and undotted cup foams.*

*Proof.* Assume we swap the arc orientation of the $k$th crossing. Then $\partial\psi_\beta(S)$ will contribute an additional factor $XY$ for every split in direction $k$. Looking case by case at Table 2.1, one checks that this change is exactly compensated by the contribution of $\psi_{\mathfrak{sl}_2}(\mathfrak{sl}(S))$. The value of $\psi_{\mathfrak{gl}_2}(S)$ does not change. Hence, changing the arc orientations does not change the value of $\partial\psi(⧉_S)$.

Assume we change the choice of undotted cup foams instead. Let $\beta$ and $\overline{\beta}$ be two choices of undotted cup foams, identical everywhere except at some vertex $W$. Denote by $\tau$ the invertible scalar such that $\overline{\beta}^W = \tau\beta^W$.[5] Then:

$$\psi_{\overline{\beta}}(\to W) = \psi_\beta(\to W)/\tau \quad \text{and} \quad \psi_{\overline{\beta}}(W \to) = \psi_\beta(W \to) \cdot \tau,$$

where $\to W$ (resp $W \to$) denotes an incoming (resp. outgoing) edge in $H_{\mathfrak{gl}_2}(S)$. In particular, $\partial\psi_\beta(S) = \partial\psi_{\overline{\beta}}(S)$. One also checks that the values of $\psi_{\mathfrak{gl}_2}(S)$ and $\psi_{\mathfrak{sl}_2}(\mathfrak{sl}(S))$ do not change. Hence, changing the choice of undotted cup foams does not change the value of $\partial\psi(⧉_S)$. □

Next we check Proposition 2.3.10 for planar isotopies. Actually, we show a bit more, namely that the value of $\partial\psi(⧉_S)$ is independent of the choice of spatial representative for $S$. A *spatial sliced representative* of a link diagram $S$ is a sliced diagram that is a representative of $S$ up to spatial isotopies (see Definition 1.2.3).

**Lemma 2.3.13.** *Let $S$ be a sliced link diagram with exactly two crossings, together with choices of arc orientations and undotted cup foams. Then the value of $\partial\psi(⧉_S)$ does not depend on the spatial sliced representative.*

*Proof.* Let $S$ and $\overline{S}$ be two sliced link diagrams in the same spatial isotopy class. Throughout the proof we use the notation $\overline{(-)}$ to distinguish features related to $\overline{S}$. By Lemma 2.3.12, we may freely choose arc orientations and undotted cup foams.

Pick a choice of arc orientations on $S$. Then there are arc orientations on $\overline{S}$ naturally associated with the one on $S$: for instance, one can use an orientation on $S$ to record arc orientations, and orientation of link diagrams is preserved by spatial isotopies. With this choice, we have $\overline{\psi}_{\mathfrak{sl}_2} = \psi_{\mathfrak{sl}_2}$.

[5]The existence of such a scalar can be deduced from Theorem 1.6.1.

Recall that if $S$ and $\overline{S}$ are planar isotopic, then by Lemma 2.1.4 there exist scalar assignments $\epsilon$ and $\overline{\epsilon}$ such that $H_{\mathfrak{gl}_2}(S;\epsilon)$ and $H_{\mathfrak{gl}_2}(\overline{S};\overline{\epsilon})$ are isomorphic. This extends to the spatial case, thanks to the following pair of isomorphisms in **GFoam**:

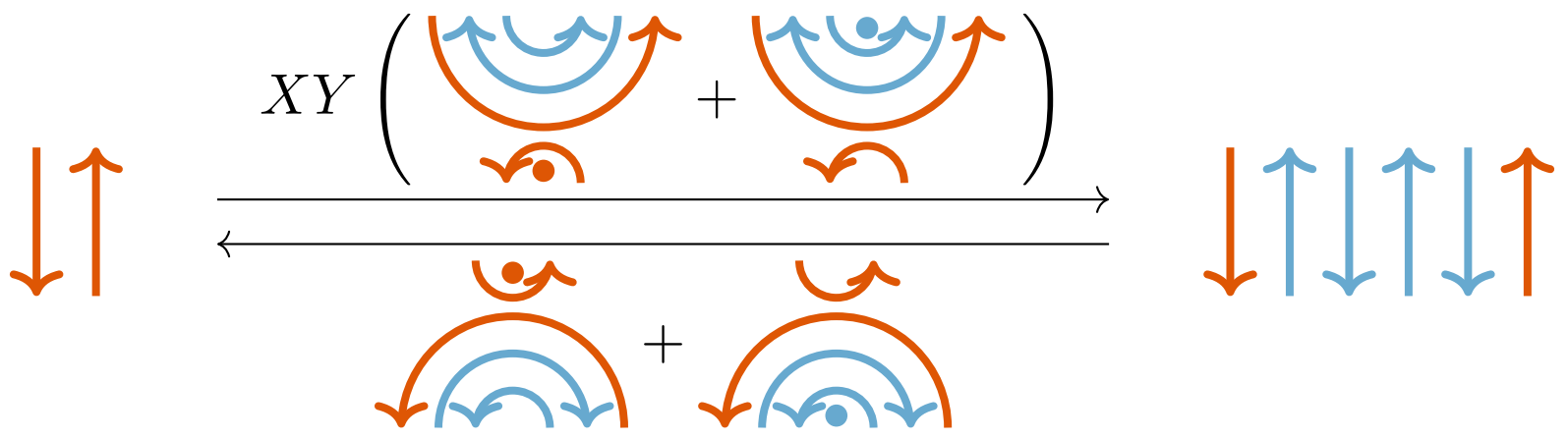

Denote $\varphi\colon H_{\mathfrak{gl}_2}(S;\epsilon) \to H_{\mathfrak{gl}_2}(\overline{S};\overline{\epsilon})$ such an isomorphism. If $\beta$ is a choice of undotted cup foams for $S$, then $\varphi \circ \beta$ is a choice of undotted cup foams for $\overline{S}$. Then, for each edge $e\colon \boldsymbol{r} \to \boldsymbol{s}$ in $\{0,1\}^2$, denoting $F_e$ and $\overline{F}_e$ the corresponding foams in $H_{\mathfrak{gl}_2}(S)$ and $H_{\mathfrak{gl}_2}(\overline{S})$ respectively, we have that:

$$\begin{aligned}\overline{\epsilon}(e)\left(\overline{F}_e \circ \overline{\beta}^{r}\right) &= \epsilon(e)\left(\varphi_{\boldsymbol{s}} \circ F_e \circ \beta^{r}\right)\\ &= \epsilon(e)\psi_\beta(e)\begin{cases}\beta^{s} & \text{if } F_e, \overline{F}_e \text{ are merges,}\\ (\beta^{s}_{i_1} + XY\beta^{s}_{i_2}) & \text{if } F_e, \overline{F}_e \text{ are splits.}\end{cases}\end{aligned}$$

Hence, we have $\psi_{\overline{\beta}}\overline{\epsilon} = \psi_\beta \epsilon$. That implies $\partial(\psi_{\overline{\beta}})\overline{\psi}_{\mathfrak{gl}_2} = \partial(\psi_\beta)\psi_{\mathfrak{gl}_2}$. This concludes the proof. ☐

Finally, we need to check that $\partial\psi(\square_S)$ only depends on the local arc presentation. Up to spatial isotopies, we can slide away all closed simple loops. The next lemma concludes the proof of Proposition 2.3.10.

**Lemma 2.3.14.** *Let $D_0$ and $D_1$ be two sliced link diagrams, such that $D_1$ has two crossings and $D_0$ none. Then $\partial\psi(\square_{D_1}) = \partial\psi(\square_{D_0 \star_0 D_1})$.*

*Proof.* If $W_0$ is the web corresponding to $D_0$ and $\beta^{W_0}$ is a choice of undotted cup foam for $W_0$, then $(\mathrm{id}_{W_0} \star_0 \beta^{W}) \star_1 \beta^{W_0}$ is an undotted cup foam for every undotted cup foam $\beta^{W}$. This defines a choice of undotted cup foams on $D_0 \star_0 D_1$, given one on $D_1$. It is also clear how to define arc orientations on $D_0 \star_0 D_1$ given such a choice on $D_1$. With these choices, $\psi_{\mathfrak{gl}_2}$, $\psi_{\mathfrak{sl}_2}$ and $\psi_\beta$ remain identical. This concludes. ☐

# 3

# A graded-categorification of the $q$-Schur algebra of level two

This chapter introduces a diagrammatic graded-2-category that categorifies the $q$-Schur algebra of level two. Then, we define a graded foamation 2-functor that relates this construction to graded $\mathfrak{gl}_2$-foams. This can be seen as a (partial) super analogue of [115] in the $\mathfrak{gl}_2$ case. See also [125] for earlier work. The casual reader is referred to the introduction, and in particular subsections i.1.8 and i.2.3, for further context.

Diagrammatic categorification of quantum groups was independently introduced by Khovanov, Lauda [107] and Rouquier [166]. A super analogue of this construction was given by Brundan and Ellis [25], building on earlier work of Kang, Kashiwara and Tsuchioka [102]. For odd $\mathfrak{sl}_2$, this was already studied in [74, 75]. See also [90, 100, 101] for related work.

In the even case, a categorification of the $q$-Schur algebra appeared in [127]. In [183], the second author defined a supercategorification of the negative half of the $q$-Schur algebra of level two. A graded version of this construction was given in [144]. In the same paper, Naisse and Putyra also defined a "1-map" from this graded version to their construction. This 1-map has similarities with our graded foamation 2-functor: we expect the two to coincide once an equivalence between [144] and our category of graded $\mathfrak{gl}_2$-foams is found.

Our presentation of the categorification of the $q$-Schur algebra of level two has similitudes with the presentation of super Kac–Moody 2-algebras in [25]. However, and contrary to the even case, it is not obtained as a quotient of their construction. For instance, the parities of cups and caps do not match. It would be interesting to find a relationship, if there exists any.

Section 3.1 review the definition of $\dot{S}_{n,d}$, the idempotented $q$-Schur algebra of level two. Its graded-categorification, that we call the *graded 2-Schur algebra* $\mathcal{GS}_{n,d}$, is introduced in section 3.2. Section 3.3 then defines the *graded foamation 2-functor* from $\mathcal{GS}_{n,d}$ into $\mathbf{GFoam}_d$, our graded-2-category of $\mathfrak{gl}_2$-foams defined in section 1.3. Finally, we show in section 3.4 that $\mathcal{GS}_{n,d}$ categorifies $\dot{S}_{n,d}$.

## 3.1 The $q$-Schur algebra of level two

**Definition 3.1.1.** *The* (idempotented) $q$-Schur algebra of level two, *or simply* Schur algebra, *is the* $\mathbb{Z}[q,q^{-1}]$*-linear category* $\dot{S}_{n,d}$ *such that:*

- *Objects are weights in the set*

$$\Lambda_{n,d} \coloneqq \{\lambda \in \{0,1,2\}^n \mid \lambda_1 + \ldots + \lambda_n = d\}.$$

- *Morphisms are* $\mathbb{Z}[q,q^{-1}]$*-linear combinations of iterated compositions of identity morphisms* $1_\lambda\colon \lambda \to \lambda$ *and generating morphisms*

$$e_i 1_\lambda\colon \lambda \to \lambda + \alpha_i \quad \text{and} \quad f_i 1_\lambda\colon \lambda \to \lambda - \alpha_i \qquad i = 1, \ldots, n-1,$$

*where* $\alpha_i \coloneqq (0,\ldots,1,-1,\ldots,0) \in \mathbb{Z}^n$ *with 1 being on the* $i$*-th coordinate. Morphisms are subject to the* Schur quotient

$$1_\lambda = 0 \qquad \text{if } \lambda \notin \Lambda_{n,d}$$

*and to the following relations:*

$$\begin{cases} (e_i f_j - f_j e_i)1_\lambda = \delta_{ij}[\overline{\lambda}_i]_q 1_\lambda, & \\ (e_i e_j - e_j e_i)1_\lambda = 0 & \text{for } |i-j| > 1, \\ (f_i f_j - f_j f_i)1_\lambda = 0 & \text{for } |i-j| > 1, \end{cases} \tag{3.1}$$

*where* $\overline{\lambda}_i := \lambda_i - \lambda_{i+1}$, *the symbol* $\delta_{ij}$ *is the Kronecker delta and*

$$[m]_q = q^{m-1} + q^{m-3} + \ldots + q^{3-m} + q^{1-m}$$

*is the* $m$*th quantum integer.*

Recall that a $\mathbb{Z}[q, q^{-1}]$-linear category is the same as a $\mathbb{Z}[q, q^{-1}]$-algebra with a distinguished set of idempotents, so that $\dot{S}_{n,d}$ is indeed an algebra. The "level two" stands for the fact that the value of coordinates is at most two. The Schur quotient implies that a morphism that factors through a weight not in $\Lambda_{n,d}$ is set to zero. In the sequel, it is understood that an expression involving a weight that does not belong to $\Lambda_{n,d}$ is set to zero.

### 3.1.1 Ladder diagrammatics

Morphisms of the Schur algebra admit a diagrammatics called *ladder diagrams*. Given some object $\lambda \in \Lambda_{n,d}$, represent the coordinates 0, 1 and 2 respectively with a dotted line ⋯⋯ , a single line — and a double line ═ . The identity $1_\lambda$ is then presented by a sequence of horizontal lines, reading coordinates from bottom to top. For instance (on the left):

$1_{(0,1,1,0,2,1)} =$ 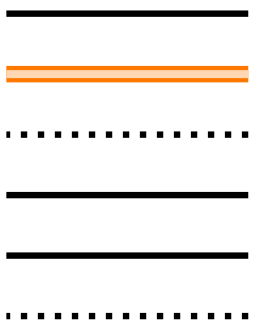 $\mathsf{F}_3(0,1,1,0,2,1) =$ 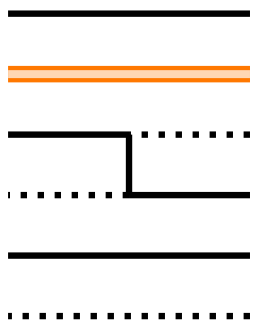

Generating morphisms $\mathsf{E}_i(\lambda)$ and $\mathsf{F}_i(\lambda)$ are pictured similarly to identities, differing only locally around coordinates $(\lambda_i, \lambda_{i+1})$ (see above on the right). Below are all the possible local forms for $\mathsf{E}$ and $\mathsf{F}$:

$\mathsf{E}(1,1) =$ 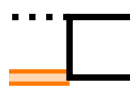 $\mathsf{E}(0,1) =$ 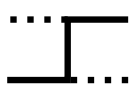 $\mathsf{E}(1,2) =$ 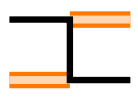 $\mathsf{E}(0,2) =$ 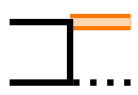

$\mathsf{F}(1,1) =$ 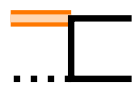 $\mathsf{F}(1,0) =$ 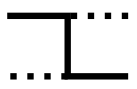 $\mathsf{F}(2,1) =$ 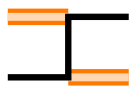 $\mathsf{F}(2,0) =$ 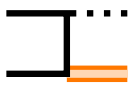

Here we read them from right to left, and composition is given by stacking horizontally. Diagrams obtained in that way are called *ladder diagrams*, and morphisms in $\dot{S}_{n,d}$ are presented by $\mathbb{Z}[q, q^{-1}]$-linear combination of ladder diagrams.

In this diagrammatics, the first relation in (3.1) when $i = j$, that is the relation

$$(\mathsf{EF} - \mathsf{FE})(\lambda) = [\overline{\lambda}]_q$$

suppressing the indices $i$, gives

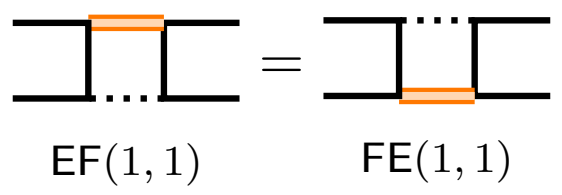

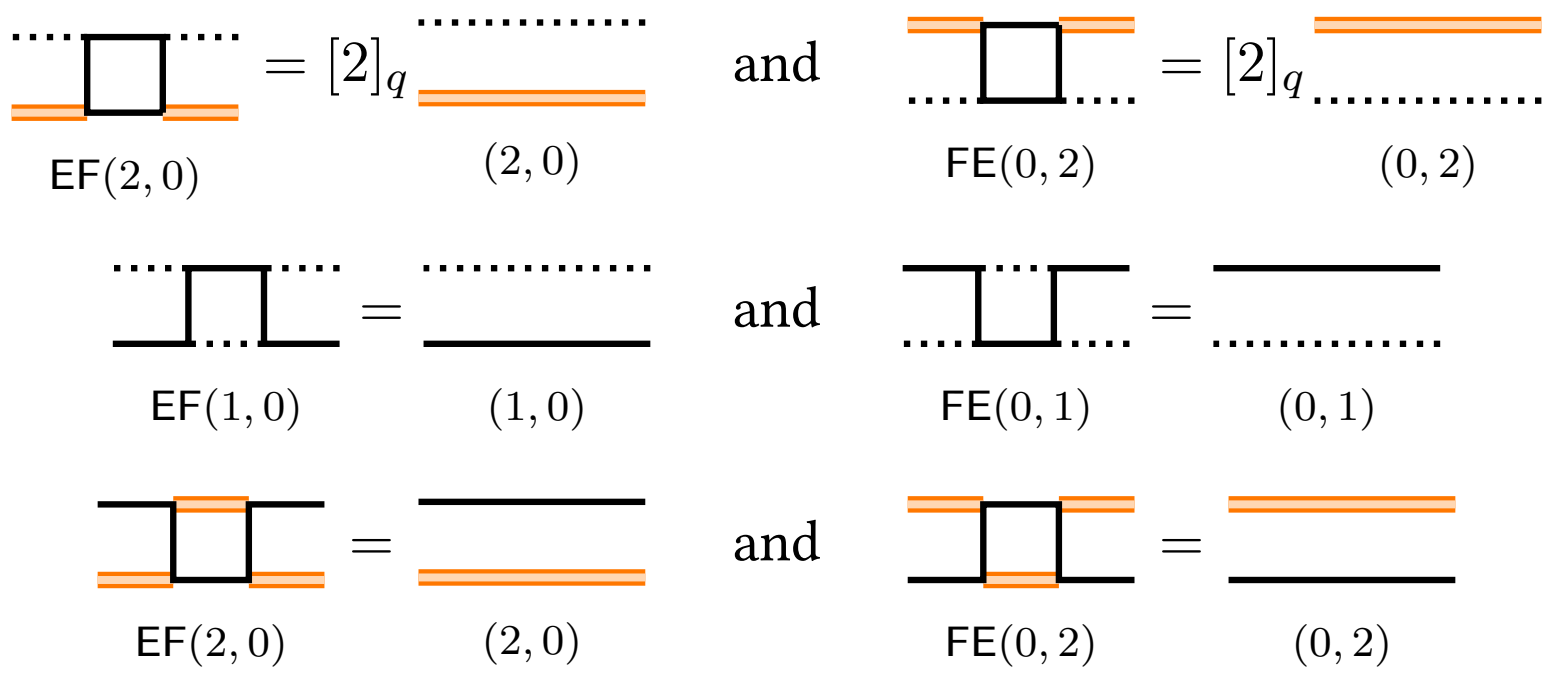

Note our notation $1_\lambda = \lambda$. When $|i - j| > 1$, the three relations in (3.1) implies that two ladder's rungs (the vertical bars) can be moved one past the other as long as they don't have any uprights (horizontal bars) in common. When $|i - j| = 1$, the first relation of (3.1) implies that E's and F's can be moved past one another, but not two E's nor two F's. For instance:

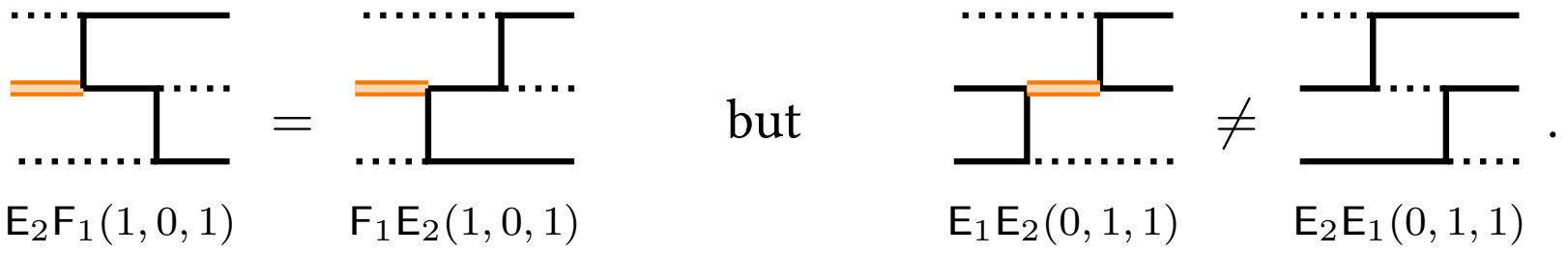

## 3.1.2 Relationship with $\mathfrak{gl}_2$-webs

The ladder diagrammatics can be understood as a "rigidification" of web diagrammatics (see section 1.2). Forgetting zero entries defines a mapping $\lambda \mapsto \underline{\lambda}$ from $\Lambda_{n,d}$ to $\underline{\Lambda}_{n,d}$, and forgetting dotted lines and smoothing out corners defines a functor $F_{n,d}$ from $\dot{S}_{n,d}$ to $\mathbf{Web}_d$.

**Lemma 3.1.2.** *The functor $F_{n,d} \colon \dot{S}_{n,d} \to \mathbf{Web}_d$ is faithful.*

*Proof.* In this proof, we want to explicitly consider the defining relations of $\dot{S}_{n,d}$ and $\mathbf{Web}_d$. For that reason, we will say that two ladders (resp. webs) are *congruent* to mean that they are equal in $\dot{S}_{n,d}$ (resp. $\mathbf{Web}_d$), and write it with the symbol $\sim$.

Both $\dot{S}_{n,d}$ and $\mathbf{Web}_d$ have scalar defining relations, in the sense that each relation is of the form $W_1 \sim \lambda W_2$ for $W_1$ and $W_2$ two diagrams (ladders or webs) and $\lambda \in \mathbb{Z}[q, q^{-1}]$. For webs, this follows directly from the definition. While the first defining relation (3.1) of the Schur algebra does not seem scalar at first glance, we have seen in the description of the ladder diagrammatics that it is, in fact, scalar. Moreover, scalars are coherently controlled by the isotopy classes of the underlying tangle diagrams, so that if $W \sim \lambda W$ either in $\dot{S}_{n,d}$ or $\mathbf{Web}_d$, then $\lambda = 1$. We deduce that both $\dot{S}_{n,d}$ and $\mathbf{Web}_d$ are free $\mathbb{Z}[q, q^{-1}]$-algebras. In particular, we can consider their respective localizations at $\mathbb{Q}(q)$, respectively denoted $(\dot{S}_{n,d})_{\mathbb{Q}(q)}$ and $(\mathbf{Web}_d)_{\mathbb{Q}(q)}$, as well as the localization of the functor $F_{n,d}$ itself:

$$(F_{n,d})_{\mathbb{Q}(q)} \colon (\dot{S}_{n,d})_{\mathbb{Q}(q)} \to (\mathbf{Web}_d)_{\mathbb{Q}(q)}$$

Showing that $F_{n,d}$ is faithful reduces to showing that $(F_{n,d})_{\mathbb{Q}(q)}$ is faithful.

Working in the localization allows to define the following ladders:

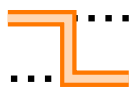 $:= \dfrac{1}{q+q^{-1}}$ 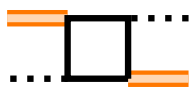 and 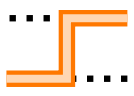 $:= \dfrac{1}{q+q^{-1}}$ 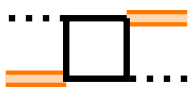

One checks that they define each other inverses in $(\dot{S}_{n,d})_{\mathbb{Q}(q)}$. Using these extra ladders together with 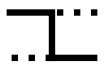 and 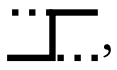, one can define an isomorphism between any pair of objects in $(\dot{S}_{n,d})_{\mathbb{Q}(q)}$ that have the same image in $(\mathbf{Web}_d)_{\mathbb{Q}(q)}$. Say that an object in $(\dot{S}_{n,d})_{\mathbb{Q}(q)}$ is *right-justified* if all zero entries are at the beginning of the tuple. Every object $\lambda \in \Lambda_{n,d}$ has an associated right-justified object $\widehat{\lambda}$, together with an isomorphism $\alpha_\lambda \colon \lambda \to \widehat{\lambda}$ in $(\dot{S}_{n,d})_{\mathbb{Q}(q)}$. To show faithfulness, we can restrict our attention to hom-spaces with right-justified source and target.

In fact, ladders themselves can be right-justified, as in the following example:

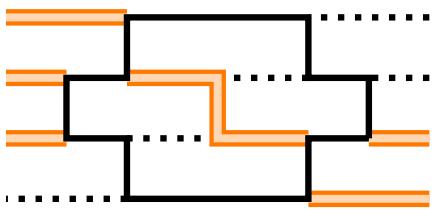 $\rightsquigarrow$ 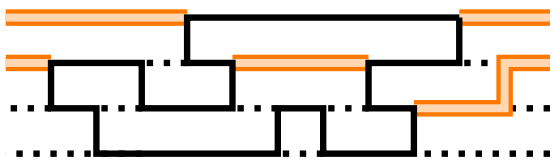

The procedure consists in right-justifying every generator of the ladder with the help of the isomorphisms $\alpha_\lambda$, reading from bottom to top. We leave the details to the reader.

Finally, webs relations lift directly to relations on right-justified ladders. For instance

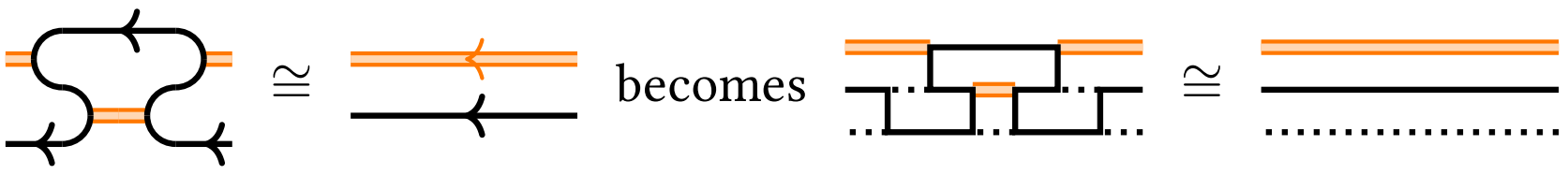

Hence, if $L_1$ and $L_2$ are two right-justified ladders and if $F_{n,d}(L_1)$ and $F_{n,d}(L_2)$ are congruent in $\mathbf{Web}_d$, then this congruence lifts to a congruence between $L_1$ and $L_2$ in $\dot{S}_{n,d}$. This concludes. □

## 3.2 The graded 2-Schur algebra

Recall the definitions of $\Bbbk$ and $\mu\colon \mathbb{Z}^2 \times \mathbb{Z}^2 \to \Bbbk^\times$ in Definition 1.3.2 and of $q\colon \mathbb{Z}^2 \to \mathbb{Z}$ in Remark 1.3.5. We use the following notation:

$$p_{ij} := ((\alpha_j)_{i+1}, -(\alpha_j)_i) = \begin{cases} (0,1) & \text{if } j = i-1, \\ (-1,-1) & \text{if } j = i, \\ (1,0) & \text{if } j = i+1, \\ (0,0) & \text{otherwise.} \end{cases}$$

**Definition 3.2.1.** *The* graded 2-Schur algebra $\mathcal{GS}_{n,d}$ *is the $\mathbb{Z}$-graded $(\mathbb{Z}^2,\mu)$-graded-2-category such that:*

- *Objects are elements $\lambda$ for $\lambda \in \Lambda_{n,d}$.*
- *1-morphisms are compositions of the generating 1-morphisms*

$$1_\lambda\colon \lambda \to \lambda, \quad \mathsf{F}_i 1_\lambda\colon \lambda \to \lambda - \alpha_i \quad \text{and} \quad \mathsf{E}_i 1_\lambda\colon \lambda \to \lambda + \alpha_i,$$

  *whenever both $\lambda$ and $\lambda - \alpha_i$ (resp. $\lambda$ and $\lambda + \alpha_i$) are objects of $\Lambda_{n,d}$. Using string diagrammatics, the identity $1_\lambda$ is not pictured, and the non-trivial 1-generators are pictured as follows:*

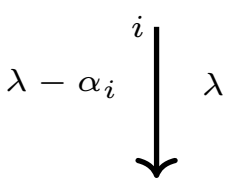

$= \mathrm{id}_{\mathsf{F}_i 1_\lambda}$

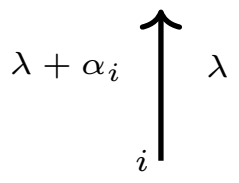

$= \mathrm{id}_{\mathsf{E}_i 1_\lambda}$

  *Note that we read from bottom to top and from right to left.*

- *2-morphisms are* $\Bbbk$*-linear combination of string diagrams generated by formal vertical and horizontal compositions of the following generating 2-morphisms:*

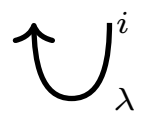 $: 1_\lambda \to \mathsf{E}_i\mathsf{F}_i 1_\lambda$ $\qquad$ 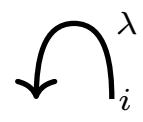 $: \mathsf{E}_i\mathsf{F}_i 1_\lambda \to 1_\lambda$

$$(\lambda_{i+1}, -\lambda_i) + (0,1) \qquad\qquad -(\lambda_{i+1}, -\lambda_i) + (1,0)$$

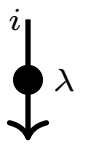 $: \mathsf{F}_i 1_\lambda \to \mathsf{F}_i 1_\lambda$ $\qquad$ 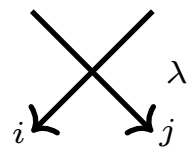 $: \mathsf{F}_i\mathsf{F}_j 1_\lambda \to \mathsf{F}_j\mathsf{F}_i 1_\lambda$

$$(1,1) \qquad\qquad p_{ij}$$

*where the* $\mathbb{Z}^2$*-degree* $\deg_{\mathbb{Z}^2}$ *is given below each generator. Such string diagrams are called* Schur diagrams.

*2-morphisms are further subject to axioms described below. The* quantum grading *on* $\mathcal{GS}_{n,d}$ *is the* $\mathbb{Z}$*-grading defined by* $\operatorname{qdeg}(D) := q(\deg_{\mathbb{Z}^2}(D))$.

One only needs to label one region is a given diagram, as it determines the label of all the other regions. If this forces a region to be labelled by a weight $\lambda$ that does not belong to $\Lambda_{n,d}$, we set this diagram to zero. In that case, we say that the diagram is zero due to the Schur quotient.

We assume that the generators in a string diagram are always in generic position, in the sense that the vertical projection defines a separative Morse function. Generating 2-morphisms are subject to the following local relations:

(1) As in any graded-2-category, we have the graded interchange law:

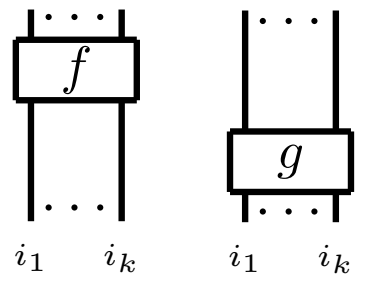 $= \mu(\deg_{\mathbb{Z}^2} f, \deg_{\mathbb{Z}^2} g)$ 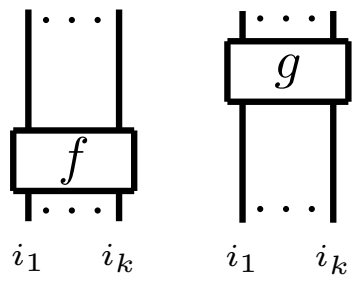

(2) Two dots annihilate:

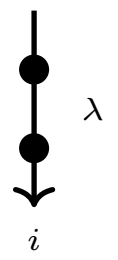 $= 0 \qquad (3.2)$

(3) Graded KLR algebra relations for downward crossings:

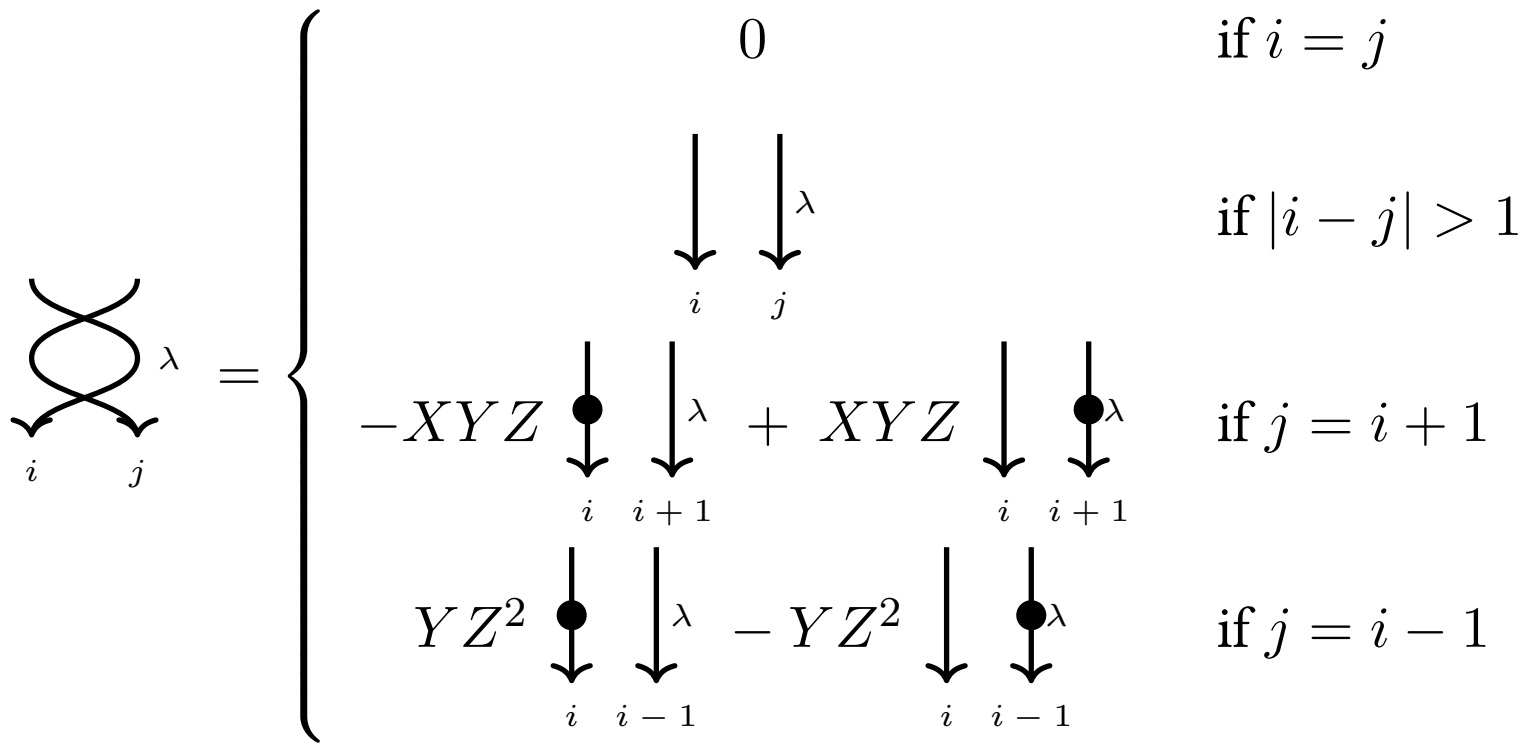

(3.3)

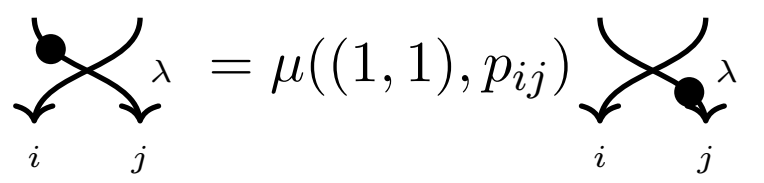

for $i \neq j$ (3.4)

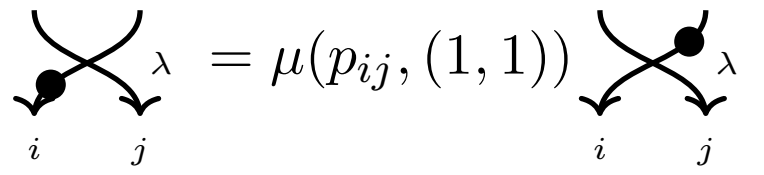

for $i \neq j$ (3.5)

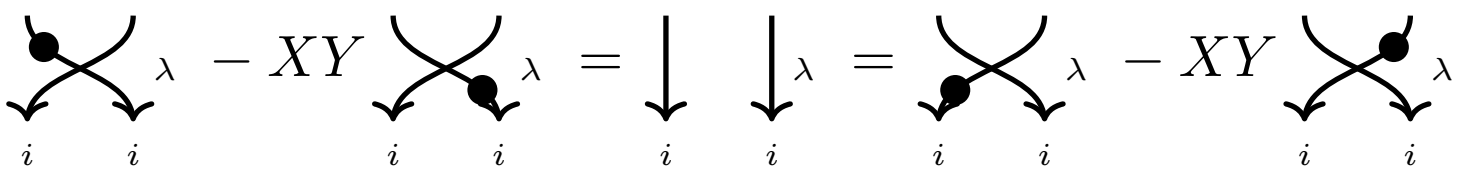

(3.6)

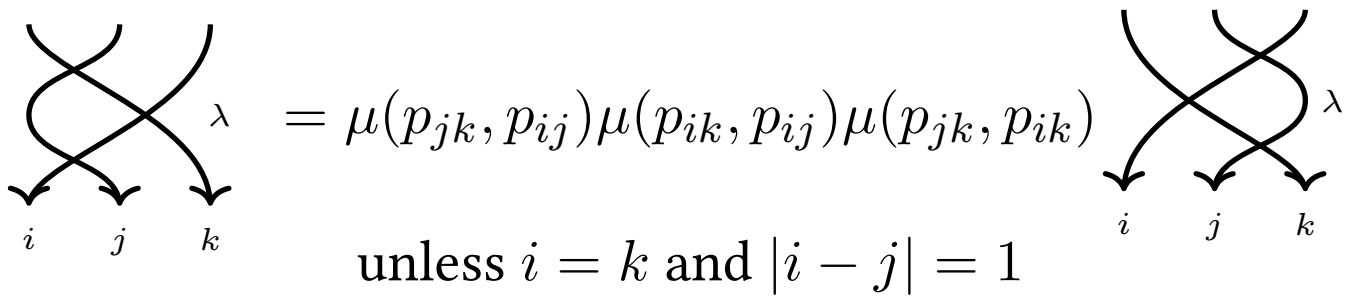

(3.7)

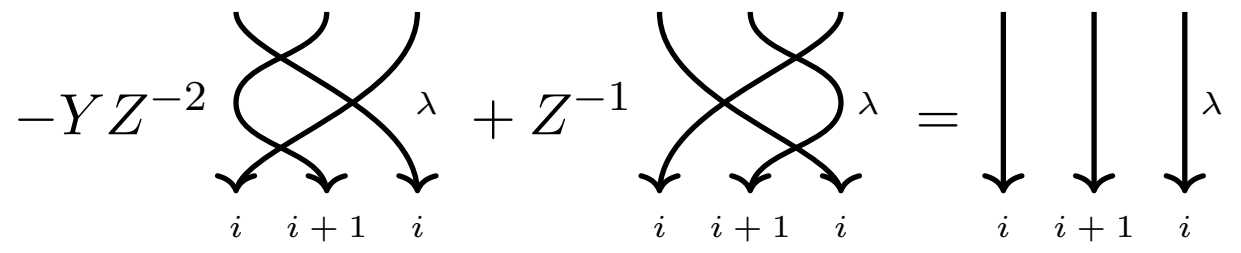

(3.8)

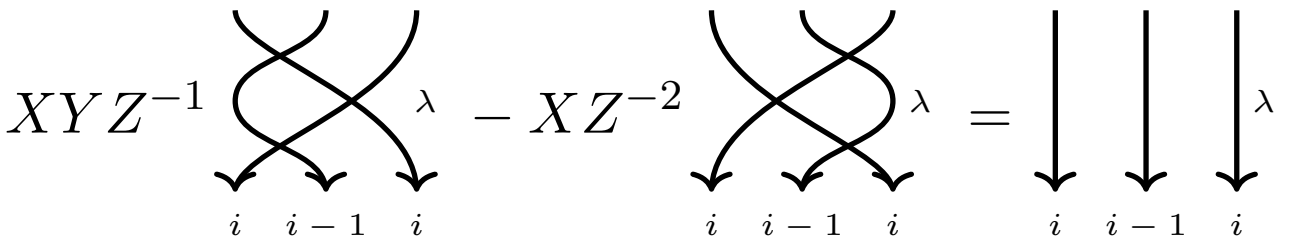

(3.9)

(4) Graded adjunction relations:

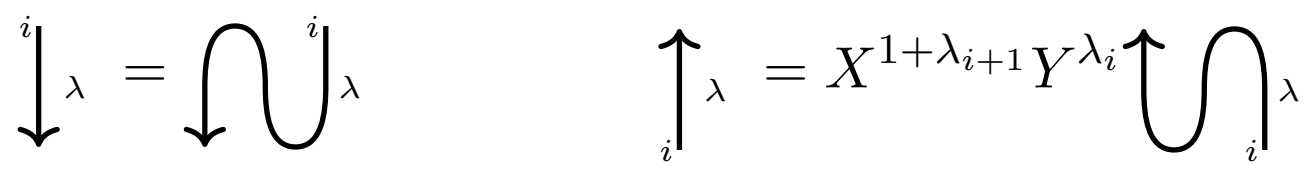

(3.10)

Finally, we require invertibility axioms. To define it, we introduce the following shorthands:

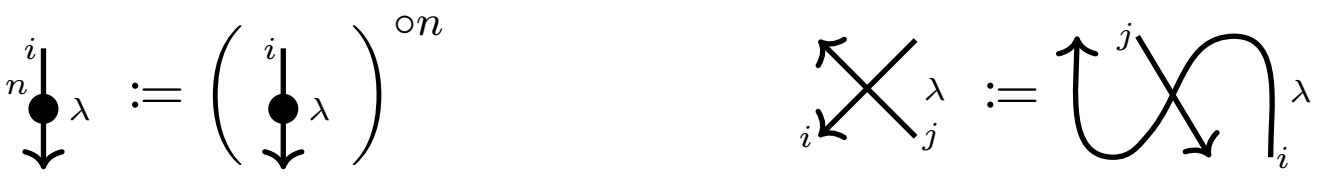

Recall the notation $\overline{\lambda}_i := \lambda_i - \lambda_{i+1}$. Then:

(5) Except if they are zero due to the Schur quotient, the following 2-morphisms are isomorphisms in the graded additive envelope of $\mathcal{GS}_{n,d}$ (see subsection 1.1.2):

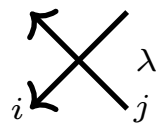 $: \mathsf{F}_i\mathsf{E}_j(\lambda) \to \mathsf{E}_j\mathsf{F}_i(\lambda)$ if $i \neq j$

(3.11)

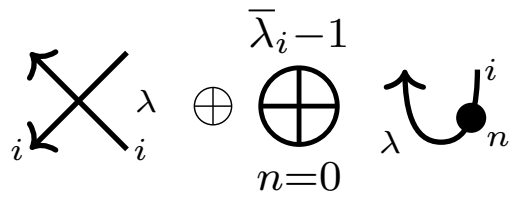 $: \mathsf{F}_i\mathsf{E}_i(\lambda) \oplus \lambda^{\oplus[\overline{\lambda}_i]} \to \mathsf{E}_i\mathsf{F}_i(\lambda)$ if $\overline{\lambda}_i \geq 0$

(3.12)

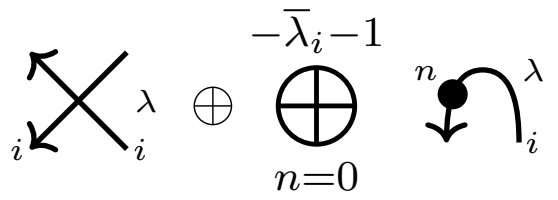 $: \mathsf{F}_i\mathsf{E}_i(\lambda) \to \mathsf{E}_i\mathsf{F}_i(\lambda) \oplus \lambda^{\oplus[-\overline{\lambda}_i]}$ if $\overline{\lambda}_i \leq 0$

(3.13)

This ends the definition of the relations on the graded 2-Schur algebra. ◇

*Remark* 3.2.2. Let us elaborate on the invertibility axioms above. They are equivalent to the existence of some unnamed generators which are entries of the inverse matrices of (3.11), (3.12) and (3.13), and some unnamed relations that precisely encompass the fact that those generators form inverse matrices. This definition follows Rouquier's approach [166] to 2-Kac–Moody algebras (categorified quantum groups) and Brundan and Ellis' approach [25] to super 2-Kac–Moody algebras. Unravelling the definition would lead to a more explicit (but heavier) definition, similar to Khovanov and Lauda's approach [107] to categorified quantum groups.

*Remark* 3.2.3. Some scalars in the relations above can be understood as artefacts of interchanging the vertical positions of the generators. For instance, the scalar in relation (3.4) is precisely the scalar that appears when interchanging a dot and a $(i, j)$-crossing. A similar reasoning applies to relations (3.5) and (3.7).

In [127, Definition 3.2], a categorification of the $q$-Schur algebra was constructed, denoted $\mathcal{S}(n, d)$. Let us write $\mathcal{S}(n, d)^\bullet$ for the $\Bbbk$-linear[1] 2-category obtained from $\mathcal{S}(n, d)$ by further imposing relation (3.2). Then:

**Proposition 3.2.4.** *Let $\underline{(\mathcal{GS}_{n,d})}_q^\oplus$ be additive shifted closure of $\mathcal{GS}_{n,d}$ (see subsection 1.1.2) with respect to the quantum grading (see subsection 1.1.2 and Definition 3.2.1). Then:*

$$\mathcal{S}(n, d)^\bullet \cong \underline{(\mathcal{GS}_{n,d})}_q^\oplus \Big|_{X=Y=Z=1}.$$

*Proof.* Explicitly defining the unnamed generators implied by the invertibility axioms and doing a relation chase exhibit the missing relations. One only need to rescale the $(i, i)$-downward crossing by $(-1)$. This exactly follows the proof from [22] showing the equivalence between Rouquier's definition of Kac–Moody 2-algebras and Khovanov-Lauda categorified quantum groups. See also [25] for a related statement in the super case. □

## 3.3 The graded foamation 2-functor

This section exhibits $\mathbf{GFoam}_d$ as a 2-representation of $\mathcal{GS}_{n,d}$. More precisely, for each $n, d \in \mathbb{N}$ there exists a $(\mathbb{Z}^2, \mu)$-graded 2-functor

$$\mathcal{F}_{n,d}\colon \mathcal{GS}_{n,d} \to \mathbf{GFoam}_d.$$

This categorifies the functor $F_{n,d}\colon \dot{S}_{n,d} \to \mathbf{Web}_d$ defined in subsection 3.1.2.

On the level of objects, the functor $F_{n,d}$ maps a weight $\lambda \in \Lambda_{n,d}$ to the weight $\underline{\lambda} \in \underline{\Lambda}_d$, obtained by forgetting all zero entries in $\lambda$. Recall the colour of a coordinate defined in section 1.2. For $i \in \{1, \ldots, n\}$ such that $\lambda_i \neq 0$, we denote $\underline{i}_\lambda$ the colour of the coordinate of the "image" of $i$ in $\underline{\lambda}$. For instance:

$$\underline{(1, 0, 1, 2, 0, 1)} = (1, 1, 2, 1) \quad \text{and} \quad \underline{1}_\lambda = 1,\ \underline{3}_\lambda = 2,\ \underline{4}_\lambda = 3 \text{ and } \underline{6}_\lambda = 5.$$

[1] The linear 2-category $\mathcal{S}(n, d)$ is defined over $\mathbb{Q}$ in [127], but it can be defined over $\mathbb{Z}$ and hence over any unital commutative ring $\Bbbk$.

In the string diagrammatics of foams, the functor $F_{n,d}$ is given by

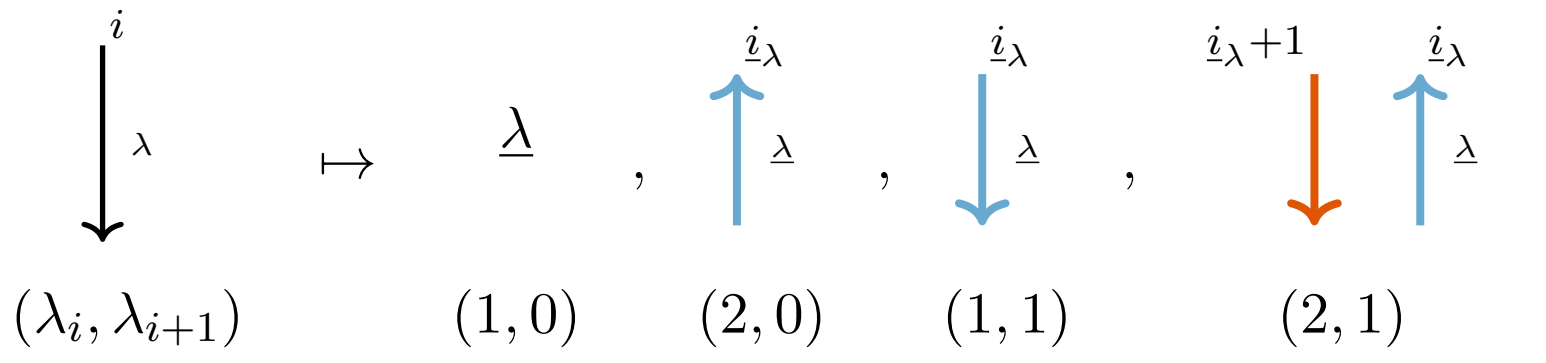

(3.14)

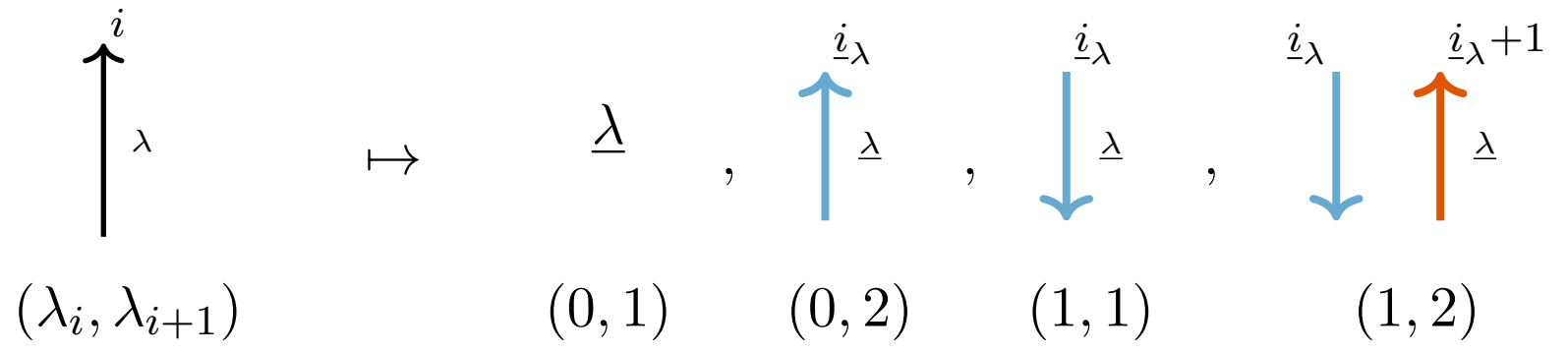

(3.15)

The local data of $(\lambda_i, \lambda_{i+1})$ is given below each case.

Following [144, p. 59], we shall use the scalar

$$\Gamma_\lambda(i) := (-XY)^{\#\{\lambda_j=1 | j \leq i\}}$$

to normalize the graded foamation 2-functor.

**Proposition 3.3.1.** *There exists a* $(\mathbb{Z}^2, \mu)$*-graded 2-functor*

$$\mathcal{F}_{n,d} \colon \mathcal{GS}_{n,d} \to \mathbf{GFoam}_d$$

*defined on generating 2-morphisms as in Fig. 3.1. We call* $\mathcal{F}_{n,d}$ *the* graded foamation 2-functor.

*Proof.* On checks that $\mathcal{F}_{n,d}$ preserves the $\mathbb{Z}^2$-grading. We need to check that the images through $\mathcal{F}_{n,d}$ of the defining relations of the graded 2-Schur algebra are relations in $\mathbf{GFoam}_d$. This is analogous to the proof of [115, Proposition 3.3]. The main additional work is to check that the scalars match. For readability, we leave implicit the label of regions for foam diagrams.

The fact that $\mathcal{F}_{n,d}$ respects relation (3.2) follows from dot annihilation in foams. For relation (3.3), the case $i = j$ follows from the evaluation of an undotted counterclockwise bubble and the case $|i-j| > 1$ follows from the Reidemeister II braid-like relation. Consider then the case $j = i-1$. For $\lambda_i = 2$, both sides of (3.3) are zero due to the Schur quotient. If $(\lambda_{i-1}, \lambda_i, \lambda_{i+1}) =$

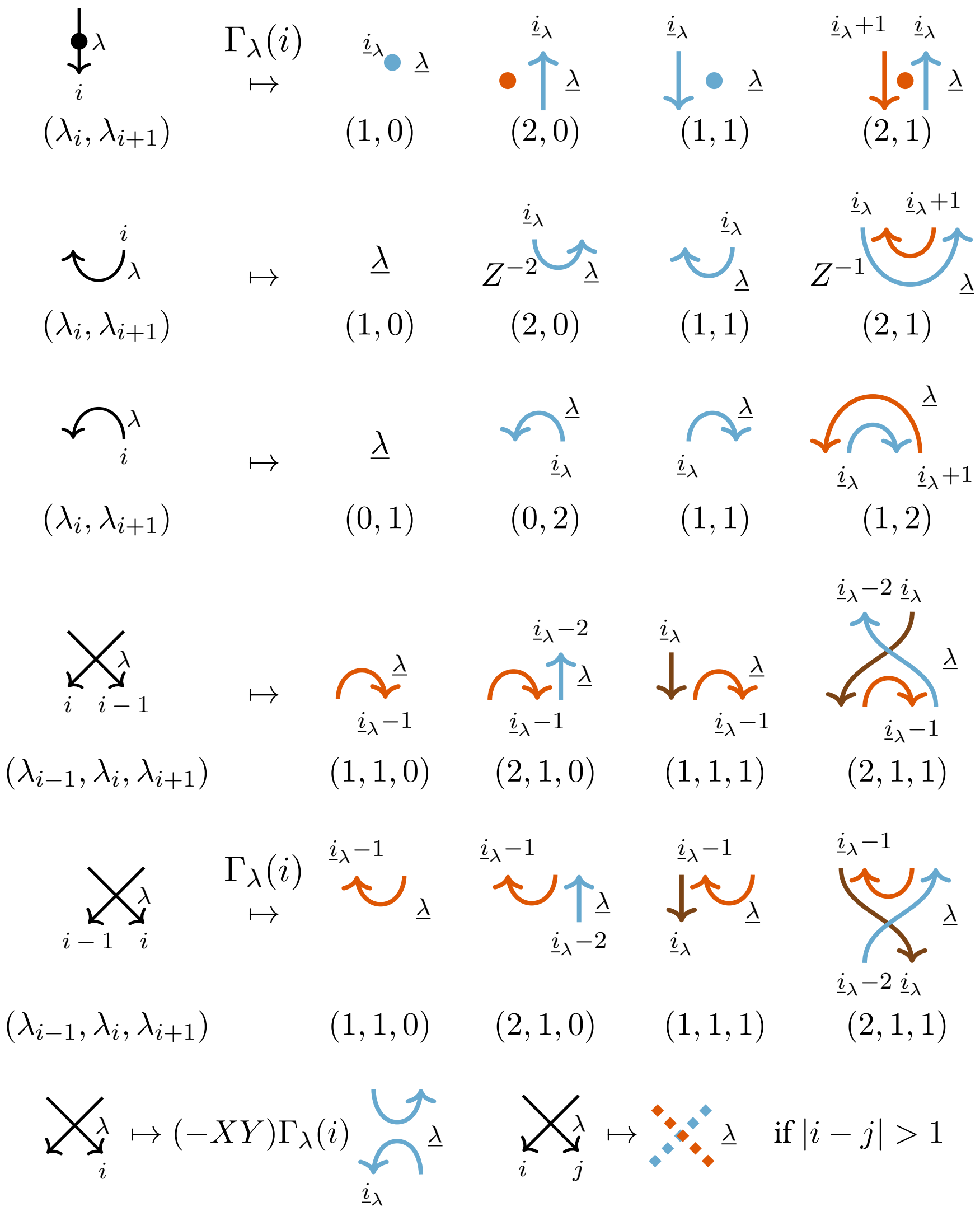

**Fig. 3.1** Definition of $\mathcal{F}_{n,d}$ on generating 2-morphisms. For dots, rightward cups and caps, and adjacent crossings, it depends on the local value of $\lambda$, which is given below each case. A symbol $\Gamma_\lambda(i)$ above a "mapsto" arrow means that the codomain should be multiplied by $\Gamma_\lambda(i)$. For distant crossings (last picture), the picture means that one should replace each strand of the Schur crossing with the corresponding strand or pair of strands as prescribed by (3.14). This defines a foam diagram consisting of one, two or four crossings.

$(1, 1, 0)$, we have:

$$\mathcal{F}_{n,d}\left(\right.$$ 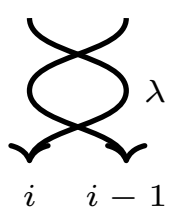 $$\left.\right) = \Gamma_\lambda(i)$$ 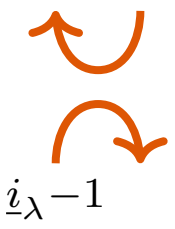

$$= \Gamma_\lambda(i)\left( YZ^2\right.$$ 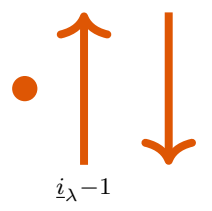 $$+ XZ^2$$ 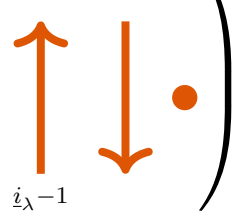 $$\left.\right)$$

$$= \mathcal{F}_{n,d}\left( YZ^2\right.$$ 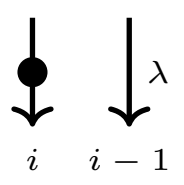 $$- YZ^2$$ 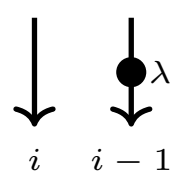 $$\left.\right)$$

where we used Lemma 1.4.4 and $\Gamma_{\lambda-\alpha_{i-1}}(i) = \Gamma_\lambda(i) = (-XY)\Gamma_\lambda(i-1)$. Other cases for which $\lambda_i = 1$ are computed similarly. If $\lambda_i = 0$, the left-hand side of (3.3) is zero due to the Schur quotient. If $(\lambda_{i-1}, \lambda_i, \lambda_{i+1}) = (2, 0, 1)$, the image of the right-hand side is

$$\Gamma_{\lambda-\alpha_{i-1}}(i)\, YZ^2$$ 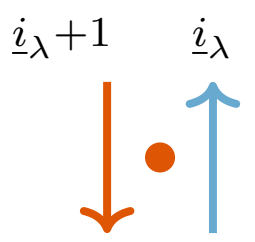 $$- \Gamma_\lambda(i-1)YZ^2$$ 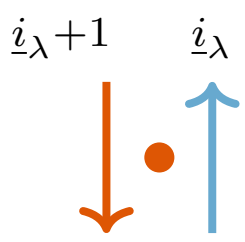 $$= 0,$$

which follows from $\Gamma_{\lambda-\alpha_{i-1}}(i) = \Gamma_\lambda(i-1)$. Other cases for which $\lambda_i = 0$ are computed similarly. When $j = i + 1$, both sides are zero due to the Schur quotient when $\lambda_j = 0$, and only the left-hand side when $\lambda_j = 2$. For the representative case $(\lambda_{j-1}, \lambda_j, \lambda_{j+1}) = (1, 1, 0)$, one gets:

$$\mathcal{F}_{n,d}\left(\right.$$ 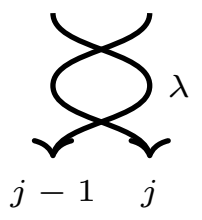 $$\left.\right) = \Gamma_\lambda(j)$$ 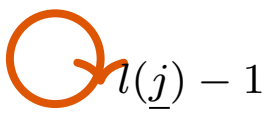

$$= \Gamma_\lambda(j)\left( Z\right.$$ 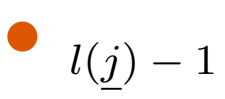 $$+ XYZ$$ 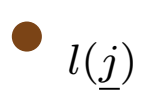 $$\left.\right)$$

$$= \mathcal{F}_{n,d}\left( -XYZ\right.$$ 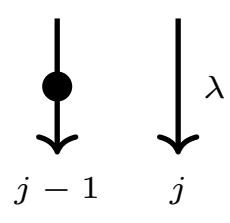 $$+ XYZ$$ 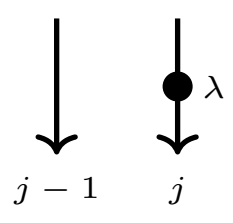 $$\left.\right)$$

which follows from $\Gamma_{\lambda-\alpha_j}(j-1) = (-XY)\Gamma_\lambda(j)$, and the latter from. For the representative case $(\lambda_{i-1}, \lambda_i, \lambda_{i+1}) = (1, 2, 0)$, the image of the right-hand

side of (3.3) gives:

$$-\Gamma_{\lambda-\alpha_j}(j-1)\,XYZ$$ 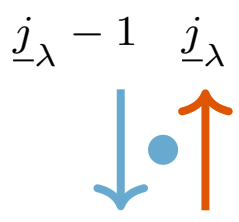

$$+\,\Gamma_\lambda(j)XYZ$$ 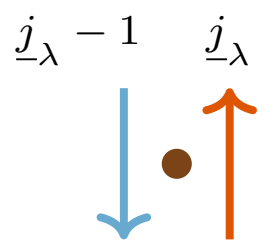

$$= 0,$$

which follows from $\Gamma_{\lambda-\alpha_j}(j-1) = \Gamma_\lambda(j)$ and dot migration.

For relations (3.4) and (3.5), it follows from the graded interchange law and dot migration. The relation (3.6) follows directly from the (vertical) neck-cutting relation.

Consider now relation (3.7). If all colours are pairwise equal or adjacent (i.e. cases $i = j = k$ and $i = j = k \pm 1$ together with permutations), then either the case is excluded by assumption, or both sides are zero by the Schur quotient. In particular, we can disregard the normalization with the Γ's (this follows from the fact that $\Gamma_{\lambda-\alpha_j}(i) = \Gamma_\lambda(i)$ whenever $i \neq j$). If a colour is distant from the two others (i.e. cases $|i-j| > 1$ and $|i-k| > 1$ together with permutations), then the image under $\mathcal{F}_{n,d}$ of relation (3.7) is a composition of graded interchanges, consisting in moving vertically the image of the only crossing with a (possibly) non-trivial $\mathbb{Z}^2$-grading. The last six cases also consist of graded interchanges, each interchanging two saddles (zip or unzip) sharing a common 1-facet. We picture below the domain and codomain of the six cases:

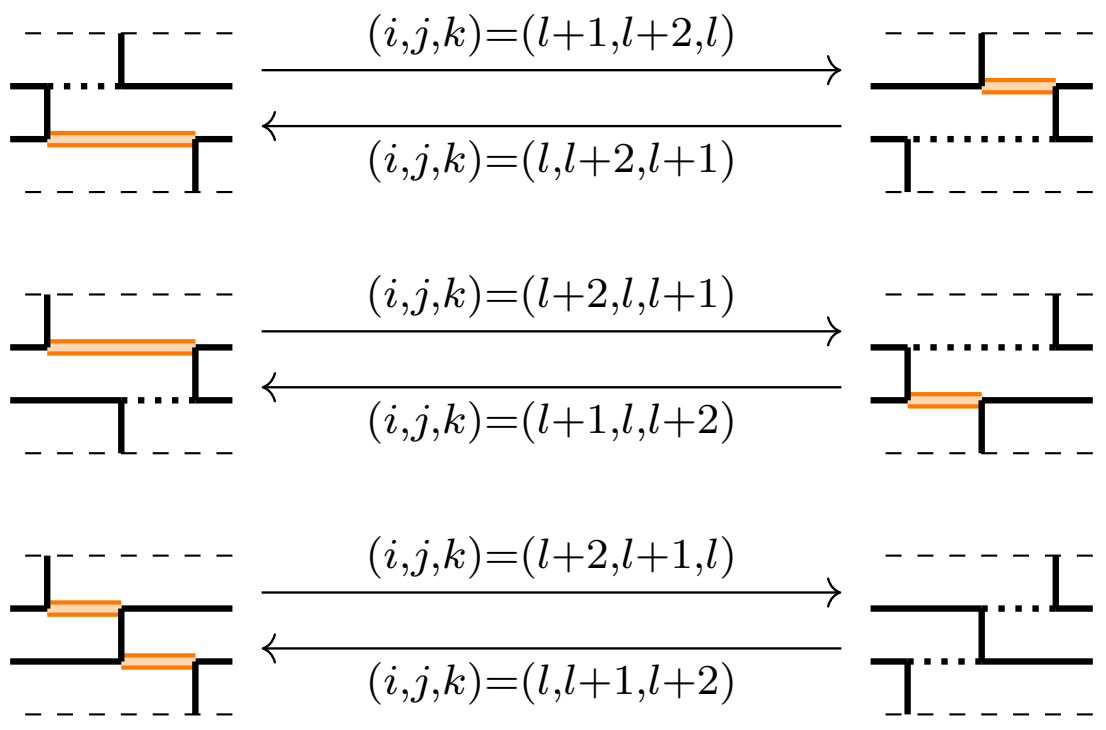

This concludes the proof for the relation (3.7).

A case-by-case analysis of the possible values for $(\lambda_i, \lambda_{i+1}, \lambda_{i+2})$ shows that except for four possibilities, the three diagrams in (3.8) are all set to zero by the Schur quotient. In the remaining cases, exactly one summand on the left-hand side is set to zero by the Schur quotient. If $(\lambda_i, \lambda_{i+1}, \lambda_{i+2}) = (2, 0, 0)$,

then

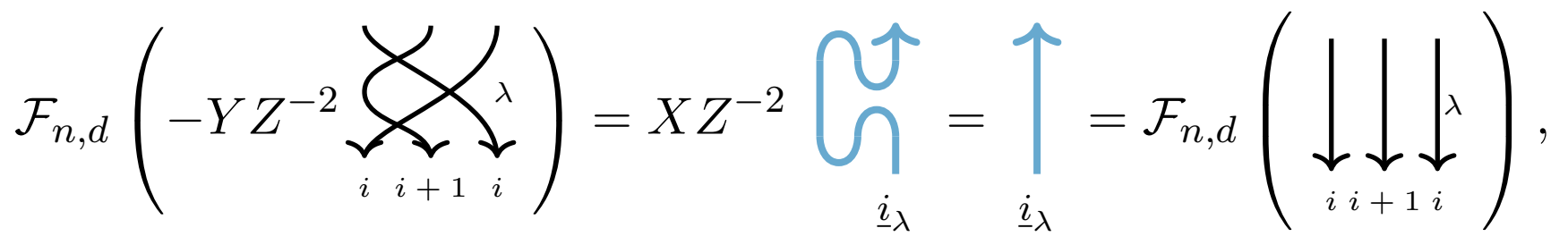

using $\Gamma_{\lambda-\alpha_i}(i+1) = \Gamma_\lambda(i)$. The case $(\lambda_i, \lambda_{i+1}, \lambda_{i+2}) = (2,0,1)$ is similar. If $(\lambda_i, \lambda_{i+1}, \lambda_{i+2}) = (2,1,0)$, then using $\Gamma_{\lambda-\alpha_{i+1}}(i) = (-XY)\Gamma_\lambda(i+1)$. The case $(\lambda_i, \lambda_{i+1}, \lambda_{i+2}) = (2,1,1)$ is similar.

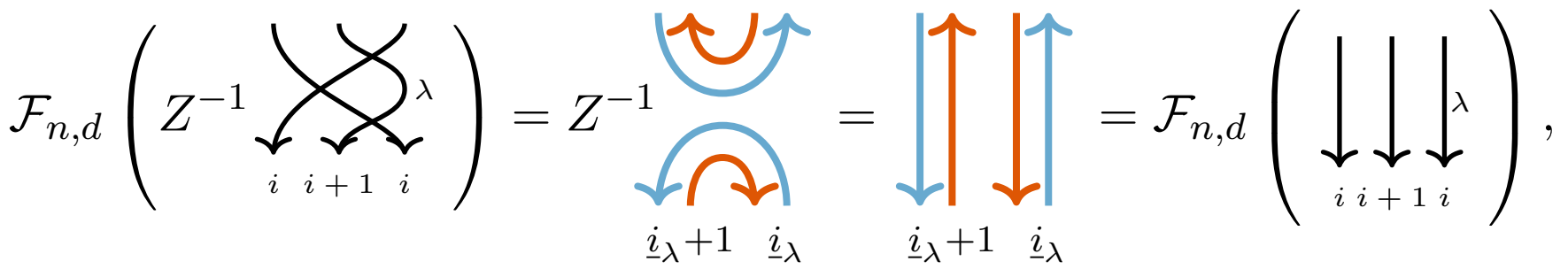

A similar analysis can be done for the relation (3.9), giving four non-trivial cases, with only two computations needed. We depict respectively the cases $(\lambda_{i-1}, \lambda_i, \lambda_{i+1}) = (1,2,0)$ and $(\lambda_{i-1}, \lambda_i, \lambda_{i+1}) = (1,1,0)$:

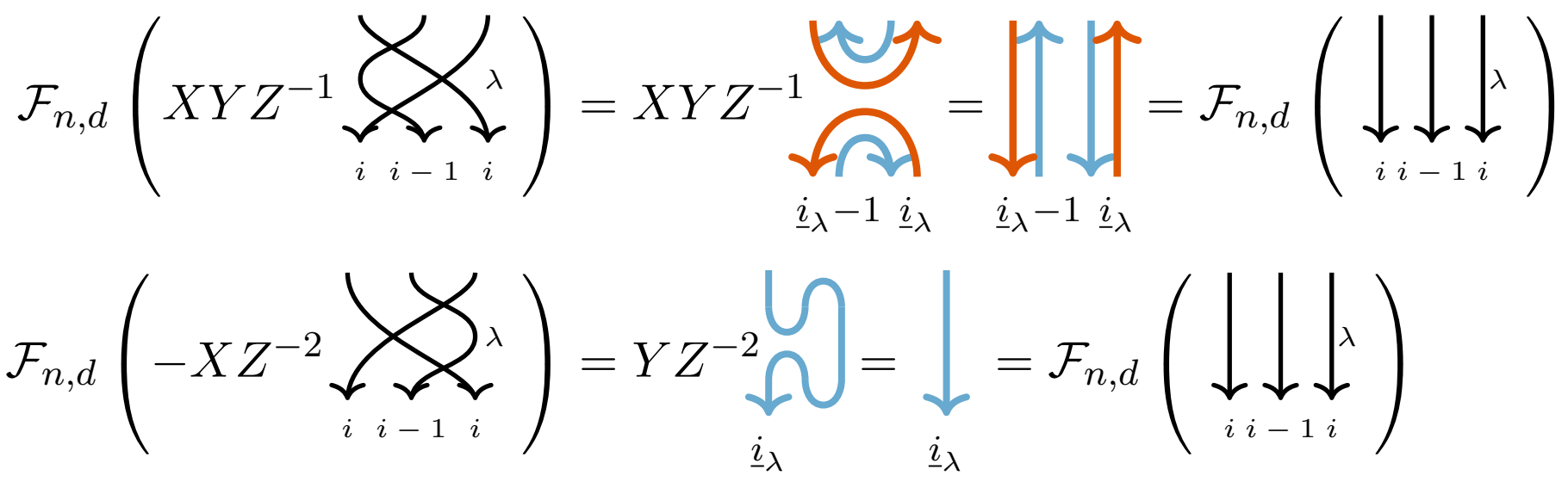

using Lemma 1.4.4, and $\Gamma_{\lambda-\alpha_i}(i) = (-XY)\Gamma_\lambda(i)$ and $\Gamma_{\lambda-\alpha_{i-1}}(i) = \Gamma_\lambda(i)$, respectively. For relation (3.10), it follows from the zigzag relations for foams:

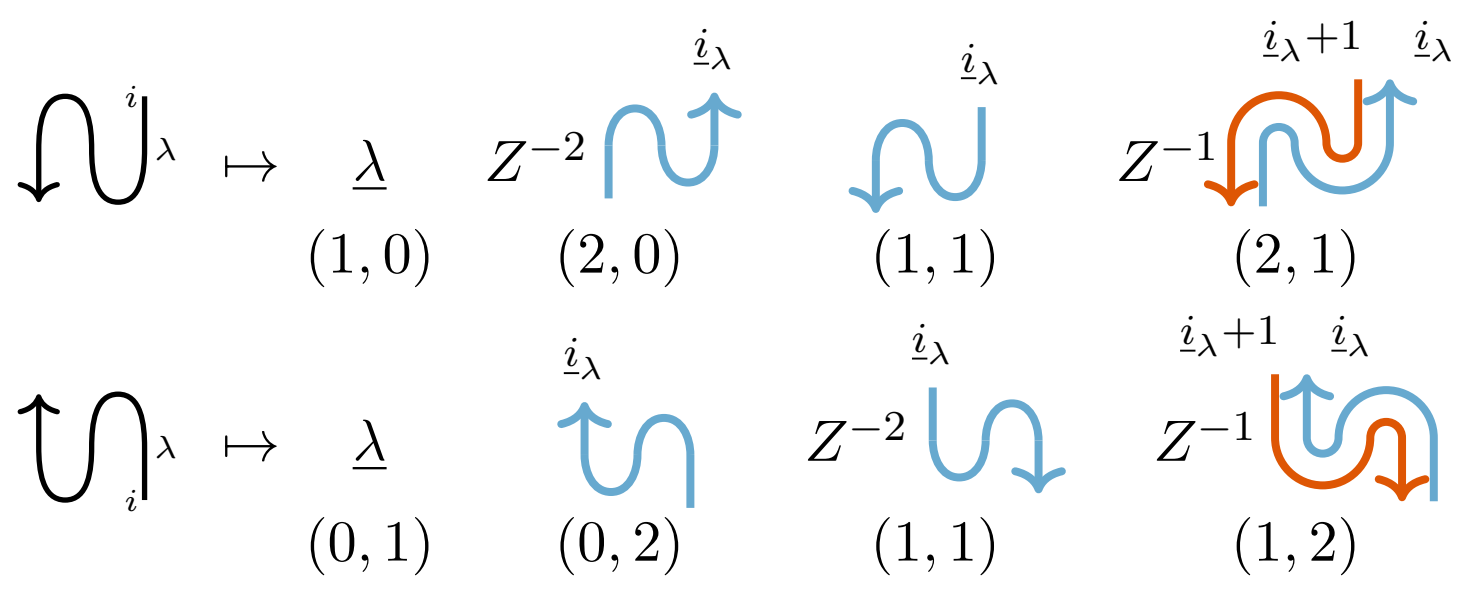

It remains to check that $\mathcal{F}_{n,d}$ preserves three inverse axioms (3.11), (3.12) and (3.13). Note that this is independent of the coefficients in the defining relations of $\mathbf{GFoam}_d$, as those inverse axioms only impose existence. For relation (3.11), it follows from graded isotopies and zigzag relations. For $\overline{\lambda}_i = 0$, one checks that the image of the leftward monochromatic crossing is the identity, up to multiplication by an invertible scalar. Otherwise, it is zero by the Schur quotient. For $\overline{\lambda}_i = \pm 1$, The image of the leftward cup (resp. cap) is invertible due to the squeezing relation. Finally, for $\overline{\lambda}_i = \pm 2$ this is precisely the (vertical) neck-cutting relation. □

In [183], the second author defined a "thick calculus" for the negative half of the super version of the graded 2-Schur $\mathcal{GS}_{n,d}$. This amounts to working in a sub-category of the Karoubi envelope. One can similarly define a thick calculus for the full graded 2-Schur, defining a graded-2-category $\check{\mathcal{GS}}_{n,d}$. Then, following the line of the proof of the analogous result in the non-graded case as given by Queffelec and Rose [156, Theorem 3.9], one can show that the foamation 2-functor factors through $\check{\mathcal{GS}}_{n,d}$. The inclusion

$$\check{\mathcal{GS}}_{n,d} \hookrightarrow \check{\mathcal{GS}}_{n+1,d}$$

is defined on objects as $(\lambda_1, \dots, \lambda_n) \mapsto (\lambda_1, \dots, \lambda_n, 0)$ and similarly for 1-morphisms and 2-morphisms. Finally, following the line of proof of [156, Proposition 3.22][2] in the $\mathfrak{gl}_2$-case, one can show the following proposition. We leave the details to the reader.

**Proposition 3.3.2.** *We have the following equivalence of $(\mathbb{Z}^2, \mu)$-graded-2-categories:*

$$\operatorname*{colim}_{n\in\mathbb{N}} \left( \dots \hookrightarrow \check{\mathcal{GS}}_{n,d} \hookrightarrow \check{\mathcal{GS}}_{n+1,d} \hookrightarrow \dots \right) \cong \mathbf{GFoam}_d.$$

## 3.4 The categorification theorem

Recall the notion of quantum Grothendieck ring from subsection 1.1.2. As claimed, the graded 2-Schur algebra categorifies the Schur algebra of level two:

[2] Or rather, the analogue of this proposition when one imposes (3.2) and dot annihilation respectively.

**Theorem 3.4.1.** *The graded 2-Schur algebra categorifies the q-Schur algebra of level two::*

$$K_0(\mathcal{GS}_{n,d})|_q \cong \dot{S}_{n,d},$$

*where the isomorphism is an isomorphism of* $\mathbb{Z}[q, q^{-1}]$*-linear categories.*

*Proof of Theorem 3.4.1.* Relations (3.1) become 2-isomorphisms in $\mathcal{GS}_{n,d}$. The first relation is categorified by the invertibility axioms (3.11), (3.12) and (3.13). The second relation is categorified by (3.3) in the case $|i-j| > 1$. Finally, the third relation is categorified by the analogue of (3.3), case $|i-j| > 1$, for upward strands. We define the *upward crossing* as:

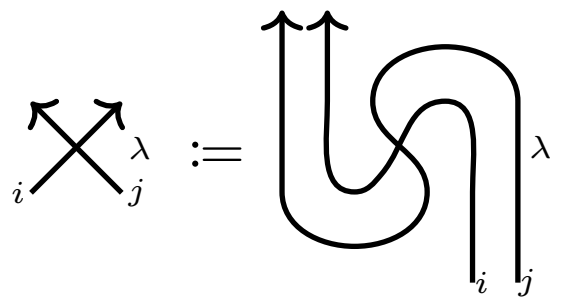

It then follows from adjunction relations (3.10) and (3.3), case $|i-j| > 1$, that:

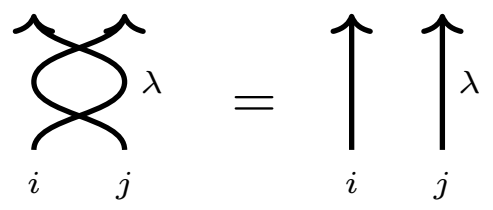

if $|i-j| > 1$.

This implies that there exists a $\mathbb{Z}[q, q^{-1}]$-linear functor $\dot{S}_{n,d} \to K_0(\mathcal{GS}^{\oplus,\mathrm{cl}}_{n,d})$, full and surjective, fitting into the following commutative diagram:

$$\begin{array}{ccc} \dot{S}_{n,d} & \longrightarrow & K_0(\mathcal{GS}^{\oplus,\mathrm{cl}}_{n,d}) \\ {\scriptstyle F_{n,d}}\downarrow & \circlearrowleft & \downarrow{\scriptstyle K_0(\mathcal{F}_{n,d})} \\ \mathbf{Web}_d & \xrightarrow{\cong} & K_0(\mathbf{GFoam}^{\oplus,\mathrm{cl}}_d) \end{array}$$

The bottom arrow is an isomorphism thanks to Theorem 1.7.1 and the left arrow is faithful thanks to Lemma 3.1.2. We conclude that the top arrow is also faithful, and hence it is an isomorphism. □

# 4

# Chain complexes in graded-monoidal categories

This chapter introduces a notion of the tensor product for chain complexes in a given graded-monoidal category (Definition 1.1.4); more precisely, we restrict our study to *homogeneous polycomplexes.* The special case of super-2-categories first appeared in the author's Master thesis [169]. Crucially, we then show that this tensor product leaves homotopy classes invariant. Although we work in the context of graded-monoidal categories for simplicity, all definitions and results extend in a straightforward way to graded-2-categories.

This chapter is not conceptually difficult, but it is technical. We refer the reader to subsection 2.1.1 for a minimalistic version, sufficient for the purpose of chapter 2.

*Notation* 4.0.1. Fix $n \in \mathbb{N}$. Recall Notation 2.1.1. For $\boldsymbol{r} \in \mathbb{Z}^n$, we write $|\boldsymbol{r}| := r_1 + \ldots + r_n$. Fix $\Bbbk$, $G$ and $\mu$ as in section 1.1, and fix a $(G,\mu)$-graded-monoidal category $\mathcal{C}$. To reduce clutter, we often abuse notation and write $f$ instead of $\deg f$, where $f$ is a morphism of $\mathcal{C}$. The distinction should be clear by the context. We also write $*$ instead of $\mu$. For instance, for $f$, $g$ and $h$ morphisms, we may write $(f+g)*h$ for $\mu(\deg f + \deg g, \deg h)$. We sometimes use the subscripted equal signs $=_{\Bbbk}$ and $=_G$ to emphasize where equality holds.

**Assumption.** To simply this already technical chapter, we assume that $\mu$ is symmetric throughout (see Definition 1.1.7).

## 4.1 A graded Koszul rule

Recall that for a graph such as $\mathbb{Z}^n$, a *cycle* is an oriented loop of edges, and a 1-cochain on a graph is said to be a *1-cocycle* if it is zero on all cycles. For $\mathbb{Z}^n$, a 1-cochain is a 1-cocycle if and only if it is zero on every square. A graph 1-cocycle is always the boundary of a graph 0-cochain.

Recall that if $V = \oplus_{g\in G} V_g$ is a $G$-graded $\Bbbk$-module, an element $v \in V$ is said to be *homogeneous* if $v \in V_g$ for some $g \in G$. If moreover $v \neq 0$, then $v$ has a well-defined degree $\deg v = g$.

**Definition 4.1.1.** *A* homogeneous $n$-polycomplex $\mathbb{A} = (A, \alpha, \psi_{\mathbb{A}})$ *is the data of:*

*(i) a family $A := (A^{\boldsymbol{r}})_{\boldsymbol{r}\in\mathbb{Z}^n}$ of objects $A^{\boldsymbol{r}} \in \mathcal{C}$,*

*(ii) a family $\alpha := (\alpha_i^{\boldsymbol{r}})_{\boldsymbol{r}\in\mathbb{Z}^n, i\in\mathcal{I}}$ of* homogeneous *morphisms $\alpha_i^{\boldsymbol{r}} \colon A^{\boldsymbol{r}} \to A^{\boldsymbol{r}+e_i}$, such that $\alpha_i^{\boldsymbol{r}+e_i} \circ \alpha_i^{\boldsymbol{r}} = 0$ for all $i \in \mathcal{I}$, and such that each square anti-commutes:*
$$\alpha_j^{\boldsymbol{r}+e_i} \circ \alpha_i^{\boldsymbol{r}} = -\alpha_i^{\boldsymbol{r}+e_j} \circ \alpha_j^{\boldsymbol{r}},$$

*(iii) a $G$-valued 1-cocycle $\psi_{\mathbb{A}}$ on $\mathbb{Z}^n$ such that $\psi_{\mathbb{A}}(\boldsymbol{r} \to \boldsymbol{r} + \boldsymbol{e}_i) = \deg \alpha_i^{\boldsymbol{r}}$ whenever $\alpha_i^{\boldsymbol{r}} \neq 0$.*

Thanks to condition (ii), we associate to $\mathbb{A}$ a chain complex, its *total complex*, denoted $(\mathrm{Tot}(\mathbb{A}), \mathrm{Tot}(\alpha))$ and given by
$$\mathrm{Tot}(\mathbb{A})_t := \bigoplus_{\boldsymbol{r}\in\mathbb{Z}^n, |\boldsymbol{r}|=t} A^{\boldsymbol{r}} \quad \text{and} \quad \mathrm{Tot}(\alpha)|_{A^{\boldsymbol{r}}} := \sum_{1\leq i\leq n} \alpha_i^{\boldsymbol{r}}.$$

In what follows, we will not distinguish between a homogeneous polycomplex and its associated total complex. In particular, the notions of chain map and chain homotopy are the usual ones.

If a square in $\mathbb{A}$ is non-zero in the sense that
$$\alpha_j^{\boldsymbol{r}+e_i} \circ \alpha_i^{\boldsymbol{r}} = -\alpha_i^{\boldsymbol{r}+e_j} \circ \alpha_j^{\boldsymbol{r}} \neq 0,$$

then each of the four maps involved are non-zero, they have a well-defined degree, and the following condition holds:
$$\deg \alpha_j^{\boldsymbol{r}+e_i} + \deg \alpha_i^{\boldsymbol{r}} = \deg \alpha_i^{\boldsymbol{r}+e_j} + \deg \alpha_j^{\boldsymbol{r}}. \tag{4.1}$$

In other words, the partially-defined 1-cochain $\deg$ is a 1-cocycle on non-zero squares. Condition (iii) states the existence of a 1-cochain $\psi_A$ that both extends $\deg$ to zero maps and the 1-cocycle condition to zero squares. In particular, a square could be zero albeit having four non-zero edges: in which case, although the equation (4.1) makes sense, it may not hold. Condition (iii) ensures that it does.

For each vertex $\boldsymbol{r} \in \mathbb{Z}^n$, let $p$ be a path in $\mathbb{Z}^n$ from $\mathbf{0}$ to $\boldsymbol{r}$. Since $\psi_{\mathbb{A}}$ is a cocycle, the value $|\alpha|(\boldsymbol{r}) := \psi_{\mathbb{A}}(p)$ does not depend on the choice of path. Most importantly, it verifies the following:

$$|\alpha|(\boldsymbol{r}+e_i) - |\alpha|(\boldsymbol{r}) =_G \psi_{\mathbb{A}}(\boldsymbol{r} \to \boldsymbol{r}+\boldsymbol{e}_i) =_G \deg \alpha_i^{\boldsymbol{r}}, \tag{4.2}$$

where the last equality holds whenever it makes sense, that is whenever $\alpha_i^{\boldsymbol{r}} \neq 0$. We now introduce the notion of tensor product of homogeneous polycomplexes, after setting some further notations.

*Notation* 4.1.2. Recall our notations from Notation 4.0.1 and Definition 4.1.1. We introduce analogous notations for a homogeneous $m$-polycomplex $\mathbb{B} = (B, \beta, \psi_{\mathbb{B}})$: we let $\mathcal{J} = \{1, \ldots, m\}$ denotes the set of indices and use generically the letter $j$ for an index $j \in \mathcal{J}$ and $\boldsymbol{s} = (s_1, \ldots, s_m)$ for a vertex in $\mathbb{Z}^m$. Finally, we generically use the letter $k$ for an index $k \in \mathcal{I} \sqcup \mathcal{J}$ and the notation $(\boldsymbol{r}, \boldsymbol{s})$ for a vertex in $\mathbb{Z}^{n+m}$.

**Definition 4.1.3.** *Let $\epsilon$ be an $\Bbbk^\times$-valued 1-cochain on $\mathbb{Z}^{n+m}$. With the notations above, the $\epsilon$-*tensor product *of $\mathbb{A}$ and $\mathbb{B}$, denoted $(\mathbb{A} \otimes \mathbb{B})(\epsilon)$, is the following data:*

*(i) family $A \otimes B = ((A \otimes B)^{(\boldsymbol{r},\boldsymbol{s})})$ of objects $(A \otimes B)^{(\boldsymbol{r},\boldsymbol{s})} := A^{\boldsymbol{r}} \otimes B^{\boldsymbol{s}}$,*

*(ii) family $\alpha \otimes \beta = ((\alpha \otimes \beta)_k^{(\boldsymbol{r},\boldsymbol{s})})$ of homogeneous morphisms*

$$(\alpha \otimes \beta)_k^{(\boldsymbol{r},\boldsymbol{s})} := \begin{cases} \epsilon(e)(\alpha_i^{\boldsymbol{r}} \otimes \mathrm{id}_{B^{\boldsymbol{s}}}) & k = i \in \mathcal{I}, \\ \epsilon(e)(\mathrm{id}_{A^{\boldsymbol{r}}} \otimes \beta_j^{\boldsymbol{s}}) & k = j \in \mathcal{J}, \end{cases}$$

*where $e$ denotes the edge $(\boldsymbol{r}, \boldsymbol{s}) \to (\boldsymbol{r}, \boldsymbol{s}) + \boldsymbol{e}_k$,*

*(iii) 1-cocycle $\psi_{\mathbb{A}} \otimes \psi_{\mathbb{B}}$ on $\mathbb{Z}^{n+m}$ given on $e = (\boldsymbol{r}, \boldsymbol{s}) \to (\boldsymbol{r}, \boldsymbol{s}) + \boldsymbol{e}_k$ by*

$$(\psi_{\mathbb{A}} \otimes \psi_{\mathbb{B}})(e) := \begin{cases} \psi_{\mathbb{A}}(\boldsymbol{r} \to \boldsymbol{r} + \boldsymbol{e}_i) & k = i \in \mathcal{I}, \\ \psi_{\mathbb{B}}(\boldsymbol{s} \to \boldsymbol{s} + \boldsymbol{e}_j) & k = j \in \mathcal{J}. \end{cases}$$

For $(\mathbb{A} \otimes \mathbb{B})(\epsilon)$ to define a homogeneous $(n+m)$-polycomplex, $\epsilon$ needs to be such that squares anti-commute. This is encapsulated in the following lemma, where $\square_{k,l}^{(\boldsymbol{r},\boldsymbol{s})}$ with $k < l$ denotes the following oriented square:

$$\begin{array}{ccc} (\boldsymbol{r},\boldsymbol{s}) & \longrightarrow & (\boldsymbol{r},\boldsymbol{s}) + \boldsymbol{e}_k \\ \big\downarrow & \circlearrowright & \big\downarrow \\ (\boldsymbol{r},\boldsymbol{s}) + \boldsymbol{e}_l & \longrightarrow & (\boldsymbol{r},\boldsymbol{s}) + \boldsymbol{e}_k + \boldsymbol{e}_l \end{array}$$

**Lemma 4.1.4.** *Say that a 1-cochain $\epsilon$ on $\mathbb{Z}^{n+m}$ is* compatible *if*

$$\partial\epsilon(\square_{k,l}^{(\boldsymbol{r},\boldsymbol{s})}) =$$
$$\begin{cases} 0 & \text{if } k,l \in \mathcal{I} \text{ or } k,l \in \mathcal{J}, \\ -\psi_{\mathbb{A}}(\boldsymbol{r} \to \boldsymbol{r} + \boldsymbol{e}_i) * \psi_{\mathbb{B}}(\boldsymbol{s} \to \boldsymbol{s} + \boldsymbol{e}_j) & \text{if } k = i \in \mathcal{I} \text{ and } l = j \in \mathcal{J}. \end{cases}$$

*If $\epsilon$ is compatible, then $(\mathbb{A} \otimes \mathbb{B})(\epsilon)$ defines a homogeneous $(n+m)$-polycomplex.*

*Proof.* The case $k,l \in \mathcal{I}$ or $k,l \in \mathcal{J}$ is clear. In the case $k = i \in \mathcal{I}$ and $l = j \in \mathcal{J}$, the square $\square_{k,l}^{(\boldsymbol{r},\boldsymbol{s})}$ has the following form:

$$\begin{array}{ccc} (\boldsymbol{r},\boldsymbol{s}) & \xrightarrow{\alpha_i^{\boldsymbol{r}} \otimes \mathrm{id}_{B^{\boldsymbol{s}}}} & (\boldsymbol{r},\boldsymbol{s}) + \boldsymbol{e}_k \\ {\scriptstyle \mathrm{id}_{A^{\boldsymbol{r}}} \otimes \beta_j^{\boldsymbol{s}}}\big\downarrow & \circlearrowright & \big\downarrow{\scriptstyle \mathrm{id}_{A^{\boldsymbol{r}+\boldsymbol{e}_k}} \otimes \beta_j^{\boldsymbol{s}}} \\ (\boldsymbol{r},\boldsymbol{s}) + \boldsymbol{e}_l & \xrightarrow[\alpha_i^{\boldsymbol{r}} \otimes \mathrm{id}_{B^{\boldsymbol{s}+\boldsymbol{e}_l}}]{} & (\boldsymbol{r},\boldsymbol{s}) + \boldsymbol{e}_k + \boldsymbol{e}_l \end{array}$$

The morphism corresponding to the path

$$(\boldsymbol{r},\boldsymbol{s}) \to (\boldsymbol{r},\boldsymbol{s}) + \boldsymbol{e}_k \to (\boldsymbol{r},\boldsymbol{s}) + \boldsymbol{e}_k + \boldsymbol{e}_l$$

is the morphism $(\beta_j^{\boldsymbol{s}} * \alpha_i^{\boldsymbol{r}})\alpha_i^{\boldsymbol{r}} \otimes \beta_j^{\boldsymbol{s}}$, while the morphism corresponding to the other path is $\alpha_i^{\boldsymbol{r}} \otimes \beta_j^{\boldsymbol{s}}$. If the square is zero, it automatically anti-commutes. Otherwise, it is sufficient to have

$$\partial\epsilon(\square_{k,l}^{(\boldsymbol{r},\boldsymbol{s})}) =_{\Bbbk} -(\beta_j^{\boldsymbol{s}} * \alpha_i^{\boldsymbol{r}})^{-1} =_{\Bbbk} -\psi_{\mathbb{A}}(\boldsymbol{r} \to \boldsymbol{r} + \boldsymbol{e}_i) * \psi_{\mathbb{B}}(\boldsymbol{s} \to \boldsymbol{s} + \boldsymbol{e}_j),$$

using symmetry of $\mu$ in the second equality. $\square$

*Remark* 4.1.5. Note that the compatibility condition is a sufficient condition to ensure that all squares anti-commute, but a priori not a necessary one. Indeed, there might be some liberty on squares that are zero. In particular, we could have $\alpha_i^{\boldsymbol{r}} \otimes \beta_j^{\boldsymbol{s}} = 0$ even if $\alpha_i^{\boldsymbol{r}} \neq 0$ and $\beta_j^{\boldsymbol{s}} \neq 0$.

Definition 4.1.3 defines a tensor product of homogeneous polycomplexes, *provided* that a compatible 1-cochain exists. Definition 4.1.6 below ensures existence.

**Definition 4.1.6** (graded Koszul rule)**.** *The* standard 1-cochain $\epsilon_{\mathbb{A}\otimes\mathbb{B}}$ *is the compatible 1-cochain on* $H^{n+m}$ *given on* $e = (\boldsymbol{r}, \boldsymbol{s}) \to (\boldsymbol{r}, \boldsymbol{s}) + \boldsymbol{e}_k$ *by*

$$\epsilon_{\mathbb{A}\otimes\mathbb{B}}(e) = \begin{cases} 0 & k = i \in \mathcal{I}, \\ (-1)^{|\boldsymbol{r}|}\, |\alpha|\,(\boldsymbol{r}) * \psi_{\mathbb{B}}(\boldsymbol{s} \to \boldsymbol{s} + \boldsymbol{e}_j) & k = j \in \mathcal{J}. \end{cases}$$

Note that if $\psi_{\mathbb{A}} = \psi_{\mathbb{B}} = 1$, then $|\alpha|\,(\boldsymbol{r}) * \psi_{\mathbb{B}}(\boldsymbol{s} \to \boldsymbol{s} + \boldsymbol{e}_j) = 1$ and we recover the usual Koszul rule. Checking that $\epsilon_{\mathbb{A}\otimes\mathbb{B}}$ is compatible is straightforward:

*Proof that the standard 1-cochain is compatible.* Consider the square $\Box_{k,l}^{(\boldsymbol{r},\boldsymbol{s})}$ as above. If $k, l \in \mathcal{I}$, the compatibility condition follows directly, and if $k, l \in \mathcal{J}$, it follows from the fact that $\psi_{\mathbb{B}}$ is a cocycle. Finally, if $k \in \mathcal{I}$ and $l \in \mathcal{J}$, we get (using property (4.2)):

$$\begin{aligned} \partial\epsilon(\Box_{k,l}^{(\boldsymbol{r},\boldsymbol{s})}) &= \Big((-1)^{|\boldsymbol{r}|+1}\, |\alpha|\,(\boldsymbol{r} + \boldsymbol{e}_k) * \psi_{\mathbb{B}}(\boldsymbol{s} \to \boldsymbol{s} + \boldsymbol{e}_j)\Big) \\ &\qquad \Big((-1)^{|\boldsymbol{r}|}\, |\alpha|\,(\boldsymbol{r}) * \psi_{\mathbb{B}}(\boldsymbol{s} \to \boldsymbol{s} + \boldsymbol{e}_j)^{-1}\Big) \\ &= -\big(\, |\alpha|\,(\boldsymbol{r} + \boldsymbol{e}_k) - |\alpha|\,(\boldsymbol{r})\big) * \psi_{\mathbb{B}}(\boldsymbol{s} \to \boldsymbol{s} + \boldsymbol{e}_j) \\ &= -\psi_{\mathbb{A}}(\boldsymbol{r} \to \boldsymbol{r} + \boldsymbol{e}_i) * \psi_{\mathbb{B}}(\boldsymbol{s} \to \boldsymbol{s} + \boldsymbol{e}_j). \end{aligned}$$ ☐

We next check that this choice of compatible 1-cochain is essentially unique, *among compatible cochains* (see Remark 4.1.5).[1]

**Lemma 4.1.7.** *Let* $\epsilon$ *and* $\epsilon'$ *be two compatible 1-cochains. Then* $(\mathbb{A} \otimes \mathbb{B})(\epsilon)$ *and* $(\mathbb{A} \otimes \mathbb{B})(\epsilon')$ *are isomorphic as chain complexes.*

*Proof.* If $\epsilon$ and $\epsilon'$ are both compatible, then $\partial(\epsilon(\epsilon')^{-1})$ is zero on all squares, and hence is a 1-cocycle. Let $\varphi$ be a 0-cochain such that $\partial\varphi = \epsilon(\epsilon')^{-1}$. Consider the map

$$\Phi\colon (\mathbb{A} \otimes \mathbb{B})(\epsilon) \to (\mathbb{A} \otimes \mathbb{B})(\epsilon')$$

[1]The proof given below follows closely the proof of Lemma 2.2 in [149].

corresponding to multiplication by $\varphi(\boldsymbol{r}, \boldsymbol{s})$ when restricted to $(\boldsymbol{r}, \boldsymbol{s})$. It is straightforward to check that this map defines an isomorphism for the associated total complexes. □

From now on, we simply call $(\mathbb{A} \otimes \mathbb{B})(\epsilon_{\mathbb{A}\otimes\mathbb{B}})$ the *tensor product* of $\mathbb{A}$ and $\mathbb{B}$ and denote it $\mathbb{A} \otimes \mathbb{B} = (A \otimes B, \alpha \otimes \beta, \psi_{\mathbb{A}} \otimes \psi_{\mathbb{B}})$.

*Example* 4.1.8. A chain complex $\mathbb{C} = (C_\bullet, \partial_\bullet)$ whose differentials $\partial_t$ are all homogeneous is a homogeneous 1-polycomplex. An $n$-fold tensor product of homogeneous 1-polycomplexes is a homogeneous $n$-polycomplex. This is the specific construction used in subsection 2.1.1 to define covering $\mathfrak{gl}_2$-Khovanov homology.

We now come to the main result of this chapter:

**Theorem 4.1.9.** *Let $\mathbb{A}_1$, $\mathbb{A}_2$, $\mathbb{B}_1$ and $\mathbb{B}_2$ be homogeneous polycomplexes. Then:*

$$\mathbb{A}_1 \simeq \mathbb{B}_1 \quad \text{and} \quad \mathbb{A}_2 \simeq \mathbb{B}_2 \qquad \text{implies} \qquad \mathbb{A}_1 \otimes \mathbb{A}_2 \simeq \mathbb{B}_1 \otimes \mathbb{B}_2,$$

*where $\simeq$ denotes homotopy equivalence (in the usual sense).*

The idea of the proof of Theorem 4.1.9 is straightforward: define induced morphisms and induced homotopies on the tensor product. Precisely, Theorem 4.1.9 holds if the following holds:

(1) Given homogeneous polycomplexes $\mathbb{A}_k, \mathbb{B}_k$ and chain maps $F_k\colon \mathbb{A}_k \to \mathbb{B}_k$ $(k = 1, 2)$, there exists a chain complex $F_1 \otimes F_2\colon \mathbb{A}_1 \otimes \mathbb{A}_2 \to \mathbb{B}_1 \otimes \mathbb{B}_2$. This definition is such that $F_1 \otimes F_2 = \mathrm{Id}_{\mathbb{A}_1\otimes\mathbb{A}_2}$ if $F_1 = \mathrm{Id}_{\mathbb{A}_1}$ and $F_2 = \mathrm{Id}_{\mathbb{A}_2}$.

(2) Given homogeneous polycomplexes $\mathbb{A}_k, \mathbb{B}_k$, chain maps $F_k, G_k\colon \mathbb{A}_k \to \mathbb{B}_k$ and homotopies $H_k\colon F_k \to G_k$ $(k = 1, 2)$, there exists a homotopy $H_1 \otimes H_2\colon F_1 \otimes F_2 \to G_1 \otimes G_2$.

This is the content of section 4.2: part (1) is shown by Proposition 4.2.4, while part (2) is shown by Proposition 4.2.6.

### 4.1.1 Further directions

As it is unnecessary for our purpose, we did not investigate the categorical properties of the tensor product. For example, it is not clear whether it is functorial, or functorial up to homotopy. Other possible questions include:

- Can we extend the definition of the tensor product to higher structures? That is, can we define induced $n$-fold homotopies on the tensor product?
- Are the definitions of induced morphisms and induced homotopies unique up to higher structures, similarly to Lemma 4.1.7?
- What is the most general case where one can define a (sensible) tensor product on chain complexes in a graded-monoidal category? In particular, how far can we weaken condition (iii) in Definition 4.1.1?

## 4.2 Induced morphisms and homotopies on the tensor product

### 4.2.1 Preliminary remarks

Recall the conventions of Notation 4.0.1 and Notation 4.1.2. We denote the graded Koszul rule defined in Definition 4.1.6 as:

$$\epsilon_{\mathbb{A}\otimes\mathbb{B}}^{\boldsymbol{r},\boldsymbol{s},k} := \begin{cases} 0 & k = i \in \mathcal{I}, \\ (-1)^{|\boldsymbol{r}|}\, |\alpha|\, (\boldsymbol{r}) * \psi_{\mathbb{B}}(\boldsymbol{s} \to \boldsymbol{s} + \boldsymbol{e}_j) & k = j \in \mathcal{J}. \end{cases}$$

Note that in Definition 4.1.3, for $k = j \in \mathcal{J}$ the scalar $\epsilon_{\mathbb{A}\otimes\mathbb{B}}^{\boldsymbol{r},\boldsymbol{s},j}$ always appears in front of an expression involving $\beta_j^{\boldsymbol{s}}$. Thus, when doing computation with $\epsilon_{\mathbb{A}\otimes\mathbb{B}}^{\boldsymbol{r},\boldsymbol{s},j}$ in this context, we can safely assume that $\psi_{\mathbb{B}}(\boldsymbol{s} \to \boldsymbol{s} + \boldsymbol{e}_j) = \deg \beta_j^{\boldsymbol{s}}$. In the proofs below, we will always encounter similar situations. To simplify the exposition, we shall act as if $\deg f$ is always well-defined, and avoid the use of $\psi$ altogether.

We extend notations using subscripts to $\mathbb{A}_k = (A_k, \alpha_k, \psi_{\mathbb{A}_k})$ and $\mathbb{B}_k = (B_k, \beta_k, \psi_{\mathbb{B}_k})$ for $k = 1, 2$. From now on, we denote $\mathbb{A} = \mathbb{A}_1 \otimes \mathbb{A}_2$ and $\mathbb{B} = \mathbb{B}_1 \otimes \mathbb{B}_2$. This includes sets of indices $\mathcal{I} = \mathcal{I}_1 \sqcup \mathcal{I}_2$ and $\mathcal{J} = \mathcal{J}_1 \sqcup \mathcal{J}_2$. We also use different shortcuts, such as $\boldsymbol{r} = (\boldsymbol{r}_1, \boldsymbol{r}_2)$ and $\boldsymbol{s} = (\boldsymbol{s}_1, \boldsymbol{s}_2)$, that should be clear from the context.

Before entering the main proofs, we define a generic $\Bbbk^\times$-valued 0-cochain, satisfying generic properties. In fact, Definition 4.2.3 and Definition 4.2.5 solely depend on this generic definition, and all "degree-wise" computations follow from those generic properties.

**Definition 4.2.1.** *For* $\lambda_k = \left(\lambda_k^{r_k,s_k}\right)_{(r_k,s_k)\in A_k\otimes B_k}$ *a pair of* $G$*-valued 0-cochains* ($k=1,2$)*, let* $\epsilon^{r,s}_{\lambda_1,\lambda_2}$ *be the following* $\Bbbk^\times$*-valued 0-cochain:*

$$\epsilon^{r,s}_{\lambda_1,\lambda_2} := \Big( \big[\lambda_1^{r_1,s_1} + |\alpha_1|\,(\boldsymbol{r}_1) - |\beta_1|\,(\boldsymbol{s}_1)\big] * |\beta_2|\,(\boldsymbol{s}_2) \Big)^{-1} \Big( |\alpha_1|\,(\boldsymbol{r}_1) * \lambda_2^{r_2,s_2} \Big).$$

**Lemma 4.2.2.** *The generic* $\Bbbk^\times$*-valued 0-cochain defined above satisfies the following:*

*(i) Assume* $\nu_1 :=_G \lambda_1^{r_1+e_{i_1},s_1} + \alpha_{i_1}^{r_1} =_G \lambda_1^{r_1,s_1-e_{j_1}} + \beta_{j_1}^{s_1-e_{j_1}}$ *holds for all* $i_1 \in \mathcal{I}_1$ *and* $j_1 \in \mathcal{J}_1$*; in particular, the* $G$*-valued 0-cochain* $\nu_1$ *is defined independently of the choice of* $i_1$ *and* $j_1$*. Then the following identities hold, for all* $i_1 \in \mathcal{I}_1$ *and* $j_1 \in \mathcal{J}_1$*:*

$$\epsilon^{r,s}_{\nu_1,\lambda_2} =_\Bbbk \epsilon^{r+e_{i_1},s}_{\lambda_1,\lambda_2} \left(\lambda_2^{r_2,s_2} * \alpha_{i_1}^{r_1}\right) =_\Bbbk \epsilon^{r,s-e_{j_1}}_{\lambda_1,\lambda_2}.$$

*(ii) Assume* $\nu_2 :=_G \lambda_2^{r_2+e_{i_2},s_2} + \alpha_{i_2}^{r_2} =_G \lambda_2^{r_2,s_2-e_{j_2}} + \beta_{j_2}^{s_2-e_{j_2}}$ *holds for all* $i_2 \in \mathcal{I}_2$ *and* $j_2 \in \mathcal{J}_2$*; in particular, the* $G$*-valued 0-cochain* $\nu_2$ *is defined independently of the choice of* $i_2$ *and* $j_2$*. Then the following identities hold, for all* $i_2 \in \mathcal{I}_2$ *and* $j_2 \in \mathcal{J}_2$*:*

$$\begin{aligned} \epsilon^{r,s}_{\lambda_1,\nu_2} &=_\Bbbk (-1)^{|r_1|} \epsilon^{r_1,r_2,j_2}_{A_1\otimes A_2} \epsilon^{r+e_{i_2},s}_{\lambda_1,\lambda_2} \\ &=_\Bbbk (-1)^{|s_1|} \epsilon^{s_1,s_2-e_{j_2},j_2}_{B_1\otimes B_2} \epsilon^{r,s-e_{j_2}}_{\lambda_1,\lambda_2} \left(\beta_{j_2}^{s_2-e_{j_2}} * \lambda_1^{r_1,s_1}\right). \end{aligned}$$

*Proof.* It follows from direct computations, using relation (4.2). Symmetry of $\mu$ gives the expressions as in the lemma.

$$\begin{aligned} \epsilon^{r+e_{i_1},s}_{\lambda_1,\lambda_2} &=_\Bbbk \Big( \big[\lambda_1^{r_1+e_{i_1},s_1} + \alpha_{i_1}^{r_1} + |\alpha_1|\,(\boldsymbol{r}_1) - |\beta_1|\,(\boldsymbol{s}_1)\big] * |\beta_2|\,(\boldsymbol{s}_2) \Big)^{-1} \\ &\qquad \cdot \Big( \big[\,|\alpha_1|\,(\boldsymbol{r}_1) + \alpha_{i_1}^{r_1}\big] * \lambda_2^{r_2,s_2} \Big) \\ &=_\Bbbk \epsilon^{r,s}_{\nu_1,\lambda_2} \left(\alpha_{i_1}^{r_1} * \lambda_2^{r_2,s_2}\right) \end{aligned}$$

$$\begin{aligned} \epsilon^{r,s-e_{j_1}}_{\lambda_1,\lambda_2} &=_\Bbbk \Big( \big[\lambda_1^{r_1,s_1-e_{j_1}} + \beta_{j_1}^{s_1-e_{j_1}} + |\alpha_1|\,(\boldsymbol{r}_1) - |\beta_1|\,(\boldsymbol{s}_1)\big] * |\beta_2|\,(\boldsymbol{s}_2) \Big)^{-1} \\ &\qquad \cdot \Big( |\alpha_1|\,(\boldsymbol{r}_1) * \lambda_2^{r_2,s_2} \Big) \\ &=_\Bbbk \epsilon^{r,s}_{\nu_1,\lambda_2} \end{aligned}$$

$$
\begin{aligned}
&(-1)^{|r_1|}\epsilon_{A_1\otimes A_2}^{r_1,r_2,j_2}\epsilon_{\lambda_1,\lambda_2}^{r+e_{i_2},s} \\
&=_{\Bbbk} \left(|\alpha_1|(r_1) * \alpha_{i_2}^{r_2,s_2}\right)\Big(\big[\lambda_1^{r_1,s_1} + |\alpha_1|(r_1) - |\beta_1|(s_1)\big] * |\beta_2|(s_2)\Big)^{-1} \\
&\qquad\qquad \cdot\Big(|\alpha_1|(r_1) * \lambda_2^{r_2+e_{i_2},s_2}\Big) \\
&=_{\Bbbk} \Big(\big[\lambda_1^{r_1,s_1} + |\alpha_1|(r_1) - |\beta_1|(s_1)\big] * |\beta_2|(s_2)\Big)^{-1} \\
&\qquad\qquad \cdot\Big(|\alpha_1|(r_1) * \big[\lambda_2^{r_2+e_{i_2},s_2} + \alpha_{i_2}^{r_2,s_2}\big]\Big) \\
&=_{\Bbbk} \epsilon_{\lambda_1,\nu_2}^{r,s}
\end{aligned}
$$

$$
\begin{aligned}
&(-1)^{|s_1|}\epsilon_{B_1\otimes B_2}^{s_1,s_2-e_{j_2},j_2}\epsilon_{\lambda_1,\lambda_2}^{r,s-e_{j_2}}\left(\lambda_1^{r_1,s_1} * \beta_{j_2}^{s_2-e_{j_2}}\right)^{-1} \\
&=_{\Bbbk} \left(|\beta_1|(s_1) * \beta_{j_2}^{r_2,s_2-e_{j_2}}\right) \\
&\qquad\qquad \cdot\Big(\big[\lambda_1^{r_1,s_1} + |\alpha_1|(r_1) - |\beta_1|(s_1)\big] * |\beta_2|(s_2 - e_{j_2})\Big)^{-1} \\
&\qquad\qquad\quad \cdot\Big(|\alpha_1|(r_1) * \lambda_2^{r_2,s_2-e_{j_2}}\Big)\Big(\lambda_1^{r_1,s_1} * \beta_{j_2}^{s_2-e_{j_2}}\Big)^{-1} \\
&=_{\Bbbk} \Big(\big[\lambda_1^{r_1,s_1} + |\alpha_1|(r_1) - |\beta_1|(s_1)\big] * |\beta_2|(s_2)\Big)^{-1} \\
&\qquad\qquad \cdot\Big(|\alpha_1|(r_1) * \big[\lambda_2^{r_2,s_2-e_{j_2}} + \beta_{j_2}^{s_2-e_{j_2}}\big]\Big) \\
&=_{\Bbbk} \epsilon_{\lambda_1,\nu_2}^{r,s}. \qquad \square
\end{aligned}
$$

### 4.2.2 Induced morphism on tensor product

Unpacking the standard[2] definition of a chain map between chain complexes, a chain map $F\colon \mathbb{A} \to \mathbb{B}$ is a set of morphisms $F^{r,s}\colon A^r \to B^s$ for all $|r| = |s|$ such that $\beta \circ F = F \circ \alpha$, that is:

$$\sum_{j\in\mathcal{J}} \beta_j^{s-e_j} \circ F^{r,s-e_j} = \sum_{i\in\mathcal{I}} F^{r+e_i,s} \circ \alpha_i^r.$$

**Definition 4.2.3.** *Let $F_k\colon \mathbb{A}_k \to \mathbb{B}_k$ be chain maps for each $k = 1, 2$. Their* induced morphism *$F = F_1 \otimes F_2$ is the chain map $F\colon \mathbb{A}_1 \otimes \mathbb{A}_2 \to \mathbb{B}_1 \otimes \mathbb{B}_2$ given by the data:*

$$F^{r,s} = \epsilon_{F_1,F_2}^{r,s} F_1^{r_1,s_1} \otimes F_2^{r_2,s_2},$$

*where $\epsilon_{F_1,F_2}^{r,s}$ is as in Definition 4.2.1 (with the abuse of notation $F_i = \deg F_i$).*

[2] The notion of chain map is unrelated to the tensor product, and hence does not require a specific definition for graded-monoidal categories.

**Proposition 4.2.4.** *Definition 4.2.3 gives a well-defined chain map. Moreover, if $F_1$ and $F_2$ are identities then $F$ is the identity.*

*Proof.* Note that if both $F_1$ and $F_2$ are identities, then

$$\epsilon^{\boldsymbol{r},\boldsymbol{s}}_{F_1,F_2} = \Big(\big[0 + |\alpha_1|\,(\boldsymbol{r}_1) + |\alpha_1|\,(\boldsymbol{r}_1)\big] * |\alpha_2|\,(\boldsymbol{s}_2)\Big)^{-1}\big(|\alpha_1|\,(\boldsymbol{r}_1) * 0\big) = 1,$$

so that two identities induce the identity on the tensor product. To show the first part of the statement, we must check that $\beta \circ F = F \circ \alpha$. First, we unfold both sides of the equation:

$$\begin{aligned}
&\sum_{j\in\mathcal{J}} \beta^{\boldsymbol{s}-e_j}_j \circ F^{\boldsymbol{r},\boldsymbol{s}-e_j} \\
&\quad= \sum_{j_1\in\mathcal{J}_1} \Big[\beta^{\boldsymbol{s}_1-e_{j_1}}_{j_1} \otimes \mathrm{Id}\Big] \circ \Big[\epsilon^{\boldsymbol{r},\boldsymbol{s}-e_{j_1}}_{F_1,F_2}\; F_1^{\boldsymbol{r}_1,\boldsymbol{s}_1-e_{j_1}} \otimes F_2^{\boldsymbol{r}_2,\boldsymbol{s}_2}\Big] \\
&\qquad+ \sum_{j_2\in\mathcal{J}_2} \Big[\epsilon^{\boldsymbol{s}_1,\boldsymbol{s}_2-e_{j_2},j_2}_{B_1\otimes B_2}\; \mathrm{Id} \otimes \beta^{\boldsymbol{s}_2-e_{j_2}}_{j_2}\Big] \\
&\qquad\qquad\qquad\qquad \circ \Big[\epsilon^{\boldsymbol{r},\boldsymbol{s}-e_{j_2}}_{F_1,F_2}\; F_1^{\boldsymbol{r}_1,\boldsymbol{s}_1} \otimes F_2^{\boldsymbol{r}_2,\boldsymbol{s}_2-e_{j_2}}\Big] \\
&\quad= \left[\sum_{j_1\in\mathcal{J}_1} \epsilon^{\boldsymbol{r},\boldsymbol{s}-e_{j_1}}_{F_1,F_2}\; \beta^{\boldsymbol{s}_1-e_{j_1}}_{j_1} \circ F_1^{\boldsymbol{r}_1,\boldsymbol{s}_1-e_{j_1}}\right] \otimes F_2^{\boldsymbol{r}_2,\boldsymbol{s}_2} && ① \\
&\qquad+ F_1^{\boldsymbol{r}_1,\boldsymbol{s}_1} \otimes \Big[\sum_{j_2\in\mathcal{J}_2} \Big(\beta^{\boldsymbol{s}_2-e_{j_2}}_{j_2} * F_1^{\boldsymbol{r}_1,\boldsymbol{s}_1}\Big) \\
&\qquad\qquad\qquad \epsilon^{\boldsymbol{s}_1,\boldsymbol{s}_2-e_{j_2},j_2}_{B_1\otimes B_2}\,\epsilon^{\boldsymbol{r},\boldsymbol{s}-e_{j_2}}_{F_1,F_2}\; \beta^{\boldsymbol{s}_2-e_{j_2}}_{j_2} \circ F_2^{\boldsymbol{r}_2,\boldsymbol{s}_2-e_{j_2}}\Big] && ②
\end{aligned}$$

$$\begin{aligned}
&\sum_{i\in\mathcal{I}} F^{\boldsymbol{r}+e_i,\boldsymbol{s}} \circ \alpha^{\boldsymbol{r}}_i \\
&\quad= \sum_{i_1\in\mathcal{I}_1} \Big[\epsilon^{\boldsymbol{r}+e_{i_1},\boldsymbol{s}}_{F_1,F_2}\; F_1^{\boldsymbol{r}_1+e_{i_1},\boldsymbol{s}_1} \otimes F_2^{\boldsymbol{r}_2,\boldsymbol{s}_2}\Big] \circ \big[\alpha^{\boldsymbol{r}_1}_{i_1} \otimes \mathrm{Id}\big] \\
&\qquad+ \sum_{i_2\in\mathcal{I}_2} \Big[\epsilon^{\boldsymbol{r}+e_{i_2},\boldsymbol{s}}_{F_1,F_2}\; F_1^{\boldsymbol{r}_1,\boldsymbol{s}_1} \otimes F_2^{\boldsymbol{r}_2+e_{i_2},\boldsymbol{s}_2}\Big] \circ \Big[\epsilon^{\boldsymbol{r}_1,\boldsymbol{r}_2,i_2}_{A_1\otimes A_2}\; \mathrm{Id} \otimes \alpha^{\boldsymbol{r}_2}_{i_2}\Big] \\
&\quad= \left[\sum_{i_1\in\mathcal{I}_1} \big(F_2^{\boldsymbol{r}_2,\boldsymbol{s}_2} * \alpha^{\boldsymbol{r}_1}_{i_1}\big)\, \epsilon^{\boldsymbol{r}+e_{i_1},\boldsymbol{s}}_{F_1,F_2}\; F_1^{\boldsymbol{r}_1+e_{i_1},\boldsymbol{s}_1} \circ \alpha^{\boldsymbol{r}_1}_{i_1}\right] \otimes F_2^{\boldsymbol{r}_2,\boldsymbol{s}_2} && ① \\
&\qquad+ F_1^{\boldsymbol{r}_1,\boldsymbol{s}_1} \otimes \left[\sum_{i_2\in\mathcal{I}_2} \epsilon^{\boldsymbol{r}+e_{i_2},\boldsymbol{s}}_{F_1,F_2}\,\epsilon^{\boldsymbol{r}_1,\boldsymbol{r}_2,i_2}_{A_1\otimes A_2}\; F_2^{\boldsymbol{r}_2+e_{i_2},\boldsymbol{s}_2} \circ \alpha^{\boldsymbol{r}_2}_{i_2}\right] && ②
\end{aligned}$$

We want to show that the two terms labelled ① (resp. ②) are equal, using the chain map relation of $F_1$ (resp. $F_2$). To do so, we only need to check that the scalars are the same. In each case, we can assume that:

① $\beta_{j_1}^{s_{j_1}-1} + F_1^{\boldsymbol{r}_1,\boldsymbol{s}_1-\boldsymbol{e}_{j_1}} =_G F_1^{\boldsymbol{r}_1+\boldsymbol{e}_{i_1},\boldsymbol{s}_1} + \alpha_{i_1}^{r_{i_1}}$ for all $i_1 \in \mathcal{I}_1$ and $j_1 \in \mathcal{J}_1$,

② $\beta_{j_2}^{s_{j_2}-1} + F_2^{\boldsymbol{r}_2,\boldsymbol{s}_2-\boldsymbol{e}_{j_2}} =_G F_2^{\boldsymbol{r}_2+\boldsymbol{e}_{i_2},\boldsymbol{s}_2} + \alpha_{i_2}^{r_{i_2}}$ for all $i_2 \in \mathcal{I}_2$ and $j_2 \in \mathcal{J}_2$, and $|\boldsymbol{r}_1| = |\boldsymbol{s}_1|$.

These are exactly the assumptions needed to apply Lemma 4.2.2. □

### 4.2.3 Induced homotopies on the tensor product

Unpacking the standard[3] definition of a homotopy between chain maps, a homotopy $H$ between chain maps $F\colon \mathbb{A} \to \mathbb{B}$ and $G\colon \mathbb{A} \to \mathbb{B}$ is a set of morphisms $H^{\boldsymbol{r},\boldsymbol{s}}\colon A^{\boldsymbol{r}} \to B^{\boldsymbol{s}}$ for all $|\boldsymbol{r}| = |\boldsymbol{s}| + 1$ such that

$$F^{\boldsymbol{r},\boldsymbol{s}} - G^{\boldsymbol{r},\boldsymbol{s}} = \sum_{i\in\mathcal{I}} H^{\boldsymbol{r}+\boldsymbol{e}_i,\boldsymbol{s}} \circ \alpha_i^{\boldsymbol{r}} + \sum_{j\in\mathcal{J}} \beta_j^{\boldsymbol{s}-\boldsymbol{e}_j} \circ H^{\boldsymbol{r},\boldsymbol{s}-\boldsymbol{e}_j}.$$

**Definition 4.2.5.** *Let $F_k\colon A_k \to B_k$ (resp. $G_k$) be chain maps and $H_k\colon F_k \to G_k$ homotopies for each $k = 1,2$. Denote $F$ (resp. $G$) the chain map induced by $F_1$ and $F_2$ (resp. $G_1$ and $G_2$). Their* induced homotopy *is the homotopy $H\colon F \to G$ given by the data:*

$$H^{\boldsymbol{r},\boldsymbol{s}} = \epsilon_{H_1,F_2}^{\boldsymbol{r},\boldsymbol{s}} H_1^{\boldsymbol{r}_1,\boldsymbol{s}_1} \otimes F_2^{\boldsymbol{r}_2,\boldsymbol{s}_2} + (-1)^{|\boldsymbol{r}_1|} \epsilon_{G_1,H_2}^{\boldsymbol{r},\boldsymbol{s}} G_1^{\boldsymbol{r}_1,\boldsymbol{s}_1} \otimes H_2^{\boldsymbol{r}_2,\boldsymbol{s}_2}$$

*where the $\epsilon$'s are defined as in Definition 4.2.1.*

**Proposition 4.2.6.** *Definition 4.2.5 gives a well-defined homotopy.*

[3]As for the notion of chain map, the notion of homotopy is unrelated to the tensor product, and hence does not require a specific definition for graded-monoidal categories.

*Proof.* We must check that $F - G = \beta \circ H + H \circ \alpha$. We unfold the two terms on the right-hand side (RHS), with first the computation of $\beta \circ H$ restricted to the paths from $A^r$ to $B^s$:

$$
\begin{aligned}
(\beta \circ H)|_{A^r}^{B^s} &= \sum_{j \in \mathcal{J}} \beta_j^{s-e_j} \circ H^{r,s-e_j} \\
&= \sum_{j_1 \in \mathcal{J}_1} \left[\beta_{j_1}^{s_1-e_{j_1}} \otimes \mathrm{Id}\right] \\
&\qquad \circ \left[\epsilon_{H_1,F_2}^{r,s-e_{j_1}} H_1^{r_1,s_1-e_{j_1}} \otimes F_2^{r_2,s_2} \right.\\
&\qquad\qquad \left. +(-1)^{|r_1|}\epsilon_{G_1,H_2}^{r,s-e_{j_1}} G_1^{r_1,s_1-e_{j_1}} \otimes H_2^{r_2,s_2}\right] \\
&\quad + \sum_{j_2 \in \mathcal{J}_2} \left[\epsilon_{B_1 \otimes B_2}^{s_1,s_2-e_{j_2},j_2} \mathrm{Id} \otimes \beta_{j_2}^{s_2-e_{j_2}}\right] \\
&\qquad \circ \left[\epsilon_{H_1,F_2}^{r,s-e_{j_2}} H_1^{r_1,s_1} \otimes F_2^{r_2,s_2-e_{j_2}} \right.\\
&\qquad\qquad \left. +(-1)^{|r_1|}\epsilon_{G_1,H_2}^{r,s-e_{j_2}} G_1^{r_1,s_1} \otimes H_2^{r_2,s_2-e_{j_2}}\right] \\
&= \sum_{j_1 \in \mathcal{J}_1} \epsilon_{H_1,F_2}^{r,s-e_{j_1}} \left[\beta_{j_1}^{s_1-e_{j_1}} \circ H_1^{r_1,s_1-e_{j_1}}\right] \otimes F_2^{r_2,s_2} && ① \\
&\qquad + (-1)^{|r_1|}\epsilon_{G_1,H_2}^{r,s-e_{j_1}} \left[\beta_{j_1}^{s_1-e_{j_1}} \circ G_1^{r_1,s_1-e_{j_1}}\right] \otimes H_2^{r_2,s_2} && ② \\
&\quad + \sum_{j_2 \in \mathcal{J}_2} \left(\beta_{j_2}^{s_2-e_{j_2}} * H_1^{r_1,s_1}\right) \epsilon_{B_1 \otimes B_2}^{s_1,s_2-e_{j_2},j_2} \epsilon_{H_1,F_2}^{r,s-e_{j_2}} \\
&\qquad\qquad H_1^{r_1,s_1} \otimes \left[\beta_{j_2}^{s_2-e_{j_2}} \circ F_2^{r_2,s_2-e_{j_2}}\right] && ③ \\
&\qquad + \left(\beta_{j_2}^{s_2-e_{j_2}} * G_1^{r_1,s_1}\right) (-1)^{|r_1|} \epsilon_{B_1 \otimes B_2}^{s_1,s_2-e_{j_2},j_2} \epsilon_{G_1,H_2}^{r,s-e_{j_2}} \\
&\qquad\qquad G_1^{r_1,s_1} \otimes \left[\beta_{j_2}^{s_2-e_{j_2}} \circ H_2^{r_2,s_2-e_{j_2}}\right] && ④
\end{aligned}
$$

The computation of $H \circ \alpha$ restricted to the paths from $A^r$ to $B^s$ gives:

$$
\begin{aligned}
&(H \circ \alpha)|_{A^r}^{B^s} \\
&= \sum_{i \in \mathcal{I}} H^{r+e_i,s} \circ \alpha_i^r \\
&= \sum_{i_1 \in \mathcal{I}_1} \Big[\epsilon_{H_1,F_2}^{r+e_{i_1},s} H_1^{r_1+e_{i_1},s_1} \otimes F_2^{r_2,s_2} \\
&\qquad\qquad +(-1)^{|r_1|+1} \epsilon_{G_1,H_2}^{r+e_{i_1},s} G_1^{r_1+e_{i_1},s_1} \otimes H_2^{r_2,s_2}\Big] \\
&\qquad \circ \Big[\alpha_{i_1}^{r_1} \otimes \mathrm{Id}\Big] \\
&\quad + \sum_{i_2 \in \mathcal{I}_2} \Big[\epsilon_{H_1,F_2}^{r+e_{i_2},s} H_1^{r_1,s_1} \otimes F_2^{r_2+e_{i_2},s_2} \\
&\qquad\qquad +(-1)^{|r_1|} \epsilon_{G_1,H_2}^{r+e_{i_2},s} G_1^{r_1,s_1} \otimes H_2^{r_2+e_{i_2},s_2}\Big] \\
&\qquad \circ \Big[\epsilon_{A_1 \otimes A_2}^{r_1,r_2,i_2} \mathrm{Id} \otimes \alpha_{i_2}^{r_2}\Big] \\
&= \sum_{i_1 \in \mathcal{I}_1} \left(F_2^{r_2,s_2} * \alpha_{i_1}^{r_1}\right) \epsilon_{H_1,F_2}^{r+e_{i_1},s} \Big[H_1^{r_1+e_{i_1},s_1} \circ \alpha_{i_1}^{r_1}\Big] \otimes F_2^{r_2,s_2} && \text{①} \\
&\qquad + \left(H_2^{r_2,s_2} * \alpha_{i_1}^{r_1}\right) (-1)^{|r_1|+1} \epsilon_{G_1,H_2}^{r+e_{i_1},s} \Big[G_1^{r_1+e_{i_1},s_1} \circ \alpha_{i_1}^{r_1}\Big] \otimes H_2^{r_2,s_2} && \text{②} \\
&\quad + \sum_{i_2 \in \mathcal{I}_2} \epsilon_{A_1 \otimes A_2}^{r_1,r_2,i_2} \epsilon_{H_1,F_2}^{r+e_{i_2},s} \; H_1^{r_1,s_1} \otimes \Big[F_2^{r_2+e_{i_2},s_2} \circ \alpha_{i_2}^{r_2}\Big] && \text{③} \\
&\qquad + (-1)^{|r_1|} \epsilon_{A_1 \otimes A_2}^{r_1,r_2,i_2} \epsilon_{G_1,H_2}^{r+e_{i_2},s} \; G_1^{r_1,s_1} \otimes \Big[H_2^{r_2+e_{i_2},s_2} \circ \alpha_{i_2}^{r_2}\Big] && \text{④}
\end{aligned}
$$

We find four different pairs of terms, labelled ① to ④, that we simplify pairwise using chain map or chain homotopy relations:

① We can assume that for all $i_1 \in \mathcal{I}_1$ and all $j_1 \in \mathcal{J}_1$:

$$\beta_{j_1}^{s_1-e_{j_1}} + H_1^{r_1,s_1-e_{j_1}} =_G H_1^{r_1+e_{i_1},s_1} + \alpha_{i_1}^{r_1} =_G F_1^{r_1,s_1} =_G G_1^{r_1,s_1}.$$

Thanks to Lemma 4.2.2 (i), the sum of the two ①-terms is:

$$\left[\epsilon_{F_1,F_2}^{r_1,s_1} F_1^{r_1,s_1} - \epsilon_{G_1,F_2}^{r_1,s_1} G_1^{r_1,s_1}\right] \otimes F_2^{r_2,s_2}.$$

The computation is similar for the case ④.

③ We can assume $|r_1| + 1 = |s_1|$ and that:

$$\beta_{j_2}^{s_2-e_{j_2}} + F_2^{r_2,s_2-e_{j_2}} =_G F_2^{r_2+e_{i_2},s_2} + \alpha_{i_2}^{r_2}.$$

Then, Lemma 4.2.2 (ii) shows that:

$$\left(H_1^{r_1,s_1} * \beta_{j_2}^{s_2-e_{j_2}}\right) \epsilon_{B_1 \otimes B_2}^{s_1,s_2-e_{j_2},j_2} \epsilon_{H_1,F_2}^{r,s-e_{j_2}} = -\epsilon_{A_1 \otimes A_2}^{r_1,r_2,i_2} \epsilon_{H_1,F_2}^{r+e_{i_2},s}$$

independently of $j_2$ and $i_2$. We conclude that the sum of the pair ③ is zero. The computation is similar for the case ②.

We conclude:

$$\begin{aligned}
\text{RHS} &= \left[\epsilon_{F_1,F_2}^{r_1,s_1} F_1^{r_1,s_1} - \epsilon_{G_1,F_2}^{r_1,s_1} G_1^{r_1,s_1}\right] \otimes F_2^{r_2,s_2} \\
&\quad + G_1^{r_1,s_1} \otimes \left[\epsilon_{G_1,F_2}^{r_1,s_1} F_2^{r_2,s_2} - \epsilon_{G_1,G_2}^{r_1,s_1} G_2^{r_2,s_2}\right] \\
&= \epsilon_{F_1,F_2}^{r_1,s_1} F_1^{r_1,s_1} \otimes F_2^{r_2,s_2} - \epsilon_{G_1,G_2}^{r_1,s_1} G_1^{r_1,s_1} \otimes G_2^{r_2,s_2} \\
&= (F-G)|_{A^r}^{B^s}.
\end{aligned}$$

□

# PART II
# Rewriting theory and higher algebra

# 5

# Linear Gray polygraphs

This chapter introduces the necessary categorical structures to present graded-2-categories and define their rewriting theory. This can be understood as a linear analogue to the work of Forest and Mimram [81, section 2 and 3]; equivalently, as a generalization of the work of Alleaume on linear $n$-polygraphs [3] to allow weak interchangers. For an introduction to the ideas of this section, see subsections ii.1.2, ii.1.6 and ii.2.1.

The notion of an $n$-sesquicategory was first defined by Street [177] in the case $n = 2$. The general case was independently introduced by Forest–Mimram [81] (following the general theory of Weber [186]) under the name of "$n$-precategories" and by Araújo [6, section 1.6] under the name of "$n$-sesquicategories". Although we shall follow Forest and Mimram's presentation, we choose Araújo's terminology to avoid confusion with already existing notions of $n$-precategories in the literature. Enriched category theory provides yet another defining approach to $n$-sesquicategories; see [81, section 2.4 and Appendix A].

To motivate the formal definitions, we start with an example in section 5.1. We expect it to be sufficient for the impatient reader. Section 5.2 review the notions of $n$-sesquicategories and $n$-sesquipolygraphs, following [81] (which they respectively call "$n$-prepolygraphs" and "$n$-precategories"). We then introduce their linear analogue in section 5.3. Finally, section 5.4 defines Gray polygraphs (still following [81]) and linear Gray polygraphs, the latter defining presentations for graded-2-categories.

In this chapter and as throughout this thesis, every categorical structure is assumed to be small.

*Notation* 5.0.1. As much as possible, we use calligraphic fonts (e.g. $\mathcal{C}$) for categories, and sans serif fonts (e.g. $\mathsf{P}$) for their presentations. If necessary, blackboard fonts (e.g. $\mathbb{P}$) refer to purely set-theoretic concepts, and typewriter fonts (e.g. $\mathtt{P}$) to purely linear concepts.

## 5.1 A motivating example

In this section, we read $\mathfrak{gl}_2$-webs from bottom to top. Recall the $\mathbb{Z}[q, q^{-1}]$-linear categories $\mathbf{Web}_d$ of $\mathfrak{gl}_2$-webs defined in section 1.2. Taking the union

$$\underline{\Lambda}_{\otimes} := \bigsqcup_{d \in \mathbb{N}} \underline{\Lambda}_d \quad \text{and} \quad \mathbf{Web}_{\otimes} := \bigsqcup_{d \in \mathbb{N}} \mathbf{Web}_d$$

defines a $\mathbb{Z}[q, q^{-1}]$-linear monoidal category $\mathbf{Web}_{\otimes}$ with objects $\underline{\Lambda}_{\otimes}$, the monoidal product being the vertical juxtaposition of $\mathfrak{gl}_2$-webs. (The reader should not confuse $\mathbf{Web}_{\otimes}$ with the category $\mathbf{Web}$ defined in Definition 2.0.1.)

View $\mathbf{Web}_{\otimes}$ as a $\mathbb{Z}[q, q^{-1}]$-linear 2-category with one object. To present higher categories, one needs to provide generators for each categorical level; in this framework, relations are simply the generators at the highest categorical level. Thus, a presentation $\mathsf{Web}$ of $\mathbf{Web}_{\otimes}$ should have four levels: object generators, 1-morphism generators, 2-morphism generators, and finally 3-morphism generators, that is, generating relations. The first three levels are given by the following sets:

$$\mathsf{Web}_0 = \{*\}, \quad \mathsf{Web}_1 = \left\{ \,|\, , \,\|\, \right\}, \quad \mathsf{Web}_2 = \left\{ W_- := \right.$$ 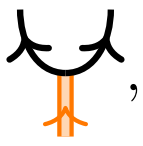 $$, W_+ :=$$ 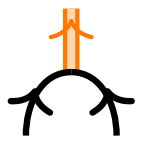 $$\left. \right\}.$$

Recall that we read $\mathfrak{gl}_2$-webs from bottom to top. Each level $\mathsf{Web}_{n+1}$ comes equipped with a source and target map into the *free $n$-sesquicategory* $\mathsf{Web}_n^*$ (subsection 5.2.4) generated by the previous level $\mathsf{Web}_n$:

$$\mathsf{Web}_n^* \underset{t_n}{\overset{s_n}{\leftleftarrows}} \mathsf{Web}_{n+1}.$$

For objects, the free 0-sesquicategory is the set itself: $\mathsf{Web}_0^* = \mathsf{Web}_0$. In our case, we have $\mathsf{Web}_0 = \{*\}$ and the maps $s_0, t_0$ are the trivial maps. The free 1-sesquicategory is given by formal horizontal juxtaposition of 1-generators, so that $\mathsf{Web}_1^* = \underline{\Lambda}_\otimes$. (If there were more than one object, we should ask the horizontal juxtaposition to be compatible with the 0-source $s_0$ and the 0-target $t_0$). The maps $s_1, t_1$ are the obvious ones, reading $s_1$ on the bottom and $t_1$ on top of a given web, respectively. The data $\mathsf{Web}_{\leq 2} = (\mathsf{Web}_0, \mathsf{Web}_1, \mathsf{Web}_2)$, together with their source and target maps, form the data of a *2-sesquipolygraph* (subsection 5.2.5):

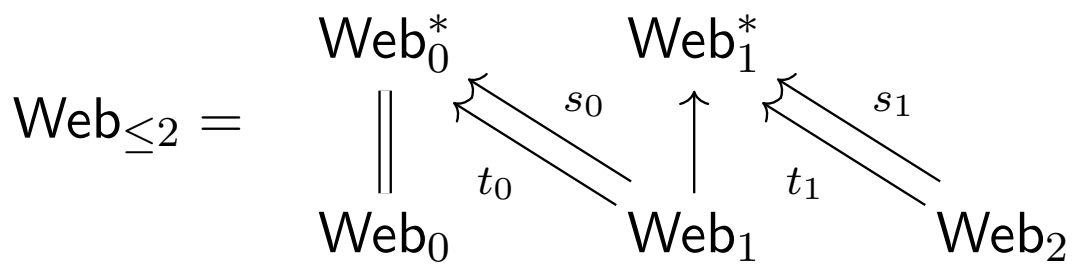

To define the highest level consisting of the defining relations, we must first describe the free 2-sesquicategory $\mathsf{Web}_2^*$ generated by $\mathsf{Web}_2$. First, juxtapose horizontally the 2-generators $W_-$ and $W_+$ with elements of $\mathsf{Web}_1^*$, both on the left and on the right. This process is called *whiskering*. This produces the webs $W_{-,i}$ and $W_{+,i}$ which we used as the generating morphisms of $\mathbf{Web}_d$ in section 1.2. Then, juxtapose vertically an arbitrary number of these whiskered generators, gluing along 1-source $s_1$ and 1-target $t_1$. This recovers generic $\mathfrak{gl}_2$-webs as defined in section 1.2, and defines the 2-morphisms of $\mathsf{Web}_2^*$. Note that we did not mod out by the interchange relation for webs: indeed, $\mathsf{Web}_2^*$ is a 2-sesquicategory, and not a strict 2-category.

Finally, the set $\mathsf{Web}_3$ consists of the following 3-generators (ignore the scalar $q + q^{-1}$ for now):

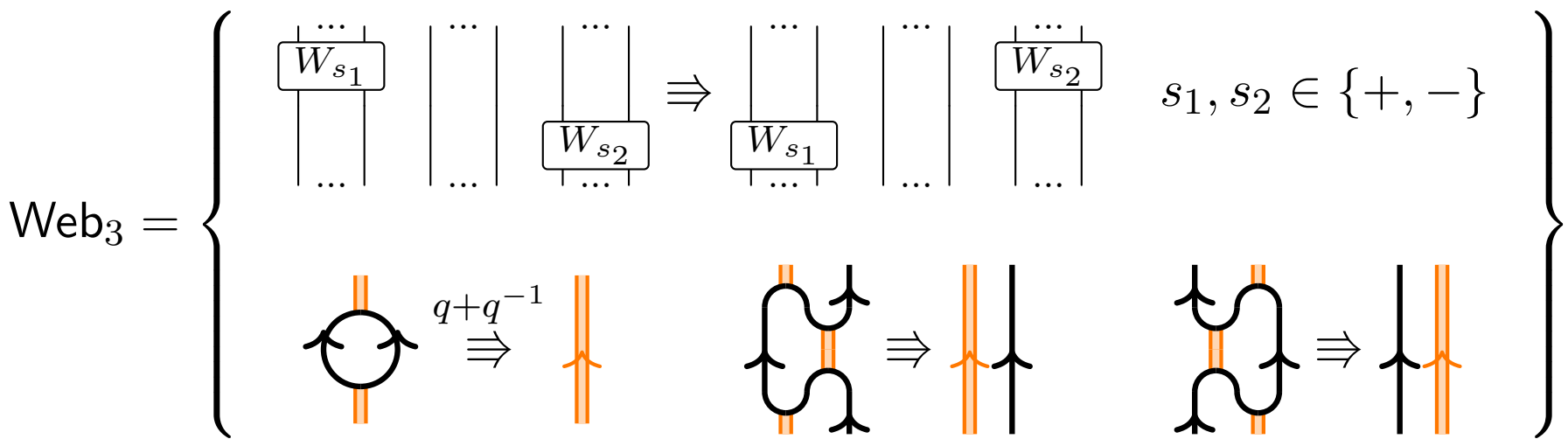

Here 3-cells are pictured as two-frame *movies* from their source $s_2$ (on the left) to their target $t_2$ (on the right). The first family of 3-generators are the interchange generators on $\mathsf{Web}_{\leq 2}$. They from the *3-sesquipolygraph of interchangers* (subsection 5.4.1), denoted $\mathsf{Web}_{\leq 2}\mathsf{Gray}$.

The set $\mathsf{Web}_3$ describes generating relations; to describe *all* relations, we must define the free 3-sesquicategory $\mathsf{Web}_3^*$ generated by $\mathsf{Web}_3$. First, we revisit our terminology: we call horizontal juxtaposition the *0-composition* (as we glue along objects) and vertical juxtaposition the *1-composition* (as we glue along 1-morphisms). We define the *contextualization* of a 3-generator as first 0-whiskering with elements of $\mathsf{Web}_1^*$, and then 1-whiskering with elements of $\mathsf{Web}_2^*$. This can be pictured as the relevant 0-compositions and 1-compositions on its source and target:

$\Gamma[A] :=$ 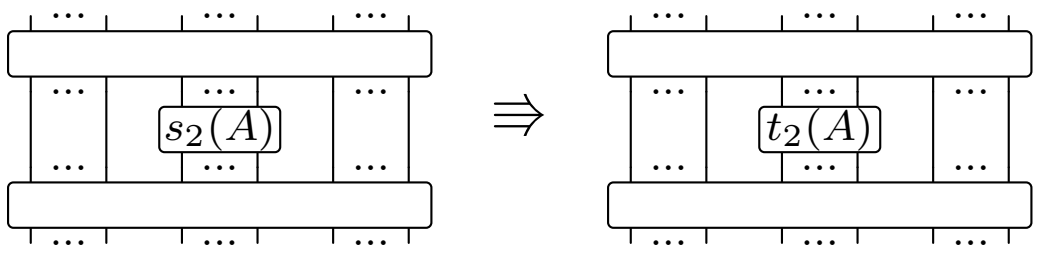 (5.1)

Here $A \in \mathsf{Web}_3$ is a 3-generator and $\Gamma$ is a *context*; that is, the data of 1-cells and 2-cells "surrounding $A$" via 0- and 1-whiskerings. We write $\mathrm{Cont}(\mathsf{Web}_3)$ the set of contextualized 3-generators. Finally, a generic 3-morphism in $\mathsf{Web}_3^*$ is a 2-composition of contextualized 3-generators, gluing along 2-source and 2-target. This is pictured as a multi-frame movie, as in the following example:

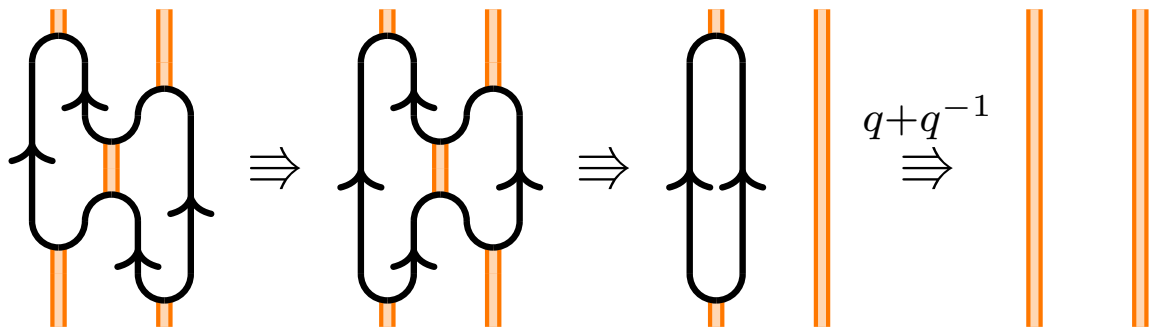

3-morphisms in $\mathsf{Web}_3^*$ have the following structure: a 2-composition (pictured as a composition of movies) and actions of respectively $\mathsf{Web}_1^*$ and $\mathsf{Web}_2^*$ via the 0- and 1-whiskerings (pictured as horizontal and vertical juxtapositions).

In general, the $(k+1)$-cells of an *$n$-sesquicategory* (subsection 5.2.2) have a $k$-composition gluing along $k$-cells, and actions of lower cells via whiskerings. *Contextualization* (subsection 5.2.3) constitutes the combined action of all whiskerings. Given a *globular extension* $\mathsf{P}$ (subsection 5.2.1), the $n$-cells of the *free $n$-sesquicategory* $\mathsf{P}^*$ (subsection 5.2.4) are given by formal $(n-1)$-compositions of contextualized $n$-generators $\mathrm{Cont}(\mathsf{P})$.

So far, we have ignored the linear structure inherent to $\mathbf{Web}_\otimes$. To encode it, each 3-generator $W_1 \Rrightarrow W_2$ in $\mathsf{Web}_3$ is equipped with the extra data of a scalar $\lambda \in \mathbb{Z}[q,q^{-1}]$, capturing the relation $W_1 = \lambda W_2$ in $\mathbf{Web}_\otimes$. This scalar is $q+q^{-1}$ for the circle evaluation, and 1 for the rest of the 3-generators. We encapsulate this data in a function $\mathrm{scl}\colon \mathsf{Web}_3 \to \mathbb{Z}[q,q^{-1}]$. The data

$\mathsf{Web} :=$ 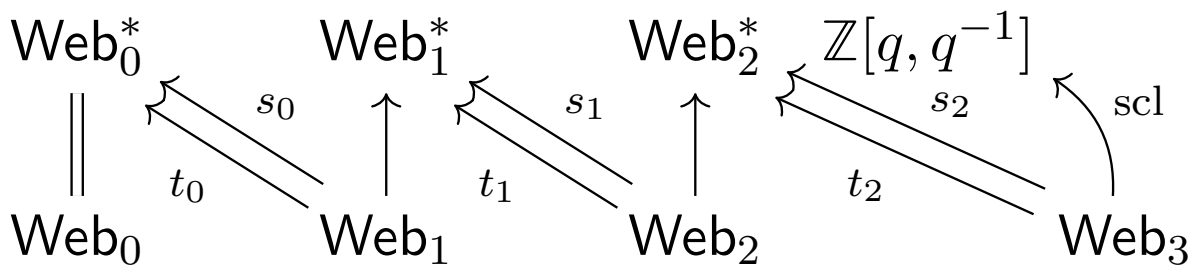

form the data of a *scalar Gray polygraph* (subsections 5.2.7 and 5.4.2). Here *scalar* stands for the fact that 3-morphisms have associated scalars, and *Gray* stands for the fact that Web contains its own 3-sesquipolygraph of interchangers, that is, $\mathsf{Web}_{\leq 2}\mathsf{Gray} \subset \mathsf{Web}_3$. It captures all the defining data of $\mathbf{Web}_{\otimes}$, and as such provides a *presentation* of $\mathbf{Web}_{\otimes}$ (Definition 5.3.4).

Equivalently, we can view Web as a *linear Gray polygraph*:

$$\mathsf{Web} := \begin{array}{ccccccc} \mathsf{Web}_0^* & & \mathsf{Web}_1^* & & \mathsf{Web}_2^l & & \\ \| & \overset{s_0}{\underset{t_0}{\nwarrow\!\nwarrow}} & \uparrow & \overset{s_1}{\underset{t_1}{\nwarrow\!\nwarrow}} & \uparrow & \overset{s_2^l}{\underset{t_2^l}{\nwarrow\!\nwarrow}} & \\ \mathsf{Web}_0 & & \mathsf{Web}_1 & & \mathsf{Web}_2 & & \mathsf{Web}_3 \end{array}$$

Here $\mathsf{Web}_2^l$ is the (set of 2-cells of the) *free linear 2-sesquicategory generated by* $\mathsf{Web}_{\leq 2}$ (subsection 5.3.2), defined as follows. Denote $\mathsf{Web}_2^*(f,g)$ the set of elements $\alpha \in \mathsf{Web}_2^*$ having $f$ (resp. $g$) as a 1-source $s_1(\alpha) = f$ (resp. as a 1-target $t_1(\alpha) = g$). Denote $\mathsf{Web}_2^l(f,g)$ the free $\mathbb{Z}[q,q^{-1}]$-module generated by $\mathsf{Web}_2^*(f,g)$. The set $\mathsf{Web}_2^l$ is the union of the $\mathsf{Web}_2^l(f,g)$s, with the obvious 1-source and 1-target maps. Finally, the set $\mathsf{Web}_3$ remains unchanged, but we set new 2-source map $s_2^l := s_2$ and 2-target map $t_2^l = \mathrm{scl} \cdot t_2$.

In this form, Web is a *linear 3-sesquipolygraph* (subsection 5.3.3). Because Web contains its own 3-sesquipolygraph of interchangers, we call it a *linear Gray polygraph* (Definition 5.4.3). The fact that it arises from a scalar 3-sesquipolygraph makes it a *monomial* linear 3-sesquipolygraph (subsection 5.3.4). In a general linear $n$-sesquipolygraph, source and target of generating 3-cells can be any linear combinations of 2-cells.

In this example, interchangers have trivial associated scalars; indeed, $\mathbf{Web}_{\otimes}$ is a genuine monoidal category. However, one can consider interchangers with arbitrary scalars using the same formalism. Given the data of $(G, \mu, \Bbbk)$ as in Notation 5.3.1 and of a *$G$-graded 2-sesquipolygraph* Q (subsection 5.2.6), one can define generating $(G,\mu)$-graded interchangers on Q. This defines the *3-sesquipolygraph of $(G,\mu)$-graded interchangers* (subsection 5.4.2), denoted QGray. In general, a *linear Gray polygraph* (Definition 5.4.3) is a linear 3-sesquipolygraph that contains its own 3-sesquipolygraph of $(G,\mu)$-graded interchangers. Linear Gray polygraphs provide a notion of presentation for graded-2-categories (Definition 5.4.4), suitable for the purpose of rewriting theory.

## 5.2 $n$-sesquicategories and their presentations

We review the notion of $n$-sesquicategories and their presentations, following the presentation given in [81]. The last subsections introduce graded and scalar variants of $n$-sesquicategories and their presentations.

### 5.2.1 $n$-globular sets

Let $n \in \mathbb{N} = \{0, 1, \ldots\}$. An *$n$-globular set* $\mathcal{C}$ is a diagram of sets and functions as follows:

$$\mathcal{C}_0 \underset{t_0}{\overset{s_0}{\leftleftarrows}} \mathcal{C}_1 \underset{t_1}{\overset{s_1}{\leftleftarrows}} \ldots \underset{t_{n-1}}{\overset{s_{n-1}}{\leftleftarrows}} \mathcal{C}_n$$

such that $s_j \circ s_{j+1} = s_j \circ t_{j+1}$ and $t_j \circ s_{j+1} = t_j \circ t_{j+1}$ for each $0 \leq j < n$. The maps $s_j$ and $t_j$ are respectively called *source maps* and *target maps*. An element $u \in \mathcal{C}_j$ is called a *$j$-cell*, with $s_{j-1}(u)$ and $t_{j-1}(u)$ respectively its *source* and *target*, which we sometimes simply denote by $s(u)$ and $t(u)$. We refer to $j$ as the *dimension* of $u$. For $0 \leq i < j$, we define the *$i$-source* of $u$ to be

$$s_i(u) = s_i \circ (\text{any suitable composition of source and target maps})\,(u),$$

where the choice in the bracket does not matter thanks to the properties of source and target maps. Note that the subscript indicates that $s_i(u)$ is an $i$-cell. We define the *$i$-target* $t_i(u)$ similarly. A *morphism of globular sets* $f\colon \mathcal{C} \to \mathcal{D}$ is a family $f_i\colon \mathcal{C}_i \to \mathcal{D}_i$ of functions that commute with the source and target maps. It is an isomorphism if each function $f_i$ is a bijection.

Given an $n$-globular set $\mathcal{C}$, a $0$-*sphere* is an ordered pair of $0$-cells in $\mathcal{C}$. For $0 < i \leq n$, an *$i$-sphere* is an ordered pair $(f, g)$ of $i$-cells such that $s(f) = s(g)$ and $t(f) = t(g)$. For $0 \leq i < k \leq n$ and an $i$-sphere $(f, g)$, we set

$$\mathcal{C}_k(f, g) = \{u \in \mathcal{C}_k \mid s_i(u) = f, t_i(u) = g\}.$$

The source and target maps restrict to maps $\mathcal{C}_k(f, g) \underset{t_k}{\overset{s_k}{\leftleftarrows}} \mathcal{C}_{k+1}(f, g)$.

A *globular extension of* $\mathcal{C}$ is an $(n+1)$-globular set $\mathsf{P}$ such that $\mathsf{P}_k = \mathcal{C}_k$ for $0 \leq k \leq n$:

$$\underbrace{\mathcal{C}_0 \overset{s_0}{\underset{t_0}{\leftleftarrows}} \mathcal{C}_1 \overset{s_1}{\underset{t_1}{\leftleftarrows}} \dots \overset{s_{n-1}}{\underset{t_{n-1}}{\leftleftarrows}} \mathcal{C}_n}_{\mathcal{C}} \overset{s_n}{\underset{t_n}{\leftleftarrows}} \mathsf{P}_{n+1}$$

For $k \in \mathbb{N}$ such that $0 \leq k \leq n$, the $k$*-restriction of* $\mathcal{C}$ is the subglobular set

$$\mathcal{C}_{\leq k} := \mathcal{C}_0 \overset{s_0}{\underset{t_0}{\leftleftarrows}} \mathcal{C}_1 \overset{s_1}{\underset{t_1}{\leftleftarrows}} \dots \overset{s_{k-1}}{\underset{t_{k-1}}{\leftleftarrows}} \mathcal{C}_k$$

### 5.2.2 $n$-sesquicategories

We review the definition of an $n$-sesquicategory (and $n$-sesquifunctor) introduced in [81, section 2.2] (with minor changes to the presentation).

An $n$*-sesquicategory* is the data of an $n$-globular set $\mathcal{C}$ together with

- identity functions $\mathrm{id}^k \colon \mathcal{C}_{k-1} \to \mathcal{C}_k$, for $0 < k \leq n$,
- composition functions $\star_{k,k} \colon \mathcal{C}_k \times_{k-1} \mathcal{C}_k \to \mathcal{C}_k$ for $0 < k \leq n$,
- left- and right-whisker functions $\star_{i,k} \colon \mathcal{C}_i \times_{i-1} \mathcal{C}_k \to \mathcal{C}_k$ and $\star_{k,i} \colon \mathcal{C}_k \times_{i-1} \mathcal{C}_i \to \mathcal{C}_k$ for $0 < i < k \leq n$, satisfying the axioms (i) and (ii) below. Hereabove we abbreviated $\times_{\mathcal{C}_k}$ with $\times_k$. We use similar notations for composition and whiskers, but one cannot confuse one for the other as they have different domain. In fact, this choice of notation emphasizes that whisker functions should be thought of compositions with identities of cells of lower dimension. As such, composition and whiskers have similar properties, and it is sometimes useful to consider $\star_{k,l}$ for *all* $0 \leq k, l \leq n$.

Note that for $0 \leq k, l \leq n$, the function $\star_{k,l}$ is defined on $(u, v) \in \mathcal{C}_k \times \mathcal{C}_l$ if and only if $s_i(u) = t_i(v)$, where $i = \min(k, l) - 1$. In that case, we say that $u$ and $v$ are $i$*-composable*, and we write $u \star_i v$, or even $u \star v$, for $u \star_{k,l} v$. While the notation $u \star_{k,l} v$ emphasizes the dimension of the respective cells, the notation $u \star_i v$ emphasizes the dimension of the compatibility condition. Also, for $u$ an $(i-1)$-cell we write $\mathrm{id}_u$ instead of $\mathrm{id}^i(u)$. The axioms of an $n$-sesquicategory are as follows:

(i) for $0 < k \leq n$, with $f \in \mathcal{C}_{k-1}$ and $\alpha, \beta, \gamma \in \mathcal{C}_k$ suitably $k$-composable:

$$t_{k-1}(\mathrm{id}_f) = f = s_{k-1}(\mathrm{id}_f)$$
$$s_{k-1}(\alpha \star_{k,k} \beta) = s_{k-1}(\beta) \qquad t_{k-1}(\alpha \star_{k,k} \beta) = t_{k-1}(\alpha)$$
$$\mathrm{id}_{t(\alpha)} \star_{k,k} \alpha = \alpha = \alpha \star_{k,k} \mathrm{id}_{s(\alpha)}$$
$$\alpha \star_{k,k} (\beta \star_{k,k} \gamma) = (\alpha \star_{k,k} \beta) \star_{k,k} \gamma$$

(ii) for $0 < i < k, k' \leq n$, with $x \in \mathcal{C}_{i-1}$, $f, g \in \mathcal{C}_i$, $\phi \in \mathcal{C}_{k-1}$ and $A \in \mathcal{C}_k$, $B \in \mathcal{C}_{k'}$ suitably composable:

$$f \star_{i,k} (g \star_{i,k} A) = (f \star_{i,i} g) \star_{i,k} A \qquad (A \star_{k,i} f) \star_{k,i} g = A \star_{k,i} (f \star_{i,i} g)$$
$$\mathrm{id}_x \star_{i,i} A = A \qquad A \star_{i,i} \mathrm{id}_x = A$$
$$(f \star_{i,k} A) \star_{k,i} g = f \star_{i,k} (A \star_{k,i} g)$$

$$s_{k-1}(f \star_{i,k} A) = f \star_{i,k-1} s_{k-1}(A) \qquad s_{k-1}(A \star_{k,i} f) = s_{k-1}(A) \star_{k-1,i} f$$
$$t_{k-1}(f \star_{i,k} A) = f \star_{i,k-1} t_{k-1}(A) \qquad t_{k-1}(A \star_{k,i} f) = t_{k-1}(A) \star_{k-1,i} f$$
$$f \star_{i,k} \mathrm{id}_\phi = \mathrm{id}_{f \star_{i,k-1} \phi} \qquad \mathrm{id}_\phi \star_{k,i} f = \mathrm{id}_{\phi \star_{k-1,i} f}$$
$$f \star_{i,\max(k,k')} (A \star_{k,k'} B) = (f \star_{i,k} A) \star_{k,k'} (f \star_{i,k'} B)$$
$$(A \star_{k,k'} B) \star_{\max(k,k'),i} f = (A \star_{k,i} f) \star_{k,k'} (B \star_{k',i} f)$$

An $n$-*sesquifunctor* between two $n$-sesquicategories is a morphism between the underlying globular sets, preserving identities and compositions as expected. It is an isomorphism if the underlying morphism of globular sets is an isomorphism. This ends the definition of an $n$-sesquicategory and of an $n$-sesquifunctor. ◇

*Remark* 5.2.1 (low-dimensional cases). A 0-sesquicategory is simply a set, and a 1-sesquicategory a category. The distinction with strict $n$-categories only appears when $n \geq 2$. Contrary to strict $n$-categories, 2-cells of an $n$-sesquicategories cannot be horizontally composed. Instead, they can be whiskered with 1-cells, understood as acting as identity 2-cells. In particular, in a 3-sesquicategory there is a priori no relationship between the two sides of the (2-dimensional) interchange law; see Fig. ii.2. We use string diagrammatics to picture 2-cells of an $n$-sesquicategory. String diagrams are equipped with a Morse function on the generators, as for graded-2-categories. 3-cells are then pictured as *movies*, i.e. paths of 2-cells, and 4-cells as *movie moves*, i.e. paths of paths of 2-cells.

*Remark* 5.2.2. With the single-index notation, the relations above become:

(i) for $0 < k \leq n$, with $f \in \mathcal{C}_{k-1}$ and $\alpha, \beta, \gamma \in \mathcal{C}_k$ suitably $k$-composable:

$$t_{k-1}(\mathrm{id}_f) = f = s_{k-1}(\mathrm{id}_f)$$
$$s_{k-1}(\alpha \star_{k-1} \beta) = s_{k-1}(\beta) \qquad t_{k-1}(\alpha \star_{k-1} \beta) = t_{k-1}(\alpha)$$
$$\mathrm{id}_{t(\alpha)} \star_{k-1} \alpha = \alpha = \alpha \star_{k-1} \mathrm{id}_{s(\alpha)}$$
$$\alpha \star_{k-1} (\beta \star_{k-1} \gamma) = (\alpha \star_{k-1} \beta) \star_{k-1} \gamma$$

(ii) for $0 < i < k, k' \leq n$, with $x \in \mathcal{C}_{i-1}$, $f, g \in \mathcal{C}_i$, $\phi \in \mathcal{C}_{k-1}$ and $A \in \mathcal{C}_k$, $B \in \mathcal{C}_{k'}$ suitably composable:

$$f \star_{i-1} (g \star_{i-1} A) = (f \star_{i-1} g) \star_{i-1} A \quad (A \star_{i-1} f) \star_{i-1} g = A \star_{i-1} (f \star_{i-1} g)$$
$$\mathrm{id}_x \star_{i-1} A = A \qquad\qquad A \star_{i-1} \mathrm{id}_x = A$$
$$(f \star_{i-1} A) \star_{i-1} g = f \star_{i-1} (A \star_{i-1} g)$$

$$s_{k-1}(f \star_{i-1} A) = f \star_{i-1} s_{k-1}(A) \qquad s_{k-1}(A \star_{i-1} f) = s_{k-1}(A) \star_{i-1} f$$
$$t_{k-1}(f \star_{i-1} A) = f \star_{i-1} t_{k-1}(A) \qquad t_{k-1}(A \star_{i-1} f) = t_{k-1}(A) \star_{i-1} f$$
$$f \star_{i-1} \mathrm{id}_\phi = \mathrm{id}_{f \star_{i-1} \phi} \qquad \mathrm{id}_\phi \star_{i-1} f = \mathrm{id}_{\phi \star_{i-1} f}$$
$$f \star_{i-1} (A \star_{\min(k,k')-1} B) = (f \star_{i-1} A) \star_{\min(k,k')-1} (f \star_{i-1} B)$$
$$(A \star_{\min(k,k')-1} B) \star_{i-1} f = (A \star_{i-1} f) \star_{\min(k,k')-1} (B \star_{i-1} f)$$

*Remark* 5.2.3. Note that if we let $i \in \mathbb{N}$ be such that $0 \leq i < n$ and $(f, g)$ be an $(i-1)$-sphere[1] in $\mathcal{C}$, the identity and composition functions restrict as follows:

$$\mathrm{id}^{i+1} \colon \mathcal{C}_i(f, g) \to \mathcal{C}_{i+1}(f, g)$$
$$\star_{i+1,i+1} \colon \mathcal{C}_{i+1}(f, g) \times_i \mathcal{C}_{i+1}(f, g) \to \mathcal{C}_{i+1}(f, g)$$

[1] For $i = 0$, we abuse notation and assume there exists a single $(-1)$-sphere $\square$, and denote $\mathcal{C}_i(\square) \coloneqq \mathcal{C}_i$.

It follows from the axioms of an $n$-sesquicategory that, equipped with the maps $\mathrm{id}^{i+1}$ and $\star_{i+1,i+1}$, the 1-globular set

$$\mathcal{C}_i(f,g) \underset{t_i}{\overset{s_i}{\longleftleftarrows}} \mathcal{C}_{i+1}(f,g)$$

defines a 1-category.

### 5.2.3 Contexts

Let $\mathcal{C}$ be an $n$-sesquicategory and let $\square = (f,g)$ an $(i-1)$-sphere for $0 < i \leq n$. A *context in* $\mathcal{C}$ *with boundary* $\square$ is a formal composition

$$\Gamma \coloneqq v_i \star_{i-1} (\ldots \star_1 (v_1 \star_0 \square \star_0 w_1) \star_1 \ldots) \star_{i-1} w_i \qquad (5.2)$$

where $v_j, w_j$ are $j$-cells suitably composable. We set $s_{i-1}(\Gamma) \coloneqq s_{i-1}(w_i)$ and $t_{i-1}(\Gamma) \coloneqq t_{i-1}(v_i)$. For instance, if $\square = (f,g)$ is a 1-sphere we picture $\Gamma$ as

$\Gamma \coloneqq$ 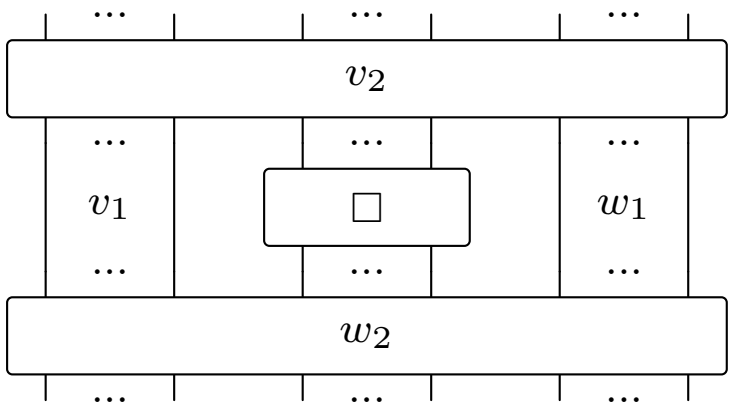

(5.3)

We write $\mathrm{Cont}(\square)$ the set of contexts in $\mathcal{C}$ with boundary $\square$.

For $k \in \mathbb{N}$ such that $i \leq k \leq n$ and for each $k$-cell $A$ in $\mathcal{C}$ such that $(s_{i-1}(A), t_{i-1}(A)) = \square$, we write

$$\Gamma[A] \coloneqq v_i \star_{i-1} (\ldots \star_1 (v_1 \star_0 A \star_0 w_1) \star_1 \ldots) \star_{i-1} w_i. \qquad (5.4)$$

We call $\Gamma[A]$ a *contextualization* of $A$. Any context $\Gamma$ defines a function

$$\Gamma \colon \mathcal{C}_k(\square) \to \mathcal{C}_k(s_{i-1}(\Gamma), t_{i-1}(\Gamma)).$$

We call this function *contextualization* with $\Gamma$.

Recall from Remark 5.2.3 how $(i-1)$-spheres in $\mathcal{C}$ defines categories of $i$-cells and $(i+1)$-cells. It follows from the axioms of an $n$-sesquicategory that contextualization suitably commutes with source, target, identity and

composition:

$$\begin{array}{ccc}
\mathcal{C}_i(\square) & \overset{s_i}{\underset{t_i}{\longleftarrow\!\!\longleftarrow}} & \mathcal{C}_{i+1}(\square) \\
\Gamma\big\downarrow & & \big\downarrow\Gamma \\
\mathcal{C}_i(s_{i-1}(\Gamma), t_{i-1}(\Gamma)) & \overset{s_i}{\underset{t_i}{\longleftarrow\!\!\longleftarrow}} & \mathcal{C}_{i+1}(s_{i-1}(\Gamma), t_{i-1}(\Gamma))
\end{array}$$

In other words, contextualization with $\Gamma$ defines a functor of categories.

A *globular extension* $\mathsf{P}$ *of an* $n$*-sesquicategory* $\mathcal{C}$ is a globular extension for the underlying globular set of $\mathcal{C}$. Given $A \in \mathsf{P}_{n+1}$ and

$$\Gamma \in \mathrm{Cont}\left((s_{n-1}(A), t_{n-1}(A))\right),$$

we define $\Gamma[A]$ as in Eq. (5.4). We write $\mathrm{Cont}(\mathsf{P})$ the set of such $(n+1)$-cells:

$$\mathrm{Cont}(\mathsf{P}) := \left\{\Gamma[A] \mid A \in \mathsf{P}_{n+1} \text{ and } \Gamma \in \mathrm{Cont}\left((s_{n-1}(A), t_{n-1}(A))\right)\right\}.$$

This defines a globular extension of $\mathcal{C}$, also denoted $\mathrm{Cont}(\mathsf{P})$ by abuse of notation, which canonically extends $\mathsf{P}$, in the sense that there is a canonical inclusion $\mathsf{P} \subset \mathrm{Cont}(\mathsf{P})$ that commutes with the source and target maps.

*Remark* 5.2.4 (low-dimensional cases). If $n = 0$ and $\mathsf{P}$ is a globular extension of a set $\mathcal{C}_0$, then $\mathrm{Cont}(\mathsf{P}) = \mathsf{P}$. If $n = 1$ and $\mathsf{P}$ is a globular extension of a category $\mathcal{C}_0 \overset{s_0}{\underset{t_0}{\longleftarrow\!\!\longleftarrow}} \mathcal{C}_1$, then

$$\mathrm{Cont}(\mathsf{P}) = \\ \left\{v \star_{1,2} A \star_{2,1} w \mid A \in \mathsf{P}, v, w \in \mathcal{C}_1, s_0(v) = t_0(A) \text{ and } s_0(A) = t_0(v)\right\}.$$

Diagrammatically:

$$v \star_{1,2} A \star_{2,1} w \;=$$ 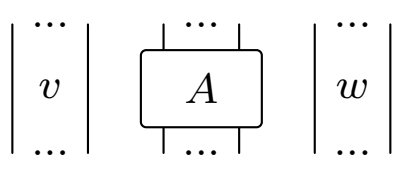

If $n = 2$, then elements of $\mathrm{Cont}(\mathsf{P})$ are as in Eq. (5.1), with $A \in \mathsf{P}$.

### 5.2.4 Free $n$-sesquicategories

Let $\mathcal{C}$ be an $n$-sesquicategory and $\mathsf{P}$ a globular extension of $\mathcal{C}$. Denote by $\mathsf{P}^*_{n+1}$ the set consisting of formal identities $\mathrm{id}_v$ for each $n$-cell $v \in \mathcal{C}_n$ and formal

compositions

$$u_1 \star_{n+1,n+1} u_2 \star_{n+1,n+1} \cdots \star_{n+1,n+1} u_d$$

with $u_i \in \mathrm{Cont}(\mathsf{P})$ and $s(u_i) = t(u_{i-1})$ for $1 \leq i \leq d$. Formal compositions of length zero (i.e. $d = 0$) are identities, and we regard the above up to the usual identity axioms. This defines an $(n+1)$-globular extension of $\mathcal{C}$, denoted $\mathsf{P}^*$, with source and target maps defined as $s(\mathrm{id}_v) = t(\mathrm{id}_v) = v$, and

$$s(u_1 \star_{n+1,n+1} u_2 \star_{n+1,n+1} \cdots \star_{n+1,n+1} u_d) = s(u_d)$$
and
$$t(u_1 \star_{n+1,n+1} u_2 \star_{n+1,n+1} \cdots \star_{n+1,n+1} u_d) = t(u_1).$$

We let $\star_{n+1,n+1} \colon \mathsf{P}^*_{n+1} \times_{\mathcal{C}_n} \mathsf{P}^*_{n+1} \to \mathsf{P}^*_{n+1}$ be the formal juxtaposition of suitably $n$-composable elements of $\mathsf{P}^*_{n+1}$, and $\star_{i,n+1} \colon \mathcal{C}_i \times_{i-1} \mathsf{P}^*_{n+1} \to \mathsf{P}^*_{n+1}$ for $0 < i < n+1$ be defined as $f \star_{i,n+1} \mathrm{id}_a = \mathrm{id}_{f \star_{i,n} a}$ and

$$\begin{aligned} &f \star_{i,n+1} (u_1 \star_{n+1,n+1} u_2 \star_{n+1,n+1} \cdots \star_{n+1,n+1} u_k) \\ &\quad = (f \star_{i,n+1} u_1) \star_{n+1,n+1} (f \star_{i,n+1} u_2) \star_{n+1,n+1} \cdots \star_{n+1,n+1} (f \star_{i,n+1} u_k). \end{aligned}$$

We similarly define $\star_{n+1,i}$. This makes $\mathsf{P}^*$ into an $(n+1)$-sesquicategory, the *free* $(n+1)$*-sesquicategory generated by* $\mathsf{P}$. We sometimes abuse notation and write $\mathsf{P}^*$ to denote the set $\mathsf{P}^*_{n+1}$.

*Remark* 5.2.5 (low-dimensional case). Recall the setting of Remark 5.2.4. If $n = 0$, then $\mathsf{P}^*$ is the free category whose morphisms are formal compositions of elements in $\mathsf{P}$. If $n = 1$, $\mathsf{P}^*$ is the free 2-sesquicategory whose 2-morphisms are formal vertical compositions, or 1-composition, of elements in $\mathrm{Cont}(\mathsf{P})$. If $n = 2$, then $\mathsf{P}^*$ is the free 3-sesquicategory whose 3-morphisms are formal 2-compositions of elements in $\mathrm{Cont}(\mathsf{P})$, which we picture as sequences of movies.

### 5.2.5 $n$-sesquipolygraphs

An $n$*-sesquipolygraph* $\mathsf{P}$ [81, section 2.5] consists of the following data:

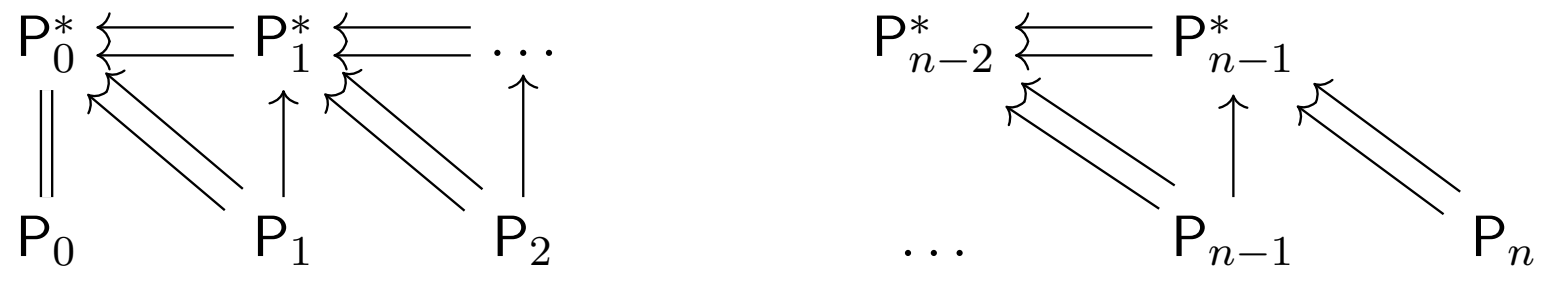

defined inductively as follows:

- $\mathsf{P}_0$ is a set, and $\mathsf{P}_0^* = \mathsf{P}_0$,
- $\mathsf{P}_{i+1}$ is a globular extension for the $i$-globular set $\mathsf{P}_0^* \leftleftarrows \ldots \leftleftarrows \mathsf{P}_i^*$, and $\mathsf{P}_{i+1}^*$ is the free $(i+1)$-sesquicategory generated by $\mathsf{P}_{i+1}$.

An $n$-sesquipolygraph provide a notion of presentation for $n$-sesquicategories, which we now describe.

Given an $n$-sesquicategory $\mathcal{C}$, an equivalence relation $\sim$ on $\mathcal{C}_n$ is said to be *higher* [81, section 2.6][2] if whenever $u \sim v$, we have

- $s(u) = s(v)$ and $t(u) = t(v)$,
- $\Gamma[u] \sim \Gamma[v]$ for each context with boundary $(s(u), t(u)) = (s(v), t(v))$.

Any $(n+1)$-sesquicategory $\mathcal{C}$ defines a higher equivalence relation on the underlying $n$-sesquicategory $\mathcal{C}_{\leq n}$, setting $\sim_{\mathcal{C}}$ to be the smallest higher equivalence relation such that $s(u) \sim_{\mathcal{C}} t(u)$ for all $u \in \mathcal{C}_{n+1}$. We write $[\mathcal{C}]_\sim$ for the $n$-sesquicategory obtained by quotienting $\mathcal{C}_{\leq n}$ with $\sim_{\mathcal{C}}$. if $\mathsf{P}$ is an $(n+1)$-sesquipolygraph, we similarly define $\sim_{\mathsf{P}}$ and $[\mathsf{P}]_\sim$ starting with the $(n+1)$-sesquicategory $\mathsf{P}^*$. Then:

**Definition 5.2.6** ([81, section 2.6]). *A* presentation *of an $n$-sesquicategory $\mathcal{C}$ is the data of an $(n+1)$-prepolygraph $\mathsf{P}$ such that $[\mathsf{P}]_\sim$ is isomorphic to $\mathcal{C}$.*

*Remark* 5.2.7 (low-dimensional cases). A 0-sesquipolygraph is a set. A 1-sesquipolygraph is the same as 1-polygraph, which is the same as a 1-globular set. A 2-sesquipolygraph is the same as a 2-polygraph (see subsection ii.1.2). For $n > 2$, $n$-sesquipolygraphs and $n$-polygraphs are distinct notions.

### 5.2.6 Graded $n$-sesquicategories and their presentations

Let $G$ be an abelian group. We extend all the above to the graded case; setting $G = \{*\}$ recovers the previously introduced notions. A set $\mathbb{P}$ is said to be *graded* if it is equipped with a degree function $\deg\colon \mathbb{P} \to G$. A function between graded sets is *homogeneous* if it preserves the degree functions.

A *graded $n$-globular set* is an $n$-globular set $\mathcal{C}$ such that $\mathcal{C}_n$ is a graded set. A *graded $n$-precategory* is an $n$-precategory whose underlying $n$-globular set is graded, such that $\star_{n,n}$ is additive with respect to the grading, and such that the action of a $k$-cell ($k < n$) on $n$-cells preserves the grading. If $\mathcal{C}$ is an $n$-sesquicategory and $\mathsf{P}$ is a graded extension of $\mathcal{C}$, then $\mathsf{P}^*$ is a graded $n$-sesquicategory, where $\mathsf{P}^*_n$ inherits a grading by setting

$$\deg\left(\Gamma_1[x_1] \star_{n+1,n+1} \ldots \star_{n+1,n+1} \Gamma_m[x_m]\right) \coloneqq \deg(x_1) + \ldots + \deg(x_m).$$

[2] In [81, section 2.6], a higher equivalence relation is called a congruence.

Here $x_k \in \mathsf{P}$ and each $\Gamma_k$ is a context in $\mathcal{C}$. A *graded $n$-sesquipolygraph* $\mathsf{P}$ is the data of a $(n-1)$-prepolygraph $\mathsf{P}_{\leq n-1}$ and a graded extension $\mathsf{P}_n$.

Given a graded $n$-globular set $\mathcal{C}$, a *homogeneous extension* is a graded globular extension $\mathsf{P}$ such that the source and target maps $s_n, t_n \colon \mathsf{P} \to \mathcal{C}_n$ are homogeneous functions. A *graded homogeneous $(n+1)$-sesquipolygraph* $\mathsf{P}$ is the data of a graded $n$-sesquipolygraph $\mathsf{P}_{\leq n}$ and a homogeneous extension $\mathsf{P}_n$. Schematically, a graded homogeneous $(n+1)$-sesquipolygraph is the following data:

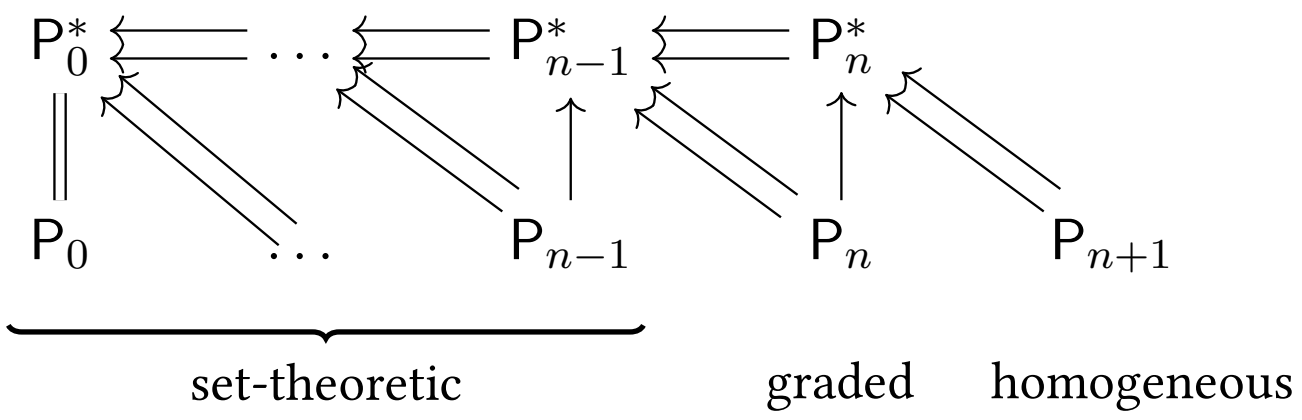

Given a graded set $\mathbb{P}$, an equivalence relation $\sim$ on $\mathbb{P}$ is *homogeneous* if

$$x \sim y \text{ implies } \deg(x) = \deg(y).$$

If a higher equivalence relation $\sim$ on an $n$-precategory $\mathcal{C}$ is homogeneous, the quotient $n$-sesquicategory $[\mathcal{C}]_\sim$ is graded. Note that the higher equivalence relation induced by a graded homogeneous $(n+1)$-sesquipolygraph on the underlying graded $n$-sesquicategory is homogeneous.

**Definition 5.2.8.** *A* presentation *of a graded $n$-sesquicategory is the data of a graded homogeneous $(n+1)$-sesquipolygraph $\mathsf{P}$ such that $[\mathsf{P}]_\sim$ is isomorphic to $\mathcal{C}$.*

### 5.2.7 Scalar $n$-sesquipolygraphs

A set is said to be *scalar* if it is $(\Bbbk, \cdot)$-graded; that is, a scalar set is a set $\mathbb{P}$ equipped with a function $\mathrm{scl} \colon \mathbb{P} \to \Bbbk$. Given a scalar set $\mathbb{P}$, we write $\sim_{\mathrm{scl}}$, or $\equiv_{\mathrm{scl}}$ depending on context, for the smallest homogeneous equivalence relation on $\mathbb{P}$; that is, we have $p \sim_{\mathrm{scl}} q$ if and only if $\mathrm{scl}(p) = \mathrm{scl}(q)$. If the image of scl consists of invertible scalars $\Bbbk^\times \subset \Bbbk$, we say that $\mathbb{P}$ is *scalar-invertible.*

The following are restatements of graded definitions introduced in the previous section.

An $n$-globular set (resp. an $n$-sesquipolygraph) $\mathsf{P}$ is *scalar* if $\mathsf{P}_n$ is a scalar set. An $n$-sesquicategory $\mathcal{C}$ is *scalar* if $\mathcal{C}_n$ is a scalar set and for every $n$-cells $\alpha, \beta$ and context $\Gamma$, we have

$$\mathrm{scl}(\Gamma[\alpha]) = \mathrm{scl}(\alpha) \quad \text{and} \quad \mathrm{scl}(\alpha \star_{n-1} \beta) = \mathrm{scl}(\alpha)\mathrm{scl}(\beta).$$

Given a scalar $n$-sesquipolygraph $\mathsf{P}$, the free $n$-sesquicategory $\mathsf{P}^*$ generated by $\mathsf{P}$ is canonically scalar. Finally, we say that $\mathsf{P}$ is a *scalar-invertible $n$-sesquipolygraph* if $\mathsf{P}_n$ is scalar-invertible.

The linear $n$-sesquicategory presented by a scalar $n$-sesquipolygraph is defined in Definition 5.3.4.

## 5.3 Linear $n$-sesquicategories and their presentations

In this section, we extend the notion of $n$-sesquicategories to the linear case. In fact, we work in the generality of graded linear structures; setting $G = \{*\}$ provides the analogous non-graded linear notions.

*Notation* 5.3.1. We fix throughout the section an abelian group $G$, a commutative ring $\Bbbk$ and a $\mathbb{Z}$-bilinear map $\mu\colon G \times G \to \Bbbk^\times$ as in section 1.1. The word “graded” always refers to $G$-graded, and “linear” to $\Bbbk$-linear. Given a homogeneous element $v$, we write $\deg(v)$ is grading.

### 5.3.1 Linear $n$-sesquicategories

A *graded linear $n$-sesquicategory* is an $n$-sesquicategory $\mathcal{C}$ where, for each $(n-1)$-sphere $(f, g)$, the set $\mathcal{C}_n(f, g)$ has the structure of a graded $\Bbbk$-module, such that the $n$-composition is bilinear and whiskering $n$-cells with a $j$-cell for $j < n$ is linear. In other words:

$$\begin{aligned}(\lambda' u' + v') \star_{n,n} (\lambda u + v) \\ = \lambda'\lambda(u' \star_{n,n} u) + \lambda'(u' \star_{n,n} v) + \lambda(v' \star_{n,n} u) + v' \star_{n,n} v,\end{aligned}$$

$$x \star_{j,n} (\lambda u + v) = \lambda(x \star_{j,n} u) + x \star_{j,n} v,$$

$$(\lambda u + v) \star_{j,n} x = \lambda(u \star_{n,j} x) + v \star_{n,j} x,$$

with scalars $\lambda, \lambda'$ in $\Bbbk$, $(n-1)$-cells $f, g, h$ in $\mathcal{C}_{n-1}$, $n$-cells $u, v$ (resp. $u', v'$) in $\mathcal{C}_n(f, g)$ (resp. $\mathcal{C}_n(g, h)$), and a $j$-cell $x$ suitably $j$-composable.

### 5.3.2 Free linear $n$-sesquicategories

Let $\mathcal{C}$ be an $n$-sesquicategory and $\mathsf{P}$ a graded globular extension of $\mathcal{C}$. The *free graded linear $n$-precategory generated by* $\mathsf{P}$ is the linear $n$-sesquicategory $\mathsf{P}^l$ such that $\mathsf{P}^l_{\leq n} := \mathcal{C}_{\leq n}$ and for each $n$-sphere $(f, g)$ in $\mathcal{C}$, $\mathsf{P}^l_{n+1}(f, g)$ is the free $\Bbbk$-module generated by $\mathsf{P}^*_{n+1}(f, g)$.

The $\Bbbk$-module $\mathsf{P}^l_{n+1}(f, g)$ inherits a grading (of $\Bbbk$-module) from the grading (of set) of the set $\mathsf{P}^*_{n+1}(f, g)$. Extending (bi)linearly the operations $\star_{n+1,n+1}$, $\star_{k,n+1}$ and $\star_{n+1,k}$ $(k \leq n)$ on $\mathsf{P}^*_{n+1}$ defines a structure of graded linear $(n+1)$-sesquicategory on $\mathsf{P}^l$.

### 5.3.3 Linear $n$-sesquipolygraphs

A *graded linear $(n+1)$-sesquipolygraph*[3] $\mathsf{P}$ is the data of a graded $n$-sesquipolygraph $\mathsf{P}_{\leq n}$, together with a homogeneous extension $\mathsf{P}_{n+1}$ of the graded linear $n$-sesquicategory $\mathsf{P}^l_{\leq n}$:

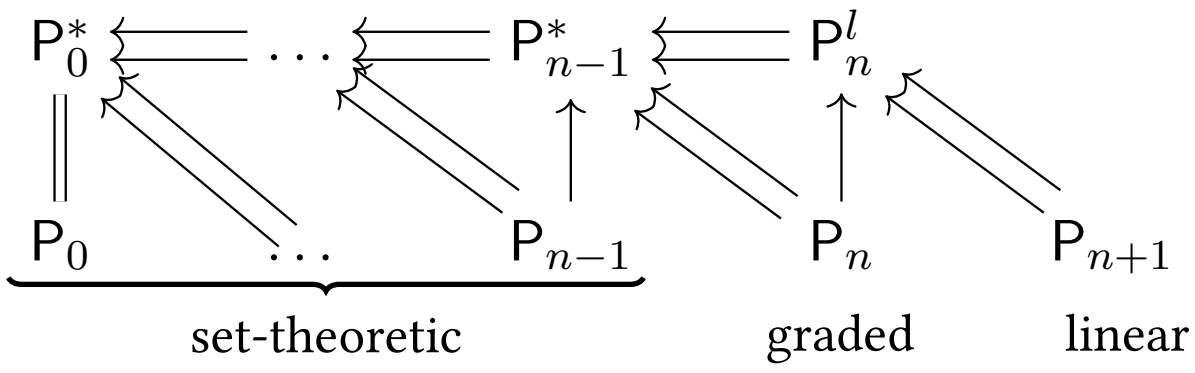

Given a graded $\Bbbk$-module $\mathtt{P}$, an equivalence relation $\sim$ on $\mathtt{P}$ is *linear homogeneous* if whenever $v \sim w$ for $v, w \in \mathtt{P}$, the following two conditions hold:

- if $v$ and $w$ are homogeneous, $\deg(v) = \deg(w)$,
- $\lambda v + u \sim \lambda w + u$ for all scalar $\lambda \in \Bbbk$ and $u \in \mathtt{P}$.

If a higher equivalence relation $\sim$ on a graded linear $n$-sesquicategory $\mathcal{C}$ is linear homogeneous, the quotient $n$-sesquicategory $[\mathcal{C}]_\sim$ is a graded linear

[3] One can more generally define *graded linear $(n, p)$-sesquipolygraphs*, adapting the approach of [3].

$n$-sesquicategory. Note that the higher equivalence relation induced by a graded linear $(n+1)$-sesquipolygraph $\mathsf{P}$ on the underlying graded linear $n$-sesquicategory $\mathsf{P}^l_{\leq n-1}$ is linear homogeneous.

**Definition 5.3.2.** *A presentation of a graded linear $n$-sesquicategory is the data of a graded linear $(n+1)$-sesquipolygraph $\mathsf{P}$ such that $[\mathsf{P}]_\sim$ is isomorphic to $\mathcal{C}$.*

*Remark* 5.3.3 (low-dimensional cases). A linear 1-sesquipolygraph is the same as a linear 1-polygraph, and a linear 2-sesquipolygraph is the same as a linear 2-polygraph. For $n > 2$, linear $n$-sesquipolygraphs and linear $n$-polygraphs [3] are distinct notions.

### 5.3.4 Monomial linear $n$-sesquipolygraphs

A graded linear $(n+1)$-sesquipolygraph $\mathsf{P}$ is *monomial* if:

$$s_n(r) \in \mathsf{P}^*_n \quad \text{and} \quad t_n(r) \in \Bbbk\mathsf{P}^*_n \text{ for all } r \in \mathsf{P}_{n+1},$$

where $\Bbbk\mathsf{P}^*_n$ is the subset of $\mathsf{P}^l_n$ consisting of vectors of the form $\lambda b$ for $\lambda \in \Bbbk$ and $b \in \mathsf{P}^*_n$. We further say it is *monomial-invertible* if $t_n(r) \in \Bbbk^\times\mathsf{P}^*_n$. We have the following canonical bijection:

$$\left\{\begin{array}{c}\text{scalar(-invertible)}\\ \text{graded homogeneous}\\ (n+1)\text{-sesquipolygraph}\end{array}\right\} \underset{\text{scl}}{\overset{\text{lin}}{\rightleftarrows}} \left\{\begin{array}{c}\text{monomial(-invertible)}\\ \text{graded linear}\\ (n+1)\text{-sesquipolygraph}\end{array}\right\}$$

Given a scalar graded homogeneous $(n+1)$-sesquipolygraph $\mathsf{P}$, its *linearization* is the monomial graded linear $(n+1)$-sesquipolygraph $\operatorname{lin}(\mathsf{P})$ defined as $\operatorname{lin}(\mathsf{P})_{\leq n} = \mathsf{P}_{\leq n}$ and $\operatorname{lin}(\mathsf{P})_{n+1} = \mathsf{P}_{n+1}$ as sets, but with the following source and target maps:

$$\operatorname{lin}(s_n)(r) = s_n(r) \quad \text{and} \quad \operatorname{lin}(t_n)(r) = \operatorname{scl}(r)t_n(r).$$

The inverse $\operatorname{scl} := \operatorname{lin}^{-1}$ is defined analogously.

**Definition 5.3.4.** *Given a scalar graded homogeneous $(n+1)$-sesquipolygraph $\mathsf{P}$, the* linear $n$-sesquicategory presented by $\mathsf{P}$ *is the linear $n$-sesquicategory presented by $\operatorname{lin}(\mathsf{P})$.*

## 5.4 Gray polygraphs and linear Gray polygraphs

This section reviews the notion of a *Gray polygraph*, as can be extracted from [81], and then introduces its linear analogue called a *linear Gray polygraph.*

As we shall not explicitly need it, the definition of a Gray category is given in an appendix Appendix A, together with a proof of coherence for interchangers.

### 5.4.1 Gray polygraphs

Let $\mathsf{Q}$ be a 2-sesquipolygraph. The *3-sesquipolygraph of interchangers* is the 3-sesquipolygraph $\mathsf{QGray}$ such that $\mathsf{QGray}_{\leq 2} = \mathsf{Q}$ and $\mathsf{QGray}_3$ consists of *interchange generators*, defined for each 0-composable $\alpha, f, \beta$ with $\alpha\colon f \Rightarrow f'$, $g \in \mathsf{P}_1^*$, $\beta\colon h \Rightarrow h' \in \mathsf{P}_2$, as the 3-cell

$$X_{\alpha,g,\beta}\colon (\alpha \star_0 g \star_0 h') \star_1 (f \star_0 g \star_0 \beta) \Rrightarrow (f' \star_0 g \star_0 \beta) \star_1 (\alpha \star_0 g \star_0 h),$$

pictured as:

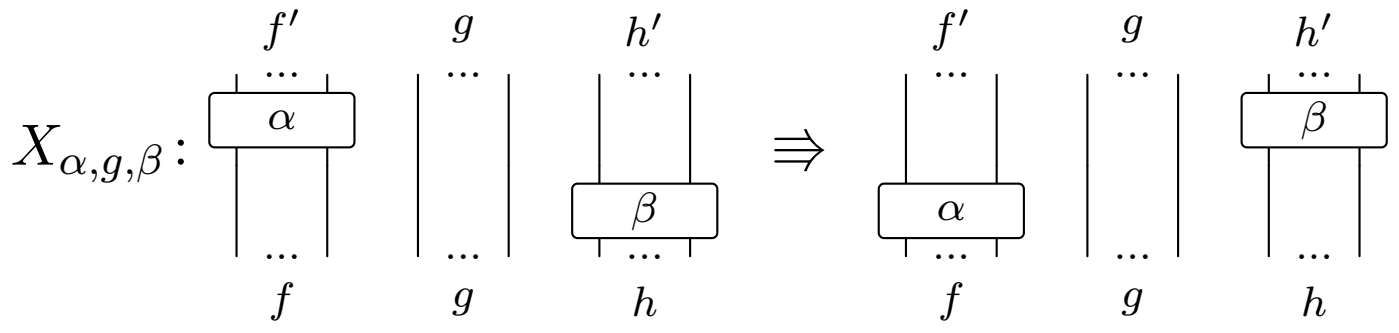

**Definition 5.4.1.** *A* Gray polygraph[4] *is a 3-sesquipolygraph* $\mathsf{P}$ *such that*

$$\mathsf{P}_{\leq 2}\mathsf{Gray} \subset \mathsf{P}.$$

*In other words,* $\mathsf{P}$ *contains its own 3-sesquipolygraph of interchangers.*

One checks that if $\mathsf{P}$ is a Gray polygraph, then $[\mathsf{P}]_\sim$ is a 2-category. This leads to the following definition:

**Definition 5.4.2.** *A presentation of a 2-category $\mathcal{C}$ is the data of a Gray polygraph* $\mathsf{P}$ *such that* $[\mathsf{P}]_\sim$ *is isomorphic to* $\mathcal{C}$.

[4]What we call a Gray polygraph is the underlying 3-sesquipolygraph (or 3-prepolygraph in their terminology) of what is called a Gray presentation in [81].

### 5.4.2 Linear Gray polygraphs

Recall Notation 5.3.1. Let $\mathsf{Q}$ be a graded 2-sesquipolygraph. The *3-sesquipolygraph of $(G,\mu)$-graded interchangers* is the scalar-invertible 3-sesquipolygraph $\mathsf{QGray}$, equipped with the function

$$\begin{aligned} \mathsf{QGray} &\to \Bbbk^{\times}, \\ X_{\alpha,g,\beta} &\mapsto \mu(\deg\alpha,\deg\beta). \end{aligned}$$

We abuse notation and similarly denote $\mathsf{QGray}$ the associated monomial-invertible linear $n$-sesquipolygraph.

**Definition 5.4.3.** *A $(G,\mu)$*-linear Gray polygraph *is a graded linear 3-sesquipolygraph* $\mathsf{P}$ *such that* $\mathsf{P}$ *contains its own monomial graded 3-sesquipolygraph of $(G,\mu)$-graded interchangers.*

Similarly, a $(G,\mu)$*-scalar Gray polygraph* is a scalar graded 3-sesquipolygraph $\mathsf{P}$ such that $\mathsf{P}$ contains its own scalar graded 3-sesquipolygraph of $(G,\mu)$-graded interchangers. In other words, a $(G,\mu)$-scalar Gray polygraph is precisely the same data as a monomial $(G,\mu)$-linear Gray polygraph.

One checks that if $\mathsf{P}$ is a $(G,\mu)$-linear Gray polygraph, then $[\mathsf{P}]_\sim$ is a $(G,\mu)$-graded-2-category. This leads to the following definition:

**Definition 5.4.4.** *A presentation of a $(G,\mu)$-graded-2-category $\mathcal{C}$ is the data of a $(G,\mu)$-linear Gray polygraph* $\mathsf{P}$ *such that* $[\mathsf{P}]_\sim$ *is isomorphic to $\mathcal{C}$.*

*Remark* 5.4.5. If $G=\{*\}$ is trivial, a $(\{*\},\mathrm{id})$-graded-2-category is just a linear 2-category. Hence, a $(\{*\},\mathrm{id})$-linear Gray polygraph, and in particular a $(\{*\},\mathrm{id})$-scalar Gray polygraph, defines a presentation of a linear 2-category. This is the case of the linear monoidal category $\mathbf{Web}_\otimes$ discussed in section 5.1, presented by the scalar $(\{*\},\mathrm{id})$-Gray polygraph $\mathsf{Web}$.

# 6

# Linear Gray rewriting modulo

This chapter introduces the rewriting theory of graded-2-categories. It has three main features: it allows *modulo*, it is *linear*, and it is *higher*. Each of these features requires its own treatment, in addition with a description of how they combine. We refer to subsections ii.1.1 and ii.1.3 to ii.1.6 for an introduction to the basic ideas, to subsection ii.2.2 for an overview of the chapter, and to subsection ii.2.3 for special cases and relationship to the literature.

Rewriting modulo aims to provide a presentation of a given object which is convergent modulo. Hence, we start in section 6.1 by describing which information can be gathered from convergence modulo. In particular, we introduce the notion of *coherence modulo*. Similar ideas have appeared in [58]: it would be interesting to better understand the relationship with our work.

Section 6.2 then introduces the fundamentals of abstract rewriting theory modulo, which underlies any rewriting theory modulo. This can be understood as an extension to modulos of classical abstract rewriting theory (see e.g. [5, chap. 1]). In addition, we introduce the notion of *tamed congruence* (see Definition 6.2.6) as a replacement for the usual notion of confluence. As far as the author is aware, this has not appeared in the literature.

In section 6.3, we extend abstract rewriting theory modulo to the linear setting. We keep the theory as general as possible, allowing the modulo data to rewrite basis element up to (not necessarily coherent) invertible scalars. Our theory can be viewed as an extension to modulos of the work of Guiraud–Hoffbeck–Malbos [86]. In addition, we emphasize the use of tamed congruence as a natural replacement for confluence (this point is not specific to the modulo context).

Section 6.4 generalizes abstract rewriting modulo to the higher setting. We emphasize that a higher rewriting system (modulo or not) is nothing else than a category of abstract rewriting systems, where morphisms are given by contexts (see subsection 5.2.3). Hence, understanding higher rewriting reduces to understanding contextualization.

Finally, section 6.5 combines all the previous sections into a theory of higher linear rewriting modulo. We emphasize the difficulty of dealing with contextualization in the linear setting, and how tamed congruence provides an answer. The latter point applies regardless whether one works modulo or not, or whether one works with graded-2-categories or genuine 2-categories.

A quick summary of the main techniques is given in section 6.6. In subsection 6.6.3, we detail the rewriting theory of $\mathfrak{gl}_2$-webs, as introduced in section 5.1.

The works of Alleaume [3] (rewriting in linear strict 2-categories) and Dupont [55, 56] (rewriting modulo in linear strict 2-categories) are important inspirations for our work. We shall point the differences with their approaches throughout the text, especially those that remain even when specifying to their respective contexts.

We remind the reader that every categorical structure is assumed to be small.

*Notation* 6.0.1. We fix the same notations as in Notation 5.3.1 throughout the chapter. As much as possible, we follow the notational conventions introduced in Notation 5.0.1. Furthermore, we use blackboard fonts (e.g. $\mathbb{P}$) for abstract rewriting systems and typewriter fonts (e.g. $\mathtt{P}$) for linear rewriting systems. Equivalence relations are denoted either with the symbol $\sim$ or the symbol $\equiv$; the former case is typically used for congruence of a rewriting system, while the latter is typically used for an equivalence relation on rewriting sequences (see section 6.2 for both notions). In particular, $\equiv$ is typically a "categorical dimension" higher than $\sim$.

## 6.1 Coherence modulo from convergence modulo

Given a groupoid $\mathcal{P}$, one may be interested in understanding its set of connected components $\pi_0(\mathcal{P})$,[1] or its $\pi_1(\mathcal{P})$, that is, its *coherence*. In this section, we describe how the property of convergence modulo answers both questions.

Let $\mathcal{S}$ be a category. The notion of higher equivalence relation on $n$-sesqui-categories defined in subsection 5.2.5 specializes to categories: an equivalence relation $\equiv$ on the set of morphisms, such that if $f \equiv f'$, then (i) $s(f) = s(f')$ and $t(f) = t(f')$, and (ii) if $g$ and $h$ are suitably composable morphisms, then $g \circ f \circ h \equiv g \circ f' \circ h$. In this context, we call such an equivalence relation an *abstract equivalence*, and say that two morphisms $f$ and $g$ are $\equiv$-equivalent if $f \equiv g$. We choose this terminology to avoid confusion with the notion of $\equiv$-congruence, used extensively in the following sections. In particular, if $\equiv = \equiv_{\text{dis.}}$ is the discrete abstract equivalence (see Definition 6.1.5), $\equiv$-equivalence reduces to equality, while $\equiv$-congruence reduces to congruence (see the next sections).

We denote $\mathcal{S}^\top$ the localisation of the category $\mathcal{S}$. Given a groupoid $\mathcal{P}$, we fix the following data for the section:

$$\begin{gathered}\text{a category } \mathcal{S} \text{ such that } \mathcal{S}^\top = \mathcal{P}, \text{ a wide}^2 \text{ subgroupoid } \mathcal{E} \subset \mathcal{S}, \\ \text{and an abstract equivalence} \equiv \text{on } \mathcal{S}.\end{gathered} \tag{$*$}$$

We use the same notation $\equiv$ for the smallest abstract equivalence on $\mathcal{S}^\top$ containing $\equiv$.[3] Morphisms in $\mathcal{S}$ are depicted with plain arrow $x \to y$ (or $x \longrightarrow y$). Morphisms in $\mathcal{E}$ are depicted with unoriented wiggly lines $x \sim y$ (or $x \rightsquigarrow y$). Unspecified morphisms $x \to_{\mathcal{S}} y$ (resp. $x \sim_{\mathcal{E}} y$) indicate the mere existence of a morphism in $\mathcal{S}$ (resp. in $\mathcal{E}$) between $x$ and $y$. Note the use of subscripts to specify the category, used extensively in what follows.

Recall from subsection 5.2.7 the notion of scalar category and of its associated abstract equivalence $\equiv_{\text{scl}}$.

[1]That is, its Grothendieck group $K_0(\mathcal{P})$.
[2]Containing all objects; equivalently, containing all identities.
[3]In particular, if $f \equiv g$ in $\mathcal{S}$, then we both have $f \equiv g$ and $f^{-1} \equiv g^{-1}$ in $\mathcal{S}^\top$.

**Definition 6.1.1.** *The category $\mathcal{P}$ is said to be* $\equiv$-coherent modulo $\mathcal{E}$ *if every endomorphism is $\equiv$-equivalent to a morphism in $\mathcal{E}$, up to conjugation. That is, for every endomorphism $f\colon x \to_{\mathcal{P}} x$, there exist morphisms $g\colon x \to_{\mathcal{P}} y$ and $e\colon y \sim_{\mathcal{E}} y$, such that $f \equiv g^{-1} \circ e \circ g$:*

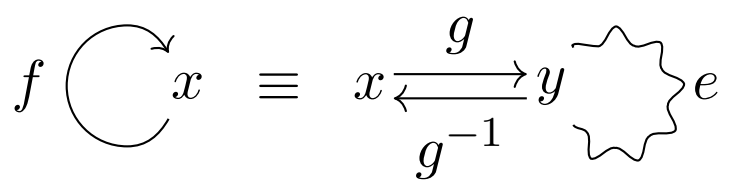

*If $\mathcal{P}$ is scalar, we say that it is* scalar-coherent modulo $\mathcal{E}$ *whenever it is $\equiv_{\mathrm{scl}}$-coherent modulo $\mathcal{E}$.*

*Remark* 6.1.2. If $\mathcal{E}$ is discrete (it does not contain any morphism; we write $\mathcal{E} = \emptyset$), then $\mathcal{P}$ is $\equiv$-coherent modulo $\mathcal{E}$ if and only if parallel morphisms are $\equiv$-equivalent, which recovers the usual notion of coherence for a category. In particular, if $\mathcal{P}$ is scalar with scalars $\mathrm{scl}\colon \mathcal{P}_1 \to \mathbb{k}$, then $\mathcal{P}$ is scalar-coherent if for every endomorphism $f\colon x \to_{\mathcal{P}} x$ we have $\mathrm{scl}(f) = 1$.

**Lemma 6.1.3** (transitivity of coherence modulo). *Let $\mathcal{C}$ be a groupoid, $\equiv$ an abstract equivalence on $\mathcal{C}$ and $\mathcal{E} \subset \mathcal{D}$ wide subgroupoids of $\mathcal{C}$. If $\mathcal{C}$ is $\equiv$-coherent modulo $\mathcal{D}$ and $\mathcal{D}$ is $\equiv$-coherent modulo $\mathcal{E}$, then $\mathcal{C}$ is $\equiv$-coherent modulo $\mathcal{E}$.* $\square$

An $\mathcal{S}$*-branching* is a pair $(f, f')$ with $f\colon x \to_{\mathcal{S}} y$ and $f'\colon x \to_{\mathcal{S}} y'$. An $\mathcal{S}$*-confluence* is a pair $(g, g')$ with $g\colon y \to_{\mathcal{S}} z$ and $g'\colon y' \to_{\mathcal{S}} z$. A branching $(f, f')$ is said to be $\equiv$*-confluent* if there exists a confluence $(g, g')$ such that $f' \circ f \equiv g' \circ g$:

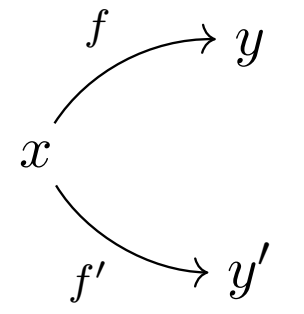

an $\mathcal{S}$-branching

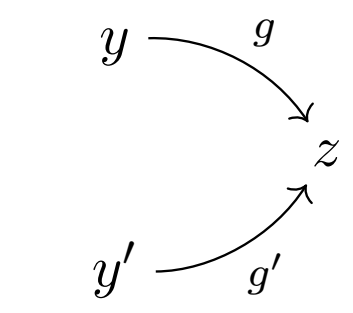

an $\mathcal{S}$-confluence

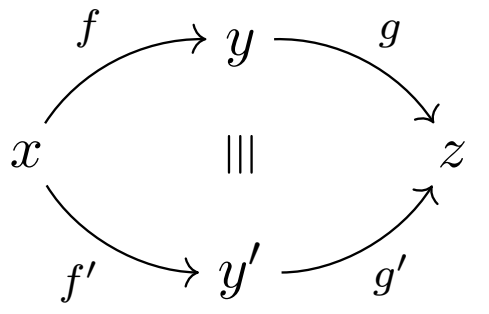

an $\mathcal{S}$-branching which is $\equiv$-confluent

We say that $\mathcal{S}$ is $\equiv$*-confluent* if every branching is $\equiv$-confluent.

We say that $\mathcal{S}$ is *terminating* if every infinite sequence $(f_n)_{n\in\mathbb{N}}$ of morphisms in $\mathcal{S}$ with $t(f_n) = s(f_{n+1})$ eventually terminates in morphisms in $\mathcal{E}$. We say that $\mathcal{S}$ is $\equiv$*-convergent* if it is both $\equiv$-confluent and terminating. An object $y \in \mathcal{S}$ for which $y \to_{\mathcal{S}} z$ implies $y \sim_{\mathcal{E}} z$ is called an $\mathcal{S}$*-normal form.* We denote $\mathrm{NF}_{\mathcal{S}}$ the set of $\mathcal{S}$-normal forms. If $x \to_{\mathcal{S}} y$ and $y$ is an $\mathcal{S}$-normal form, we say that $y$ is an $\mathcal{S}$-normal form for $x$. Note that the notions of termination, convergence and normal form all depend on the data of $\mathcal{E}$. To stress this dependency, we sometimes add *modulo* $\mathcal{E}$, as in *terminating modulo* $\mathcal{E}$.

The *Church–Rosser property* is a classical property of $\equiv$-confluent categories (see for instance [5, p. 1.3.14]):

**Lemma 6.1.4** (Church–Rosser property). *Let $\mathcal{S}$, $\equiv$ and $\mathcal{E}$ as in $(*)$. If $\mathcal{S}$ is $\equiv$-confluent, then every morphism $f \in \mathcal{P} = \mathcal{S}^\top$ can be decomposed as $f \equiv h^{-1} \circ g$ for some $g, h \in \mathcal{S}$.*

*Proof.* By definition, $f$ decomposes as $f = h_1 \circ g_1^{-1} \circ h_2 \circ g_2^{-1} \circ \ldots \circ h_k \circ g_k^{-1}$ for some $k \in \mathbb{N}$ and $g_i, h_i \in \mathcal{S}$. Inductively applying $\equiv$-confluence leads to the desired decomposition. □

**Definition 6.1.5.** *The* discrete abstract equivalence *on $\mathcal{S}$ is the abstract equivalence $\equiv_{\text{dis.}}$ such that every pair of parallel morphisms is equivalent with respect to $\equiv_{\text{dis.}}$.*

We simply denote $\equiv_{\text{dis.}}$-confluence as *confluence*, and $\equiv_{\text{dis.}}$-convergence modulo as *convergence modulo*. Note that $\equiv$-confluence (resp. $\equiv$-convergence modulo) for any abstract equivalence $\equiv$ implies confluence (resp. convergence modulo).

Finally, the following two propositions explain how convergence modulo leads to a description of $\pi_0(\mathcal{P})$ and $\pi_1(\mathcal{P})$, respectively.

**Proposition 6.1.6.** *Let $\mathcal{P}$, $\mathcal{S}$ and $\mathcal{E}$ as in $(*)$. If $\mathcal{S}$ is convergent modulo $\mathcal{E}$, then the canonical mapping $\mathrm{NF}_\mathcal{S} \to \pi_0(\mathcal{P})$ sending an $\mathcal{S}$-normal form to its connected component induces a bijection*

$$\mathrm{NF}_\mathcal{S}/\pi_0(\mathcal{E}) \xrightarrow{\sim} \pi_0(\mathcal{P}),$$

*where the quotient identifies objects belonging to the same connected component of $\mathcal{E}$.*

*Proof.* If we have $f\colon x \to_\mathcal{S} y$, then $x$ and $y$ belong to the same connected component of $\mathcal{P}$. In particular, since $\mathcal{S}$ is terminating the mapping $\mathrm{NF}_\mathcal{S} \to \pi_0(\mathcal{P})$ is surjective. Assume then that $x$ and $y$ are $\mathcal{S}$-normal forms belonging to the same connected component of $\mathcal{P}$. By the Church–Rosser property (Lemma 6.1.4), there exists $f\colon x \to_\mathcal{S} z$ and $g\colon y \to_\mathcal{S} z$. Since $x$ and $y$ are $\mathcal{S}$-normal forms, $f$ and $g$ must be in $\mathcal{E}$, so that $x$ and $y$ belong to the same connected component of $\mathcal{E}$. □

**Proposition 6.1.7.** *Let* $\mathcal{S}$, $\equiv$ *and* $\mathcal{E}$ *as in* $(*)$. *If* $\mathcal{S}$ *is* $\equiv$*-convergent modulo* $\mathcal{E}$*, then* $\mathcal{P} = \mathcal{S}^\top$ *is* $\equiv$*-coherent modulo* $\mathcal{E}$*.*

*Proof.* Thanks to the Church–Rosser property (Lemma 6.1.4), for any endomorphism $f\colon x \to_{\mathcal{P}} x$ of $\mathcal{P}$ there exist $f_1, f_2\colon x \to_{\mathcal{S}} z$ morphisms in $\mathcal{S}$ such that $f \equiv f_1^{-1} \circ f_2$. Since $\mathcal{S}$ is terminating, there exists an $\mathcal{S}$-normal form $y$ and a morphism $h\colon z \to_{\mathcal{S}} y$. Because $\mathcal{S}$ is $\equiv$-confluent, the branching $(h \circ f_1, h \circ f_2)$ admits a $\equiv$-confluence. Since $y$ is an $\mathcal{S}$-normal form, this confluence is in $\mathcal{E}$, and there exists $e \in \mathcal{E}$ such that $e \circ h \circ f_1 \equiv h \circ f_2$. Setting $g := h \circ f_1$ gives $e \circ g \equiv g \circ f$, which concludes. □

## 6.2 Abstract rewriting modulo

Recall that 1-polygraph, 1-sesquipolygraph and 1-globular set are identical notions (Remark 5.2.7). Recall also the notion of the free category $\mathbb{P}^*$ generated by a 1-polygraph $\mathbb{P}$ (subsection 5.2.4). We denote $\mathbb{P}^\top$ the localisation of $\mathbb{P}^*$.

**Definition 6.2.1.** *An* abstract rewriting system (abstract RS) *is the data* $(\mathbb{P}; \equiv)$ *of a 1-polygraph* $\mathbb{P}$ *together with an abstract equivalence* $\equiv$ *on* $\mathbb{P}^*$*.*

Unpacking the definition, an abstract RS is the data of a set $X$ of 0-cells, called *elements*, and a set $\mathbb{P}$ of generating 1-cells, called *rewriting steps*, equipped with *source* and *target maps*:

$$X \overset{s}{\underset{t}{\leftleftarrows}} \mathbb{P},$$

together with an abstract equivalence $\equiv$ on the free category $\mathbb{P}^*$ generated by $\mathbb{P}$. Note that we abuse notation and denote $\mathbb{P}$ both the 1-polygraph and the set of relations. A morphism in $\mathbb{P}^*$ (resp. in $\mathbb{P}^\top$) is called a *rewriting sequence* (resp. a *congruence*). If a rewriting sequence (resp. a congruence) decomposes in $n$ $\mathbb{P}$-rewriting steps, we call the number $n$ its *length.* We write $x \to_{\mathbb{P}} y$ to refer to an unspecified rewriting step with source $x$ and target $y$ (note that there could be more than one rewriting step between given source and target), or to indicate the existence of such a rewriting step. Similarly, we write $x \overset{*}{\to}_{\mathbb{P}} y$ (resp. $x \sim_{\mathbb{P}} y$) to denote a rewriting sequence (resp. a congruence), and say that $x$ *rewrites into* (resp. *is congruent to*) $y$.

| $x \to_{\mathbb{P}} y$ | $x \overset{*}{\to}_{\mathbb{P}} y$ | $x \sim_{\mathbb{P}} y$ |
|---|---|---|
| rewriting step | rewriting sequence | congruence |

Note that following these notations, $x \xrightarrow{*}_{\mathbb{P}} y$ coincides with $x \to_{\mathbb{P}^*} y$, and $x \sim_{\mathbb{P}} y$ coincides with $x \to_{\mathbb{P}^\top} y$.

Let $\mathbb{R} = (X, \mathbb{R})$ and $\mathbb{E} = (X, \mathbb{E})$ be two 1-polygraphs with the same underlying set of elements. We define the following 1-polygraph:

$$_{\mathbb{E}}\mathbb{R}_{\mathbb{E}} := \mathbb{E}^\top \times_X \mathbb{R} \times_X \mathbb{E}^\top.$$

In other words, a rewriting step in $_{\mathbb{E}}\mathbb{R}_{\mathbb{E}}$ is a triple $(e, r, e') \in \mathbb{E}^\top \times \mathbb{R} \times \mathbb{E}^\top$ with $t(e) = s(r)$ and $t(r) = s(e')$. The source and target maps are defined as $s(e, r, e') = s(e)$ and $t(e, r, e') = t(e')$.

**Definition 6.2.2.** *An* abstract rewriting system modulo (abstract RSM) *is the data* $\mathbb{S} = (\mathbb{R}, \mathbb{E}; \equiv)$ *of two 1-polygraphs* $\mathbb{R} := (X, \mathbb{R})$ *and* $\mathbb{E} = (X, \mathbb{E})$*, together with an abstract equivalence* $\equiv$ *on* $(_{\mathbb{E}}\mathbb{R}_{\mathbb{E}})^* \cup \mathbb{E}^\top$*. In that case, we say that* $\mathbb{S}$ *is an* abstract RSM modulo $\mathbb{E}$ on the set $X$.[4]

An $\mathbb{S}$*-rewriting sequence* is either an $\mathbb{E}$-congruence (in which case it has *length zero*) or an $_{\mathbb{E}}\mathbb{R}_{\mathbb{E}}$-rewriting sequence (in which case it has the same length as the $_{\mathbb{E}}\mathbb{R}_{\mathbb{E}}$-rewriting sequence).[5] We denote $\mathbb{S}^*$ the set of $\mathbb{S}$-rewriting sequences. Note that $(_{\mathbb{E}}\mathbb{R}_{\mathbb{E}})^* \cup \mathbb{E}^\top = \mathbb{S}^*$. We similarly denote $x \to_{\mathbb{S}} y$ (resp. $x \xrightarrow{*}_{\mathbb{S}} y$, resp. $x \sim_{\mathbb{S}} y$) an $\mathbb{S}$-rewriting step (resp. an $\mathbb{S}$-rewriting sequence, resp. an $\mathbb{S}$-congruence).

An abstract RSM inherits all the notions and results introduced in the previous section, setting $(*)$ as

$$\mathcal{S} := \mathbb{S}^* \quad \text{and} \quad \mathcal{E} := \mathbb{E}^\top.$$

For instance, we say that $\mathbb{S}$ is *terminating* (or *terminating modulo* $\mathbb{E}$) if $\mathcal{S}$ is terminating modulo $\mathcal{E}$; equivalently, $\mathbb{S}$ is terminating if there is no infinite sequence $\{f_n\}_{n\in\mathbb{N}}$ such that $f_n$ is an $\mathbb{S}$-rewriting step and $t(f_n) = s(f_{n+1})$. (Note that with our conventions, one *cannot* replace “step” by “sequence” in the above characterization.)

However, in contrast with the previous section, an abstract RSM provides a notion of *locality*. A *local* $\mathbb{S}$*-branching* is a pair $(f, g)$ with $f\colon x \to_{\mathbb{S}} y$ and

[4] In [56, 58], the authors allow for a more general definition, where $\mathbb{S}$ is any abstract RS such that $\mathbb{R} \subset \mathbb{S} \subset {}_{\mathbb{E}}\mathbb{R}_{\mathbb{E}}$. We do not work in this generality. Moreover, we impose $\mathbb{R}$ and $\mathbb{E}$ to be defined on the same set of elements; see Footnote 16 for further comments on this.

[5] Our terminology differs from [55, 56], where $\mathbb{S}$-rewriting sequences coincide with $_{\mathbb{E}}\mathbb{R}_{\mathbb{E}}$-rewriting sequences. In particular, our notions of branching and confluence do not explicitly depend on the modulo data, while the notion of termination does; this is the converse of [55, 56].

$f\colon x \to_{\mathbb{S}} y'$ two rewriting steps in $\mathbb{S}$:

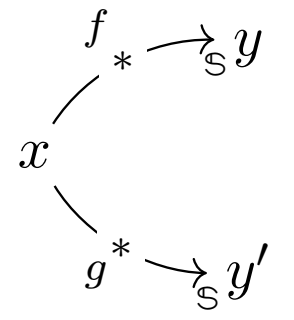

an $\mathbb{S}$-branching

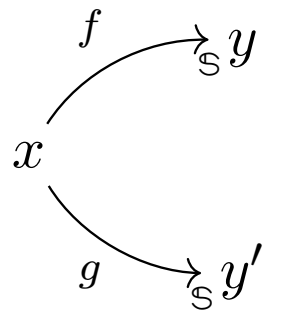

a local $\mathbb{S}$-branching

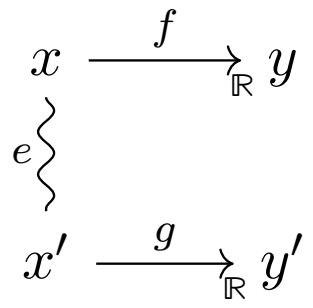

an $\mathbb{S}$-local triple

We say that $\mathbb{S}$ is *locally* $\equiv$*-confluent* if all local branchings are $\equiv$-confluent. An $\mathbb{S}$*-local triple* is a triple $[f, e, g]$ with $\mathbb{R}$-rewriting steps $f$ and $g$ and $\mathbb{E}$-congruence $e$, such that $s(f) = s(e)$ and $t(e) = s(g)$. Note that every $\mathbb{S}$-local triple defines a local $\mathbb{S}$-branching $(f, g \circ e)$, and that if every $\mathbb{S}$-local triple is $\equiv$-confluent, then $\mathbb{S}$ is $\equiv$-confluent. Local $\mathbb{S}$-branchings and $\mathbb{S}$-local triples should be thought as essentially identical notions, one being more suited for general statements while the other being better suited for explicit computations.

Given an abstract RSM $\mathbb{S} = (\mathbb{R}, \mathbb{E}, \equiv)$, an abstract RSM $\mathbb{T} = (\mathbb{R}_1, \mathbb{E}, \equiv_1)$ is a *sub-abstract RSM of* $\mathbb{S}$ if $\mathbb{R}_1 \subset \mathbb{R}$ and for all $f, g \in \mathbb{T}^*$, $f \equiv_1 g$ implies that $f \equiv g$. If an $\mathbb{S}$-branching $(f, g)$ admits a $\mathbb{T}$-congruence $(f', g')$ such that $f' \circ f \equiv g' \circ g$, we say that $(f, g)$ is $(\mathbb{T}, \equiv)$*-confluent*. We say that an $\mathbb{S}$-congruence $h$ is $(\mathbb{T}, \equiv)$*-confluent* (resp. $(\mathbb{T}, \equiv)$*-congruent*) if there exists a $\mathbb{T}$-confluence $(f, g)$ with $s(f) = s(h)$, $s(g) = t(h)$, such that $f \equiv g \circ h$ (resp. a $\mathbb{T}$-congruence $h'$ such that $h \equiv h'$). We use similar notations to indicate that a given notion related to $\mathbb{S}$ restricts to a notion related to $\mathbb{T}$.

*Remark* 6.2.3 (scalar abstract RSM). Recall from subsection 5.2.7 the notion of scalar 1-polygraph, scalar relation and free scalar 1-category. An abstract RSM $\mathbb{S} = (\mathbb{R}, \mathbb{E}, \equiv)$ is said to be *scalar* if $\mathbb{R}$ is scalar, $\mathbb{E}$ is scalar-invertible and $\equiv \, = \, \equiv_{\mathrm{scl}}$ (that is, $r \equiv s$ if and only if $\mathrm{scl}(r) = \mathrm{scl}(s)$).

*Remark* 6.2.4 (quotient). Given an abstract RSM $\mathbb{S} = (X; \mathbb{R}, \mathbb{E}; \equiv)$, we can define its *quotient* abstract RS $[\mathbb{S}]^{\mathbb{E}} := ([X]^{\mathbb{E}}; [\mathsf{R}]^{\mathbb{E}}; [\equiv]^{\mathbb{E}})$, where $[X]^{\mathbb{E}}$ is the set of $\mathbb{E}$-congruence classes and $[\mathsf{R}]^{\mathbb{E}}$ and $[\equiv]^{\mathbb{E}}$ are defined analogously. In that case, $\mathbb{S}$ is terminating modulo $\mathbb{E}$ if and only if $[\mathbb{S}]^{\mathbb{E}}$ is terminating, and if $\mathbb{E}$ is $\equiv$-coherent, then $\mathbb{S}$ is $\equiv$-confluent if and only if $[\mathbb{S}]^{\mathbb{E}}$ is $[\equiv]^{\mathbb{E}}$-confluent. In other words, *when* $\mathbb{E}$ *is coherent*, the theory of rewriting modulo reduces to rewriting (without modulo) on its quotient.

However, as we will see in section 6.5 (see also subsection ii.2.2), explicitly working with the modulo data enlights some of the difficulties of higher linear rewriting. Moreover, one does not always have the luxury of a coherent modulo data; indeed, the modulo data for graded $\mathfrak{gl}_2$-foams (chapter 7) is not coherent.

The next subsections describe how confluence can be achieved from a local analysis.

### 6.2.1 Tamed Newmann's lemma

The classical Newmann's lemma states that under termination, confluence follows from local confluence. It readily extends to modulos:

**Lemma 6.2.5** (Newmann's lemma). *Let $\mathbb{S} = (\mathbb{R}, \mathbb{E}; \equiv)$ be an abstract RSM. If $\mathbb{S}$ is a terminating and locally $\equiv$-confluent, then it is $\equiv$-confluent.* □

However, local confluence turns out to be too restrictive for some purposes, especially in the linear context. In this subsection, we introduce the weaker notion of *$\succ$-tamed $\equiv$-congruence*, for which an analogue of Newmann's lemma still holds.

Recall that a binary transitive relation on $X$ is called a *preorder*. Let $\mathbb{E} = (X, \mathbb{E})$ be an abstract RS. A preorder $\succ$ is said to be *$\mathbb{E}$-invariant* if ($x' \sim_{\mathbb{E}} x$ and $x \succ y$ and $y \sim_{\mathbb{E}} y'$) implies ($x' \succ y'$). We always implicitly assume that preorders on $X$ are $\mathbb{E}$-invariant. If $M \subset \mathcal{P}(X)$ is a set of elements in $X$, we write $x \succ M$ if $x \succ y$ for all $y \in M$. If $f = f_n \circ \ldots \circ f_1$ is a sequence of composable arrows on $X$, we write $x \succ f$ to mean $x \succ \{s(f_1), t(f_1), \ldots, t(f_n)\}$.

**Definition 6.2.6.** *Let $\mathbb{S} = (\mathbb{R}, \mathbb{E}; \equiv)$ be an abstract RSM on the set $X$ and $\succ$ a preorder on $X$. An $\mathbb{S}$-branching $(f, g)$ of source $\bullet$ is said to be $\succ$-*tamely $\equiv$-congruent *(resp. $\succ$-*tamely $\equiv$-confluent*) if there exists a $\equiv$-congruence $h$ (resp. $\equiv$-confluence $(f', g')$) such that $\bullet \succ h$ (resp. $\bullet \succ f'^{-1} \circ g'$).*

In particular, $\succ$-tameness implies $\bullet \succ t(f)$ and $\bullet \succ t(g)$. Here is a schematic for a $\succ$-tamed congruence, where horizontal positions are used to suggest relative orderings with respect to $\succ$:

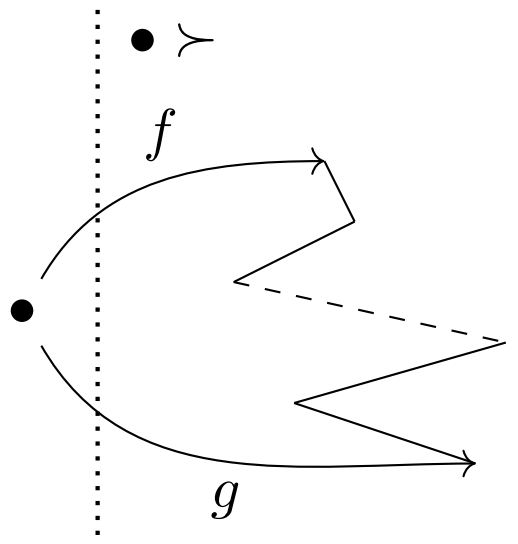

Our notion of $\succ$-tameness is reminiscent of the notion of *confluence by decreasingness* as introduced by van Oostrom [182], and as appearing in [3, 55, 56] in the context of higher linear rewriting. Indeed, $\mathbb{S}$-rewriting steps inherit a preorder from $\succ$ by stating that $h_1 \succ h_2$ if and only if $s(h_1) \succ s(h_2)$ and $s(h_1) \succ t(h_2)$. With this choice, a $\succ$-tame $\equiv$-confluence is decreasing in the sense of [182]. However, moving the order from rewriting steps to elements allows a meaningful weakening of the notion to congruence. As far as we are aware, this has not appeared in the literature.

**Definition 6.2.7.** *Let $\mathbb{S}$ be an abstract RSM. A preorder $\succ$ on $X$ is said to be* compatible with $\mathbb{S}$ *if $x \to_{\mathbb{S}} y$ implies $x \succ y$.*

We denote $\succ_{\mathbb{S}}$ the minimal preorder compatible with $\mathbb{S}$. Note that if $\succ$ is compatible with $\mathbb{S}$, then minimal elements for $\succ$ are $\mathbb{S}$-normal forms. In particular, if $\succ$ is well-founded then $\mathbb{S}$ is terminating. The converse holds if $\succ = \succ_{\mathbb{S}}$.

Note also that if $\succ$ is compatible with $\mathbb{S}$, $\equiv$-confluence implies $\succ$-tame $\equiv$-confluence, and that irrespective of whether $\succ$ is compatible, $\succ$-tame $\equiv$-confluence always implies $\succ$-tame $\equiv$-congruence.

**Lemma 6.2.8** (tamed Newmann's lemma). *Let $\mathbb{S} = (\mathbb{R}, \mathbb{E}; \equiv)$ be an abstract RSM and $\succ$ a preorder compatible with $\mathbb{S}$. If $\succ$ is terminating and every local $\mathbb{S}$-branching is $\succ$-tamely $\equiv$-congruent, then $\mathbb{S}$ is $\equiv$-confluent.*

*Proof.* Since $\succ$ is terminating, we can proceed by induction on $\succ$ (see e.g. [5, section 1.3.9]) to show that the following property holds for every $x \in X$:

> $P(x)$: every $\mathbb{S}$-branching with source $y$ such that $x \succ y$ is $\mathbb{S}$-confluent.

If $x$ is minimal for $\succ$, then in particular $x$ is an $\mathbb{S}$-normal form, and so $P(x)$ automatically holds. Consider then $x$ generic and assume that $P(y)$ holds whenever $x \succ y$.

If $h$ is an $\mathbb{S}$-congruence such that $x \succ h$, one can use the induction hypothesis to show that $h$ is $\mathbb{S}$-confluent. (Recall that $h$ being $\mathbb{S}$-confluent means that there exists a confluence $(f, g)$ with $s(f) = s(h)$, $s(g) = t(h)$, such that $f \equiv g \circ h$.) In particular, every local $\mathbb{S}$-branching with source $x$ is $\mathbb{S}$ is confluent.

Consider a (not necessarily local) $\mathbb{S}$-branching $(f, g)$. Decomposing $f$ and $g$ into $\mathbb{S}$-rewriting steps $f_m \circ \ldots \circ f_1$ and $g_n \circ \ldots \circ g_1$ respectively gives a local $\mathbb{S}$-branching $(f_1, g_1)$, $\mathbb{S}$-confluent by the previous paragraph. The rest of the confluence can be completed using induction on $\succ$, as shown in the following diagram:

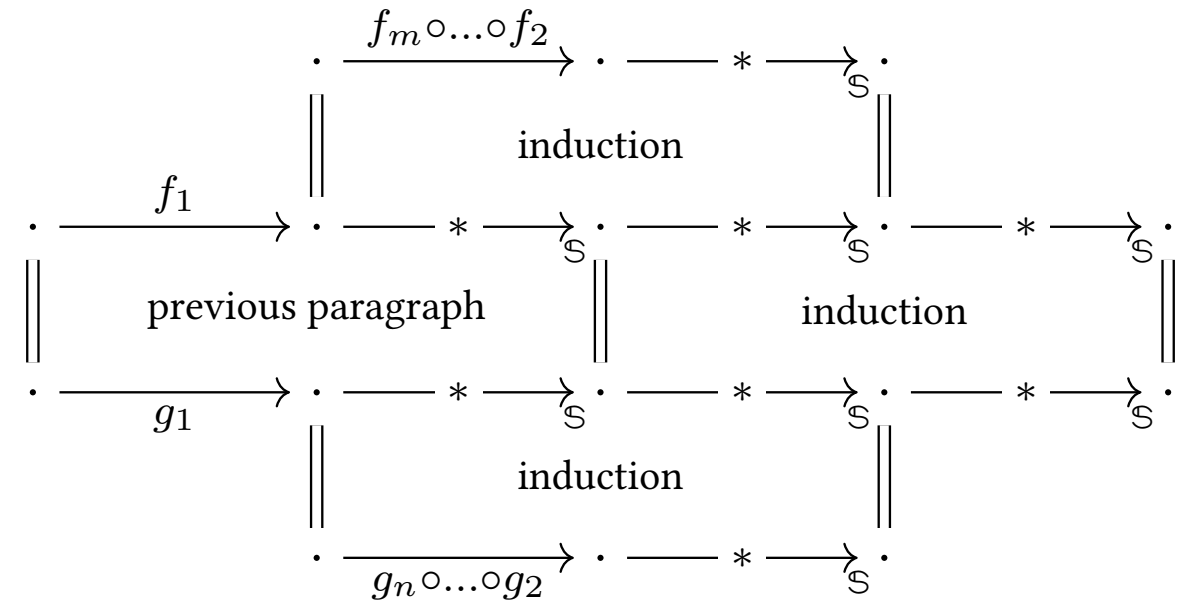

This concludes. □

### 6.2.2 Branchwise $\mathbb{E}$-congruence and Newmann's lemma

Following the (tamed) Newmann's lemma, we wish to study confluence (or tamed congruence) of local branchings. In principle, working modulo makes it a difficult task. Indeed, given that the length of the $\mathbb{E}$-congruence in a local triple is not limited, the number of local branchings is in general infinite. To circumvent this problem, rewriting steps need to be understood *up to* $\mathbb{E}$*-congruence*, similarly to how elements in $X$ are understood up to $\mathbb{E}$-congruence. In practice, one has canonical ways to do so, coming from naturality axioms with regard to the modulo: naturality of interchangers when working modulo interchange (see subsection 6.4.6), naturality of the pivotal structure when working up to isotopies (see section 7.1), and so on. This subsection formalizes this situation.

Let $\mathbb{S} = (\mathbb{R}, \mathbb{E}; \equiv)$ be an abstract RSM. Two $\mathbb{S}$-rewriting sequences $f$ and $g$ are said to be $(\mathbb{E}, \equiv)$*-congruent* if there exist $\mathbb{E}$-congruences $e_s\colon s(f) \sim_{\mathbb{E}} s(g)$ and $e_t\colon t(f) \sim_{\mathbb{E}} t(g)$ such that $e_t \circ f \equiv g \circ e_s$:

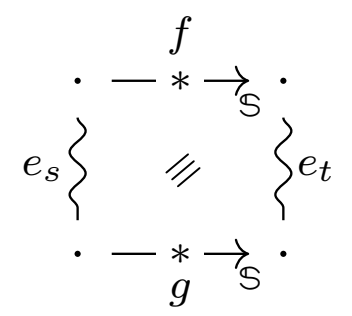

two $(\mathbb{E}, \equiv)$-congruent $\mathbb{S}$-rewriting sequences

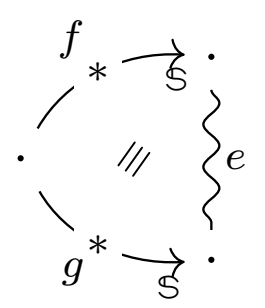

a $(\mathbb{E}, \equiv)$-congruent $\mathbb{S}$-branching

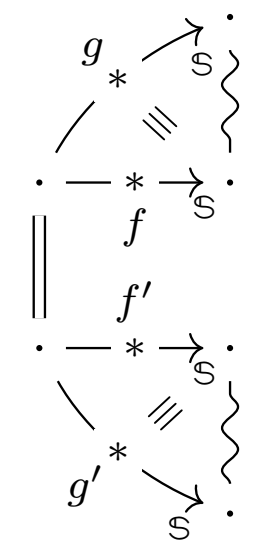

two branchwise $(\mathbb{E}, \equiv)$-congruent $\mathbb{S}$-branchings

An $\mathbb{S}$-branching $(f, g)$ is said to be $(\mathbb{E}, \equiv)$-*congruent* if there exists an $\mathbb{E}$-congruent $e\colon t(f) \sim_{\mathbb{E}} t(g)$ such that $e \circ f \equiv g$. Two $\mathbb{S}$-branchings $(f, f')$ and $(g, g')$ are said to be *branchwise* $(\mathbb{E}, \equiv)$-*congruent* if they have the same source and $(f, g)$ (resp. $(f', g')$) is $(\mathbb{E}, \equiv)$-congruent.

The following lemma states that "confluence is preserved under branchwise $\mathbb{E}$-congruence":

**Lemma 6.2.9.** *Let $\mathbb{S} = (\mathbb{R}, \mathbb{E}; \equiv)$ be an abstract RSM. If $(f, g)$ and $(f', g')$ are branchwise $(\mathbb{E}, \equiv)$-congruent $\mathbb{S}$-branchings, then $(f, g)$ is $\equiv$-confluent if and only if $(f', g')$ is.* □

The proof of the above lemma fits into one picture:

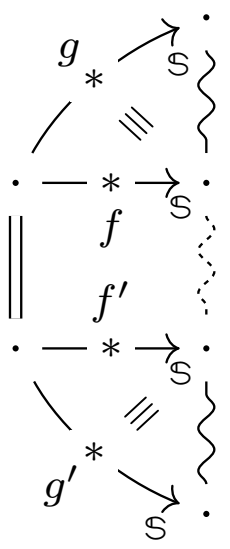

In practice, one works with $\mathbb{S}$-local triples. Two $\mathbb{S}$-local triples $[f, e, g]$ and $[f', e', g']$ are said to be $(\mathbb{E}, \equiv)$-*congruent* if $f$ (resp. $g$) is $(\mathbb{E}, \equiv)$-congruent to $f'$ (resp. $g'$) such that the relevant square of $\mathbb{E}$-congruences is $\equiv$-equivalent:

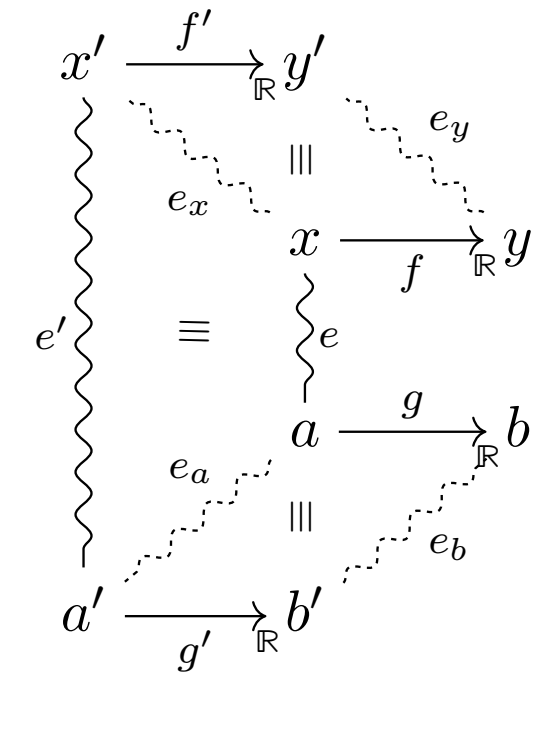

In other words, $[f, e, g]$ and $[f', e', g']$ represent branchwise $(\mathbb{E}, \equiv)$-congruent $\mathbb{S}$-branchings. For local triples, Lemma 6.2.9 says that $[f, e, g]$ is $(\mathbb{S}, \equiv)$-confluent if and only if $[f', e', g']$ is.

### 6.2.3 Branchwise confluence and tamed Newmann's lemma

Replacing $\mathbb{E}$-congruence with $\equiv$-confluence (resp. $\succ$-tamed $\equiv$-congruence) defines the following analogous branchwise notions:

**Definition 6.2.10.** *Let $\mathbb{S} = (\mathbb{R}, \mathbb{E})$ be an abstract RSM and $\succ$ a preorder on $X$. Two $\mathbb{S}$-branchings $(f, g)$ and $(f', g')$ are* branchwise $\equiv$-confluent *(resp.* branchwise $\succ$-tamely $\equiv$-congruent*) if they have the same source and the branchings $(f, f')$ and $(g, g')$ are respectively $\equiv$-confluent (resp. $\succ$-tamely $\equiv$-congruent).*

Contrary to branchwise $\mathbb{E}$-congruence, working up to branchwise confluence does not preserve confluence. However, the following lemma states that "tamed congruence is preserved under branchwise tamed congruence". Having applications in mind, we state it with respect to a sub-abstract RSM.

**Lemma 6.2.11** (Branchwise Tamed Congruence Lemma). *Let $\mathbb{S} = (\mathbb{R}, \mathbb{E})$ be an abstract RSM and $\succ$ a preorder on $X$ compatible with $\mathbb{S}$. Let also $\mathbb{T} \subset \mathbb{S}$ be a sub-abstract RSM. If $(f, g)$ and $(f', g')$ are branchwise $\succ$-tamely $(\mathbb{T}, \equiv)$-congruent $\mathbb{S}$-branchings, then $(f, g)$ is $(\mathbb{T}, \equiv)$-tamely $(\mathbb{T}, \equiv)$-congruent if and only if $(f', g')$ is.*

The proof of the above lemma fits into one picture:

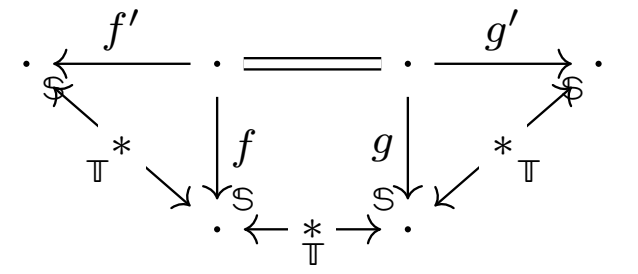

The Branchwise Tamed Congruence Lemma 6.2.11 is an important practical tool, as it can greatly simplify the study of tamed congruence. We will use it heavily in chapter 7.

### 6.2.4 Branchwise rewriting

As we can rewrite elements, we can similarly rewrite a rewriting sequence into another rewriting sequence, or a branching $(f, g)$ into another branching $(f', g')$. The latter appears as a special case of branchwise $\equiv$-confluence, with $(f, g)$ branchwise $\equiv$-confluent to a trivial branching.

**Definition 6.2.12.** *Let $\mathbb{S}$ be an abstract RSM. We say that an $\mathbb{S}$-branching $(f, g)$* rewrites *into a branching $(f', g')$ if there exist $\mathbb{S}$-rewriting sequences connecting*

*$t(f)$ with $t(f')$, $s(f) = s(g)$ with $s(f') = s(g')$ and $t(g)$ with $t(g')$, such that the relevant squares are $\equiv$-equivalent:*

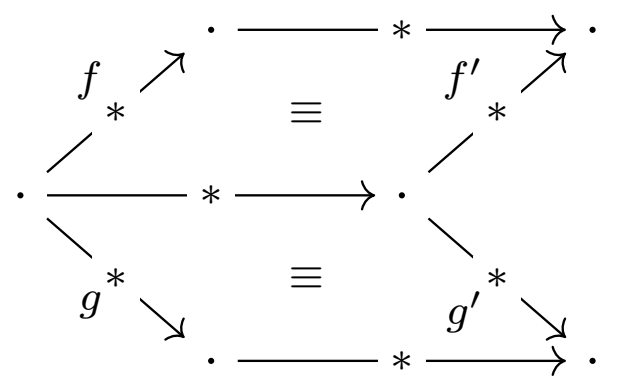

**Lemma 6.2.13.** *Let $\mathbb{S}$ be an abstract RSM and $\succ$ a preorder compatible with $\mathbb{S}$. Let also $(f, g)$ be an $\mathbb{S}$-branching that rewrites into another $\mathbb{S}$-branching $(f', g')$. If $(f', g')$ is $\succ$-tamed $\equiv$-congruent, then so is $(f, g)$.* □

Branchwise rewriting is another practical tool used in chapter 7. We shall come back to it in subsection 6.4.5 and subsection 6.5.4.

## 6.3 Linear rewriting modulo

This section generalizes linear rewriting theory to the modulo setting. In its polygraphic framework, linear rewriting theory was introduced By Guiraud, Hoffbeck and Malbos in [86], with application to associative algebras. This was extended to linear strict $n$-categories by Alleaume [3], and further extended modulo by Dupont [56].

As in Notation 5.3.1, we fix $\Bbbk$ a commutative ring. Given a set $\mathcal{B}$, we write $\langle \mathcal{B} \rangle_{\Bbbk}$ the free $\Bbbk$-module generated by $\mathcal{B}$.

### 6.3.1 Linear rewriting systems modulo

Recall that linear 1-polygraph and linear 1-sesquipolygraph are identical notions. Unpacking the definition, a linear 1-polygraph is the data of a set $\mathcal{B}$ of 0-cells, called *monomials*,[6] and a set $\mathtt{P}$ of generating 1-cells, called *relations*, equipped with *source* and *target maps* $s$ and $t$:

$$\langle \mathcal{B} \rangle_{\Bbbk} \underset{t}{\overset{s}{\leftleftarrows}} \mathtt{P}.$$

[6]This terminology is historic: in the context of rewriting in commutative algebras (i.e. polynomial algebras), elements of $\mathcal{B}$ are monomials.

Elements in $V := \langle\mathcal{B}\rangle_{\Bbbk}$ are called *vectors*. We fix a choice of monomials $\mathcal{B}$ (and hence vectors) for the reminder of the section. The $\Bbbk$-*module presented by* $(\mathcal{B}, \mathtt{P})$ is the module $\langle\mathcal{B}\rangle_{\Bbbk}/\langle\mathtt{P}\rangle_{\Bbbk}$. Note that $\langle\mathtt{P}\rangle_{\Bbbk}$ can be viewed as a globular extension of $\langle\mathcal{B}\rangle_{\Bbbk}$, extending $s$ and $t$ linearly. If $\mathtt{P}^{=}$ denotes the reflexive closure of $\mathtt{P}$, we write $\mathtt{P}^{l} := \langle\mathtt{P}^{=}\rangle_{\Bbbk}$, viewed as a globular extension:

$$\langle\mathcal{B}\rangle_{\Bbbk} \underset{t}{\overset{s}{\leftleftarrows}} \mathtt{P}^{l}.$$

Explicitly, relations in $\mathtt{P}^{l}$ are of the form $\sum_i \lambda_i r_i + v$ for $\lambda_i \in \Bbbk$, $r_i \in \mathtt{P}$ and $v \in V$, with source $\sum_i \lambda_i s(r_i) + v$ and target $\sum_i \lambda_i t(r_i) + v$. As is explained in subsection 6.3.2, $\mathtt{P}^{l}$ should be thought of the set of congruences associated to $\mathtt{P}$, analogous to $\mathbb{P}^{\top}$ in the abstract case.

We say that $\mathtt{R}$ is *left-monomial* if for all $r \in \mathtt{R}$, we have $s(r) \in \mathcal{B}$; in other words, each $r$ is of the form $b \xrightarrow{r}_{\mathtt{R}} \sum_i \lambda_i b_i$, with $\lambda_i \in \Bbbk$ and $b, b_i \in \mathcal{B}$. We say that $\mathtt{R}$ is *adapted* if it is left-monomial and we have $s(r) \notin \operatorname{supp}(t(r))$ for every $r \in \mathtt{R}$. Then:

**Definition 6.3.1.** *A* linear rewriting system (linear RS) *is the data* $(\mathcal{B}; \mathtt{P})$ *of an adapted linear 1-polygraph* $\mathtt{P}$ *on a set* $\mathcal{B}$.

We sometimes leave $\mathcal{B}$ implicit, and call $\mathtt{P}$ a linear RS. The adaptedness condition is not an important restriction. Indeed, if $\mu b + \sum_i \mu_i b_i = 0$ is some relation in a $\Bbbk$-module presentation, we can rewrite it as $b = -\mu^{-1} \sum_i \mu_i b_i$ (provided that $\mu$ is invertible). Doing so with every relation gives an adapted linear 1-polygraph presenting the given $\Bbbk$-module (provided we can always find such an invertible scalar $\mu$).

We now extend to modulos the notion of linear RS. Let $\mathtt{E}$ be another set of linear relations on $\mathcal{B}$:

$$\langle\mathcal{B}\rangle_{\Bbbk} \underset{t}{\overset{s}{\leftleftarrows}} \mathtt{E}.$$

Denote $\Bbbk^{\times}\mathcal{B}$ the subset of $\langle\mathcal{B}\rangle_{\Bbbk}$ consisting of vectors of the form $\lambda b$ for $\lambda \in \Bbbk^{\times}$ and $b \in \mathcal{B}$. We say that $\mathtt{E}$ is *monomial-invertible* if it is left-monomial and for all $e \in \mathtt{E}$, we have $t(e) \in \Bbbk^{\times}\mathcal{B}$; in other words, each $e$ is of the form $b \xrightarrow{e}_{\mathtt{E}} \lambda b'$, with $\lambda \in \Bbbk^{\times}$ and $b, b' \in \mathcal{B}$. This coincides with the notion of monomial-invertible linear 1-polygraph defined in subsection 5.3.4. If we drop the condition that scalars are invertible, we simply say that $\mathtt{E}$ is *monomial*. We will always assume $\mathtt{E}$ to be monomial-invertible.

Similarly to abstract RSM, we write $u \sim_{\mathtt{E}} v$ if there exists $e \in \mathtt{E}^{l}$ such that $s(e) = u$ and $t(e) = v$, and in that case we say that $u$ and $v$ are $\mathtt{E}$-*congruent*.

In the linear context, we furthermore have a notion of *projective* $\mathtt{E}$*-congruence*, an equivalence relation on the set of monomials $\mathcal{B}$, defined as $b \mathrel{\dot{\sim}_{\mathtt{E}}} b'$ if and only if there exists $\lambda \in \Bbbk^{\times}$ such that $b \sim_{\mathtt{E}} \lambda b'$. For $v$ a vector in $V$, we set

$$\widetilde{\dot{\mathrm{supp}}}_{\mathtt{E}}(v) := \big\{ b \in \mathcal{B} \mid b \mathrel{\dot{\sim}_{\mathtt{E}}} b' \text{ for some } b' \in \mathrm{supp}(v) \big\}$$

and call $\widetilde{\dot{\mathrm{supp}}}_{\mathtt{E}}(v)$ the $\mathtt{E}$*-projective support of* $v$.

A set $\mathtt{R}$ of linear relations on $\mathcal{B}$ is said to be $\mathtt{E}$*-adapted* if it is left-monomial and we have $s(r) \notin \widetilde{\dot{\mathrm{supp}}}_{\mathtt{E}}(t(r))$ for every $r \in \mathtt{R}$.

**Definition 6.3.2.** *A* linear rewriting system modulo (linear $\mathsf{RSM}$) *is the data* $\mathtt{S} = (\mathcal{B}; \mathtt{R}, \mathtt{E})$ *of two linear 1-polygraphs* $\mathtt{R} := (\mathcal{B}; \mathtt{R})$ *and* $\mathtt{E} := (\mathcal{B}; \mathtt{E})$*, such that* $\mathtt{E}$ *is monomial-invertible and* $\mathtt{R}$ *is* $\mathtt{E}$*-adapted.*[7]

We sometimes leave $\mathcal{B}$ implicit, and call $(\mathtt{R}, \mathtt{E})$ a linear $\mathsf{RSM}$. Note that if $\mathtt{E} = \emptyset$ (i.e. the set of relations in $\mathtt{E}$ is empty), we recover the notion of linear rewriting system. The *module presented by* $\mathtt{S}$ is the module $\langle \mathcal{B} \rangle_{\Bbbk} / \langle \mathtt{R} \sqcup \mathtt{E} \rangle_{\Bbbk}$. We write:

$$[-]_{\mathtt{S}} \colon \langle \mathcal{B} \rangle_{\Bbbk} \to \langle \mathcal{B} \rangle_{\Bbbk} / \langle \mathtt{R} \sqcup \mathtt{E} \rangle_{\Bbbk}$$

the associated quotient map. Finally, we write $\mathtt{S}^{l} := (\mathtt{R} \sqcup \mathtt{E})^{l}$.

We conclude this subsection with a few remarks; some appeal to concepts only defined in the next subsections.

*Remark* 6.3.3. The ($\mathtt{E}$-)adaptedness condition prevents "obvious" obstructions to termination. Indeed, if $s(r) \in \widetilde{\dot{\mathrm{supp}}}_{\mathtt{E}}(t(r))$ then there exists an infinite sequence of positive $\mathtt{S}$-rewriting steps, all of type $r$. Moreover, without this assumption the Church–Rosser property for positive rewritings steps (Lemma 6.3.10) does not hold; see Remark 6.3.11. Given how fundamental this result is, we choose to enforce the adaptedness condition in the definition, in contrast with the abstract case. This is only a choice of presentation; indeed, in practice one eventually wishes to work with terminating rewriting systems, which implies adaptedness. In [86], "left-monomial" encompasses both our "left-monomial" and "adapted" notions. However, the adaptedness condition is dropped in [3] (in the context of strict higher linear rewriting) and in [55, 56] (in the context of strict higher linear rewriting modulo).[8] This motivates our change of terminology, hoping to avoid further confusion. Note that while the adaptedness condition is easy to check in the non-modulo context, it can be more involved in the modulo context, as one needs to scan through projective E-congruence classes.

[7] Our definition differs from that in [55, 56], as already commented in Footnotes 4 and 5 for the abstract case; see Remark 6.3.3 for a further difference, specific to the linear case.

[8] In particular, Lemma 4.2.9 in [3] and Lemma 1.1.5 in [56] are not correct as stated. However, in both cases they eventually impose the condition of "exponentiation freedom", which implies adaptedness. (This condition is missing in [3], but this is corrected in [4].)

*Remark* 6.3.4. Recall the notion of scalar abstract RSM from Remark 6.2.3. A linear RSM $(\mathcal{B};\mathtt{R},\mathtt{E})$ is *monomial* if R is monomial. The bijection between scalar 1-polygraphs and monomial linear 1-polygraphs (see subsection 5.3.4) extends to rewriting systems:

$$\left\{\begin{array}{c}\text{scalar}\\ \text{abstract RSMs}\end{array}\right\} \underset{\text{scl}}{\overset{\text{lin}}{\rightleftarrows}} \left\{\begin{array}{c}\text{monomial}\\ \text{linear RSMs}\end{array}\right\},$$

sending $\mathbb{S} = (X;\mathbb{R},\mathbb{E},\equiv)$ to $\mathtt{S} = (\mathcal{B};\mathtt{R},\mathtt{E})$, with $\mathcal{B} = \langle X\rangle_{\Bbbk}$, $\mathtt{R} = \mathrm{lin}(\mathbb{R})$ and $\mathtt{E} = \mathrm{lin}(\mathbb{E})$. (Caveat: we did not impose the adaptedness condition in the abstract case; see also Remark 6.3.3.) It is canonical in the sense that a property holds for $\mathbb{S}$ if and only the suitable analogue holds for S; for instance, $\mathbb{S}$ is scalar-confluent if and only if $\mathtt{S}^+$ is confluent.

*Remark* 6.3.5. Recall the notion of quotient for a given abstract RSM. If E is scalar-coherent in the sense that $\mathrm{scl}(\mathtt{E})^\top$ is scalar-coherent, we can similarly define a quotient for a linear RSM $\mathtt{S} = (\mathcal{B};\mathtt{R},\mathtt{E})$, working on the $[\mathcal{B}]^{\mathtt{E}}$ of *projective* E-congruence classes. This applies for instance when E is *scalar-free*, in the sense that $s(e), t(e) \in \mathcal{B}$ for all $e \in \mathtt{E}$, which is the setting of rewriting modulo interchange and pivotality in strict 2-categories [56].

Note however that coherence of E is necessary to define the quotient, as otherwise the module $\langle\mathcal{B}\rangle_{\Bbbk}/\langle\mathtt{E}\rangle_{\Bbbk}$ is not free. The case of graded $\mathfrak{gl}_2$-foams, described chapter 7, provides an example where E is not fully coherent: interchanging two identical dots gives a scalar $XY \in \Bbbk$ (for $\Bbbk$ as defined in Definition 1.3.2).

*Remark* 6.3.6. We have chosen not to equip a linear RSM with the extra data of an abstract equivalence $\equiv$, as we did in the abstract case. In other words, a linear RSM is implicitly equipped with the discrete abstract equivalence $\equiv_{\text{disc.}}$, identifying any two parallel relations. While this is unnecessary for our purpose, there is no obstruction in generalizing linear RSMs to arbitrary equivalences on relations.

### 6.3.2 Positive rewriting steps

Compared to the abstract case, the linear case requires a specific notion of rewriting step. This is due to the fact that linearizing a set of relations $\mathtt{P}$ makes it reflexive, transitive and symmetric; that is, $\mathtt{P}^l = (\mathtt{P}^l)^\top$. Indeed:

*Proof that* $\mathtt{P}^l = (\mathtt{P}^l)^\top$. By definition, $\mathtt{P}^l := \langle \mathtt{P}^= \rangle_{\mathbb{k}}$ is reflexive. If $r_1, r_2 \in \mathtt{P}^l$ are relations such that $t(r_1) = s(r_2)$, then $r_1 + r_2 - t(r_1) = r_1 + r_2 - s(r_2)$ is a relation in $\mathtt{P}^l$ whose source is $s(r_1)$ and target is $t(r_2)$; hence, $\mathtt{P}^l$ is transitive. If $\alpha \in \mathtt{P}^l$ then $\alpha^{-1} = s(\alpha) + t(\alpha) - \alpha \in \langle \mathtt{P} \rangle_{\mathbb{k}}$; hence, $\mathtt{P}^l$ is symmetric. □

In particular, $\mathtt{P}^l$ cannot terminate. To tackle this issue, we work with a subset $\mathtt{P}^+ \subset \langle \mathtt{P} \rangle_{\mathbb{k}}$, the subset of *positive* $\mathtt{P}$*-rewriting steps*, that avoids this formal obstruction to termination while describing the same congruence, i.e. $(\mathtt{P}^+)^\top = \mathtt{P}^l$.

**Definition 6.3.7.** *Let* $\mathtt{P}$ *be a linear RS. A* $\mathtt{P}$-rewriting step *is an element* $\alpha \in \langle \mathtt{P} \rangle_{\mathbb{k}}$ *of the form*

$$\alpha = \lambda r + v, \quad \lambda \in \mathbb{k} \setminus \{0\}, r \in \mathtt{P}, v \in \langle \mathcal{B} \rangle_{\mathbb{k}}.$$

*In that case, we say that* $\alpha$ *is* of type $r$. *The* $\mathtt{P}$*-rewriting step* $\alpha$ *is said to be* positive*;*[9] *if* $s(r) \notin \operatorname{supp}(v)$. *We write a* $\mathtt{P}$*-rewriting step as* $\alpha\colon s(\alpha) \dashrightarrow_{\mathtt{P}} t(\alpha)$*, and a positive* $\mathtt{P}$*-rewriting step as* $\alpha\colon s(\alpha) \to_{\mathtt{P}} t(\alpha)$*. The set of (resp. positive)* $\mathtt{P}$*-rewriting steps is denoted* $\mathtt{P}^{\mathrm{st}}$ *(resp.* $\mathtt{P}^+$*).*

*If* $v = 0$*, we say that* $\alpha$ *is* monomial*; in that case,* $\alpha$ *is necessarily positive.*

Note that while $\mathtt{P}^{\mathrm{st}}$ is a priori not transitive, it is symmetric; indeed, if $\alpha$ is a $\mathtt{P}$-rewriting step then so is its inverse

$$\alpha^{-1} = (-\lambda) r + \lambda(s(r) + t(r)) + v,$$

as defined in the proof above. Hence, $\mathtt{P}^{\mathrm{st}}$ does not terminate and cannot provide a suitable reduction algorithm. However, if $\alpha$ is positive then $\alpha^{-1}$ *cannot* be positive; indeed, the assumptions $s(r) \notin \operatorname{supp}(t(r))$ (adaptedness of $\mathtt{P}$) and $s(r) \notin \operatorname{supp}(v)$ (positiveness of $\alpha$) imply that

$$s(r) \in \operatorname{supp}(\lambda(s(r) + t(r)) + v).$$

This makes positive $\mathtt{P}$-rewriting steps suitable candidates to define a terminating reduction algorithm.

[9] This terminology appears e.g. in [35] other references [3, 56, 86] use the terminology *elementary relation* for rewriting step, and rewriting step for positive rewriting step.

These notions are readily extended to modulos:

**Definition 6.3.8.** *Let* $\mathtt{S} = (\mathtt{R}, \mathtt{E})$ *be a linear RSM. An* $\mathtt{S}$-rewriting step *is a composition*

$$u \;\sim_{\mathtt{E}}\; \lambda s(r) + v \;\dashrightarrow_{\mathtt{R}}\; \lambda t(r) + v \;\sim_{\mathtt{E}}\; w,$$

*where the middle arrow is a* $\mathtt{R}$*-rewriting step. This* $\mathtt{S}$*-rewriting step is said to be* positive *if* $s(r) \notin \widetilde{\dot{\mathrm{supp}}}_{\mathtt{E}}(v)$*. We write an* $\mathtt{S}$*-rewriting step as* $\dashrightarrow_{\mathtt{S}}$*, and a positive* $\mathtt{P}$*-rewriting step as* $\rightarrow_{\mathtt{S}}$*:*

| $s(\alpha) \dashrightarrow_{\mathtt{S}} t(\alpha)$ | $s(\alpha) \rightarrow_{\mathtt{S}} t(\alpha)$ |
|---|---|
| *(not necessarily positive) rewriting step* | *positive rewriting step* |

*The set of (resp. positive)* $\mathtt{S}$*-rewriting steps is denoted* $\mathtt{S}^{\mathrm{st}}$ *(resp.* $\mathtt{S}^{+}$*).*

*If* $v = 0$*, we say that* $\alpha$ *is* monomial*; in that case,* $\alpha$ *is necessarily positive.*

As in the non-modulo context, one can argue that $\mathtt{S}^{l} = (\mathtt{S}^{l})^{\top}$, with positive $\mathtt{S}$-rewriting steps providing an answer to this formal obstruction to termination.

Note that the positivity condition is $s(r) \notin \widetilde{\dot{\mathrm{supp}}}_{\mathtt{E}}(v)$, and not $s(r) \notin \mathrm{supp}(v)$. In other words, a positive $\mathtt{S}$-rewriting step is *not* a composition as above such that the middle arrow is a positive $\mathtt{R}$-rewriting step. The positivity condition is stronger, and depends on $\mathtt{E}$. Otherwise, positive $\mathtt{S}$-rewriting steps do not provide a suitable solution to termination.[10] Note that for the positivity condition to make sense, $\mathtt{E}$ must be monomial. Extending the theory to non-monomial modulo rules remains a non-trivial question.[11]

Both $\mathtt{S}^{\mathrm{st}}$ and $\mathtt{S}^{+}$ provide an abstract RSM associated to $\mathtt{S}$, namely respectively $(\langle\mathcal{B}\rangle_{\Bbbk}; \mathtt{S}^{\mathrm{st}}, \mathtt{E}^{\mathrm{st}})$ and $(\langle\mathcal{B}\rangle_{\Bbbk}; \mathtt{S}^{+}, \mathtt{E}^{\mathrm{st}})$. As such, an $\mathtt{S}$-rewriting step (resp. a positive $\mathtt{S}$-rewriting step) in the sense of Definition 6.3.8 is the same as an $\mathtt{S}^{\mathrm{st}}$-rewriting step (resp. an $\mathtt{S}^{+}$-rewriting step), and we shall use the two terminologies interchangeably. Similarly, a *positive* $\mathtt{S}$*-rewriting sequence* denotes an $\mathtt{S}^{+}$-rewriting sequence.

The following lemma shows that both $\mathtt{S}$-rewriting steps and positive $\mathtt{S}$-rewriting steps are suited to study $\mathtt{S}$-congruence:

[10] For instance, consider the linear RSM $\mathtt{S} = (\mathcal{B}; \mathtt{R}, \mathtt{E})$ with $\mathcal{B} = \{a, a', b\}$, $\mathtt{R} = \{a \to b\}$ and $\mathtt{E} = \{a \sim a'\}$, and the $\mathtt{S}$-rewriting sequence $0 = a - a \sim_{\mathtt{E}} a - a' \dashrightarrow_{\mathtt{R}} b - a' \sim_{\mathtt{E}} b - a \dashrightarrow_{\mathtt{R}} b - b = 0$.

[11] In [55, 56], a positive $\mathtt{S}$-rewriting step *is* defined as a composition as above such that the middle arrow is a positive $\mathtt{R}$-rewriting step, and $\mathtt{E}$ is not constrained to be monomial. In [57], Dupont's work is applied with a non-monomial $\mathtt{E}$.

**Lemma 6.3.9.** *Let* $\mathtt{S} = (\mathcal{B}; \mathtt{R}, \mathtt{E})$ *be a linear RSM,* $u, v \in \langle \mathcal{B} \rangle_{\Bbbk}$ *and recall the notation* $[u]_{\mathtt{S}}, [v]_{\mathtt{S}}$ *for their respective image in* $\langle \mathcal{B} \rangle_{\Bbbk} / \langle \mathtt{R} \sqcup \mathtt{E} \rangle_{\Bbbk}$*. The following statements are equivalent:*

*(i)* $[u]_{\mathtt{S}} = [v]_{\mathtt{S}}$ *in* $\langle \mathcal{B} \rangle_{\Bbbk} / \langle \mathtt{R} \sqcup \mathtt{E} \rangle_{\Bbbk}$*;*

*(ii)* $u$ *and* $v$ *are* $\mathtt{S}^{\mathrm{st}}$*-congruent;*

*(iii)* $u$ *and* $v$ *are* $\mathtt{S}^{+}$*-congruent.*

*In particular,* $\mathtt{S}^{l} = (\mathtt{S}^{+})^{\top} = (\mathtt{S}^{\mathrm{st}})^{\top}$.

The lemma implies that there is no distinction between the properties of $\mathtt{S}^{+}$-congruence, an $\mathtt{S}^{\mathrm{st}}$-congruence or an $\mathtt{S}^{l}$-congruence; note however that a *given* $\mathtt{S}^{\mathrm{st}}$-congruence may not be positive. In order to prove the lemma, we need the following Church–Rosser property for positive rewriting steps, which generalizes modulo Lemma 3.1.2 in [86]. The last statement is explained and used in subsection 6.3.4 and can be ignored for the purpose of this subsection.

**Lemma 6.3.10.** *Let* $\mathtt{S} = (\mathcal{B}; \mathtt{R}, \mathtt{E})$ *be a linear RSM. If* $f$ *is an* $\mathtt{S}$*-rewriting step, then there exist positive* $\mathtt{S}$*-rewriting steps* $g, h$ *of length at most one such that* $f = h^{-1} \circ g$*:*

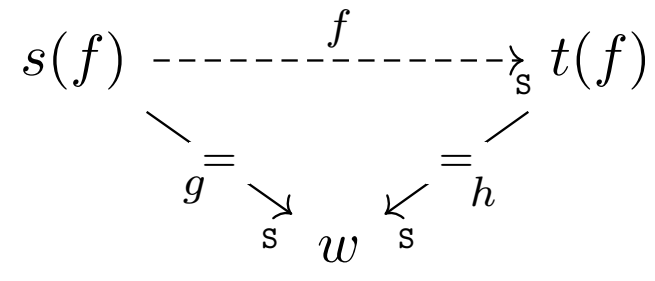

*Moreover, we have* $f \succcurlyeq^{\mathrm{rel}} w$ *for any* $\mathtt{E}$*-invariant linear preorder* $\succ$ *on* $\mathcal{B}$ *(see Definition 6.3.23).*

*Proof.* Let $\lambda r + v$ be the $\mathtt{R}$-rewriting step associated to $f$, with $r \in \mathtt{R}$, $\lambda \in \Bbbk \backslash \{0\}$ and $v \in \langle \mathcal{B} \rangle_{\Bbbk}$. Extracting $s(r)$ from the decomposition of $v$, we let $\mu \in \Bbbk$ and $v' \in \langle \mathcal{B} \rangle_{\Bbbk}$ such that $v \sim_{\mathtt{E}} \mu s(r) + v'$ and $s(r) \notin \widetilde{\mathrm{supp}}_{\mathtt{E}}(v')$. Set

$$w := (\lambda + \mu) t(r) + v'.$$

Then $s(f) \sim_{\mathtt{E}} (\lambda + \mu) s(r) + v'$ (resp. $t(f) \sim_{\mathtt{E}} \mu s(r) + (\lambda t(r) + v')$) is either equal to $w$ if $\mu = -\lambda$ (resp. if $\mu = 0$), or there exists a positive rewriting step $(\lambda + \mu) s(r) + v' \to_{\mathtt{R}} w$ (resp. $\mu s(r) + (\lambda t(r) + v') \to_{\mathtt{R}} w$):

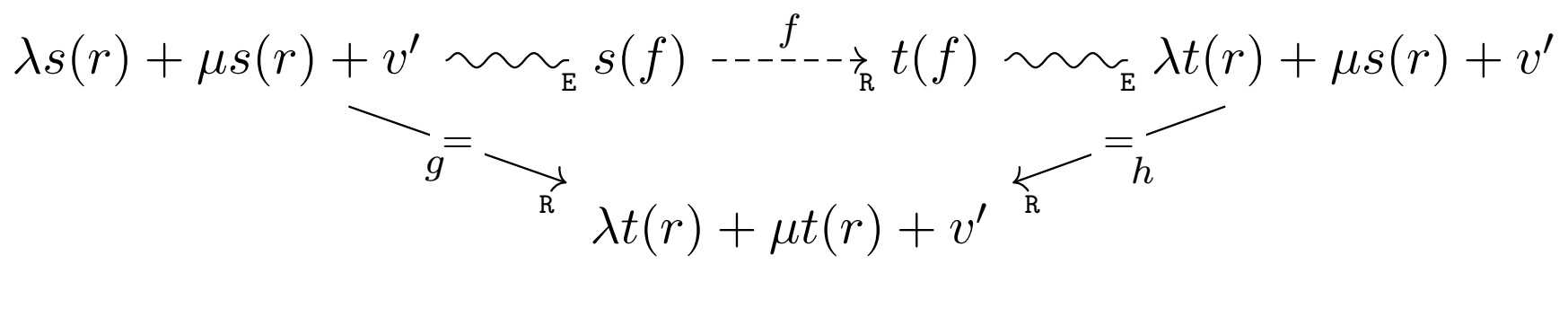

Here we use the fact that $s(r) \notin \widetilde{\dot{\mathrm{supp}}}_{\mathtt{E}}(t(r))$ (E-adaptedness) to ensure that $s(r) \notin \widetilde{\dot{\mathrm{supp}}}_{\mathtt{E}}(\lambda t(r) + v')$. Finally, it follows from

$$\widetilde{\dot{\mathrm{supp}}}_{\mathtt{E}}\left(\lambda s(r) + \mu s(r) + v'\right) \cup \widetilde{\dot{\mathrm{supp}}}_{\mathtt{E}}\left(\lambda t(r) + \mu s(r) + v'\right) \supset \widetilde{\dot{\mathrm{supp}}}_{\mathtt{E}}\left((\lambda + \mu)t(r) + v'\right)$$

that $f \not\succcurlyeq^{\mathrm{rel}} w$ for any E-invariant linear preorder $\succ$ on $\mathcal{B}$. □

*Proof of Lemma 6.3.9.* (iii) $\Leftrightarrow$ (ii) is given by Lemma 6.3.10. If $u \dashrightarrow_{\mathtt{S}} v$, then $[u] = [v]$, so (ii) $\Rightarrow$ (i). To show (i) $\Rightarrow$ (ii), assume that $[u] = [v]$. In that case, $u \sim_{\mathtt{E}} v + \sum_{i=0}^{n} \lambda_i (s(r_i) - t(r_i))$ for some scalars $\lambda_i \in \Bbbk$ and relations $r_i \in \mathtt{R}$. Write $u_j = \sum_{i \leq j} \lambda_i (t(r_i) - s(r_i))$ and $\alpha_j \colon u_j \to u_{j-1}$ the obvious S-rewriting step. Successively applying the $\alpha_j$'s defines an S-rewriting sequence

$$u \sim v + u_n \overset{\alpha_n}{\dashrightarrow}_R v + u_{n-1} \overset{\alpha_n}{\dashrightarrow}_R \ldots \overset{\alpha_n}{\dashrightarrow}_R v.$$

This concludes. □

*Remark* 6.3.11. Without the adaptedness condition, Lemma 6.3.10 does not hold: for instance, one can consider the linear RS $\mathtt{P} = \{a \to 2a\}$ and the non-positive rewriting step $a + a \to 2a + a$ as a counterexample.

### 6.3.3 Basis from convergence modulo

Fix $\mathtt{S} = (\mathcal{B}; \mathtt{R}, \mathtt{E})$ a linear RSM. In this subsection, we explain how convergence modulo of $\mathtt{S}^+$ can provide a basis for the underlying module. Unsurprisingly, finding a basis is closely related to understanding congruence. The following is a direct consequence of Lemma 6.3.9:

**Corollary 6.3.12.** *Let* $\mathtt{S} = (\mathcal{B}; \mathtt{R}, \mathtt{E})$ *be a linear RSM and* $B \subset \mathcal{B}$ *a subset of* $\mathcal{B}$. *Then* $[B]_{\mathtt{S}}$ *is a basis of* $\langle \mathcal{B} \rangle_{\Bbbk} / \langle \mathtt{R} \sqcup \mathtt{E} \rangle_{\Bbbk}$ *if and only if* $\langle B \rangle_{\Bbbk}$ *is a set of unique* $\mathtt{S}^+$*-congruence representatives.* □

*Remark* 6.3.13. If $(\mathcal{B}; \mathtt{E})$ is a monomial linear RS, then the above corollary implies that $[\mathcal{B}]_{\mathtt{E}}$ is free if and only if $\mathrm{scl}(\mathtt{E})^{\top}$ is scalar-coherent (see Remark 6.3.4). In that case, $[B]_{\mathtt{E}}$ is a basis if and only if $B$ is a unique set of monomial E-congruence representatives.

As in the abstract setting, normal forms are prime candidates for congruence representatives. To get a candidate basis, we look at monomial normal forms:

**Definition 6.3.14.** *Let* $\mathtt{S} = (\mathcal{B}; \mathtt{R}, \mathtt{E})$ *be a linear RSM. A* monomial $\mathtt{S}^+$-normal form *is a monomial* $b \in \mathcal{B}$ *which is a normal form for* $\mathtt{S}^+$*; that is, we have* $b \not\to_{\mathtt{E}} s(r)$ *for all* $r \in \mathtt{R}$*. We denote* $\mathcal{B}\mathrm{NF}_{\mathtt{S}}$ *the set of monomial* $\mathtt{S}^+$*-normal forms.*

*Remark* 6.3.15. A linear combination of $\mathtt{S}^+$-normal forms is an $\mathtt{S}^+$-normal form, and monomials in the support of an $\mathtt{S}^+$-normal form are $\mathtt{S}^+$-normal forms. The zero vector 0 is always an $\mathtt{S}^+$-normal form. In particular, $\mathrm{NF}_{\mathtt{S}} = \langle \mathcal{B}\mathrm{NF}_{\mathtt{S}} \rangle_{\Bbbk}$.

Given that convergence implies existence and uniqueness of normal form representatives, we get the following proposition:

**Proposition 6.3.16.** *Let* $\mathtt{S} = (\mathcal{B}; \mathtt{R}, \mathtt{E})$ *be a linear RSM and* $B \subset \mathcal{B}\mathrm{NF}_{\mathtt{S}}$ *such that* $[B]_{\mathtt{E}}$ *is a basis for the module* $\mathrm{NF}_{\mathtt{S}}/\langle E \rangle_{\Bbbk}$*. If* $\mathtt{S}^+$ *is convergent modulo, then* $[B]_{\mathtt{S}}$ *is a basis for the module* $\langle \mathcal{B} \rangle_{\Bbbk}/\langle \mathtt{R} \sqcup \mathtt{E} \rangle_{\Bbbk}$ *presented by* $\mathtt{S}$.

*Proof.* By Corollary 6.3.12, we must check that $\langle B \rangle_{\Bbbk}$ is a set of unique $\mathtt{S}^+$-congruence representatives. By Proposition 6.1.6, this is equivalent to $\langle B \rangle_{\Bbbk}$ being a set of unique $\mathtt{E}$-congruence representatives on the set $\mathrm{NF}_{\mathtt{S}}$. Finally, by Corollary 6.3.12 again this is equivalent to $[B]_{\mathtt{E}}$ being a basis for the module $\mathrm{NF}_{\mathtt{S}}/\langle E \rangle_{\Bbbk}$. □

Setting $\mathtt{E} = \emptyset$, Proposition 6.3.16 becomes:

**Corollary 6.3.17.** *Let* $(\mathcal{B}; \mathtt{P})$ *be a linear RS. If* $\mathtt{P}^+$ *is convergent, then* $\mathcal{B}\mathrm{NF}_{\mathtt{P}}$ *is a basis for the module* $\langle \mathcal{B} \rangle_{\Bbbk}/\langle \mathtt{P} \rangle_{\Bbbk}$ *presented by* $\mathtt{P}$. □

Combining Proposition 6.3.16 and Corollary 6.3.17, we arrive at the following theorem:

**Theorem 6.3.18** (Basis From Convergence Theorem)**.** *Let* $\mathtt{S} = (\mathcal{B}; \mathtt{R}, \mathtt{E})$ *be a linear RSM. If* $\mathtt{S}^+$ *is convergent and if* $\mathrm{scl}(\mathtt{E})^\top$ *is scalar-coherent on the set* $\mathcal{B}\mathrm{NF}_{\mathtt{S}}$ *of monomial* $\mathtt{S}^+$*-normal forms, then the module* $\langle \mathcal{B} \rangle_{\Bbbk}/\langle \mathtt{R} \sqcup \mathtt{E} \rangle_{\Bbbk}$ *presented by* $\mathtt{S}$ *is free, and any choice of projective* $\mathtt{E}$*-congruence representatives on* $\mathcal{B}\mathrm{NF}_{\mathtt{S}}$ *defines a basis.* □

The Basis From Convergence Theorem 6.3.18 is used in chapter 7 to provide a basis for graded $\mathfrak{gl}_2$-foams.

### 6.3.4 Termination order

This section extends to the linear case the notion of abstract compatible preorder (subsection 6.2.1).

Given a linear RSM $\mathtt{S} = (\mathcal{B}; \mathtt{R}, \mathtt{E})$, we say that a relation $\succ$ on $\mathcal{B}$ is $\mathtt{E}$-*invariant* if it is invariant with respect to projective $\mathtt{E}$-congruence in the sense of subsection 6.2.1; that is, if

$$(a' \dot{\sim}_{\mathtt{E}} a \text{ and } a \succ b \text{ and } b \dot{\sim}_{\mathtt{E}} b') \text{ implies } (a' \succ b').$$

In the presence of a linear RSM, We shall always implicitly assume that relations are $\mathtt{E}$-invariant.

**Definition 6.3.19.** *Let* $\mathtt{S} = (\mathcal{B}; \mathtt{R}, \mathtt{E})$ *be a linear RSM and* $\succ$ *a preorder on* $\mathcal{B}$. *We say that* $\succ$ *is* compatible with $\mathtt{S}$ *if*

$$s(r) \succ b \qquad \textit{for all } r \in \mathtt{R} \textit{ and } b \in \operatorname{supp}(t(r)).$$

*We denote* $\succ_{\mathtt{S}}$ *the smallest* $\mathtt{E}$*-invariant preorder on* $\mathcal{B}$ *compatible with* $\mathtt{S}$.

Recall from subsection 5.3.3 the notion of a linear relation. An ($\mathtt{E}$-invariant) preorder $\succ$ on $\mathcal{B}$ induces an ($\mathtt{E}^l$-invariant) linear relation $\succ^+$ on $\langle\mathcal{B}\rangle_{\mathbb{k}}$, setting $u \succ^+ v$ whenever the following two conditions hold:

(a) $\overset{\cdot}{\widetilde{\mathrm{supp}}}_{\mathtt{E}}(u) \neq \overset{\cdot}{\widetilde{\mathrm{supp}}}_{\mathtt{E}}(v)$;

(b) for every $a$ in $\overset{\cdot}{\widetilde{\mathrm{supp}}}_{\mathtt{E}}(v) \setminus \overset{\cdot}{\widetilde{\mathrm{supp}}}_{\mathtt{E}}(u)$, there exists $b \in \overset{\cdot}{\widetilde{\mathrm{supp}}}_{\mathtt{E}}(u) \setminus \overset{\cdot}{\widetilde{\mathrm{supp}}}_{\mathtt{E}}(v)$ such that $b \succ a$.

This definition corresponds to setting $u \succ^+ v$ if and only if $\overset{\cdot}{\widetilde{\mathrm{supp}}}_{\mathtt{E}}(u) \succ^{\mathrm{set}} \overset{\cdot}{\widetilde{\mathrm{supp}}}_{\mathtt{E}}(v)$ where $\succ^{\mathrm{set}}$ is the multi-set relation induced by $\succ$,[12] defined on the power set $\mathcal{P}(\mathcal{B})$ as:

$$M \succ^{\mathrm{set}} N \quad \Leftrightarrow \quad M \neq N \text{ and } \forall y \in N \setminus M, \exists x \in M \setminus N \text{ such that } x \succ y.$$

Equivalently, $M \succ^{\mathrm{set}} N$ if and only if one can go from $M$ to $N$ by a sequence of moves consisting in removing an element $b$ and adding elements $a_i$ with $b \succ a_i$. This last interpretation motivates the definition of $\succ^+$, designed precisely such that the following holds:

**Lemma 6.3.20.** *Let* $\mathtt{S} = (\mathcal{B}; \mathtt{R}, \mathtt{E})$ *be a linear RSM and* $\succ$ *a relation on* $\mathcal{B}$. *If* $\succ$ *is a strict*[13] *preorder on* $\mathcal{B}$ *compatible with* $\mathtt{S}$ *in the linear sense of Definition 6.3.19, then* $\succ^+$ *is a strict preorder on* $\langle\mathcal{B}\rangle_{\mathbb{k}}$ *compatible with* $\mathtt{S}^+$ *in the abstract sense of Definition 6.2.7.*

[12]See e.g. [7, Definition 2.5.3]; we only use the special case of sets, and therefore denote it $\succ^{\mathrm{set}}$.

[13]A preorder $\succ$ is *strict* if $b \not\succ b$ for all $b \in \mathcal{B}$.

*Proof.* It is shown in [7, lemma 2.5.4] that if $\succ$ is a strict preorder, so is $\succ^{\text{set}}$; hence if $\succ$ is a strict preorder, so is $\succ^{+}$. Note that without the strictness condition, $\succ^{+}$ may not even be transitive.

Consider then $r\colon s(r) \to_{\mathtt{R}} t(r)$. Strictness (or E-adaptedness) implies that

$$\widetilde{\dot{\operatorname{supp}}}_{\mathtt{E}}(s(r)) \setminus \widetilde{\dot{\operatorname{supp}}}_{\mathtt{E}}(t(r)) = \widetilde{\dot{\operatorname{supp}}}_{\mathtt{E}}(s(r)),$$

so that $s(r) \succ^{+} t(r)$. The general case follows from the following lemma. □

**Lemma 6.3.21.** *Let* $\mathtt{S} = (\mathcal{B}; \mathtt{R}, \mathtt{E})$ *be a linear RSM and* $\succ$ *a relation on* $\mathcal{B}$*. For vectors* $u, v \in \langle\mathcal{B}\rangle_{\Bbbk}$*, we have:*

$$u \succ^{+} v \quad \Rightarrow \quad \lambda u + w \succ^{+} \lambda v + w$$

*for all* $\lambda \in \Bbbk \setminus \{0\}$ *and* $w \in \langle\mathcal{B}\rangle_{\Bbbk}$ *such that* $\widetilde{\dot{\operatorname{supp}}}_{\mathtt{E}}(u) \cap \widetilde{\dot{\operatorname{supp}}}_{\mathtt{E}}(w) = \emptyset$.

*Proof.* Since $\widetilde{\dot{\operatorname{supp}}}_{\mathtt{E}}(\lambda u + w) = \widetilde{\dot{\operatorname{supp}}}_{\mathtt{E}}(u) \sqcup \widetilde{\dot{\operatorname{supp}}}_{\mathtt{E}}(w)$, we have:

$$\widetilde{\dot{\operatorname{supp}}}_{\mathtt{E}}(\lambda u + w) \setminus \widetilde{\dot{\operatorname{supp}}}_{\mathtt{E}}(\lambda v + w) \supset \widetilde{\dot{\operatorname{supp}}}_{\mathtt{E}}(u) \setminus \widetilde{\dot{\operatorname{supp}}}_{\mathtt{E}}(v),$$
$$\widetilde{\dot{\operatorname{supp}}}_{\mathtt{E}}(\lambda v + w) \setminus \widetilde{\dot{\operatorname{supp}}}_{\mathtt{E}}(\lambda u + w) \subset \widetilde{\dot{\operatorname{supp}}}_{\mathtt{E}}(v) \setminus \widetilde{\dot{\operatorname{supp}}}_{\mathtt{E}}(u).$$ □

We now relate to termination:[14]

**Lemma 6.3.22.** *Let* $\mathtt{S} = (\mathcal{B}; \mathtt{R}, \mathtt{E})$ *be a linear RSM and* $\succ$ *a preorder on* $\mathcal{B}$*. If* $\succ$ *is compatible with* $\mathtt{S}$*, we have the following implications:*

$$\succ \textit{ is terminating on } \mathcal{B} \quad \Leftrightarrow \quad \succ^{+} \textit{ is terminating on } \langle\mathcal{B}\rangle_{\Bbbk} \quad \Rightarrow \quad \mathtt{S}^{+} \textit{ terminates.}$$

*Moreover, if* $\succ = \succ_{\mathtt{S}}$ *(see Definition 6.3.19) then the converse of the last implication holds.*

*Proof.* Note that if a preorder terminates, it is necessarily strict. It is shown in [7, theorem 2.5.5] that $\succ$ is terminating if and only if $\succ^{\text{set}}$ is terminating. In that case, Lemma 6.3.20 implies that $\mathtt{S}^{+}$ terminates.

It remains to show that if $\mathtt{S}^{+}$ terminates, then $\succ_{\mathtt{S}}$ is terminating on $\mathcal{B}$. We proceed by contraposition and assume that there exists an infinite sequence

$$b_0 \succ_{\mathtt{S}} b_1 \succ_{\mathtt{S}} b_2 \succ_{\mathtt{S}} \dots b_n \succ_{\mathtt{S}} b_{n+1} \succ_{\mathtt{S}} \cdots$$

[14] We remind the reader that $\mathtt{S}^{+}$ terminates whenever there is no infinite sequence of consecutive rewriting steps.

in $\mathcal{B}$. Let $r_n \in \mathtt{R}$ such that $s(r_n) \dot{\sim}_{\mathtt{E}} b_n$ and $b_{n+1} \in \widetilde{\mathrm{supp}}_{\mathtt{E}}(t(r_n))$. We construct a sequence of $\mathtt{S}^+$-rewriting steps starting with $v_0 = b_0$ and defining $v_{n+1}$ recursively by applying $r_m$ on $v_n$ (possibly after an E-congruence), for $m$ the biggest index possible. The assumption on $m$ ensures that we always have $b_{m+1} \in \mathrm{supp}(v_{n+1})$, so that this process does not end. Hence, $\mathtt{S}^+$ does not terminate. □

We conclude the subsection with the following notion, which appears in the statement of Lemma 6.3.10:

**Definition 6.3.23.** *Let $\succ$ be a relation on a set $\mathcal{B}$. the* relative relation $\succcurlyeq^{\mathrm{rel}}$ induced by $\succ$ *is the following relation on the power set $\mathcal{P}(\mathcal{B})$:*

$$M \succcurlyeq^{\mathrm{rel}} N \quad \Leftrightarrow \quad \forall b \in \mathcal{B}, (b \succ M) \Rightarrow (b \succ N).$$

*If $\mathtt{S} = (\mathcal{B}; \mathtt{R}, \mathtt{E})$ is a linear RSM and $\succ$ an E-invariant preorder on $\mathcal{B}$, we set $u \succcurlyeq^{\mathrm{rel}} v$ if and only if $\widetilde{\mathrm{supp}}_{\mathtt{E}}(u) \succcurlyeq^{\mathrm{rel}} \widetilde{\mathrm{supp}}_{\mathtt{E}}(v)$.*

Note that equivalently, $M \succcurlyeq^{\mathrm{rel}} N$ if and only if $\forall L \in \mathcal{P}(\mathcal{B})$, $(L \succ M)$ implies $(L \succ N)$.

## 6.3.5 Tamed linear Newmann's lemma

The linear structure of a linear RSM $\mathtt{S} = (\mathcal{B}; \mathtt{R}, \mathtt{E})$ induces canonical types of $\mathtt{S}^+$-local triples. Let $[f, e, g]$ be a $\mathtt{S}^+$-local triple such that $f$ and $g$ is of type $r_1$ and $r_2$, respectively. Then:

- if $s(r_1) \not\dot{\sim}_{\mathtt{E}} s(r_2)$, one says that $[f, e, g]$ is *additive*,
- if $s(r_1) \dot{\sim}_{\mathtt{E}} s(r_2)$, one says that $[f, e, g]$ is *intersecting*.

Recall the notion of branchwise E-congruence for branchings (subsection 6.2.2). A $\mathtt{S}^+$-local triple $[f, e, g]$ is additive if, up to branchwise E-congruence, it has of the following form:

$$[f, e, g] = [\lambda_1 r_1 + \lambda_2 s(r_2) + v, \mathrm{id}, \lambda_1 s(r_1) + \lambda_2 r_2 + v],$$

An additive branching has a canonical $\mathtt{S}^{\mathrm{st}}$-confluence:

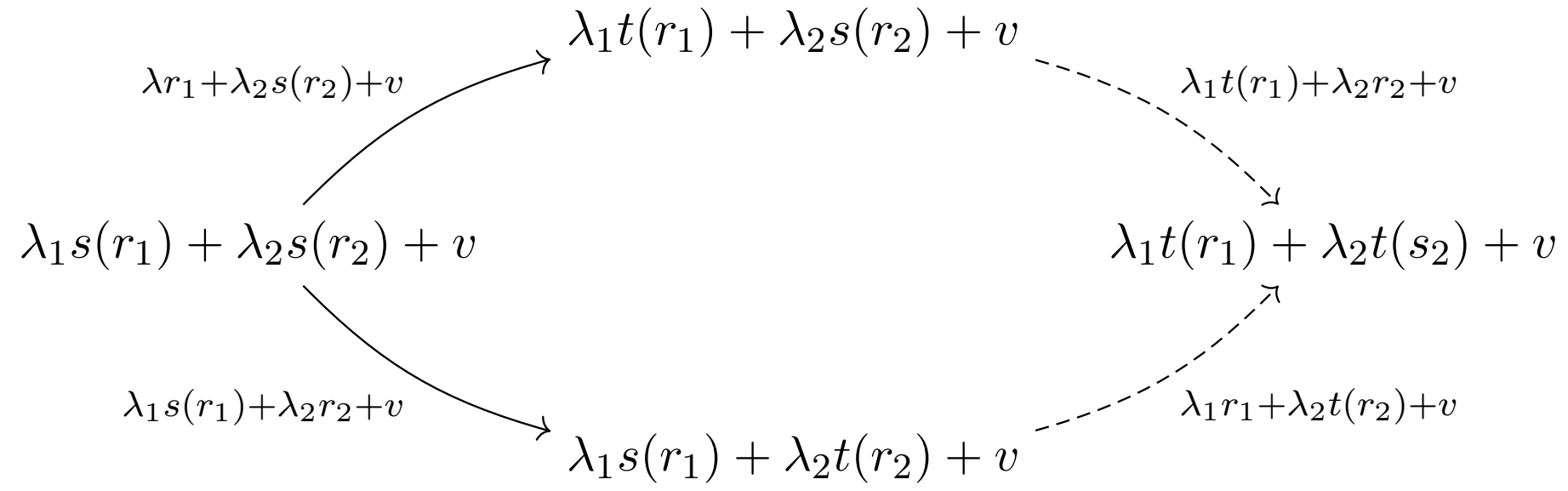

On the other hand, $[f, e, g]$ is intersecting if, up to branchwise $\mathtt{E}$-congruence, it has the following form:

$$[f, e, g] = [\lambda_1 r_1 + w, e_r + w, \lambda_2 r_2 + w],$$

where $e_r \colon \lambda_1 s(r_1) \sim_{\mathtt{E}} \lambda_2 s(r_2)$. We say that $[f, e, g]$ is *monomial* if it is intersecting and if, up to branchwise $\mathtt{E}$-congruence, it has the form $\lambda_1 = 1$ and $w = 0$ with the notations above. That is, $[f, e, g]$ is monomial if, up to branchwise $\mathtt{E}$-congruence, it has the following form:

$$[f, e, g] = [r_1, e_r, \lambda r_2].$$

where $e_r \colon s(r_1) \sim_{\mathtt{E}} \lambda s(r_2)$.

Note that both branches of a monomial branching are monomial rewriting steps. By definition, every intersecting branching is of the form $\lambda_1[f, e, g] + w$ for $[f, e, g]$ a monomial branching. Still by definition, additive, intersecting and monomial branchings are positive branchings.

We end the subsection with a linear analogue of the tamed Newmann's lemma (Lemma 6.2.5). Given a linear RSM $\mathtt{S} = (\mathcal{B}; \mathtt{R}, \mathtt{E})$ and a preorder $\succ$ on $\mathcal{B}$, we say"$\succ$-tameness" to refer to $\succ^+$-tameness.

**Theorem 6.3.24** (Tamed Linear Newmann's Lemma). *Let* $\mathtt{S} = (\mathcal{B}; \mathtt{R}, \mathtt{E})$ *be a linear RSM and* $\succ$ *a preorder on* $\mathcal{B}$ *compatible with* $\mathtt{S}$*. If* $\succ$ *is terminating and every monomial local* $\mathtt{S}^+$*-branching is* $\succ$*-tamely* $\mathtt{S}^{\mathrm{st}}$*-congruent, then* $\mathtt{S}^+$ *is convergent.*

The Tamed Linear Newmann's Lemma 6.3.24 reduces the study of convergence to tamed congruence of monomial local branchings. The Tamed Linear Newmann's Lemma 6.3.24 and the Basis From Convergence Theorem 6.3.18 are the two main results of linear rewriting modulo theory, allowing one to deduce bases from a local analysis. Given that monomial local branchings are the critical branchings (see subsection ii.1.1) of linear rewriting theory, the Tamed Linear Newmann's Lemma 6.3.24 could also be called the "tamed linear critical lemma".

Our proof can be understood as a generalization of [86, Theorem 4.2.1, part (ii)], working with tamed $\mathtt{S}^{\mathrm{st}}$-congruence instead of $\mathtt{S}^{+}$-confluence.

*Proof.* It follows from hypothesis and Lemma 6.3.21 that every intersecting $\mathtt{S}^{+}$-branching is $\succ$-tamely $\mathtt{S}^{\mathrm{st}}$-congruent. Consider then an additive $\mathtt{S}^{+}$-local triple as pictured above. We wish to show that its canonical $\mathtt{S}^{\mathrm{st}}$-confluence is tamed, that is:

$$\lambda_1 s(r_1) + \lambda_2 s(r_2) + w \succ^{+} \lambda_1 t(r_1) + \lambda_2 t(s_2) + w.$$

This happens precisely if $s(r_2) \notin \widetilde{\dot{\mathrm{supp}}}_{\mathtt{E}}(t(r_1))$ or if $s(r_1) \notin \widetilde{\dot{\mathrm{supp}}}_{\mathtt{E}}(t(r_2))$. (These two cases correspond respectively to the upper branch or lower branch being positive.) Assuming otherwise, the compatibility of $\succ$ implies that $s(r_2) \succ s(r_1)$ and $s(r_1) \succ s(r_2)$. This contradicts the assumption that $\succ$ terminates.[15] We conclude that every local $\mathtt{S}^{+}$-branching is $\succ$-tamely $\mathtt{S}^{\mathrm{st}}$-congruent.

Note then that:

**Lemma 6.3.25.** *A branching is $\succ$-tamely $\mathtt{S}^{\mathrm{st}}$-congruent if and only if it is $\succ$-tamely $\mathtt{S}^{+}$-congruent.*

*Proof.* This follows from Lemma 6.3.10 and the fact that for $M, N_1, N_2, L \in \mathcal{P}(\mathcal{B})$, if $M \succcurlyeq^{\mathrm{rel}} N_1 \cup N_2$ and $N_1 \succcurlyeq^{\mathrm{rel}} L$, then $M \succcurlyeq^{\mathrm{rel}} N_1 \cup N_2$. □

It follows that every local $\mathtt{S}^{+}$-branching is $\succ$-tamely $\mathtt{S}^{+}$-congruent. We can now apply the tamed (abstract) Newmann's lemma (Lemma 6.2.8) and conclude. □

## 6.4 Higher rewriting modulo

This section introduces (weak) higher rewriting modulo, extending the work of Forest and Mimram [81], who studied the non-modulo and context-agnostic (see below) setting. Although we restrict to 3-dimensional rewriting, many notions are not specific to the 3-dimensional case; in particular, this section could be adapted to the 2-dimensional case (rewriting modulo in categories, including monoids).

[15] Note that the termination is necessary hypothesis to have $\succ$-tameness. This fixes a small gap in the proof of Theorem 4.2.1 in [86]. Similar gaps appear in [3, lemma 4.2.12] and [56, theorem 2.2.9].

To some extent, it serves as a blueprint for the next section, which deals with higher *linear* rewriting modulo.

### 6.4.1 Higher rewriting system modulo

Recall the notion of 3-sesquicategory (subsection 5.2.2), 3-sesquipolygraph and higher equivalence (subsection 5.2.5). For $\mathcal{C}$ a 3-sesquicategory and $\equiv$ a higher equivalence on $\mathcal{C}$, we say that $\equiv$ satisfies the independence axiom if the following holds:

> *independence axiom:* for every pair of 1-composable 3-cells $\mathsf{f}\colon \phi \to \phi'$ and $\mathsf{g}\colon \psi \to \psi'$ in $\mathsf{P}^*$, the 3-cells $(\phi' \star_1 \mathsf{g}) \star_2 (\mathsf{f} \star_1 \psi)$ and $(\mathsf{f} \star_1 \psi') \star_2 (\phi \star_1 \mathsf{g})$ are $\equiv$-equivalence:

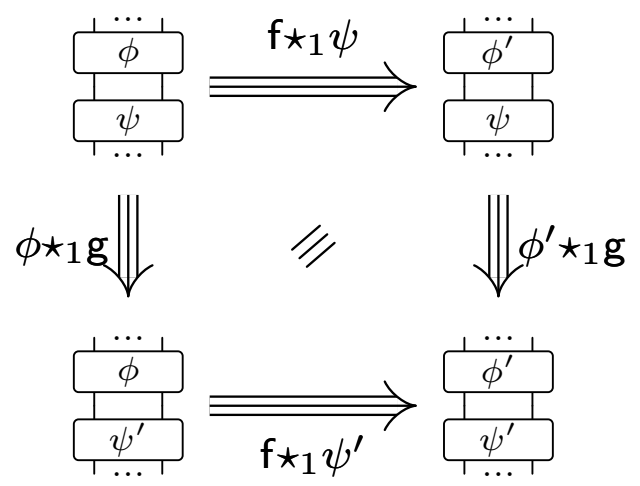

The independence axiom captures the interchange of 3-cells in $\mathsf{P}^*$.

**Definition 6.4.1.** *A* higher rewriting system *(higher* RS*)* $(\mathsf{P};\equiv)$ *is the data of a 3-sesquipolygraph* $\mathsf{P}$*, together with a higher equivalence* $\equiv$ *on* $\mathsf{P}^*$ *(see subsection 5.2.5) satisfying the independence axiom. A* higher rewriting system modulo *(higher* RSM*)* $\mathsf{S} = (\mathsf{R},\mathsf{E};\equiv)$ *is the data of two 3-sesquipolygraphs* $\mathsf{R}$ *and* $\mathsf{E}$ *with the same underlying 2-sesquipolygraph* $\mathsf{R}_2 = \mathsf{E}_2$*,*[16] *together with a higher equivalence* $\equiv$ *on* $\mathsf{R}^* \cup \mathsf{E}^\top$ *satisfying the independence axiom.*

Recall the notion of contexts (subsection 5.2.3). Given a higher RSM $\mathsf{S} = (\mathsf{R},\mathsf{E};\equiv)$ every choice of 1-sphere $\square$ in $\mathsf{R}_1^*$ defines an abstract RSM $\mathbb{S}(\square) =$

[16] In [55, 56], Dupont allows the more general definition $\mathsf{E}_2 \subset \mathsf{R}_2$ in the context of linear 3-polygraph. We prefer to avoid this generality, as it leads to ambiguity for the statement "$\mathsf{E}$ is convergent". For instance, for the linear 3-polygraphs $E$ and $R$ in [55, 56], we have that $E$ is convergent when viewed on $E_2$, but not when viewed on $R_2$. Similar issues appear in [57].

$(\mathbb{R}(\square), \mathbb{E}(\square); \equiv)$ on the underlying set $X(\square)$, where

$$X(\square) = \mathsf{R}_2^*(\square) = \mathsf{E}_2^*(\square),$$
$$\mathbb{R}(\square) = \mathrm{Cont}(\mathsf{R}_3)(\square) \quad \text{and} \quad \mathbb{E}(\square) = \mathrm{Cont}(\mathsf{E}_3)(\square).$$

We abuse notation and write $\equiv$ for its restriction on $(\mathsf{R}^* \cup \mathsf{E}^\top)(\square)$. Moreover, every context $\Gamma$ on $\square$ defines a morphism of abstract RSM

$$\Gamma \colon \mathbb{S}(\square) \to \mathbb{S}\big(s_1(\Gamma), t_1(\Gamma)\big). \tag{6.1}$$

(A morphism of abstract RSM is a pair of morphisms of 1-globular sets, the latter denoting a pair of functions commuting with the source and target maps; see subsection 5.2.1.) We say that a branching $(f, g)$ is a *contextualization* of another branching $(f', g')$ whenever $\Gamma[f', g'] = (f, g)$.

Hence, we think of a higher RSM as a category of abstract RSMs, related by contextualization. In principle, one could deal with each abstract RSM $\mathbb{S}(\square)$ independently, using the tools of abstract rewriting theory modulo (section 6.2). In practice, one would like to leverage the fact that these abstract RSMs gather together in a category, and relate to one another via contextualization. In that sense, higher rewriting theory (modulo) is nothing else that the combination of abstract rewriting theory (modulo) with the idea of contextualization.

In practice, termination can often only be obtained with context-dependent termination rule, where $r$ being a rewriting rule does not imply that $\Gamma[r]$ is. For that reason, we introduce the following notion:

**Definition 6.4.2.** *Let* $\mathsf{S} = (\mathsf{R}, \mathsf{E}; \equiv)$ *be a higher RSM. An* abstract sub-system $\mathbb{T}$ *is the data of a family of subsets* $\mathbb{T}(\square) \subset \mathsf{R}_3^*(\square)$ *for each 1-sphere* $\square$ *in* $\mathsf{R}_1^*$*, with the property that if* $\Gamma[r] \in \mathbb{T}(s(\Gamma), t(\Gamma))$*, then* $r \in \mathbb{T}(s_1(r), t_1(r))$.

We think of $\mathbb{T}$ as a family of sub-abstract RSMs $\mathbb{T}(\square) \subset \mathbb{S}(\square)$ modulo $\mathsf{E}(\square)$ and on the set $X(\square)$, and write $\mathbb{T} \subset \mathsf{S}$ to emphasize this point. We say that $\mathbb{T}$ is $\equiv$*-confluent* (resp. *terminating*, or *terminating modulo* $\mathsf{E}$) if every sub-abstract RSM $\mathbb{T}(\square)$ is $\equiv$-confluent (resp. terminating modulo $\mathsf{E}(\square)$), and similarly for other notions of abstract rewriting theory.

One can always consider $\mathsf{S}$ as its own abstract sub-system, so that the notion of abstract sub-system generalizes the notion of higher RSM. We shall often refer to the former as being *context-dependent*, and to the latter as being *context-agnostic*.

### 6.4.2 Compatibility and contextualization

The notion of compatible abstract preorder generalizes verbatim to the higher setting:

**Definition 6.4.3.** *Let* $\mathsf{S} = (\mathsf{R}, \mathsf{E}; \equiv)$ *be a higher RSM and* $\mathbb{T} \subset \mathsf{S}$ *an abstract sub-system. A preorder on* $\mathsf{R}_2^*$ *is said to be* compatible with $\mathbb{T}$ *if it is compatible with every sub-abstract RSM* $\mathbb{T}(\square)$*, in the sense of Definition 6.2.7. We denote* $\succ_{\mathbb{T}}$ *the smallest preorder compatible with* $\mathbb{T}$*.*

As customary, $\mathsf{E}$-invariance of $\succ$ is left implicit.

We now describe contextualization:

**Lemma 6.4.4.** *Let* $\mathsf{S} = (\mathsf{R}, \mathsf{E}; \equiv)$ *be a higher RSM and* $\mathbb{T} \subset \mathsf{S}$ *an abstract sub-system. Let* $(f, g)$ *be a local* $\mathbb{T}$*-branching and* $\Gamma$ *a context. Assume* $(f, g)$ *admits a* $(\mathbb{T}, \equiv)$*-confluence* $(f', g')$ *such that* $\Gamma[f', g']$ *belongs to* $\mathbb{T}$*. Then* $\Gamma[f, g]$ *is* $\mathbb{T}$*-confluent.*

*Proof.* Since $\equiv$ is a higher equivalence, $f' \circ f \equiv g' \circ g$ implies that $\Gamma[f' \circ f] \equiv \Gamma[g' \circ g]$. □

In particular, in the context-agnostic case:

**Lemma 6.4.5.** *Let* $\mathsf{S} = (\mathsf{R}, \mathsf{E}; \equiv)$ *be a higher RSM. Let* $(f, g)$ *be a local* $\mathsf{S}$*-branching and* $\Gamma$ *a context. If* $(f, g)$ *is* $(\mathsf{S}, \equiv)$*-confluent, so is* $\Gamma[f, g]$*.* □

### 6.4.3 Independent branchings

The higher structure of a higher RS $(\mathsf{P}; \equiv)$ induces canonical types of local $\mathsf{P}$-branchings. A local $\mathsf{P}$-branching is *independent* if it is of the form $(f, g) = (\mathsf{f} \star_1 \psi, \phi \star_1 \mathsf{g})$, for $\mathsf{P}$-rewriting steps $\mathsf{f} \colon \phi \to \phi'$ and $\mathsf{g} \colon \psi \to \psi'$:

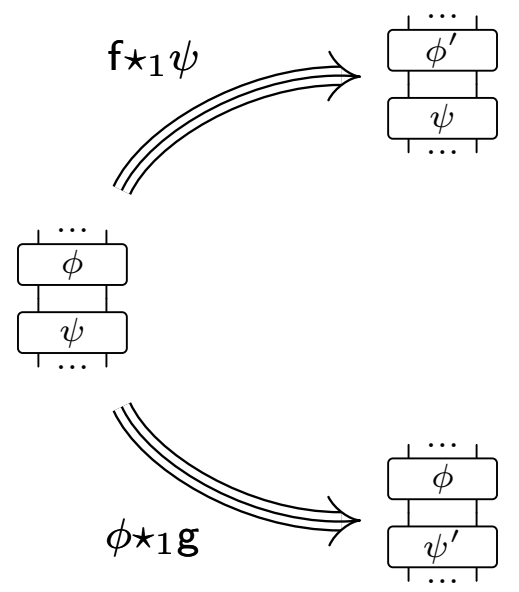

When working modulo, that is with a higher RSM $\mathsf{S} = (\mathsf{R}, \mathsf{E}; \equiv)$, independent branchings refer to independent R-branchings. Note that the property of being an independent branching is preserved by context: if $(f, g)$ is an independent branching, so is $\Gamma[f, g]$. A local S-branching (or a S-local triple) is said to be *overlapping* if it is not branchwise E-congruent to an independent R-branching.

An independent branching always have a canonical S-confluence, given by

$$(g', f') := (\phi' \star_1 \mathsf{g}, \mathsf{f} \star_1 \psi').$$

The following is tautological, thanks to the independence axiom:

**Lemma 6.4.6.** *Let* $\mathsf{S} = (\mathsf{R}, \mathsf{E}; \equiv)$ *be a higher RSM,* $\mathbb{T} \subset \mathsf{S}$ *an abstract subsystem and* $\succ$ *a preorder on* $\mathsf{R}_2^*$ *compatible with* $\mathbb{T}$*. Let* $(f, g)$ *is an independent* $\mathbb{T}$*-branching. If its canonical* $\mathbb{T}$*-confluence* $(g', f')$ *is in* $\mathbb{T}$*, then it defines a* $(\mathbb{T}, \equiv)$*-confluence for* $(f, g)$*.* ☐

In particular, in the context-agnostic case:

**Lemma 6.4.7.** *Let* $\mathsf{S} = (\mathsf{R}, \mathsf{E}; \equiv)$ *be a higher RSM. Every independent* S*-branching is* $(\mathsf{S}, \equiv)$*-confluent.* ☐

### 6.4.4 Higher Newmann's lemma

In the context-agnostic case, the above results can be gathered in a single black-box:

**Lemma 6.4.8** (higher Newmann's lemma)**.** *Let* $\mathsf{S} = (\mathsf{R}, \mathsf{E}; \equiv)$ *be a higher RSM. If* S *is terminating and if every overlapping local* S*-branching is, up to branchwise* E*-congruence, a contextualization of a* $(\mathsf{S}, \equiv)$*-confluent branching, then* S *is convergent.*

*Proof.* This follows from the (abstract) Newmann's lemma (Lemma 6.2.5), in combination with the branchwise confluence lemma (Lemma 6.2.9), confluence of contextualizations (Lemma 6.4.5) and confluence of independent branchings (Lemma 6.4.7). ☐

Given that overlapping local S-branchings are the critical branchings of higher rewriting theory, the higher Newmann's lemma could also be called the "higher critical lemma".

### 6.4.5 Independent rewriting

Recall the abstract notion of rewriting branchings as introduced in subsection 6.2.4. In the higher setting, one often wants to rewrite part of a diagram away from a given branching. This gives the notion of *independent rewriting*. Given rewriting steps $\mathsf{f}\colon \phi \to \phi_1$, $\mathsf{g}\colon \phi \to \phi_2$ and $\mathsf{h}\colon \psi \to \psi'$, setting

$$f = \mathsf{f} \star_1 \psi, \quad g = \mathsf{g} \star_1 \psi \quad \text{and} \quad h = \phi \star_1 \mathsf{h}$$

defines three pairs of branchings $(f, g)$, $(f, h)$ and $(h, g)$, the latter two being independent branchings. Define also

$$f' = \mathsf{f} \star_1 \psi', \quad g' = \mathsf{g} \star_1 \psi', \quad h_1 = \phi_1 \star_1 \mathsf{h} \quad \text{and} \quad h_2 = \phi_2 \star_1 \mathsf{h}.$$

We say that the branching $(f, g)$ *rewrites* into the branching $(f', g')$ via the triple $(h_1, h, h_2)$ (and similarly if the vertical positions of $f, g$ and $h$ are swapped), as pictured below:

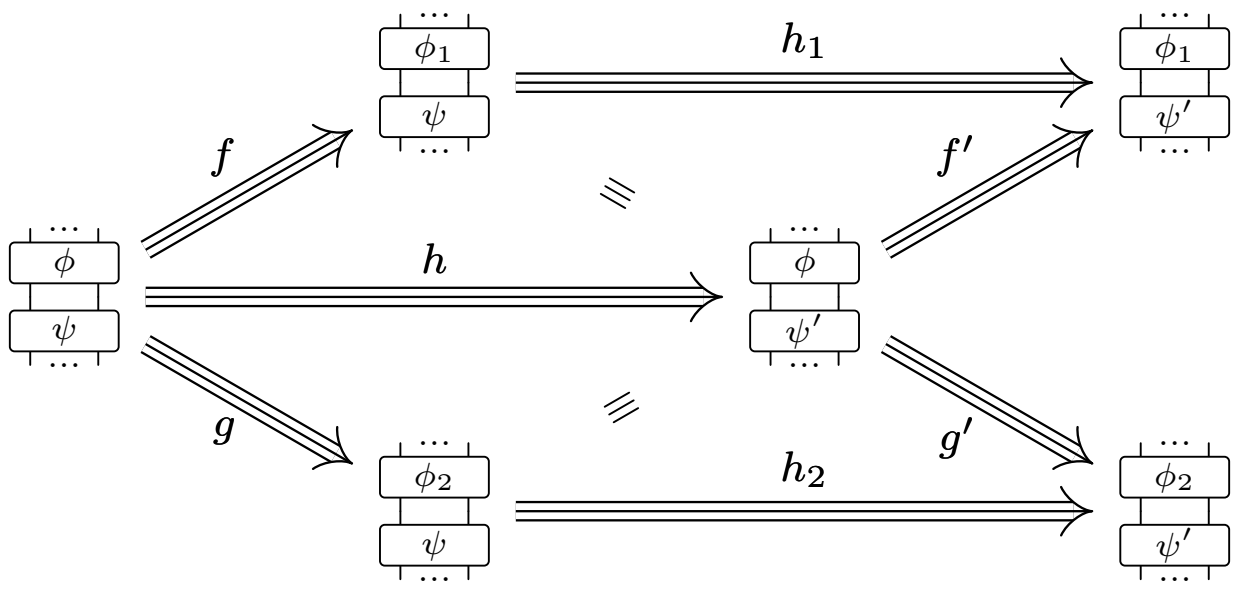

We have the following lemma:

**Lemma 6.4.9.** *Let* $\mathsf{S} = (\mathsf{R}, \mathsf{E}; \equiv)$ *be a higher RSM,* $\mathbb{T} \subset \mathsf{S}$ *an abstract sub-system and* $\succ$ *a preorder on* $\mathsf{R}_2^*$ *compatible with* $\mathbb{T}$*. In the situation above, assume* $(f, g)$ *is a* $\mathbb{T}$*-branching. If* $h_1$ *and* $h_2$ *are in* $\mathbb{T}$ *and if* $(f', g')$ *is* $\succ$*-tamed* $(\mathbb{T}, \equiv)$*-congruent, then so is* $(f, g)$*.* □

In particular, in the context-agnostic case:

**Lemma 6.4.10.** *Let* $\mathsf{S} = (\mathsf{R}, \mathsf{E}; \equiv)$ *be a higher RSM and* $\succ$ *a preorder on* $\mathsf{R}_2^*$ *compatible with* $\mathsf{S}$*. In the situation above, if* $(f', g')$ *is* $\succ$*-tamed* $(\mathsf{S}, \equiv)$*-congruent, then so is* $(f, g)$*.* □

### 6.4.6 Gray rewriting system modulo

A *Gray rewriting system* $(\mathsf{P}; \equiv)$ is the data of a higher RS $(\mathsf{P}, \equiv)$ such that $\mathsf{P}$ is a Gray polygraph and $\equiv$ verifies the interchange naturality axiom:

*interchange naturality axioms:* for all 0-composable $A, g, \beta$ with $A\colon \phi \Rrightarrow \phi'\colon f \Rightarrow f' \in \mathsf{P}_3$, $g \in \mathsf{P}_1^*$ and $\beta\colon h \Rightarrow h' \in \mathsf{P}_2$, and for all 0-composable $\alpha, g, B$ with $\alpha\colon f \Rightarrow f' \in \mathsf{P}_2$, $g \in \mathsf{P}_1^*$ and $B\colon \psi \Rrightarrow \psi'\colon h \Rightarrow h' \in \mathsf{P}_3$, we have

$$X_{\phi',g,\beta} \star_2 \Big( (A \star_0 g \star_0 h') \star_1 (f \star_0 g \star_0 \beta) \Big) \equiv \Big( (f' \star_0 g \star_0 \beta) \star_1 (A \star_0 g \star_0 h) \Big) \star_2 X_{\phi,g,\beta}$$

$$\text{and} \quad X_{\alpha,g,\psi'} \star_2 \Big( (\alpha \star_0 g \star_0 h') \star_1 (f \star_0 g \star_0 B) \Big) \equiv \Big( (f' \star_0 g \star_0 B) \star_1 (\alpha \star_0 g \star_0 h) \Big) \star_2 X_{\alpha,g,\psi},$$

where $X_{\phi',g,\beta}$ denotes a composition of interchange generators, interchanging 2-cells in $\psi$ with the 2-cell $\beta$, and similarly for the other $X$s. These interchange naturality axioms can be pictured as:

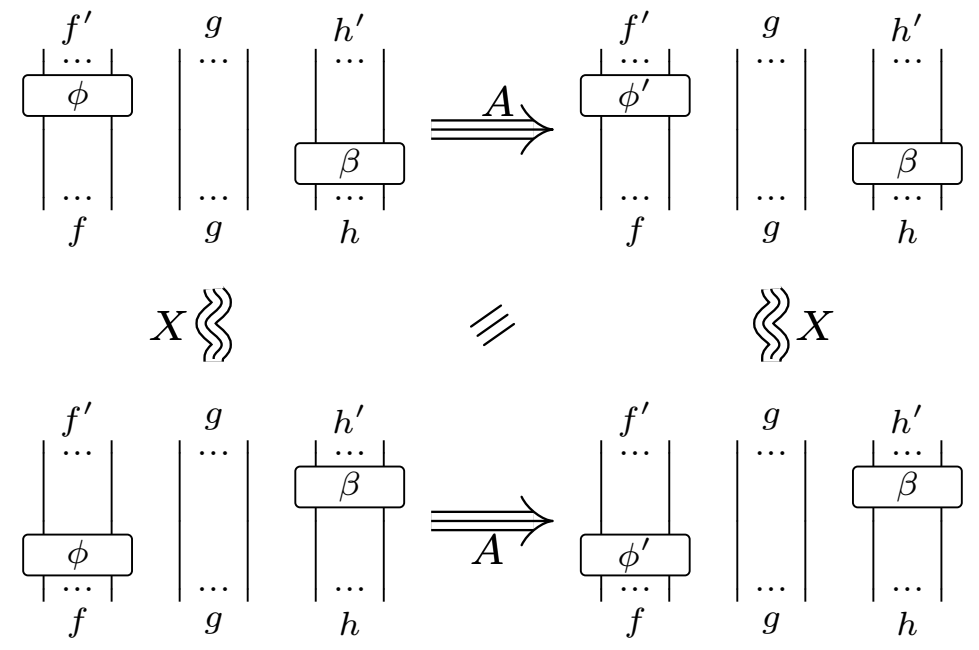

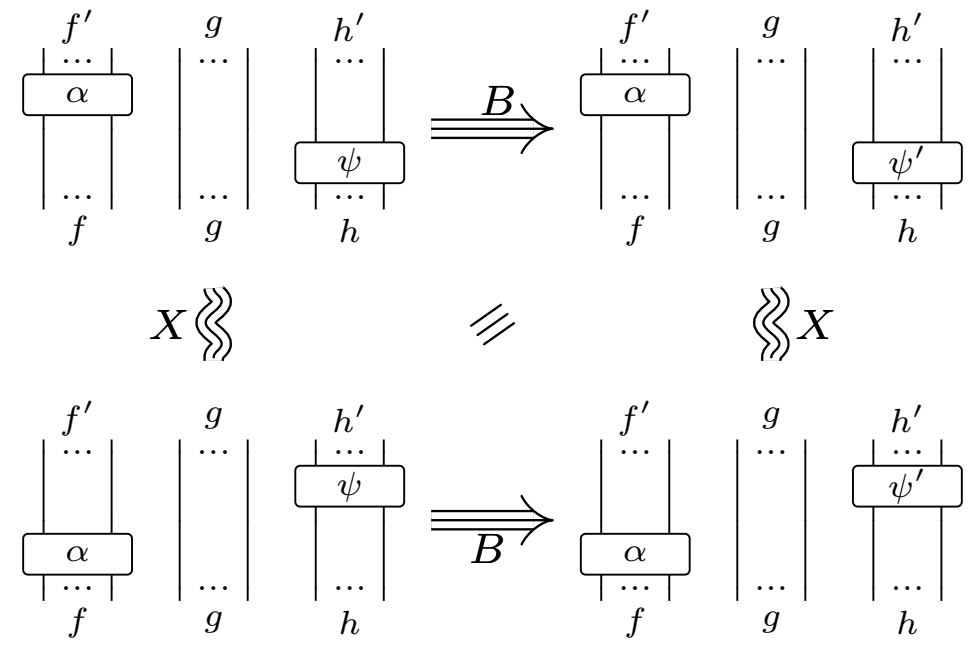

A *Gray rewriting system modulo* (Gray RSM) $\mathsf{S} = (\mathsf{R}, \mathsf{E}; \equiv)$ is the data of a higher RSM such that $\mathsf{E}$ is a Gray polygraph.

## 6.5 Higher linear rewriting modulo

This section finally introduces (weak) higher linear rewriting modulo, building on all the previous sections. In particular, linear Gray rewriting modulo is defined in subsection 6.5.5. Similarly to section 6.4, while we restrict to 3-dimensional rewriting, many notions are not specific to the 3-dimensional case; in particular, this could be adapted to the 2-dimensional case (rewriting modulo in linear categories, including algebras).

### 6.5.1 Higher linear rewriting system modulo

Recall the notions introduced in section 5.3.

**Definition 6.5.1.** *A* higher linear rewriting system *(higher linear RS) is the same data as a linear 3-sesquipolygraph* $\mathsf{P}$. *A* higher linear rewriting system modulo *(higher linear RSM)* $\mathsf{S} = (\mathsf{R}, \mathsf{E})$ *is the data of a left-monomial 3-sesquipolygraph* $\mathsf{R}$ *and a monomial-invertible 3-sesquipolygraph* $\mathsf{E}$*, such that* $\mathsf{R}_{\leq 2} = \mathsf{E}_{\leq 2}$*.*

As in the abstract higher case, a higher linear RSM $\mathsf{S} = (\mathsf{R}, \mathsf{E})$ defines a graph of linear 1-polygraphs $\mathtt{S}(\square) = (\mathcal{B}(\square), \mathtt{R}(\square), \mathtt{E}(\square))$, setting

$$\mathcal{B}(\square) := \mathsf{R}_2^*(\square) = \mathsf{E}_2^*(\square),$$
$$\mathtt{R}(\square) := \mathrm{Cont}(\mathsf{R})(\square) \quad \text{and} \quad \mathtt{E}(\square) := \mathrm{Cont}(\mathsf{E})(\square),$$

for each 1-sphere $\square$ in $\mathsf{R}_1^* = \mathsf{E}_1^*$. Each context $\Gamma$ on $\square$ defines a morphism of linear 1-polygraphs:

$$\Gamma \colon \mathtt{S}(\square) \to \mathtt{S}(s_1(\Gamma), t_1(\Gamma)).$$

Contrary to the abstract higher case however, $\mathsf{S}(\square)$ *needs not be a linear RSM*, as we did not impose the adaptedness condition appearing in Definition 6.3.2.

As for the abstract higher case, we define a notion of sub-systems to properly speak of context-dependent rewriting rules, imposing now the adaptedness condition:

**Definition 6.5.2.** *Let* $\mathsf{S} = (\mathsf{R}, \mathsf{E})$ *be a higher linear RSM. A* linear sub-system $\mathtt{T}$ *is the data of a family of sub-sets* $\mathtt{T}(\square) \subset \mathsf{R}_3^*(\square)$ *for each 1-sphere* $\square$ *in* $\mathsf{R}_1^*$*, such that if* $\Gamma[r] \in \mathtt{T}(s_1(\Gamma), t_1(\Gamma))$*, then* $r \in \mathtt{T}(s_1(r), t_1(r))$*. Moreover, we assume the* adaptedness condition*:*

$$s(r) \notin \widehat{\dot{\mathrm{supp}}}_{\mathtt{E}}(t(r)) \quad \forall r \in \mathtt{T}.$$

*We say that* $\mathsf{S}$ *is* adapted *if* $\mathsf{S}$ *is a linear sub-system of itself.*

As in the higher abstract case, we think of $\mathtt{T}$ as defining a family of sub-linear RSM $\mathtt{T}(\square) \subset \mathsf{S}(\square)$ modulo $\mathsf{E}(\square)$, and write $\mathtt{T} \subset \mathsf{S}$ to emphasize that point. We say that $\mathtt{T}$ is $\equiv$*-confluent* (resp. *terminating*, or *terminating modulo* $\mathsf{E}$) if every sub-linear RSM $\mathtt{T}(\square)$ is confluent (resp. terminating modulo $\mathsf{E}(\square)$), and similarly for other notions of linear rewriting theory.

Note that adaptedness *must* be defined via a linear sub-system, since in general, it depends on context: if $r$ is adapted, $\Gamma[r]$ needs not be.

A $\mathtt{T}$-rewriting step has the following general form (compare with Definition 6.3.8):

$$\alpha \colon u \quad \sim_{\mathtt{E}} \quad \lambda\Gamma[s(r)] + v \quad \dashrightarrow_{\mathtt{T}} \quad \lambda\Gamma[t(r)] + v \quad \sim_{\mathtt{E}} \quad w.$$

We say that $\alpha$ is *of type* $r$.[17] As before, $\alpha$ is positive if $\Gamma[s(r)] \notin \widehat{\dot{\mathrm{supp}}}_{\mathtt{E}}(v)$. Note that *contextualization needs not preserve positivity*: we may have

$$\Gamma'[\Gamma[s(r)]] \in \widehat{\dot{\mathrm{supp}}}_{\mathtt{E}}(\Gamma'[v]), \qquad \text{even if} \qquad \Gamma[s(r)] \notin \widehat{\dot{\mathrm{supp}}}_{\mathtt{E}}(v).$$

However, it is true that if $\alpha$ is monomial (that is, if $v = 0$) then $\Gamma[\alpha]$ is monomial.

*Remark* 6.5.3. In [3, 4], [55, 56] and [57], the authors used quasi-terminating normal forms to propose a basis theorem suitable for their setting. Under the hood, choosing quasi-normal forms amounts to choosing a context-dependent termination. This perspective avoids the use of quasi-normal forms, which are typically badly behaved.[18]

[17] This is a slight abuse of notation: as belonging to a linear RSM, $\alpha$ is of type $\Gamma[r]$ according to Definition 6.3.7.

[18] For instance, a monomial in the support of a quasi-normal form needs not be a quasi-normal form, and a linear combination of quasi-normal forms needs not be a quasi-normal form; compare with Remark 6.3.15.

### 6.5.2 Strong compatibility and contextualization

**Definition 6.5.4.** *Let* $\mathsf{S} = (\mathsf{R}, \mathsf{E})$ *be a higher linear* RSM *and* $\mathtt{T}$ *a linear-sub system. A preorder* $\succ$ *on* $\mathsf{R}_2^*$ *is said to be* strongly compatible with $\mathtt{T}$ *if:*

$$\Gamma[s(r)] \succ b, \quad \textit{for all } r \in \mathtt{T}, \textit{context } \Gamma \textit{ and } b \in \Gamma\big[\dot{\widetilde{\mathrm{supp}}}_{\mathtt{E}}(t(r))\big].$$

*We denote* $\succ_{\mathtt{T}}$ *the smallest preorder strongly compatible with* $\mathtt{T}$.

As customary, $\mathsf{E}$-invariance of $\succ$ is left implicit. For any vector $w$, we have

$$\Gamma\big[\dot{\widetilde{\mathrm{supp}}}_{\mathtt{E}}(w)\big] \supseteq \dot{\widetilde{\mathrm{supp}}}_{\mathtt{E}}(\Gamma[w]).$$

This inclusion needs not be an equality! Recall that in general, contextualization does not act freely (see subsection ii.2.2): there may exist monomials $b_1, b_2$ such that $b_1 \not\dot{\sim}_{\mathtt{E}} b_2$ but $\Gamma[b_1] \dot{\sim}_{\mathtt{E}} \Gamma[b_2]$. In particular, if $w = \lambda_1 b_1 + \lambda_2 b_2$ and $\lambda_1, \lambda_2$ have well-chosen values, we have $\Gamma[w] = 0$ and thus $\Gamma[b_1] \notin \dot{\widetilde{\mathrm{supp}}}_{\mathtt{E}}(\Gamma[w])$.

Recall that for a linear RSM $\mathtt{S} = (\mathcal{B}; \mathtt{R}, \mathtt{E})$, a preorder $\succ$ is compatible (Definition 6.3.19) if for all $r \in \mathtt{R}$, $s(r) \succ b$ for all $b \in \widetilde{\mathrm{supp}}_{\mathtt{E}}(t(r))$. In particular, strong compatibility implies compatibility with each linear RSM $\mathtt{T}(\square)$. The converse does not hold: by the paragraph above, having that $\Gamma[s(r)] \succ b$ for all $b \in \dot{\widetilde{\mathrm{supp}}}_{\mathtt{E}}(\Gamma[t(r)])$ is not sufficient to have strong compatibility.

*Remark* 6.5.5. In practice, one checks strong compatibility by first checking that $\succ$ is $\mathsf{E}$-invariant, and then that for each generating relation $r \in \mathtt{T} \cap \mathsf{R}$, one has: (i) $s(r) \succ b$ for each $b \in \dot{\widetilde{\mathrm{supp}}}_{\mathtt{E}}(t(r))$, and (ii) if $\Gamma$ is a context such that $\Gamma[r] \in \mathtt{T}$, then $\Gamma[s(r)] \succ b$.

**Lemma 6.5.6.** *Let* $\mathsf{S} = (\mathsf{R}, \mathsf{E})$ *be a higher linear RSM,* $\mathtt{T} \subset \mathsf{S}$ *a linear sub-system and* $\succ$ *a preorder on* $\mathsf{R}_2^*$ *strongly compatible with* $\mathtt{T}$. *Consider a composition* $g \star_2 f$*, where* $f$ *is a monomial* $\mathtt{T}^+$*-rewriting step and* $g$ *is a* $\mathtt{T}^+$*-rewriting sequence. If* $\Gamma$ *is a context such that both* $\Gamma[f]$ *and* $\Gamma[g]$ *are in* $\mathtt{T}$*, then* $\Gamma[g]$ *is* $\succ$*-tamed by* $\Gamma[s(f)]$.

*Proof.* Decomposing $g$ as a composition of $\mathtt{T}^+$-rewriting steps, the situation is:

$$\begin{array}{ccccccccc} s(f) & \xrightarrow[\mathtt{T}]{f} & v_0 & \xrightarrow[\mathtt{T}]{g_0} & v_1 & \xrightarrow[\mathtt{T}]{g_1} & \cdots & \xrightarrow[\mathtt{T}]{g_{n-1}} & v_n \\ & & & & & & & & \downarrow \Gamma \\ \Gamma[s(f)] & \xrightarrow[\mathtt{T}]{\Gamma[f]} & \Gamma[v_0] & \dashrightarrow_{\mathtt{T}}^{\Gamma[g_0]} & \Gamma[v_1] & \dashrightarrow_{\mathtt{T}}^{\Gamma[g_1]} & \cdots & \dashrightarrow_{\mathtt{T}}^{\Gamma[g_{n-1}]} & \Gamma[v_n] \end{array}$$

We wish to show that $\Gamma[s(f)] \succ \dot{\widetilde{\mathrm{supp}}}_{\mathsf{E}}(\Gamma[v_i])$ for all $0 \leq i \leq n$. Rather, given the inclusion $\Gamma\left[\dot{\widetilde{\mathrm{supp}}}_{\mathsf{E}}(w)\right] \supseteq \dot{\widetilde{\mathrm{supp}}}_{\mathsf{E}}(\Gamma[w])$ which holds for generic vector $w$, it suffices to show

$$\Gamma[s(f)] \succ \Gamma\left[\dot{\widetilde{\mathrm{supp}}}_{\mathsf{E}}(v_i)\right].$$

For $i = 0$, this follows from strong compatibility. It remains to show that

$$\Gamma\left[\dot{\widetilde{\mathrm{supp}}}_{\mathsf{E}}(v_i)\right] \succcurlyeq^{\mathrm{rel}} \Gamma\left[\dot{\widetilde{\mathrm{supp}}}_{\mathsf{E}}(v_{i+1})\right].$$

(Recall the definition of relative relation from Definition 6.3.23.) Write $g_i = \lambda r + v$, its canonical decomposition as a positive $\mathtt{T}$-rewriting step. On the one hand:

$$\Gamma\left[\dot{\widetilde{\mathrm{supp}}}_{\mathsf{E}}(t(r)+v)\right] \subset \Gamma\left[\dot{\widetilde{\mathrm{supp}}}_{\mathsf{E}}(t(r))\right] \cup \Gamma\left[\dot{\widetilde{\mathrm{supp}}}_{\mathsf{E}}(v)\right].$$

On the other hand, thanks to the positivity of $g_i$, we have $\dot{\widetilde{\mathrm{supp}}}_{\mathsf{E}}(s(r)+v) = \{s(r)\} \sqcup \dot{\widetilde{\mathrm{supp}}}_{\mathsf{E}}(v)$, and hence

$$\Gamma\left[\dot{\widetilde{\mathrm{supp}}}_{\mathsf{E}}(s(r)+v)\right] = \{\Gamma[s(r)]\} \cup \Gamma\left[\dot{\widetilde{\mathrm{supp}}}_{\mathsf{E}}(v)\right].$$

Strong compatibility concludes. □

**Corollary 6.5.7** (Contextualization Lemma). *Let* $\mathsf{S} = (\mathsf{R}, \mathsf{E})$ *be a higher linear RSM,* $\mathtt{T} \subset \mathsf{S}$ *a linear sub-system and* $\succ$ *a preorder on* $\mathsf{R}_2^*$ *strongly compatible with* $\mathtt{T}$. *Let* $(f, g)$ *be a monomial local* $\mathtt{T}^+$*-branching and* $\Gamma$ *a context. Assume* $(f, g)$ *admits a* $\mathtt{T}^+$*-confluence* $(f', g')$ *such that* $\Gamma[f', g']$ *belongs to* $\mathtt{T}$. *Then* $\Gamma[f, g]$ *is* $\succ$*-tamely* $\mathtt{T}^{\mathrm{st}}$*-congruent.*

In general, $\Gamma[f', g']$ needs not be positive,[19] as the following example, taken from chapter 7 (Lemma 7.3.9), illustrates. Recall that plain arrows denote positive rewriting steps, while dashed arrows denote not-necessarily-positive rewriting steps. In the example below, the graded interchange law is part of the modulo. On the left-hand side, the branches of the confluence (labelled “dd”) evaluate 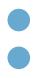 to zero, so that only the other diagram remains. The two diagrams on the top (resp. on the bottom) are *not* projectively congruent modulo, so that both branches of the confluence are positive. The same confluent branching is pictured on the right-hand side, contextualized with a cap and cup. Because the graded interchange law is part of the modulo, dots can now freely move between the top and the bottom. In particular, the two diagrams on the top (resp. on the bottom) are now projectively congruent modulo. In particular, the branches of the confluence are not positive anymore.

[19] In particular, there is are related gaps in the proofs of [3, lemma 4.2.12] and [56, theorem 2.2.9]. These gaps are distinct from the ones mentioned in Footnote 15.

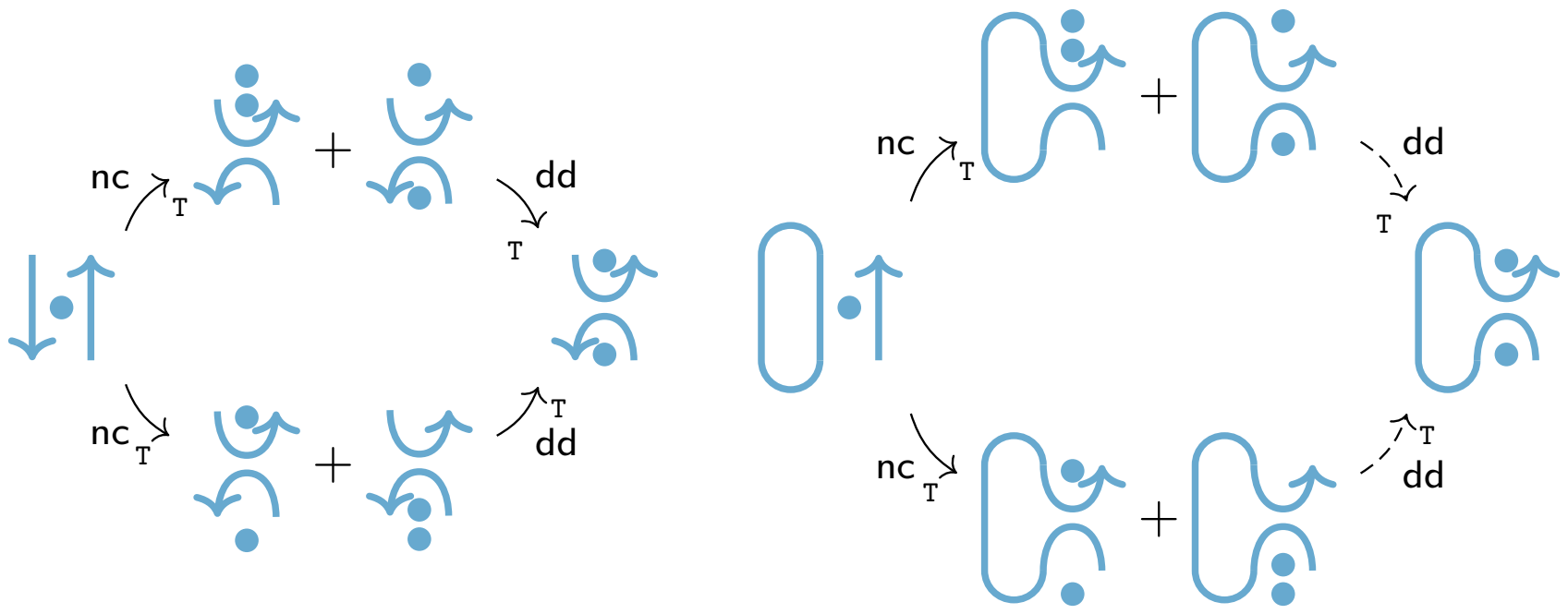

In the context-agnostic setting, the CONTEXTUALIZATION LEMMA 6.5.7 reduces to the following:

**Corollary 6.5.8.** *Let* $\mathsf{S} = (\mathsf{R}, \mathsf{E})$ *be an adapted higher linear RSM and* $\succ$ *a preorder on* $\mathsf{R}_2^*$ *strongly compatible with* $\mathsf{S}$*. Let* $(f, g)$ *be a monomial local* $\mathsf{S}^+$*-branching and* $\Gamma$ *a context. If* $(f, g)$ *is* $\mathsf{S}^+$*-confluent, then* $\Gamma[f, g]$ *is* $\succ$*-tamely* $\mathsf{S}$*-congruent.*

## 6.5.3 Independent branchings

As we have seen, for higher RSM the statement that independent branchings are $\equiv$-confluent is tautological (Lemma 6.4.7). As we now explain, the linear case is more subtle.

Fix $\mathsf{P}$ a higher linear RS and $(f, g)$ an independent local $\mathsf{P}^+$-branching. Without loss of generality, this means that $f = \mathsf{f} \star_1 s(\mathsf{g})$ and $g = s(\mathsf{f}) \star_1 \mathsf{g}$, for some monomial $\mathsf{P}^+$-rewriting steps $\mathsf{f}$ and $\mathsf{g}$. The $\mathsf{P}$-confluence $(g', f') = (t(\mathsf{f}) \star_1 \mathsf{g}, \mathsf{f} \star_1 t(\mathsf{g}))$ defines a canonical $\mathsf{P}$-confluence for $(f, g)$. Write:

$$\mathsf{f} \colon s(\mathsf{f}) \to \lambda_1 x_1 + \ldots + \lambda_m x_m \quad \text{and} \quad \mathsf{g} \colon s(\mathsf{g}) \to \mu_1 y_1 + \ldots + \mu_n y_n,$$

where the $x_i$'s (resp. $y_j$'s) are monomials with pairwise distinct $\mathsf{E}$-projective classes. Then the $\mathsf{P}$-confluence $(g', f')$ can be explicitly decomposed as a $\mathsf{P}^{\mathrm{st}}$-confluence, setting $f' = f'_n \star_2 \ldots \star_2 f'_1$ and $g' = g'_m \star_2 \ldots \star_2 g'_1$, where:

$$f'_j = \sum_{k<j} \mu_k (t(\mathsf{f}) \star_1 y_k) + \mu_j (\mathsf{f} \star_1 y_j) + \sum_{j<k} \mu_k (s(\mathsf{f}) \star_1 y_k)$$

$$\text{and } g'_i = \sum_{k<i} \lambda_k(x_k \star_1 t(\mathrm{g})) + \lambda_i(x_i \star_1 \mathrm{g}) + \sum_{i<k} \lambda_k(x_k \star_1 s(\mathrm{g})).$$

**Lemma 6.5.9** (independent branchings, context-dependent case)**.** *Let* $\mathsf{S} = (\mathsf{R}, \mathsf{E})$ *be a higher linear RSM,* $\mathtt{T} \subset \mathsf{S}$ *a linear sub-system and* $\succ$ *a linear preorder strongly compatible with* $\mathtt{T}$*. Let* $(f, g)$ *be an independent local* $\mathtt{T}^+$*-branchings. If for each* $1 \leq j \leq n$ *(resp.* $1 \leq i \leq m$*) and with the notations above, the* $\mathsf{S}^+$*-rewriting step* $\mathrm{f} \star_1 y_j$ *(resp.* $x_i \star_1 \mathrm{g}$*) is in* $\mathtt{T}$*, then the canonical* $\mathtt{T}^{\mathrm{st}}$*-confluence* $(g', f')$ *is* $\succ$*-tamed by* $s(f) = s(g)$*.*

*Proof.* For each $1 \leq j \leq n$, the supports of $s(f'_j)$ and $t(f'_j)$ are contained in the following set:

$$M \coloneqq \bigcup_{1\leq j\leq n} \widetilde{\mathrm{supp}}_{\mathtt{E}}\big(t(f) \star_1 y_j\big) \cup \bigcup_{1\leq j\leq n} \widetilde{\mathrm{supp}}\big(s(f) \star_1 y_j\big).$$

Since $\mathrm{f}\star_1 y_i$ belongs to $\mathtt{T}$, it follows from strong compatibility that $s(f)\star_1 s(g) \succ s(f) \star_1 y_j$ and $s(f) \star_1 y_j \succ t(f) \star_1 y_j$, so that $s(f) \star_1 s(g) \succ_{\mathsf{S}} M$. In other words, $f'$ is $\succ$-tamed by $s(f) \star_1 s(g)$. An analogous statement holds for the branch $g'$. This concludes. □

Note the necessity of strong compatibility in the proof of Lemma 6.5.9, even in the context-agnostic case: $s(f) \star_1 s(g) \succ s(f) \star_1 y_j$ does *not* follow from the fact that $s(f) \star_1 g$ is a $\mathsf{S}^+$-rewriting step, because $s(f) \star_1 y_j$ may not be in the support of $s(f) \star_1 t(g)$.

The proof of Lemma 6.5.9 resembles the proof of Theorem 4.2.1, part (iii) in [86], with tameness appearing implicitly. In our terminology, their setting (associative algebra) restrict to context-agnostic systems, and as contextualization acts freely, there is no distinction between compatibility and strong compatibility.

The following is a direct corollary of Lemma 6.5.9:

**Lemma 6.5.10** (independent branchings, context-agnostic case)**.** *Let* $\mathsf{S} = (\mathsf{R}, \mathsf{E})$ *be an adapted higher linear RSM and* $\succ$ *a preorder on* $\mathsf{R}_2^*$ *strongly compatible with* $\mathsf{S}$*. Every independent* $\mathsf{S}^+$*-branching is* $\succ$*-tamely* $\mathsf{S}^{\mathrm{st}}$*-congruent.* □

In the context-dependent case, one cannot rely on such a general statement. However:

**Lemma 6.5.11.** *Let* $\mathsf{S} = (\mathsf{R}, \mathsf{E})$ *be a higher linear RSM,* $\mathtt{T} \subset \mathsf{S}$ *a linear sub-system and* $\succ$ *a preorder on* $\mathsf{R}_2^*$ *strongly compatible with* $\mathtt{T}$*. Assume given another linear sub-system* $\mathtt{B} \subset \mathtt{T}$*, such that* $\mathtt{B}^+$ *is convergent. If every* $(\mathsf{S}\backslash\mathtt{T})$*-rewriting step is* $\mathtt{B}$*-congruent, then every independent* $\mathtt{T}$*-branching is* $\succ$*-tamely* $\mathtt{T}^{\mathrm{st}}$*-congruent.*

This will be the situation when discussing the rewriting theory of graded $\mathfrak{gl}_2$-foams in the next chapter. Note that in the above situation, two 2-cells are $\mathsf{S}$-congruent if and only if they are $\mathtt{T}$-congruent, so that $\mathsf{S}$ and $\mathtt{T}$ present the same underlying family of modules.

*Proof.* We adapt the proof of Lemma 6.5.9, and borrow its notations. Assume, without loss of generality, that $\mathsf{f} \star_1 y_j$ is in $\mathtt{T}$ for $1 \leq j \leq k$ and that $\mathsf{f} \star_1 y_j$ is *not* in $\mathtt{T}$ for $k < j \leq m$. By the same argument as in the proof of Lemma 6.5.9, $f'_k \star_2 \ldots \star_2 f'_1$ is a $\mathtt{T}^{\mathrm{st}}$-rewriting sequence $\succ$-tamed by $s(\mathsf{f}) \star_1 s(\mathsf{g})$. On the other hand, by the hypothesis of the lemma, $f'_m \star_2 \ldots \star_2 f'_{k+1}$ is $\mathtt{B}$-congruent. Arguing similarly for $g'$, we get a $\mathtt{T}$-congruence for $(f, g)$ given as a $\mathtt{B}$-congruence "sandwiched" between two $\succ$-tamed $\mathtt{T}^{\mathrm{st}}$-congruence. Since $\mathtt{B}^+$ is convergent, we can replace the $\mathtt{B}$-congruence by a $\mathtt{B}^+$-confluence. This concludes. □

### 6.5.4 Independent rewriting

We consider the linear analogue of abstract independent rewriting defined in subsection 6.4.5. Let us use the same notations with $\mathsf{f}$, $\mathsf{g}$ and $\mathsf{h}$ monomial rewriting steps, which decompose as:

$$\begin{aligned} \mathsf{f}\colon \phi &\to_{\mathsf{R}} \lambda_1 x_1 + \ldots + \lambda_l x_l, \\ \mathsf{g}\colon \phi &\to_{\mathsf{R}} \mu_1 y_1 + \ldots + \mu_m y_m, \\ \text{and}\quad \mathsf{h}\colon \psi &\to_{\mathsf{R}} \nu_1 z_1 + \ldots + \nu_n z_n, \end{aligned}$$

where $x_i$'s (resp. $y_j$'s and $z'_k s$) are monomials belonging to distinct $\mathsf{E}$-projective classes. In these decompositions, the independent rewriting of $(f, g)$ into $(f', g')$ via the triple $(h, h_1, h_2)$ is pictured as follows:

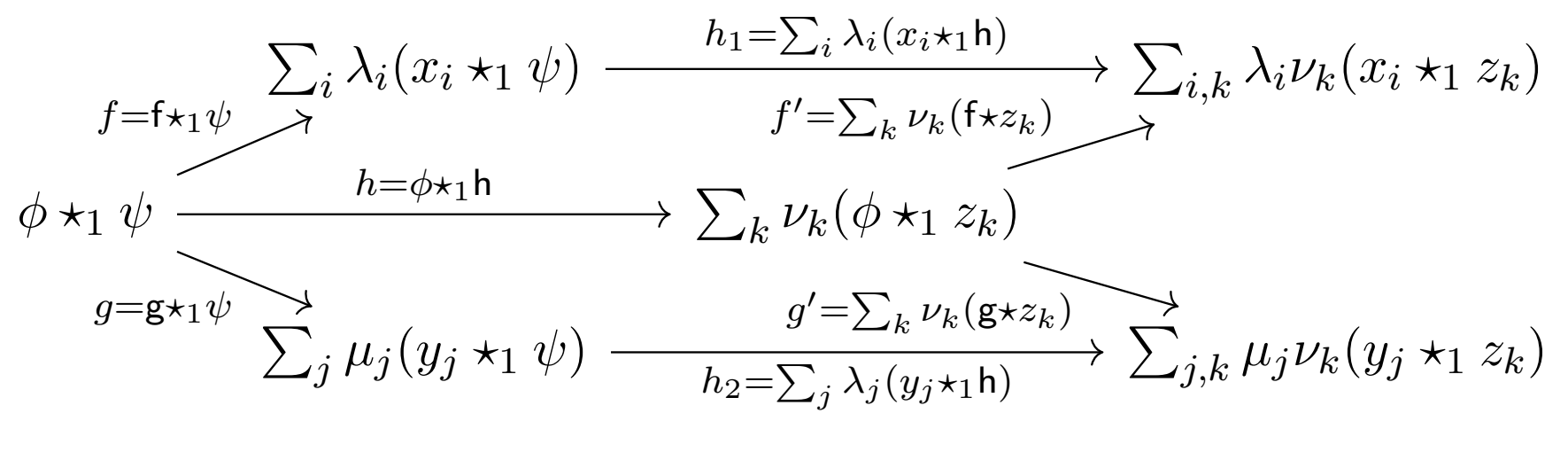

**Lemma 6.5.12** (Independent Rewriting Lemma). *Let* $\mathsf{S} = (\mathsf{R}, \mathsf{E})$ *be a higher linear RSM,* $\mathsf{T} \subset \mathsf{S}$ *a linear sub-system and* $\succ$ *a preorder on* $\mathsf{R}_2^*$ *strongly compatible with* $\mathsf{T}$*. Let* $\mathsf{f}$, $\mathsf{g}$ *and* $\mathsf{h}$ *monomial* $\mathsf{R}^+$*-rewriting steps as above, such that:*

*(i)* $f = \mathsf{f} \star_1 \psi$, $g = \mathsf{g} \star_1 \psi$ *and* $h = \phi \star_1 \mathsf{h}$ *are in* $\mathsf{T}$*, and for each* $1 \leq i \leq l$ *(resp.* $1 \leq j \leq m$*, resp.* $1 \leq k \leq n$*), the rewriting step* $x_i \star_1 \mathsf{h}$ *(resp.* $y_j \star_1 \mathsf{h}$*, resp.* $\mathsf{f} \star_1 z_k$ *and* $\mathsf{g} \star_1 z_k$*) is in* $\mathsf{T}$*,*

*(ii) for each* $1 \leq k \leq n$*, the branching* $(\mathsf{f} \star_1 z_k, \mathsf{g} \star_1 z_k)$ *is* $\succ$*-tamely* $\mathsf{T}^{\mathrm{st}}$*-congruent,*

*then* $(f, g)$ *is* $\succ$*-tamely* $\mathsf{T}^{\mathrm{st}}$*-congruent.*

Condition (i) simply says that all monomial $\mathsf{R}^+$-rewriting steps involved are in $\mathsf{T}$ (note that $f$, $g$ and $h$ in $\mathsf{T}$ imply that $\mathsf{f}$, $\mathsf{g}$ and $\mathsf{h}$ are in $\mathsf{T}$).

*Proof.* By strong compatibility, we have $\phi \star_1 \psi \succ \phi \star_1 z_k$. Let $(f'_k, g'_k)$ be a $\mathsf{T}^{\mathrm{st}}$-congruence for $(\mathsf{f} \star_1 z_k, \mathsf{g} \star_1 z_k)$ tamed by $\phi \star_1 z_k$, and hence by $\phi \star_1 \psi$. The $\mathsf{T}^{\mathrm{st}}$-congruence $\sum_k \nu_k (f'_k, g'_k)$ defines a $\mathsf{T}^{\mathrm{st}}$-congruence for $\sum_k \nu_k (\mathsf{f} \star_1 z_k, \mathsf{g} \star_1 z_k)$, tamed by $\phi \star_1 \psi$. On the other hand, strong compatibility implies that $\phi \star_1 \psi \succ x_i \star_1 \psi$ and $x_i \star_1 \psi \succ x_i \star_1 z_k$, for all suitable $i$ and $k$, so that $h_1$ admits a decomposition as a $\mathsf{T}^{\mathrm{st}}$-congruence tamed by $\phi \star_1 \psi$. The same statement holds for $h_2$, which concludes. □

**Lemma 6.5.13** (independent rewriting lemma, context-agnostic case). *Let* $\mathsf{S} = (\mathsf{R}, \mathsf{E})$ *be a higher linear RSM and* $\succ$ *a preorder on* $\mathsf{R}_2^*$ *strongly compatible with* $\mathsf{S}$*. If for each* $1 \leq k \leq n$*, the branching* $(\mathsf{f} \star_1 z_k, \mathsf{g} \star_1 z_k)$ *is* $\succ$*-tamely* $\mathsf{T}^{\mathrm{st}}$*-congruent, then* $(f, g)$ *is* $\succ$*-tamely* $\mathsf{T}^{\mathrm{st}}$*-congruent.*

### 6.5.5 Linear Gray rewriting system modulo

Linear Gray rewriting modulo is simply higher linear rewriting modulo where the modulo contains graded interchangers. Below are the formal definitions.

Recall the notations $G$ and $\mu$ from Notation 5.3.1:

**Definition 6.5.14.** *A* $(G, \mu)$*-graded linear Gray rewriting system modulo is a higher linear RSM* $\mathsf{S} = (\mathsf{R}, \mathsf{E})$ *such that* $\mathsf{E}$ *is a* $(G, \mu)$*-graded linear Gray polygraph.*

In particular, we can specialize to the scalar case:

**Definition 6.5.15.** *A* $(G,\mu)$-scalar Gray rewriting system modulo $\mathsf{S} := (\mathsf{R},\mathsf{E};\mathrm{scl})$ *is the data of two scalar* $G$*-graded 3-sesquipolygraphs* $(\mathsf{R};\mathrm{scl})$ *and* $(\mathsf{E};\mathrm{scl})$ *with* $\mathsf{R}_2^* = \mathsf{E}_2^*$*, such that* $\mathsf{E}$ *is a* $(G,\mu)$*-scalar Gray polygraph.*

*Remark* 6.5.16. Gray rewriting modulo corresponds to rewriting modulo in strict 2-categories, and linear Gray rewriting corresponds to rewriting in graded-2-categories. In particular, linear Gray rewriting where $\mu$ is trivial corresponds to rewriting in linear strict 2-categories.

## 6.6 Summary

Given a presented linear 2-sesquicategory, how can one use rewriting theory to find a basis? This section summarizes the main tools developed throughout this chapter. While we focus on the linear case (which is more involved), the same ideas apply for higher (non-linear) rewriting modulo.

### 6.6.1 The setup

Let $\mathcal{C}$ be a linear 2-sesquicategory presented by a linear 3-sesquipolygraph $\mathsf{P}$ in the sense of Definition 5.3.2. Assume the following choices of data has been made:

(a) a splitting of $\mathsf{P}$ as $\mathsf{P} = \mathsf{R} \sqcup \mathsf{E}$, defining a higher linear RSM $\mathsf{S} = (\mathsf{R},\mathsf{E})$ presenting $\mathcal{C}$ (see Definition 6.5.1),

(b) a linear sub-system $\mathtt{T} \subset \mathsf{S}$ (see Definition 6.5.2),

(c) a preorder $\succ$ strongly compatible with $\mathtt{T}$ (see Definition 6.5.4 and Remark 6.5.5).

Typically, $\mathcal{C}$ would be a graded-2-category and $\mathsf{E}$ would contain the 3-sesquipolygraph of graded-interchangers, making the data $\mathsf{S} = (\mathsf{R},\mathsf{E})$ a linear Gray rewriting system. To find a basis using the Basis From Convergence Theorem 6.3.18, the following needs to be checked:

($\alpha$) $\mathsf{S}$ *and* $\mathtt{T}$ *present the same underlying (family of) module(s)*: it suffices to show that every $(\mathsf{S}\setminus\mathtt{T})^+$-rewriting step is $\mathtt{T}$-congruent, so that two vectors are $\mathsf{S}$-congruent if and only if they are $\mathtt{T}$-congruent (see Lemma 6.3.9).

($\beta$) $\mathrm{scl}(\mathsf{E})^{\top}$ *is scalar-coherent on* $\mathcal{B}\mathrm{NF}_{\mathtt{T}}$: recall that $\mathcal{B}\mathrm{NF}_{\mathtt{T}}$ denotes the set of monomial $\mathtt{T}^{+}$-normal forms (Definition 6.3.14), and that $\mathrm{scl}(\mathsf{E})^{\top}$ is scalar-coherent on $\mathcal{B}\mathrm{NF}_{\mathtt{T}}$ (Definition 6.1.1 and Remark 6.1.2) if for all $b \in \mathcal{B}\mathrm{NF}_{\mathtt{T}}$, the existence of an $\mathsf{E}$-congruence $b \sim_{\mathsf{E}} \lambda b$ for some scalar $\lambda \in \Bbbk$ implies that $\lambda = 1$.

Showing scalar-coherence of $\mathrm{scl}(\mathsf{E})^{\top}$ on $\mathcal{B}\mathrm{NF}_{\mathtt{T}}$ can be done either using ad-hoc arguments, or by apply higher rewriting modulo theory to $\mathrm{scl}(\mathsf{E})^{\top}$, together with Proposition 6.1.7 (coherence from convergence).

($\gamma$) $\mathtt{T}^{+}$ *is convergent:* we use the Tamed Linear Newmann's Lemma 6.3.24, showing on one hand the $\succ$ is terminating,[20] and on the other that every monomial local $\mathtt{T}^{+}$-branching is $\succ$-tamely $\mathtt{T}^{\mathrm{st}}$-congruent.

The analysis of monomial local $\mathtt{T}^{+}$-branching is the hardest task, which we now discuss. Of course, in the context-agnostic case $\mathtt{T} = \mathsf{S}$, ($\alpha$) is automatic.

### 6.6.2 How to classify monomial local branchings?

Working modulo makes classifying monomial local branchings difficult, as it considerably increases the number of rewriting steps. However, in the context of diagrammatic algebra, the modulo data typically has a topological interpretation, which can be leveraged: e.g. rectilinear isotopies when working modulo interchangers, or planar isotopies when working modulo a pivotal structure. We describe here a general strategy to classify local branchings, to be adapted depending on the example at hand:

(i) *Understand coherence of the modulo.* In other words, provide a topological or combinatorial description of when two monomial 2-cells are projectively $\mathsf{E}$-congruent. this can be done either using ad-hoc arguments, or by applying higher rewriting modulo theory to $\mathrm{scl}(\mathsf{E})$. In practice, ($\beta$) in the previous section comes as a byproduct.

(ii) *Describe naturalities of the modulo.* The modulo typically captures some underlying categorical structure, which should come with natural compatibilities with the rewriting steps of the rewriting system. For instance, a Gray RSM comes with interchange naturalities (see subsection 6.4.6). We call these $\mathsf{E}$-*naturalities.*

[20] While this will not be the case in this thesis, defining a terminating order may be hard in general. For this problem, the method of derivation of Guiraud and Malbos may be useful [89].

(iii) *Characterize rewriting steps modulo.* In other words, provide a topological or combinatorial description of when two rewriting steps are $\mathsf{E}$-congruent. (One needs not show that this characterization captures *all* $\mathsf{E}$-congruences.) To do so, use the coherence statement to express isotopies in terms of $\mathsf{E}$-naturalities. In fact, given the Branchwise Tamed Congruence Lemma 6.2.11, one can instead characterize when two rewriting steps are $\succ$-tamed $\mathsf{T}^+$-congruent.

(iv) *Classify monomial local branchings.* Using the characterization above, provide a list of monomial local $\mathsf{T}^+$-branchings, called *critical branchings*, such that every monomial local $\mathsf{T}^+$-branching rewrites into (a linear combination of) either an independent branching or a branchwise $\succ$-tamely $\mathsf{T}^{\mathrm{st}}$-congruent to a contextualization of a critical branching.

In the last point (iv), three main tools are at play: the Independent Rewriting Lemma 6.5.12, the Branchwise Tamed Congruence Lemma 6.2.11 and the Contextualization Lemma 6.5.7. The results of subsection 6.5.3 also help to deal with independent branchings. Hopefully, this leaves only a few critical branchings for which an explicit computation is needed, and one can conclude that every monomial local $\mathsf{T}^+$-branching is $\succ$-tamely $\mathsf{T}^{\mathrm{st}}$-congruent.

*Remark* 6.6.1. If $\mathsf{E}$ is scalar-free (Remark 6.3.5), the strategy described above greatly simplifies. For instance, this is the case of linear Gray rewriting modulo scalar-free (or non-graded) interchangers, which coincides with rewriting in linear strict 2-categories [3]. However, even in that case it is useful to explicitly consider the interchange law as part of the modulo, given the caveats related to freeness of contextualization, described in subsection ii.2.2.[21]

### 6.6.3 Example: $\mathfrak{gl}_2$-webs

Recall the scalar Gray polygraph $\mathsf{Web}$ from section 5.1. We solve its rewriting theory following the above strategy. In fact, as $\mathsf{Web}$ is scalar, one can use higher (non-linear) rewriting modulo. The approach given is largely overkilled for the example at hand, but helps illustrate the main ideas.

We define the following data:

(a) Let $\mathsf{E} = \mathsf{Web}_{\leq 2}\mathsf{Gray}$ the 3-prepolygraph of interchangers on $\mathsf{Web}_{\leq 2}$ and $\mathsf{R}$ the 3-prepolygraph on $\mathsf{Web}_{\leq 2}$ whose 3-cells are $\mathsf{R}_3 = \mathsf{Web}_3 \setminus \mathsf{E}_3$.

[21] See for instance Footnote 19.

This defines a splitting $\mathsf{Web} = \mathsf{R} \sqcup \mathsf{E}$ and an associated scalar higher RSM denoted $\mathsf{S} = (\mathsf{R}, \mathsf{E})$.

(b) Web does not need a context-dependent termination: we take $\mathsf{S}$ as its own abstract sub-system.

(c) We let $\succ$ be the order induced by counting the number of 2-generators $W_-$ and $W_+$ in a $\mathfrak{gl}_2$-web. Clearly, $\succ$ is $\mathsf{E}$-invariant as the number of 2-generators is preserved by interchangers, and $\succ$ is compatible with $\mathsf{S}$ as each $\mathsf{R}$-rewriting step strictly reduces the number of 2-generators, irrespective of contextualization.

In this setting, $(\alpha)$ and $(\beta)$ are automatically given. The order $\succ$ is terminating by definition, and in this simplified setting, the higher Newmann's lemma (Lemma 6.4.8) will be sufficient. Hence, following $(\gamma)$ we are reduced to show that:

**Lemma 6.6.2** (step (iv))**.** *Every overlapping* $\mathsf{S}$*-branching is, up to branchwise* $\mathsf{E}$*-congruence, a contextualization of one of the two* $(\mathsf{S}, \equiv_{\mathrm{scl}})$*-confluent branchings:*

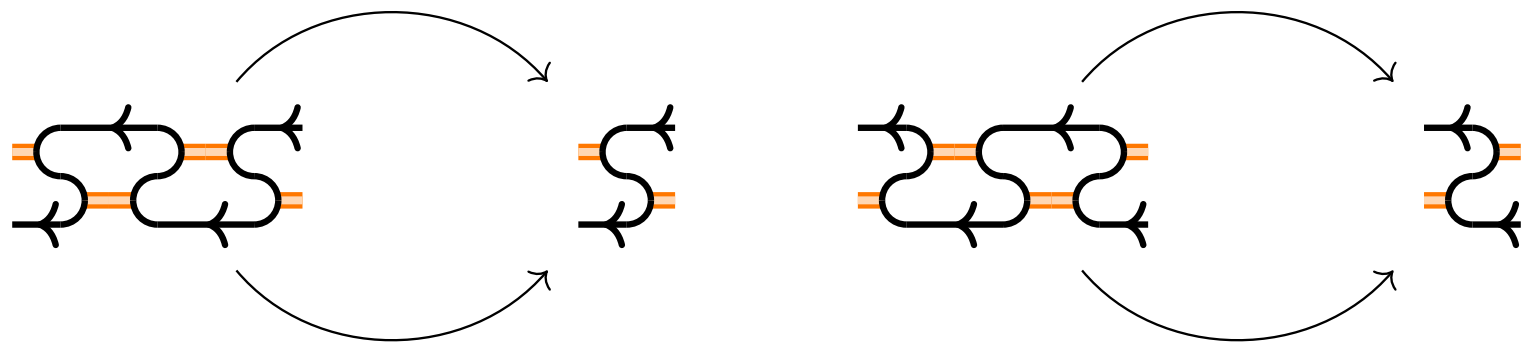

Following the strategy given in subsection 6.6.2, the first step (i) is to understand coherence of the modulo data. For interchangers, the topological interpretation is given by rectilinear isotopies. The following gives coherence:

**Lemma 6.6.3** (step (i))**.** $\mathsf{E}^\top$ *is scalar-coherent.*

*Proof.* $\mathsf{E}_3$ is scalar-free, and hence $\mathsf{E}^\top$ is scalar-coherent. However, we can also argue that $\mathsf{E}_2^*$ does not contain bubbles in the sense of Appendix A and hence, following the coherence theorem for interchangers (Theorem A.2.4), $\mathsf{E}^\top$ is scalar-coherent. The latter argument does not use the fact that $\mathsf{E}_3$ is scalar-free. □

In particular, this gives another proof of $(\beta)$. Rewriting modulo interchangers comes with canonical naturalities, namely the interchange naturalities given in subsection 6.4.6. Step (ii) consists in checking that $\mathsf{R}$-rewriting steps verify these naturalities:

**Lemma 6.6.4** (step (ii))**.** *The 3-cells in* $\mathsf{R}_3$ *verify the interchange naturalities given in subsection 6.4.6.*

*Proof.* Again, as $\mathsf{E}_3$ is scalar-free, this is obvious. However, if we were working with graded interchangers, it would suffice to check that the 3-cells in $\mathsf{R}_3$ preserve the grading. □

In some obvious sense, the data of 2-generators is preserved by rectilinear isotopies. given $\Gamma[r]$ a R-rewriting step with $r \in \mathsf{R}_3$, we call the 2-generators in $s(r)$ the *combinatorial data* associated to $\Gamma[r]$. As this combinatorial data is preserved under E-congruence, it makes sense to talk of *the* combinatorial data associated to an E-congruence class of rewriting steps.

**Lemma 6.6.5** (characterization of rewriting steps; step (iii))**.** *Let* $[f, e, g]$ *be a* S*-local triple. If* $f$ *and* $g$ *have the same combinatorial data, then* $[f, e, g]$ *is* E*-congruent.*

*Proof.* If $f$ and $g$ have the same combinatorial data, there exists a sequence of interchange naturalities $e'$ from $s(f)$ to $s(g)$ that identifies their respective 2-generators. Thanks to the scalar-coherence of $\mathsf{E}^{\top}$ given by step (i), we have $e \equiv_{\text{scl}} e'$. We conclude with step (ii). □

We now have all the ingredients to prove step (iv):

*Proof of Lemma 6.6.2.* Let $[f, e, g] = (\Gamma[p], e, \Gamma[q])$ be a S-local triple, where we have $p, q \in \mathsf{R}_3$. If $s(p)$ and $s(q)$ do not share any 2-generator in $s(f) \sim_{\mathsf{E}} s(g)$, we can find a 2-cell $x$ such that

$$s(f) \sim_{\mathsf{E}} x \sim_{\mathsf{E}} s(g)$$

and an independent R-branching $(f', g')$ on $x$, such that $f'$ (resp. $g'$) has the same combinatorial data as $f$ (resp. $g$). By step (iii), the characterization of rewriting steps, the S-local triple $[f, e, g]$ is branchwise E-congruent to $(f', g')$.

Hence, classifying overlapping S-branchings comes down to classifying how the combinatorial data of $f$ and $g$ can overlap; that is, how their 2-generators can overlap modulo interchangers. This leads to the two critical branchings given in the statement of the lemma. By the characterization of rewriting steps, if the combinatorial data of $f$ and $g$ overlap then $[f, e, g]$ is branchwise E-congruent to a contextualization of one of these two critical branchings. □

Finally, we conclude with the higher Newmann's lemma (Lemma 6.4.8):

**Proposition 6.6.6.** $\mathsf{S}$ *is convergent modulo* $\mathsf{E}$. □

In particular, we have the following corollaries:

**Corollary 6.6.7.** $\mathsf{S}^{\top}$ *is scalar-coherent.*

*Proof.* This follows from the coherence-from-convergence theorem (Proposition 6.1.7), the transitivity of coherence modulo (Lemma 6.1.3) and the scalar-coherence of $\mathsf{E}^{\top}$ (step (i)). □

**Corollary 6.6.8.** *For each 1-sphere* $(x, y)$*, the module* $\mathbf{Web}_{\otimes}(x, y)$ *is free. A basis is given by any choice of family of webs with exactly one web for each isotopy class of embedded 1-manifolds without closed components.*

*Proof.* This is a consequence of the Basis From Convergence Theorem 6.3.18. □

*Proof of Lemma 1.2.4.* Topological considerations show that there exists a $\mathsf{S}$-congruence between $W$ and $W'$ if and only if $c(W)$ and $c(W')$ are isotopic, one removed all there closed components. The result then follow from scalar-coherence of $\mathsf{S}^{\top}$. (One could also appeal to Theorem 6.3.18.) □

### 6.6.4 Addendum: graded deformations

Of course, we greatly overshooted by using rewriting theory to prove the above corollary. However, the rewriting approach has the advantage of being fairly independent of the specific scalars appearing in the presentation of $\mathbf{Web}_{\otimes}$. Indeed, we only used them to check scalar-confluence of the two critical $\mathsf{S}^{+}$-branchings in Lemma 6.6.2.

In particular, it informs us on the existence of *graded deformations* of $\mathbf{Web}_{\otimes}$, in the following sense:

**Definition 6.6.9.** *Let* $G$ *be an abelian group,* $\Bbbk^{\mathrm{ev}}$ *a commutative ring (the "evaluation ring") and* $\mathcal{C}$ *a* $G$*-graded* $\Bbbk$*-linear 2-category. A* graded deformation *of* $\mathcal{C}$ *is the data of a commutative ring* $\Bbbk$*, a ring morphism* $\Bbbk \to \Bbbk^{\mathrm{ev}}$*, a* $\mathbb{Z}$*-bilinear map* $\mu\colon G \times G \to \Bbbk^{\times}$ *and a* $(G, \mu)$*-graded-2-category* $\mathcal{GC}$ *such that*

$$\mathcal{GC} \otimes_{\Bbbk} \Bbbk^{\mathrm{ev}} \simeq \mathcal{C},$$

*where the equivalence is an equivalence of* $G$*-graded* $\Bbbk^{\mathrm{ev}}$*-linear 2-categories.*

To illustrate, we construct a graded deformation $\mathbf{GWeb}_\otimes$ of $\mathbf{Web}_\otimes$; in particular, $\Bbbk^{\mathrm{ev}} = \mathbb{Z}[q, q^{-1}]$. To begin with, we investigate gradings. Looking at the relations, we see that $\mathbf{Web}_\otimes$ admits a $\mathbb{Z}$-grading, unique up to symmetry, setting $\deg(W_-) = 1$ and $\deg(W_+) = -1$. From now on, we view $\mathbf{Web}_\otimes$ as a $\mathbb{Z}$-graded $\mathbb{Z}[q, q^{-1}]$-linear monoidal category.

For $G = \mathbb{Z}$, a $\mathbb{Z}$-bilinear map $\mu\colon G \times G \to \Bbbk^\times$ is determined by the value $X := \mu(1,1)$, with $\mu(n,m) = X^{nm}$, so that we expect $\Bbbk$ to contain such a formal variable $X$.

Looking at the two critical branchings, one sees that adding formal variables $A$ and $B$ to the defining relations of $\mathbf{Web}_\otimes$ as given below preserves the scalar-confluence of the two critical branchings:

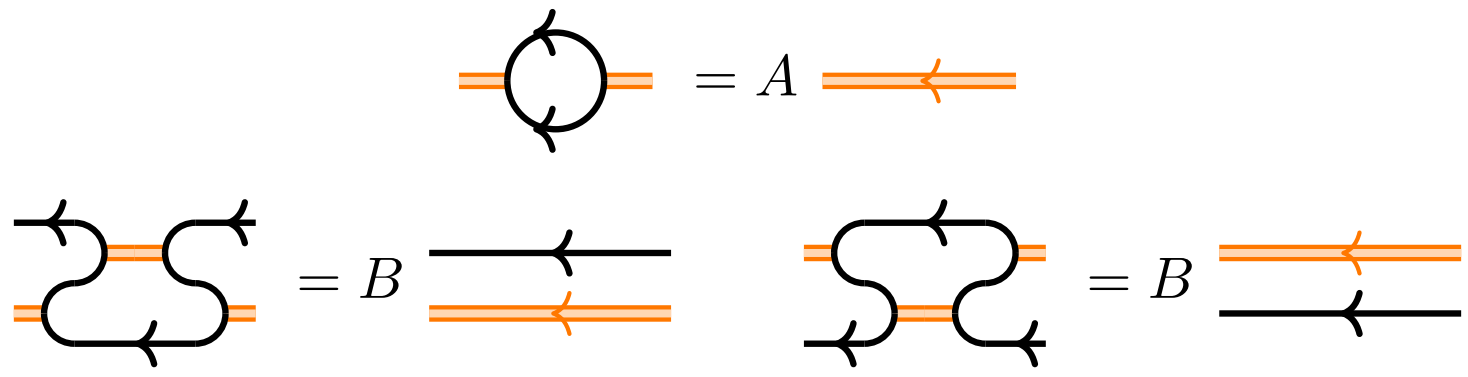

Hence we define $\Bbbk = \mathbb{Z}[X, A, B]$, $\mu\colon G \times G \to \Bbbk^\times$ to be the $\mathbb{Z}$-bilinear map defined by $\mu(1,1) = X$, and $\mathbf{GWeb}_\otimes$ to be the $(\mathbb{Z}, \mu)$-graded-2-category presented by the $(\mathbb{Z}, \mu)$-scalar Gray polygraph GWeb, where $\mathsf{GWeb}_{\leq 2} = \mathsf{Web}_{\leq 2}$ and $\mathsf{GWeb}_3$ consists of graded interchangers and the relations above.

By construction, the analogue of Lemma 6.6.2 holds for $\mathbf{GWeb}_\otimes$, and hence so does the analogue of Corollary 6.6.8. It follows that $\mathbf{GWeb}_\otimes$ is a graded deformation of $\mathbf{Web}_\otimes$, with $\Bbbk \to \Bbbk^{\mathrm{ev}}$ defined as $X \mapsto 1$, $A \mapsto q + q^{-1}$ and $B \mapsto 1$.

Note that if we decategorify $\mathbf{GFoam}_\otimes$ with respect to the full grading $G = \mathbb{Z} \times \mathbb{Z}$, one gets the deformation $X = 1$, $A = u + v^{-1}$ and $B = uv^{-1}$ for some formal variables $u$ and $v$. This suggests that the existence of graded deformations at the decategorified level informs on the possible existence of graded categorifications at the categorified level.

# 7

# A basis for graded $\mathfrak{gl}_2$-foams via rewriting theory

In this chapter, we apply linear Gray rewriting modulo as developed in chapter 6 to graded $\mathfrak{gl}_2$-foams, and show the non-degeneracy theorem (Corollary 1.6.4). The strategy follows the blueprint given in subsection 6.6.2. The modulo data is given by foams isotopies; as such, we begin section 7.1 with a careful study of coherence of foam isotopies. Section 7.2 describes all the necessary working data and deals with termination. Finally, section 7.3 studies tamed congruence to conclude the proof.

Throughout we will assume a good familiarity with chapter 1, and more particularly with the string diagrammatics (section 1.4) and the shading diagrammatics (section 1.5).

## 7.1 Coherence of foam isotopies

Subsection 7.1.1 review some terminologies. Subsection 7.1.2 then gives a polygraphic presentation of foam isotopies and state coherence of foam isotopies (Proposition 7.1.2). Finally, subsection 7.1.3 proves coherence of foam isotopies using the higher Newmann's lemma (Lemma 6.4.8). In that regard, the strategy of proof is closer to the one described in subsection 6.6.3.

### 7.1.1 Terminology

We recall some terminology from chapter 1. A strand labelled by a colour $i$ is called an $i$-*strand*. We similarly define $i$-*dots*. Two colours $i$ and $j$ are said to be *distant* if $|i-j|>1$. Given a colour $i$ and a diagram $\psi$, we say that $i$ *is distant from* $\psi$ if for each $j$-strand (resp. $j$-dot) in $\psi$, we have $|i-j|>1$ (resp. $j\neq i,i+1$).

### 7.1.2 A scalar Gray rewriting system for foam isotopies

**Definition 7.1.1.** *The graded 2-polygraph* $\mathsf{GFoam}_d$ *has the set* $\underline{\Lambda}_d$ *as its underlying set of objects,* $\{W_{i,-},W_{i,+}\mid 1\leq i\leq d-1\}$ *as its generating set of 1-morphisms, and the generators in Definition 1.4.1 as its generating set of graded 2-morphisms.*

We define different 3-sesquipolygraphs with $\mathsf{GFoam}_d$ as their underlying 2-polygraph.

Recall the definition of $\Bbbk$ and $\mu$ from Definition 1.3.2. Let $\mathsf{FGray}$ be the $(\mathbb{Z}^2,\mu)$-scalar 3-sesquipolygraph of graded interchangers on $\mathsf{GFoam}_d$ (see subsection 5.4.2). Let $\mathsf{X}$ and $\mathsf{Z}$ be two 3-sesquipolygraphs with $\mathsf{X}_{\leq 2}=\mathsf{Z}_{\leq 2}=\mathsf{GFoam}_d$, such that $\mathsf{X}_3$ is the scalar globular extension consisting of the following 3-generators:

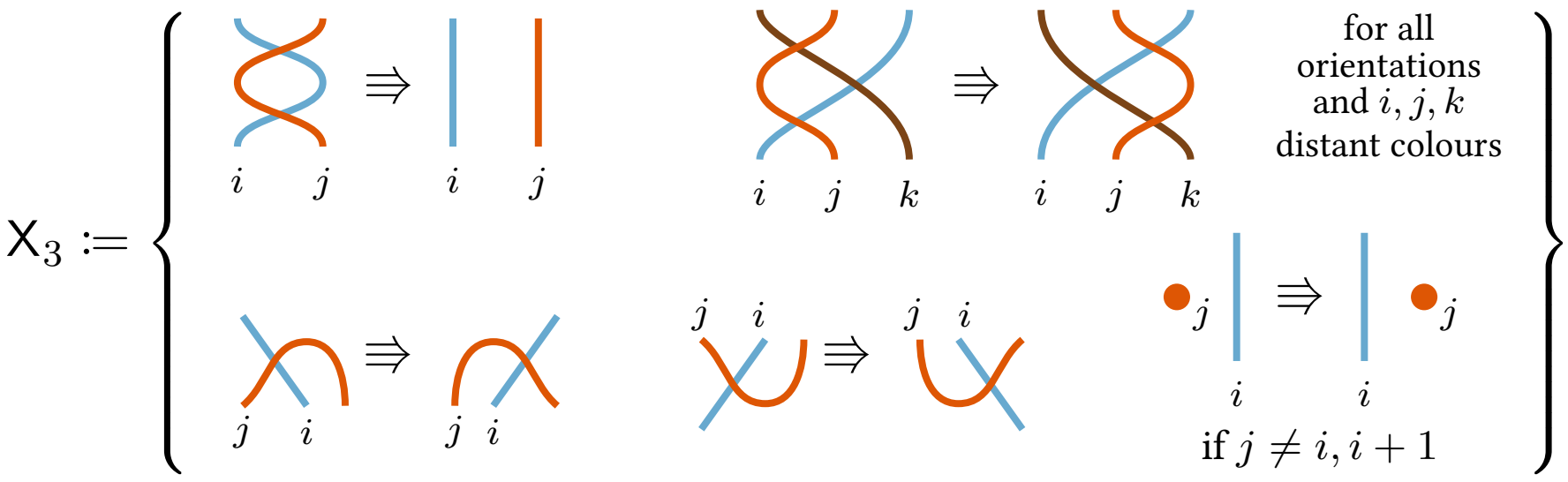

and similarly $\mathsf{Z}_3$ is the scalar globular extension consisting of the following 3-generators:

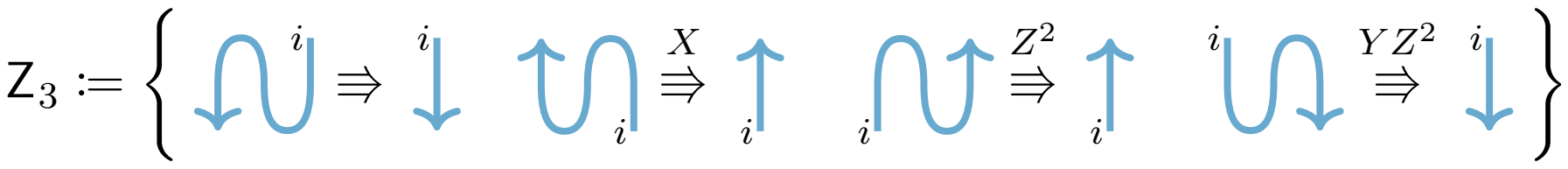

Finally, we define $\mathsf{E}:=\mathsf{FGray}\sqcup\mathsf{X}\sqcup\mathsf{Z}$; in other words, we have $\mathsf{E}_{\leq 2}=\mathsf{GFoam}_d$ and

$$\mathsf{E}_3=\mathsf{FGray}_3\sqcup\mathsf{X}_3\sqcup\mathsf{Z}_3.$$

Obviously, each element of $\mathsf{E}_3$ corresponds to a generating isotopy of foams (see Fig. 1.3): Morse-preserving isotopies for $\mathsf{FGray}_3$, Morse-regular isotopies for $\mathsf{X}_3$ with respectively braided-like relations, pitchfork relations and dot slide, and birth-death isotopies or zigzag relations for $\mathsf{Z}_3$. We use the same terminology to refer to the corresponding element of $\mathsf{E}_3$.

Note that the data of dots and strands in a given foam diagram is preserved under foam isotopies; that is, if $e\colon D_0 \to D_1$ is an isotopy, then there is a canonical bijection between the set of dots and strands of $D_0$, with the set of dots and strands of $D_1$. This makes sense of the following statement, which we show in the next subsection:

**Proposition 7.1.2** (coherence of foam isotopies)**.** *If two parallel morphisms in $\mathsf{E}^\top$ define the same bijection on dots and strands, then they have the same associated scalar.*

Note the following corollary:

**Corollary 7.1.3.** *$\mathsf{E}^\top$ is scalar-coherent on reduced foams.*

*Proof.* A reduced foam does not have any closed strands (Corollary 1.6.8), nor isotopic dots with the same colour. Hence, foam isotopies on reduced foams all induce the same bijection on dots and strands. The statement then follows from Proposition 7.1.2. □

### 7.1.3 Proof of coherence via rewriting theory

Let

$$\overline{\mathsf{E}} = \mathsf{FGray} \sqcup \mathsf{X}.$$

To show Proposition 7.1.2, we use the Gray RSM $(\mathsf{Z}, \overline{\mathsf{E}})$. We follow the general strategy sketched in subsection 6.6.2.

We start by providing the analogous coherence result for the modulo data:

**Lemma 7.1.4** (coherence of foam isotopies, excepted birth-death isotopies)**.** *If two parallel morphisms in $\overline{\mathsf{E}}^\top$ define the same bijection on dots and strands, then they have the same associated scalar.*

*Proof.* Since $\mu$ is symmetric, the scalar of a relation in $\mathsf{FGray}^\top$ only depends on how generators vertically permute (see also Appendix A). In fact, it only depends on how generators *with a non-trivial grading* vertically permute. Relations in $\mathsf{X}^\top$ have trivial associated scalar. While they do not preserve the data of the set of generators, they do preserve the data of the set of *non-trivially graded* generators. That is, if $e\colon D_0 \to D_1$ is a generating relation in $\mathsf{X}^\top$, there is a canonical way to identify the non-trivially graded generators of $D_0$ with those of $D_1$. In other words, the permutation of non-trivially graded generators is a well-defined data associated to any relation in $\overline{\mathsf{E}}^\top$, and this data determines the associated scalar. Finally, we check that this permutation data is itself determined by the bijection on dots and strands. □

The modulo $\overline{\mathsf{E}}^\top$ can be understood as encompassing a braided-like structure. Indeed, $\overline{\mathsf{E}}^\top$ is equivalently generated by the braided-like relations and the following 3-morphisms:

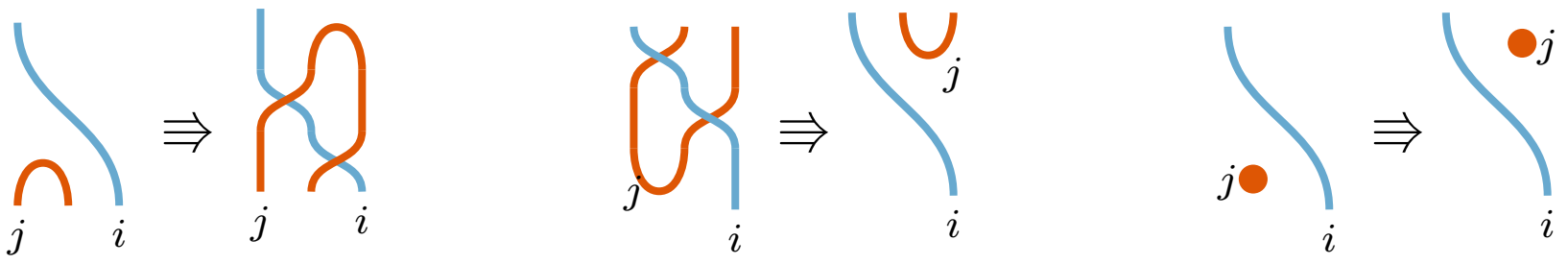

if $|i-j|>1$ if $|i-j|>1$ if $j \neq i, i+1$

which capture the naturality of the braiding with respect to the cap, the cup and the dot. We say “braided-like” to emphasize the restriction on the labels. This is of course no coincidence, as the crossings are interchangers witnessing the interchange law on webs (see Remark 1.3.4).

Recall that any Gray RSM has canonical interchange naturalities (subsection 6.4.6). For $(\mathsf{Z}, \overline{\mathsf{E}})$, we further have the following braided-like naturalities:

**Lemma 7.1.5** (braided-like $\overline{\mathsf{E}}$-naturalities). *For any relation $f\colon \psi_0 \to \psi_1$ in $\mathsf{Z}$ and $i$ any colour distant from $\psi_0$ (see subsection 7.1.1), we have the following $\overline{\mathsf{E}}$-congruences:*

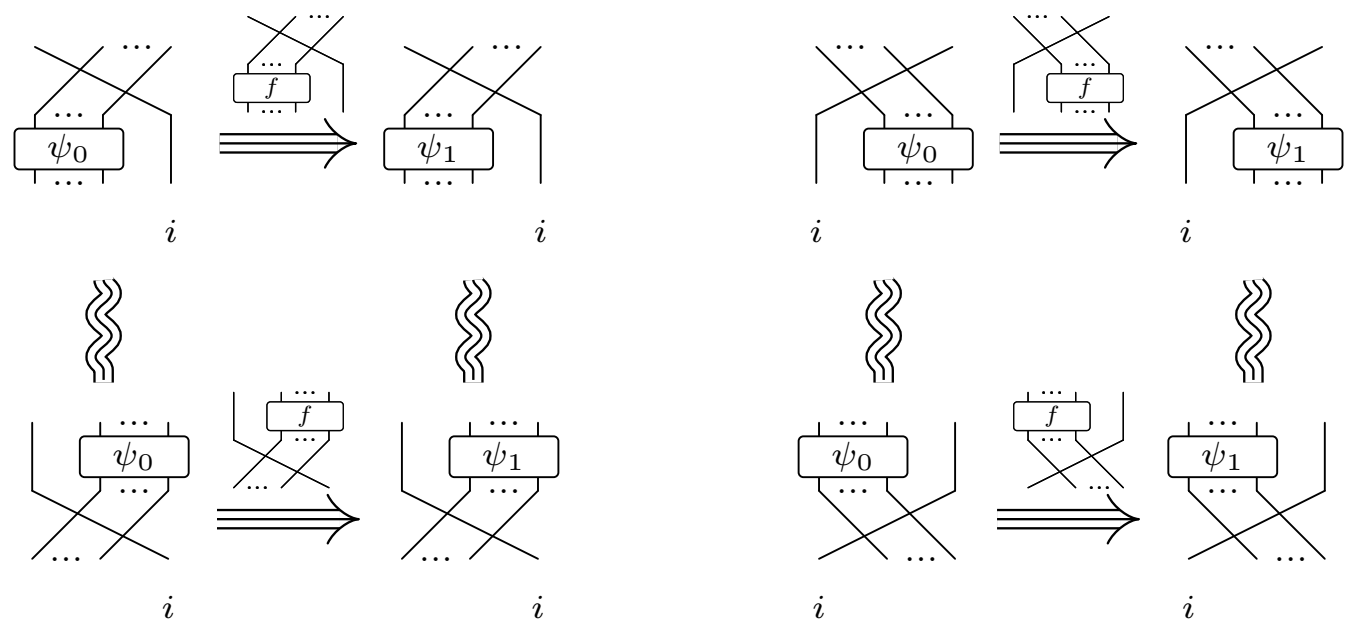

□

Our next step is to upgrade our knowledge of $\overline{\mathsf{E}}$-naturalities to a general characterization of $\overline{\mathsf{E}}$-congruence classes, leveraging our understanding of coherence in $\overline{\mathsf{E}}$:

**Lemma 7.1.6** (characterization of $\overline{\mathsf{E}}$-congruence classes for $\mathsf{Z}$). *If $(f,g)$ is a local $\mathsf{Z}^{+}$-branching such that $f$ and $g$ are of the same type with identical cup and cap, then $(f,g)$ is $\overline{\mathsf{E}}$-scalar-congruent.*

*Proof.* Let $[f,e,g]$ be a $\mathsf{Z}$-local triple with $f$ and $g$ as in the statement of the lemma. Thanks to coherence of the modulo (Lemma 7.1.4), we can choose the isotopy $e\colon s(f)\to s(g)$ such that each step either does not overlap $s(f)$, or consists of a braided-like $\overline{\mathsf{E}}$-naturality (Lemma 7.1.5). This concludes. □

With this characterization at hand, we can classify local $\mathsf{Z}$-branchings and reduce scalar-confluence to four explicit branchings, the critical branchings:

**Lemma 7.1.7.** *Every local $\mathsf{Z}$-branching is, up to branchwise $\overline{\mathsf{E}}$-congruence, either an independent branching or a contextualization of one the following scalar-confluent branchings:*

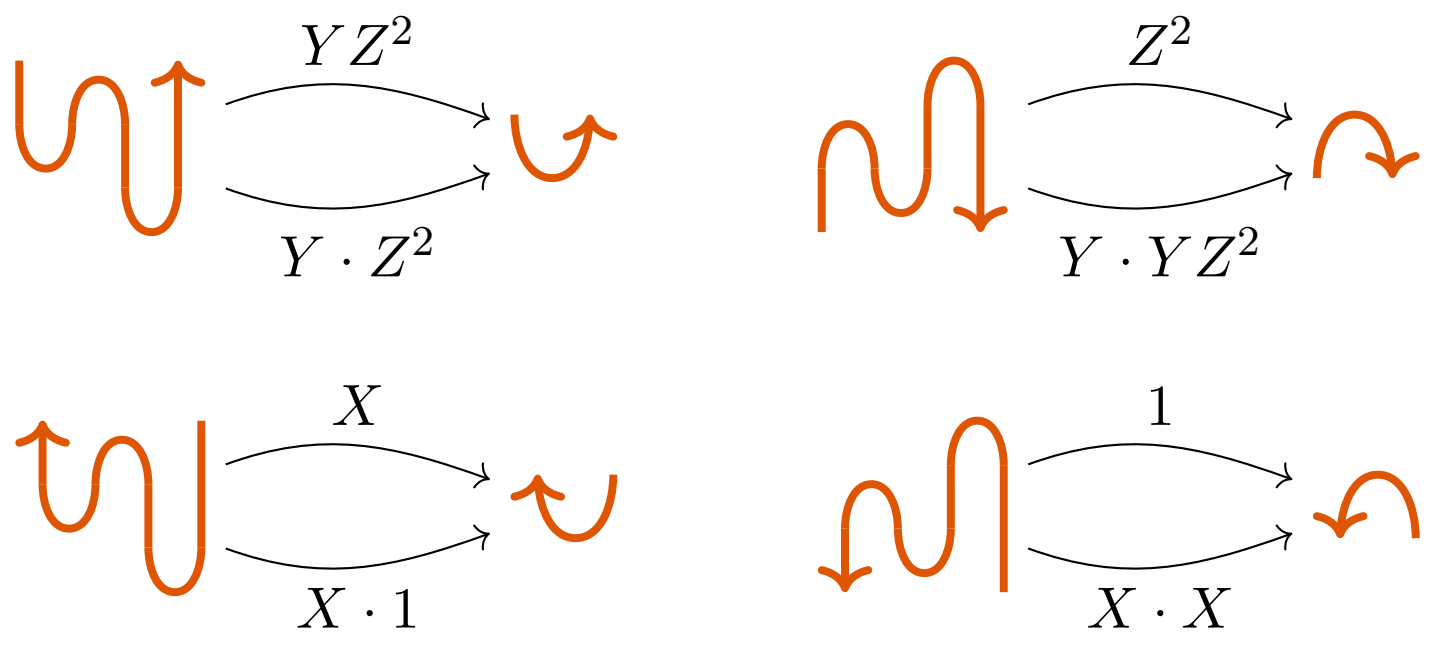

*Proof.* Let $[f,e,g]=[\Gamma_f[r_f],e,\Gamma_f[r_g]]$ be a $\mathsf{Z}$-local triple, with $r_f,r_g\in\mathsf{Z}_3$. We wish to simplify $[f,e,g]$ by replacing it with a branchwise $\overline{\mathsf{E}}$-congruent branching using the characterization given in Lemma 7.1.6. Given a cup $\cup$ and a cap $\cap$ belonging to the same strand $S$, we refer to a "strand between the $\cup$ and $\cap$" to mean a strand crossing $S$ on the piece of strand connecting $\cup$ to $\cap$.

Consider the cup and cap of $s(r_g)$, respectively denoted $\cup_g$ and $\cap_g$, as sitting in $s(f)$. They necessarily belong to the same strand $S$, possibly with strands between $\cup_g$ and $\cap_g$. Using isotopies in $\overline{\mathsf{E}}$, we can slide these strands away, so that no strand lies between $\cup_g$ and $\cap_g$. Moreover, we can do so without adding strands between $\cup_f$ and $\cap_f$, the cup and cap of $s(r_f)$ (viewed as sitting in $s(f)$). Thanks to the characterization given in Lemma 7.1.6, this does not change the branchwise $\overline{\mathsf{E}}$-congruence class of $[f,e,g]$.

If the cups and caps of $s(r_f)$ and $s(r_g)$ are disjoint, we can further use interchanges to move the cup and cap of $s(r_g)$ one below the other, again without affecting $s(r_f)$. This shows that $[f,e,g]$ is branchwise $\overline{\mathsf{E}}$-congruent to an independent $\mathsf{Z}$-branching. If $s(r_f)$ and $s(r_g)$ have precisely one cap or one cup in common, we find the four critical branchings given in the lemma. □

Moreover:

**Lemma 7.1.8.** $\mathsf{Z}$ *terminates modulo* $\overline{\mathsf{E}}$.

*Proof.* $\mathsf{Z}$ strictly diminishes the number of caps and cups, which is kept constant by $\overline{\mathsf{E}}$. □

We now have all the ingredients to show Proposition 7.1.2:

*Proof of Proposition 7.1.2.* Using the higher Newmann's lemma (Lemma 6.4.8), Lemma 7.1.7 and termination (Lemma 7.1.8), we conclude that $\mathsf{Z}$ is convergent $\overline{\mathsf{E}}$. Hence, $\mathsf{Z}^\top$ is scalar-coherent modulo $\overline{\mathsf{E}}^\top$ (coherence modulo from convergence modulo; Proposition 6.1.7).

To show the proposition, it suffices to show that if a loop $e$ in $\mathsf{E}^\top$ defines the identity bijection of dots and strands, then its associated scalar is one. Since $\mathsf{Z}^\top$ is scalar-coherent modulo $\overline{\mathsf{E}}^\top$, we can write $e$ as $e = g \circ h \circ g^{-1}$ with $h$ in $\overline{\mathsf{E}}^\top$. We conclude using coherence of $\overline{\mathsf{E}}^\top$ (Lemma 7.1.4). □

## 7.2 A convergent presentation of graded $\mathfrak{gl}_2$-foams

This section sets up the "working data" for the rewriting theory of graded $\mathfrak{gl}_2$-foams. This is the data (a), (b) and (c) as described in subsection 6.6.1. Subsection 7.2.1 describes a suitable polygraphic presentation of graded $\mathfrak{gl}_2$-foams. In subsection 7.2.2, we define a linear sub-system $\mathtt{T}$ and a strongly compatible and terminating order $\succ$. Finally, subsection 7.2.3 states the basis theorem and review what remain to be shown; this is $(\alpha)$ and $(\gamma)$ in the general strategy described in subsection 6.6.1. The rest of the proof, and in particular the tamed congruence of monomial branching (Proposition 7.2.6), is the subject of section 7.3, the last section of the chapter.

### 7.2.1 A linear Gray RSM for graded $\mathfrak{gl}_2$-foams

Let $\mathsf{R}$ be the linear 3-sesquipolygraph with $\mathsf{R}_{\leq 2} = \mathsf{GFoam}_d$ and $\mathsf{R}_3$ the linear globular extension of $\langle \mathsf{GFoam}_d \rangle_{\Bbbk}$ defined in Fig. 7.1. Obviously, each element of $\mathsf{R}_3$ corresponds to a generating relation in $\mathbf{GFoam}_d$ (see Fig. 1.3): dot annihilation (dd), dot migration (dm), evaluation of counter-clockwise ($\mathsf{bb}_{\circlearrowleft}$) or clockwise ($\mathsf{bb}_{\circlearrowright}$) bubbles, neck-cutting relation (nc), or squeezing relation (sq). As for $\mathsf{E}_3$, we use the same terminology to refer to the correspond element of $\mathsf{R}_3$.

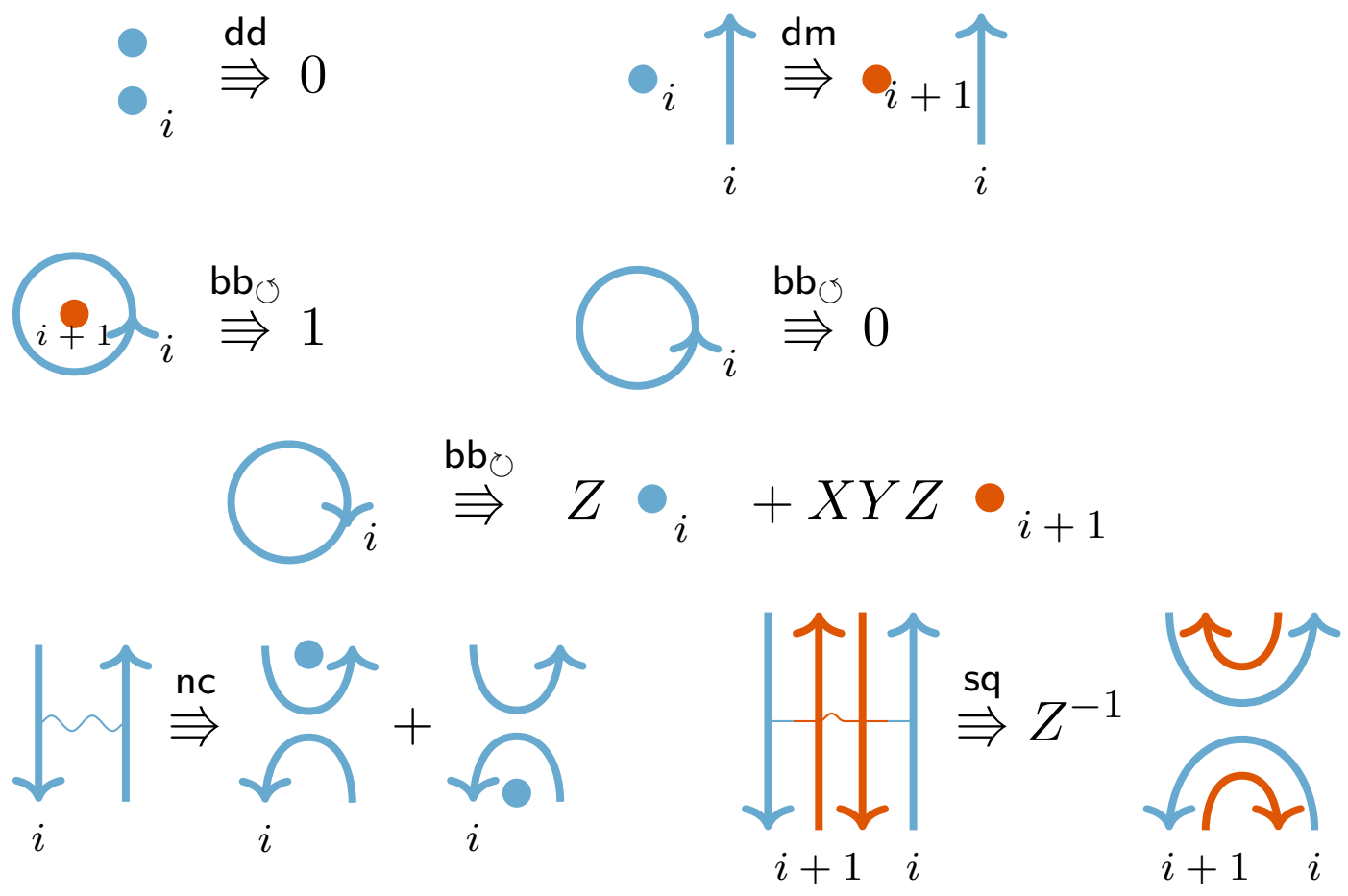

**Fig. 7.1** 3-cells in $\mathsf{R}_3$. In the last two cases, the wiggly lines are only visual aids, and are not part of the data.

**Lemma 7.2.1.** *The linear Gray RSM* $\mathsf{S} := (\mathsf{R}, \mathsf{E})$ *presents the graded-2-category* $\mathbf{GFoam}_d$. ☐

We further define the linear sub-3-sesquipolygraph

$$\mathsf{B}_3 := \{\mathsf{dd}, \mathsf{dm}, \mathsf{bb}_{\circlearrowleft}, \mathsf{bb}_{\circlearrowright}\}.$$

We write $\mathtt{R} := \mathrm{Cont}(\mathsf{R})$, $\mathtt{E} := \mathrm{Cont}(\mathsf{E})$ and $\mathtt{B} := \mathrm{Cont}(\mathsf{B})$ the associated linear RSs, and $\mathtt{S} := (\mathtt{R}, \mathtt{E})$ the associated linear RSM.

**Lemma 7.2.2.** *The following are* $\mathtt{S}$*-rewriting steps:*

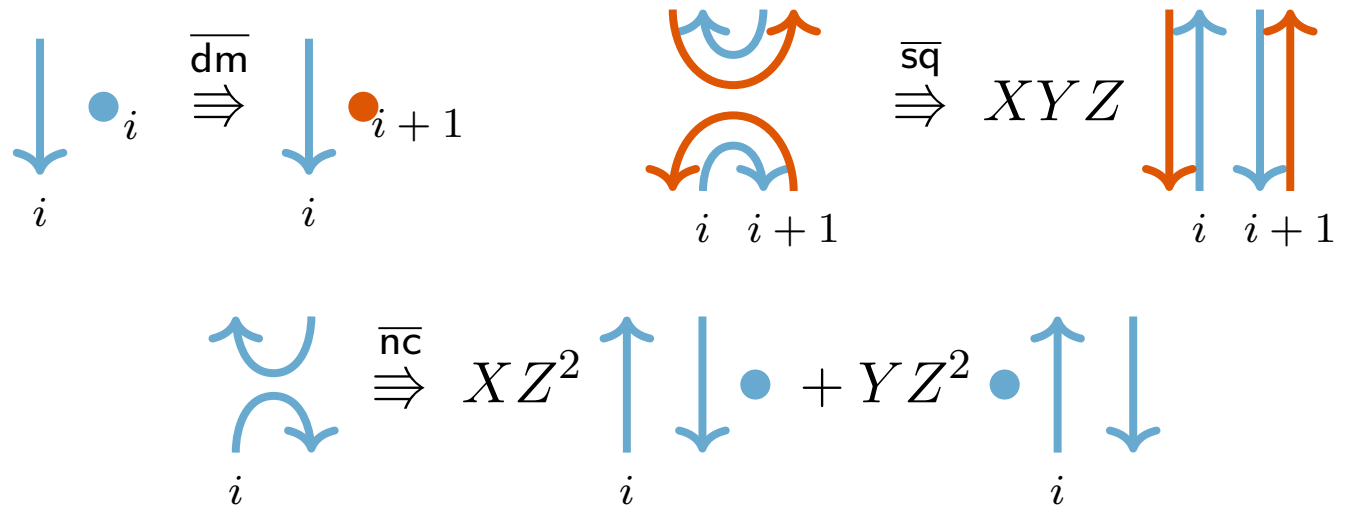

### 7.2.2 A context-dependent linear sub-system

The linear Gray RSM S is not terminating, and hence not suited for a reduction algorithm. For instance, we have the following infinite sequence:

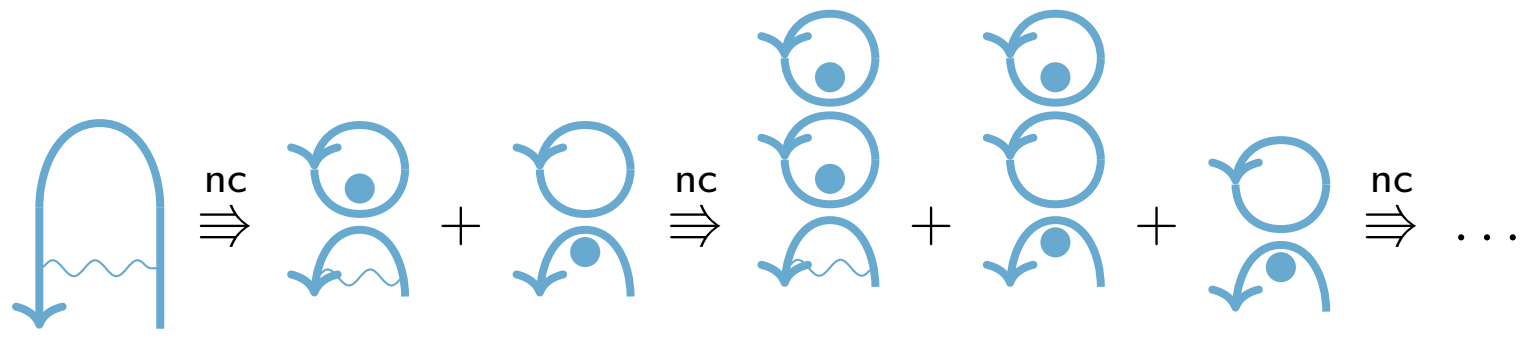

To avoid this obstruction to termination, we define a linear sub-system of S that prevents one from applying a neck-cutting on two pieces of the same strand. Taking into account the analogue obstruction for sq leads to the following definition:

**Definition 7.2.3.** *The (family of) linear RSM(s)* T *is the linear sub-system* T $\subset$ S *such that for each* $r \in \mathsf{R}_3$, $\Gamma[r] \in$ T *if and only if either* $r \in$ B, *or* $r =$ nc *(resp.* sq*) of colour* $i$ *(resp.* $(i, i+1)$*), and the two pieces of* $i$*-strands belong to distinct strands in* $\Gamma[r]$.

The reduction algorithm defined by $\mathtt{T}^+$ is precisely the eponymous algorithm defined in chapter 1. Hence, a monomial $\mathtt{T}^+$-normal form is the same as a reduced foam (subsection 1.6.1 and Proposition 1.6.5), and we shall use the later terminology. Note that T is *not* of form Cont(T) for some sub-Gray RSM T $\subset$ S.

Given a foam $F$, we associate to $F$ the following data:

$$\#\mathrm{sh}_i(F) := \text{number of } i\text{-shadings in the shading diagrammatics of } F,$$
$$\#\mathrm{cl}_i(F) := \text{number of closed } i\text{-strands in the string diagrammatics of } F,$$
$$\#\mathrm{d}_i(F) := \text{number of } i\text{-dots in } F.$$

Each of these data defines a partial order on foams. Let $\succ$ be the lexicographic order induced by these partial orders:

$$\succ := (\#\mathrm{sh}_1, \ldots, \#\mathrm{sh}_{d-1}, \#\mathrm{cl}_1, \ldots, \#\mathrm{cl}_{d-1}, \#\mathrm{d}_1, \ldots, \#\mathrm{d}_d).$$

By definition, $\succ$ is well-founded.

**Proposition 7.2.4.** *The preorder* $\succ$ *is strongly compatible with* T*. In particular,* $\mathtt{T}^+$ *terminates.* □

### 7.2.3 A basis for graded $\mathfrak{gl}_2$-foams

We aim to show the following basis theorem using $\mathtt{T}^+$ and Theorem 6.3.18:

**Theorem 7.2.5.** *Let $W, W' \colon \mu \to \lambda$ be two parallel webs. If $B$ is a set of unique isotopy representative for reduced foams, then $B$ is a basis of the free $\Bbbk$-module* $\mathrm{Hom}_{\mathbf{GFoam}_d}(W, W')$.

We have already seen that $\mathsf{E}$ is coherent on reduced foams (Corollary 7.1.3). Thanks to the Basis From Convergence Theorem 6.3.18 and the Tamed Linear Newmann's Lemma 6.3.24, to prove Theorem 7.2.5 it remains to check that, on the one hand:

**Proposition 7.2.6.** *Every monomial local $\mathtt{T}^+$-branching is $\succ$-tamely $\mathtt{T}$-congruent.*

And on the other hand:

**Lemma 7.2.7.** *The linear RSM $\mathtt{T}$ and $\mathtt{S}$ present the same underlying module.*

Both statements are shown in the next and last section.

## 7.3 Tamed congruence modulo isotopies

We denote $\mathtt{A} := \mathtt{T} \setminus \mathtt{B}$; that is, $\mathtt{A}$ consists of the rewriting step of type $\mathsf{nc}$ and $\mathsf{sq}$ which are in $\mathtt{T}$.

This section is devoted to the proof of Proposition 7.2.6. The study of monomial local $\mathtt{T}^+$-branchings is roughly divided in three parts: both branches are of type $\mathsf{B}$ (Proposition 7.3.3), one branch is of type $\mathsf{B}$ and the other is of type $\mathsf{A}$ (Lemma 7.3.12), and both branches are of type $\mathsf{A}$ (Lemma 7.3.13). The general strategy is, with some variations, the one described in subsection 6.6.2: give naturalities of the modulo, characterize rewriting steps modulo, and finally enumerate monomial local branchings.

### 7.3.1 Foam isotopy naturalities

**Lemma 7.3.1** (braided-like $\mathsf{E}$-naturalities). *The relations in* $\mathsf{R}_3 \sqcup \{\overline{\mathsf{dm}}, \overline{\mathsf{nc}}, \overline{\mathsf{sq}}\}$ *satisfy the same braided-like naturalities as pictured in Lemma 7.1.5.*

*Proof.* This follows from the fact that if $i$ is distant from $\psi_0$, then $i$ is distant from $\psi_1$ (using the same notations as in Lemma 7.1.5). This is straightforward in most cases; we only detail type dm. Let $j$ be the colour of the migrating dot in $\psi_0$. By assumption, the associated $j$-strand crosses an $i$-strand, so that we have $|i-j|>1$. In particular $j+1 \neq i, i+1$ and the $(j+1)$-dot can still slide through the $i$-strand once it has migrated. □

For $A\colon \psi_0 \to \psi_1$ a 3-cell in $\mathsf{R}_3 \sqcup \{\overline{\mathsf{dm}}, \overline{\mathsf{nc}}, \overline{\mathsf{sq}}\}$, we write $\rotatebox[origin=c]{180}{$A$}\colon \rotatebox[origin=c]{180}{$\psi_0$} \to \rotatebox[origin=c]{180}{$\psi_1$}$ the 3-cell obtained by rotating each diagram by a half-turn. Most 3-cells rotate to themselves, except for dm (resp. $\overline{\mathsf{dm}}$), which rotate to $\overline{\mathsf{dm}}$ (resp. dm).

**Lemma 7.3.2** (pivotal E-naturalities)**.** *Let $A\colon \psi_0 \to \psi_1$ be a 3-cell in* $\mathsf{R}_3 \sqcup \{\overline{\mathsf{dm}}, \overline{\mathsf{nc}}, \overline{\mathsf{sq}}\}$ *and $i$ be any color distant from $\psi_0$. Then the following squares commute:*

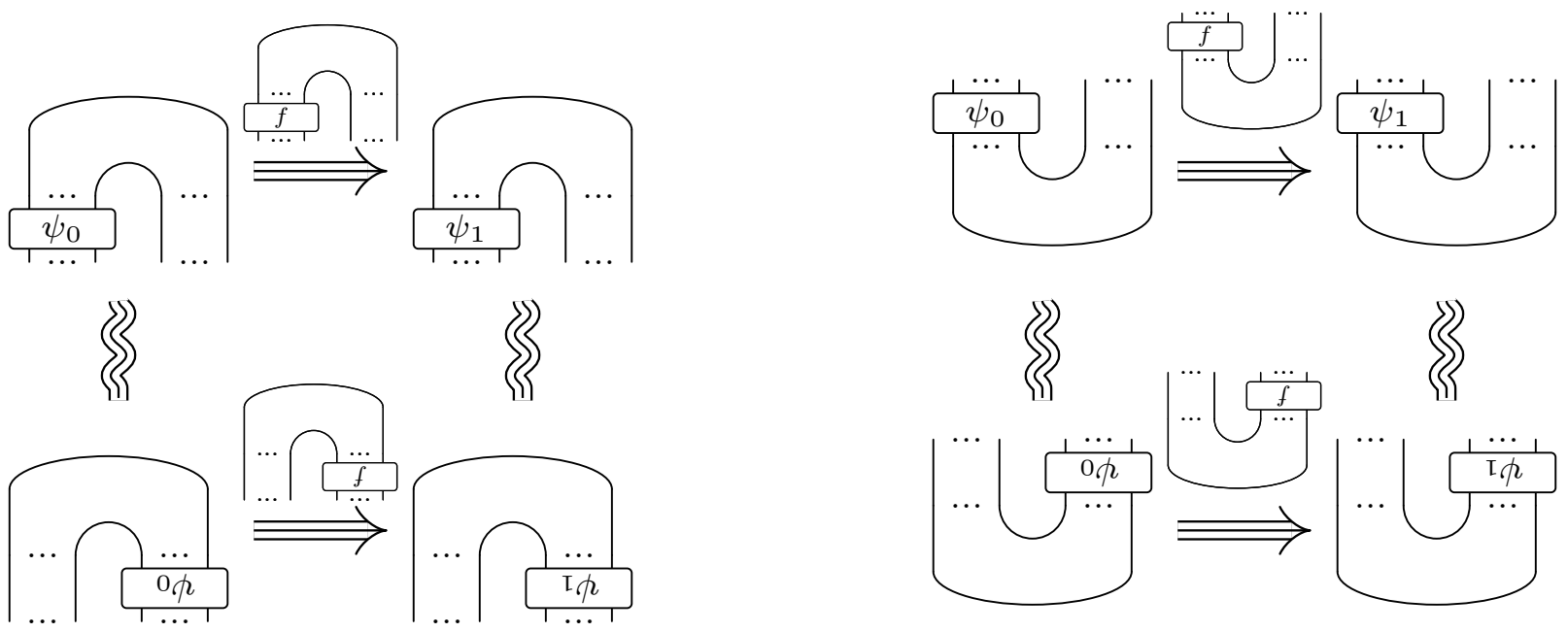

*Proof.* The statement is trivial for types dd, $\mathsf{bb}_{\circlearrowleft}$ and $\mathsf{bb}_{\circlearrowright}$, and follows from graded interchange for types dm. For types nc and sq, it comes down to the following computations:

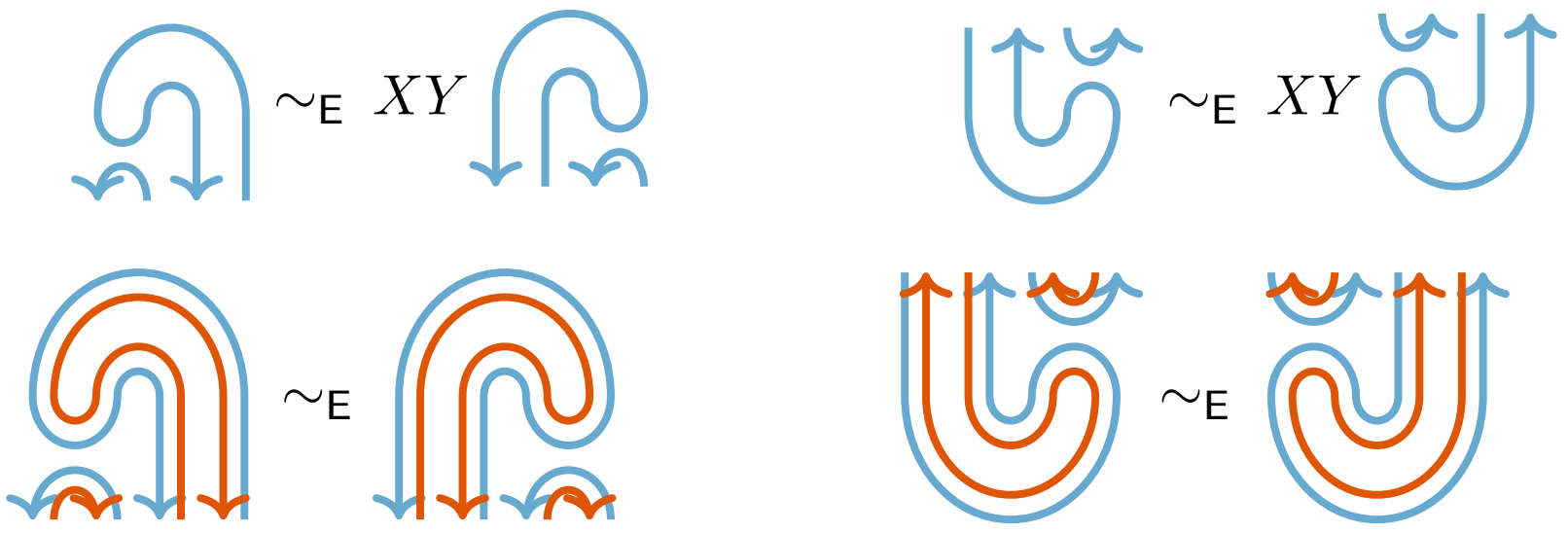

This also implies the lemma for the overlined types $\overline{\mathsf{dm}}$, $\overline{\mathsf{nc}}$ and $\overline{\mathsf{sq}}$, using coherence of E. □

### 7.3.2 Confluence of monomial local $\mathsf{B}^+$-branchings

We begin our study of confluence with positive branchings in the context-agnostic linear sub-system $\mathtt{B} \subset \mathtt{T}$, which derives from a linear Gray RSM: $\mathtt{B} = \mathrm{Cont}(\mathsf{B})$. In fact, we show that in that case, confluences can also be taken in $\mathsf{B}$:

**Proposition 7.3.3.** *Every monomial local $\mathsf{B}^+$-branching is $\succ$-tamely $\mathsf{B}^{\mathrm{st}}$-congruent.*

Note that since $\succ$ is strongly compatible with $\mathtt{T}$ (Proposition 7.2.4), it is also strongly compatible with $\mathsf{B}$. As a preliminary step, we give a topological characterization of $\mathsf{E}$-congruence classes:

**Lemma 7.3.4** (characterization of $\mathsf{E}$-congruence classes for $\mathsf{B}$). *Let $(f,g)$ be a monomial local $\mathsf{B}^+$-branching, with $f$ and $g$ of the same type. The following holds:*

*(a)* type dd: *$(f,g)$ is $\mathsf{E}$-congruent;*

*(b)* type dm: *If $f$ and $g$ have isotopic $i$-dot and $i$-strand, then $(f,g)$ is $\mathsf{E}$-congruent;*

*(c)* type $\mathsf{bb}_{\circlearrowleft}$ and $\mathsf{bb}_{\circlearrowright}$: *If $f$ and $g$ have isotopic $i$-strand, then $(f,g)$ is $\mathsf{E}$-congruent.*

*We call the data associated to each type its* combinatorial data, *and say that this combinatorial data* characterizes *the $\mathsf{E}$-congruence class of the type.*

*Proof.* Denote $\psi$ the local picture of the rewriting step $f$; that is, $\psi = s(r)$ for $r \in \mathtt{B}$ and $f = \Gamma[r]$ for some context $\Gamma$. In each type, we use coherence of $\mathsf{E}$ (Proposition 7.1.2) to present the isotopy $e$ as a composition of $\mathsf{E}$-naturalities, as described in Lemma 7.3.1 and Lemma 7.3.2:

(a) Trivial, since both $f$ and $g$ rewrites to zero.

(b) Through the isotopy $e$, the $i$-dot starts and ends next to the $i$-strand. Hence, we can choose $e$ such that the $i$-dot always remains next to the $i$-strands. In that case, the only isotopies that overlap $\psi$ are braided-like naturalities.

(c) Given that the bubble $\psi$ starts and ends without any strand crossing it, we can choose the isotopy $e$ such $D$ always crosses strands "at once". In that case, the only isotopies that overlap $\psi$ are braided-like isotopies. □

As $\mathsf{B}$ is context-agnostic, the case of independent $\mathsf{B}^+$-branchings comes for free (Lemma 6.5.10):

**Lemma 7.3.5.** *Each independent $\mathsf{B}^+$-branchings is $\succ$-tamely $\mathsf{B}$-congruent.* □

We can now prove the proposition:

*Proof of Proposition 7.3.3.* We use the Branchwise Tamed Congruence Lemma 6.2.11 and the fact that independent $\mathsf{B}^+$-branchings are $\succ$-tamely $\mathsf{B}$-congruent (Lemma 7.3.5) without further mention.

Let $[f,e,g]$ be a monomial $\mathsf{B}$-local triple. If $f$ and $g$ are of type $\mathsf{dd}$, $\mathsf{bb}_{\circlearrowleft}$ or $\mathsf{bb}_{\circlearrowright}$ (not necessarily both of the same type), then either they have the same combinatorial data and hence are $\mathsf{E}$-congruent, or their combinatorial data are disjoint. We can isotope them away from each other, and $[f,e,g]$ is branchwise $\mathsf{E}$-congruent to an independent branching. A similar argument applies if $f$ is of type $\mathsf{dm}$ and $g$ is of type $\mathsf{bb}_{\circlearrowright}$.

Assume $f$ is of type $\mathsf{dm}$ and $g$ is of type $\mathsf{dd}$. If their combinatorial data are disjoint, $[f,e,g]$ is branchwise $\mathsf{E}$-congruent to an independent branching. Otherwise, they share the data of an $i$-dot. In that case, $[f,e,g]$ is $\mathsf{E}$-congruent to a contextualization of the following local $\mathsf{B}^+$-branching, shown to be $\mathsf{B}^+$-confluent:

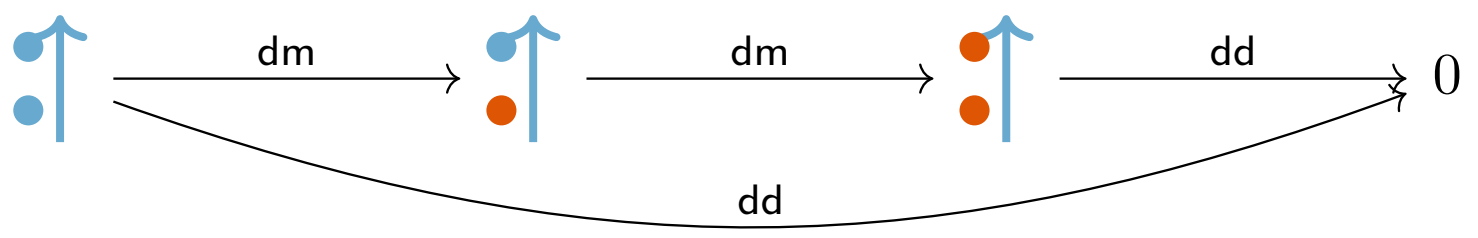

It follows from the Contextualization Lemma 6.5.7 that $[f,e,g]$ is $\succ$-tamely $\mathsf{B}$-confluent.

Assume finally that $f$ is of type $\mathsf{dm}$ and $g$ is of type $\mathsf{bb}_{\circlearrowleft}$. A similar reasoning reduces the statement to the following $\mathsf{B}^+$-confluence:

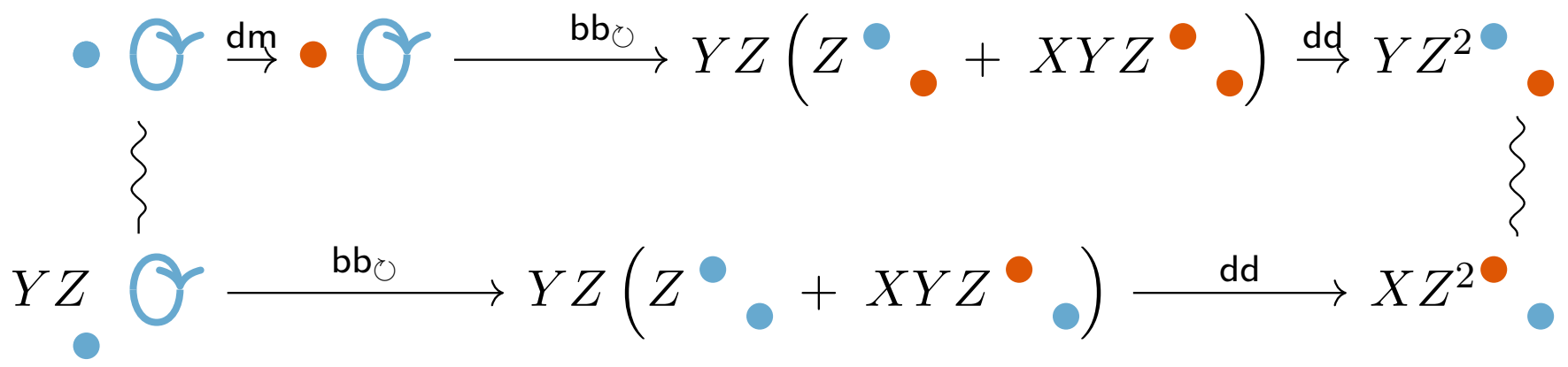

This concludes. □

Finally, it follows from the Tamed Linear Newmann's Lemma 6.3.24 that:

**Corollary 7.3.6.** $\mathsf{B}^+$ *is convergent.* □

A *bubble* is an endomorphism of an identity 1-morphism, possibly viewed inside a bigger diagram. Recall from Definition 1.6.6 the notion of generalized bubble, and from Lemma 1.6.7 the evaluation of a generalized bubble. Looking at the proof of this lemma, it is clear that this evaluation only involves $\mathsf{B}^+$-rewriting steps. Proposition 7.3.3 shows that this evaluation is uniquely defined up to $\mathsf{E}$-congruence, so that we can speak of the $\mathsf{B}^+$-rewriting sequence evaluating a generalized bubble. In particular, we can do so for bubbles:

**Definition 7.3.7.** *For each bubble $\phi$, its* bubble evaluation *is the* $\mathsf{B}^+$*-rewriting sequence, well-defined up to* $\mathsf{E}$*-congruence, which rewrites $\phi$ into a sum of dots:*

$$\mathsf{bb}^*: \quad \textcircled{\phi} \xrightarrow{*}_{\mathsf{B}} \sum \textit{dots}.$$

*A* $\mathsf{B}$*-rewriting sequence defined as a contextualized* $\mathsf{bb}^*$ *is said to be* of type $\mathsf{bb}^*$.

## 7.3.3 Characterizing neck-cutting and squeezing relation up to B-confluence

In order to study rewriting steps of type $\mathsf{nc}$ and $\mathsf{sq}$, we would like to give a topological characterization of their $\mathsf{E}$-congruence classes akin to the one given for $\mathsf{B}$ in Lemma 7.3.4. In fact, it will be easier to characterize their branchwise $\mathsf{B}^+$-confluence classes.

**Proposition 7.3.8** (characterization of branchwise $\mathsf{B}$-confluence classes for $\mathsf{nc}$ and $\mathsf{sq}$)**.** *Let $(f,g)$ be a monomial local* $\mathsf{S}^+$*-branching with $f$ and $g$ of type* $\mathsf{nc}$ *and label $i$ (resp. of type* $\mathsf{sq}$ *and label $(i, i+1)$). If $f$ and $g$ apply to the same $i$-strand(s), then $(f,g)$ is* $\mathsf{B}^+$*-confluent.*

We call the data of the $i$-strand(s) the *combinatorial data* of type $\mathsf{nc}$ (resp. $\mathsf{sq}$), and say that it *characterizes* their $\mathsf{B}^+$-confluence class. Before proving Proposition 7.3.8, we show the following elementary $\mathsf{B}^+$-confluences:

**Lemma 7.3.9** (spatial-like $\mathsf{B}^+$-confluence)**.** *The following branchings are* $\mathsf{B}^+$*-confluent, where $\phi$ is an arbitrary bubble and the dotted wiggly line denotes either a neck-cutting or a squeezing relation:*

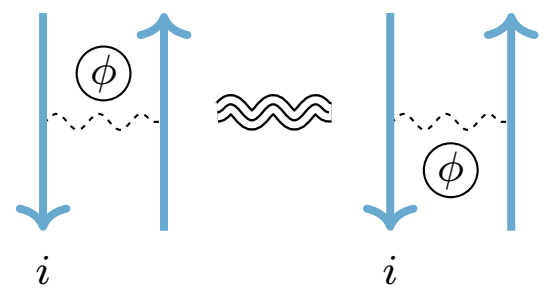

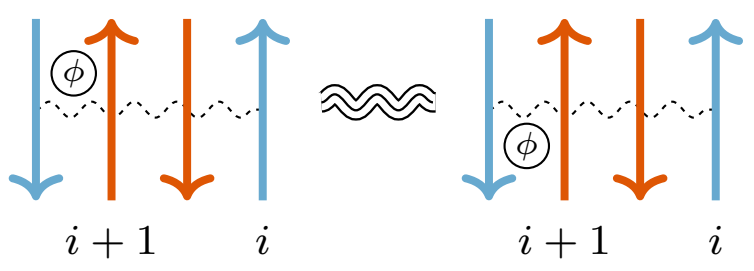

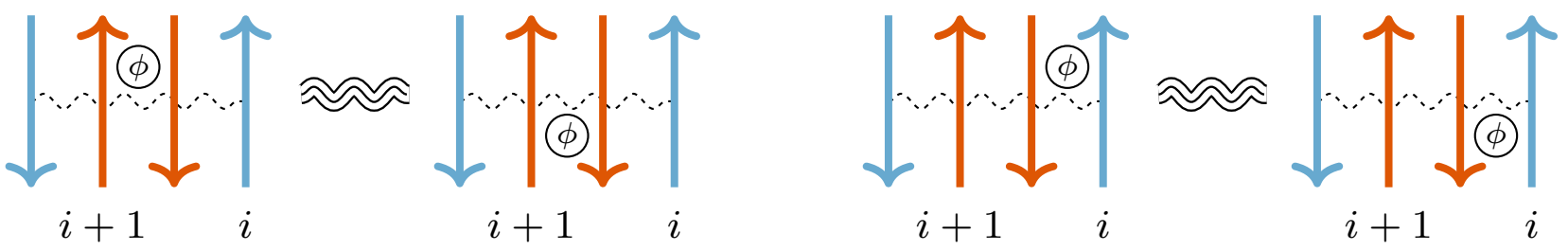

*Proof.* We shall see that in each case, it suffices to evaluate the bubble $\phi$ (see Definition 7.3.7) and apply some additional dot migrations to achieve $\mathsf{B}^+$-confluence. Denote

$$\mathsf{bb}^*\colon \quad \textcircled{\scriptsize $\phi$} \xrightarrow{*}_{\mathsf{B}} \sum_{\delta} \textcircled{\scriptsize $\delta$}$$

the bubble evaluation of $\phi$. Denote by $\delta$ a generic union of dots appearing in this bubble evaluation.

Consider the first branching. We compare

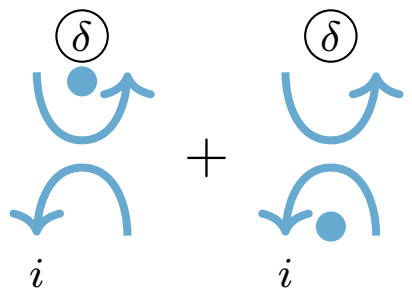

with

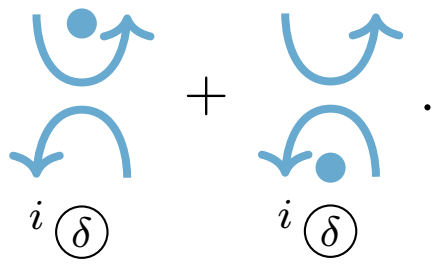

If $\delta$ only consists of $j$-dots with $j \neq i, i+1$, then $\delta$ can slide from top to bottom, without additional scalar. If $\delta$ contains at least two $j$-dots with $j = i, i+1$, then both sides rewrite to zero, possibly migrating an $i$-dot into a $(i+1)$-dot first. Finally, If $\delta$ contains exactly one $j$-dot with $j = i, i+1$, then both sides rewrites to a single diagram consisting of a $(i+1)$-dot on top and a $(i+1)$-dot on the bottom (again using dot migrations).

The other branchings are treated similarly. For the second and fourth branchings, one can use dot migrations to rewrite $i$- and $(i+1)$-dots into $(i+2)$-dots, allowing them to cross the $i$-strands. For the third branching, shading conditions (Lemma 1.5.1) prevent $(i+1)$- and $(i+2)$-dots, and $i$-dots can be treated as before, sliding then first across the $(i+1)$-strand. □

We can now prove the proposition:

*Proof of Proposition 7.3.8.* Let $[f, e, g]$ a local triple as in the proposition. Since $\mathsf{B}^+$ is convergent (Corollary 7.3.6), it suffices to show that $[f, e, g]$ is $\mathsf{B}$-congruent. As in Lemma 7.3.4, denote $\psi$ the local picture of the rewriting step $f$. That is, $\psi = s(\mathsf{nc})$ is two vertical pieces of $i$-strands (resp. $\psi = s(\mathsf{sq})$ is four vertical pieces of strands), with $f = \Gamma[\mathsf{nc}]$ (resp. $f = \Gamma[\mathsf{sq}]$) for some context $\Gamma$. The main idea is to treat $\psi$ as an extra formal generator, performing only isotopies independent of $\psi$, or E-naturalities and B-confluences as described in Lemma 7.3.1 (braided-like E-naturalities), Lemma 7.3.2 (pivotal E-naturalities), and Lemma 7.3.9 (spatial-like $B^+$-confluences). This does not change whether $[f, e, g]$ is B-congruent, thanks to the Branchwise Tamed Congruence Lemma 6.2.11.

We describe the procedure in more details for type nc; the type sq is analogous. An example is given below, picturing only the two $i$-strands and some bubble $\phi$:

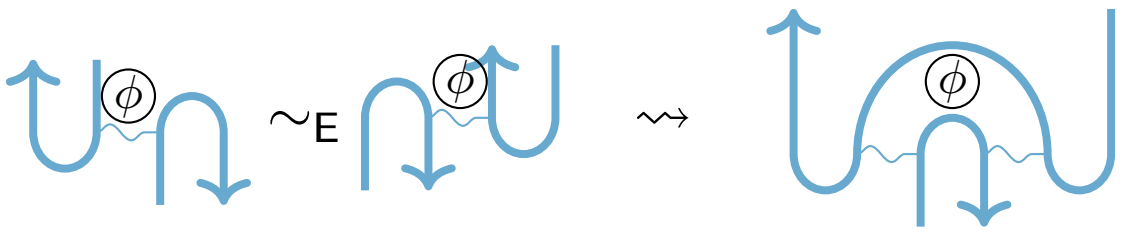

First, we move the two $i$-strands in $s(f)$ closer to one another, so that they remain parallel throughout, except possibly close to their endpoints. This procedure may require rewriting $[f, e, g]$ using independent B-rewriting steps, in order to ensure that no bubble prevents braided-like relations, similarly to the proof of Lemma 1.6.7. Doing so does not affect what we need to demonstrate, thanks to the Independent Rewriting Lemma 6.5.12.

Then, $\psi$ can be slid along the two parallel $i$-strands using pivotal E-naturalities (Lemma 7.3.2). Doing so, it may cross distant strands or bubbles: both go through thanks to braided-like E-naturalities (Lemma 7.3.1) and spatial-like $B^+$-confluences (Lemma 7.3.9). Applying the same procedure to $s(g)$ eventually leads to the same diagram, up to isotopies independent of $\psi$. ☐

**Lemma 7.3.10.** *Every monomial* $(\mathtt{S} \setminus \mathtt{T})^+$*-rewriting step is* E*-congruent to a* B*-congruence.*

*Proof.* By Proposition 7.3.8, a monomial $(\mathtt{S} \setminus \mathtt{T})^+$-rewriting step $f$ of type nc is branchwise $B^+$-confluent to a rewriting step which, up to contextualization, has the following form:

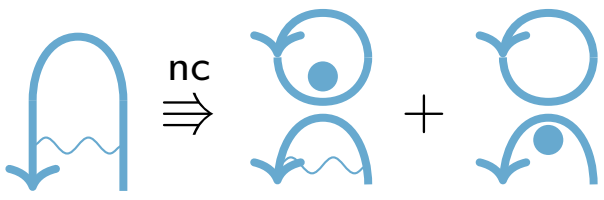

which is readily seen to be $B^+$-confluent. The Contextualization Lemma 6.5.7 together with the Branchwise Tamed Congruence Lemma 6.2.11 implies that $f$ is $\succ$-tamely $B^{st}$-congruent. A similar argument holds for $(\mathtt{S} \setminus \mathtt{T})^+$-rewriting step of type sq, reducing the argument to the following rewriting step, readily seen to be $\mathsf{B}^+$-confluent (see e.g. Lemma 1.4.4):

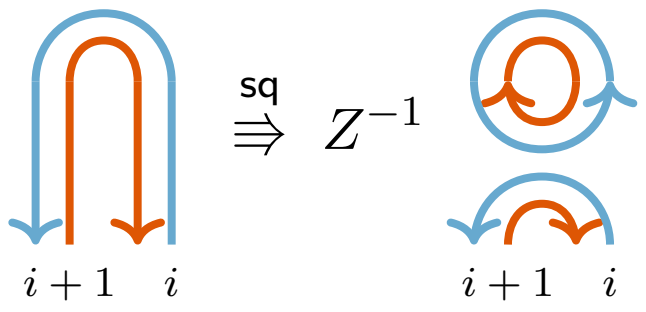

This concludes. ☐

As a direct corollary, we have that $\mathtt{T}$ and $\mathtt{S}$ present the same underlying module (Lemma 7.2.7). Moreover:

**Corollary 7.3.11.** *Every independent local* $\mathtt{T}^+$*-branching is* $\succ$*-tamely* $\mathtt{T}^{\mathrm{st}}$*-congruent.*

*Proof.* This follows from Lemma 6.5.11. ☐

### 7.3.4 Confluence of monomial local $(\mathtt{A}^+, \mathtt{B}^+)$-branchings

Recall the notation $\mathtt{A} = \mathtt{T} \setminus \mathtt{B}$. We use the terminology $(\mathtt{A}^+, \mathtt{B}^+)$*-branching* to refer to a branching which has one branch in $\mathtt{A}^+$ and the other in $\mathtt{B}^+$.

**Lemma 7.3.12.** *Every monomial local* $(\mathtt{A}^+, \mathtt{B}^+)$*-branching is* $\succ$*-tamely* $\mathtt{T}^{\mathrm{st}}$*-congruent.*

*Proof.* Throughout we use that independent local $\mathtt{T}^+$-branchings are $\succ$-tamely $\mathtt{T}^{\mathrm{st}}$-congruent (Corollary 7.3.11) without further mention. We also stop explicitly mentioning the Branchwise Tamed Congruence Lemma 6.2.11.

Let $[f, e, g]$ be a monomial $(\mathtt{A}^+, \mathtt{B}^+)$-local triple. As a result of the characterization of $\mathsf{B}$ (Lemma 7.3.4) and the characterization of $\{\mathsf{nc}, \mathsf{sq}\}$ (Proposition 7.3.8), we can freely choose the combinatorial representatives of $f$ and $g$. For instance, if $f$ and $g$ have distinct combinatorial data then we can choose the combinatorial representatives of $f$ and $g$ so that $[f, e, g]$ is an independent branching and hence $\succ$-tamely $\mathtt{T}^{\mathrm{st}}$-congruent. In particular, if $g$ is of type $\mathsf{dd}$ then $[f, e, g]$ is automatically $\succ$-tamely $\mathtt{T}^{\mathrm{st}}$-congruent. A similar reasoning applies if $g$ is of type $\mathsf{dm}$, choosing to apply the dot migration away from where the neck-cutting or the squeezing relation happens.

Assume that $f$ is of type $\mathsf{nc}$ and $g$ is of type $\mathsf{bb}_{\circlearrowleft}$, such that the strand in its combinatorial data (call it $s_1$) is one of the two $i$-strands in the combinatorial data of $f$ (call them $s_1$ and $s_2$). The fact that $s_1$ is closed forces $s := s_1 = s_2$, so that $f$ is in fact not a $\mathtt{T}^+$-rewriting step to start with. If instead $g$ is of type $\mathsf{bb}_{\circlearrowleft}$ sharing part of its combinatorial data with $f$, then $[f, e, g]$ is branchwise $\succ$-tamely $\mathtt{T}^{\mathrm{st}}$-congruent to a contextualization of the $\mathtt{T}^{+}$-confluent branching pictured in Fig. 7.2. Since the given $\mathtt{T}^{+}$-confluence in $\mathsf{B}$, the Contextualization Lemma 6.5.7 applies.

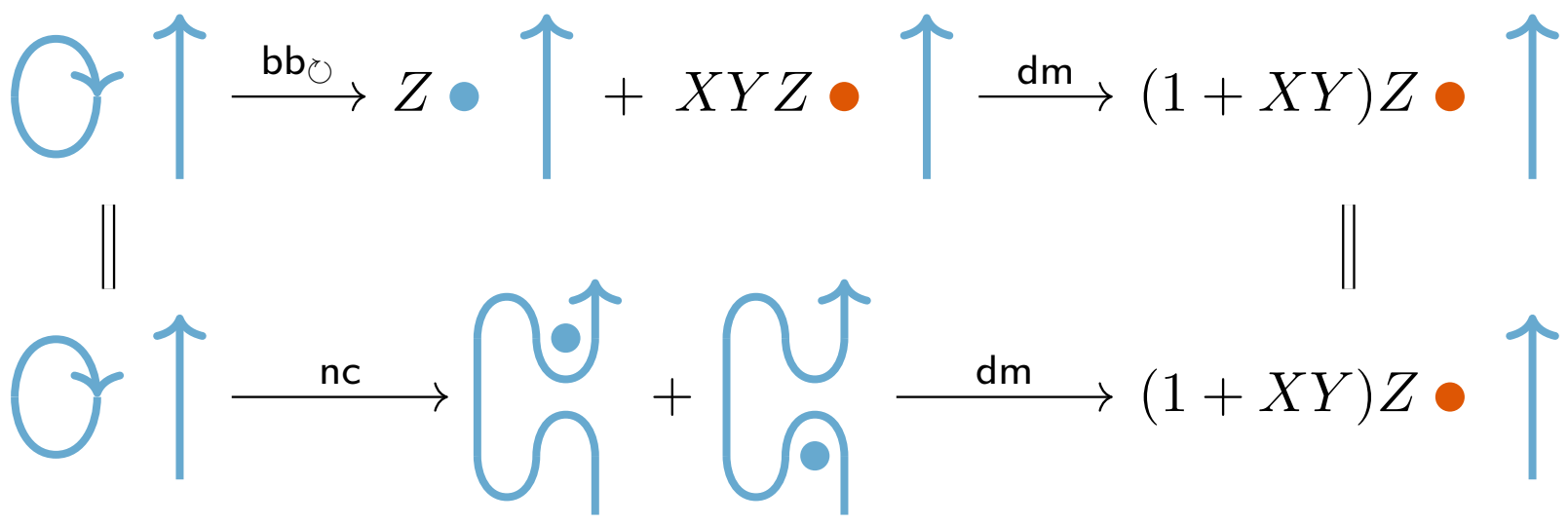

**Fig. 7.2** Critical branching between types B and type nc.

Similar arguments can be given if $f$ is of type $\mathsf{nc}$ and $g$ is of type $\mathsf{bb}_{\circlearrowleft}$ (resp. $\mathsf{bb}_{\circlearrowright}$), showing that $[f, e, g]$ is branchwise $\succ$-tamely $\mathtt{T}^{\mathrm{st}}$-congruent to a contextualization of the $\mathtt{T}^{+}$-confluent branching pictured in Fig. 7.3. In that case however, the $\mathtt{T}^{+}$-confluences are not in $\mathsf{B}$, and more care is needed to apply the Contextualization Lemma 6.5.7; in fact, it *does not* apply to the first of the two critical branchings in Fig. 7.3. One can give yet another critical branching to deal with this case. Instead, we describe another argument that does not require further computation and works for all critical branchings. This argument is essentially the same used in Lemma 6.5.11.

Note that each confluence in Fig. 7.2 and Fig. 7.3 contains at most one $\mathtt{A}^{+}$-rewriting step; called it $r$. Assume that for some context $\Gamma$, one of these confluences fails to verify the hypothesis of the Contextualization Lemma 6.5.7. That means that $\Gamma[r]$ is not in $\mathtt{A}$, and hence by Lemma 7.3.10, $[f, e, g]$ is $\mathsf{B}$-congruent. We conclude that $[f, e, g]$ is $\mathsf{B}^{+}$-confluent, thanks to the convergence of $\mathsf{B}$ (Corollary 7.3.6). ☐

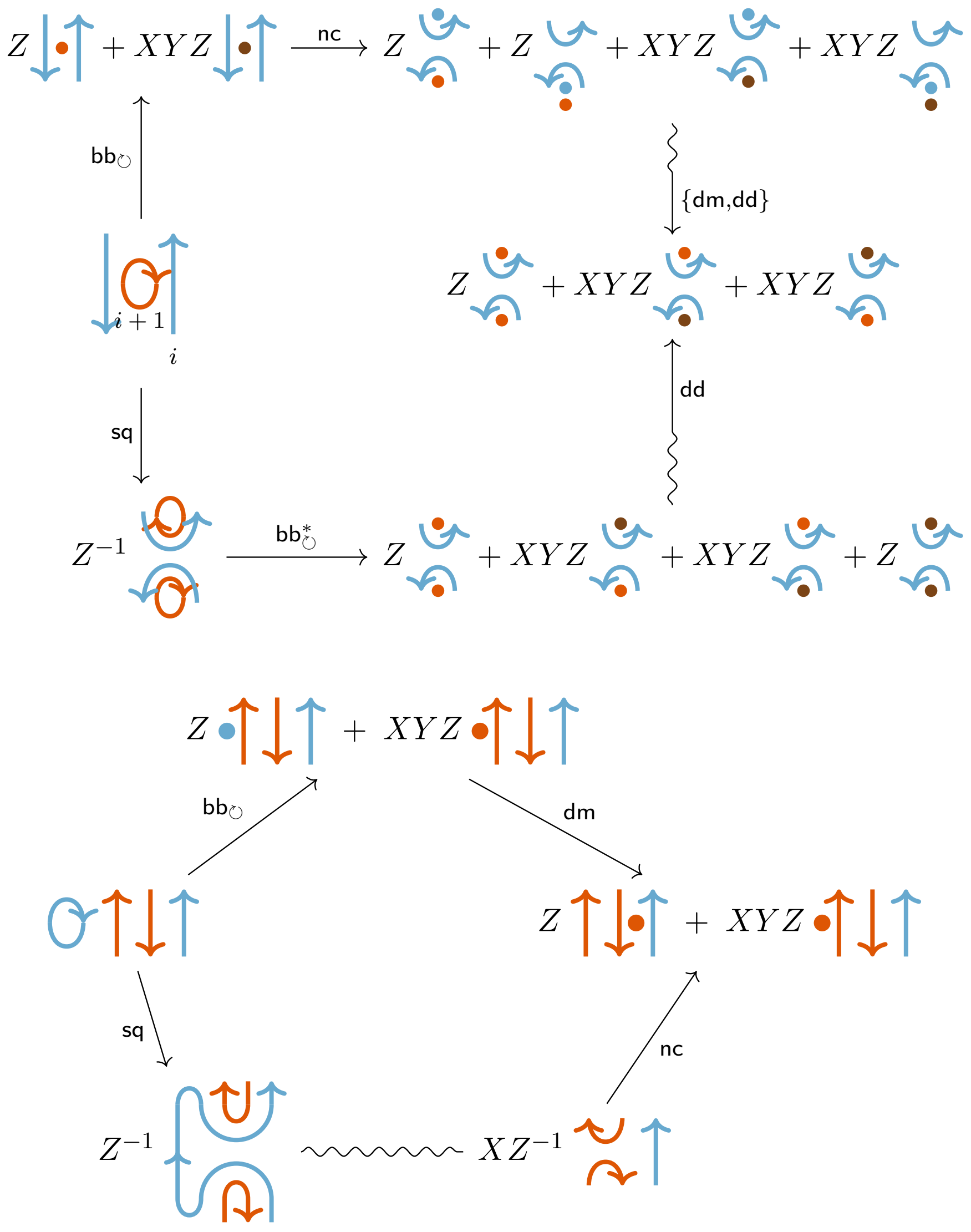

**Fig. 7.3** Critical branchings between types B and type sq.

### 7.3.5 Confluence of monomial local $\mathtt{A}^+$-branchings

**Lemma 7.3.13.** *Every monomial local $\mathtt{A}^+$-branching is $\succ$-tamely $\mathtt{T}^{\mathrm{st}}$-congruent.*

*Proof.* In what follows, we use that independent local $\mathtt{T}^+$-branchings are $\succ$-tamely $\mathtt{T}^{\mathrm{st}}$-congruent (Corollary 7.3.11), the Branchwise Tamed Congruence Lemma 6.2.11 and the Contextualization Lemma 6.5.7 without explicit mentions.

Let $[f, e, g]$ be a monomial $\mathtt{A}^+$-local triple. As in the proof of Lemma 7.3.12, the characterization of $\{\mathsf{nc}, \mathsf{sq}\}$ Proposition 7.3.8 together with Lemma 6.2.11 implies that we can freely choose the combinatorial representatives of $f$ and $g$.

Contrary to the proof of Lemma 7.3.12 however, even if $f$ and $g$ have distinct combinatorial data, there may not exist combinatorial representatives for which $f$ and $g$ are independent. On the other hand, we shall see that in most cases, we can choose the combinatorial representatives of $f$ and $g$ such that $[f, e, g]$ *rewrites* into an independent branching, in the sense of subsection 6.5.4. Only two cases will not follow this scheme: they are the critical branchings given in Fig. 7.4.

Consider first the case where both $f$ and $g$ are of type $\mathsf{nc}$. If their respective colour are $i$ and $j$, we can assume that either $j = i$ or $j = i + 1$ (otherwise, we can choose combinatorial representatives such that $[f, e, g]$ is an independent branching). In these cases, we can choose combinatorial representatives such that $[f, e, g] = \Gamma[f', e', g']$ for some context $\Gamma$ and $[f', e', g']$ as pictured below:

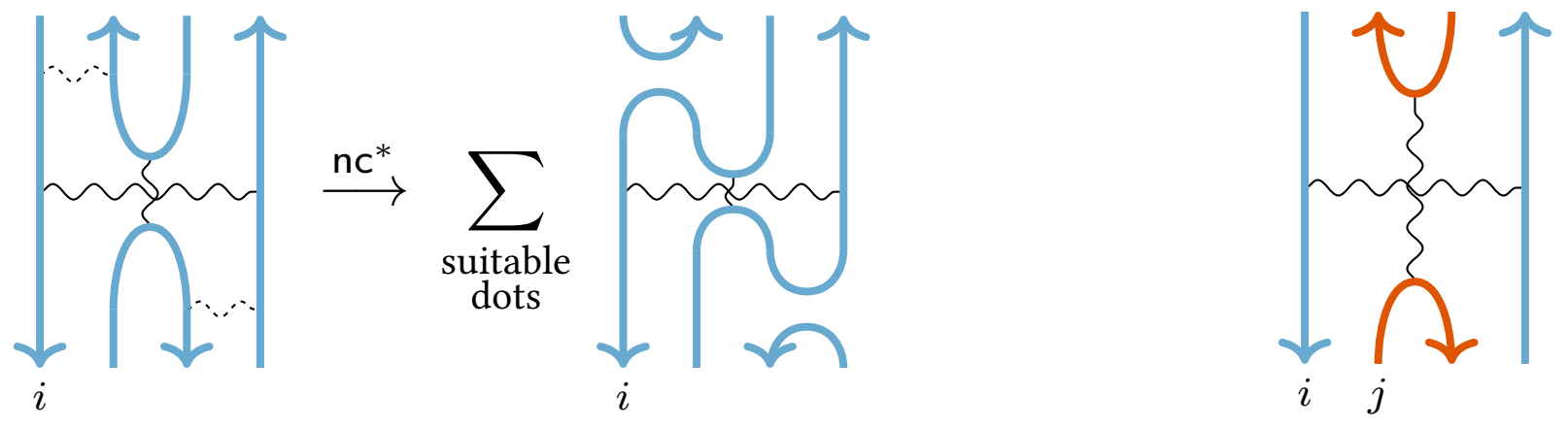

If $j = i$, then $\Gamma[f', e', g']$ rewrites via two neck-cuttings into a linear combination of branchings $\Gamma[f'', e'', g'']$ with $f''$ and $g''$ associated to the same combinatorial data, as pictured in the schematic above (we only picture the $i$-strands, leaving the dots implicit). If $\Gamma$ does not connect any of the $i$-strands involved, then all monomial rewriting steps involved are in $\mathtt{T}$, and we can apply the Independent Rewriting Lemma 6.5.12. Otherwise, it means that only one, or none, neck-cutting is necessary to get to the same situation, and the Independent Rewriting Lemma 6.5.12 is still applicable. In any case, $[f, e, g]$ is $\succ$-tamely $\mathtt{T}^{\mathrm{st}}$-congruent.

If $j = i+1$, the $\mathtt{T}^+$-confluence of $[f', e', g']$ is the first critical branching pictured in Fig. 7.4. As this confluence only uses one rewriting step of type $\mathtt{A}$, we conclude as in Lemma 7.3.12 that $[f, e, g] = \Gamma[f', e', g']$ is $\succ$-tamely $\mathtt{T}^{\mathrm{st}}$-congruent.

Consider now the case where $f$ is of type $\mathsf{sq}$ and $g$ is of type $\mathsf{nc}$. Denote $(i, i+1)$ and $j$ their respective colour. As before, we can assume that $j = i-1$, $i$, $i+1$ or $i+2$. Then, up to choice of combinatorial representatives and vertical symmetry, $[f, e, g]$ is equal to $\Gamma[f', e', g']$ for some context $\Gamma$ and $[f', e', g']$ as pictured below:

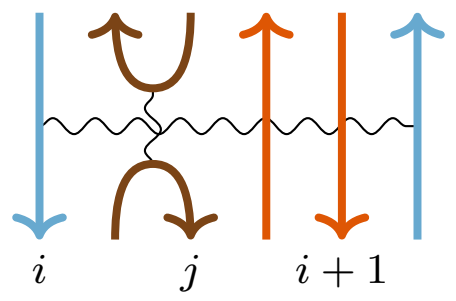

If $j = i+2$, we can further isotope the branching to get an independent branching. If $j = i$ or if $j = i+1$, we can rewrite it using neck-cuttings and conclude similarly as above.

The last case $j = i-1$ leads to the second critical branching pictured in Fig. 7.4. Note that flipping everything vertically leads to the same critical branching because $j$ and $i+1$ are distant colours. The confluence uses two rewriting steps of type $\mathtt{A}$, one of type $\mathsf{sq}$ and the other of type $\mathsf{nc}$; we refer to them simply as "$\mathsf{sq}$" and "$\mathsf{nc}$". If $\Gamma[\text{"}\mathsf{sq}\text{"}]$ does not belong to $\mathtt{T}$, then $f$ does not belong to $\mathtt{T}$ either; hence if $[f, e, g] = \Gamma[f', e', g']$ is a $\mathtt{T}^+$-branching, $\Gamma[\text{"}\mathsf{sq}\text{"}]$ necessarily belong to $\mathtt{T}$. On the other hand, if $\Gamma[\text{"}\mathsf{nc}\text{"}]$ does not belong to $\mathtt{T}$, we can replace it with a $\mathsf{B}^+$-confluence. We conclude that $[f, e, g]$ is $\succ$-tamely $\mathtt{T}^{\mathrm{st}}$-congruent.

Finally, consider the case where both $f$ and $g$ are of type $\mathsf{sq}$, with respective colours $(i, i+1)$ and $(j, j+1)$. Without loss of generality we can assume that $j \geq i$, and furthermore that either $j = i$, $j = i+1$ or $j = i+2$. Each case will allow different choices of combinatorial representatives. To help the exposition, we fix the positions of the $i$ and $(i+1)$-strands associated to $f$, simply called the *$i$- and $(i+1)$-strands* below, and discuss how the $j$- and $(j+1)$-strands associated to $g$, simply called the *$j$- and $(j+1)$-strands* below, can be isotoped with regard to the $i$- and $(i+1)$-strands.

If $j = i+2$, then $j$ is adjacent to $i+1$. In particular, the two $j$-strands cannot be isotoped through the $(i+1)$-strands. As there exists an isotopy joining the $j$ and $(j+1)$-strands, the four strands must lie on one side of

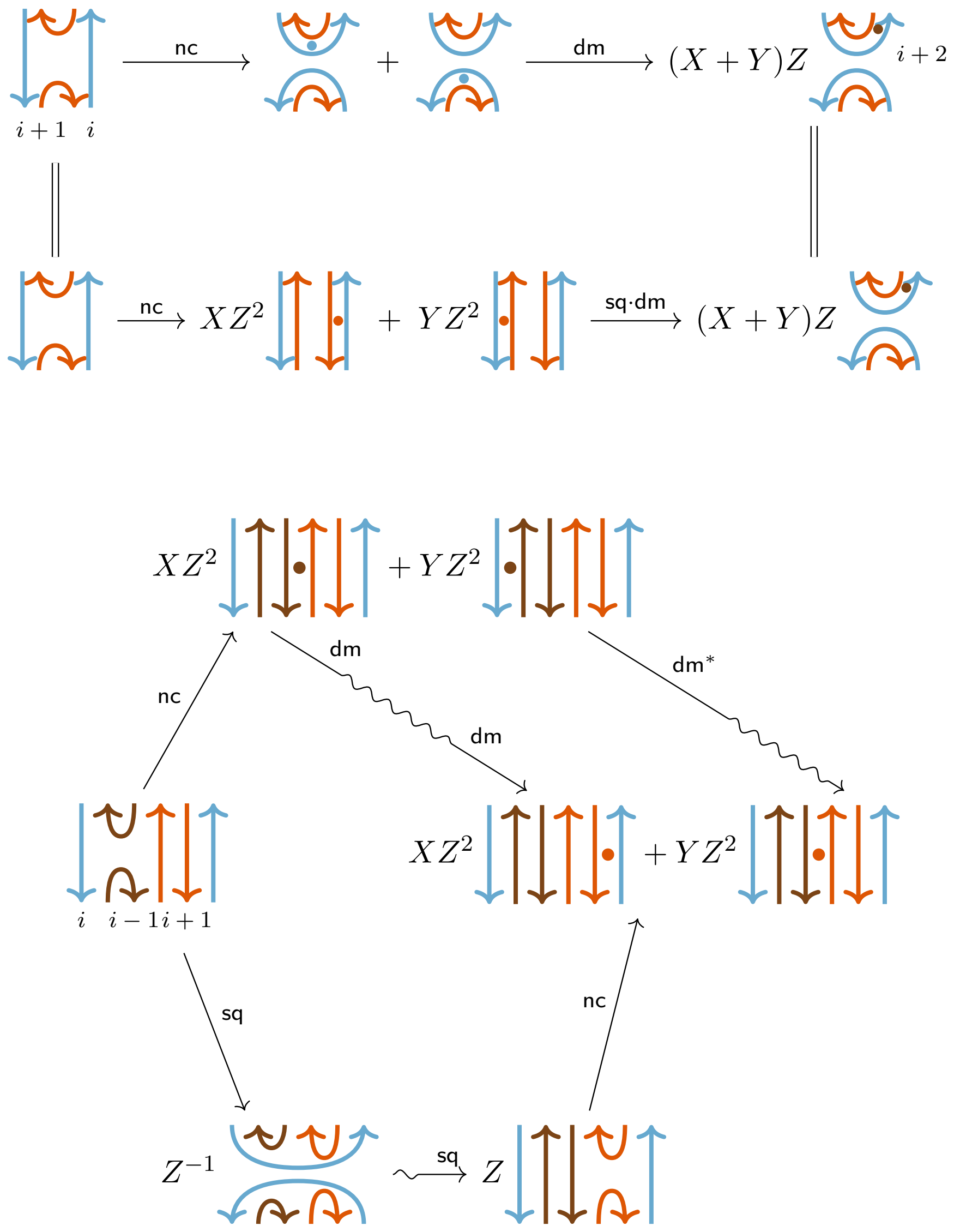

Fig. 7.4 Critical branchings in types {nc, sq}

the $(i+1)$-strands (at least partially for the $(j+1)$-strands), as pictured in the following schematic (replacing wiggly lines with straight lines to avoid clutter):

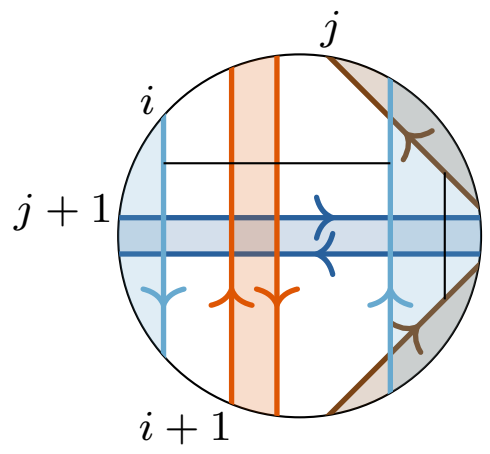

On that side, there is one $i$-strand. Given that $i$ is distant from $j$ and $j+1$, it can be isotoped through the $j$ and $(j+1)$-strands. In that way, $(f,g)$ is branchwise E-congruent to an independent branching.

If $j = i+1$, then the $(i+1)$-strands and $(j+1)$-strands cannot intersect. This leads to three possible schematics, up to symmetries:

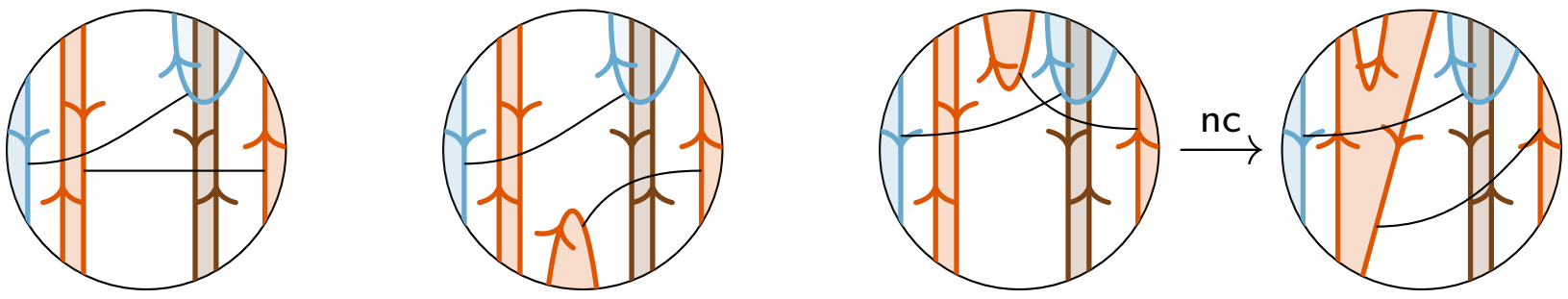

In the first schematic, one of the $(i+1)$-strands coincides with one of the $j$-strands. The first two schematic are independent branchings on the nose. The last one rewrites into an independent branching, as pictured.

Finally, if $j = i$, either have the $(i+1)$-strands and $(j+1)$-strands do not intersect, or they coincide. Below we only picture the two schematics for which rewriting the branching is necessary:

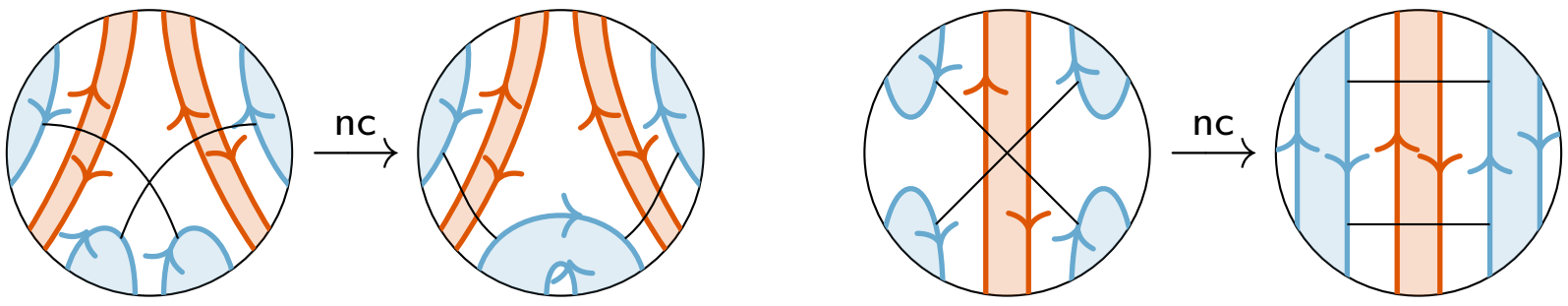

This concludes. □

## 7.4 Addendum: another deformation of $\mathfrak{gl}_2$-foams

As we have seen in the case of $\mathfrak{gl}_2$-webs (subsection 6.6.3), a rewriting approach allows one to explore graded analogues. We apply the same procedure to graded $\mathfrak{gl}_2$-foams, up to some simplifying assumptions.

The structure of graded-2-category on $\mathbf{GFoam}_d$ is already somehow universal. Indeed, the abelian group $\mathbb{Z}^2$ is isomorphic to the abelian group presented by generators $D, \cup^L, \cup^R, \cap^L, \cap^R$ and relations

$$\begin{aligned} \cap^R + \cup^R &= 0 \\ \cap^L + \cup^L &= 0 \\ \cup^R + \cap^L + D &= 0 \\ \cup^L + \cap^R &= D \\ \cup^L + \cup^R + \cap^L + \cap^R &= 0 \end{aligned}$$

obtained from taking the “abelianization” of the defining local relations in Fig. 1.2. Hence, if we assume that the grading is independent on the colours of the generators, then the $\mathbb{Z}^2$-grading on $\mathbf{GFoam}_d$ is the most general one. Moreover, a symmetric bilinear map $\mu\colon \mathbb{Z}^2 \times \mathbb{Z}^2 \to \Bbbk^\times$ is determined by its values on the generators of $\mathbb{Z}^2 \times \mathbb{Z}^2$, with relations:

$$\begin{aligned} \mu((1,0),(1,0))^2 = \mu((0,1),(0,1))^2 = 1 \\ \text{and } \mu((1,0),(0,1))\mu((0,1),(1,0)) = 1. \end{aligned}$$

This gives the parameters $X$, $Y$ and $Z$.

Let us now look how the defining relations could be deformed. For simplicity, we assume that the braid-like relations, pitchfork relations and dot slide (i.e. relations captured by $\mathsf{X}$) remain scalar-free. Up to normalization, we can further assume that the scalars of the zigzag relations, the dot migration ($\mathsf{dm}$) and the evaluation of counter-clockwise bubble ($\mathsf{bb}_{\circlearrowleft}$) keep the same scalars.

Going over all the critical branchings, one finds exactly one extra possibility where nc remains as it is, and we have:

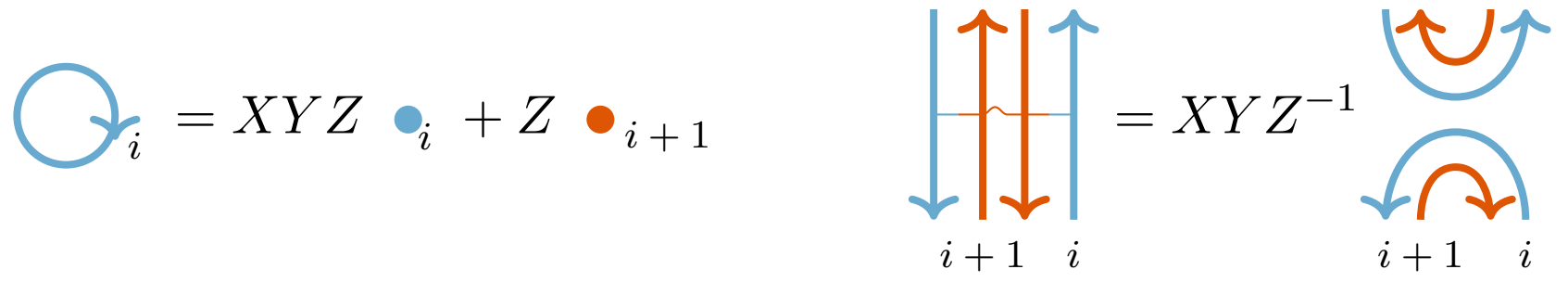

This defines an a priori distinct $(\mathbb{Z}^2, \mu)$-graded-2-category $\mathbf{GFoam}'_d$. Renormalizing the rightward cap, we can also define $\mathbf{GFoam}'_d$ by only modifying the zigzag relations as follows:

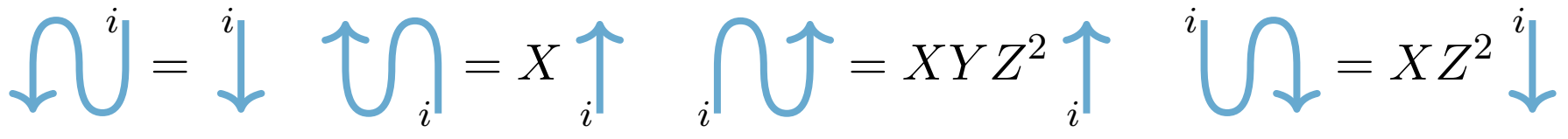

This other graded deformation has exactly the same properties as $\mathbf{GFoam}_d$. In particular, it also categorifies $\mathbf{Web}_d$ in the same sense as in Theorem 1.7.1 and can be used to define a tangle invariant generalizing odd Khovanov homology.

*Remark* 7.4.1. As we noted in Remark 2.3.11, $\mathbf{GFoam}_d$ defines "type X" odd Khovanov homology, while $\mathbf{GFoam}'_d$ defines "type Y" odd Khovanov homology (see Remark 2.2.2), although these two variants are isomorphic. The existence of these variants is explained by so-called "ladybug squares", certain squares in the hypercube that compose to zero. Interestingly, the same ladybug squares are at the heart of *Khovanov homotopy type* [121], a stable homotopy refinement of Khovanov homology. Most of the refinement is canonical, except on ladybug squares, for which a choice has to be made. This failure of canonicity is what allows Khovanov homotopy type to be a strictly stronger invariant than Khovanov homology [122, 171].

# A

# Coherence for Gray categories

In this appendix, we state a coherence result for interchangers in Gray categories, similar to the coherence theorem for braided monoidal categories. We expect this result to be known to experts, although we could not find any reference in the literature.

## A.1 Gray categories

We reproduce the definition given in [81, section 3.2]. A *Gray category* is a 3-sesquicategory together with a family of distinguished 3-cells, the *interchangers*

$$X_{\phi,\psi}\colon (\phi \star_0 g') \star_1 (f \star_0 \psi) \Rrightarrow (f' \star_0 \psi) \star_1 (\phi \star_0 g),$$

defined for each pair of 0-composable 2-cells $\phi\colon f \Rightarrow f'$ and $\psi\colon g \Rightarrow g'$. An interchanger can be pictured as:

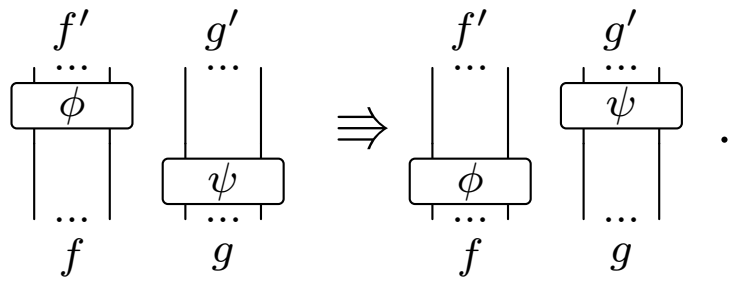

Furthermore, interchangers are required to satisfy the following set of axioms:

(i) *compatibility with compositions and identities*: for 2-cells $\phi\colon f \Rightarrow f'$, $\phi'\colon f' \Rightarrow f''$, $\psi\colon g \Rightarrow g'$ and $\psi'\colon g' \Rightarrow g''$ and 1-cells $e, h$ suitably composable:

$$X_{\mathrm{id}_f,\psi} = \mathrm{id}_{f\star_0\psi} \qquad X_{\phi,\mathrm{id}_g} = \mathrm{id}_{\phi\star_0 g}$$
$$X_{\phi'\star_1\phi,\psi} = (X_{\phi',\psi} \star_1 (\phi \star_0 g)) \star_2 ((\phi' \star_0 g') \star_1 X_{\phi,\psi})$$
$$X_{\phi,\psi'\star_1\psi} = ((f' \star_0 \psi') \star_1 X_{\phi,\psi'}) \star_2 (X_{\phi,\psi'} \star_1 (f \star_0 \psi))$$

and

$$X_{e\star_0\phi,\psi} = e \star_0 X_{\phi,\psi} \qquad X_{\phi,\psi\star_0 h} = X_{\phi,\psi} \star_0 h$$

Moreover, for 2-cells $\phi, \psi$ and 1-cell $f$ suitably 0-composable:

$$X_{\phi\star_0 f,\psi} = X_{\phi,f\star_0\psi}$$

(ii) *naturality of interchangers*: For 3-cells

$$A\colon \phi \Rrightarrow \phi'\colon u \Rightarrow u' \quad \text{and} \quad B\colon \psi \Rrightarrow \psi'\colon v \Rightarrow v',$$

such that $A$ and $B$ are 0-composable:

$$((A \star_0 v) \star_1 (u' \star_0 \psi)) \star_2 X_{\phi',\psi} = X_{\phi,\psi} \star_2 ((u \star_0 \psi) \star_1 (A \star_0 v'))$$
$$((\phi \star_0 v) \star_1 (u' \star_0 B)) \star_2 X_{\phi,\psi'} = X_{\phi,\psi} \star_2 ((u \star_0 B) \star_1 (\phi \star_0 v'))$$

Furthermore, the interchange law for 3-cells holds strictly in a Gray category. That is, for all 1-composable 3-cells $A\colon \phi \Rrightarrow \phi'$ and $B\colon \psi \Rrightarrow \psi'$, we have:

$$\big(A \star_1 t(B)\big) \star_2 \big(s(A) \star_1 B\big) = \big(t(A) \star_1 B\big) \star_2 \big(A \star_1 s(B)\big).$$

This ends the definition of a Gray category. ◇

*Remark* A.1.1. To understand the statement that "not every tricategory is equivalent to a strict 3-category", one can specialize to tricategories with only one object and only one morphism, which are precisely braided monoidal categories [36]. Similarly, one-object one-morphism Gray categories are precisely braided monoidal categories with strict associator and unitors, that is, braided *strict* monoidal categories[1]. We recover the known fact that any braided monoidal category is braided equivalent to a braided strict monoidal category. On the other hand, a one-object one-morphism strict 3-category is a strict monoidal category equipped with a trivial braiding; in other words, each component of the braiding $\beta_{X,Y} \colon X \otimes Y \to Y \otimes X$ is an identity, which in particular implies that $X \otimes Y = Y \otimes X$[2]. This latter condition is certainly very restrictive, making such monoidal categories rather degenerate instances of braided monoidal categories. Similarly, one should think of strict 3-categories as rather degenerate instances of tricategories.

## A.2 Coherence of interchangers

Recall that a braided strict monoidal category $\mathcal{C}$ is a strict monoidal category equipped with a natural family of isomorphisms $\{\beta_{\phi,\psi} \colon \phi \otimes \psi \to \psi \otimes \phi\}_{\phi,\psi \in \mathrm{ob}(\mathcal{C})}$ satisfying a certain hexagon relation. The string diagrammatics associates to the braiding an actual crossing:

$\beta_{\phi,\psi} \mapsto$ 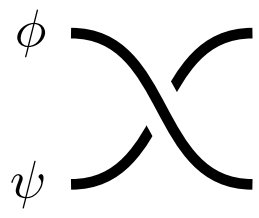

,  $\beta^{-1}_{\phi,\psi} \mapsto$ 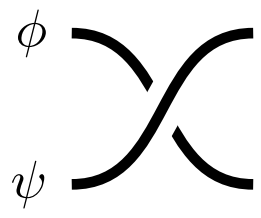

.

This assigns a braid word to any composition of braiding isomorphisms, labelled with the relevant objects. We say that the composition *represents* the corresponding (labelled) braid.

**Theorem A.2.1** (coherence for braided monoidal categories [97])**.** *Let* $\mathsf{Q}_2$ *be a set of objects and denote* $\mathsf{Q}_2\mathcal{B}raid$ *the free braided strict monoidal category with* $\mathsf{Q}_2$ *as set of objects. Two parallel morphisms in* $\mathsf{Q}_2\mathcal{B}raid$ *are equal if and only if they represent the same braid.*

Let $\mathsf{Q}$ be a 2-polygraph and denote $\mathsf{Q}\mathcal{G}ray$ the free Gray category generated by $\mathsf{Q}$.[3] Recall that a braided strict monoidal category is a one-object one-morphism Gray category. Under this equivalence, $\mathsf{Q}_2\mathcal{B}raid$ is precisely the free

[1] There are more commonly called strict braided monoidal categories, but this terminology is slightly misleading, given that the braiding need not be strict.

[2] For that reason, such monoidal categories are sometimes called *commutative* (see [10] and this MathOverflow question).

[3] $\mathsf{Q}\mathcal{G}ray$ can be more precisely defined as $\mathsf{Q}\mathcal{G}ray = ([\mathsf{QGray}]^{\equiv})^{\top}$, where $\equiv$ is the minimal higher congruence making $\mathsf{QGray}$ a Gray rewriting system (see subsection 6.4.6) and $\top$ denotes the localisation.

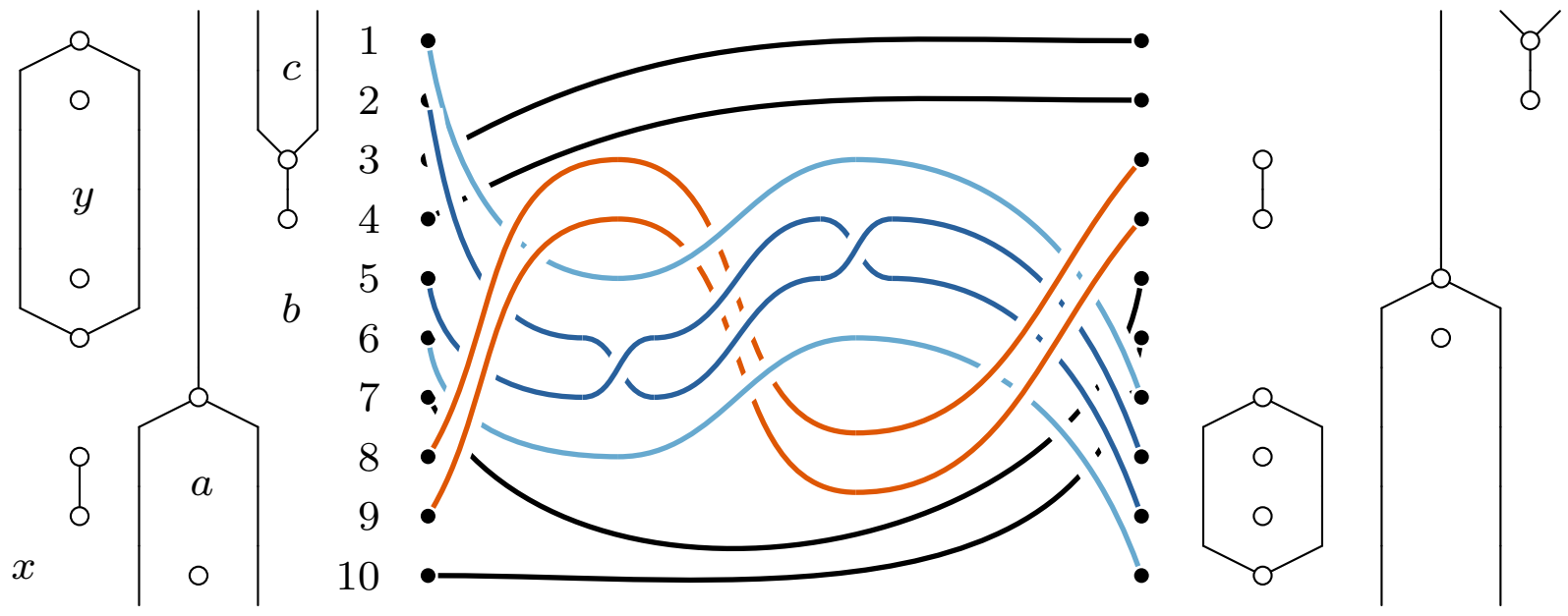

**Fig. A.1** A 3-morphism in $\mathsf{Q}\mathcal{G}ray$ and the braid it represents. The regions of the source are labelled $x$, $y$, $a$, $b$ and $c$.

Gray category associated to the 2-sesquipolygraph $\{*\} \Longleftarrow \{*\} \Longleftarrow \mathsf{Q}_2$, with $\{*\}$ denoting the one-object set, each source and target map being trivial. This leads to a faithful Gray functor $\mathsf{Q}\mathcal{G}ray \to \mathsf{Q}_2\mathcal{B}raid$. In particular, there is a notion of braids represented by 3-morphisms in $\mathsf{Q}\mathcal{G}ray$. See Fig. A.1 for an example.

The following is a direct corollary of Theorem A.2.1:

**Corollary A.2.2.** *Let* $\mathsf{Q}$ *be a 2-sesquipolygraph and* $\mathsf{Q}\mathcal{G}ray$ *the free Gray category generated by* $\mathsf{Q}$*. Two parallel 3-morphisms in* $\mathsf{Q}\mathcal{G}ray$ *are equal if and only if they represent the same braid.* □

Therefore, the problem of coherence for Gray categories reduces to understanding when two parallel 3-morphisms represent the same braid, taking the data of objects and 1-morphisms back into account. Coherence typically fails when we have two 2-endomorphisms, say $\phi$ and $\psi$, of the same identity 1-morphism:

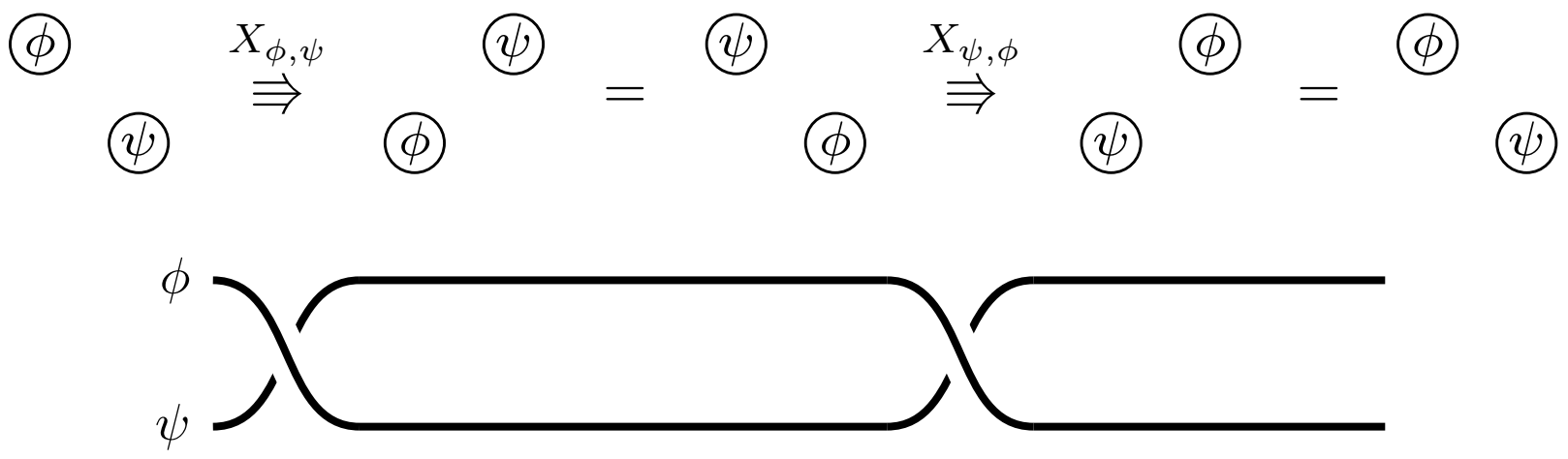

The composition $X_{\psi,\phi} \star_2 X_{\phi,\psi}$, and similarly $X_{\phi,\phi} \colon \phi \star_1 \phi \Rrightarrow \phi \star_1 \phi$, is not an identity in $\mathsf{Q}\mathcal{G}ray$. We shall see that in fact, this is the only obstruction to coherence.

To formalize this idea, we introduce some terminology. A *diagram* $D$ is the string diagrammatics of a 2-cell in $\mathsf{Q}_2^*$; in what follows, every notion is defined with respect to the string diagrammatics. A *generator of* $D$ is a generating 2-cell in $\mathsf{Q}_2$ appearing in $D$. If the same 2-cell appears multiple times in $D$, we view these copies as distinct generators. Note that the set of generators is invariant under interchangers. A *connected component* is a union of generators whose string diagrammatics correspond to a connected component of $D$. It is *closed* if it is not connected to the boundary.

A *bubble* is a union of generators that coincide with a closed connected component and its interior. Note that a bubble can have other bubbles in its interior. Two bubbles are *adjacent* if they are distinct and belong to the same region (recall that different copies of the same generator are viewed as distinct). For instance, in Fig. A.1 $\{2\}$ and $\{5\}$ are adjacent bubbles in region $y$, and so are $\{1, 2, 5, 6\}$ and $\{8, 9\}$ in region $x$.

Let $\mathfrak{I}$ be a 3-morphism in $\mathsf{Q}\mathcal{G}ray$ and denote $\beta_{\mathfrak{I}}$ the braid it represents. Let $x$ be a region of $D$ and denote $\phi_1, \ldots, \phi_n$ the set of adjacent bubbles in $x$. Pick a generator $\alpha_k$ in each bubble $\phi_k$. The *bubble braid associated to* $x$ is the labelled braid obtained from $\beta_{\mathfrak{I}}$ by picking out the strands associated to the $\alpha_k$s, labelled with $\phi_k$, and regarded up to conjugation. For instance, in Fig. A.1 the only regions with a non-trivial bubble braid are $x$ and $y$, and the bubble braid is respectively given by (irrespective of the choice of generators):

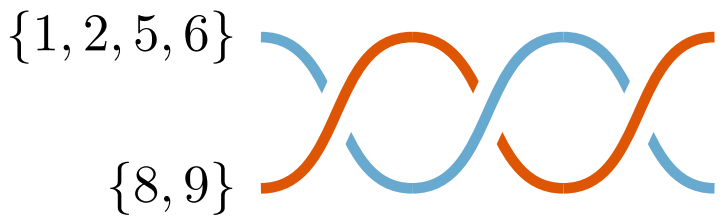

bubble braid associated to $x$

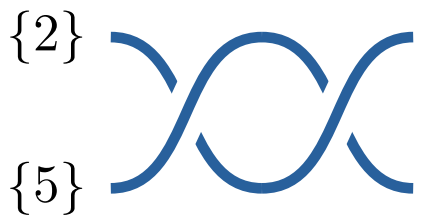

bubble braid associated to $y$

**Lemma A.2.3.** *If $\mathfrak{I}$ is an endomorphism, the definition of a bubble braid is independent of the choice of generators $\alpha_k$s.*

*Proof.* We give a definition of the bubble braid readily independent of the choice of generators $\alpha_k$s. Let $i_k$ be the number of generators in $\phi_k$. Let $\beta_x$ be the braid obtained from $\beta_{\mathfrak{I}}$ by picking out the strands associated to all the generators of the $\phi_k$s. The braid $\beta_x$ defines a loop $l_x$, or more precisely an element of the $\pi_1$, in the configuration space of $i_1 + \ldots + i_n$ points in the plane.

Given the definition of $\beta_x$, the loop $l_x$ extends to a loop in the configuration space of bubbles $\phi_k$ in the plane, that is, the space of planar embeddings of the planar graph $\phi_1 \sqcup \ldots \sqcup \phi_n$. Contracting each bubble to a point gives a path $[l_x]$ in the configuration space of $n$ points in the plane. Choosing an ordering on these $n$ points turns this loop back into a braid. Choosing a different ordering leads to the same braid, up to conjugation. Ordering the $n$ points in the order the $\alpha_k$s appear recovers the previous definition of the bubble braid. □

As described in [49] regions and bubbles in a diagram can be organized in a tree. For instance, here is the tree for the diagram in the example, on the left-hand side:

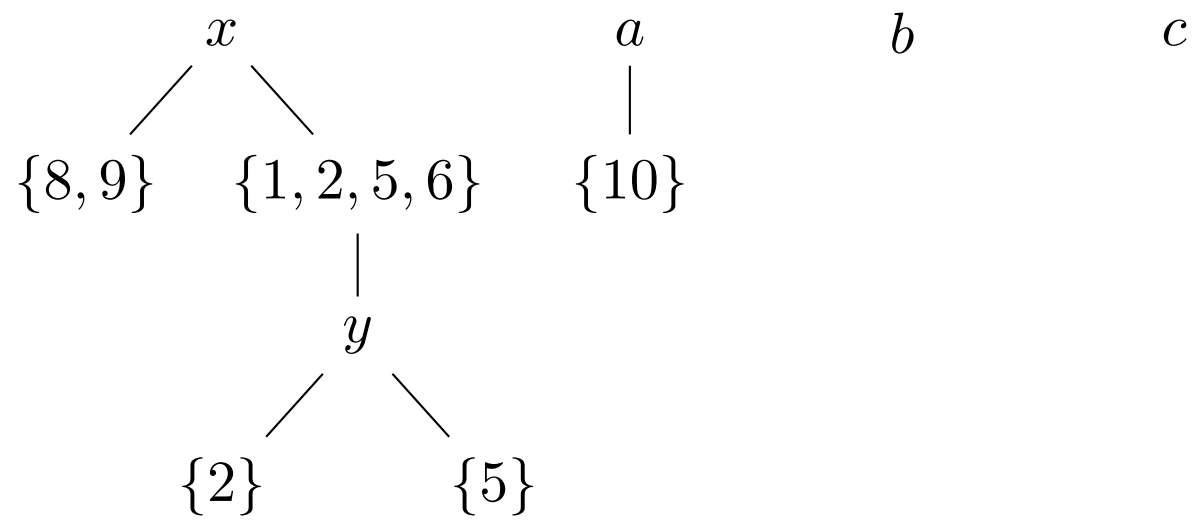

Following this tree and the data of a bubble braid for each region, one can inductively recover the braid presented by the corresponding 3-endomorphism. We leave the details to the reader. This leads to the following coherence result for Gray categories:

**Theorem A.2.4** (coherence for Gray categories)**.** *Let* Q *be a 2-sesquipolygraph and* Q$\mathcal{G}ray$ *the free Gray category generated by* Q*. Two parallel 3-endomorphisms in* Q$\mathcal{G}ray$ *are equal if and only if they represent the same family of bubble braids.* □

In particular, if the source (or target) of two parallel 3-endomorphisms in Q$\mathcal{G}ray$ does not contain any bubble, then the two endomorphisms are equal. More generally, if Q$\mathcal{G}ray$ does not contain any bubble, then Q$\mathcal{G}ray$ is automatically coherent. For instance, this is the case of Q $=$ Web$_{\leq 2}$ defined in section 5.1. This is also the idea behind the notion of positivity of a Gray presentation as introduced in [81, p. 606].

## A.3 Applications to graded-2-categories

Let $(G, \mu, \Bbbk)$ as in Notation 5.3.1. Let Q be a scalar $G$-graded 2-sesquipolygraph Q. Recall that in that case, QGray is scalar. Since the axioms of a Gray category are linear homogeneous with respect to the scalar structure, the free Gray category $\mathsf{Q}\mathcal{G}ray$ is also scalar.

The data $(G, \mu)$ associates a scalar in $\Bbbk^{\times}$ to any bubble braid. Given a 3-endomorphism of $\mathsf{Q}\mathcal{G}ray$, its *bubble braid scalar* is the product of all the scalars associated to its bubble braids. As a direct corollary of Theorem A.2.4, we have:

**Theorem A.3.1** (coherence of graded-interchangers)**.** *Let* Q *be a* $G$*-graded 2-sesquipolygraph and* $\mathsf{Q}\mathcal{G}ray$ *the free* $(G, \mu)$*-scalar Gray category generated by* Q*. Two parallel 3-endomorphisms in* $\mathsf{Q}\mathcal{G}ray$ *have the same associated scalar if and only if they have the same associated bubble braid scalar.* □

# Bibliography

[1] D. Aasen, E. Lake, and K. Walker, Fermion condensation and super pivotal categories, *J. Math. Phys.*, vol. 60, no. 12 (2019), pp. 121901, 111 (cit. on p. 19).

[2] M. Aganagic, Knot categorification from mirror symmetry, part I: coherent sheaves, preprint 2020, arXiv: 2004.14518 [hep-th] (cit. on p. 12).

[3] C. Alleaume, Rewriting in higher dimensional linear categories and application to the affine oriented Brauer category, *J. Pure Appl. Algebra*, vol. 222, no. 3 (2018), pp. 636–673 (cit. on pp. 35, 37, 39, 40, 42, 44, 141, 156, 157, 162, 170, 174, 176, 178, 187, 195, 197, 204).

[4] C. Alleaume, Higher-dimensional linear rewriting and coherence in categorification and representation theory, Ph.D. thesis, Université de Lyon, 2018 (cit. on pp. 39, 176, 195).

[5] D. Ara, A. Burroni, Y. Guiraud, P. Malbos, F. Métayer, and S. Mimram, Polygraphs: From Rewriting to Higher Categories, preprint 2023, arXiv: 2312.00429 [cs, math] (cit. on pp. 34, 161, 165, 170).

[6] M. Araújo, "Simple string diagrams and n-sesquicategories", *Theory Appl. Categ.*, vol. 38 (2022), Paper No. 34, 1284–1325 (cit. on pp. 37, 141).

[7] F. Baader and T. Nipkow, Term rewriting and all that, Cambridge University Press, Cambridge, 1998, xii+301 (cit. on pp. 183, 184).

[8] J. C. Baez and L. Langford, Higher-dimensional algebra IV: 2-tangles, *Adv. Math.*, vol. 180, no. 2 (2003), pp. 705–764 (cit. on p. 16).

[9] J. C. Baez and A. D. Lauda, "A prehistory of n-categorical physics", in: *Deep Beauty*, Cambridge Univ. Press, Cambridge, 2011, pp. 13–128 (cit. on p. 7).

[10] J. C. Baez and J. Master, Open Petri nets, *Math. Structures Comput. Sci.*, vol. 30, no. 3 (2020), pp. 314–341 (cit. on p. 235).

[11] K. Bar, A. Kissinger, and J. Vicary, Globular: an online proof assistant for higher-dimensional rewriting, *Log. Methods Comput. Sci.*, vol. 14, no. 1 (2018), Paper No. 8, 16 (cit. on pp. 38, 45).

[12] K. Bar and J. Vicary, "Data structures for quasistrict higher categories", in: *2017 32nd Annual ACM/IEEE Symposium on Logic in Computer Science (LICS)*, IEEE, [Piscataway], NJ, 2017, p. 12 (cit. on pp. 38, 45).

[13] D. Bar-Natan, Fast Khovanov homology computations, *J. Knot Theory Ramifications*, vol. 16, no. 3 (2007), pp. 243–255 (cit. on p. 85).

[14] D. Bar-Natan, Khovanov's homology for tangles and cobordisms, *Geom. Topol.*, vol. 9, no. 3 (2005), pp. 1443–1499 (cit. on pp. 12, 14, 85–87).

[15] S. Barbier, Diagram categories of Brauer type, preprint 2024, arXiv: 2406.18436 [math] (cit. on p. 46).

[16] S. Beier, An integral lift, starting in odd Khovanov homology, of Szabó's spectral sequence, preprint 2012, arXiv: 1205.2256 [math] (cit. on p. 21).

[17] A. Beliakova, M. Hogancamp, K. K. Putyra, and S. M. Wehrli, On the functoriality of $\mathfrak{sl}_2$ tangle homology, *Algebr. Geom. Topol.*, vol. 23, no. 3 (2023), pp. 1303–1361 (cit. on pp. 47, 69, 73).

[18] G. M. Bergman, The diamond lemma for ring theory, *Adv. in Math.*, vol. 29, no. 2 (1978), pp. 178–218 (cit. on p. 32).

[19] C. Blanchet, An oriented model for Khovanov homology, *J. Knot Theory Ramifications*, vol. 19, no. 2 (2010), pp. 291–312 (cit. on pp. 12, 15, 22, 47).

[20] L. A. Bokut', "Imbeddings into simple associative algebras", *Algebra i Logika*, vol. 15, no. 2 (1976), pp. 117–142, 245 (cit. on p. 32).

[21] F. Borceux, "Handbook of categorical algebra. 2", vol. 51, Encyclopedia of Mathematics and Its Applications, Cambridge University Press, Cambridge, 1994, xviii+443 (cit. on p. 50).

[22] J. Brundan, On the definition of Kac–Moody 2-category, *Math. Ann.*, vol. 364, no. 1-2 (2016), pp. 353–372 (cit. on pp. 10, 116).

[23] J. Brundan and N. Davidson, Categorical actions and crystals, in: *Categorification and Higher Representation Theory*, vol. 683, Contemp. Math. Amer. Math. Soc., Providence, RI, 2017, pp. 105–147 (cit. on p. 47).

[24] J. Brundan and A. P. Ellis, Monoidal supercategories, *Comm. Math. Phys.*, vol. 351, no. 3 (2017), pp. 1045–1089 (cit. on pp. 19, 48, 52, 53).

[25] J. Brundan and A. P. Ellis, Super Kac–Moody 2-categories, *Proc. London Math. Soc.*, vol. 115, no. 5 (2017), pp. 925–973 (cit. on pp. 19, 20, 26, 47, 48, 107, 108, 115, 116).

[26] J. Brundan and A. Kleshchev, Odd Grassmannian bimodules and derived equivalences for spin symmetric groups, preprint 2022, arXiv: 2203.14149 [math] (cit. on pp. 19, 20, 39).

[27] B. Buchberger, An algorithm for finding the basis elements of the residue class ring of a zero dimensional polynomial ideal, *J. Symbolic Comput.*, vol. 41, no. 3-4 (2006), pp. 475–511 (cit. on p. 32).

[28] A. Burroni, Higher-dimensional word problems with applications to equational logic, in: *Theoretical Computer Science*, vol. 115, 1, 1993, pp. 43–62 (cit. on p. 33).

[29] C. L. Caprau, $sl(2)$ tangle homology with a parameter and singular cobordisms, *Algebr. Geom. Topol.*, vol. 8, no. 2 (2008), pp. 729–756 (cit. on pp. 12, 15).

[30] J. S. Carter and M. Saito, Knotted surfaces and their diagrams, vol. 55, Mathematical Surveys and Monographs, American Mathematical Society, Providence, RI, 1998, xii+258 (cit. on p. 16).

[31] S. Cautis, Clasp technology to knot homology via the affine Grassmannian, *Math. Ann.*, vol. 363, no. 3 (2015), pp. 1053–1115 (cit. on p. 12).

[32] S. Cautis and J. Kamnitzer, Knot homology via derived categories of coherent sheaves. I. The $\mathfrak{sl}(2)$-case, *Duke Math. J.*, vol. 142, no. 3 (2008), pp. 511–588 (cit. on p. 12).

[33] S. Cautis, J. Kamnitzer, and S. Morrison, Webs and quantum skew Howe duality, *Math. Ann.*, vol. 360, no. 1-2 (2014), pp. 351–390 (cit. on p. 17).

[34] J. Cerf, “Sur les difféomorphismes de la sphère de dimension trois ($\Gamma_4 = 0$)”, vol. No. 53, Lecture Notes in Mathematics, Springer-Verlag, Berlin-New York, 1968, xii+133 (cit. on p. 61).

[35] C. Chenavier, B. Dupont, and P. Malbos, Confluence of algebraic rewriting systems, *Math. Structures Comput. Sci.*, vol. 32, no. 7 (2022), pp. 870–897 (cit. on pp. 36, 178).

[36] E. Cheng and N. Gurski, “The periodic table of n-categories II: Degenerate tricategories”, *Cah. Topol. Géom. Différ. Catég.*, vol. 52, no. 2 (2011), pp. 82–125 (cit. on p. 234).

[37] J. Chuang and R. Rouquier, Derived equivalences for symmetric groups and $\mathfrak{sl}_2$-categorification, *Ann. of Math. (2)*, vol. 167, no. 1 (2008), pp. 245–298 (cit. on pp. 10, 11).

[38] C. Chung, T. Sale, and W. Wang, Quantum supergroups VI: roots of 1, *Lett Math Phys*, vol. 109, no. 12 (2019), pp. 2753–2777 (cit. on pp. 20, 26).

[39] D. Clark, S. Morrison, and K. Walker, Fixing the functoriality of Khovanov homology, *Geom. Topol.*, vol. 13, no. 3 (2009), pp. 1499–1582 (cit. on pp. 12, 15).

[40] S. Clark, Odd knot invariants from quantum covering groups, *Algebr. Geom. Topol.*, vol. 17, no. 5 (2017), pp. 2961–3005 (cit. on p. 26).

[41] S. Clark, Quantum supergroups IV: the modified form, *Math. Z.*, vol. 278, no. 1-2 (2014), pp. 493–528 (cit. on pp. 20, 26).

[42] S. Clark, Z. Fan, Y. Li, and W. Wang, Quantum supergroups III. Twistors, *Comm. Math. Phys.*, vol. 332, no. 1 (2014), pp. 415–436 (cit. on pp. 20, 26).

[43] S. Clark and D. Hill, Quantum supergroups V. Braid group action, *Comm. Math. Phys.*, vol. 344, no. 1 (2016), pp. 25–65 (cit. on pp. 20, 26).

[44] S. Clark, D. Hill, and W. Wang, Quantum supergroups I. Foundations, *Transformation Groups*, vol. 18, no. 4 (2013), pp. 1019–1053 (cit. on pp. 20, 26).

[45] S. Clark, D. Hill, and W. Wang, Quantum supergroups II. Canonical basis, *Represent. Theory*, vol. 18, no. 9 (2014), pp. 278–309 (cit. on pp. 20, 26).

[46] R. L. Cohen, J. D. S. Jones, and G. B. Segal, "Floer's infinite-dimensional Morse theory and homotopy theory", in: *The Floer Memorial Volume*, vol. 133, Progr. Math. Birkhäuser, Basel, 1995, pp. 297–325 (cit. on p. 46).

[47] L. Crane, Clock and category: is quantum gravity algebraic?, *J. Math. Phys.*, vol. 36, no. 11 (1995), pp. 6180–6193 (cit. on p. 6).

[48] L. Crane and I. B. Frenkel, Four-dimensional topological quantum field theory, Hopf categories, and the canonical bases, *J. Math. Phys.*, vol. 35 (1994), pp. 5136–5154 (cit. on p. 6).

[49] A. Delpeuch and J. Vicary, Normalization for planar string diagrams and a quadratic equivalence algorithm, *Log. Methods Comput. Sci.*, vol. 18, no. 1 (2022), Paper No. 10, 38 (cit. on p. 238).

[50] S. K. Donaldson, An application of gauge theory to four-dimensional topology, *J. Differential Geom.*, vol. 18, no. 2 (1983), pp. 279–315 (cit. on p. 5).

[51] C. Dorn, Associative $n$-categories, preprint 2023, arXiv: 1812.10586 [math] (cit. on pp. 38, 45).

[52] A. Dranowski, M. Guo, A. Lauda, and A. Manion, Spectral 2-actions, foams, and frames in the spectrification of Khovanov arc algebras, preprint 2024, arXiv: 2402.11368 [math] (cit. on p. 46).

[53] V. G. Drinfel'd, "Quantum groups", in: *Proceedings of the International Congress of Mathematicians, Vol. 1, 2 (Berkeley, Calif., 1986)*, Amer. Math. Soc., Providence, RI, 1987, pp. 798–820 (cit. on p. 4).

[54] B. Dupont, Réécriture modulo dans les catégories diagrammatiques, Ph.D. thesis, Université de Lyon, 2020 (cit. on p. 39).

[55] B. Dupont, Rewriting modulo isotopies in Khovanov-Lauda-Rouquier's categorification of quantum groups, *Adv. Math.*, vol. 378 (2021), Paper No. 107524, 75 (cit. on pp. 37, 39, 40, 162, 167, 170, 176, 179, 188, 195).

[56] B. Dupont, Rewriting modulo isotopies in pivotal linear (2,2)-categories, *J. Algebra*, vol. 601 (2022), pp. 1–53 (cit. on pp. 37, 39, 40, 44, 162, 167, 170, 174, 176–179, 187, 188, 195, 197).

[57] B. Dupont, M. Ebert, and A. D. Lauda, Super rewriting theory and nondegeneracy of odd categorified sl(2), 2021, arXiv: 2102.00276 [math] (cit. on pp. 20, 37, 40, 44, 179, 188, 195).

[58] B. Dupont and P. Malbos, Coherent confluence modulo relations and double groupoids, *J. Pure Appl. Algebra*, vol. 226, no. 10 (2022), Paper No. 107037, 57 (cit. on pp. 36, 161, 167).

[59] J. N. Eberhardt, G. Naisse, and A. Wilbert, Real Springer fibers and odd arc algebras, *J. London Math. Soc.*, vol. 103, no. 4 (2021), pp. 1415–1452 (cit. on pp. 19, 21).

[60] M. Ebert, A new presentation of the osp(1|2)-polynomial link invariant and categorification, preprint 2022, arXiv: 2210.09583 [math] (cit. on p. 26).

[61] M. Ebert, A. D. Lauda, and L. Vera, Derived superequivalences for spin symmetric groups and odd sl(2)-categorifications, preprint 2022, arXiv: 2203.14153 [math] (cit. on p. 19).

[62] I. Egilmez and A. D. Lauda, DG structures on odd categorified quantum $sl(2)$, *Quantum Topology*, vol. 11, no. 2 (2020), pp. 227–294 (cit. on p. 19).

[63] M. Ehrig, C. Stroppel, and D. Tubbenhauer, Generic $\mathfrak{gl}_2$-foams, web and arc algebras, preprint 2020, arXiv: 1601.08010v3 (cit. on p. 47).

[64] M. Ehrig, C. Stroppel, and D. Tubbenhauer, The Blanchet-Khovanov algebras, in: *Categorification and Higher Representation Theory*, vol. 683, Contemp. Math. Amer. Math. Soc., Providence, RI, 2017, pp. 183–226 (cit. on p. 47).

[65] M. Ehrig, D. Tubbenhauer, and P. Wedrich, Functoriality of colored link homologies, *Proc. London Math. Soc.*, vol. 117, no. 5 (2018), pp. 996–1040 (cit. on p. 15).

[66] B. Elias, A diamond lemma for Hecke-type algebras, *Trans. Amer. Math. Soc.*, vol. 375, no. 3 (2022), pp. 1883–1915 (cit. on p. 37).

[67] B. Elias, The two-color Soergel calculus, *Compos. Math.*, vol. 152, no. 2 (2016), pp. 327–398 (cit. on p. 11).

[68] B. Elias and M. Khovanov, Diagrammatics for Soergel categories, *Int. J. Math. Math. Sci.* (2010), Art. ID 978635, 58 (cit. on p. 11).

[69] B. Elias, S. Makisumi, U. Thiel, and G. Williamson, Introduction to Soergel bimodules, RSME Springer Series volume 5, Cham, Switzerland: Springer, 2020, 588 pp. (cit. on pp. 27, 53).

[70] B. Elias and G. Williamson, Soergel calculus, *Represent. Theory*, vol. 20 (2016), pp. 295–374 (cit. on p. 11).

[71] B. Elias and G. Williamson, The Hodge theory of Soergel bimodules, *Ann. of Math. (2)*, vol. 180, no. 3 (2014), pp. 1089–1136 (cit. on p. 11).

[72] A. P. Ellis, The odd Littlewood-Richardson rule, *J. Algebraic Combin.*, vol. 37, no. 4 (2013), pp. 777–799 (cit. on p. 20).

[73] A. P. Ellis and M. Khovanov, The Hopf algebra of odd symmetric functions, *Adv. Math.*, vol. 231, no. 2 (2012), pp. 965–999 (cit. on p. 20).

[74] A. P. Ellis, M. Khovanov, and A. D. Lauda, The odd nilHecke algebra and its diagrammatics, *Int. Math. Res. Not. IMRN*, no. 4 (2014), pp. 991–1062 (cit. on pp. 19, 20, 107).

[75] A. P. Ellis and A. D. Lauda, An odd categorification of $U_q(\mathfrak{sl}_2)$, *Quantum Topol.*, vol. 7, no. 2 (2016), pp. 329–433 (cit. on pp. 19, 20, 26, 107).

[76] A. P. Ellis and Y. Qi, The differential graded odd nilHecke algebra, *Comm. Math. Phys.*, vol. 344, no. 1 (2016), pp. 275–331 (cit. on p. 19).

[77] P. Etingof, D. Nikshych, and V. Ostrik, On fusion categories, *Ann. of Math. (2)*, vol. 162, no. 2 (2005), pp. 581–642 (cit. on p. 11).

[78] A. Floer, An instanton-invariant for 3-manifolds, *Comm. Math. Phys.*, vol. 118, no. 2 (1988), pp. 215–240 (cit. on p. 5).

[79] A. Floer, Morse theory for Lagrangian intersections, *J. Differential Geom.*, vol. 28, no. 3 (1988), pp. 513–547 (cit. on p. 5).

[80] A. Floer, Symplectic fixed points and holomorphic spheres, *Comm. Math. Phys.*, vol. 120, no. 4 (1989), pp. 575–611 (cit. on p. 5).

[81] S. Forest and S. Mimram, Rewriting in Gray categories with applications to coherence, *Math. Structures Comput. Sci.*, vol. 32, no. 5 (2022), pp. 574–647 (cit. on pp. 37–39, 43, 141, 146, 147, 152, 153, 158, 187, 233, 238).

[82] I. Frenkel, M. Khovanov, and C. Stroppel, A categorification of finite-dimensional irreducible representations of quantum $\mathfrak{sl}_2$ and their tensor products, *Selecta Math. (N.S.)*, vol. 12, no. 3-4 (2006), pp. 379–431 (cit. on p. 10).

[83] R. Gordon, A. J. Power, and R. Street, Coherence for tricategories, *Mem. Amer. Math. Soc.*, vol. 117, no. 558 (1995), pp. vi+81 (cit. on p. 37).

[84] L. Guetta, Homology of categories via polygraphic resolutions, *J. Pure Appl. Algebra*, vol. 225, no. 10 (2021), Paper No. 106688, 33 (cit. on p. 34).

[85] Y. Guiraud, Termination orders for 3-polygraphs, *C. R. Math. Acad. Sci. Paris*, vol. 342, no. 4 (2006), pp. 219–222 (cit. on p. 33).

[86] Y. Guiraud, E. Hoffbeck, and P. Malbos, Convergent presentations and polygraphic resolutions of associative algebras, *Math. Z.*, vol. 293, no. 1-2 (2019), pp. 113–179 (cit. on pp. 32, 35, 42, 44, 46, 161, 174, 176, 178, 180, 187, 199).

[87] Y. Guiraud and P. Malbos, "Higher-dimensional categories with finite derivation type", *Theory Appl. Categ.*, vol. 22 (2009), No. 18, 420–478 (cit. on p. 33).

[88] Y. Guiraud and P. Malbos, "Identities among relations for higher-dimensional rewriting systems", in: *OPERADS 2009*, vol. 26, Sémin. Congr. Soc. Math. France, Paris, 2013, pp. 145–161 (cit. on p. 33).

[89] Y. Guiraud and P. Malbos, Polygraphs of finite derivation type, *Math. Structures Comput. Sci.*, vol. 28, no. 2 (2018), pp. 155–201 (cit. on pp. 33, 203).

[90] D. Hill and W. Wang, Categorification of quantum Kac-Moody superalgebras, *Trans. Amer. Math. Soc.*, vol. 367, no. 2 (2015), pp. 1183–1216 (cit. on pp. 26, 107).

[91] P. Hu, D. Kriz, and I. Kriz, "Field theories, stable homotopy theory, and Khovanov homology", *Topology Proc.*, vol. 48 (2016), pp. 327–360 (cit. on p. 46).

[92] G. Huet, "Confluent reductions: abstract properties and applications to term rewriting systems", in: *18th Annual Symposium on Foundations of Computer Science (Providence, R.I., 1977)*, IEEE, Long Beach, CA, 1977, pp. 30–45 (cit. on p. 36).

[93] M. Jimbo, A q-difference analogue of $U(\mathfrak{g})$ and the Yang-Baxter equation, *Lett. Math. Phys.*, vol. 10, no. 1 (1985), pp. 63–69 (cit. on p. 4).

[94] V. F. R. Jones, A polynomial invariant for knots via von Neumann algebras, *Bull. Amer. Math. Soc. (N.S.)*, vol. 12, no. 1 (1985), pp. 103–111 (cit. on p. 2).

[95] J.-P. Jouannaud and H. Kirchner, Completion of a set of rules modulo a set of equations, *SIAM J. Comput.*, vol. 15, no. 4 (1986), pp. 1155–1194 (cit. on p. 36).

[96] J.-P. Jouannaud and J. Li, "Church-Rosser properties of normal rewriting", in: *Computer Science Logic 2012*, vol. 16, LIPIcs. Leibniz Int. Proc. Inform. Schloss Dagstuhl. Leibniz-Zent. Inform., Wadern, 2012, pp. 350–365 (cit. on p. 36).

[97] A. Joyal and R. Street, Braided tensor categories, *Adv. Math.*, vol. 102, no. 1 (1993), pp. 20–78 (cit. on p. 235).

[98] A. Joyal and R. Street, The geometry of tensor calculus. I, *Adv. Math.*, vol. 88, no. 1 (1991), pp. 55–112 (cit. on p. 7).

[99] S.-J. Kang and M. Kashiwara, Categorification of highest weight modules via Khovanov-Lauda-Rouquier algebras, *Invent. Math.*, vol. 190, no. 3 (2012), pp. 699–742 (cit. on p. 39).

[100] S.-J. Kang, M. Kashiwara, and S.-j. Oh, Supercategorification of quantum Kac-Moody algebras, *Adv. Math.*, vol. 242 (2013), pp. 116–162 (cit. on pp. 19, 20, 107).

[101] S.-J. Kang, M. Kashiwara, and S.-j. Oh, Supercategorification of quantum Kac–Moody algebras II, *Adv. Math.*, vol. 265 (2014), pp. 169–240 (cit. on pp. 19, 20, 107).

[102] S.-J. Kang, M. Kashiwara, and S. Tsuchioka, Quiver Hecke superalgebras, *J. Reine Angew. Math.*, vol. 711 (2016), pp. 1–54 (cit. on pp. 19, 20, 107).

[103] D. Kazhdan and G. Lusztig, Representations of Coxeter groups and Hecke algebras, *Invent. Math.*, vol. 53, no. 2 (1979), pp. 165–184 (cit. on p. 11).

[104] M. Khovanov, Crossingless matchings and the cohomology of (n,n) Springer varieties, *Commun. Contemp. Math.*, vol. 6, no. 4 (2004), pp. 561–577 (cit. on p. 21).

[105] M. Khovanov, sl(3) link homology, *Algebr. Geom. Topol.*, vol. 4, no. 2 (2004), pp. 1045–1081 (cit. on pp. 12, 15).

[106] M. Khovanov and A. Lauda, A categorification of quantum $sl(n)$, *Quantum Topol.* (2010), pp. 1–92 (cit. on pp. 10, 39).

[107] M. Khovanov and A. Lauda, A diagrammatic approach to categorification of quantum groups I, *Represent. Theory*, vol. 13, no. 14 (2009), pp. 309–347 (cit. on pp. 24, 73, 107, 115).

[108] M. Khovanov and R. Lipshitz, Categorical lifting of the Jones polynomial: a survey, *Bull. Amer. Math. Soc. (N.S.)*, vol. 60, no. 4 (2023), pp. 483–506 (cit. on pp. 12, 21).

[109] M. Khovanov, K. Putyra, and P. Vaz, Odd two-variable Soergel bimodules and Rouquier complexes, in: *Algebraic and Topological Aspects of Representation Theory*, vol. 791, Contemp. Math. Amer. Math. Soc., [Providence], RI, 2024, pp. 205–227 (cit. on p. 19).

[110] M. Khovanov and L. Rozansky, Matrix factorizations and link homology, *Fundamenta Mathematicae*, vol. 199 (2008), pp. 1–91 (cit. on p. 12).

[111] P. B. Kronheimer and T. S. Mrowka, Khovanov homology is an unknot-detector, *Publ. Math. Inst. Hautes Études Sci.*, no. 113 (2011), pp. 97–208 (cit. on p. 12).

[112] Y. Lafont, Towards an algebraic theory of Boolean circuits, *J. Pure Appl. Algebra*, vol. 184, no. 2-3 (2003), pp. 257–310 (cit. on p. 33).

[113] Y. Lafont and F. Métayer, Polygraphic resolutions and homology of monoids, *J. Pure Appl. Algebra*, vol. 213, no. 6 (2009), pp. 947–968 (cit. on p. 34).

[114] Y. Lafont, F. Métayer, and K. Worytkiewicz, A folk model structure on omega-cat, *Adv. Math.*, vol. 224, no. 3 (2010), pp. 1183–1231 (cit. on p. 34).

[115] A. D. Lauda, H. Queffelec, and D. E. V. Rose, Khovanov homology is a skew Howe 2–representation of categorified quantum $\mathfrak{sl}_m$, *Algebr. Geom. Topol.*, vol. 15, no. 5 (2015), pp. 2517–2608 (cit. on pp. 12, 15, 18, 25, 78, 107, 117).

[116] A. D. Lauda, A categorification of quantum $sl(2)$, *Adv. Math.*, vol. 225, no. 6 (2010), pp. 3327–3424 (cit. on p. 10).

[117] A. D. Lauda and H. M. Russell, Oddification of the cohomology of type A Springer varieties, *Int. Math. Res. Not. IMRN*, no. 17 (2014), pp. 4822–4854 (cit. on pp. 19, 21).

[118] A. D. Lauda and J. Sussan, An invitation to categorification, *Notices Amer. Math. Soc.*, vol. 69, no. 1 (2022), pp. 11–21 (cit. on p. 8).

[119] T. Lawson, R. Lipshitz, and S. Sarkar, Khovanov homotopy type, Burnside category and products, *Geom. Topol.*, vol. 24, no. 2 (2020), pp. 623–745 (cit. on p. 46).

[120] R. Lipshitz, P. S. Ozsvath, and D. P. Thurston, Bordered Heegaard Floer homology, *Mem. Amer. Math. Soc.*, vol. 254, no. 1216 (2018), pp. viii+279 (cit. on p. 46).

[121] R. Lipshitz and S. Sarkar, A Khovanov stable homotopy type, *J. Amer. Math. Soc.*, vol. 27, no. 4 (2014), pp. 983–1042 (cit. on pp. 46, 92, 232).

[122] R. Lipshitz and S. Sarkar, A Steenrod square on Khovanov homology, *J. Topol.*, vol. 7, no. 3 (2014), pp. 817–848 (cit. on p. 232).

[123] Y. L. Liu, Braiding on complex oriented Soergel bimodules, preprint 2024, arXiv: 2407.04891 [math] (cit. on p. 46).

[124] Y. L. Liu, A. Mazel-Gee, D. Reutter, C. Stroppel, and P. Wedrich, A braided monoidal $(\infty, 2)$-category of Soergel bimodules, preprint 2024, arXiv: 2401.02956 [math] (cit. on p. 46).

[125] M. Mackaay, sl(3)-Foams and the Khovanov-Lauda categorification of quantum sl(k), preprint 2009, arXiv: 0905.2059 [math] (cit. on p. 107).

[126] M. Mackaay, M. Stošić, and P. Vaz, $\mathfrak{sl}(N)$–link homology ($N \geq 4$) using foams and the Kapustin–Li formula, *Geom. Topol.*, vol. 13, no. 2 (2009), pp. 1075–1128 (cit. on pp. 12, 15).

[127] M. Mackaay, M. Stošić, and P. Vaz, A diagrammatic categorification of the $q$-Schur algebra, *Quantum Topol.*, vol. 4, no. 1 (2013), pp. 1–75 (cit. on pp. 18, 24, 107, 116).

[128] M. Mackaay and B. Webster, Categorified skew Howe duality and comparison of knot homologies, *Adv. Math.*, vol. 330 (2018), pp. 876–945 (cit. on p. 12).

[129] A. Manion and R. Rouquier, Higher representations and cornered Heegaard Floer homology, preprint 2020, arXiv: 2009.09627 [math] (cit. on pp. 16, 46).

[130] C. Marché, "Normalized rewriting: a unified view of Knuth-Bendix completion and Gröbner bases computation", in: *Symbolic Rewriting Techniques (Ascona, 1995)*, vol. 15, Progr. Comput. Sci. Appl. Logic, Birkhäuser, Basel, 1998, pp. 193–208 (cit. on p. 36).

[131] A. Markoff, "On the impossibility of certain algorithms in the theory of associative systems", *C. R. (Doklady) Acad. Sci. URSS (N.S.)*, vol. 55 (1947), pp. 583–586 (cit. on p. 29).

[132] V. Mazorchuk and C. Stroppel, A combinatorial approach to functorial quantum $\mathfrak{sl}_k$ knot invariants, *Amer. J. Math.*, vol. 131, no. 6 (2009), pp. 1679–1713 (cit. on p. 12).

[133] V. Mazorchuk and V. Miemietz, Additive versus abelian 2-representations of fiat 2-categories, *Mosc. Math. J.*, vol. 14, no. 3 (2014), pp. 595–615, 642 (cit. on p. 11).

[134] V. Mazorchuk and V. Miemietz, Cell 2-representations of finitary 2-categories, *Compos. Math.*, vol. 147, no. 5 (2011), pp. 1519–1545 (cit. on p. 11).

[135] V. Mazorchuk and V. Miemietz, Endomorphisms of cell 2-representations, *Int. Math. Res. Not. IMRN*, no. 24 (2016), pp. 7471–7498 (cit. on p. 11).

[136] V. Mazorchuk and V. Miemietz, Isotypic faithful 2-representations of $\mathcal{J}$-simple fiat 2-categories, *Math. Z.*, vol. 282, no. 1-2 (2016), pp. 411–434 (cit. on p. 11).

[137] V. Mazorchuk and V. Miemietz, Morita theory for finitary 2-categories, *Quantum Topol.*, vol. 7, no. 1 (2016), pp. 1–28 (cit. on p. 11).

[138] V. Mazorchuk and V. Miemietz, Transitive 2-representations of finitary 2-categories, *Trans. Amer. Math. Soc.*, vol. 368, no. 11 (2016), pp. 7623–7644 (cit. on p. 11).

[139] F. Métayer, "Resolutions by polygraphs", *Theory Appl. Categ.*, vol. 11 (2003), No. 7, 148–184 (cit. on p. 34).

[140] V. Mikhaylov and E. Witten, Branes and supergroups, *Comm. Math. Phys.*, vol. 340, no. 2 (2015), pp. 699–832 (cit. on pp. 19, 26).

[141] S. Mimram, "Computing critical pairs in 2-dimensional rewriting systems", in: *RTA 2010: Proceedings of the 21st International Conference on Rewriting Techniques and Applications*, vol. 6, LIPIcs. Leibniz Int. Proc. Inform. Schloss Dagstuhl. Leibniz-Zent. Inform., Wadern, 2010, pp. 227–241 (cit. on p. 34).

[142] S. Mimram, Towards 3-dimensional rewriting theory, *Log. Methods Comput. Sci.*, vol. 10, no. 2 (2014), 2:1, 47 (cit. on p. 33).

[143] S. Morrison, K. Walker, and P. Wedrich, Invariants of 4-manifolds from Khovanov-Rozansky link homology, preprint 2021, arXiv: 1907.12194 [math] (cit. on p. 12).

[144] G. Naisse and K. Putyra, Odd Khovanov homology for tangles, preprint 2020, arXiv: 2003.14290 (cit. on pp. 20, 21, 24, 25, 77, 107, 117).

[145] G. Naisse and P. Vaz, Odd Khovanov's arc algebra, *Fund. Math.*, vol. 241, no. 2 (2018), pp. 143–178 (cit. on pp. 19, 20, 25).

[146] T. Ohtsuki, "Quantum invariants: a study of knots, 3-manifolds, and their sets", K & E Series on Knots and Everything v. 29, Singapore ; River Edge, NJ: World Scientific, 2002, 489 pp. (cit. on p. 84).

[147] V. Ostrik, Module categories, weak Hopf algebras and modular invariants, *Transform. Groups*, vol. 8, no. 2 (2003), pp. 177–206 (cit. on p. 11).

[148] P. Ozsváth and Z. Szabó, On the Heegaard Floer homology of branched double-covers, *Adv. Math.*, vol. 194, no. 1 (2005), pp. 1–33 (cit. on p. 20).

[149] P. S. Ozsváth, J. Rasmussen, and Z. Szabó, Odd Khovanov homology, *Algebr. Geom. Topol.*, vol. 13, no. 3 (2013), pp. 1465–1488 (cit. on pp. 18, 19, 88, 91, 92, 129).

[150] D. J. Peifer, "Reinforcement Learning in Buchberger's Algorithm", ProQuest LLC, Ann Arbor, MI, 2021, 170 pp. (cit. on p. 45).

[151] R. Penrose, "Applications of negative dimensional tensors", in: *Combinatorial Mathematics and Its Applications (Proc. Conf., Oxford, 1969)*, Academic Press, London-New York, 1971, pp. 221–244 (cit. on p. 7).

[152] G. E. Peterson and M. E. Stickel, Complete sets of reductions for some equational theories, *J. Assoc. Comput. Mach.*, vol. 28, no. 2 (1981), pp. 233–264 (cit. on p. 36).

[153] L. Piccirillo, The Conway knot is not slice, *Ann. of Math. (2)*, vol. 191, no. 2 (2020), pp. 581–591 (cit. on p. 12).

[154] E. L. Post, Recursive unsolvability of a problem of Thue, *J. Symbolic Logic*, vol. 12 (1947), pp. 1–11 (cit. on p. 29).

[155] K. K. Putyra, A 2-category of chronological cobordisms and odd Khovanov homology, in: *Knots in Poland III. Part III*, vol. 103, Banach Center Publ. Polish Acad. Sci. Inst. Math., Warsaw, 2014, pp. 291–355 (cit. on pp. 20, 21, 24, 25, 47, 59, 61, 77, 88, 91, 92).

[156] H. Queffelec and D. E. V. Rose, The $\mathfrak{sl}_n$ foam 2-category: A combinatorial formulation of Khovanov–Rozansky homology via categorical skew Howe duality, *Adv. Math.*, vol. 302 (2016), pp. 1251–1339 (cit. on pp. 12, 18, 47, 56, 61, 122).

[157] H. Queffelec, D. E. V. Rose, and A. Sartori, Annular Evaluation and Link Homology, 2018, arXiv: 1802.04131 [math] (cit. on p. 18).

[158] J. Rasmussen, Khovanov homology and the slice genus, *Invent. math.*, vol. 182, no. 2 (2010), pp. 419–447 (cit. on p. 12).

[159] Q. Ren and M. Willis, Khovanov homology and exotic 4-manifolds, preprint 2024, arXiv: 2402.10452 [math] (cit. on p. 12).

[160] N. Reshetikhin and V. G. Turaev, Invariants of 3-manifolds via link polynomials and quantum groups, *Invent. Math.*, vol. 103, no. 3 (1991), pp. 547–597 (cit. on pp. 2, 6).

[161] N. Yu. Reshetikhin and V. G. Turaev, Ribbon graphs and their invariants derived from quantum groups, *Comm. Math. Phys.*, vol. 127, no. 1 (1990), pp. 1–26 (cit. on p. 2).

[162] D. Reutter and J. Vicary, "High-level methods for homotopy construction in associative n-categories", in: *2019 34th Annual ACM/IEEE Symposium on Logic in Computer Science (LICS)*, IEEE, [Piscataway], NJ, 2019, [13 pp.] (Cit. on pp. 38, 45).

[163] L.-H. Robert, A new way to evaluate MOY graphs, preprint 2015, arXiv: 1512.02370 [math] (cit. on p. 69).

[164] L.-H. Robert and E. Wagner, A closed formula for the evaluation of $\mathfrak{sl}_N$-foams, 2018, arXiv: 1702.04140 [math] (cit. on p. 27).

[165] L.-H. Robert and E. Wagner, Symmetric Khovanov-Rozansky link homologies, *J. Éc. polytech. Math.*, vol. 7 (2020), pp. 573–651 (cit. on p. 18).

[166] R. Rouquier, 2-Kac–Moody algebras, preprint 2008, arXiv: 0812.5023 (cit. on pp. 10, 24, 107, 115).

[167] T. Sano, Fixing the functoriality of Khovanov homology: a simple approach, *J. Knot Theory Ramifications*, vol. 30, no. 11 (2021), Paper No. 2150074, 12 (cit. on pp. 12, 15).

[168] S. Sarkar, C. Scaduto, and M. Stoffregen, An odd Khovanov homotopy type, *Adv. Math.*, vol. 367 (2020), p. 107112 (cit. on p. 46).

[169] L. Schelstraete, Supercategorification and Khovanov-like tangle invariants, master's thesis, UCLouvain, 2020 (cit. on pp. 23, 125).

[170] L. Schelstraete and P. Vaz, Odd Khovanov homology and higher representation theory, preprint 2023, arXiv: 2311.14394 [math] (cit. on p. 21).

[171] C. Seed, Computations of the Lipshitz-Sarkar Steenrod Square on Khovanov Homology, preprint 2012, arXiv: 1210.1882 [math] (cit. on p. 232).

[172] P. Seidel and I. Smith, A link invariant from the symplectic geometry of nilpotent slices, *Duke Math. J.*, vol. 134, no. 3 (2006), pp. 453–514 (cit. on p. 12).

[173] P. Selinger, A survey of graphical languages for monoidal categories, in: *New Structures for Physics*, vol. 813, Lecture Notes in Phys. Springer, Heidelberg, 2011, pp. 289–355 (cit. on pp. 7, 56).

[174] A. I. Shirshov, Some algorithmic problems for Lie algebras, in: *Selected Works of A.I. Shirshov*, ed. by L. Bokut, I. Shestakov, V. Latyshev, and E. Zelmanov, Basel: Birkhäuser, 2009, pp. 125–130 (cit. on p. 32).

[175] A. N. Shumakovitch, Patterns in odd Khovanov homology, *J. Knot Theory Ramifications*, vol. 20, no. 1 (2011), pp. 203–222 (cit. on p. 18).

[176] W. Soergel, The combinatorics of Harish-Chandra bimodules, *J. Reine Angew. Math.*, vol. 429 (1992), pp. 49–74 (cit. on p. 11).

[177] R. Street, Categorical structures, in: *Handbook of Algebra, Vol. 1*, vol. 1, Handb. Algebr. Elsevier/North-Holland, Amsterdam, 1996, pp. 529–577 (cit. on pp. 37, 141).

[178] R. Street, Limits indexed by category-valued 2-functors, *Journal of Pure and Applied Algebra*, vol. 8, no. 2 (1976), pp. 149–181 (cit. on p. 33).

[179] C. Stroppel, "Categorification: tangle invariants and TQFTs", in: *ICM—International Congress of Mathematicians. Vol. 2. Plenary Lectures*, EMS Press, Berlin, 2023, pp. 1312–1353 (cit. on p. 16).

[180] C. Stroppel and B. Webster, 2-block Springer fibers: convolution algebras and coherent sheaves, *Comment. Math. Helv.*, vol. 87, no. 2 (2012), pp. 477–520 (cit. on p. 21).

[181] R. Usher, Fermionic 6j-symbols in superfusion categories, *J. Algebra*, vol. 503 (2018), pp. 453–473 (cit. on p. 19).

[182] V. van Oostrom, Confluence by decreasing diagrams, *Theoret. Comput. Sci.*, vol. 126, no. 2 (1994), pp. 259–280 (cit. on p. 170).

[183] P. Vaz, Not even Khovanov homology, *Pacific J. Math.*, vol. 308, no. 1 (2020), pp. 223–256 (cit. on pp. 24, 25, 77, 107, 122).

[184] P. Viry, Rewriting modulo a rewrite system, Pisa, IT: Università di Pisa, 1995 (cit. on p. 36).

[185] P. Vogel, Functoriality of Khovanov homology, *J. Knot Theory Ramifications*, vol. 29, no. 04 (2020), p. 2050020 (cit. on p. 15).

[186] M. Weber, "Free products of higher operad algebras", *Theory Appl. Categ.*, vol. 28 (2013), No. 2, 24–65 (cit. on p. 141).

[187] B. Webster, Knot invariants and higher representation theory, *Mem. Amer. Math. Soc.*, vol. 250, no. 1191 (2017), pp. v+141 (cit. on pp. 12, 16, 39).

[188] B. Webster, Unfurling Khovanov-Lauda-Rouquier algebras, preprint 2018, arXiv: 1603.06311 [math] (cit. on p. 39).

[189] E. Witten, Khovanov homology and gauge theory, in: *Proceedings of the Freedman Fest*, vol. 18, Geom. Topol. Monogr. Geom. Topol. Publ., Coventry, 2012, pp. 291–308 (cit. on p. 12).

[190] E. Witten, Quantum field theory and the Jones polynomial, *Comm. Math. Phys.*, vol. 121, no. 3 (1989), pp. 351–399 (cit. on p. 6).